\documentclass[10pt]{amsart}

\usepackage{amscd,epsfig,caption,mathrsfs,pgf,rotating,stmaryrd,tikz,MnSymbol,stmaryrd}
\usepackage{calc}
\usepackage[all]{xy} 
\usetikzlibrary{decorations.pathreplacing}
\usepackage[unicode=true,pdfborder={0 0 0},bookmarksdepth=-1
]{hyperref}

\usepackage{stackengine,scalerel}
\stackMath
\newcommand\tsup[2][2]{\ThisStyle{%
 \def\useanchorwidth{T}%
  \ifnum#1>1%
    \stackon[-5.1\LMpt]{\SavedStyle\tsup[\numexpr#1-1\relax]{#2}}{\SavedStyle\mathchar`~}%
  \else%
    \stackon[-4.5\LMpt]{\SavedStyle#2}{\SavedStyle\mathchar`~}%
  \fi%
}}

\makeatletter
  \def\@wrindex#1{%
    \protected@write\@indexfile{}%
      {\string\indexentry{#1}{ \S\thesubsection (p.\thepage)}}
    \endgroup
  \@esphack}

\makeatother

\makeindex

\newcommand{\DR}{\mathrm{DR}}
\newcommand{\incl}[1][r]
  {\ar@<-0.2pc>@{^(-}[#1] \ar@<+0.2pc>@{-}[#1]}
 
\newcommand{\Add}[1]{\textcolor{black}{#1}} 

\author{Benjamin Enriquez}
\author{Hidekazu Furusho}

\address{Institut de Recherche Math\'{e}matique Avanc\'{e}e, UMR 7501, 
Universit\'{e} de Strasbourg et CNRS, 7 rue Ren\'{e} Descartes, 67000 Strasbourg, France}
\email{enriquez@math.unistra.fr}

\address{Graduate School of Mathematics, Nagoya University, 
Furo-cho, Chikusa-ku, Nagoya, 464-8602, Japan}
\email{furusho@math.nagoya-u.ac.jp}

\date{February 12, 2026} 

\newtheorem{thm}{Theorem}[section]
\newtheorem{lem}[thm]{Lemma}

\newtheorem{cor}[thm]{Corollary}
\newtheorem{prop}[thm]{Proposition}

{\theoremstyle{definition} \newtheorem{rem}[thm]{Remark}}
{\theoremstyle{definition} \newtheorem{notation}[thm]{Notation}}
{\theoremstyle{definition} \newtheorem{defn}[thm]{Definition}}

{\theoremstyle{remark} }

\numberwithin{equation}{subsection}

\begin{document}

\title[Double shuffle Lie algebra and special derivations]{Double shuffle Lie algebra and special derivations}

\begin{abstract}
Racinet's double shuffle Lie algebra $\mathfrak{dmr}_0$ is a Lie subalgebra of the Lie algebra $\mathfrak{tder}$ of 
tangential derivations of the free Lie algebra with generators $x_0,x_1$, i.e. of derivations such that $x_1\mapsto 0$ and 
$x_0\mapsto [a,x_0]$ for some element $a$. We prove: (a) $\mathfrak{dmr}_0$ is contained in the Lie subalgebra 
$\mathfrak{sder}$ of $\mathfrak{tder}$ of special derivations, i.e. satisfying the additional condition that 
$x_\infty\mapsto [b,x_\infty]$ for some element $b$, where $x_\infty:=x_1-x_0$; (b)  $\mathfrak{dmr}_0$ is stable 
under the involution of $\mathfrak{sder}$ induced by the exchange of $x_0$ and $x_\infty$. 
The first statement 
(a) says that any element of $\mathfrak{dmr}_0$ 
satisfies the "senary relation" (a fact announced without proof 
by Ecalle in 2011), and implies the inclusion $\mathfrak{dmr}_0\subset \mathfrak{krv}_2$ (which was proved by Schneps in 
2012 only conditionally to the truth of (a)). We also derive the analogues of statements (a) and 
(b) respective to Racinet's ``double shuffle schemes'' $\mathsf{DMR}_\mu(\mathbf k)$
and to the Betti double shuffle group $\mathsf{DMR}^{\mathrm{B}}(\mathbf k)$ introduced in our
earlier work. 
\end{abstract}
\maketitle
\setcounter{tocdepth}{2}
{\footnotesize \tableofcontents}

\section*{Introduction}

\subsection{Context and motivation}\label{sect:context:motivation}

The multiple zeta values (MZVs) are the real numbers defined by 
$$
\zeta(n_1,\ldots,n_s):=\sum_{0<k_1<\ldots<k_s}\frac{1}{
k_1^{n_1}\cdots k_s^{n_s}} 
$$
for $n_1\geq1,\ldots,n_{s-1}\geq1,n_s\geq2$. They arise naturally in several areas of  mathematics
and theoretical physics. It is natural to study the algebraic relations satisfied by these numbers, 
especially in connection with the theory of motives (\cite{De,DeG}).

Three main families of such relations are known: 

(a) the relations arising from associator theory (\cite{Dr,LM}), which are equivalent to the confluence relations 
(\cite{HS,Fu:confl}), 

(b) the regularized double shuffle relations (\cite{IKZ,R,Ec1}),

(c) the relations from Kashiwara-Vergne theory (\cite{AT,AET}). 

Remarkably, the schemes arising from all three systems are torsors under algebraic groups over $\mathbb Q$, which 
are semidirect products by $\mathbb G_m$ of prounipotent algebraic groups corresponding to graded Lie algebras, 
which are denoted $\mathfrak{grt}_1,\mathfrak{dmr}_0$ and $\mathfrak{krv}_2$ in the cases (a), (b)  and (c), 
and are all Lie subalgebras of one and the same Lie algebra $(\mathfrak G,\langle,\rangle)$ 
(see Lem. \ref{lem06:1109}), where $\mathfrak G$ is the free Lie algebra over two generators, equipped with the 
``Ihara bracket'' $\langle,\rangle$. The logical relations between the three systems
(a),(b),(c) are then equivalent to inclusion relations between the corresponding torsors, which are
in their turn equivalent to inclusion relations between the Lie algebras 
$\mathfrak{grt}_1,\mathfrak{dmr}_0$ and $\mathfrak{krv}_2$. 

The implication of relations $(\mathrm{a})\Rightarrow(\mathrm{b})$, which corresponds to the Lie algebra inclusion 
$\mathfrak{grt}_1\subset\mathfrak{dmr}_0$, was proved in \cite{Fu:assandDS}. An alternative proof 
is given in \cite{EF1,EF2} based on the ideas of \cite{DT,T}; another alternative proof can be derived from the
combination of \cite{HS} (confluence relations imply double shuffle relations) and \cite{Fu:confl} 
(associator relations are equivalent to confluence relations). 

The inclusion $\mathfrak{grt}_1\subset\mathfrak{krv}_2$ was proved in \cite{AT}, and the corresponding 
implication  $(\mathrm{a})\Rightarrow(\mathrm{c})$ in \cite{AET}. 

The relations (c) consist of two parts, the ``speciality'' and ``trace'' conditions. The scheme arising from the 
speciality condition is a torsor, and the corresponding Lie algebra, denoted $(\mathfrak G_{\mathrm{inert}},\langle,\rangle)$
(see Lem. \ref{lem06:1109}), is a Lie subalgebra of $(\mathfrak G,\langle,\rangle)$. The statement that 
relations (b) imply the speciality relations corresponds to 
the Lie algebra inclusion $\mathfrak{dmr}_0\subset (\mathfrak G_{\mathrm{inert}},\langle,\rangle)$; it 
was stated in \cite{Ec2} without proof. It was proved in \cite{Sch1} that relations (b) imply the ``trace''
conditions, so that the Lie algebra inclusion $\mathfrak{dmr}_0\subset (\mathfrak G_{\mathrm{inert}},\langle,\rangle)$
both implies the inclusion $\mathfrak{dmr}_0\subset \mathfrak{krv}_2$ and the implication 
$(\mathrm{b})\Rightarrow(\mathrm{c})$ (see also \cite{FK}).  

The purpose of this paper is threefold: 

(i) to prove the inclusion $\mathfrak{dmr}_0\subset (\mathfrak G_{\mathrm{inert}},\langle,\rangle)$
(a result which was recently established independently in \cite{Sch2,Kaw}), as
well as the stability of $\mathfrak{dmr}_0$ 
under an involution $\mathrm{Lie}\Theta$ of $\mathfrak G_{\mathrm{inert}}$ 
(see Lem. \ref{lem02:2506}(d)); 

(ii) to derive the consequences of statements (i) and (ii) concerning the scheme of double 
shuffle relations 
$\mathsf{DMR}_\mu(\mathbf k)$ from \cite{R} and the Betti group 
$\mathsf{DMR}^{\mathrm B}(\mathbf k)$ (see \cite{EF3}); 

(iii) we observe that   $\mathfrak{dmr}_0\subset (\mathfrak G_{\mathrm{inert}},\langle,\rangle)$
from (i) implies both $(\mathrm{b})\Rightarrow(\mathrm{c})$ and the inclusion 
$\mathfrak{dmr}_0\subset \mathfrak{krv}_2$, from where one derives that the sequence of 
inclusions 
$$
\mathfrak{grt}_1\subset \mathfrak{dmr}_0\subset \mathfrak{krv}_2\subset (\mathfrak G_{\mathrm{inert}},\langle,\rangle)
$$ 
holds unconditionally.

The sequel of this Introduction is organized as follows. 
In \S\S\ref{sect:0:2} and \ref{sect:0:3}, we introduce the framework necessary for the 
formulation of the results of the text. In \S\ref{sect:0:2}, we introduce the 
bitorsor $(\mathcal G\rtimes\mathbf k^\times,\mathcal G\times\mathbf k^\times,
\mathcal G^{\mathrm B}\rtimes\mathbf k^\times)$, its subbitorsor 
$(\mathcal G_{\mathrm{inert}} \rtimes \mathbf k^\times,
\sqcup_{\mu\in\mathbf k^\times}\mathcal G_{\mathrm{inert}}^\mu,
(\mathcal G^{\mathrm B} \rtimes \mathbf k^\times)_{\mathrm{inert}})$ and the involution  
$(\Theta\rtimes id,\sqcup_{\mu\in\mathbf k^\times}\Theta^\mu,\Theta^{\mathrm B}\rtimes id)$ of the latter. 
In \S\ref{sect:0:3}, we introduce the bitorsor $(\mathsf{DMR}_0(\mathbf k)\rtimes\mathbf k^\times,
\sqcup_{\mu\in\mathbf k^\times}\mathsf{DMR}_\mu(\mathbf k),
\mathsf{DMR}^{\mathrm B}(\mathbf k))$, which is known to be a subbitorsor of $(\mathcal G\rtimes\mathbf k^\times,\mathcal G\times\mathbf k^\times,
\mathcal G^{\mathrm B}\rtimes\mathbf k^\times)$, its subbitorsor 
$(\mathcal G_{\mathrm{inert}} \rtimes \mathbf k^\times,
\sqcup_{\mu\in\mathbf k^\times}\mathcal G_{\mathrm{inert}}^\mu,
(\mathcal G^{\mathrm B} \rtimes \mathbf k^\times)_{\mathrm{inert}})$ by \cite{EF3}. 
We then formulate the main results in \S\ref{subsect:summary}, namely the inclusion of $(\mathsf{DMR}_0(\mathbf k)\rtimes\mathbf k^\times,
\sqcup_{\mu\in\mathbf k^\times}\mathsf{DMR}_\mu(\mathbf k),
\mathsf{DMR}^{\mathrm B}(\mathbf k))$ in $(\mathcal G_{\mathrm{inert}} \rtimes \mathbf k^\times,
\sqcup_{\mu\in\mathbf k^\times}\mathcal G_{\mathrm{inert}}^\mu,
(\mathcal G^{\mathrm B} \rtimes \mathbf k^\times)_{\mathrm{inert}})$ and
its stability under 
$(\Theta\rtimes id,\sqcup_{\mu\in\mathbf k^\times}\Theta^\mu,\Theta^{\mathrm B}\rtimes id)$, corresponding to the above i) and ii). We also formulate precisely iii).

Throughout the paper, we denote by $\mathbf k$ a $\mathbb Q$-algebra. 

\subsection{The bitorsor $(\mathcal G_{\mathrm{inert}} \rtimes \mathbf k^\times,
\sqcup_{\mu\in\mathbf k^\times}\mathcal G_{\mathrm{inert}}^\mu,
(\mathcal G^{\mathrm B} \rtimes \mathbf k^\times)_{\mathrm{inert}})$}\label{sect:0:2}

\subsubsection{The group $(\mathcal G,\circledast)$, its subgroup 
 $\mathcal G_{\mathrm{inert}}$ and the involution $\Theta$}

For $\mathfrak g$ a complete graded Lie algebra, let $(U\mathfrak g)^\wedge$ be 
the completion of its enveloping algebra $U\mathfrak g$ for the topology 
defined by the powers of its augmentation ideal; 
this is a topological Hopf algebra. Its group of group-like elements
is $\mathrm{exp}(\mathfrak g)$, where $\mathrm{exp} : (U\mathfrak g)^\wedge_0\to 
(U\mathfrak g)^\wedge$ is the exponential map, and the source is the augmentation ideal of 
the target. The assignment $\mathfrak g\mapsto\mathrm{exp}(\mathfrak g)$ is then functorial. 

Define $\mathfrak{lie}_{\{0,1\}}$ to be\footnote{the Lie algebras $\mathfrak{lie}_{\{0,1\}},
\mathfrak{lie}_{\{0,1\}}^\wedge$ are denoted $\mathfrak f_2,\hat{\mathfrak f}_2$ 
in \cite{EF1}} the free $\mathbf k$-Lie algebra over generators $e_0,e_1$, and 
$\mathfrak{lie}_{\{0,1\}}^\wedge$ its degree completion, where $e_0,e_1$ have 
degree 1. Then $\mathfrak{lie}_{\{0,1\}}^\wedge$ is the set of primitive elements of the 
topologically free algebra $\mathbf k\langle\langle e_0,e_1\rangle\rangle$, 
equipped with the coproduct for which $e_0,e_1$ are primitive. 
We denote by $\mathfrak{lie}_{\{0,1\}}^\wedge\to\mathbf k\overline e_0\oplus
\mathbf k\overline e_1$ the abelianization morphism of $\mathfrak{lie}_{\{0,1\}}^\wedge$.

\begin{defn}
    $\mathcal G$ is the kernel of the group morphism 
    $\mathrm{exp}(\mathfrak{lie}_{\{0,1\}}^\wedge)\to 
    \mathrm{exp}(\mathbf k\overline e_0\oplus\mathbf k\overline e_1)$ 
    induced by abelianization; this is the subset of 
    $\mathrm{exp}(\mathfrak{lie}_{\{0,1\}}^\wedge)$ of elements of 
    $\mathrm{exp}(\mathrm{ker}(\mathfrak{lie}_{\{0,1\}}^\wedge\to\mathbf k\overline e_0\oplus\mathbf k\overline e_1))$.  
\end{defn}

\begin{lem}\label{def:mathcalG} (Lem. \ref{def:mathcalG:w:proof}) 
The product defined by 
$$
(g\circledast h)(e_0,e_1):=h(g(e_0,e_1)\cdot e_0\cdot g(e_0,e_1)^{-1},e_1)\cdot g(e_0,e_1). 
$$
equips $\mathcal G$ with a group structure. 
\end{lem}

\begin{lem}\label{lem:actions:ktimes:a} (Lem. \ref{lem:actions:ktimes:a:w:proof})
(a)   The action of $\mathbf k^\times$ on $\mathbf k\langle\langle e_0,e_1\rangle\rangle$,
where $\lambda\in\mathbf k^\times$ acts by taking $e_0,e_1$ to $\lambda e_0,\lambda e_1$, 
induces an action of $\mathbf k^\times$ on $\mathrm{exp}(\mathfrak{lie}_{\{0,1\}}^\wedge)$, 
which restricts to an action on $\mathcal G$. 

(b) This action is compatible with the product $\circledast$, hence is an action by 
automorphisms of the group $(\mathcal G,\circledast)$. 
\end{lem}

Set $e_\infty:=-e_0-e_1$. 
\begin{lem}\label{lem02:2506:wo:proof} (Lem. \ref{lem02:2506})
    (a) If $g\in \mathcal G$ is such that there exists $h\in \mathcal G$ such that 
$\mathrm{Ad}_g(e_0)+e_1+\mathrm{Ad}_h(e_\infty)=0$ (equality in $\mathfrak{lie}_{\{0,1\}}^\wedge$), then 
$h$ is unique; it will be denoted $h_g$. 

(b) The subset $\mathcal G_{\mathrm{inert}}\subset  \mathcal G$ of all elements 
$g$ as in (a) is a subgroup of 
$(\mathcal G,\circledast)$.

(c) There is a unique automorphism $s_{(0,\infty)}$ of 
$\mathfrak{lie}_{\{0,1\}}^\wedge$, such that 
$e_1\mapsto e_1$ and $e_0\leftrightarrow e_\infty$; it is an involution. 

(d) The map $\Theta : g\mapsto s_{(0,\infty)}(h_g)$ is an involutive automorphism of $(\mathcal G_{\mathrm{inert}},\circledast)$.
\end{lem}

\begin{lem}\label{lem:actions:ktimes:b}
 (a) (Lem. \ref{lem:15:14}(a))
 The action of $\mathbf k^\times$ on $(\mathcal G,\circledast)$ from Lem. \ref{lem:actions:ktimes:a}
 preserves the subgroup $\mathcal G_{\mathrm{inert}}$, hence is an action by automorphisms of the 
 group $(\mathcal G_{\mathrm{inert}},\circledast)$. 

(b) (Lem. \ref{lem:15:15}(a))
The latter action commutes with the involution $\Theta$.     
\end{lem}

\begin{defn}
   $\mathcal G\rtimes\mathbf k^\times$ (resp. 
   $\mathcal G_{\mathrm{inert}}\rtimes\mathbf k^\times$) is the semidirect 
   product\footnote{
Recall that the semidirect product $A\rtimes G$ arising from a pair of groups $A,G$ and an 
action $*$ of $G$ on $A$, is the group explicitly 
defined as the set $A\times G$, equipped with $(a,g)\cdot (a',g')=(a\cdot (g*a'),g\cdot g')$. 
It is equipped with group morphisms $i,j$ from 
$A,G$ to it, satisfying the identity $j(g)i(a)j(g)^{-1}=i(g*a)$; namely, $i(a)=(a,1)$ and 
$j(g)=(1,g)$; the elements
$i(a)$ and $j(g)$ by sometimes be simply denoted $a$ and $g$. If $\Gamma$ is a group, there 
is a bijection between the set of group morphisms $A\rtimes G\to\Gamma$ and the set of pairs 
of group morphisms 
$\alpha : A\to \Gamma$, $\gamma : G\to \Gamma$, satisfying the identity $\phi(a*g)
=\alpha(a)\phi(g)\alpha(a)^{-1}$.}  
   of the group $(\mathcal G,\circledast)$ (resp. $(\mathcal G_{\mathrm{inert}},\circledast)$)
   with the action of $\mathbf k^\times$ from Lem. \ref{lem:actions:ktimes:a}; both group 
   structures are given by 
$$
(g(e_0,e_1),\lambda)\circledast (h(e_0,e_1),\mu)
:=(h(g(e_0,e_1)\lambda e_0 g(e_0,e_1)^{-1},\lambda e_1)  \cdot g(e_0,e_1),\lambda\mu). 
$$
Lem. \ref{lem:actions:ktimes:b} then gives rise to a group inclusion
   $$
   \mathcal G_{\mathrm{inert}}\rtimes\mathbf k^\times \subset 
\mathcal G\rtimes\mathbf k^\times 
   $$ 
as well as to an involution $\Theta\rtimes id$ of its source.  
\end{defn}

\begin{lem}\label{lem06:1109}
 (a) The Lie algebra of the group functor $\mathbf k\mapsto (\mathrm{exp}(\mathfrak{lie}_{\{0,1\}}^\wedge),\circledast)$ is 
$(\mathfrak{lie}_{\{0,1\}}^\wedge,\langle,\rangle)$, with 
$$
\langle a,b\rangle:=[a,b]+d_a(b)-d_b(a), \quad \textrm{where} \quad d_a : e_0\mapsto 0,e_1\mapsto[e_1,a].  
$$

(b) The Lie algebra of the group functor $\mathbf k\mapsto (\mathcal G,\circledast)$ is 
$(\mathfrak G,\langle,\rangle)$, where 
$\mathfrak G\subset \mathfrak{lie}_{\{0,1\}}^\wedge$ is the subspace of 
elements with vanishing degree $1$ part 
 (the degree being defined by the fact that both $e_0,e_1$ are of degree $1$). 

(c) The Lie algebra of the group functor $\mathbf k\mapsto (\mathcal G_{\mathrm{inert}},\circledast)$ is then 
$(\mathfrak G_{\mathrm{inert}},\langle,\rangle)$, where $\mathfrak G_{\mathrm{inert}}\subset\mathfrak G$ is defined as 
$\{a\in\mathfrak G|\exists b\in\mathfrak G,[a,e_0]+[b,e_\infty]=0\}$; for $a\in \mathfrak G_{\mathrm{inert}}$, 
the element $b$ such that $[a,e_0]+[b,e_\infty]=0$ is unique and will be denoted $b_a$. 

(d) The Lie algebra involution $\mathrm{Lie}(\Theta)$ of $(\mathcal G_{\mathrm{inert}},\circledast)$
induced by $\Theta$ is given by $a\mapsto s_{(0,\infty)}(b_a)$. 
\end{lem}

\begin{proof} These results are obtained by the usual linearization procedure.  
\end{proof}

\subsubsection{The Betti group $(\mathcal G^{{\mathrm B}},\circledast)$, the group 
$((\mathcal G^{{\mathrm B}}\rtimes\mathbf k^\times)_{\mathrm{inert}},\circledast)$ 
and its involution $\Theta^{\mathrm B}$} 


For $\Gamma$ a discrete group, let $(\mathbf k\Gamma)^\wedge$ be 
the completion of its group algebra $\mathbf k\Gamma$ for the topology 
defined by the powers of its augmentation ideal; 
this is is a topological Hopf algebra. 
Denote by $\Gamma(\mathbf k)$ its group of group-like elements. The assignment 
$\Gamma\mapsto \Gamma(\mathbf k)$ is then functorial. 


\begin{defn}
    Define $\mathcal G^{\mathrm B}:=\mathrm{ker}(F_2(\mathbf k)\to\mathbb Z^2(\mathbf k)
    =\mathbf k^2)$, 
    where the morphism $F_2\to\mathbb Z^2$ is the abelianization morphism of the free group $F_2$ 
    with generators $X_0,X_1$.
\end{defn}

\begin{lem}\label{lem:mathcalG:betti} (\cite{EF3}, §2.1.3)
A group structure $\circledast$ is defined on $\mathcal G^{\mathrm B}$ by 
$$
g(X_0,X_1)\circledast h(X_0,X_1)
:=h(g(X_0,X_1)X_0 g(X_0,X_1)^{-1},X_1)  \cdot g(X_0,X_1).
$$
The group $\mathbf k^\times$ acts on $(\mathcal G^{\mathrm B},\circledast)$ by 
$\lambda\bullet g(X_0,X_1):=g(X_0^\lambda,X_1^\lambda)$. The resulting semidirect 
product group $\mathcal G^{\mathrm B}\rtimes \mathbf k^\times$ is the set 
$\mathcal G^{\mathrm B}\times \mathbf k^\times$, equipped with the product 
$$
(g(X_0,X_1),\lambda)\circledast (h(X_0,X_1),\mu)
:=(h(g(X_0,X_1)X_0^\lambda g(X_0,X_1)^{-1},X_1^\lambda)  \cdot g(X_0,X_1),\lambda\mu). 
$$
\end{lem}

\begin{defn}
    $(\mathcal G^{\mathrm B}\rtimes\mathbf k^\times)_{\mathrm{inert}}$ is the subset of 
    $\mathcal G^{\mathrm B}\rtimes\mathbf k^\times$
    of all pairs $(g,\lambda)$ such that there exists 
    $h\in F_2(\mathbf k)$ such that 
\begin{equation}\label{def:dual:g:h}
    X_1^\lambda\mathrm{Ad}_{g}(X_0^\lambda)=\mathrm{Ad}_{h}((X_1X_0)^\lambda). 
\end{equation} 
\end{defn}

\begin{lem}\label{lem:semi:B:wo:proof} (see Lem. \ref{lem:semi:B})
$((\mathcal G^{\mathrm B}\rtimes\mathbf k^\times)_{\mathrm{inert}},\circledast)$
    is a subgroup of $(\mathcal G^{\mathrm B}\rtimes\mathbf k^\times,\circledast)$.  
\end{lem}

\begin{lem}\label{lem:15:12:new} (Lem. \ref{lem:15:12:new:w:proof})
(a) There is an involutive automorphism $\sigma$ of $F_2(\mathbf k)$, determined by 
$$
\sigma : X_0\mapsto X_1^{-1/2}X_\infty X_1^{1/2},\quad X_\infty\mapsto X_1^{1/2}X_0X_1^{-1/2},
\quad X_1\mapsto X_1.  
$$

(b) There is a unique involution $\Theta^{\mathrm B}$ of 
$((\mathcal G^{\mathrm B}\rtimes\mathbf k^\times)_{\mathrm{inert}},\circledast)$
such that 
$$
\forall (g,\lambda)\in (\mathcal G^{\mathrm B}\rtimes\mathbf k^\times)_{\mathrm{inert}}, \quad 
\Theta^{\mathrm B}(g,\lambda)=(X_1^{-\lambda/2} \sigma(h_g)X_1^{1/2},\lambda). 
$$
\end{lem}

\subsubsection{ The bitorsor $(\mathcal G \rtimes \mathbf k^\times,
\mathcal G \times\mathbf k^\times,
\mathcal G^{\mathrm B} \rtimes \mathbf k^\times)$}

Recall from \cite{Gi}, 
Chap. III, Def. 1.5.3 that a bitorsor is $(G,X,H)$ is the data of a nonempty set $X$, of 
groups $G$ and $H$, and of commuting left and right actions of $G$ and $H$ on $X$, which are
transitive and with trivial stabilizers. A subbitorsor $(G',X',H')$ is then the data of
subgroups $G',H'$ of $G,H$ and of a nonempty subset $X'$ of $X$, stable under the actions of $G'$ and
$H'$ and such that $(G',X',H')$ is a torsor. An involution of the bitorsor 
$(G,X,H)$ is a triple of group involutions of $G$ and $H$ and of a set involution of 
$X$, which are compatible with the actions. 


\begin{lem} \label{lem:GGG:bitorsor}
(\cite{EF3}, §2.3)
(a) The maps $((g(e_0,e_1),\lambda),(\phi(e_0,e_1),\mu))\mapsto (g(e_0,e_1),\lambda)\circledast(\phi(e_0,e_1),\mu)$ and 
$$
((\phi(e_0,e_1),\mu),(k(X_0,X_1),\nu))\mapsto 
(\phi(e_0,e_1),\mu)\bullet(k(X_0,X_1),\nu)
:=(k(\phi(e_0,e_1)e^{\mu e_0}\phi^{-1}(e_0,e_1),e^{\mu e_1})\phi(e_0,e_1),\mu\nu)
$$
respectively define a left action of 
$\mathcal G\rtimes\mathbf k^\times$ and a right action of $\mathcal G^{\mathrm B}
\rtimes\mathbf k^\times$ on the set 
$\mathcal G\times\mathbf k^\times$.  

(b) $(\mathcal G\rtimes\mathbf k^\times,\mathcal G\times\mathbf k^\times,
\mathcal G^{\mathrm B}\rtimes\mathbf k^\times)$ is a bitorsor. 
\end{lem}

\subsubsection{The 
subbitorsor $(\mathcal G_{\mathrm{inert}} \rtimes \mathbf k^\times,
\sqcup_{\mu\in\mathbf k^\times}\mathcal G_{\mathrm{inert}}^\mu,
(\mathcal G^{\mathrm B} \rtimes \mathbf k^\times)_{\mathrm{inert}})$}

\begin{defn}
    For $x,y$ free noncommutative variables, let $\mathrm{cbh}(x,y)\in\mathfrak{lie}_{x,y}^\wedge$
    be the Campbell-Baker-Hausdorff series defined by $\mathrm{cbh}(x,y):=\mathrm{log}e^xe^y$. 
    Let $\mathrm{cbh}(x,y)=\sum_{k\geq1}\mathrm{cbh}_k(x,y)$ be its degree decomposition. 
    
    For $\mu\in\mathbf k$, set $\mathrm{cbh}_\mu(x,y):=\sum_{k\geq1}\mu^{k-1}\mathrm{cbh}_k(x,y)
    \in (\mathfrak{lie}_{x,y}\otimes\mathbf k)^\wedge$. 
    Then $\mathrm{cbh}_\mu(x,y)=x+y+(\mu/2)[x,y]+\cdots$. 
    We also use the notation $x*_\mu y:=\mathrm{cbh}_\mu(x,y)$. 
\end{defn}

Note that for $\mu\in\mathbf k^\times$, one has $x*_\mu y=\mu^{-1}\mathrm{log}(e^{\mu x}e^{\mu y})$. 

\begin{defn}\label{def:G:mu:inert}
    For $\mu\in\mathbf k$, define $\mathcal G_{\mathrm{inert}}^\mu$ to be the subset of $\mathcal G$ of all 
    elements $g$ such that for some $h\in\mathrm{exp}(\mathfrak{lie}_{\{0,1\}}^\wedge)$, one has 
     $\mathrm{Ad}_ge_0*_\mu \mathrm{Ad}_he_\infty
 =e_0+e_\infty
 $ (equality in $\mathfrak{lie}_{\{0,1\}}^\wedge$). 
\end{defn}

Note that $\mathcal G_{\mathrm{inert}}^0=\mathcal G_{\mathrm{inert}}$.  

\begin{lem}\label{lem:G:mu} 
(Lem. \ref{lem:15:14}(c))
 Set $\sqcup_{\mu\in\mathbf k^\times}\mathcal G_{\mathrm{inert}}^\mu:=\{(\mu,g)\in
 \mathbf k^\times\times \mathcal G|g\in \mathcal G_{\mathrm{inert}}^\mu\}$. Then 
 $$
 (\mathcal G_{\mathrm{inert}} \rtimes \mathbf k^\times,
\sqcup_{\mu\in\mathbf k^\times}\mathcal G_{\mathrm{inert}}^\mu,
(\mathcal G^{\mathrm B} \rtimes \mathbf k^\times)_{\mathrm{inert}})
$$ 
is a subbitorsor of $(\mathcal G\rtimes\mathbf k^\times,\mathcal G\times
    \mathbf k^\times,\mathcal G^{\mathrm B}\rtimes\mathbf k^\times)$. 
\end{lem}

\subsubsection{The bitorsor involution 
$(\Theta\rtimes id,\sqcup_{\mu\in\mathbf k^\times}\Theta^\mu,\Theta^{\mathrm B}\rtimes id)$}

\begin{lem}\label{def:Theta:mu} (see Lem. \ref{def:Theta:mu:w:proofs})
Let $\mu\in\mathbf k$. 

    (a) If $g\in \mathcal G_{\mathrm{inert}}^\mu$, then there exists a unique 
    $h\in\mathrm{exp}(\mathfrak{lie}_{\{0,1\}}^\wedge)$ such that 
 $\mathrm{Ad}_ge_0*_\mu \mathrm{Ad}_he_\infty
 =e_0+e_\infty
 $ and $\mathrm{log}h\equiv (\mu/2)e_1$ mod $F^2\mathfrak{lie}_{\{0,1\}}^\wedge$; it 
    will be denoted $h_g$.

    (b) If $g\in \mathcal G_{\mathrm{inert}}^\mu$, then 
    $e^{-\mu e_1/2}(s_{(0,\infty)}(h_g))\in \mathcal G_{\mathrm{inert}}^\mu$. 

    (c) The map $\Theta^\mu : g\mapsto e^{-\mu e_1/2}(s_{(0,\infty)}(h_g))$ is an involution of the set $\mathcal G_{\mathrm{inert}}^\mu$. 
\end{lem}

\begin{lem}\label{lem:invol:bitorsor} (see Lem. \ref{lem:15:15}(b))
    The triple
$(\Theta\rtimes id,\sqcup_{\mu\in\mathbf k^\times}\Theta^\mu,\Theta^{\mathrm B}\rtimes id)$
is an involution of the bitorsor 
$(\mathcal G_{\mathrm{inert}} \rtimes \mathbf k^\times,
\sqcup_{\mu\in\mathbf k^\times}\mathcal G_{\mathrm{inert}}^\mu,
(\mathcal G^{\mathrm B} \rtimes \mathbf k^\times)_{\mathrm{inert}})$ (see Lem. \ref{lem:G:mu}). 
\end{lem}

\subsection{The double shuffle bitorsor}\label{sect:0:3}

\subsubsection{The group $\mathsf{DMR}_0(\mathbf k)$}

\begin{defn}
    For $g\in\mathcal G$, set $\Gamma_g^{-1}(-e_1):=\mathrm{exp}(\sum_{n\geq1}(g|e_0^{n-1}e_1)e_1^n/n)\in\mathbf k[[e_1]]$. 
\end{defn}


\begin{defn} \label{def:mrs:1608:BIS:BIS} (\cite{EF1},\S1.1)
(a) Set $\hat{\mathcal V}:=\mathbf k\langle\langle e_0,e_1\rangle\rangle$. 

(b) $\hat{\mathcal W}$ is the complete $\mathbb Z_{\geq 0}$-graded subalgebra of 
$\hat{\mathcal V}$ given by 
$\hat{\mathcal W}:=\mathbf k\oplus \hat{\mathcal V}e_1$. 

(c) $\hat{\mathcal M}$ is the left $\hat{\mathcal W}$-module given by 
$\hat{\mathcal M}:=\hat{\mathcal V}/\hat{\mathcal V}e_0$; $1_{\mathcal M}\in \hat{\mathcal M}$ is the class of $1$. 
\end{defn}
(The notation used in \cite{EF1} is $\mathcal V^\DR,\mathcal W^\DR,\mathcal M^\DR,\hat\Delta^{\mathcal W,\DR},1_\DR$.)

\begin{lem} (see \cite{EF1}, \S1.2)
(a) There is a unique topological $\mathbf  k$-algebra morphism $\hat\Delta^{\mathcal W} : 
\hat{\mathcal W}\to \hat{\mathcal W}^{\hat\otimes 2}$, such that  
\begin{equation}\label{pties:Delta:W}
   \textstyle \forall n\geq 1,\quad \hat\Delta^{\mathcal W}(e_0^{n-1}e_1)=e_0^{n-1}e_1\otimes1+1\otimes e_0^{n-1}e_1-\sum_{n',n''>0,n'+n''=n}e_0^{n'-1}e_1
    \otimes e_0^{n''-1}e_1.  
\end{equation}

(b) There is a unique topological $\mathbf k$-module morphism $\hat\Delta^{\mathcal M} : 
\hat{\mathcal M}\to \hat{\mathcal M}^{\hat\otimes 2}$, such that  
$$
\forall w\in\hat{\mathcal W},\quad \hat\Delta^{\mathcal M}(w\cdot 1_{\mathcal M})=\hat\Delta^{\mathcal W}(w)\cdot 1_{\mathcal M}^{\otimes2}. 
$$
\end{lem}

Then $\hat{\mathcal W},\hat{\mathcal M}$ are complete graded and
$\hat\Delta^{\mathcal W},\hat\Delta^{\mathcal M}$ have degree 0. 
Define $\mathcal G(\hat{\mathcal M}):=\{m\in \hat{\mathcal M}|\hat\Delta^{\mathcal M}(m)=m\otimes m\}$. 

\begin{defn} (see \cite{R}, \S3.2.1 and \cite{EF2}, \S2.4)
For $\mu\in\mathbf k$, define\footnote{The letters `DMR' stand for the French 
`double mélange et régularisation'.} $\mathsf{DMR}_\mu(\mathbf k)$ as the set of elements 
$\Phi\in \mathcal G$ such that 
$$
\Gamma_\Phi^{-1}(-e_1)\Phi\cdot 1_{\mathcal M}\in\mathcal G(\hat{\mathcal M}),\quad (\Phi|e_0e_1)=\mu^2/24. 
$$ 
\end{defn}

\begin{thm} \label{thm:racinet}
(see \cite{R}, Thm. I and \cite{EF2}, \S2.4)
    (a) $\mathsf{DMR}_0(\mathbf k)$ is a subgroup of $(\mathcal G,\circledast)$. 

    (b) For $\mu\in\mathbf k$, the left action of $(\mathcal G,\circledast)$ on itself restricts 
    to an action of the group $\mathsf{DMR}_0(\mathbf k)$ on the set $\mathsf{DMR}_\mu(\mathbf k)$, 
    which is then a torsor over this group. 
\end{thm}

\begin{lem} (see \cite{R}, \S3.3.8)
    The assignment $\mathbf k\mapsto \mathsf{DMR}_0(\mathbf k)$ is a prounipotent subgroup scheme, and its Lie subalgebra 
    is $\mathfrak{dmr}_0:=\{a\in\mathfrak G|(x+\sum_{n\geq1}(a|e_0^{n-1}e_1)e_1^n/n)\cdot 1_{\mathcal M}\in
    \mathcal P(\hat{\mathcal M})$ 
    and $(x|e_0e_1)=0\}\}$ where $\mathcal P(\hat{\mathcal M})
    =\{m\in \hat{\mathcal M}|\hat\Delta^{\mathcal M}(m)=m\otimes1_{\mathcal M}+1_{\mathcal M}\otimes m$, which is a Lie subalgebra
    of $(\mathfrak G,\langle,\rangle)$. 
\end{lem}
 
\subsubsection{The group $\mathsf{DMR}^{\mathrm B}(\mathbf k)$}

Define $\hat{\mathcal W}^{\mathrm B}:=\mathbf k+(\mathbf kF_2)^\wedge(X_1-1)\subset 
(\mathbf kF_2)^\wedge$; there is a unique continuous $\mathbf k$-algebra morphism 
$\hat\Delta^{\mathcal W,\mathrm B} : \hat{\mathcal W}^{\mathrm B}\to 
\hat{\mathcal W}^{\mathrm B}\hat\otimes
\hat{\mathcal W}^{\mathrm B}$, such that $X_1^{\pm1}\mapsto X_1^{\pm1}\otimes X_1^{\pm1}$ 
and $X_0^k(1-X_1)\mapsto X_0^k(1-X_1)\otimes1
+1\otimes X_0^k(1-X_1)+\sum_{i=1}^{k-1}X_0^i(1-X_1)\otimes X_0^{k-i}(1-X_1)$ for $k\in\mathbb Z$
(with $\sum_{i=1}^{k-1}f(i)$ being defined as $0$ for $k=1$ and 
as $-f(0)-f(-1)\cdots -f(k)$ for $k\leq0$). 
Let 
$\hat{\mathcal M}^{\mathrm B}:=(\mathbf kF_2)^\wedge/(\mathbf kF_2)^\wedge\cdot (X_0-1)$, 
and denote by $x\mapsto x\cdot 1_{\mathrm B}$ 
the natural projection $(\mathbf kF_2)^\wedge\to \hat{\mathcal M}^{\mathrm B}$. 
Then the map $\hat{\mathcal W}^{\mathrm B}\to \hat{\mathcal M}^{\mathrm B}$, $x\mapsto 
x\cdot 1_{\mathrm B}$ is an isomorphism. Let $\hat\Delta^{\mathcal M,\mathrm B} : 
\hat{\mathcal M}^{\mathrm B}\to 
\hat{\mathcal M}^{\mathrm B}\hat\otimes\hat{\mathcal M}^{\mathrm B}$ be the map such that 
$\hat\Delta^{\mathcal M,\mathrm B}(w\cdot 1_{\mathrm B})
=\hat\Delta^{\mathcal W,\mathrm B}(w)\cdot (1_{\mathrm B}\otimes 1_{\mathrm B})$
for any $w\in \hat{\mathcal W}^{\mathrm B}$
and let $\mathcal G(\hat{\mathcal M}^{\mathrm B})$ be set of group-like elements of 
$(\hat{\mathcal M}^{\mathrm B},\hat\Delta^{\mathcal M,\mathrm B})$. 

For $g\in F_2(\mathbf k)$, define $\Gamma_g$ 
by 
\begin{equation}\label{def:gamma:betti}
\textstyle\Gamma_g(t):=\mathrm{exp}(\sum_{n\geq 1}(-1)^{n+1}(g|(\mathrm{log}X_0)^{n-1}
\mathrm{log}X_1)t^n/n), 
\end{equation}
where the coordinates are taken with respect to the topological 
basis of $(\mathbf kF_2)^\wedge$ formed by the set of words in $\mathrm{log}X_0,\mathrm{log}X_1$. 

\begin{defn}\label{def:DMLR:betti} (see \cite{EF3}, Lem.-Def. 3.10 and Thm. 4.5)
The subset $\mathsf{DMR}^{\mathrm B}(\mathbf k)\subset \mathcal G^{\mathrm B}
\rtimes\mathbf k^\times$ is defined as the set of pairs $(g,\lambda)$ such that 

(a) $(\Gamma_g^{-1}(-\mathrm{log}X_0) \cdot g) 
\cdot 1_{\mathrm B} \in \mathcal G(\hat{\mathcal M}^{\mathrm B})$, where $\Gamma_g$ is as in 
\eqref{def:gamma:betti};

(b) $(g|\mathrm{log}X_0) = (g|\mathrm{log}X_1) = 0$, $\lambda^2 = 1 + 24(g|\mathrm{log}X_0\mathrm{log}X_1)$.
\end{defn}

\subsubsection{The bitorsor $(\mathsf{DMR}_0(\mathbf k)\rtimes\mathbf k^\times,
\sqcup_{\mu\in\mathbf k^\times}\mathsf{DMR}_\mu(\mathbf k),
\mathsf{DMR}^{\mathrm B}(\mathbf k))$}

\begin{lem}\label{lem:15:5} (see \cite{R} and \cite{EF3}, Lem. 3.11)
(a) The subgroup $\mathsf{DMR}_0(\mathbf k)$ of $\mathcal G$ is preserved by the 
action of $\mathbf k^\times$. 

(b) $(\mathsf{DMR}_0(\mathbf k)\rtimes\mathbf k^\times,
\sqcup_{\mu\in\mathbf k^\times}\mathsf{DMR}_\mu(\mathbf k),
  \mathsf{DMR}^{\mathrm B}(\mathbf k))$ is a subbitorsor of  
  $(\mathcal G\rtimes\mathbf k^\times,\mathcal G\times\mathbf k^\times,
\mathcal G^{\mathrm B}\rtimes\mathbf k^\times)$
(the inclusion being given by $\sqcup_{\mu\in\mathbf k^\times}\mathsf{DMR}_\mu(\mathbf k)
=\{(g,\mu)|\mu\in\mathbf k^\times,g\in \mathsf{DMR}_\mu(\mathbf k)\}
\subset \mathcal G\times\mathbf k^\times$).  
\end{lem}

\subsection{Summary of the results}\label{subsect:summary}

\subsubsection{Results on $\mathsf{DMR}_0(\mathbf k)$}


In Def. \ref{defn:2:7:korea} and Lem. \ref{lem:43:0401:TER}(b), we introduce a set $\mathrm{Hom}_{\mathcal C\text{-alg}}(\hat{\mathcal W},\hat V)$, 
an element $\Delta^{\mathcal W}_{r,l}$ in this set, and an action of the group $\mathcal G\times\mathbf k[[u,v]]^\times$ on the same set. 
This gives rise (see Lem. \ref{lem:3101009:BIS}) to an action of $\mathcal G$ on the quotient set 
$\mathbf k[[u,v]]^\times\backslash\mathrm{Hom}_{\mathcal C\text{-alg}}
(\hat{\mathcal W},\hat V)$, and denoting by 
$\mathbf k[[u,v]]^\times\bullet\Delta^{\mathcal W}_{r,l}$ 
the class of $\Delta^{\mathcal W}_{r,l}$ in this quotient set, to a 
subgroup\footnote{For $(X,x)$ a pointed set with action of a group $G$, we denote by $\mathrm{Stab}_G(x)$ (or sometimes $\mathrm{Stab}_G(X,x)$)
the stabilizer group of $x$, which is a subgroup of $G$.} 
$\mathrm{Stab}_{\mathcal G}(\mathbf k[[u,v]]^\times\bullet\Delta^{\mathcal W}_{r,l}
)\subset\mathcal G$. The main result of Part 1 is: 
\begin{thm} \label{thm:013}(Lem. \ref{lem:12:2209} and Thm. \ref{thm:main}) 
One has the equality 
$$
\mathsf{DMR}_0(\mathbf k)=\mathrm{Stab}_{\mathcal G}(\mathbf k[[u,v]]^\times\bullet\Delta^{\mathcal W}_{r,l}
) 
$$
of subgroups of $\mathcal G$. 
\end{thm}

In §\ref{sect:4:2012}, we introduce a set $\mathrm{Hom}_{\mathcal C\text{-alg}}(\hat{\mathcal V},M_3\hat V)$, 
an element $\rho_{\mathrm{DT}}$ in this set, and compatible actions of the groups $\mathcal G$ and $\mathrm{GL}_3\hat V$
on $\mathrm{Hom}_{\mathcal C\text{-alg}}(\hat{\mathcal V},M_3\hat V)$. 
 This gives rise to an action of $\mathcal G$ on the quotient set $\mathrm{GL}_3\hat V\backslash
\mathrm{Hom}_{\mathcal C\text{-alg}}(\hat{\mathcal V},M_3\hat V)$, and denoting by 
$\mathrm{GL}_3\hat V\bullet \rho_{\mathrm{DT}}
$ the class 
of $\rho_{\mathrm{DT}}$ in this quotient set, to a 
subgroup $\mathrm{Stab}_{\mathcal G}(\mathrm{GL}_3\hat V\bullet \rho_{\mathrm{DT}}
)\subset\mathcal G$. The 
main result of Part 2 is: 

\begin{thm} \label{thm:015}(Thm. \ref{thm:5:31:3103}) 
One has the inclusion
$$
\mathrm{Stab}_{\mathcal G}(\mathrm{GL}_3\hat V\bullet \rho_{\mathrm{DT}}
)\subset 
\mathrm{Stab}_{\mathcal G}(\mathbf k[[u,v]]^\times\bullet\Delta^{\mathcal W}_{r,l}
).
$$
of subgroups of $\mathcal G$.
\end{thm}

The main result of Part 3 is: 
\begin{thm} \label{thm:016} (Cor. \ref{cor:11:2:17apr}) 
One has the equality
$$
\mathrm{Stab}_{\mathcal G}(\mathrm{GL}_3\hat V\bullet \rho_{\mathrm{DT}}
)=
\mathrm{Stab}_{\mathcal G}(\mathbf k[[u,v]]^\times\bullet\Delta^{\mathcal W}_{r,l}
).
$$
of subgroups of $\mathcal G$.
\end{thm}

Thm. \ref{thm:016} supersedes Thm. \ref{thm:015}; however, the proof of Thm. \ref{thm:016} (Part 3) 
being considerably more involved than that of Thm. \ref{thm:015} (Part 2), it seemed to us appropriate 
to give a separate proof of this weaker result in Part 2. 

The main result of Part 4 is: 

\begin{thm} \label{thm:014}(Thms. \ref{thm:13:22:0205} and \ref{thm:13:36:17apr}) 
(a) One has the inclusion
$\mathrm{Stab}_{\mathcal G}(\mathrm{GL}_3\hat V\bullet \rho_{\mathrm{DT}}
)\subset \mathcal G_{\mathrm{inert}}$ 
of subgroups of $\mathcal G$. 

(b) The subgroup $\mathrm{Stab}_{\mathcal G}(\mathrm{GL}_3\hat V\bullet \rho_{\mathrm{DT}}
)$ is stable under the involution 
$\Theta$ of $\mathcal G_{\mathrm{inert}}$ (see Lem. \ref{lem02:2506}(d)).
\end{thm}

Let us now show how this implies one of the announced results (see i) of \S\ref{sect:context:motivation}).  

\begin{cor}\label{main:cor}
(a)  One has the inclusion 
$$
\mathsf{DMR}_0(\mathbf k)\hookrightarrow\mathcal{G}_{\mathrm{inert}} 
$$
of subgroups of $\mathcal G$. 

(b) $\mathsf{DMR}_0(\mathbf k)$ is stable under the involution $\Theta$ of $\mathcal{G}_{\mathrm{inert}}$
(see Lem. \ref{lem02:2506}(d)). 

(c) One has the inclusion 
$$
\mathfrak{dmr}_0\subset\mathfrak G_{\mathrm{inert}}, 
$$
of Lie subalgebras of $(\mathfrak G,\langle,\rangle)$ and $\mathfrak{dmr}_0$ is stable
under the involution $\mathrm{Lie}\Theta$ of $\mathfrak G_{\mathrm{inert}}$.   
 \end{cor}  

\begin{proof}
(a) follows from 
$$
\mathcal G_{\mathrm{inert}}\hookleftarrow
\mathrm{Stab}_{\mathcal G}(\mathrm{GL}_3\hat V\bullet \rho_{\mathrm{DT}})
=\mathrm{Stab}_{\mathcal G}(\mathbf k[[u,v]]^\times\bullet\Delta^{\mathcal W}_{r,l})
=\mathsf{DMR}_0(\mathbf k), 
$$
where the inclusion follows from Thm. \ref{thm:014}(a), 
the first equality follows from Thm. \ref{thm:016}, and the last
equality follows from Thm. \ref{thm:013}. 
(b) follows from (a) and from Thm. \ref{thm:014}(b). 
(c) follows from (a) and (b) by taking Lie algebras.  
\end{proof}

The conclusions of Cor. \ref{main:cor}(c) can be made explicit using the following result.  

\begin{lem} (explicit description of 
$\mathfrak G_{\mathrm{inert}}$ and $\mathrm{Lie}\Theta$)
\\
 (a) The Lie subalgebra $\mathfrak G_{\mathrm{inert}}$ of $(\mathfrak G,\langle,\rangle)$ is 
 $\mathfrak G\cap\mathrm{ker}(push-id)$, where $push$ is the linear endomorphism of $U\mathfrak{lie}_{\{0,1\}}$ 
 given by $e_\infty^{a_0}e_0e_\infty^{a_1}e_0\cdots e_0e_\infty^{a_r}\mapsto 
 e_\infty^{a_r}e_0e_\infty^{a_0}e_0\cdots e_0e_\infty^{a_{r-1}}$ for $r\geq0$,
 $a_0,\ldots,a_r\geq0$. 
\\
(b) The linear endomorphism of $U\mathfrak{lie}_{\{0,1\}}$ given by 
$$
\textstyle a\mapsto s_{(0,\infty)}(\sum_{i\geq0}((-1)^i/i!)e_\infty^ie_0\partial_\infty^i(a_0))  
$$
where $a\mapsto a_\infty,a_0$ are the endomorphisms of 
$U\mathfrak{lie}_{\{0,1\}}$ defined by the identity 
$a=a_\infty e_\infty+a_0e_0$ and $\partial_\infty$ is the derivation of $U\mathfrak{lie}_{\{0,1\}}$ such that 
$e_\infty\mapsto 1,e_0\mapsto 0$ induces 
the involution $\mathrm{Lie}\Theta$ of $(\mathfrak G_{\mathrm{inert}},\langle,\rangle)$. 
\\
\end{lem}

\begin{proof}
Thm. 2.1 from \cite{Sch1} (with $x:=e_\infty,y:=e_0$) implies both (a) and that 
the linear endomorphism of $U\mathfrak{lie}_{\{0,1\}}$ given by 
 $\textstyle a\mapsto b_a:=\sum_{i\geq0}((-1)^i/i!)e_\infty^ie_0\partial_\infty^i(a_0)$
induces the map $\mathfrak G_{\mathrm{inert}}\to \mathfrak G$ 
from Lem. \ref{lem06:1109}(c). (b) then follows from 
$\mathrm{Lie}(\Theta)(a)=s_{(0,\infty)}(b_a)$ (see Lem. \ref{lem02:2506:wo:proof}(d)). 
\end{proof}

\subsubsection{Bitorsor consequences}

\begin{thm}\label{thm:0:31:25jan} (Thm. \ref{thm:DMRmu})
  Let $\mu\in\mathbf k$. Then (with the notation of Def. \ref{def:G:mu:inert} 
  and Lem. \ref{def:Theta:mu}):
  
  (a) the inclusion $\mathsf{DMR}_\mu(\mathbf k)\subset\mathcal G^\mu_{\mathrm{inert}}$ holds
  (inclusion of sets); 
  
  (b) the subset $\mathsf{DMR}_\mu(\mathbf k)$ of $\mathcal G^\mu_{\mathrm{inert}}$ is stable 
  under the involution $\Theta^\mu$ of this set. 
\end{thm}
 
\begin{thm} \label{thm:0:33} (Thms. \ref{thm:betti:a} and \ref{thm:betti:b})
The subgroup $\mathsf{DMR}^{\mathrm B}(\mathbf k)$ of $\mathcal G^{\mathrm B}\rtimes\mathbf k^\times$
(see Defs. \ref{lem:mathcalG:betti} and \ref{def:DMLR:betti})
is contained in $(\mathcal G^{\mathrm B}\rtimes\mathbf k^\times)_{\mathrm{inert}}$ and is 
stable under the involution $\Theta^{\mathrm B}$ of this group. 
\end{thm}
It follows that the bitorsor $(\mathsf{DMR}_0(\mathbf k)\rtimes\mathbf k^\times,
\sqcup_{\mu\in\mathbf k^\times}\mathsf{DMR}_\mu(\mathbf k),
  \mathsf{DMR}^{\mathrm B}(\mathbf k))$ is a subtorsor of the torsor 
$(\mathcal G_{\mathrm{inert}} \rtimes \mathbf k^\times,
\sqcup_{\mu\in\mathbf k^\times}\mathcal G_{\mathrm{inert}}^\mu,
(\mathcal G^{\mathrm B} \rtimes \mathbf k^\times)_{\mathrm{inert}})$ and is stable 
under the torsor involution $(\Theta\rtimes id,\sqcup_{\mu\in\mathbf k^\times}\Theta^\mu,
\Theta^{\mathrm B}\rtimes id)$, see Lem. \ref{lem:invol:bitorsor}.

\subsubsection{Relationship with the Kashiwara-Vergne Lie algebra}\label{sect:announcement:krv}


Let $\mathbb L$ be the Lie algebra freely generated by $x,y$
($\mathbb L,x,y$ are denoted $\mathfrak{Lib}(X),x_0,x_1$ in the specialization to $\Gamma=1$ of \cite{R}). 
Let $i : \mathbb L\to\mathfrak{lie}_{\{0,1\}}$  be the Lie algebra isomorphism induced by $x\mapsto e_0$, $y\mapsto -e_1$. 
Then $i$ induces a Lie algebra isomorphism $\mathbb L\stackrel{\sim}{\to}\mathfrak{lie}_{\{0,1\}}$, the bracket on 
$\mathbb L$ being given by 
$$
\langle a,b\rangle:=[a,b]+d_a^x(b)-d_b^x(a), \quad \textrm{where} \quad d^x_a : x\mapsto 0,y\mapsto[y,a]. 
$$

\begin{defn}\label{def:ds}
  $\mathfrak{ds}$ 
is the preimage by $i$ of 
$\mathfrak{dmr}_0$ 
(this Lie algebra is denoted  $\mathfrak{dmr}_0$ in 
 \cite{R}, §3.3.8 and $\mathfrak{dmr}$ in \cite{FK}).   
\end{defn}

\begin{defn} (see  \cite{FK} for (a,b,c))
 (a) $\mathfrak{der}$ is the Lie algebra of derivations of $\mathbb L$; 

(b) $\nu : \mathbb L\to \mathfrak{der}$ is the linear map $\tilde f\mapsto d_F^{\mathrm{FK}}$, where 
(see (0.2) in {\it loc. cit.}), where $f(x,y):=\tilde f(x,-y)$, $F(x,y):=f(-x-y,y)$ (see (2.2) in {\it loc. cit.}), 
and $d_F^{\mathrm{FK}}$ is defined by $y\mapsto [y,F(x,y)]$, $x+y\mapsto 0$. 

(c) $\mathfrak{sder}\subset\mathfrak{der}$ is the set of derivations $D$ such that there 
exist $u,v\in \mathbb L$ with $D : x\mapsto [x,u],y\mapsto [y,v], x+y\mapsto 0$. 
\end{defn}

\begin{lem}\label{main:cor:c} (see Lem. \ref{main:cor:c:w:proof})
If $\tilde f \in\mathfrak{ds}$, then the mould $M:=\mathrm{ma}_{\tilde f}$ 
defined in \cite{FK} satisfies the "senary relation" (3.64) from \cite{Ec2}
(see also (1) in \cite{FK}). 
\end{lem}

The paper \cite{Ec2} contains the statement that any "alternal//alternil mould" satisfies (3.64), 
which is unproven there, but has been proved in \cite{Kaw}. 
A specialization of this statement is equivalent to Lem. \ref{main:cor:c} of the present paper, and is proved
in \cite{Sch2} on the basis of independent techniques.

In \cite{AT}, one defines a Lie subalgebra $\mathfrak{krv}\subset\mathfrak{sder}$ (see Def. 3.2 in \cite{FK}). 

\begin{thm} \label{thm:schneps:nu} (see Thm. \ref{thm:schneps:nu:w:proof})
(see \cite{Sch1,Sch2}) The map $\nu : \mathbb L\to\mathfrak{der}$ induces an injection of Lie algebras $\mathfrak{ds}\hookrightarrow\mathfrak{krv}$. 
\end{thm}

In \S\ref{sect:comments:involutions}, we discuss the behavior of the terms of the sequence of inclusions 
$$
\mathfrak{grt}_1\subset \mathfrak{dmr}_0\subset \mathfrak{krv}
$$
with respect to the action of 
$\mathfrak S_3$ on $\mathfrak{krv}$.

\subsection{Organization of the text} 

This text is divided into six parts: \\
\textbullet\ Parts 1 to 4 are devoted to the proof of Cor. \ref{main:cor}, with Parts 1 (resp. 2, 3, 4) 
proving Thms. \ref{thm:013} (resp. \ref{thm:015}, \ref{thm:016}, \ref{thm:014});  \\
\textbullet\ Part \ref{part 5} is devoted to the proof of Thms. \ref{thm:0:31:25jan} and \ref{thm:0:33}, 
which is the bitorsor counterpart of Cor. \ref{main:cor}; \\
\textbullet\ Part \ref{part:kv} is devoted to the proof of Lem. \ref{main:cor:c} and Thm. \ref{thm:schneps:nu}, 
which are the applications of the Cor. \ref{main:cor} to the relation between the double shuffle and 
Kashiwara-Vergne Lie algebras. \\
Each of the parts starts with its own introduction, in which its contents are described.  


\begin{notation}\label{NOTATION}
(a) For $S$ a set of formal variables, $\mathbf k[[S]]$ is the algebra of commutative formal series over $S$.   

(b) For $A$ an associative $\mathbb Q$-algebra, $\mathrm{Der}_{\mathbb Q\operatorname{-alg}}(A)$ is the Lie algebra
of derivations of $A$, and $\mathrm{ad}_a\in \mathrm{Der}_{\mathbb Q\operatorname{-alg}}(A)$ is the inner derivation 
associated to an element $a\in A$. 

(c) For $X$ an object in a category and $Y$ a subobject, $i_{Y,X} : Y\to X$ is
the corresponding morphism. 

(d) We denote by $\mathrm{pr}_{\mathcal G}$ the projections $\mathcal G\times H\to\mathcal G$ and $H\rtimes\mathcal G\to \mathcal G$, 
where $H$ is any group or group equipped with an action of $\mathcal G$. 

(e) For $A$ an associative unital algebra, $A^\times$ is the group of its units (invertible elements). 
\end{notation}

\paragraph{\textbf{Acknowledgements}} 
The collaboration of both authors was partially supported the following grants:  
ANR project HighAGT ANR20-CE40-0016, JSPS KAKENHI JP18H01110, JP 24K00520,  JP 24K21510. 
Both authors 
also thank Anton Alekseev and Leila Schneps for discussions.

\newpage

\part{The group equality $\mathrm{Stab}_{\mathcal G}(\mathbf k[[u,v]]^\times\bullet\Delta^{\mathcal W}_{r,l}
)=\mathsf{DMR}_0(\mathbf k)$}\label{part 1}

The objective of Part \ref{part 1} is the proof of Thm. \ref{thm:013} (Lem. \ref{lem:12:2209} and 
Thm. \ref{thm:main}).  
In §\ref{sect:1}, we formulate these results, prove Lem. \ref{lem:12:2209} and give a 
plan for the proof of Thm. 
\ref{thm:main}, which is carried out in the remainder of Part \ref{part 1}. In 
§\ref{sect:2} and §\ref{sect:3:1009}, we set up diagrams of 
pointed sets with group actions.
In §\ref{sect:4:korea}, we prove that the corresponding morphisms of stabilizer groups are isomorphisms, thus 
proving Thm. \ref{thm:main}.

\section{Statement of the group equality  $\mathrm{Stab}_{\mathcal G}(\mathbf k[[u,v]]^\times\bullet\Delta^{\mathcal W}_{r,l}
)=\mathsf{DMR}_0(\mathbf k)$ and plan of proof}\label{sect:1}

The purpose of this section is to formulate Lem. \ref{lem:12:2209} and Thm. \ref{thm:main}, to prove 
Lem. \ref{lem:12:2209} and to give a sketch of the proof of Thm. \ref{thm:main}. The three first steps are 
carried out in §\ref{sect:1:1}. §\ref{sect:1:2} contains categorical material, which is necessary for 
the proof of Thm. \ref{thm:main}, as well as for other results of this paper. The sketch of proof of 
Thm. \ref{thm:main} is given in §\ref{sect:1:3}: there, we formulate Props. \ref{prop:1:7}, \ref{prop:1:8}, 
\ref{prop:1:9} and Prop. \ref{lem:final:c:2909}(a)(b) (whose proofs are carried out in §§\ref{sect:3:1009} 
and \ref{sect:4:korea}), and show that these statements imply Thm. \ref{thm:main} (see Prop. 
\ref{lem:final:c:2909}(c)).

\subsection{Proofs of Lems. \ref{def:mathcalG} and \ref{lem:actions:ktimes:a}}

We now give the proofs of Lems. \ref{def:mathcalG} and \ref{lem:actions:ktimes:a}, 
the statements of which we recall. 
\begin{lem}\label{def:mathcalG:w:proof} (see Lem. \ref{def:mathcalG})
The product defined by 
\begin{equation}\label{FORMULA:1912}
(g\circledast h)(e_0,e_1):=h(g(e_0,e_1)\cdot e_0\cdot g(e_0,e_1)^{-1},e_1)\cdot g(e_0,e_1). 
\end{equation}
equips $\mathcal G$ with a group structure. 
\end{lem}
\begin{proof}
The fact that $\mathrm{exp}(\mathfrak{lie}_{\{0,1\}}^\wedge)$ equipped with 
the product \eqref{FORMULA:1912} is a group, follows from \cite{R}, Prop. 3.1.6. 
One then checks that the map 
\begin{equation}\label{MORPH}
(\mathrm{exp}(\mathfrak{lie}_{\{0,1\}}^\wedge),\circledast)\to 
(\mathbf k^2,+),\quad g\mapsto((g|e_0),(g|e_1))
\end{equation}
defines a group morphism, which implies the statement.
\end{proof}

\begin{lem}\label{lem:actions:ktimes:a:w:proof} (see Lem. \ref{lem:actions:ktimes:a})
(a)   The action of $\mathbf k^\times$ on $\mathbf k\langle\langle e_0,e_1\rangle\rangle$,
where $\lambda\in\mathbf k^\times$ acts by taking $e_0,e_1$ to $\lambda e_0,\lambda e_1$, 
induces an action of $\mathbf k^\times$ on $\mathrm{exp}(\mathfrak{lie}_{\{0,1\}}^\wedge)$, 
which restricts to an action on $\mathcal G$. 

(b) This action is compatible with the product $\circledast$, hence is an action by 
automorphisms of the group $(\mathcal G,\circledast)$.  
\end{lem}

\begin{proof}
 (a) and (b) follow from \cite{EF2}, §1.6.3 and from the equivariance of the group morphism \eqref{MORPH} under 
 the action of $\mathbf k^\times$. 
\end{proof}

\subsection{Statement of group equality $\mathrm{Stab}_{\mathcal G}(\mathbf k[[u,v]]^\times\bullet\Delta^{\mathcal W}_{r,l}
)
=\mathsf{DMR}_0(\mathbf k)$}\label{sect:1:1}

The action of $\mathbf k^\times$ on $(\mathrm{exp}(\mathfrak{lie}_{\{0,1\}}^\wedge),\circledast)$
preserves the subgroup $\mathsf{DMR}_0(\mathbf k)$; this results in the inclusion of the corresponding 
semidirect products
$$
\mathsf{DMR}_0(\mathbf k)\rtimes\mathbf k^\times
\subset \mathrm{exp}(\mathfrak{lie}_{\{0,1\}}^\wedge)
\rtimes\mathbf k^\times
$$
(denoted $\mathsf{DMR}^\DR(\mathbf k)\subset \mathsf G^\DR(\mathbf k)$ in {\it loc. cit.}). 
In \cite{EF2}, one constructs subgroups 
$\mathsf{Stab}(\hat\Delta^{\mathcal M,\DR})(\mathbf k)$ and $\mathsf{Stab}(\hat\Delta^{\mathcal W,\DR})(\mathbf k)$ 
of $\mathrm{exp}(\mathfrak{lie}_{\{0,1\}}^\wedge)\rtimes\mathbf k^\times$. The intersections 
of these groups with $(\mathrm{exp}(\mathfrak{lie}_{\{0,1\}}^\wedge),\circledast)$ are stable under the action of 
$\mathbf k^\times$, and
these groups are the corresponding semidirect products. 

\begin{lem}\label{lem:12:2209}
One has 
$$
\mathsf{DMR}_0(\mathbf k)=\mathsf{Stab}(\hat\Delta^{\mathcal M,\DR})(\mathbf k)\cap \mathcal G
$$
(equality of subgroups of $(\mathrm{exp}(\mathfrak{lie}_{\{0,1\}}^\wedge),\circledast)$). 
\end{lem}

\begin{proof}
The translation of \cite{EF0}, Thm. 1.2 in the language of the present paper was done in the proof of 
Lem. 3.7 in \cite{EF2} and yields the equality 
$$
\{e^{\alpha e_1}|\alpha\in\mathbf k\}\cdot\mathsf{DMR}_0(\mathbf k)\cdot\{e^{\beta e_0}|\beta\in\mathbf k\}
=\mathsf{Stab}(\hat\Delta^{\mathcal M,\DR})(\mathbf k)\cap \mathrm{exp}(\mathfrak{lie}_{\{0,1\}}^\wedge)
$$
of subsets of $\mathrm{exp}(\mathfrak{lie}_{\{0,1\}}^\wedge)$. Since $\mathsf{DMR}_0(\mathbf k)
\subset\mathcal G$, the intersection of the left-hand side of this equality with $\mathcal G$ is 
$\mathsf{DMR}_0(\mathbf k)$. The intersection of its right-and side is  
$\mathsf{Stab}(\hat\Delta^{\mathcal M,\DR})(\mathbf k)\cap \mathcal G$. The result therefore follows from 
the intersection of the said equality with $\mathcal G$.  
\end{proof}

In Lem.  \ref{lem:3101009:BIS}, one defines a pointed set 
$(\mathbf k[[u,v]]^\times\backslash\mathrm{Hom}_{\hat{\mathcal C}_{\mathbf k}\operatorname{-alg}}(\hat{\mathcal W},\hat V),
\mathbf k[[u,v]]^\times\bullet\Delta^{\mathcal W}_{r,l}
)$, and equips it with an action of the group $(\mathcal G,\circledast)$ (see Lem. \ref{def:mathcalG}. 
The corresponding stabilizer group, denoted  $\mathrm{Stab}_{\mathcal G}(
\mathbf k[[u,v]]^\times\bullet\Delta^{\mathcal W}_{r,l}
)$, is a subgroup of this group.

The main result of Part 1 is
\begin{thm}\label{thm:main}
One has the equality 
\begin{equation*}
\mathsf{Stab}(\hat\Delta^{\mathcal M,\DR})(\mathbf k)\cap \mathcal G
=\mathrm{Stab}_{\mathcal G}(
\mathbf k[[u,v]]^\times\bullet\Delta^{\mathcal W}_{r,l}
). 
\end{equation*}
of subgroups of $(\mathcal G,\circledast)$. 
\end{thm}

\subsection{Categorical material}\label{sect:1:2}

\begin{defn}
A pointed set is a pair $(X,x_0)$, where $X$ is a set and $x_0\in X$. If $(X,x_0)$ and $(Y,y_0)$ are pointed sets, a morphism 
of pointed sets $f : (X,x_0)\to (Y,y_0)$ is a set map $f : X\to Y$ such that $f(x_0)=y_0$. We denote by $\mathbf{PS}$ the category of pointed 
sets.  
\end{defn}

\begin{defn}\label{def:psga}
(a) A {\it $\mathcal G$-equivariant pointed set with group action ($\mathcal G$-PSGA)} is a tuple $(X,x_0,A,\bullet,*)$, 
there $(X,x_0)$ is a pointed set, $A$ is a group, 
$\bullet$ is an action of $A$ on $X$, and $*$ is a pair of actions of $\mathcal G$ on the group 
$A$ and on the set $X$, such that $g*(a\bullet x)=(g*a)\bullet(g*x)$. 

(b) If $(X,x_0,A,\bullet,*)$ and  $(Y,y_0,B,\bullet,*)$ are $\mathcal G$-PSGAs, a $\mathcal G$-PSGA morphism 
$(X,x_0,A,\bullet,*)\to(Y,y_0,B,\bullet,*)$
 is a pair $(f,\alpha)$, where $f : (X,x_0)\to(Y,y_0)$ is a morphism of pointed sets and $\alpha : A\to B$ is a group morphism, satisfying the identities
 $f(a\bullet x)=\alpha(a)\bullet f(x)$, $\alpha(g*a)=g*\alpha(a)$, $f(g*x)=g*f(x)$. 
\end{defn}

\begin{defn}
(a) A {\it split extension (SE) of $\mathcal G$} is a triple $(\Sigma,p,i)$ where $\Sigma$ is a group and 
$p : \Sigma\to\mathcal G$ and $i : \mathcal G\to \Sigma$ are group morphisms such that 
$p\circ i=id_{\mathcal G}$. 

(b) If $(\Sigma,p,i)$ and $(\Upsilon,q,j)$ are split extensions of $\mathcal G$, a morphism 
$(\Sigma,p,i)\to(\Upsilon,q,j)$ is a group morphism $\phi : \Sigma\to\Upsilon$, such that 
$q\circ\phi=p$ and $\phi\circ i=j$. 
\end{defn}

\begin{defn}
    (a) A {\it pointed set with action of a group with split morphism to $\mathcal G$ ($\mathcal G$-PSAGSM)} is a tuple
$(X,x_0,\Sigma,p,i,*)$ where $(X,x_0)$ is a pointed set, $(\Sigma,p,i)$ is a SE of $\mathcal G$, and $*$ is an action of 
$\Sigma$ on $X$. 

(b) A morphism $(X,x_0,\Sigma,p,i,*)\to(Y,y_0,\Upsilon,q,j,*)$ of $\mathcal G$-PSAGSMs is a pair $(f,\sigma)$ where $f : (X,x_0)\to (Y,y_0)$
is a morphism of pointed sets and $\phi : (\Sigma,p,i)\to(\Upsilon,q,j)$ is a morphism of SEs of $\mathcal G$, satisfying the identity
$f(\sigma*x)=\phi(\sigma)*f(x)$. 
\end{defn}

\begin{lem}
(a) When equipped with the obvious identities and compositions of morphisms, both the class of $\mathcal G$-PSGAs and of 
$\mathcal G$-PSAGSMs form categories, which we denote $\mathcal G\operatorname{-}\mathbf{PSGA}$ and $\mathbf{PSGA}_{\mathcal G}^\rtimes$. 

(b) There is a functor 
 $\mathcal G\operatorname{-}\mathbf{PSGA}\to\mathbf{PSGA}_{\mathcal G}^\rtimes$ 
given by $(X,x_0,A,\bullet,*)\mapsto (X,x_0,A\rtimes\mathcal G,p,i,\odot)$, where 
$(A\rtimes\mathcal G,p,i)$ is the SE of $\mathcal G$ (semidirect product) corresponding to the action of 
$\mathcal G$ on $A$, and $\odot$ is the action of $A\rtimes\mathcal G$ derived from the actions of 
$\mathcal G$ and $A$ on $X$, and from their compatibility with the action of $\mathcal G$ on $A$. 

(c) There is a functor
 $\mathbf{PSGA}_{\mathcal G}^\rtimes\to\mathcal G\operatorname{-}\mathbf{PSGA}$ induced by 
$(X,x_0,\Sigma,p,i,\odot)\mapsto(X,x_0,A,\bullet,*)$, where $A:=\mathrm{ker}(p)\subset\Sigma$, its action on $X$ is the restriction of that of $\Sigma$, 
the action of $\mathcal G$ on $\Sigma$ is induced by the SE structure $(\Sigma,p,i)$ of $\Sigma$, and its action on $X$ is the pullback of that of $\Sigma$ by 
$i : \mathcal G\to\Sigma$.  

(d) The functors from (b) and (c) are quasi-inverse of one another, and build up a category equivalence 
 $\mathcal G\operatorname{-}\mathbf{PSGA}\simeq\mathbf{PSGA}_{\mathcal G}^\rtimes$.
 \end{lem}

\begin{proof}
One first observes that both the classes of SEs of $\mathcal G$ and of groups with an action of $\mathcal G$ form categories, and that 
there are functors $\{$groups with an action of $\mathcal G\}\to\{$SEs of $\mathcal G\}$ given by $A\mapsto (A\rtimes\mathcal G,p,i)$, where 
$p : A\rtimes\mathcal G\to\mathcal G$ and $i : \mathcal G\to A\rtimes\mathcal G$ are the natural morphisms, and 
$\{$SEs of $\mathcal G\}\to\{$groups with an action of $\mathcal G\}$ taking $(\Sigma,p,i)$ to the group $\mathrm{ker}(p)$, which is 
a normal subgroup of $\Sigma$, so that the action of $\Sigma$ by conjugation restricts to an action on it, equipped with the action of $\mathcal G$ defined as the 
pull-back by $i : \mathcal G\to\Sigma$ of this action. The statements then follow from the equivalence of an   
action of a semidirect product group $A\rtimes\mathcal G$ on a set with a pair of actions of $A$ and $\mathcal G$, which also satisfy the compatibility condition 
from Def. \ref{def:psga}(a). 
\end{proof}

\begin{defn}
(a) A {\it group with morphism to $\mathcal G$} is a pair $(\Sigma,p)$ where $\Sigma$ is a group and 
$p : \Sigma\to\mathcal G$ is a group morphism. 

(b) If $(\Sigma,p)$ and $(\Upsilon,q)$ are groups with morphism $\mathcal G$, a morphism 
$(\Sigma,p)\to(\Upsilon,q)$ is a group morphism $\phi : \Sigma\to\Upsilon$, such that 
$q\circ\phi=p$. 
\end{defn}
We denote by $\mathbf{Gp}_{\mathcal G}$ the category of groups with morphism to $\mathcal G$. 

\begin{defn}
    (a) A pointed set with action of a group with morphism to $\mathcal G$ ($\mathcal G$-PSAGM) is a tuple
$(X,x_0,\Sigma,p,*)$ where $(X,x_0)$ is a pointed set, $(\Sigma,p)$ is a group with morphism to 
$\mathcal G$, and $*$ is an action of 
$\Sigma$ on $X$. 

(b) A morphism $(X,x_0,\Sigma,p,*)\to(Y,y_0,\Upsilon,q,*)$ of $\mathcal G$-PSAGMs is a pair $(f,\phi)$ where $f : (X,x_0)\to (Y,y_0)$
is a morphism of pointed sets and $\phi : (\Sigma,p)\to(\Upsilon,q)$ is a morphism of groups with morphisms to
$\mathcal G$, satisfying the identity $f(\sigma*x)=\phi(\sigma)*f(x)$. 

(c) We denote by $\mathbf{PSGA}_{\mathcal G}$  the category of $\mathcal G$-PSAGMs. 
\end{defn}

\begin{lem}
The assignment $(X,x_0,\Sigma,p,i,*)\mapsto (X,x_0,\Sigma,p,*)$ defines a
functor 
$$
\mathbf f :\mathbf{PSGA}_{\mathcal G}^\rtimes\to\mathbf{PSGA}_{\mathcal G}.
$$ 
\end{lem}

\begin{proof}
    Obvious. 
\end{proof}

\begin{defn}
    A {\it pointed set with action of $\mathcal G$} is a triple $(X,x_0,\bullet)$, where $(X,x_0)$ is a pointed set and 
    $\bullet$ is an action of $\mathcal G$ on $X$. A morphism $(X,x_0,\bullet)\to(Y,y_0,\bullet)$ of pointed sets with action of $\mathcal G$
     is a morphism $f:(X,x_0)\to(Y,y_0)$ of pointed sets, which is $\mathcal G$-equivariant.  
\end{defn}

The class of pointed sets with action of $\mathcal G$ then forms a category, which will be denoted $\mathbf{PS}_{\mathcal G}$. 

\begin{lem}\label{lem:foncteur:q}
For $(X,x_0,A,\bullet,*)$ a $\mathcal G$-PSGA, there is an action $*$ of $\mathcal G$ on $A\backslash X$, 
uniquely defined by the condition that the projection $X\to A\backslash X$ is 
$\mathcal G$-equivariant. The assignment $(X,x_0,A,\bullet,*)\mapsto (A\backslash X,A\bullet x_0,*)$ 
defines a functor $$\mathbf q : \mathcal G\operatorname{-}\mathbf{PSGA}\to\mathbf{PS}_{\mathcal G}.$$ 
\end{lem}

\begin{proof}
The proof is straightforward.
\end{proof}

\begin{defn}
$\mathbf{Subgp}_{\mathcal G}$ is the category where the objects are the subgroups of $\mathcal G$, and the morphisms are the inclusions of 
subgroups of $\mathcal G$. 
\end{defn}

\begin{lem}\label{lem:def:stab:ulstab}\label{BASIC}
(a) If $(\varphi,f) : (X,x,\Gamma,\bullet)\to(X',x',\Gamma',\bullet')$ is a morphism of pointed
sets with group actions, then $\varphi(\mathrm{Stab}_\Gamma(x))\subset \mathrm{Stab}_{\Gamma'}(x')$.

(b) The assignment $(X,x_0,\bullet)\mapsto \mathrm{Stab}_{\mathcal G}(x_0)$ 
defines a functor $$\mathrm{Stab} : \mathbf{PS}_{\mathcal G}\to\mathbf{Subgp}_{\mathcal G}.$$ 

(c)    The assignment $(X,x_0,\Sigma,*,p)\mapsto (\mathrm{Stab}_\Sigma(x_0),p|_{\mathrm{Stab}_\Sigma(x_0)})$ defines a 
functor $\underline{\mathrm{Stab}} : \mathbf{PSGA}_{\mathcal G}\to
\mathbf{Gp}_{\mathcal G}$. 
\end{lem}

\begin{proof}
  Immediate. 
\end{proof}

There are functors 
$$
can : \mathbf{Subgp}_{\mathcal G}\to\mathbf{Gp}_{\mathcal G},
$$
taking a subgroup $H$ of $\mathcal G$
to the pair $(H,H\hookrightarrow\mathcal G)$ and 
$\mathbf c : \mathbf{PS}_{\mathcal G}\to\mathcal G\operatorname{-}\mathbf{PSGA}$ given by $(X,x_0,*)\mapsto (X,x_0,\{1\},\bullet,*)$, 
where $\bullet$ is the trivial action of the trivial group $\{1\}$ on $X$. 
The categories and functors of this section are summarized in the following diagram
\begin{equation}\label{categorical:diagram}
\xymatrix{
\mathcal G\operatorname{-}\mathbf{PSGA}\ar^{\simeq}[r]\ar@/^1pc/^{\mathbf q}[d]&
\mathbf{PSGA}^\rtimes_{\mathcal G}\ar^{\mathbf f}[r]&
\mathbf{PSGA}_{\mathcal G}\ar^{\underline{\mathrm{Stab}}}[r]&
\mathbf{Gp}_{\mathcal G}\\
\mathbf{PS}_{\mathcal G}\ar@/^1pc/^{\mathbf c}[u]\ar_{\mathrm{Stab}}[rrr]
&&&
\mathbf{Subgp}_{\mathcal G}\ar_{can}[u]}    
\end{equation}
Moreover, one checks that the composed functors $can\circ \mathrm{Stab}$ and 
$\underline{\mathrm{Stab}}\circ \mathbf f \circ \mathbf c: \mathbf{PS}_{\mathcal G}\to \mathbf{Gp}_{\mathcal G}$ are equivalent.

\begin{lem}\label{lem:tranfo:nat}
    The assignment $(X,x_0,A,\bullet,*)\mapsto ((p_{X,A},1) : (X,x_0,A,\bullet,*)\to(A\backslash X,A\bullet x_0,\{1\},\underline\bullet,\underline *))$, where 
    $\underline\bullet,\underline *$ are the natural actions of $\{1\},\mathcal G$ on $A\backslash X$, and $p_{X,A} : X\to A\backslash X$ and $1 : A\to\{1\}$ are the canonical morphisms, 
    defines a natural transformation relating the 
    endofunctors $id$ and $\mathbf c\circ \mathbf q$ of $\mathcal G\operatorname{-}\mathbf{PSGA}$.  
\end{lem}

\begin{proof}
    Obvious. 
\end{proof}

\subsection{Sketch of the proof of Thm. \ref{thm:main}}\label{sect:1:3}

The statement of Thm. \ref{thm:main} is obtained in Prop. \ref{lem:final:c:2909}(c). 
The proofs of Props. \ref{prop:1:7}, \ref{prop:1:8}, \ref{prop:1:9}, and Prop. \ref{lem:final:c:2909}(a)(b) are 
done in §§\ref{sect:a:diag:in:PS}, \ref{sect:3:1009} and \ref{sect:4:korea}. 

\begin{prop}\label{prop:1:7}
    (a) The pairs 
$$
\mathbf M:=(\mathrm{Cop}_{\hat{\mathcal C}}(\hat{\mathcal M}),\hat\Delta^{\mathcal M}), \quad 
\mathbf{WM}:=(\mathrm{Cop}_{{\mathcal C}\operatorname{-alg-mod}}(\hat{\mathcal W},\hat{\mathcal M}),
(\hat\Delta^{\mathcal W},\hat\Delta^{\mathcal M})), \quad 
\mathbf W:=(\mathrm{Cop}_{\hat{\mathcal C}\operatorname{-alg}}(\hat{\mathcal W}),\hat\Delta^{\mathcal W}),
$$
$$
\mathbf E:=(\mathrm{Cop}_{e_1}(\hat{\mathcal W}),\hat\Delta^{\mathcal W}), \quad 
\mathbf E':=((1+t^2\mathbf k[[t]])\backslash\mathrm{Cop}_{e_1}(\hat{\mathcal W}),[\hat\Delta^{\mathcal W}]),
$$
$$
\mathbf E'':=(\mathbf k[[u,v]]^\times\backslash\mathrm{Cop}_{e_1}(\hat{\mathcal W}),[\!\![\hat\Delta^{\mathcal W}]\!\!]),\quad 
\mathbf H'':=(\mathbf k[[u,v]]^\times\backslash\mathrm{Hom}_{\hat{\mathcal C}\operatorname{-alg}}(\hat{\mathcal W},\hat V),
\mathbf k[[u,v]]^\times\bullet\Delta^{\mathcal W}_{r,l}
)
$$
where the right-hand sides are as in Defs. \ref{def:M:W:MW}, \ref{def:E:1001}, \ref{def:E':E'':H''} are objects in $\mathbf{PS}$. 

(b) The morphisms from Defs.  \ref{def:M:W:MW}, \ref{def:i:EW}, Lem. \ref{lem:def:p_UV} set up the diagram 
\begin{align*}
\xymatrix{
\mathbf M&
\ar_{p_{\mathbf M}}[l]\mathbf{WM}\ar^{p_{\mathbf W}}[r]&
\mathbf{W}&
\ar_{i_{\mathbf{EW}}}[l]\mathbf{E}\ar^{p_{\mathbf{EE'}}}[r]&
\mathbf{E}'\ar^{p_{\mathbf{E'E''}}}[r]&
\mathbf{E}''\ar^{i_{\mathbf{E''H''}}}[r]&
\mathbf{H}''
}   
\end{align*}
of morphisms in $\mathbf{PS}$. 
\end{prop}

\begin{proof}
The statements follow from the contents of \S\S\ref{sect:a:diag:in:PS}, \ref{sect:def:iEH}, \ref{sect:diag:XXXX}, \ref{sect:TRANSITION}.  
\end{proof}

Following \cite{EF1}, (1.4.1), define the map 
$\Gamma : \mathcal G\to 1+t^2\mathbf k[[t]]$, which takes 
 $g\in \mathcal G$ to 
\begin{equation}\label{def:Gamma:Phi:BIS}
\Gamma_g(t):=\mathrm{exp}(\sum_{n \geq 1}((-1)^{n+1}/n)(g|e_0^{n-1}e_1)t^n) 
\in 1+t^2\mathbf k[[t]]. 
\end{equation} 
The relation $\Gamma_g(t)\in1+t^2\mathbf k[[t]]$ for $g\in\mathcal G$ follows from $(g|e_1)=0$.

\begin{defn}\label{defn:sigma:BIS}
Set 
$$
\tilde{\mathcal G}:=\mathcal G \times (1+t^2\mathbf k[[t]]), 
$$
and define $\sigma  : \mathcal G\to\tilde{\mathcal G}$ as the map 
$g\mapsto (g,\Gamma_g(-t)^{-1})$. 
\end{defn}

\begin{lem}\label{lem:sigma:3107:BIS}
The map $\sigma$ is a group morphism.  
\end{lem}

\begin{proof}
Lemma 1.13, (a) in \cite{EF2} implies that the map $g\mapsto \Gamma_g(-t)^{-1}$ is a group morphism 
$(\mathrm{exp}(\mathfrak{lie}_{\{0,1\}}^\wedge),\circledast)\to 1+t\mathbf k[[t]]$. 
It follows that so is $\sigma$. 
\end{proof}

\begin{prop}\label{prop:1:8}
The pairs
$$
(\mathcal G,id),\quad (\tilde{\mathcal G},\mathrm{pr}_{\mathcal G}), 
$$
are objects in $\mathbf{Gp}_{\mathcal G}$, and 
$$
(\mathcal G,id)\stackrel{\sigma}{\to} (\tilde{\mathcal G},\mathrm{pr}_{\mathcal G}) 
$$ 
is a morphism in $\mathbf{Gp}_{\mathcal G}$, where $\sigma : \mathcal G\to G$ 
is in Lem. \ref{lem:sigma:3107:BIS}. 
\end{prop}

\begin{proof}
Immediate.  
\end{proof}

\begin{prop}\label{prop:1:9}
(a) The actions $\bullet_\Gamma$ and $\odot$ being given by Def. \ref{def32:1009} and Lem. \ref{lem:3101009:BIS}, 
and $\mathbf M$ and $\mathbf H''$ being as in Defs. \ref{def:M:W:MW} and \ref{def:E':E'':H''}, 
the pairs $(\mathbf M,\bullet_\Gamma)$ and $(\mathbf{H''},\odot)$ are objects in $\mathbf{PS}_{\mathcal G}$. 

(b) The actions $\bullet_{\mathcal M}$, $\bullet_{\mathcal W\mathcal M}$, $\bullet_{\mathcal W}$, $\odot$
being given by Def. \ref{def:35:1009} and Lem. \ref{lem:3101009:BIS}, the triples
$$
\mathbf X_{\mathcal G,\mathbf M}:=((\mathcal G,id),\mathbf M,\bullet_\Gamma), \quad
\mathbf X_{\tilde{\mathcal G},\mathbf M}:=((\tilde{\mathcal G},\mathrm{pr}_{\mathcal G}),\mathbf M,\bullet_{\mathcal M}), \quad 
\mathbf X_{\tilde{\mathcal G},\mathbf{WM}}:=((\tilde{\mathcal G},\mathrm{pr}_{\mathcal G}),\mathbf{WM},\bullet_{\mathcal{WM}}), 
$$
$$
\mathbf X_{\tilde{\mathcal G},\mathbf{W}}:=((\tilde{\mathcal G},\mathrm{pr}_{\mathcal G}),\mathbf{W},\bullet_{\mathcal{W}}), \quad 
\mathbf X_{\tilde{\mathcal G},\mathbf{E}}:=((\tilde{\mathcal G},\mathrm{pr}_{\mathcal G}),\mathbf{E},\odot), 
$$
$$
\mathbf X_{\mathcal G,\mathbf{E'},}:=((\mathcal G,id),\mathbf{E'},\odot), \quad 
\mathbf X_{\mathcal G,\mathbf{E''}}:=((\mathcal G,id),\mathbf{E''},\odot), \quad 
\mathbf X_{\mathcal G,\mathbf{H''}}:=((\mathcal G,id),\mathbf{H''},\odot).  
$$
are objects in $\mathbf{PSGA}_{\mathcal G}$. Moreover, $\mathbf X_{\mathcal G,\mathbf M}$ and
$\mathbf X_{\mathcal G,\mathbf{H''}}$ are the images of 
$(\mathbf M,\bullet_\Gamma)$ and $(\mathbf{H''},\odot)$ by the functor
$\mathbf f\mathbf c : \mathbf{PS}_{\mathcal G}\to\mathbf{PSGA}_{\mathcal G}$. 

(c) The following
\begin{equation}\label{diag:psga:2909}
 \xymatrix{\mathbf X_{\mathcal G,\mathbf M}\ar^{(\sigma,id)}[r]&
\mathbf X_{\tilde{\mathcal G},\mathbf M}&
\ar_{(id,p_{\mathbf M})}[l]\mathbf X_{\tilde{\mathcal G},\mathbf{WM}}\ar^{(id,p_{\mathbf W})}[r]&
\mathbf X_{\tilde{\mathcal G},\mathbf{W}}&
\ar_{(id,i_{\mathbf{EW}})}[l]\mathbf X_{\tilde{\mathcal G},\mathbf{E}}\ar^{(pr_1,p_{\mathbf{EE'}})}[r]&
\mathbf X_{\mathcal G,\mathbf{E}'}\ar^{(id,pr_{\mathbf{E'E''}})}[r]&
\mathbf X_{\mathcal G,\mathbf{E}''}\ar^{(id,i_{\mathbf{E''H''}})}[r]&
\mathbf X_{\mathcal G,\mathbf{H}''}
}   
\end{equation}
is a diagram of morphisms in $\mathbf{PSGA}_{\mathcal G}$. 
\end{prop}

\begin{proof}
(a) follows from \S\S\ref{sect:3:2} and \ref{sect:diag:XXXX}. The first part of (b) follows from \S\S\ref{sect:XXXX}, \ref{sect:TRANSITION} and 
\ref{sect:diag:XXXX}; the statement on 
 $\mathbf X_{\mathcal G,\mathbf M}$ and $\mathbf X_{\mathcal G,\mathbf{H''}}$ follows from the definition of $\mathbf f\mathbf c$. (c) follows from \S\S\ref{sect:XXXX}, \ref{sect:TRANSITION} 
and \ref{sect:diag:XXXX}. 
\end{proof}

\begin{prop}\label{lem:final:c:2909}
(a)  In the diagram 
$$
\tiny{ \xymatrix{\underline{\mathrm{Stab}}\mathbf X_{\mathcal G,\mathbf M}\ar[r]&
\underline{\mathrm{Stab}}\mathbf X_{\tilde{\mathcal G},\mathbf M}&
\ar[l]\underline{\mathrm{Stab}}\mathbf X_{\tilde{\mathcal G},\mathbf{WM}}\ar[r]&
\underline{\mathrm{Stab}}\mathbf X_{\tilde{\mathcal G},\mathbf{W}}&
\ar[l]\underline{\mathrm{Stab}}\mathbf X_{\tilde{\mathcal G},\mathbf{E}}\ar[d]&
&
\\
&&&&\underline{\mathrm{Stab}}\mathbf X_{\mathcal G,\mathbf{E}'}\ar[r]&
\underline{\mathrm{Stab}}\mathbf X_{\mathcal G,\mathbf{E}''}\ar[r]&
\underline{\mathrm{Stab}}\mathbf X_{\mathcal G,\mathbf{H}''}
}}   
$$
in $\mathbf{Gp}_{\mathcal G}$ obtained by applying $\underline{\mathrm{Stab}}$ to \eqref{diag:psga:2909} (see its definition in 
Lem. \ref{BASIC}), all the morphisms are isomorphisms. 

(b) The objects $\mathsf{Stab}(\hat\Delta^{\mathcal M,\DR})(\mathbf k)\cap \mathcal G$ 
and $\mathrm{Stab}_{\mathcal G}(\mathbf k[[u,v]]^\times\bullet\Delta^{\mathcal W}_{r,l}
)$ of 
$\mathbf{Subgp}_{\mathcal G}$ are the images of the objects 
$(\mathbf M,\bullet_\Gamma)$ and $(\mathbf{H''},\odot)$ by the functor 
$\mathrm{Stab} : \mathbf{PS}_{\mathcal G}\to \mathbf{Subgp}_{\mathcal G}$.

(c) One has $\mathsf{Stab}(\hat\Delta^{\mathcal M,\DR})(\mathbf k)\cap \mathcal G=
\mathrm{Stab}_{\mathcal G}(\mathbf k[[u,v]]^\times\bullet\Delta^{\mathcal W}_{r,l}
)$ (equality of subgroups of $\mathcal G$). 
\end{prop}

\begin{proof}
(a) This follows from Cor. \ref{cor:4:6:nov}, Prop. \ref{prop:4:10:nov}, Cor. \ref{cor:4:20:nov}, Prop. \ref{prop:4:22:nov}, Prop. \ref{prop:4:25:1110},
Prop. \ref{prop:4:34:nov}, and Prop. \ref{prop:missing:iso}. 
(b) It follows from §\ref{sect:3:2} that the stabilizer group of
$(\mathcal G,\mathrm{Cop}_{\mathcal C}(\hat{\mathcal M}),\bullet_\Gamma)$
is $\mathsf{Stab}(\hat\Delta^{\mathcal M,\DR})(\mathbf k)\cap \mathcal G$, which implies the first statement. The second 
statement follows from the definition of $(\mathbf{H''},\odot)$.      
(c) one has the following isomorphisms in $\mathbf{Gp}_{\mathcal G}$ (see \eqref{categorical:diagram}): 
\begin{align*}
    &
can(\mathsf{Stab}(\hat\Delta^{\mathcal M,\DR})(\mathbf k)\cap \mathcal G)
\simeq can\circ \mathrm{Stab}(\mathbf M,\bullet_\Gamma)\simeq
\underline{\mathrm{Stab}}\circ\mathbf f\mathbf c(\mathbf M,\bullet_\Gamma)
\simeq\underline{\mathrm{Stab}}\mathbf X_{\mathcal G,(\mathbf M}
\simeq \underline{\mathrm{Stab}}\mathbf X_{\mathcal G,\mathbf{H''}}
    \\&
\simeq \underline{\mathrm{Stab}}\circ\mathbf f\mathbf c(\mathbf{H''},\odot)    
\simeq can \circ \mathrm{Stab}(\mathbf{H''},\odot)\simeq can(\mathrm{Stab}_{\mathcal G}(
\mathbf k[[u,v]]^\times\bullet\Delta^{\mathcal W}_{r,l}
), 
\end{align*}
where the first and seventh isomorphisms follow from (b), the second and sixth isomorphisms follow from 
the equivalence of the functors $can\circ \mathrm{Stab}$ and 
$\underline{\mathrm{Stab}}\circ \mathbf f\mathbf c$, the third and fifth isomorphisms follow from the second statement in Prop. \ref{prop:1:9}(b), 
and the fourth isomorphism follows from Prop. \ref{lem:final:c:2909}(a). (c) then follows from the
fact that if $H,K\subset\mathcal G$ are subgroups such that $can(H)$ and $can(K)$ are isomorphic
objects of $\mathbf{Gp}_{\mathcal G}$, then $H=K$. 
\end{proof}

\section{Diagrams in $\mathcal G\operatorname{-}\mathbf{PSGA}$ and $\mathbf{PSGA}_{\mathcal G}$}
\label{sect:2}

Based on the preliminary material introduced in §§\ref{213:toto:BIS:5786} and \ref{sect:def:iEH}, we set up in 
§§\ref{sect:2:3:5786} and \ref{sect:diag:XXXX} diagrams of pointed sets with group actions (Cor. \ref{cor:HH} 
and Lem. \ref{lem:def:p_UV}).

\subsection{Background}\label{213:toto:BIS:5786}

\subsubsection{The sets $\mathrm{Cop}_{\mathcal C}(X)$}

For $\mathcal A$ a symmetric tensor category (STC; see \cite{McL}), define $\mathcal A$-alg to be the 
category of algebras in $\mathcal A$: an object in $\mathcal A$-alg is an
associative unital algebra in $\mathcal A$. This is again a STC. 
If $\mathcal A$ is a STC and $X$ is an object in $\mathcal A$, define $\mathrm{Cop}_{\mathcal A}(X):=
\mathrm{Hom}_{\mathcal A}(X,X^{\otimes2})$
(the set of coproducts of $X$).

\subsubsection{}

In \cite{EF1}, Part 1, we introduced the category $\mathcal C:=\mathbf k\operatorname{-mod}_{top}$ of topological $\mathbf k$-modules, which are 
$\mathbf k$-modules equipped with a decreasing map $\mathbb Z\to\{\mathbf k$-submodules of $M\}$, $i\mapsto F^iM$, such that 
there exists $i(M)$ with $F^{i(M)}M=M$ and $M\to \varprojlim_i M/F^iM$ is a $\mathbf k$-module isomorphism, and where 
the morphisms are the $\mathbf k$-module morphisms compatible with the filtrations. This is a STC with tensor product 
denoted $\hat\otimes$. A tensor functor from the category of graded  $\mathbf k$-modules to $\mathcal C$ is defined by taking 
$M=\oplus_d M_d$ to $\hat M:=\varprojlim_i M/F^iM$, where $F^iM:=\oplus_{j\geq i}M_j$, and $\hat M$ is then called a . 
complete $\mathbb Z_{\geq0}$-graded $\mathbf k$-module. 

\subsubsection{}
\label{sect:2:3:titi:BIS}

One checks that  
$$
\hat\Delta^{\mathcal W}\in\mathrm{Cop}_{\mathcal C\operatorname{-alg}}(\hat{\mathcal W}). 
$$
where $\hat\Delta^{\mathcal W}$ is given by \eqref{pties:Delta:W}.

\subsection{The morphism $i_{\mathbf E\mathbf H} : \mathbf E\to\mathbf H$ in $\mathbf{PS}$}\label{sect:def:iEH}

\begin{defn}
 $\mathrm{Cop}_{e_1}(\hat{\mathcal W})\subset \mathrm{Cop}_{{\mathcal C}\operatorname{-alg}}(\hat{\mathcal W})$ 
 is the set of all $\Delta$ such that $\Delta(e_1)=e_1 \otimes1+1\otimes e_1$.   
\end{defn}

It follows from \eqref{pties:Delta:W} that $\hat\Delta^{\mathcal W} \in \mathrm{Cop}_{e_1}(\hat{\mathcal W})$. 

\begin{defn}\label{def:E:1001}
 $\mathbf E$ is the pointed set $(\mathrm{Cop}_{e_1}(\hat{\mathcal W}),\hat\Delta^{\mathcal W})$.
\end{defn}

\begin{defn}
Set $\hat V:=\hat{\mathcal V}^{\hat\otimes2}$. Then $\hat V$ is an object in $\mathcal C\operatorname{-alg}$, with grading such that 
$e_i\otimes1,1\otimes e_i$ have degree $1$ for $i=0,1$.    
\end{defn}

\begin{lem}\label{lem:toto:0930}
(a) $\hat{\mathcal W}_r=\mathbf k\oplus e_1\hat{\mathcal V}$ is a complete graded subalgebra of $\hat{\mathcal V}$. 

(b) There is a unique map $\mathrm{Ad}_{e_1} : \hat{\mathcal W}\to\hat{\mathcal W}_r$, such that 
$\lambda+ve_1\mapsto \lambda+e_1v$ for any $(\lambda,v)\in\mathbf k\times\hat{\mathcal V}$. It is an 
isomorphism of complete graded algebras. 
   
(c) There is a unique map $\mathrm{Cop}_{e_1}(\hat{\mathcal W})\to 
\mathrm{Hom}_{\mathcal C\operatorname{-alg}}(\hat{\mathcal W},\hat V)$ given by 
$\Delta\mapsto(i_{\mathcal W_r,\mathcal V}\circ\mathrm{Ad}_{e_1})^{\otimes2}\circ\Delta$, 
$i_{\mathcal W_r,\mathcal V} : \hat{\mathcal W}_r\hookrightarrow\hat{\mathcal V}$ being the 
canonical inclusion. 
\end{lem}

\begin{proof}
    (a) is immediate. The first statement of (b) follows from the fact that the maps 
    $\mathbf k\times\hat{\mathcal V}\to \hat{\mathcal W}$, $(\lambda,v)\mapsto \lambda+ve_1$ and 
$\mathbf k\times\hat{\mathcal V}\to \hat{\mathcal W}_r$, $(\lambda,v)\mapsto \lambda+e_1v$
are both algebra isomorphisms, $\mathbf k\times\hat{\mathcal V}$ being equipped with the product 
$(\lambda,v)\cdot (\lambda',v')=(\lambda\lambda',\lambda v'_v\lambda'+ve_1v')$. 
(c) follows from the fact that both $\mathrm{Ad}_{e_1}$ and $i_{\mathcal W_r,\mathcal V}$
are algebra morphisms. 
\end{proof}

\begin{defn}\label{defn:2:7:korea}
(a)    $\Delta^{\mathcal W}_{r,l}:=(i_{\mathcal W_r,\mathcal V}\circ
    \mathrm{Ad}_{e_1})^{\hat\otimes 2}\circ\hat\Delta^{\mathcal W}\in
\mathrm{Hom}_{\hat{\mathcal C}\operatorname{-alg}}(\hat{\mathcal W},
\hat V)$ is the image of $\hat\Delta^{\mathcal W}\in\mathrm{Cop}_{e_1}(\hat{\mathcal W})$ 
by the map from Lem. \ref{lem:toto:0930}(c). 

(b) $\mathbf H$ is the pointed set 
$(\mathrm{Hom}_{\hat{\mathcal C}\operatorname{-alg}}(\hat{\mathcal W},
\hat V),\Delta^{\mathcal W}_{r,l})$. 

(c) $i_{\mathbf E\mathbf H} : \mathbf E\to\mathbf H$ is the morphism in $\mathbf{PS}$
induced by the map from Lem. \ref{lem:toto:0930}(c). 
\end{defn}

\subsection{The diagram $(\mathbf E,1+t^2\mathbf k[[t]],\bullet,*)\to(\mathbf E,\mathbf k[[u,v]]^\times,\bullet,*)\to
(\mathbf H,\mathbf k[[u,v]]^\times,\bullet,*)$ in $\mathcal G\operatorname{-}\mathbf{PSGA}$}\label{sect:2:3:5786}

Set 
$$
e_1:=e_1\otimes1\in\hat{\mathcal W}^{\hat\otimes2}\subset \hat V,\quad f_1:=1\otimes e_1\in\hat{\mathcal W}^{\hat\otimes2}\subset \hat V, 
$$
and 
$$
e_0:=e_0\otimes1\in\hat V,\quad f_0:=1\otimes e_0\in\hat V.  
$$

\begin{lem}\label{lem28:1001}
(a) There is a unique group morphism $\theta : (1+t^2\mathbf k[[t]],\cdot)\to
(\mathbf k[[u,v]]^\times,\cdot)$ such that $f(t)\mapsto f(u)f(v)/f(u+v)$.    

(b) The map $\mathbf k[[u,v]]^\times\times \mathrm{Cop}_{e_1}(\hat{\mathcal W})
\to \mathrm{Cop}_{e_1}(\hat{\mathcal W})$, $(f,\Delta)\mapsto 
f\bullet \Delta:=\mathrm{Ad}_{f(e_1,f_1)}\circ\Delta$ defines an action of 
$\mathbf k[[u,v]]^\times$ on $\mathrm{Cop}_{e_1}(\hat{\mathcal W})$.

(c) The map $\mathbf k[[u,v]]^\times\times 
\mathrm{Hom}_{\mathcal C\operatorname{-alg}}(\hat{\mathcal W},\hat V)
\to \mathrm{Hom}_{\mathcal C\operatorname{-alg}}(\hat{\mathcal W},\hat V)$, 
$(f,\Delta)\mapsto f\bullet \Delta:=
\mathrm{Ad}_{f(e_1,f_1)}\circ\Delta$ defines an action of 
$\mathbf k[[u,v]]^\times$ on 
$\mathrm{Hom}_{\mathcal C\operatorname{-alg}}(\hat{\mathcal W},\hat V)$. 

(d) The map from Lem. \ref{lem:toto:0930}(c) is equivariant with respect to the actions of $\mathbf k[[u,v]]^\times$.  
\end{lem}
    
\begin{proof}
(a) follows from the fact that $f(t)\mapsto f(u),f(v)$ and $f(u+v)$ are all group morphisms, and from the commutativity of 
$\mathbf k[[u,v]]^\times$. (b) follows from 
$$
\mathrm{Ad}_{f(e_1,f_1)}\circ\Delta(e_1)=f(e_1,f_1)(e_1+f_1)f(e_1,f_1)^{-1}=e_1+f_1
$$
for any $\Delta\in\mathrm{Cop}_{e_1}(\hat{\mathcal W})$, where the second equality follows
from the commutativity of $e_1$ and $f_1$. (c) is immediate. 
(d) follows from the fact that the algebra morphism
$(i_{\mathcal W,\mathcal V} \circ \mathrm{Ad}_{e_1})^{\otimes2} : \hat{\mathcal W}^{\hat\otimes2}\to\hat V$
intertwines the automorphisms $\mathrm{Ad}_{f(e_1,f_1)}$ of its source and target, which is itself a consequence of 
$(i_{\mathcal W,\mathcal V} \circ \mathrm{Ad}_{e_1})^{\otimes2}(f(e_1,f_1))=f(e_1,f_1)$.   
\end{proof}

\begin{lem}\label{lem:2:8:jan25}
(a) For $g\in\mathcal G$, there is a unique algebra automorphism $\mathrm{aut}_{g}^{\mathcal V}$ of $\hat{\mathcal V}$
such that 
\begin{equation}\label{def:aut:g:V:(1)}
\mathrm{aut}_g^{\mathcal V} : e_1 \mapsto e_1,\quad e_0 \mapsto g(e_0,e_1)\cdot e_0\cdot g(e_0,e_1)^{-1}. 
\end{equation}
The map 
$\mathcal G\to\mathrm{Aut}_{\mathcal C\operatorname{-alg}}(\hat{\mathcal V})$, $g\mapsto \mathrm{aut}_{g}^{\mathcal V}$
is a group morphism. 

(b) For any $g\in\mathcal G$, the automorphism $\mathrm{aut}_{g}^{\mathcal V}$ restricts to an automorphism of the 
subalgebra $\hat{\mathcal W}$; the corresponding restriction will be denoted $\mathrm{aut}_{g}^{\mathcal W}$. 
Then $\mathrm{aut}_{g}^{\mathcal W}$ is such that 
$1\mapsto 1$ and $ve_1\mapsto \mathrm{aut}_{g}^{\mathcal V}(v)e_1$ for any $v\in\hat{\mathcal V}$, and 
$\mathcal G\to\mathrm{Aut}_{\mathcal C\operatorname{-alg}}(\hat{\mathcal W})$, $g\mapsto \mathrm{aut}_{g}^{\mathcal W}$
is a group morphism. 

(c) The map $\mathcal G\to\mathrm{Aut}_{\mathcal C\operatorname{-alg}}(\hat V)$, $g\mapsto 
\mathrm{aut}_g^V:=(\mathrm{aut}_{g}^{\mathcal V})^{\otimes2}$
is a group morphism. 
\end{lem}

\begin{proof}
Immediate.  
\end{proof}

\begin{lem}\label{lem:43:0401:TER}
(a) There is an action $*$ of $\mathcal G$ on $\mathrm{Cop}_{e_1}(\hat{\mathcal W})$, 
defined by $g*\Delta:=(\mathrm{aut}_{g}^{\mathcal W})^{\otimes2}\circ
\Delta\circ (\mathrm{aut}_{g}^{\mathcal W})^{-1}$ 
for any 
$\Delta\in\mathrm{Cop}_{e_1}(\hat{\mathcal W})$, 
$g\in\mathcal G$. This action commutes with the action of $\mathbf k[[u,v]]^\times$ 
from Lem. \ref{lem28:1001}(b). 

(b) There is an action $*$ of $\mathcal G$ on 
$\mathrm{Hom}_{\mathcal C\operatorname{-alg}}(\hat{\mathcal W},\hat V)$, 
defined by $g*\Delta:=
\mathrm{aut}_{g}^V\circ\Delta\circ (\mathrm{aut}_{g}^{\mathcal W})^{-1}$ 
for any 
$\Delta\in\mathrm{Hom}_{\mathcal C\operatorname{-alg}}(\hat{\mathcal W},\hat V)$, 
$g\in\mathcal G$. 
This action commutes with the action of $\mathbf k[[u,v]]^\times$ from Lem. \ref{lem28:1001}(c). 

(c) The map from Lem. \ref{lem:toto:0930}(c) is 
equivariant with respect to the actions of $\mathcal G$ from (a) and (b). 
\end{lem}

\begin{proof}
(a) Since $g\mapsto \mathrm{aut}_g^{\mathcal W}$ defines an action of $\mathcal G$ on the 
algebra $\hat{\mathcal W}$, 
the said formula defines an action of $\mathcal G$ on 
$\mathrm{Cop}_{\mathcal C\operatorname{-alg}}(\hat{\mathcal W})$. 
Since $\mathrm{aut}_{g}^{\mathcal W}(e_1)=e_1$ for any $g\in\mathcal G$, 
this action preserves $\mathrm{Cop}_{e_1}(\hat{\mathcal W})$. 
The said commutativity follows from the equalities
$$
g*(f\bullet\Delta)=(\mathrm{aut}_{g}^{\mathcal W})^{\otimes2}\circ
\mathrm{Ad}_{f(e_1,f_1)}\circ \Delta\circ (\mathrm{aut}_{g}^{\mathcal W})^{-1}
=\mathrm{Ad}_{f(e_1,f_1)}\circ(\mathrm{aut}_{g}^{\mathcal W})^{\otimes2}\circ
 \Delta\circ (\mathrm{aut}_{g}^{\mathcal W})^{-1}
=f\bullet(g*\Delta) 
$$
for any $f\in \mathbf k[[u,v]]^\times$, $g\in \mathcal G$, $\Delta\in 
\mathrm{Cop}_{\mathcal C\operatorname{-alg}}(\hat{\mathcal W})$, where the second equality follows from 
$(\mathrm{aut}_{g}^{\mathcal W})^{\otimes2}(f(e_1,f_1))=f(e_1,f_1)$, 
which follows from the invariance by 
$(\mathrm{aut}_{g}^{\mathcal W})^{\otimes2}$ of $e_1$ and $f_1$. 

    (b) The proof is similar to (a), using that $g\mapsto \mathrm{aut}_g^{\mathcal W}$ and 
    $g\mapsto \mathrm{aut}_g^{V}$ define actions of  
    $\mathcal G$ on the algebras $\hat{\mathcal W}$ and $\hat V$, and that $\mathrm{aut}_g^{V}$ leaves 
    invariant $e_1$ and $f_1$.

(c) follows from the fact that the algebra morphism 
$(i_{\mathcal W_r,\mathcal V}^{\otimes2}\circ \mathrm{Ad}_{e_1})^{\otimes2} : 
\hat{\mathcal W}^{\hat\otimes2}\to\hat V$ intertwines the algebra automorphisms
$\mathrm{Ad}_{h(e_1,f_1)}$ on the source and on the target for any $h\in\mathbf k[[u,v]]^\times$ 
one the one hand, and  $(\mathrm{aut}_{g}^{\mathcal W})^{\otimes2}$
and $\mathrm{aut}_g^V$ for any $g\in\mathcal G$ on the other hand. 
\end{proof}

\begin{cor}\label{cor:HH}
(a) The quadruples $(\mathbf E,1+t^2\mathbf k[[t]],\bullet,*)$, $(\mathbf E,\mathbf k[[u,v]]^\times,\bullet,*)$, and $
(\mathbf H,\mathbf k[[u,v]]^\times,\bullet,*)$ are objects of $\mathcal G\operatorname{-}\mathbf{PSGA}$, where in the first quadruple
the action is the pull-back of the action of $\mathbf k[[u,v]]^\times$ on $\mathbf E$ by $\theta$.  

(b) The pairs $(id_{\mathbf E},\theta) : (\mathbf E,1+t^2\mathbf k[[t]],\bullet,*)\to(\mathbf E,\mathbf k[[u,v]]^\times,\bullet,*)$ and  
$(i_{\mathbf E\mathbf H},id) : (\mathbf E,\mathbf k[[u,v]]^\times,\bullet,*)\to
(\mathbf H,\mathbf k[[u,v]]^\times,\bullet,*)$ are morphisms in $\mathcal G\operatorname{-}\mathbf{PSGA}$.  
\end{cor}

\begin{proof}
(a) The statements on $(\mathbf E,\mathbf k[[u,v]]^\times,\bullet,*)$, and $
(\mathbf H,\mathbf k[[u,v]]^\times,\bullet,*)$ follow from Lem. \ref{lem28:1001}(b),(c) and 
Lem. \ref{lem:43:0401:TER}(a),(b). The statements on $(\mathbf E,1+t^2\mathbf k[[t]],\bullet,*)$ is 
an example of the fact that an object $(X,x_0,A,\bullet,*)$ in $\mathcal G\operatorname{-}\mathbf{PSGA}$
and a group morphism $\phi : B\to A$ give rise to a pulled back object 
$(X,x_0,B,\bullet,*)$ in $\mathcal G\operatorname{-}\mathbf{PSGA}$. 

(b) The first statement follows from the fact that in the situation evoked in (a), 
$(id_X,\phi) : (X,x_0,B,\bullet,*)\to (X,x_0,A,\bullet,*)$ is a morphism in 
$\mathcal G\operatorname{-}\mathbf{PSGA}$. The second statement follows from Lem. \ref{lem28:1001}(d) and 
Lem. \ref{lem:43:0401:TER}(c).
\end{proof}

\subsection{The diagram $\mathbf X_{\tilde{\mathcal G},\mathbf{E}}\stackrel{}{\to}
\mathbf X_{\mathcal G,\mathbf{E}'}\stackrel{}{\to}
\mathbf X_{\mathcal G,\mathbf{E}''}\stackrel{}{\to}
\mathbf X_{\mathcal G,\mathbf{H}''}$ in $\mathbf{PSGA}_{\mathcal G}$}\label{sect:diag:XXXX}

\begin{defn}
(a) Define $\mathbf X_{\tilde{\mathcal G},\mathbf{E}}$, $\mathbf X_{\tsup{\mathcal G},\mathbf E}$, $\mathbf X_{\tsup{\mathcal G},\mathbf H}$
to be the images of the objects  $(\mathbf E,1+t^2\mathbf k[[t]],\bullet,*)$, 
$(\mathbf E,\mathbf k[[u,v]]^\times,\bullet,*)$, $(\mathbf H,\mathbf k[[u,v]]^\times,\bullet,*)$
by the functor $\mathbf f : \mathcal G\operatorname{-}\mathbf{PSGA}\to \mathbf{PSGA}_{\mathcal G}$. 

(b) Define $\mathbf X_{\mathcal G,\mathbf{E}'}$, $
\mathbf X_{\mathcal G,\mathbf{E}''}$, $\mathbf X_{\mathcal G,\mathbf{H}''}$
to be the images of the objects  $(\mathbf E,1+t^2\mathbf k[[t]],\bullet,*)$, $(\mathbf E,\mathbf k[[u,v]]^\times,\bullet,*)$, and $
(\mathbf H,\mathbf k[[u,v]]^\times,\bullet,*)$ by the functor $\mathbf f \mathbf c \mathbf q : \mathcal G\operatorname{-}\mathbf{PSGA}\to \mathbf{PSGA}_{\mathcal G}$. 
\end{defn}

\begin{defn}\label{def:doubletildeG}
 $\tsup{\mathcal G}$ is the group $\mathcal G\times\mathbf k[[u,v]]^\times$. 
\end{defn}

\begin{defn}\label{def:previous}
    (a) $\odot$ is the action of the group $\tsup{\mathcal G}$ on the set $\mathrm{Cop}_{e_1}(\hat{\mathcal W})$ 
    obtained by combining the commuting actions $*$ and $\bullet$ of the groups $\mathcal G$ and $\mathbf k[[u,v]]^\times$ 
    on this set (see Lem. \ref{lem:43:0401:TER}(a)). 

    (b) $\odot$ is the action of the group $\tilde{\mathcal G}$ on the set $\mathrm{Cop}_{e_1}(\hat{\mathcal W})$ obtained from the action of (a) by pull-back by 
    the group morphism $id_{\mathcal G}\times\theta : \tilde{\mathcal G}=\mathcal G\times(1+t^2\mathbf k[[t]])\to\mathcal G\times\mathbf k[[u,v]]^\times=\tsup{\mathcal G}$. 
    
  (c)   $\odot$ is the action of the group $\tsup{\mathcal G}$ on the set $\mathrm{Hom}_{\mathcal C\operatorname{-alg}}(\hat{\mathcal W},\hat V)$ obtained by combining 
  the commuting actions $*$ and $\bullet$ of the groups $\mathcal G$ and $\mathbf k[[u,v]]^\times$  on this set (see Lem. \ref{lem:43:0401:TER}(b)). 
\end{defn}

\begin{lem}\label{lem:2:15:toto}
(a) One has
$$
\mathbf X_{\tilde{\mathcal G},\mathbf E}=((\tilde{\mathcal G},\mathrm{pr}_{\mathcal G}),\mathbf{E},\odot), \quad
\mathbf X_{\tsup{\mathcal G},\mathbf E}=((\tsup{\mathcal G},\mathrm{pr}_{\mathcal G}),\mathbf{E},\odot),\quad
\mathbf X_{\tsup{\mathcal G},\mathbf H}=((\tsup{\mathcal G},\mathrm{pr}_{\mathcal G}),\mathbf{H},\odot), 
$$
where the actions are respectively as in Def. \ref{def:previous}(b),(a),(c). 

(b) The morphisms $\mathbf X_{\tilde{\mathcal G},\mathbf E}\to\mathbf X_{\tsup{\mathcal G},\mathbf E}$ and $\mathbf X_{\tsup{\mathcal G},\mathbf E}\to
\mathbf X_{\tsup{\mathcal G},\mathbf H}$ in $\mathbf{PSGA}_{\mathcal G}$ given by the images by $\mathbf f$ of the morphisms 
$(\mathbf E,1+t^2\mathbf k[[t]],\bullet,*)\stackrel{(id_{\mathbf E},\theta)}\to (\mathbf E,\mathbf k[[u,v]]^\times,\bullet,*)$
and $(\mathbf E,\mathbf k[[u,v]]^\times,\bullet,*)\stackrel{(i_{\mathbf E\mathbf H},id)}{\to}(\mathbf H,\mathbf k[[u,v]]^\times,\bullet,*)$ in 
$\mathcal G\operatorname{-}\mathbf{PSGA}$ are respectively given by the pairs $(id_{\mathcal G}\times\theta,id_{\mathbf E})$ and $(id_{\mathcal G},i_{\mathbf E\mathbf H})$. 
\end{lem}

\begin{proof}
    Immediate. 
\end{proof}

\begin{defn}\label{def:E':E'':H''}
Define pointed sets as follows: 

(a) $\mathbf E''$ is the pair $(\mathbf k[[u,v]]^\times\backslash\mathrm{Cop}_{e_1}(\hat{\mathcal W}),[\!\![\hat\Delta^{\mathcal W}]\!\!])$, where 
the quotient is with respect to the action $\bullet$ of Lem. \ref{lem28:1001}(b) and $[\!\![\hat\Delta^{\mathcal W}]\!\!]$ is the image of 
$\hat\Delta^{\mathcal W}$ in this quotient; 

(b) $\mathbf E'$ is the pair $((1+t^2\mathbf k[[t]])\backslash\mathrm{Cop}_{e_1}(\hat{\mathcal W}),[\hat\Delta^{\mathcal W}])$, where 
the quotient is with respect to the pull-back by $\theta$ of the action $\bullet$ of Lem. \ref{lem28:1001}(b) and $[\hat\Delta^{\mathcal W}]$ is the image of 
$\hat\Delta^{\mathcal W}$ in this quotient; 

(c) $\mathbf H''$ is the pair $(\mathbf k[[u,v]]^\times\backslash\mathrm{Hom}_{\mathcal C\operatorname{-alg}}(\hat{\mathcal W},\hat V),\mathbf k[[u,v]]^\times\bullet\Delta^{\mathcal W}_{r,l}
)$, where  
the quotient is with respect to the action $\bullet$ of Lem. \ref{lem28:1001}(c) and 
$\mathbf k[[u,v]]^\times\bullet\Delta^{\mathcal W}_{r,l}
$ is the image of 
$\Delta^{\mathcal W}_{r,l}$ in this quotient. 
\end{defn}

\begin{lem}\label{lem:3101009:BIS}
One has 
$$
\mathbf X_{\mathcal G,\mathbf{E'}}=((\mathcal G,id),\mathbf{E'},\odot), \quad 
 \mathbf X_{\mathcal G,\mathbf{E''}}=((\mathcal G,id),\mathbf{E''},\odot), \quad 
\mathbf X_{\mathcal G,\mathbf{H''}}=((\mathcal G,id),\mathbf{H''},\odot).  
$$
where the actions of $\mathcal G$ are as follows: 
in the first case and second case, its actions on $(1+t^2\mathbf k[[t]])\backslash\mathrm{Cop}_{e_1}(\hat{\mathcal W})$
and $\mathbf k[[u,v]]^\times\backslash\mathrm{Cop}_{e_1}(\hat{\mathcal W})$ are induced by its action $*$ on $\mathrm{Cop}_{e_1}(\hat{\mathcal W})$; 
in the third case, its action on $\mathbf k[[u,v]]^\times\backslash\mathrm{Hom}_{\mathcal C\operatorname{-alg}}(\hat{\mathcal W},\hat V)$
is induced by its action $*$ on $\mathrm{Hom}_{\mathcal C\operatorname{-alg}}(\hat{\mathcal W},\hat V)$. 
\end{lem}

\begin{proof}
Direct verification. 
\end{proof}

\begin{lem}\label{lem:def:p_UV}
The diagram 
$$
\mathbf X_{\tilde{\mathcal G},\mathbf{E}}\stackrel{}{\to}
\mathbf X_{\mathcal G,\mathbf{E}'}\stackrel{}{\to}
\mathbf X_{\mathcal G,\mathbf{E}''}\stackrel{}{\to}
\mathbf X_{\mathcal G,\mathbf{H}''}
$$ 
in $\mathbf{PSGA}_{\mathcal G}$
extracted from the image by $\mathbf f$ of the commutative diagram 
arising from the diagram 
\begin{equation}\label{diag:E:E:H}
(\mathbf E,1+t^2\mathbf k[[t]],\bullet,*)\stackrel{(id_{\mathbf E},\theta)}{\to}(\mathbf E,\mathbf k[[u,v]]^\times,\bullet,*)
\stackrel{(i_{\mathbf E\mathbf H},id)}{\to}(\mathbf H,\mathbf k[[u,v]]^\times,\bullet,*)    
\end{equation}
in $\mathcal G\operatorname{-}\mathbf{PSGA}$ and the natural transformation 
$id\to \mathbf c\mathbf q$ (see Lem. \ref{lem:tranfo:nat}) is given as follows: the morphism  $\mathbf X_{\tilde{\mathcal G},\mathbf{E}}\stackrel{}{\to}
\mathbf X_{\mathcal G,\mathbf{E}'}$ is given by the pair formed by $\mathrm{pr}_{\mathcal G}$ and the morphism $p_{\mathbf E\mathbf E'} : 
\mathbf E\to\mathbf E'$ induced by the projection 
$\mathrm{Cop}_{e_1}(\hat{\mathcal W})\to(1+t^2\mathbf k[[t]])\backslash \mathrm{Cop}_{e_1}(\hat{\mathcal W})$, the morphism 
$\mathbf X_{\mathcal G,\mathbf{E}'}\stackrel{}{\to}\mathbf X_{\mathcal G,\mathbf{E}''}$ is given by the pair formed by 
$id_{\mathcal G}$ and by the morphism $p_{\mathbf E'\mathbf E''} : \mathbf E'\to\mathbf E''$ induced by the projection
$(1+t^2\mathbf k[[t]])\backslash \mathrm{Cop}_{e_1}(\hat{\mathcal W})
\to\mathbf k[[u,v]]^\times\backslash \mathrm{Cop}_{e_1}(\hat{\mathcal W})$, 
and the morphism $\mathbf X_{\mathcal G,\mathbf{E}''}\stackrel{}{\to}
\mathbf X_{\mathcal G,\mathbf{H}''}$ is given by the pair formed by $id_{\mathcal G}$ and by the morphism 
$i_{\mathbf E''\mathbf H''} : \mathbf E''\to\mathbf H''$ induced by
$\mathbf k[[u,v]]^\times\backslash \mathrm{Cop}_{e_1}(\hat{\mathcal W})\to
\mathbf k[[u,v]]^\times\backslash \mathrm{Hom}_{\mathcal C\operatorname{-alg}}(\hat{\mathcal W},\hat V)$. 
\end{lem}

\begin{proof}
The diagram \eqref{diag:E:E:H} in $\mathcal G\operatorname{-}\mathbf{PSGA}$, combined with the natural transformation 
$id\to cq$, gives rise to the commutative diagram 
$$
\xymatrix{(\mathbf E,1+t^2\mathbf k[[t]],\bullet,*)\ar^{(id_{\mathbf E},\theta)}[r]\ar[d]&
(\mathbf E,\mathbf k[[u,v]]^\times,\bullet,*)\ar^{(i_{\mathbf E\mathbf H},id)}[r]\ar[d]&
(\mathbf H,\mathbf k[[u,v]]^\times,\bullet,*)\ar[d]\\
 \mathbf c\mathbf q(\mathbf E,1+t^2\mathbf k[[t]],\bullet,*)\ar^{ \mathbf c\mathbf q(id_{\mathbf E},\theta)}[r]&
 \mathbf c\mathbf q(\mathbf E,\mathbf k[[u,v]]^\times,\bullet,*)\ar^{ \mathbf c\mathbf q(i_{\mathbf E\mathbf H},id)}[r]&
 \mathbf c\mathbf q(\mathbf H,\mathbf k[[u,v]]^\times,\bullet,*)
}
$$
in $\mathcal G\operatorname{-}\mathbf{PSGA}$. Applying $\mathbf f$ to it and using Lems. \ref{lem:2:15:toto} and \ref{lem:3101009:BIS}, one obtains
the following commutative diagram 
\begin{equation*}
\xymatrix{\mathbf X_{\tilde{\mathcal G},\mathbf{E}}\ar^{\mathbf f(id_{\mathbf E},\theta)}[r]\ar[d]&
\mathbf X_{\tsup{\mathcal G},\mathbf E}
\ar^{\mathbf f(i_{\mathbf E\mathbf H},id)}[r]\ar[d]&
\mathbf X_{\tsup{\mathcal G},\mathbf H}
\ar[d]\\
\mathbf X_{\mathcal G,\mathbf{E}'}\ar^{ \mathbf f\mathbf c\mathbf q(id_{\mathbf E},\theta)}[r]&
\mathbf X_{\mathcal G,\mathbf{E}''}\ar^{\mathbf f\mathbf c\mathbf q(i_{\mathbf E\mathbf H},id)}[r]&
\mathbf X_{\mathcal G,\mathbf{H}''}}
\end{equation*}
in $\mathbf{PSGA}_{\mathcal G}$. One checks that the various morphisms $\mathbf X_a\to\mathbf X_b$ are given by the announced formulas.    
\end{proof}

\section{A diagram in $\mathbf{PSGA}_{\mathcal G}$}\label{sect:3:1009}

Based on a discussion of coproducts in §\ref{sect:3:1:5786}, we introduce in §\ref{sect:2:3:titi} a diagram of 
pointed sets (Def. \ref{def:M:W:MW}). Together with the constructions of actions in §\ref{sect:3:2} and 
\ref{sect:3:4:5786}, this enables us to construct in §§\ref{sect:XXXX} and \ref{sect:TRANSITION} 
diagrams of pointed sets with group actions (Lems. \ref{lem:37:1009} and \ref{lem:3111009}).

\subsection{Sets of coproducts}\label{sect:3:1:5786}

For $\mathcal A$ a symmetric tensor category, define $\mathcal A$-alg-mod to be 
the category of algebra-modules in $\mathcal A$: an object in $\mathcal A$-alg-mod is a pair $(W,M)$ of an object $W$
of $\mathcal A$-alg and a left module $M$ over $W$. This is again a STC and 
there is a diagram 
$$
\mathcal A\leftarrow\mathcal A\operatorname{-alg-mod}\to\mathcal A\operatorname{-alg}
$$
of forgetful tensor functors given by $M\leftmapsto (W,M)\mapsto W$. 
In particular, an object $(W,M)$ in $\mathcal A$-alg-mod gives rise to a
diagram of sets
\begin{equation}\label{diag:com:gal}
\mathrm{Cop}_{\mathcal A}(M)
\leftarrow\mathrm{Cop}_{\mathcal A\operatorname{-alg-mod}}(W,M)
\to\mathrm{Cop}_{\mathcal A\operatorname{-alg}}(W). 
\end{equation}

\begin{lem}\label{lem:cat:act}
Let $\mathcal A$ is a STC and $X$ be an object in $\mathcal A$. 

(a) The group $\mathrm{Aut}_{\mathcal A}(X)$ acts on $\mathrm{Cop}_{\mathcal A}(X)$ by $g\diamond c:=g^{\otimes 2}\circ c\circ g^{-1}$. 

(b) A functor $F : \mathcal A\to\mathcal B$ of STCs gives rise to a map $\mathrm{Cop}_{\mathcal A}(X)\to\mathrm{Cop}_{\mathcal B}(FX)$ and to a 
group morphism $\mathrm{Aut}_{\mathcal A}(X)\to\mathrm{Aut}_{\mathcal B}(FX)$, which are compatible with the actions on both sides. 
\end{lem}

\begin{proof} Both statements are immediate.  
\end{proof}

\subsection{The diagram $\mathbf M\stackrel{p_{\mathbf M}}{\leftarrow}\mathbf{WM}\stackrel{p_{\mathbf W}}{\to} 
\mathbf{W}$ in $\mathbf{PS}$}
\label{sect:2:3:titi}\label{sect:a:diag:in:PS}



It follows from definitions that  
\begin{equation}\label{elts:2909}
\hat\Delta^{\mathcal M}\in\mathrm{Cop}_{\mathcal C}(\hat{\mathcal M}), \quad 
\hat\Delta^{\mathcal W}\in\mathrm{Cop}_{\mathcal C\operatorname{-alg}}(\hat{\mathcal W}), \quad 
(\hat\Delta^{\mathcal W},\hat\Delta^{\mathcal M})\in\mathrm{Cop}_{\mathcal C\operatorname{-alg-mod}}
(\hat{\mathcal W},\hat{\mathcal M}).     
\end{equation}
The diagram 
\begin{equation}\label{maps:2909}
\mathrm{Cop}_{\mathcal C}(\hat{\mathcal M})\leftarrow 
\mathrm{Cop}_{\mathcal C\operatorname{-alg-mod}}(\hat{\mathcal W},\hat{\mathcal M})\to 
\mathrm{Cop}_{\mathcal C\operatorname{-alg}}(\hat{\mathcal W})
\end{equation}
arising from \eqref{diag:com:gal} is compatible with the elements \eqref{elts:2909}, therefore
it induces a diagram in $\mathbf{PS}$.

\begin{defn}\label{def:M:W:MW}
     $\mathbf M\stackrel{p_{\mathbf M}}{\leftarrow}\mathbf{WM}\stackrel{p_{\mathbf W}}{\to} 
\mathbf{W}$ is the diagram in $\mathbf{PS}$ arising from \eqref{elts:2909}, \eqref{maps:2909}. 
\end{defn}

\subsection{The action $\bullet_\Gamma$}\label{sect:3:2}

In \cite{EF2}, Lem. 1.18, one defines a group morphism 
$$
\mathsf G^\DR(\mathbf k)\to 
\mathrm{Aut}_{\mathcal C\operatorname{-alg-mod}}(\hat{\mathcal W},\hat{\mathcal M}), \quad 
(\mu,g)\mapsto
(\ ^\Gamma\!\mathrm{aut}^{\mathcal W,(1)}_{(\mu,g)},\ ^\Gamma\!\mathrm{aut}^{\mathcal M,(10)}_{(\mu,g)}). 
$$
Post-composing it with the group morphism 
$\mathrm{Aut}_{\mathcal C\operatorname{-alg-mod}}(\hat{\mathcal W},\hat{\mathcal M})
\to \mathrm{Aut}_{\mathcal C}(\hat{\mathcal M})$, one obtains a group morphism 
$$
\mathsf G^\DR(\mathbf k)\to 
\mathrm{Aut}_{\mathcal C\operatorname{-mod}}(\hat{\mathcal M}), \quad 
(\mu,g)\mapsto
\ ^\Gamma\!\mathrm{aut}^{\mathcal M,(10)}_{(\mu,g)}. 
$$
The subgroup $\mathsf{Stab}(\hat\Delta^{\mathcal M})$ of $\mathsf G^\DR(\mathbf k)=(\mathrm{exp}(\mathfrak{lie}_{\{0,1\}}^\wedge),\circledast) 
\rtimes \mathbf k^\times$ is defined in \cite{EF2}, Def. 2.20 
to be the stabilizer group of $(\mathsf G^\DR(\mathbf k),\mathrm{Cop}_{\hat{\mathcal C}_{\mathbf k}}(\hat{\mathcal M}),\hat\Delta^{\mathcal M},
(\mu,g)\mapsto ^\Gamma\mathrm{aut}^{\mathcal M,(10)}_{(\mu,g)}\diamond-)$. 

Pre-composing the latter group morphism with the group inclusion $\mathcal G\hookrightarrow \mathsf G^\DR(\mathbf k)$, 
$g\mapsto (1,g)$, one obtains a group morphism 
\begin{equation}\label{def;mor:Gamma:2312:BIS}
\mathcal G\to \mathrm{Aut}_{\mathcal C}(\hat{\mathcal M}), 
\quad g\mapsto\ ^\Gamma\!\mathrm{aut}^{\mathcal M,(10)}_{(1,g)}.     
\end{equation}
It follows from {\it loc. cit.}, (1.6.13) that for $g\in\mathcal G$, one has 
\begin{equation}\label{id:1512:BIS}
 ^\Gamma\mathrm{aut}_g^{\mathcal M,(10)}
=\ell_{\Gamma_g(-e_1)^{-1}}\circ \mathrm{aut}_g^{\mathcal M},  
\end{equation}
where for $a\in\hat{\mathcal V}$, $\ell_a$ is the endomorphism $m\mapsto am$ of $\hat{\mathcal M}$, and where 
$\mathrm{aut}_g^{\mathcal M}$ is the automorphism of $\hat{\mathcal M}$ such that for any $v\in\hat{\mathcal V}$, 
one has $\mathrm{aut}_g^{\mathcal M}(v\cdot 1_{\mathcal M})=\mathrm{aut}_g^{\mathcal V}(v)g\cdot 1_{\mathcal M}$
(where $\mathrm{aut}_g^{\mathcal V}$ is the automorphism of $\hat{\mathcal V}$ induced by 
the automorphism $\mathrm{aut}_g^{\mathcal V}$ of $\mathfrak{lie}_{\{0,1\}}^\wedge$).

\begin{defn}\label{def32:1009}
The action $\bullet_\Gamma$ of $\mathcal G$ on $\mathrm{Cop}_{\mathcal C}(\hat{\mathcal M})$ is the pull-back by the morphism  
\eqref{def;mor:Gamma:2312:BIS} of the action of 
$\mathrm{Aut}_{\mathcal C}(\hat{\mathcal M})$ on this set arising from Lem. 
\ref{lem:cat:act} (where $(\mathcal A,X)=(\mathcal C,\hat{\mathcal M})$).  
\end{defn}

\subsection{The actions $\bullet_{\mathcal{WM}},\bullet_{\mathcal W}$ and $\bullet_{\mathcal M}$}
\label{sect:3:4:5786}

Let $\hat{\mathcal V}_1\subset \hat{\mathcal V}$ be the subset of elements with constant term equal to 1. Define a law 
$\circledast$ on $\hat{\mathcal V}_1$ by extending \eqref{FORMULA:1912}; then $(\hat{\mathcal V}_1,\circledast)$ is a group. 
For $g\in \hat{\mathcal V}_1$, define $\mathrm{aut}^{\mathcal V}_g\in\mathrm{Aut}_{\mathcal C\operatorname{-alg}}(\hat{\mathcal V})$ by 
$e_1\mapsto e_1$ and $e_0\mapsto g\cdot e_0\cdot g^{-1}$. Then $\mathrm{aut}^{\mathcal V}_g$ induces an automorphism of the subalgebra 
$\hat{\mathcal W}$, which will be denoted $\mathrm{aut}^{\mathcal W}_g$. There is also an automorphism $\mathrm{aut}^{\mathcal M}_g$
of the $\mathbf k$-module $\hat{\mathcal M}$ such that $v\cdot 1_{\mathcal M}\mapsto \mathrm{aut}^{\mathcal V}_g(v)g
\cdot 1_{\mathcal M}$ for any $v\in\hat{\mathcal V}$. 
It follows from \cite{EF2}, Lem. 1.11, that the map 
\begin{equation}\label{gp:mor:mathcalV:2312:BIS}
(\hat{\mathcal V}_1,\circledast)\to
\mathrm{Aut}_{\mathcal C\operatorname{-alg-mod}}(\hat{\mathcal W},\hat{\mathcal M}), 
\quad 
g\mapsto (\mathrm{aut}^{\mathcal W}_g,\mathrm{aut}^{\mathcal M}_g)   
\end{equation}
is a group morphism (this map is denoted $g\mapsto (\mathrm{aut}^{\mathcal W,\DR,(1)}_g,\mathrm{aut}^{\mathcal M,\DR,(10)}_g)$ in 
{\it loc. cit.}), and $g\mapsto \mathrm{aut}_g^{\mathcal V}$, $g\mapsto \mathrm{aut}_g^{\mathcal W}$
are extensions of the maps from Lem. \ref{lem:2:8:jan25}. 

\begin{lem}\label{lem:can:BIS}
There is a unique group morphism $\tilde{\mathcal G}\to(\hat{\mathcal V}_1,\circledast)$, 
whose restriction to $\mathcal G$ coincides with inclusion of this group in $(\hat{\mathcal V}_1,\circledast)$ 
(see \cite{EF2}, \S1.6.3, end of third paragraph) and whose restriction to $1+t^2\mathbf k[[t]]$ is given by $h(t)\mapsto 
h(e_1)$.  
\end{lem}

\begin{proof}
The canonical injection and $h(t)\mapsto h(e_1)$ are respectively group morphisms from 
$\mathcal G$ and $1+t^2\mathbf k[[t]]$ to $(\hat{\mathcal V}_1,\circledast)$. The images of these
morphisms commute: indeed, for $g\in \mathcal G$ and $h\in 1+t^2\mathbf k[[t]]$, one has 
$g\circledast h(e_1)=g\cdot h(e_1)=h(e_1)\circledast g$ (see \eqref{FORMULA:1912}). It follows that there is a 
unique group morphism 
\begin{equation}\label{pre:pre:mor:BIS}
\tilde{\mathcal G}=\mathcal G\times(1+t^2\mathbf k[[t]])\to(\hat{\mathcal V}_1,\circledast),
\end{equation}
whose restrictions to the two factors of the source are respectively the canonical inclusion and 
the morphism $h(t)\mapsto h(e_1)$. 
\end{proof}

\begin{defn}\label{def:gp:mor:BIS}
The group morphism 
\begin{equation}\label{mor:G:Aut:alg:mod}
  \tilde{\mathcal G}\to\mathrm{Aut}_{\mathcal C\operatorname{-alg-mod}}(\hat{\mathcal W},\hat{\mathcal M})  
\end{equation}
obtained as the composition of the morphism $(\hat{\mathcal V}_1,\circledast)\to
\mathrm{Aut}_{\mathcal C\operatorname{-alg-mod}}(\hat{\mathcal W},\hat{\mathcal M})$
from \eqref{gp:mor:mathcalV:2312:BIS} with the morphism  $\tilde{\mathcal G}\to (\hat{\mathcal V}_1,\circledast)$ (see Lem. \ref{lem:can:BIS})
is denoted $(g,h)\mapsto(\mathrm{aut}_{(g,h)}^{\mathcal W},\mathrm{aut}_{(g,h)}^{\mathcal M})$. 
Its composition with the group morphisms from its target to $\mathrm{Aut}_{\mathcal C\operatorname{-alg}}(\hat{\mathcal W})$
and $\mathrm{Aut}_{\mathcal C}(\hat{\mathcal M})$ are then the group morphisms 
\begin{equation}\label{mor:G:Aut:alg}
 \tilde{\mathcal G}\to\mathrm{Aut}_{\mathcal C\operatorname{-alg}}(\hat{\mathcal W}),\quad   
 (g,h)\mapsto\mathrm{aut}_{(g,h)}^{\mathcal W}=\mathrm{Ad}_{h(e_1)} \circ \mathrm{aut}_{g}^{\mathcal W}
\end{equation}
and
\begin{equation}\label{mor:G:Aut:mod}
 \tilde{\mathcal G}\to\mathrm{Aut}_{\mathcal C}(\hat{\mathcal M}),\quad   (g,h)\mapsto\mathrm{aut}_{(g,h)}^{\mathcal M}
 =\ell_{h(e_1)} \circ \mathrm{aut}_{g}^{\mathcal M}. 
\end{equation}
\end{defn}

\begin{defn}\label{def:35:1009}
The action $\bullet_{\mathcal{WM}}$ (resp. $\bullet_{\mathcal{W}}$, $\bullet_{\mathcal{M}}$) of $\tilde{\mathcal G}$ on 
$\mathrm{Cop}_{\mathcal C\operatorname{-alg-mod}}(\hat{\mathcal W},\hat{\mathcal M})$ 
(resp. $\mathrm{Cop}_{\mathcal C\operatorname{-alg}}(\hat{\mathcal W})$, $\mathrm{Cop}_{\mathcal C}(\hat{\mathcal M})$)
is the pull-back by the morphism  
\eqref{mor:G:Aut:alg:mod} of the action of 
$\mathrm{Aut}_{\mathcal C\operatorname{-alg-mod}}(\hat{\mathcal W},\hat{\mathcal M})$ (resp. 
$\mathrm{Aut}_{\mathcal C\operatorname{-alg}}(\hat{\mathcal W})$, $\mathrm{Aut}_{\mathcal C}(\hat{\mathcal M})$) on this set arising from Lem. 
\ref{lem:cat:act}, with $(\mathcal A,X)$ equal to $(\mathcal C\operatorname{-alg-mod},(\hat{\mathcal W},\hat{\mathcal M}))$
(resp. $(\mathcal C\operatorname{-alg},\hat{\mathcal W})$, $(\mathcal C,\hat{\mathcal M})$). 
\end{defn}

\subsection{The diagram $\mathbf X_{\mathcal G,\mathbf M}\stackrel{(\sigma,id)}{\to}
\mathbf X_{\tilde{\mathcal G},\mathbf M}\stackrel{(id,p_{\mathbf M})}{\leftarrow}\mathbf X_{\tilde{\mathcal G},\mathbf{WM}}\stackrel{(id,p_{\mathbf W})}{\to}
\mathbf X_{\tilde{\mathcal G},\mathbf{W}}$ in $\mathbf{PSGA}_{\mathcal G}$}\label{sect:XXXX}

\begin{lem}\label{lem:mor:A'':A:BIS}
The pair $(\sigma,\mathrm{id})$ gives rise to a morphism 
$\mathbf X_{\mathcal G,\mathbf M}\to 
\mathbf X_{\tilde{\mathcal G},\mathbf M}$ in $\mathbf{PSGA}_{\mathcal G}$, where the objects are as in Prop. \ref{prop:1:9}(b).  
\end{lem}

\begin{proof}
For $g\in\mathcal G$, one has 
$$
 ^\Gamma\mathrm{aut}_g^{\mathcal M,(10)}
=\ell_{\Gamma_g(-e_1)^{-1}}\circ \mathrm{aut}_g^{\mathcal M}
=\mathrm{aut}_{(g,\Gamma_g(-t)^{-1})}^{\mathcal M}
=\mathrm{aut}_{\sigma(g)}^{\mathcal M}, 
$$
where the first (resp. second, third) equality follows from \eqref{id:1512:BIS} (resp. 
\eqref{mor:G:Aut:mod}, Def. \ref{defn:sigma:BIS}). 
It follows that for $g\in\mathcal G$, the permutations 
$^\Gamma\mathrm{aut}_g^{\mathcal M,(10)}\diamond-$ and $\mathrm{aut}_{\sigma(g)}^{\mathcal M}\diamond-$
of $\mathrm{Cop}_{\mathcal C}(\hat{\mathcal M})$ coincide, therefore that 
the actions $\bullet_\Gamma$ and $\bullet_{\mathcal M}$ of $\mathcal G$ and $\tilde{\mathcal G}$ 
on $\mathrm{Cop}_{\mathcal C}(\hat{\mathcal M})$ are related by 
$\sigma(g)\bullet_{\mathcal M}x=x\bullet_{\Gamma}x$ for any $(g,x)\in \mathcal G\times 
\mathrm{Cop}_{\mathcal C}(\hat{\mathcal M})$, which implies the claim. 
\end{proof}

\begin{lem}\label{lem:37:1009}
The pairs $(id,p_{\mathbf M})$ and $(id,p_{\mathbf W})$ give rise to morphisms $\mathbf X_{\tilde{\mathcal G},\mathbf{WM}}\to
\mathbf X_{\tilde{\mathcal G},\mathbf{M}}$ and $\mathbf X_{\tilde{\mathcal G},\mathbf{WM}}\to\mathbf X_{\tilde{\mathcal G},\mathbf{W}}$ 
in $\mathbf{PSGA}_{\mathcal G}$,  where the objects are as in Prop. \ref{prop:1:9}(b).  
\end{lem}

\begin{proof}
Applying Lem. \ref{lem:cat:act} to the forgetful functor $\mathcal C\operatorname{-alg-mod}\to\mathcal C\operatorname{-alg}$
and $X:=(\hat{\mathcal W},\hat{\mathcal M})$ one sees that 
$(\mathrm{Aut}_{\mathcal C\operatorname{-alg-mod}}(\hat{\mathcal W},\hat{\mathcal M}),
\mathrm{Cop}_{\mathcal C\operatorname{-alg-mod}}(\hat{\mathcal W},\hat{\mathcal M}),\bullet_{\mathcal W\mathcal M})
\to(\mathrm{Aut}_{\mathcal C\operatorname{-alg}}(\hat{\mathcal W}),
\mathrm{Cop}_{\mathcal C\operatorname{-alg}}(\hat{\mathcal W}),\bullet_{\mathcal W})$
is a morphism of sets with group actions. One shows similarly the same statement on
$$
(\mathrm{Aut}_{\mathcal C\operatorname{-alg-mod}}(\hat{\mathcal W},\hat{\mathcal M}),
\mathrm{Cop}_{\mathcal C\operatorname{-alg-mod}}(\hat{\mathcal W},\hat{\mathcal M}),\bullet_{\mathcal W\mathcal M})
\to(\mathrm{Aut}_{\mathcal C}(\hat{\mathcal M}),\mathrm{Cop}_{\mathcal C}(\hat{\mathcal M}),\bullet_{\mathcal M}).
$$
The statement then follows from the commutativity of the triangles in 
$$
\xymatrix{
\mathrm{Aut}_{\mathcal C}(\hat{\mathcal M})& 
\mathrm{Aut}_{\mathcal C\operatorname{-alg-mod}}(\hat{\mathcal W},\hat{\mathcal M})\ar[r]\ar[l]& 
\mathrm{Aut}_{\mathcal C\operatorname{-alg}}(\hat{\mathcal W})\\ & \tilde{\mathcal G}\ar[u]\ar[ur]\ar[ul]& }
$$
which follows from Defs. \ref{mor:G:Aut:alg}, \ref{mor:G:Aut:mod}. 
\end{proof}

\subsection{The morphism $(id,i_{\mathbf E\mathbf W}) : \mathbf X_{\tilde{\mathcal G},\mathbf E}\to\mathbf X_{\tilde{\mathcal G},\mathbf W}$ in 
$\mathbf{PSGA}_{\mathcal G}$}\label{sect:TRANSITION}

\begin{defn}\label{def:i:EW}
     We denote by $i_{\mathbf E,\mathbf W} : \mathbf E\to\mathbf W$ 
the morphism in $\mathbf{PS}$  induced by the canonical injection  $\mathrm{Cop}_{e_1}(\hat{\mathcal W})\subset \mathrm{Cop}_{{\mathcal C}\operatorname{-alg}}(\hat{\mathcal W})$. 
\end{defn}

\begin{lem}\label{lem:3111009}
(a) The canonical inclusion $\mathrm{Cop}_{e_1}(\hat{\mathcal W})\hookrightarrow
\mathrm{Cop}_{{\mathcal C}\operatorname{-alg}}(\hat{\mathcal W})$ is equivariant with respect to the actions of $\tilde{\mathcal G}$ on its source 
by $\odot$ and on its target by $\bullet_{\mathcal W}$. 

(b) $(id,i_{\mathbf E\mathbf W})$ is a morphism 
$\mathbf X_{\tilde{\mathcal G},\mathbf E}\to\mathbf X_{\tilde{\mathcal G},\mathbf W}$ in $\mathbf{PSGA}_{\mathcal G}$. 
\end{lem}

\begin{proof}
Let $g\in\mathcal G$, $h\in1+t^2\mathbf k[[t]]$ and $\Delta\in\mathrm{Cop}_{e_1}(\hat{\mathcal W})$. Then 
$$
(g,1)\odot\Delta=(\mathrm{aut}_{(g,1)}^{\mathcal W})^{\otimes2}\circ\Delta\circ (\mathrm{aut}_{(g,1)}^{\mathcal W})^{-1}=
\mathrm{aut}_{(g,1)}^{\mathcal W}\diamond\Delta=(g,1)\bullet_{\mathcal W}\Delta
$$
and
$$
(1,h)\bullet_{\mathcal W}\Delta=\mathrm{aut}_{(1,h)}^{\mathcal W}\diamond\Delta
=\mathrm{Ad}_{h(e_1)}\diamond\Delta
=\mathrm{Ad}_{h(e_1)}^{\otimes2}\circ\Delta\circ(\mathrm{Ad}_{h(e_1)})^{-1}
=\mathrm{Ad}_{h(e_1)h(f_1)/h(e_1+f_1)}\circ\Delta=(1,h)\odot\Delta, 
$$
where the second equality follows from 
$\mathrm{aut}_{(1,h)}^{\mathcal W}=(e_1\mapsto e_1,e_0\mapsto h(e_1)e_0h(e_1)^{-1})
=\mathrm{Ad}_{h(e_1)}$ and the fourth equality follows from $\Delta(e_1)=e_1+f_1$. 
This implies (a). (b) follows. 
\end{proof}

\section{Isomorphisms in $\mathbf{Gp}_{\mathcal G}$}\label{sect:4:korea}

This section is the final stage of the proof of Thm. \ref{thm:main}. It consists essentially in 
proving that the the group morphisms arising from the diagrams of pointed sets introduced 
in §§\ref{sect:2} and \ref{sect:3:1009} are isomorphisms. Such isomorphisms are established in 
§\ref{sect:4:1:5786} (based on a involved algebraic argument), §\ref{sect:4:2:5786} 
(based on an argument involving structures of free rank one modules over algebras), 
§\ref{sect:4:3:5786} (based on results of \cite{EF4}), §\ref{sect:4:4:5786} (based on 
local injectivity arguments), §\ref{sect:4:5:5786} (based on computations of centralizers). 
In §\ref{sect:4:6:5786}, we establish various 
algebraic results, one of which being related with Hochschild homology, and in §\ref{sect:G:H:2402:BIS}, 
we  establish the surjectivity of a group morphism; based on these results, we prove the 
remaining group isomorphisms in §§\ref{sect:H:I:2402:BIS} and \ref{sect:4:9:5786}.

\subsection{Isomorphism status of $\underline{\mathrm{Stab}}\mathbf X_{\mathcal G,\mathbf M}\to
\underline{\mathrm{Stab}}\mathbf X_{\tilde{\mathcal G},\mathbf M}$}\label{sect:4:1:5786}

\begin{lem}
Set $\mathcal G(\hat{\mathcal M}):=\{m\in\hat{\mathcal M}|m$ is a generator of $\hat{\mathcal M}$ over $\hat{\mathcal W}$ and 
$\hat\Delta^{\mathcal M}(m)=m\otimes m\}$. Then 
\begin{equation}\label{FIRSTSTEP:BIS}
    \mathrm{Stab}_{\tilde{\mathcal G}}(\hat\Delta^{\mathcal M}) \subset \{(g,h)\in \tilde{\mathcal G}=\mathcal G\times(1+t^2\mathbf k[[t]]) | 
    h(-e_1)g\cdot 1_{\mathcal M}\in\mathcal G(\hat{\mathcal M})\},  
\end{equation}
where the left-hand side is the subgroup of $\tilde{\mathcal G}$ associated with the action of this group on the pointed set $\mathbf M$ 
by $\bullet_\Gamma$ (see Defs. \ref{def:M:W:MW} and \ref{def32:1009}).  
\end{lem}

\begin{proof}
If $(g,h)\in\mathcal G\times(1+t^2\mathbf k[[t]])$, then 
\begin{equation}\label{EQN:BIS}
    \mathrm{aut}^{\mathcal M}_{(g,h)}(1_{\mathcal M})=\mathrm{aut}^{\mathcal M}_{(1,h)}\circ 
    \mathrm{aut}^{\mathcal M}_{(g,1)}(1_{\mathcal M})=\ell_{h(e_1)}(g\cdot 1_{\mathcal M})=h(e_1)g\cdot 1_{\mathcal M}. 
\end{equation} 
As $g \in \mathrm{exp}(\mathfrak{lie}_{\{0,1\}}^\wedge)$ and $h\in 1+t^2\mathbf k[[t]]$, one has $h(e_1)g\in\hat{\mathcal V}^\times$ therefore 
$h(e_1)g\cdot 1_{\mathcal M}$ 
is a generator of $\hat{\mathcal M}$ as a $\hat{\mathcal W}$-module. Assume now that 
$(g,h)\in\mathrm{Stab}_{\tilde{\mathcal G}}(\hat\Delta^{\mathcal M})$. 
Then $(\mathrm{aut}^{\mathcal M}_{(g,h)})^{\otimes 2}\circ \hat\Delta^{\mathcal M}=\hat\Delta^{\mathcal M}\circ
\mathrm{aut}^{\mathcal M}_{(g,h)}$. Applying this equality to $1_{\mathcal M}$, using $\hat\Delta^{\mathcal M}(1_{\mathcal M})
=1_{\mathcal M}^{\otimes 2}$ and \eqref{EQN:BIS}, one obtains $\hat\Delta^{\mathcal M}(h(e_1)g\cdot 1_{\mathcal M})=
(h(e_1)g\cdot 1_{\mathcal M})^{\otimes 2}$. 
\end{proof}

Recall that $\mathcal W$ is the graded $\mathbf k$-subalgebra of $\mathcal V$ defined by 
$\mathcal W:=\mathbf k\oplus\mathcal V e_1$ (see \cite{EF1}, \S1.1). Then $\mathcal W_+:=\mathcal V e_1$
is a positively graded $\mathbf k$-subalgebra of $\mathcal W$ (without unit). 

\begin{lem}\label{lemma:basics:BIS}
Define graded $\mathbf k$-submodules  $U_0,U_1,\mathcal A$ of $\mathcal W_+$ by 
$U_1:=\mathbf k[e_1]e_1$, $U_0:=\mathbf k[e_0]e_1$, 
$\mathcal A:=(\mathcal Ve_0\mathcal Ve_1\mathcal V+\mathcal Ve_1\mathcal Ve_0\mathcal V)e_1$.

(a) There holds the direct sum decomposition $\mathcal W_+=U_1\oplus U_0\oplus\mathcal A$. Let $U:=U_0\oplus U_1$ 
and for $x\in \mathcal W_+$, let $x=x_U+x_{\mathcal A}$ be the 
decomposition corresponding to $\mathcal W_+=U\oplus\mathcal A$. 

(b) Let $(u_n)_{n\geq 1}$ be the family of elements of $\mathcal V$ defined by $\sum_{n\geq 1}u_nt^n=
\mathrm{log}(1-\sum_{i\geq 1}t^ie_0^{i-1}e_1)$, then $u_n\in\mathcal W_+$ for any $n\geq1$. One has $(u_1)_U=-e_1$ and 
$(u_n)_U=-e_0^{n-1}e_1-e_1^n/n$ for any $n>1$. 

(c) One has $[\mathcal W_+,\mathcal W_+]\subset\mathcal A$. 

(d) For $v_0\in \hat U_0$, $v_1\in \hat U_1$, $a\in\hat{\mathcal A}$ one has $\mathrm{exp}(v_0+v_1+a)\in 
\mathrm{exp}(v_1)+v_0+\hat{\mathcal A}$ (equality in $\hat{\mathcal W}$; recall that the degree completion functor is 
denoted $X\mapsto\hat X$). 
\end{lem}

\begin{proof}
(a) Recall that $\mathcal V$ is a $\mathbb Z_{\geq 0}^2$-graded algebra, the grading being such that $e_0$, $e_1$ are of degrees $(1,0)$ and 
$(0,1)$. The partition $\mathbb Z_{\geq 0}^2=(\{0\}\times\mathbb Z_{\geq0})\sqcup(\mathbb Z_{>0}\times\{0\})\sqcup\mathbb Z_{>0}^2$ gives 
rise to the decomposition $\mathcal V=\mathbf k[e_1]\oplus e_0\mathbf k[e_0]\oplus(\mathcal Ve_0\mathcal Ve_1\mathcal V\oplus
\mathcal Ve_1\mathcal Ve_0\mathcal V)$. Applying the linear isomorphism $\mathcal V\to\mathcal W_+$ given by $x\mapsto xe_1$ gives the claimed 
decomposition. 

(b) For $n\geq1$, one has $u_n=-\sum_{(k,n_1,\ldots,n_k)|k,n_1,\ldots,n_k \geq 1, 
n_1+\ldots+n_k=n} e_0^{n_1-1}e_1\cdots e_0^{n_k-1}e_1/k$ therefore $u_n\in\mathcal W_+$. 
Assume that $n>1$. If the index $(k,n_1,\ldots,n_k)$ is such that $k\geq 2$ and $n_i>1$ for some index $i\in[\!\![1,k]\!\!]$, then 
the corresponding component in the expression of $u_n$ belongs to $\mathcal A$ ; exactly two indices do not fall in this class, 
namely $(n,{\underbrace{1,1,\cdots,1}_{n}})$ and $(1,n)$, for which the contributions are respectively $-e_1^n/n$ and $-e_0^{n-1}e_1$, 
which belong to $U$, giving the announced expression of $(u_n)_U$.  
The announced expression of $(u_1)_U$  follows from $u_1=-e_1$.

(c) Set $\mathcal B:=U_0\oplus\mathcal A$, there is a decomposition $\mathcal W_+=U_1\oplus\mathcal B$ and an inclusion
$\mathcal W_+\supset \mathcal A$. Moreover, $\mathcal W_+$ is a 
$\mathbb Z_{\geq 0}\times\mathbb Z_{>0}$-graded algebra, the decomposition corresponds to the partition $\mathbb Z_{\geq 0}\times
\mathbb Z_{>0}=(\{0\}\times\mathbb Z_{>0})\sqcup\mathbb Z_{>0}^2$ and inclusion corresponds to the inclusion 
$\mathbb Z_{\geq 0}\times\mathbb Z_{>0}\supset\mathbb Z_{>0}\times\mathbb Z_{>1}$. Then $(\{0\}\times\mathbb Z_{>0}) + \mathbb Z_{>0}^2 
\subset \mathbb Z_{>0}\times\mathbb Z_{>1}$ and $\mathbb Z_{>0}^2+\mathbb Z_{>0}^2 \subset \mathbb Z_{>0}\times\mathbb Z_{>1}$, which implies 
$U_1\cdot\mathcal B\subset\mathcal A$, $\mathcal B\cdot U_1\subset\mathcal A$ and $\mathcal B\cdot\mathcal B\subset\mathcal A$, therefore 
$[\mathcal B,U_1]\subset\mathcal A$ and $[\mathcal B,\mathcal B]\subset\mathcal A$. The result follows from the conjunction of these inclusions and 
of the equality $[U_1,U_1]=0$.   

(d) Set $b:=v_0+a$, so that $b \in\hat{\mathcal B}$ (see (c)). One has 
\begin{equation}\label{1503:1007:BIS}
\mathrm{exp}(v_0+v_1+a)=\mathrm{exp}(v_1+b)=\sum_{n\geq 0}\sum_{s_0,\ldots,s_n\geq 0}v_1^{s_0}b\cdot\cdots\cdot bv_1^{s_n}/(s_0+...+s_n+n)!
\end{equation}
Let $\tilde U_1:=\mathbf k\oplus \hat U_1$, then it follows from (c) that $\tilde U_1\cdot\hat{\mathcal B}\subset\hat{\mathcal B}$, 
$\hat{\mathcal B}\cdot\tilde U_1\subset\hat{\mathcal B}$. Since $b^s\in \tilde U_1$ for any $s\geq 0$, the summand in the right-hand side of 
\eqref{1503:1007:BIS} associated to $n$ is contained in $\hat{\mathcal B}\cdot\cdots\cdot 
\hat{\mathcal B}$ ($n$ times), which by (c) is contained in $\hat{\mathcal A}$ if $n\geq 2$. The contribution for $n=1$ is 
$\sum_{s_0,s_1\geq 0} v_1^{s_0}bv_1^{s_1}/(s_0+s_1+1)!$. The inclusions $\hat U_1\cdot\hat{\mathcal B}\subset\hat{\mathcal A}$, 
$\hat{\mathcal B}\cdot\hat U_1\subset\hat{\mathcal A}$ from (c) imply that the contributions for $n=1$ and $(s_0,s_1)\neq (0,0)$ belong 
to $\hat{\mathcal A}$. It follows that $\mathrm{exp}(v_0+v_1+a)\in\mathrm{exp}(v_1)+b+\hat{\mathcal A}=\mathrm{exp}(v_1)+v_0
+\hat{\mathcal A}$, as claimed. 
\end{proof}

\begin{lem}
There holds the inclusion 
\begin{equation}\label{SECONDSTEP:BIS}
\{(g,h)\in \mathcal G\times(1+t^2\mathbf k[[t]]) | h(e_1)g\cdot 1_{\mathcal M}\in\mathcal G(\hat{\mathcal M})\}
\subset\mathrm{im}(\sigma),  
\end{equation}
where $\sigma$ as in Def. \ref{defn:sigma:BIS}. 
\end{lem}

\begin{proof}
Let $(g,h) \in \mathcal G\times(1+t^2\mathbf k[[t]])$ be such that $h(e_1)g\cdot 1_{\mathcal M}\in
\mathcal G(\hat{\mathcal M})$. The map $\hat{\mathcal W}\to\hat{\mathcal M}$, $x\mapsto x\cdot 1_{\mathcal M}$, sets up a 
bijection $\mathcal G(\hat{\mathcal W})\stackrel{\sim}{\to}\mathcal G(\hat{\mathcal M})$, where $\mathcal G(\hat{\mathcal W})$ is 
the group of group-like elements of the topological Hopf algebra $(\hat{\mathcal W},\hat\Delta^{\mathcal W})$. The condition 
$h(e_1)g\cdot 1_{\mathcal M} \in \mathcal G(\hat{\mathcal M})$ is therefore equivalent to $\pi(h(e_1)g)\in\mathcal G(\hat{\mathcal W})$, 
where $\pi : \hat{\mathcal V}\to\hat{\mathcal W}$ is the composition of the projection map 
$\hat{\mathcal V}\to\hat{\mathcal M}$, $x\mapsto x\cdot 1_{\mathcal M}$ with the inverse of the isomorphism 
$\hat{\mathcal W}\to\hat{\mathcal M}$, $x\mapsto x\cdot 1_{\mathcal M}$. Let $g_0,g_1\in\hat{\mathcal V}$ be the elements such that 
$g=1+g_0e_0+g_1e_1$, then $\pi(h(e_1)g)=h(e_1)(1+g_1e_1)$, so the latter condition is equivalent to 
$h(e_1)(1+g_1e_1)\in\mathcal G(\hat{\mathcal W})$. 

Recall that $(\hat{\mathcal W},\hat\Delta^{\mathcal W})$ is isomorphic to the enveloping algebra of the topologically free Lie algebra 
over the family $(u_n)_{n\geq 1}$ of generators defined in Lem. \ref{lemma:basics:BIS}(b) (see \cite{R}, \S 2.3.6, also \cite{EF1}, Rem. 1.1). 
It follows that the set $\mathcal P(\hat{\mathcal W})$ of primitive elements of $(\hat{\mathcal W},\hat\Delta^{\mathcal W})$ is the 
topological Lie subalgebra of $\hat{\mathcal W}_+$ generated by $(u_n)_{n\geq 1}$, and that 
$\mathcal G(\hat{\mathcal W})=\mathrm{exp}(\mathcal P(\hat{\mathcal W}))$. In particular, 
there exists a Lie series $P(u_1,u_2,\ldots)$ in $\mathcal P(\hat{\mathcal W})$ such that 
\begin{equation}\label{EQUALITy:BIS}
h(e_1)(1+g_1e_1)=\mathrm{exp}(P(u_1,u_2,\ldots)). 
\end{equation}

Let us denote by $(\lambda_i)_{i\geq 1}$ the coefficients such that the linear part of 
$P(u_1,u_2,\ldots)$ is $\sum_{i\geq 1}\lambda_i u_i$. Then 
\begin{equation}\label{2226:0607:BIS}
 P(u_1,u_2,\ldots)\in
\sum_{i\geq 1}\lambda_i u_i+[\mathcal P(\mathcal W),\mathcal P(\mathcal W)]^\wedge
\subset \sum_{i\geq 1}\lambda_i u_i+[\mathcal W_+,\mathcal W_+]^\wedge \subset 
\sum_{i\geq 1}\lambda_i u_i+\hat{\mathcal A},    
\end{equation}
where $(-)^\wedge$ means the degree completion of a graded vector space, 
the first inclusion follows from $\mathcal P(\hat{\mathcal W})\subset\hat{\mathcal W}_+$ and the last inclusion follows from 
Lem. \ref{lemma:basics:BIS}(c). \eqref{2226:0607:BIS} and Lem. \ref{lemma:basics:BIS}(b) then 
imply 
$$
 P(u_1,u_2,\ldots)_U=-\sum_{i\geq 1}\lambda_ie_1^i/i-\sum_{i\geq 2}\lambda_ie_0^{i-1}e_1,    
$$
where the left-hand side has the meaning explained in Lem. \ref{lemma:basics:BIS}(a).
The latter equation enables one to apply Lem. \ref{lemma:basics:BIS} (d) with $u_1:=-\sum_{i\geq 1}\lambda_ie_1^i/i$, 
$u_0=-\sum_{i\geq 2}\lambda_ie_0^{i-1}e_1$, $a:=P(u_1,u_2,\ldots)-u_0-u_1$, which yields
$$
\mathrm{exp}(P(u_1,u_2,\ldots))\in\mathrm{exp}(-\sum_{i\geq 1}\lambda_ie_1^i/i)-\sum_{i\geq 2}\lambda_ie_0^{i-1}e_1+\hat{\mathcal A}
$$
which together with \eqref{EQUALITy:BIS} implies
\begin{equation}\label{above:id:BIS}
h(e_1)(1+g_1e_1)\in\mathrm{exp}(-\sum_{i\geq 1}\lambda_ie_1^i/i)-\sum_{i\geq 2}\lambda_ie_0^{i-1}e_1+\hat{\mathcal A}. 
\end{equation}
Equip $\mathbf k[[e_1]]$ with the coproduct for which $e_1$ is primitive and $\hat{\mathcal V}$ with the coproduct 
$\hat\Delta^{\mathcal V}$. There is unique topological algebra morphism $\varpi : \hat{\mathcal V}\to\mathbf k[[e_1]]$ induced 
by $e_0\mapsto 0$, $e_1\mapsto e_1$, which is also a topological Hopf algebra morphism. Since $g\in\mathrm{exp}(\mathfrak{lie}_{\{0,1\}}^\wedge)$, 
it follows that $\varpi(g)\in\mathcal G(\mathbf k[[e_1]])$, which implies that there exists $\alpha\in\mathbf k$ such that 
$\varpi(g)=\mathrm{exp}(\alpha e_1)$, and since $(g|e_1)=0$ since $g\in\mathcal G$, one has $\alpha=0$, therefore 
$\varpi(g)=1$. The image of \eqref{above:id:BIS} by $\varpi$ yields 
$h(e_1)=\mathrm{exp}(-\sum_{i\geq 1}\lambda_ie_1^i/i)$ (equality in $\mathbf k[[e_1]]$) therefore 
\begin{equation}\label{2242:0607:BIS}
h(t)=\mathrm{exp}(-\sum_{i\geq 1}\lambda_it^i/i). 
\end{equation}
(equality in $\mathbf k[[t]]$). On the other hand, \eqref{above:id:BIS} also implies that for $i\geq2$, one has $\lambda_i=-(g|e_0^{i-1}e_1)$. 
Combining this with \eqref{2242:0607:BIS}, one obtains  
$$
h(t)=\mathrm{exp}(-\lambda_1t+\sum_{i\geq 2}(g|e_0^{i-1}e_1)t^i/i). 
$$
Together with $h\in 1+t^2\mathbf k[[t]]$, this implies $\lambda_1=0$, therefore  
$$
h(t)=\mathrm{exp}(\sum_{i\geq 2}(g|e_0^{i-1}e_1)t^i/i)=\tilde\Gamma_g^{-1}(-t). 
$$
Therefore $(g,h)=\sigma(g)$, hence $(g,h)\in\mathrm{im}(\sigma)$. 
\end{proof}

\begin{prop}\label{prop:incl:1510}
One has the inclusion 
$$
\mathrm{Stab}_{\tilde{\mathcal G}}(\hat\Delta^{\mathcal M})\subset \mathrm{im}(\sigma)
$$
of subgroups of $\tilde{\mathcal G}=\mathcal G\times (1+t^2\mathbf k[[t]])$. 
\end{prop}

\begin{proof}
This is a consequence of \eqref{FIRSTSTEP:BIS} and \eqref{SECONDSTEP:BIS}.     
\end{proof}

\begin{lem}\label{lem:general:nonsense}
    Assume that $(X,x)$ is a pointed set, $A$ is a group, $\bullet$ is an action of $\mathcal G\times A$ on $X$, and $\tilde\alpha : \mathcal G\to A$
    is a group morphism, and let $\alpha:=(id,\tilde\alpha) : \mathcal G\to\mathcal G\times A$. 

    Let $\tilde\bullet$ be the action of 
    $\mathcal G$ on $X$ obtained by pulling back $\bullet$ by the group morphism $\alpha$; then $(\alpha,id) : ((\mathcal G,id),(X,x),\tilde\bullet)
    \to ((\mathcal G\times A,\mathrm{pr}_{\mathcal G}),(X,x),\bullet)$ is a morphism in  $\mathbf{PSGA}_{\mathcal G}$.

If $\mathrm{Stab}_{\mathcal G\times A}(x)\subset \mathrm{im}(\alpha)$, then the morphism 
\begin{equation}\label{mor:stab:stab}
\underline{\mathrm{Stab}}((\mathcal G,id),(X,x),\tilde\bullet)
\to\underline{\mathrm{Stab}}((\mathcal G\times A,\mathrm{pr}_{\mathcal G}),(X,x),\bullet)
\end{equation}
is an isomorphism is $\mathbf{Gp}_{\mathcal G}$.     
\end{lem}

\begin{proof}
The group morphism $\alpha$ restricts and corestricts to a morphism 
$\underline{\alpha} : \mathrm{Stab}_{\mathcal G}(x)\to \mathrm{Stab}_{\mathcal G\times A}(x)$. 
The morphism \eqref{mor:stab:stab} is given by the diagram 
$$
\xymatrix{\mathrm{Stab}_{\mathcal G}(x)\ar_{i_{\mathrm{Stab}_{\mathcal G}(x),\mathcal G}}[rd]\ar^{\underline{\alpha}}[rr]& & 
\mathrm{Stab}_{\mathcal G\times A}(x)
\ar^{\mathrm{pr}_{\mathcal G}\circ i_{\mathrm{Stab}_{\mathcal G\times A}(x),\mathcal G\times A}}[ld]
\\ & \mathcal G& }
$$
Since $\alpha$ is injective, so is $\underline{\alpha}$. The assumption 
$\mathrm{Stab}_{\mathcal G\times A}(x)\subset \mathrm{im}(\alpha)$ means that for any 
$(g,a)\in\mathcal G\times A$, $(g,a)\bullet x=x$ implies $a=\tilde\alpha(g)$. 
Therefore 
$$
\mathrm{Stab}_{\mathcal G\times A}(x)=\{(g,a)\in\mathcal G\times A|(g,a)\bullet x=x\}
=\{(g,a)\in\mathcal G\times A|a=\tilde\alpha(g) \text{ and }(g,\tilde\alpha(g))\bullet x=x\}
=\alpha(\mathrm{Stab}_{\mathcal G}(x))
$$
which implies the surjectivity of $\underline{\alpha}$. It follows that $\underline{\alpha}$
is an isomorphism, which implies the statement. 
\end{proof}

\begin{cor}\label{cor:4:6:nov}
    The morphism $\underline{\mathrm{Stab}}\mathbf X_{\mathcal G,\mathbf M}\to
\underline{\mathrm{Stab}}\mathbf X_{\tilde{\mathcal G},\mathbf M}$ is an isomorphism in $\mathbf{Gp}_{\mathcal G}$. 
\end{cor}

\begin{proof}
This follows from Prop. \ref{prop:incl:1510} and Lem. \ref{lem:general:nonsense}, applied to the morphism
$(\sigma,id) : \mathbf X_{\mathcal G,\mathbf M}\to\mathbf X_{\tilde{\mathcal G},\mathbf M}$. 
\end{proof}

\subsection{Isomorphism status of $\underline{\mathrm{Stab}}\mathbf X_{\tilde{\mathcal G},\mathbf M}\to\underline{\mathrm{Stab}}\mathbf X_{\tilde{\mathcal G},\mathbf{WM}}$}\label{sect:4:2:5786}

\begin{defn}\label{def:local:injectivity:1018}
    A morphism $f : (X,x)\to (X',x')$ in $\mathbf{PS}$ will be called {\it locally injective} if $f^{-1}(x')=\{x\}$. 
\end{defn}

\begin{lem}\label{lem:3010:1129}
    Let $((\Gamma,h),(X,x),\bullet)$ and $((\Gamma,h),(X',x'),\bullet')$
     be objects of $\mathbf{PSGA}_{\mathcal G}$ and let $f : (X,x)\to(X',x')$ be a morphism in 
     $\mathbf{PS}$ such that $(id,f) : ((\Gamma,h),(X,x),\bullet)\to ((\Gamma,h),(X',x'),\bullet')$ is a 
     morphism in $\mathbf{PSGA}_{\mathcal G}$ and that $f$ is locally injective. Then 
     \begin{equation}\label{eq:subgps:Gamma}
      \mathrm{Stab}_\Gamma(X,x)=\mathrm{Stab}_\Gamma(X',x')   
     \end{equation}
   (equality of subgroups of $\Gamma$) and the morphism 
     $$ \underline{\mathrm{Stab}}((\Gamma,h),(X,x),\bullet)\to\underline{\mathrm{Stab}}((\Gamma,h),(X',x'),\bullet')$$  
     is an isomorphism in $\mathbf{Gp}_{\mathcal G}$. 
\end{lem}

\begin{proof}
The groups $\mathrm{Stab}_\Gamma(X,x)$ and $\mathrm{Stab}_\Gamma(X',x')$ are subgroups of $\Gamma$. 
For $\gamma\in\Gamma$, 
$$
(\gamma\in\mathrm{Stab}_\Gamma(X',x'))\iff(\gamma\bullet'x'=x')\iff(f(\gamma\bullet x)=f(x))
\iff (\gamma\bullet x=x)\iff(\gamma\in\mathrm{Stab}_\Gamma(X,x)),
$$ 
where the third equivalence uses the assumption on $f$. This implies \eqref{eq:subgps:Gamma}, 
which by the definition of $\underline{\mathrm{Stab}}$ (see Lem. \ref{lem:def:stab:ulstab}) implies the statement. 
\end{proof}

\begin{lem}\label{lem:3010:1129a}
    The map $p_{\mathbf M} : \mathbf{WM}\to\mathbf M$ is locally injective. 
\end{lem}

\begin{proof}
Recall that this map is given by $\mathrm{Cop}_{\mathcal C\operatorname{-alg-mod}}
(\hat{\mathcal W},\hat{\mathcal M})\to \mathrm{Cop}_{\mathcal C}(\hat{\mathcal M})$,  
$(\Delta_W,\Delta_M)\mapsto \Delta_M$. Let then $\Delta_W\in\mathrm{Cop}_{\mathcal C\operatorname{-alg}}(\hat{\mathcal W})$
such that $(\Delta_W,\hat\Delta^{\mathcal M})$ belongs to $\mathrm{Cop}_{\mathcal C\operatorname{-alg-mod}}
(\hat{\mathcal W},\hat{\mathcal M})$. Then for any $w\in\hat{\mathcal W}$, one has
$$
\hat\Delta^{\mathcal M}(w\cdot 1_{\mathcal M})=\Delta_W(w)\cdot \hat\Delta^{\mathcal M}(1_{\mathcal M})
=\Delta_W(w)\cdot 1_{\mathcal M}^{\otimes2},  
$$
which since $(\hat\Delta^{\mathcal W},\hat\Delta^{\mathcal M})$ belongs to $\mathrm{Cop}_{\mathcal C\operatorname{-alg-mod}}
(\hat{\mathcal W},\hat{\mathcal M})$ implies
$$
\hat\Delta^{\mathcal W}(w)\cdot 1_{\mathcal M}^{\otimes2}=\Delta_W(w)\cdot 1_{\mathcal M}^{\otimes2},
$$
which since $\hat{\mathcal M}^{\hat\otimes2}$ is freely generated by $1_{\mathcal M}^{\otimes2}$ over $\hat{\mathcal W}^{\hat\otimes2}$
implies $\Delta_W=\hat\Delta^{\mathcal W}$. 
\end{proof}

\begin{prop}\label{prop:4:10:nov}
One has 
\begin{equation}\label{eq:stab:G:0411}
   \mathrm{Stab}_{\tilde{\mathcal G}}((\hat\Delta^{\mathcal W},\hat\Delta^{\mathcal M}))
\subset\mathrm{Stab}_{\tilde{\mathcal G}}(\hat\Delta^{\mathcal M}) 
\end{equation}
(equality of subgroups of $\tilde{\mathcal G}$), where the left-hand side is 
the subgroup of $\tilde{\mathcal G}$ associated with the action of this group on the pointed set $\mathbf W\mathbf M$ by 
$\bullet_{\mathcal W,\mathcal M}$ (see Defs. \ref{def:M:W:MW} and \ref{def32:1009}) and 
the morphism $\underline{\mathrm{Stab}}\mathbf X_{\tilde{\mathcal G},\mathbf M}\to\underline{\mathrm{Stab}}\mathbf X_{\tilde{\mathcal G},\mathbf{WM}}$
    is an isomorphism in $\mathbf{Gp}_{\mathcal G}$. 
\end{prop}

\begin{proof}
  Follows from Lems. \ref{lem:3010:1129} and \ref{lem:3010:1129a}. 
\end{proof}

\subsection{Isomorphism status of $\underline{\mathrm{Stab}}\mathbf X_{\tilde{\mathcal G},\mathbf{WM}}\to
\underline{\mathrm{Stab}}\mathbf X_{\tilde{\mathcal G},\mathbf W}$}\label{sect:4:3:5786}

\begin{lem}\label{lem:4:11:411}
There is an inclusion $\mathrm{Stab}_{\tilde{\mathcal G}}(\hat\Delta^{\mathcal M})
\subset\mathrm{Stab}_{\tilde{\mathcal G}}(\hat\Delta^{\mathcal W})$ of subgroups of $\tilde{\mathcal G}$, where
the right-hand side is the subgroup of $\tilde{\mathcal G}$ associated with the action of this group on the pointed set $\mathbf W$ by 
$\bullet_{\mathcal W}$ (see Defs. \ref{def:M:W:MW} and \ref{def32:1009}).
\end{lem}

\begin{proof}
The morphism $\mathbf{WM}\to\mathbf{W}$ is $\tilde{\mathcal G}$-equivariant, which induces an inclusion 
$\mathrm{Stab}_{\tilde{\mathcal G}}((\hat\Delta^{\mathcal W},\hat\Delta^{\mathcal M}))
\subset\mathrm{Stab}_{{\tilde{\mathcal G}}}(\hat\Delta^{\mathcal W})$ of subgroups of ${\tilde{\mathcal G}}$. 
The result then follows from \eqref{eq:stab:G:0411}. 
\end{proof}

The assignment $\mathbf k\mapsto {\tilde{\mathcal G}}=\mathcal G\times (1+t^2\mathbf k\mathbf [[t]],\cdot)$ is a functor from the 
category of $\mathbb Q$-algebras to that of groups. It corresponds to a prounipotent $\mathbb Q$-group scheme and 
$\mathbf k\mapsto \mathrm{Stab}_{{\tilde{\mathcal G}}}(\hat\Delta^{\mathcal M})$, 
$\mathbf k\mapsto \mathrm{Stab}_{{\tilde{\mathcal G}}}(\hat\Delta^{\mathcal W})$ 
are the subfunctors corresponding to prounipotent subgroup schemes. Lem. \ref{lem:4:11:411} gives rise to a sequence of inclusions 
of $\mathbb Q$-Lie algebras 
\begin{equation}\label{inclusions:1312:BIS}
\mathfrak{stab}_{\tilde{\mathfrak G}}(\hat\Delta^{\mathcal M})\subset \mathfrak{stab}_{\tilde{\mathfrak G}}(\hat\Delta^{\mathcal W})
\subset \tilde{\mathfrak G},  
\end{equation}
where $\tilde{\mathfrak G}$ (resp. $\mathfrak{stab}_{\tilde{\mathfrak G}}(\hat\Delta^{\mathcal M})$, 
$\mathfrak{stab}_{\tilde{\mathfrak G}}(\hat\Delta^{\mathcal W})$) is the Lie algebra of the group scheme corresponding to $\mathbf k\mapsto {\tilde{\mathcal G}}$ 
(resp. $\mathbf k\mapsto \mathrm{Stab}_{{\tilde{\mathcal G}}}(\hat\Delta^{\mathcal M})$, 
$\mathbf k\mapsto \mathrm{Stab}_{{\tilde{\mathcal G}}}(\hat\Delta^{\mathcal W})$). One has 
$\tilde{\mathfrak G}=\mathfrak G\times t^2\mathbb Q[[t]]$, where $\mathfrak G$ is as in Lem. \ref{lem06:1109} and 
$t^2\mathbb Q[[t]]$ is viewed as a graded abelian Lie algebra, with $t$ of degree 1. 

For $\mathcal X\in\{\mathcal W,\mathcal M\}$, denote by $\hat{\mathcal X}_{\mathbb Q}$ the 
specialization of $\hat{\mathcal X}$ for $\mathbf k=\mathbb Q$.

\begin{lem}\label{lem:comp:stabs:BIS}
For $\mathcal X\in\{\mathcal W,\mathcal M\}$, one has 
 \begin{align*}
\mathfrak{stab}_{\tilde{\mathfrak G}}(\hat\Delta^{\mathcal X})=\{(x,h)\in \tilde{\mathfrak G}|
(\mathrm{der}_{(x,h)}^{\mathcal X}\otimes \mathrm{id}+\mathrm{id}\otimes \mathrm{der}_{(x,h)}^{\mathcal X})\circ\hat\Delta^{\mathcal 
X}=\hat\Delta^{\mathcal X}\circ \mathrm{der}_{(x,h)}^{\mathcal X}
\text{(equality in $\mathrm{Cop}_{\hat{\mathcal C}_{\mathbb Q}}((F)\hat{\mathcal X}_{\mathbb Q})$)}\}
\end{align*}
where $\mathrm{der}_{(x,h)}^{\mathcal W}:=\mathrm{der}_{x}^{\mathcal W}+\mathrm{ad}_{h(e_1)}$, 
$\mathrm{der}_{(x,h)}^{\mathcal M}:=\mathrm{der}_{x}^{\mathcal M}+\ell_{h(e_1)}$ and 
for $x\in\mathfrak{lie}_{\{e_0,e_1\},\mathbb Q}^\wedge$, we denote by $\mathrm{der}_x^{\mathcal W}$, 
$\mathrm{der}_x^{\mathcal M}$ the elements 
$\mathrm{der}_{(0,x)}^{\mathcal W,(1),\DR}$, $\mathrm{der}_{(0,x)}^{\mathcal M,(10),\DR}$ from \cite{EF2}, Lem. 3.11; 
$F$ is the forgetful functor $\hat{\mathcal C}_{\mathbb Q}\operatorname{-alg}\to\hat{\mathcal C}_{\mathbb Q}$
and the parenthesis means this it is applied only if $\mathcal X=\mathcal W$.  
\end{lem}

\begin{proof} Assume $\mathcal X=\mathcal M$. The group morphism $\mathcal G\times(1+t^2\mathbf k[[t]])
={\tilde{\mathcal G}}\to \mathrm{Aut}_{\hat{\mathcal C}_{\mathbf k}}(\hat{\mathcal M})$ given by $(g,h)\mapsto \mathrm{aut}^{\mathcal M}_{(g,h)}$ 
satisfies $\mathrm{aut}^{\mathcal M}_{(g,h)}=\mathrm{aut}^{\mathcal M,(10)}_g\circ \ell_{h(e_1)}$. The underlying Lie algebra 
morphism is then $\tilde{\mathfrak G}=\mathfrak G\times t^2\mathbb Q[[t]]\to 
\mathrm{End}_{\hat{\mathcal C}_{\mathbb Q}}(\hat{\mathcal M}_{\mathbb Q})$ given by 
$(x,h)\mapsto \mathrm{der}_x^{\mathcal M}+\ell_{h(e_1)}=\mathrm{der}_{(x,h)}^{\mathcal M}$. 
This enables to identify the right-hand side with the stabilizer Lie algebra of 
$\hat\Delta^{\mathcal M}$ for the action of $\tilde{\mathfrak G}$ on 
$\mathrm{Cop}_{\hat{\mathcal C}_{\mathbb Q}}(\hat{\mathcal M}_{\mathbb Q})$, which is 
$\mathfrak{stab}_{\tilde{\mathfrak G}}(\hat\Delta^{\mathcal M})$. 

Assume now $\mathcal X=\mathcal W$. The group morphism $\mathcal G\times(1+t^2\mathbf k[[t]])={\tilde{\mathcal G}}\to 
\mathrm{Aut}_{\hat{\mathcal C}_{\mathbf k}\operatorname{-alg}}(\hat{\mathcal W})$ given by 
$(g,h)\mapsto \mathrm{aut}^{\mathcal W}_{(g,h)}$ is given by 
$\mathrm{aut}^{\mathcal W}_{(g,h)}=\mathrm{aut}^{\mathcal W}_g\circ \mathrm{Ad}_{h(e_1)}$. The underlying Lie algebra 
morphism $\tilde{\mathfrak G}\to \mathrm{End}_{\hat{\mathcal C}_{\mathbb Q}}(\hat{\mathcal W}_{\mathbb Q})\cap 
\mathrm{Der}_{\mathbb Q\operatorname{-alg}}(\hat{\mathcal W}_{\mathbb Q})$ is therefore given by 
$(x,h)\mapsto \mathrm{der}_x^{\mathcal W}+\mathrm{ad}_{h(e_1)}=\mathrm{der}_{(x,h)}^{\mathcal W}$. 
This enables to identify the right-hand side with the stabilizer Lie algebra of 
$\hat\Delta^{\mathcal W}$ for the action of $\tilde{\mathfrak G}$ on 
$\mathrm{Cop}_{\hat{\mathcal C}_{\mathbb Q}\operatorname{-alg}}(\hat{\mathcal W}_{\mathbb Q})$, 
which is $\mathfrak{stab}_{\tilde{\mathfrak G}}(\hat\Delta^{\mathcal X})$. 
\end{proof}

Let $\mathrm{Der}_{\hat\Delta^{\mathcal W}}(\hat{\mathcal W}_{\mathbb Q},\hat{\mathcal W}_{\mathbb Q}^{\hat\otimes 2})$ be the set of
morphisms $a\in \mathrm{Cop}_{\hat{\mathcal C}_{\mathbb Q}}(F\hat{\mathcal W}_{\mathbb Q})$ 
satisfying the identity $a(ww')=a(w)\hat\Delta^{\mathcal W}(w')+\hat\Delta^{\mathcal W}(w)a(w')$ for $w,w'\in\hat{\mathcal W}_{\mathbb Q}$.

\begin{lem}\label{toy:lemma:BIS}
Let $a,b\in \mathrm{Der}_{\hat\Delta^{\mathcal W}}(\hat{\mathcal W}_{\mathbb Q},\hat{\mathcal W}_{\mathbb Q}^{\hat\otimes 2})$ and let $\alpha,\beta\in \mathrm{Cop}_{\hat{\mathcal C}_{\mathbb Q}}(\hat{\mathcal M}_{\mathbb Q})$ 
satisfy the identities $\alpha(wm)=a(w)\hat\Delta^{\mathcal M}(m)+\hat\Delta^{\mathcal W}(w)\alpha(m)$, 
$\beta(wm)=b(w)\hat\Delta^{\mathcal M}(m)+\hat\Delta^{\mathcal W}(w)\beta(m)$ for $w\in\hat{\mathcal W}_{\mathbb Q}$, 
$m\in\hat{\mathcal M}_{\mathbb Q}$. Then $\alpha=\beta$ iff $a=b$ and $\alpha(1_{\mathcal M})=\beta(1_{\mathcal M})$. 
\end{lem}

\begin{proof} Assume that $\alpha=\beta$. Then $\alpha(1_{\mathcal M})=\beta(1_{\mathcal M})$.
Moreover, for any $w\in\hat{\mathcal W}_{\mathbb Q}$, one has 
\begin{align*}
& a(w)1_{\mathcal M}^{\otimes2}+\hat\Delta^{\mathcal W}(w)\alpha(1_{\mathcal M})
=a(w)\hat\Delta^{\mathcal M}(1_{\mathcal M})+\hat\Delta^{\mathcal W}(w)\alpha(1_{\mathcal M})
=\alpha(w\cdot 1_{\mathcal M})=\beta(w\cdot 1_{\mathcal M})
\\ & =b(w)\hat\Delta^{\mathcal M}(1_{\mathcal M})+\hat\Delta^{\mathcal W}(w)\beta(1_{\mathcal M})
=b(w)1_{\mathcal M}^{\otimes2}+\hat\Delta^{\mathcal W}(w)\beta(1_{\mathcal M})
\end{align*}
therefore $a(w)1_{\mathcal M}^{\otimes2}+\hat\Delta^{\mathcal W}(w)\alpha(1_{\mathcal M})
=b(w)1_{\mathcal M}^{\otimes2}+\hat\Delta^{\mathcal W}(w)\beta(1_{\mathcal M})$, which, since 
$\alpha(1_{\mathcal M})=\beta(1_{\mathcal M})$, implies 
$a(w)1_{\mathcal M}^{\otimes2}=b(w)1_{\mathcal M}^{\otimes2}$. Since $\hat{\mathcal M}_{\mathbb Q}^{\hat\otimes 2}$
is a free rank one module over $\hat{\mathcal W}_{\mathbb Q}^{\hat\otimes 2}$ with basis $1_{\mathcal M}^{\otimes2}$, the 
latter equality implies $a(w)=b(w)$, therefore $a=b$. 

Assume that $a=b$ and $\alpha(1_{\mathcal M})=\beta(1_{\mathcal M})$. Then for any $m\in \hat{\mathcal M}_{\mathbb Q}$, there exists 
$w\in\hat{\mathcal W}_{\mathbb Q}$
such that $m=w\cdot 1_{\mathcal M}$. Then $\alpha(m)=\alpha(w\cdot 1_{\mathcal M})
=a(w)\hat\Delta^{\mathcal M}(1_{\mathcal M})+\hat\Delta^{\mathcal W}(w)\alpha(1_{\mathcal M})
=a(w)1_{\mathcal M}^{\otimes2}+\hat\Delta^{\mathcal W}(w)\alpha(1_{\mathcal M})
=b(w)1_{\mathcal M}^{\otimes2}+\hat\Delta^{\mathcal W}(w)\beta(1_{\mathcal M})
=b(w)\hat\Delta^{\mathcal M}(1_{\mathcal M})+\hat\Delta^{\mathcal W}(w)\beta(1_{\mathcal M})
=\beta(w\cdot 1_{\mathcal M})=\beta(m)$, so $\alpha=\beta$. 
\end{proof}

\begin{rem}
If $a,b\in \mathrm{Cop}_{\hat{\mathcal C}_{\mathbb Q}}(\hat{\mathcal W}_{\mathbb Q})$
and $\alpha,\beta\in \mathrm{Cop}_{\hat{\mathcal C}_{\mathbb Q}}(\hat{\mathcal M}_{\mathbb Q})$
satisfy the relation from Lem. \ref{toy:lemma:BIS}, then one can prove using the freeness of rank 1 of $\hat{\mathcal M}_{\mathbb Q}$
over $\hat{\mathcal W}_{\mathbb Q}$ that 
$a,b\in \mathrm{Der}_{\hat\Delta^{\mathcal W}}(\hat{\mathcal W}_{\mathbb Q},\hat{\mathcal W}_{\mathbb Q}^{\hat\otimes 2})$. 
\end{rem}

Define $\mathcal P(\hat{\mathcal M}_{\mathbb Q})$ as the kernel of the morphism $\hat{\mathcal M}_{\mathbb Q}
\to \hat{\mathcal M}_{\mathbb Q}^{\hat\otimes2}$ in $\hat{\mathcal C}_{\mathbb Q}$ 
given by hat $\hat\Delta^{\mathcal M}-id \otimes 1_{\mathcal M}-1_{\mathcal M}\otimes id$, where $id \otimes 1_{\mathcal M}$ 
(resp. $1_{\mathcal M}\otimes id$) is the composed morphism 
$\hat{\mathcal M}_{\mathbb Q} \simeq \mathbb Q\otimes \hat{\mathcal M}_{\mathbb Q}
\stackrel{1_{\mathcal M}\otimes id}{\to}\hat{\mathcal M}_{\mathbb Q}^{\hat\otimes2}$ 
(resp. $\hat{\mathcal M}_{\mathbb Q} \simeq \hat{\mathcal M}_{\mathbb Q}\otimes\mathbb Q 
\stackrel{id\otimes 1_{\mathcal M}}{\to}\hat{\mathcal M}_{\mathbb Q}^{\hat\otimes2}$), so  
$\mathcal P(\hat{\mathcal M}_{\mathbb Q})$ can be identified with $\{x\in \mathcal P(\hat{\mathcal M}_{\mathbb Q})|
\hat\Delta^{\mathcal M}(x)=x\otimes 1_{\mathcal M}+1_{\mathcal M}\otimes x\}$. 

\begin{lem}\label{lem:3:11:3105:BIS}
One has 
$$
\{(x,h)\in \tilde{\mathfrak G}|(x+h(e_1))\cdot 1_{\mathcal M}
\in\mathcal P(\hat{\mathcal M}_{\mathbb Q})\}\cap \mathfrak{stab}_{\tilde{\mathfrak G}}(\hat\Delta^{\mathcal W})=
\mathfrak{stab}_{\tilde{\mathfrak G}}(\hat\Delta^{\mathcal M})
$$
\end{lem}

\begin{proof} Let $(x,h)\in \tilde{\mathfrak G}$. Set $a_{(x,h)}:=\hat\Delta^{\mathcal W}_{\mathbb Q}\circ\mathrm{der}^{\mathcal W}_{(x,h)}$, 
$b_{(x,h)}:=(\mathrm{der}^{\mathcal W}_{(x,h)}\otimes id+id\otimes \mathrm{der}^{\mathcal W}_{(x,h)})\circ
\hat\Delta^{\mathcal W}_{\mathbb Q}$,
$\alpha_{(x,h)}:=\hat\Delta^{\mathcal M}_{\mathbb Q}\circ\mathrm{der}^{\mathcal M}_{(x,h)}$, 
$\beta_{(x,h)}:=(\mathrm{der}^{\mathcal M}_{(x,h)}\otimes id+id\otimes \mathrm{der}^{\mathcal M}_{(x,h)})\circ
\hat\Delta^{\mathcal M}_{\mathbb Q}$.

By Lemma \ref{lem:comp:stabs:BIS}, 
$\mathfrak{stab}_{\tilde{\mathfrak G}}(\hat\Delta^{\mathcal M})=\{(x,h)\in\tilde{\mathfrak G}|\alpha_{(x,h)}=\beta_{(x,h)}\}$ and 
$\mathfrak{stab}_{\tilde{\mathfrak G}}(\hat\Delta^{\mathcal W})=\{(x,h)\in\tilde{\mathfrak G}|a_{(x,h)}=b_{(x,h)}\}$. Then 
$\alpha_{(x,h)}(1_{\mathcal M})=\hat\Delta^{\mathcal M}((x+h(e_1))\cdot 1_{\mathcal M})$, $\beta_{(x,h)}(1_{\mathcal M})=
(x+h(e_1))\cdot 1_{\mathcal M}\otimes 1_{\mathcal M}+1_{\mathcal M}\otimes (x+h(e_1))\cdot 1_{\mathcal M}$, therefore 
$\{(x,h)\in\tilde{\mathfrak G}|\alpha_{(x,h)}(1_{\mathcal M})=\beta_{(x,h)}(1_{\mathcal M})\}
=\{(x,h)\in\tilde{\mathfrak G}|(x+h(e_1))\cdot 1_{\mathcal M}\in\mathcal P(\hat{\mathcal M}_{\mathbb Q})\}$. 

Then $a_{(x,h)},b_{(x,h)},\alpha_{(x,h)},\beta_{(x,h)}$ satisfy the assumptions of Lemma \ref{toy:lemma:BIS}. This lemma then implies the 
statement. 
\end{proof}

\begin{lem}\label{lem:3:12:3105:BIS}
One has 
$$
\mathfrak{stab}_{\tilde{\mathfrak G}}(\hat\Delta^{\mathcal W})\subset \{(x,h)\in \tilde{\mathfrak G}|(x+h(e_1))\cdot 1_{\mathcal M}
\in\mathcal P(\mathcal M)\}
$$
\end{lem}

\begin{proof}
Let $(\hat{\mathcal V}_{\mathbb Q})_0$ be the kernel of the augmentation morphism $\hat{\mathcal V}_{\mathbb Q}\to\mathbb Q$. 
By \cite{EF4}, Lem. 2.9 and Prop. 6.1, 
$$
\{y\in(\hat{\mathcal V}_{\mathbb Q})_0|(\mathrm{der}_y^{\mathcal W}\otimes id+id\otimes \mathrm{der}_y^{\mathcal W})\circ
\hat\Delta^{\mathcal W}
=\hat\Delta^{\mathcal W}\circ\mathrm{der}_y^{\mathcal W}\}\subset\{y\in(\hat{\mathcal V}_{\mathbb Q})_0|y\cdot 1_\mathcal M\in\mathcal P(\hat{\mathcal M}_{\mathbb Q})\}
$$
where for $y\in(\hat{\mathcal V}_{\mathbb Q})_0$, $\mathrm{der}_y^{\mathcal W}$ is the restriction to $\hat{\mathcal W}$ 
of the algebra derivation of $\hat{\mathcal V}$ given by $e_0\mapsto[y,e_0]$, $e_1\mapsto0$.   

One checks the identity $\mathrm{der}^{\mathcal W}_{(x,h(e_1))}=\mathrm{der}^{\mathcal W}_{x+h(e_1)}$, which implies that the preimage by 
the linear map $\tilde{\mathfrak G}\to(\hat{\mathcal V}_{\mathbb Q})_0$ given by $(x,h)\mapsto x+h(e_1)$ of the left-hand (resp. right-hand) 
side of this inclusion is the left-hand (resp. right-hand) side of the stated inclusion, which implies the said inclusion. 
\end{proof}

\begin{lem}\label{lem:eq:1329:3105:BIS}
There holds the equality $\mathfrak{stab}_{\tilde{\mathfrak G}}(\hat\Delta^{\mathcal M})
=\mathfrak{stab}_{\tilde{\mathfrak G}}(\hat\Delta^{\mathcal W})$ of Lie subalgebras of $\tilde{\mathfrak G}$. 
\end{lem}

\begin{proof}
Lemma \ref{lem:3:12:3105:BIS} implies the sequence of inclusions
$\mathfrak{stab}_{\tilde{\mathfrak G}}(\hat\Delta^{\mathcal W})\subset \{(x,h)\in \tilde{\mathfrak G}|(x+h(e_1))\cdot 1_{\mathcal M}
\in\mathcal P(\hat{\mathcal M}_{\mathbb Q})\}\cap \mathfrak{stab}_{\tilde{\mathfrak G}}(\hat\Delta^{\mathcal W})$, which together with Lemma
\ref{lem:3:11:3105:BIS} implies the inclusion $\mathfrak{stab}_{\tilde{\mathfrak G}}(\hat\Delta^{\mathcal W})\subset\mathfrak{stab}_{\tilde{\mathfrak G}}
(\hat\Delta^{\mathcal M})$. The opposite inclusion follows from \eqref{inclusions:1312:BIS}. 
\end{proof}

\begin{lem}\label{lem:3:9:1328:BIS} 
There holds the equality of subgroups 
$\mathrm{Stab}_{\tilde{\mathcal G}}(\hat\Delta^{\mathcal M})=\mathrm{Stab}_{\tilde{\mathcal G}}(\hat\Delta^{\mathcal W})$ 
of $\tilde{\mathcal G}$
for any $\mathbf k$. 
\end{lem}

\begin{proof}
The exponential map sets up a bijection $\tilde{\mathfrak G}(\mathbf k)\to {\tilde{\mathcal G}}$, which is restricted to bijections
$\mathfrak{stab}_{\tilde{\mathfrak G}}(\hat\Delta^{\mathcal M})(\mathbf k)\to \mathrm{Stab}_{{\tilde{\mathcal G}}}(\hat\Delta^{\mathcal M})$
and $\mathfrak{stab}_{\tilde{\mathfrak G}}(\hat\Delta^{\mathcal W})(\mathbf k)\to \mathrm{Stab}_{{\tilde{\mathcal G}}}(\hat\Delta^{\mathcal W})$. 
The equality of Lie subalgebras 
$\mathfrak{stab}_{\tilde{\mathfrak G}}(\hat\Delta^{\mathcal M})=\mathfrak{stab}_{\tilde{\mathfrak G}}(\hat\Delta^{\mathcal W})$ of $\tilde{\mathfrak G}$ (see Lem. \ref{lem:eq:1329:3105:BIS})
implies
then implies the equality 
$\mathrm{Stab}_{{\tilde{\mathcal G}}}(\hat\Delta^{\mathcal M})=\mathrm{Stab}_{{\tilde{\mathcal G}}}(\hat\Delta^{\mathcal W})$. 
\end{proof}

\begin{prop}\label{prop:eq:411}
    One has 
    $$\mathrm{Stab}_{\tilde{\mathcal G}}((\hat\Delta^{\mathcal W},\hat\Delta^{\mathcal M}))
=\mathrm{Stab}_{\tilde{\mathcal G}}(\hat\Delta^{\mathcal W})
    $$
    (equality of subgroups of ${\tilde{\mathcal G}}$). 
\end{prop}

\begin{proof}
This follows from Lems. \ref{lem:3:9:1328:BIS} and \ref{lem:eq:1329:3105:BIS}. 
\end{proof}

\begin{cor}\label{cor:4:20:nov}
The morphism $\underline{\mathrm{Stab}}\mathbf X_{{\tilde{\mathcal G}},\mathbf{WM}}\to
\underline{\mathrm{Stab}}\mathbf X_{{\tilde{\mathcal G}},\mathbf W}$ is an isomorphism in  
 $\mathbf{Gp}_{\mathcal G}$. 
\end{cor}

\begin{proof}
This follows from Prop. \ref{prop:eq:411} and from 
the equalities of $\underline{\mathrm{Stab}}\mathbf X_{{\tilde{\mathcal G}},\mathbf{WM}}$ and $
\underline{\mathrm{Stab}}\mathbf X_{{\tilde{\mathcal G}},\mathbf W}$ with the pairs formed by 
$\mathrm{Stab}_{\tilde{\mathcal G}}((\hat\Delta^{\mathcal W},\hat\Delta^{\mathcal M}))
$ and $\mathrm{Stab}_{\tilde{\mathcal G}}(\hat\Delta^{\mathcal W})$ with their canonical morphisms to 
$\mathcal G$ (see the definition of $\underline{\mathrm{Stab}}$ in Lem. \ref{lem:def:stab:ulstab}). 
\end{proof}

\subsection{Isomorphism status of $\underline{\mathrm{Stab}}\mathbf X_{{\tilde{\mathcal G}},\mathbf E}\to
\underline{\mathrm{Stab}}\mathbf X_{{\tilde{\mathcal G}},\mathbf{W}}$ and $\underline{\mathrm{Stab}}\mathbf X_{\mathcal G,\mathbf E''}\to
\underline{\mathrm{Stab}}\mathbf X_{\mathcal G,\mathbf H''}$}\label{sect:4:4:5786}

\begin{lem}\label{lem:421:0511}
(a) The morphism $i_{\mathbf E,\mathbf{W}} : \mathbf E\to\mathbf{W}$ in $\mathbf{PS}$ is locally injective. 

(b) The morphism $i_{\mathbf E'',\mathbf H''} : \mathbf E''\to\mathbf H''$ is in $\mathbf{PS}$ is locally injective. 
\end{lem}

\begin{proof}
    (a) follows from the injectivity of the map $\mathrm{Cop}_{e_1}(\hat{\mathcal W})\to 
    \mathrm{Cop}_{\mathcal C\operatorname{-alg}}(\hat{\mathcal W})$ underlying $i_{\mathbf E,\mathbf{WM}}$ (see Def. \ref{def:E:1001}). 

    (b) By Lem. \ref{lem:def:p_UV}, the map underlying $i_{\mathbf E'',\mathbf H''}$ is the map
\begin{equation}\label{underlying:1105}
 \mathbf k[[u,v]]^\times\backslash\mathrm{Cop}_{e_1}(\hat{\mathcal W})\to 
\mathbf k[[u,v]]^\times\backslash\mathrm{Hom}_{\mathcal C\operatorname{-alg}}(\hat{\mathcal W},\hat V)   
\end{equation}    
arising from the 
$\mathbf k[[u,v]]^\times$-equivariance of the map  $\mathrm{Cop}_{e_1}(\hat{\mathcal W})\to 
\mathrm{Hom}_{\mathcal C\operatorname{-alg}}(\hat{\mathcal W},\hat V)$ given by 
$\Delta\mapsto(i_{\mathcal W_r,\mathcal V}\circ\mathrm{Ad}_{e_1})^{\otimes2}\circ\Delta$ (see Lem. \ref{lem:toto:0930}(c) 
and Lem. \ref{lem28:1001}(d)). One checks that $i_{\mathcal W_r,\mathcal V}\circ\mathrm{Ad}_{e_1}$ is injective, which implies the
injectivity of the latter map and therefore of \eqref{underlying:1105}.   
\end{proof}

\begin{prop}\label{prop:4:22:nov}
The morphisms $\underline{\mathrm{Stab}}\mathbf X_{{\tilde{\mathcal G}},\mathbf E}\to
\underline{\mathrm{Stab}}\mathbf X_{{\tilde{\mathcal G}},\mathbf{W}}$ and $\underline{\mathrm{Stab}}\mathbf X_{\mathcal G,\mathbf E''}\to
\underline{\mathrm{Stab}}\mathbf X_{\mathcal G,\mathbf H''}$ are isomorphisms in  
 $\mathbf{Gp}_{\mathcal G}$. 
\end{prop}

\begin{proof}
 Follows from Lems. \ref{lem:421:0511} and \ref{lem:3010:1129}. 
\end{proof}

\subsection{Isomorphism status of $\underline{\mathrm{Stab}}\mathbf X_{{\tilde{\mathcal G}},\mathbf E}\to
\underline{\mathrm{Stab}}\mathbf X_{\mathcal G,\mathbf E'}$}\label{sect:4:5:5786}

\begin{lem}\label{lem:general:0611}
    Let $(X,x)$ be a pointed set, $A$ be a group, $\bullet$ be an action of $\mathcal G\times A$ on $X$. 
    The coset space $A\backslash X$ is then equipped with an action $\overline\bullet$ of $\mathcal G$, uniquely determined by the 
    condition that $(\mathcal G\times A,X,\bullet)\to(\mathcal G,A\backslash X,\overline\bullet)$ is a morphism of 
    sets with group actions, and 
    $$
    \mathrm{Stab}_{\mathcal G}(A\bullet x)=\mathrm{pr}_{\mathcal G}(\mathrm{Stab}_{\mathcal G\times A}(x)). 
    $$
    (see Notation \ref{NOTATION}(d)). 
\end{lem}

\begin{proof}
The first statement is obvious, and the second statement follows from 
\begin{align*}
&  \mathrm{Stab}_{\mathcal G}(A\bullet x)
=\{g\in\mathcal G|g\overline\bullet(A\bullet x)=A\bullet x\}
=\{g\in\mathcal G|A\bullet ((g,1)\bullet x)=A\bullet x\}
\\ & =\{g\in\mathcal G|\exists a\in A, (g,1)\bullet x=(1,a^{-1})\bullet x\}
=\{g\in\mathcal G|\exists a\in A, (g,a)\bullet x=x\}
=\mathrm{pr}_{\mathcal G}(\mathrm{Stab}_{\mathcal G\times A}(x)).   
\end{align*}
\end{proof}

\begin{lem}\label{lem:part:0611}
    One has
    $$
    \mathrm{Stab}_{1+t^2\mathbf k[[t]]}(\hat\Delta^{\mathcal W})=1, 
    $$
where the left-hand side is the stabilizer group corresponding to 
the action of the group $1+t^2\mathbf k[[t]]$ on the pointed set 
$\mathbf E=(\mathrm{Cop}_{e_1}(\hat{\mathcal W}),\hat\Delta^{\mathcal W})$ by $\odot$ (see Lem. \ref{lem:3101009:BIS}).   
\end{lem}

\begin{proof} The left-hand side is equal to 
$\{h\in 1+t^2\mathbf k[[t]] | \mathrm{Ad}_{h(e_1)h(f_1)/h(e_1+f_1)} \circ\hat\Delta^{\mathcal W}
=\hat\Delta^{\mathcal W}\}$. The exponential map $\mathrm{exp} : t^2\mathbf k[[t]]\to 1+t^2\mathbf k[[t]]$ sets up a bijection of the latter set with the Lie subalgebra 
$\{h\in t^2\mathbf k[[t]] | \mathrm{ad}_{h(e_1+f_1)-h(e_1)-h(f_1)} \circ\hat\Delta^{\mathcal W}=0\}$ of the abelian Lie algebra $t^2\mathbf k[[t]]$, where for $x\in
\hat{\mathcal W}^{\hat\otimes2}$ one defines $\mathrm{ad}_x$ to be the endomorphism of $\hat{\mathcal W}^{\hat\otimes2}$ given by $y\mapsto [x,y]$. 
If $h$ belongs to this Lie subalgebra, one has in particular  $\mathrm{ad}_{h(e_1+f_1)-h(e_1)-h(f_1)} \circ\hat\Delta^{\mathcal W}(e_0e_1)=0$, which since 
$\hat\Delta^{\mathcal W}(e_0e_1)=e_0e_1+f_0f_1-e_1f_1$ and using the commutativity of $e_1$ and $f_1$ implies the relation 
\begin{equation}\label{sncf:BIS}
[h(e_1+f_1)-h(e_1)-h(f_1),e_0e_1+f_0f_1]=0\text{ in }\hat{\mathcal W}^{\hat\otimes 2}.
\end{equation} Let $\sum_{k\geq 2}h_kt^k$ be the expression of $h(t)$, where $h_k\in\mathbf k$. 
For each $k\geq2$, the component of degree $k+2$ of \eqref{sncf:BIS} is $h_k[(e_1+f_1)^k-(e_1)^k-(f_1)^k,e_0e_1+f_0f_1]$, which therefore vanishes. 
Recall that $\hat{\mathcal W}$ has a topological basis given by the set all the words in $e_0,e_1$ not ending in $e_0$; its tensor square is then a topological 
basis of $\hat{\mathcal W}^{\hat\otimes 2}$. The coefficient of the element $[(e_1+f_1)^k-(e_1)^k-(f_1)^k,e_0e_1+f_0f_1]\in 
\hat{\mathcal W}^{\hat\otimes 2}$ with respect to the element $e_0e_1^kf_1$ of the latter basis 
is $-k$. The coefficient of this element in the left-hand side of \eqref{sncf:BIS} is then $-kh_k$, therefore \eqref{sncf:BIS} implies $-kh_k=0$, which  
implies $h_k=0$ as $\mathbf k$ contains $\mathbb Q$. It follows that $h=0$.    
\end{proof}

\begin{prop}\label{prop:4:25:1110}
(a) The morphism  $\mathrm{pr}_{1+t^2\mathbf k[[t]]} : {\tilde{\mathcal G}}\to\mathcal G$ induces an isomorphism 
$ \mathrm{Stab}_{\tilde{\mathcal G}}(\hat\Delta^{\mathcal W})\to
\mathrm{Stab}_{\mathcal G}((1+t^2\mathbf k[[t]])\odot \hat\Delta^{\mathcal W})$, where the right-hand side is the subgroup of $\mathcal G$ 
associated with the action of this group on the pointed set $\mathbf E'$ by $\odot$ (see Lem. \ref{lem:3101009:BIS}). 

(b) The morphism $\underline{\mathrm{Stab}}\mathbf X_{{\tilde{\mathcal G}},\mathbf E}\to
\underline{\mathrm{Stab}}\mathbf X_{\mathcal G,\mathbf E'}$ in 
$\mathbf{PS}_{\mathcal G}$ is an isomorphism. 
\end{prop}

\begin{proof}
It follows from the compatibility of actions that the morphism  $\mathrm{pr}_{1+t^2\mathbf k[[t]]}$ induces an 
morphism $\mathrm{Stab}_{\tilde{\mathcal G}}(\hat\Delta^{\mathcal W})\to
\mathrm{Stab}_{\mathcal G}((1+t^2\mathbf k[[t]])\odot \hat\Delta^{\mathcal W})$. 
 It follows from Lem. \ref{lem:general:0611}
applied to $(A,X,x,\bullet):=(1+t^2\mathbf k[[t]],\mathrm{Cop}_{e_1}(\hat{\mathcal W}),\hat\Delta^{\mathcal W},\odot)$ 
that this morphism is surjective, while Lem. \ref{lem:part:0611} implies that it is injective. This implies (a). 
The source and target of the morphism from (b) are the pairs formed 
respectively by $\mathrm{Stab}_{\tilde{\mathcal G}}(\hat\Delta^{\mathcal W})$ and 
$\mathrm{Stab}_{\mathcal G}((1+t^2\mathbf k[[t]])\odot \hat\Delta^{\mathcal W})$, 
together with their canonical 
morphisms to $\mathcal G$. The morphism itself is induced by the group 
morphism  $\mathrm{Stab}_{\tilde{\mathcal G}}(\hat\Delta^{\mathcal W})\to
\mathrm{Stab}_{\mathcal G}((1+t^2\mathbf k[[t]])\odot \hat\Delta^{\mathcal W})$ induced by 
$\mathrm{pr}_{1+t^2\mathbf k[[t]]}$; (a) implies that it is an 
isomorphism, which implies (b). 
\end{proof}

\subsection{Material for §§\ref{sect:G:H:2402:BIS} and \ref{sect:H:I:2402:BIS}}\label{sect:4:6:5786}

\begin{lem}\label{lem:filtration:BIS}
Let $A,B$ be $\mathbf k$-modules equipped with increasing $\mathbf k$-module filtrations $(F_iA)_{i\geq-1}$, $(F_iB)_{i\geq-1}$ with $A_{-1}=0$, 
$\cup_i F_iA=A$ and $\cup_i F_i B=B$, and let $f:A\to B$ be a $\mathbf k$-module morphism such that $f(F_iA)\subset F_iB$ 
for any $i$. If the associated graded morphism $\mathrm{gr}(f):\mathrm{gr}(A)\to\mathrm{gr}(B)$ is injective, then $f$ is injective.  
\end{lem}

\proof Let $a\in \mathrm{Ker}(f)$. Let $i$ be an integer $>0$ such that $a\in F_i(A)$. As 
$f(a)=0$, the image by $\mathrm{gr}_i(f)$ of the class of $a$ in $\mathrm{gr}_i(A)=F_i(A)/F_{i-1}(A)$
is zero, which implies that this class is zero, and therefore that $a\in F_{i-1}(A)$. One finally 
obtains $a\in F_{-1}(A)=0$. \hfill\qed\medskip 

For $M$ a $\mathbf k$-module, $m\in M$, $n\geq 1$ and $i\in[\!\![1,n]\!\!]$, set $m^{(i)}:=1^{\otimes i-1}\otimes m\otimes 1^{\otimes n-i}$.  

\begin{lem}\label{lemma:comm:0304:BIS}
For any $n\geq 1$, the centralizer of $e_1^{(1)}+\cdots+e_1^{(n)}$ in $\mathcal V^{\otimes n}$ (resp. $\hat{\mathcal V}^{\hat\otimes n}$) 
is equal to the $\mathbf k$-subalgebra (resp. topological $\mathbf k$-subalgebra) generated by $e_1^{(1)},\ldots,e_1^{(n)}$.   
\end{lem}

\proof 
Recall that $\mathcal V$ is $\mathbb Z_{\geq 0}^2$-graded, with $e_0,e_1$ of degrees $(1,0),(0,1)$. 
The direct sum of components with degree contained in 
$\{0\}\times\mathbb Z_{\geq 0}$ (resp. $\mathbb Z_{>0}\times\mathbb Z_{\geq 0}$) is $\mathcal V_0:=\mathbf k[e_1]$ 
(resp. $\mathcal V_1:=\mathcal Ve_0\mathcal V$), which implies the direct sum decomposition $\mathcal V=
\mathcal V_0\oplus\mathcal V_1$. The endomorphism $[e_1,-]$ of $\mathcal V$ is compatible with this decomposition, where we use the 
notation $[a,-]$ for the map $x\mapsto[a,x]$. Observe that the map $\mathbf k[x]\otimes\mathcal V\to \mathcal Ve_0\mathcal V$, 
$P(x)\otimes a\mapsto P(e_1)e_0a$ sets up an isomorphism of $\mathbf k$-modules. 

The decomposition of $\mathcal V$ leads to 
the decomposition $\mathcal V^{\otimes n}=\oplus_{a\in\{0,1\}^n}(\mathcal V^{\otimes n})_a$, where 
$(\mathcal V^{\otimes n})_a:=\mathcal V_{a(1)}\otimes\cdots\otimes \mathcal V_{a(n)}$. It follows from the compatibility of 
$[e_1,-]$ with the decomposition of $\mathcal V$ that the endomorphism 
$[e_1^{(1)}+\cdots+e_1^{(n)},-]$ of $\mathcal V^{\otimes n}$ is compatible with the decomposition of this $\mathbf k$-module. 
Its restriction to $(\mathcal V^{\otimes n})_0=\mathbf k[e_1]^{\otimes n}$ is zero. Let us show that its restriction to
$(\mathcal V^{\otimes n})_a$ is injective if $a\neq 0$. By tensoring a collection of 
isomorphisms $\mathbf k[x]\otimes\mathcal V\to \mathcal Ve_0\mathcal V$ indexed by 
$a^{-1}(1)$ with a collection of identity automorphisms of $\mathbf k[x]$ indexed by 
$a^{-1}(0)$, one obtains an isomorphism of $\mathbf k$-modules
$(\mathcal V^{\otimes n})_a\to\mathbf k[x_1,\ldots,x_n]\otimes(\otimes_{i\in a^{-1}(1)}\mathcal V)$ with inverse 
$P(x_1,\ldots,x_n)\otimes (\otimes_{i\in a^{-1}(1)}v_i)\mapsto 
P(e_1^{(1)},\ldots,e_1^{(n)})e_0^{(1)}\cdots e_0^{(n)}\prod_{i\in a^{-1}(1)}v_i^{(i)}$. 
The conjugation of the restriction of  $[e_1^{(1)}+\cdots+e_1^{(n)},-]$ to $(\mathcal V^{\otimes n})_a$ by this isomorphism
is the endomorphism of $\mathbf k[x_1,\ldots,x_n]\otimes(\otimes_{i\in a^{-1}(1)}\mathcal V)$ given by  
\begin{equation}\label{transport:endo:BIS}
f:=((\sum_{i\in a^{-1}(1)}x_i)\cdot-)\otimes\mathrm{id}-\mathrm{id}\otimes(-\cdot (\sum_{i\in a^{-1}(1)}e_1^{(i)})),
\end{equation} 
where $a\cdot -$ (resp. $-\cdot a$) denotes the left (resp. right) multiplication by $a$. 

Let $A$ be the $\mathbf k$-module $\mathbf k[x_1,\ldots,x_n]\otimes(\otimes_{i\in a^{-1}(1)}\mathcal V)$ equipped with the filtration defined by 
$F_i(A):=\mathbf k[x_1,\ldots,x_n]_{\leq i}\otimes(\otimes_{i\in a^{-1}(1)}\mathcal V)$ for $i\geq -1$, where $\mathbf k[x_1,\ldots,x_n]_{\leq i}$ 
is the space of polynomials of degree $\leq i$. Set $B:=A$ and set $F'_i(B):=F_{i+1}(A)$. 
Then $f$ given by \eqref{transport:endo:BIS} is a morphism $A\to B$ of filtered $\mathbf k$-modules. The associated 
graded modules can both be identified with $\mathbf k[x_1,\ldots,x_n]\otimes(\otimes_{i\in a^{-1}(1)}\mathcal V)$ (with natural grading for 
the source, and with shifted grading for the target), and $\mathrm{gr}(f)$ can then be identified with the endomorphism 
$((\sum_{i\in a^{-1}(1)}x_i)\cdot-)\otimes\mathrm{id}$
of $\mathbf k[x_1,\ldots,x_n]\otimes(\otimes_{i\in a^{-1}(1)}\mathcal V)$, which is 
injective. Lemma \ref{lem:filtration:BIS} then implies that \eqref{transport:endo:BIS} is injective, therefore that the 
endomorphism $[e_1^{(1)}+\cdots+e_1^{(n)},-]$ of $(\mathcal V^{\otimes n})_a$ is injective. 

This proves that statement on the centralizer of $e_1^{(1)}+\cdots+e_1^{(n)}$ in $\mathcal V^{\otimes n}$, which implies the similar statement 
about $\hat{\mathcal V}^{\hat\otimes n}$. 
\hfill\qed\medskip 

Recall the notation $e_i,f_i$ for the elements $e_i\otimes1$, $1\otimes e_i$ of $\mathcal W^{\otimes2}$ 
($i\in\{0,1\}$); we introduce the notation $e_i,f_i,g_i$ for the elements $e_i\otimes1^{\otimes2}$, $1\otimes e_i\otimes1$, 
$1^{\otimes2}\otimes e_i$.  We also denote by $\Delta^{\mathcal W,(2)}$ the morphism 
$(\Delta^{\mathcal W}\otimes \mathrm{id})\circ\Delta^{\mathcal W}:\mathcal W\to
\mathcal W^{\otimes3}$ and by $\hat\Delta^{\mathcal W,(2)}:
\hat{\mathcal W}\to\hat{\mathcal W}^{\hat\otimes3}$ its completion.  

\begin{lem}\label{lem:comm:Delta:BIS}
If $C\in\mathbf k[e_1,f_1,g_1]\subset\mathcal W^{\otimes3}$ (resp. 
$C\in\mathbf k[[e_1,f_1,g_1]]\subset\mathcal W^{\hat\otimes3}$) commutes with the image of 
$\Delta^{\mathcal W,(2)}$ (of $\hat\Delta^{\mathcal W,(2)}$), then $C$ belongs to $\mathbf k$. 
\end{lem}

\proof 
Let $C\in\mathbf k[e_1,f_1,g_1]$. 
Decompose $C$ as $\sum_{a,b,c\geq0}C_{abc}e_1^af_1^bg_1^c$. 
There is a direct sum decomposition 
\begin{equation}\label{decomp:1030:BIS}
\mathcal V^{\otimes3}=\mathbf k\otimes\mathcal V^{\otimes2}\oplus e_0\mathcal V^{\otimes3}\oplus
e_1\mathcal V^{\otimes3}. 
\end{equation}
Set $C_+:=\sum_{a>0,b,c\geq0}C_{abc}e_1^af_1^bg_1^c$. The commutator 
$[e_0e_1+f_0f_1+g_0g_1,C]$ can be decomposed as 
$e_0e_1 C_+-C_+e_0e_1+[f_0f_1,C]+[g_0g_1,C]$, 
where $e_0e_1 C_+$ belongs to the second summand of \eqref{decomp:1030:BIS}, 
$C_+e_0e_1$ belongs to the third summand, and $[f_0f_1,C]$ and $[g_0g_1,C]$
both belong to the direct sum of the first and third summands. 
It follows that the projection of $[e_0e_1+f_0f_1+g_0g_1,C]$ on the second summand 
is $e_0e_1 C_+$. 

Assume that $C$ commutes with the image of $\Delta^{\mathcal W,(2)}$. 
Since $\Delta^{\mathcal W,(2)}(e_0e_1)=e_0e_1+f_0f_1+g_0g_1-e_1f_1-e_1g_1-f_1g_1$, $C$ commutes with 
$e_0e_1+f_0f_1+g_0g_1$. The vanishing of $[e_0e_1+f_0f_1+g_0g_1,C]$ implies the vanishing of its projection on the
second summand of \eqref{decomp:1030:BIS}, therefore of $e_0e_1 C_+$. It follows that $C_{abc}=0$ 
for any $(a,b,c)$ with $a>0$. One proves similarly the vanishing of $C_{abc}$ for any $(a,b,c)$ such that $b>0$ or 
$c>0$. It follows that $C\in\mathbf k1$. The statement in the completed case then follows from the fact that  
$\Delta^{\mathcal W,(2)}$ is graded and $\hat\Delta^{\mathcal W,(2)}$ is its graded completion. \hfill\qed\medskip 

\begin{lem}\label{lemma:acyclicity:BIS}
The sequence of multiplicative group morphisms 
\begin{equation}\label{wanted:complex:BIS}
1+u^2\mathbf k[[u]]\stackrel{\theta}{\to}\mathbf k[[u,v]]^\times\stackrel{d_2}{\to}\mathbf k[[u,v,w]]^\times, 
\end{equation}
given by $(d_2\psi)(u,v,w):=
\psi(u,v)\psi(u+v,w)/(\psi(v,w)\psi(u,v+w))$, is such that $\mathrm{ker}(d_2)=\mathbf k^\times\cdot \mathrm{im}(\theta)$.  
\end{lem}

\proof Let $E$ be a finite dimensional $\mathbb Q$-vector space. Let $S(E)$ be the symmetric algebra over $E$ and let 
$\Delta:S(E)\to S(E)^{\otimes2}$ be the coproduct for which the elements of $E$ are primitive. The co-Hochschild complex
is $0\to \mathbb Q\to S(E)\to S(E)^{\otimes2}\to\cdots$, where the map $S(E)^{\otimes n}\to S(E)^{\otimes n+1}$ is given by 
$x\mapsto \sum_{k=1}^n(-1)^{k+1}(\mathrm{id}^{\otimes k-1}\otimes\Delta\otimes\mathrm{id}^{\otimes n-k})(x)
-1\otimes x+(-1)^nx\otimes 1$. It is quasi-isomorphic to the complex $0\to E\to \Lambda^2(E)\to\cdots$ with zero 
differerential, the morphism $\Lambda^\bullet(E)\to S(E)^{\otimes\bullet}$ being given by $\Lambda^n(E)\subset E^{\otimes n}
\to S(E)^{\otimes n}$ (see \cite{Dr}, Proposition 2.2). The complex $S(E)^{\otimes n}$ is naturally graded by the condition 
that $E$ has degree 1. One therefore obtains a new complex by restricting to positive degrees, taking the tensor product with 
$\mathbf k$, and then the completion with respect to degree, which is quasi-isomorphic to the positive degree part of 
$\Lambda^\bullet(E)\otimes\mathbf k$. When $E=\mathbb Q$, one obtains that the complex 
$0\to\mathbf k[[u]]_+\stackrel{d'_1}{\to}\mathbf k[[u,v]]_+\stackrel{d'_2}{\to}
\mathbf k[[u,v,w]]_+\to\cdots$  (the indices $+$ denote the augmentation ideals) with differentials 
$(d'_1f)(u,v):=f(u)+f(v)-f(u+v)$, $(d'_2g)(u,v,w):=
g(u+v,w)-g(u,v+w)-g(v,w)+g(u,v)$ is quasi-isomorphic to the complex $0\to \mathbf k u\to 0\to 0\to 0\to\cdots$. 
It follows that the complex
$0\to u^2\mathbf k[[u]]\stackrel{d'_1}{\to}\mathbf k[[u,v]]_+\stackrel{d'_2}{\to}
\mathbf k[[u,v,w]]_+\to\cdots$ is acyclic. 

Taking exponentials, this implies the acyclicity of the complex of multiplicative groups 
\begin{equation}\label{complex:1+:BIS}
1+u^2\mathbf k[[u]]\stackrel{\theta}{\to}1+\mathbf k[[u,v]]_+\stackrel{d_2}{\to}1+\mathbf k[[u,v,w]]_+,
\end{equation} 
therefore 
\begin{equation}\label{toto:1110}
\mathrm{ker}(1+\mathbf k[[u,v]]_+\stackrel{d_2}{\to}1+\mathbf k[[u,v,w]]_+)
=\theta(1+u^2\mathbf k[[u]]). 
\end{equation} 
The complex \eqref{wanted:complex:BIS} of multiplicative groups
is the 
product of the complexes $1\to\mathbf k^\times\stackrel{x\mapsto 1}{\to}\mathbf k^\times$,  
and
\eqref{complex:1+:BIS}. It follows that 
$\mathrm{ker}(1+\mathbf k[[u,v]]_+\stackrel{d_2}{\to}1+\mathbf k[[u,v,w]]_+)$
is the product of $\mathrm{ker}(1+\mathbf k[[u,v]]_+\stackrel{d_2}{\to}1+\mathbf k[[u,v,w]]_+)$
with $\mathrm{ker}(\mathbf k^\times\stackrel{x\mapsto 1}{\to}\mathbf k^\times)$, which is 
$\mathbf k^\times$. The result then follows from \eqref{toto:1110}. 
\hfill\qed\medskip

\subsection{Surjectivity of the group morphism $\mathrm{Stab}_{\tilde{\mathcal G}}(\hat\Delta^{\mathcal W})
\to\mathrm{Stab}_{\tsup{\mathcal G}}(\hat\Delta^{\mathcal W})$}\label{sect:G:H:2402:BIS}

\begin{defn}
  $\beta : {\tilde{\mathcal G}}\to\tsup{\mathcal G}$ 
is the group morphism given by $id_{\mathcal G}\times\theta$ (see Lem. \ref{lem28:1001} and Def. \ref{def:doubletildeG}). 
\end{defn}

\begin{lem}\label{lem:4:32:1110} 
One has $\beta(\mathrm{Stab}_{{\tilde{\mathcal G}}}(\hat\Delta^{\mathcal W})) \subset 
\mathrm{Stab}_{\tsup{\mathcal G}}(\hat\Delta^{\mathcal W})$.  
\end{lem}

\begin{proof}
This follows from the fact that the action of $\tilde{\mathcal G}$ on 
$\mathrm{Cop}_{e_1}(\hat{\mathcal W})$ is the pull-back by 
$\beta : \tilde{\mathcal G}\to \tsup{\mathcal G}$ of the action of $\tsup{\mathcal G}$ on this set. 
\end{proof}

\begin{prop}\label{prop:4:32:1110}
    (a) One has $\mathrm{Stab}_{\tsup{\mathcal G}}
(\hat\Delta^{\mathcal W})\subset\mathrm{im}(\beta)
    \cdot (1\times\mathbf k^\times)$. 

    (b) One has $\mathrm{Stab}_{\tsup{\mathcal G}}(\hat\Delta^{\mathcal W})=\beta(\mathrm{Stab}_{\tilde{\mathcal G}}(\hat\Delta^{\mathcal W}))
    \cdot (1\times\mathbf k^\times)$ (equality in $\tsup{\mathcal G}=\mathcal G \times \mathbf k[[u,v]]^\times$). 
\end{prop}

\begin{proof} 
(a) 
Let $(g,a)\in\mathrm{Stab}_{\tsup{\mathcal G}}(\hat\Delta^{\mathcal W})$. 
Then $(g,a)\odot\hat\Delta^{\mathcal W}=\hat\Delta^{\mathcal W}$, so 
\begin{equation}\label{13500104BIS:BIS}
(g,1)\odot\hat\Delta^{\mathcal W}
=(1,a)^{-1}\odot\hat\Delta^{\mathcal W}    
\end{equation}
(equality in $\mathrm{Cop}_{\mathcal C\operatorname{-alg}}(\hat{\mathcal W})$).

Call an element $\Delta\in\mathrm{Cop}_{\mathcal C\operatorname{-alg}}(\hat{\mathcal W})
=\mathrm{Hom}_{\mathcal C\operatorname{-alg}}(\hat{\mathcal W},
\hat{\mathcal W}^{\hat\otimes2})$ coassociative if $(\Delta\otimes \mathrm{id})\circ\Delta
=(\mathrm{id}\otimes\Delta)\circ\Delta$ (equality in 
$\mathrm{Hom}_{\mathcal C\operatorname{-alg}}(\hat{\mathcal W},
\hat{\mathcal W}^{\hat\otimes3})$); 
the subset of $\mathrm{Cop}_{\mathcal C\operatorname{-alg}}(\hat{\mathcal W})$ 
of coassociative elements is preserved by the action of $\mathcal G$. 
Indeed, if $\Delta$ is coassociative and $g\in\mathcal G$, then 
\begin{align*}
& (((g,1)\odot\Delta)\otimes\mathrm{id})\circ ((g,1)\odot\Delta)
=(\mathrm{aut}_g^{\mathcal W})^{\otimes 3}\circ (\Delta\otimes\mathrm{id})\circ\Delta\circ
(\mathrm{aut}_g^{\mathcal W})^{-1}
\\ & 
=(\mathrm{aut}_g^{\mathcal W})^{\otimes 3}\circ (\mathrm{id}\otimes\Delta)\circ\Delta\circ
(\mathrm{aut}_g^{\mathcal W})^{-1}
=(\mathrm{id}\otimes ((g,1)\odot\Delta))\circ ((g,1)\odot\Delta). 
\end{align*}
It was noted in \cite{EF1}, \S1.2 that $\hat\Delta^{\mathcal W}$ is coassociative, 
therefore the left-hand side of \eqref{13500104BIS:BIS} is coassociative. 
Equation \eqref{13500104BIS:BIS} then implies that the right-hand side of this equation is coassociative, so that 
\begin{equation}\label{eq:coass:a:Delta:BIS}
(((1,a)^{-1}\odot\hat\Delta^{\mathcal W})\otimes\mathrm{id})\circ ((1,a)^{-1}\odot\hat\Delta^{\mathcal W})
=(\mathrm{id}\otimes((1,a)^{-1}\odot\hat\Delta^{\mathcal W}))\circ ((1,a)^{-1}\odot\hat\Delta^{\mathcal W}). 
\end{equation}
Define maps $d'_2,d''_2 :\mathbf k[[u,v]]^\times\to\mathbf k[[u,v,w]]^\times$ by  $d'_2(a)(u,v,w):=a(u,v)a(u+v,w)$ and $d''_2(a)(u,v,w):=a(v,w)a(u,v+w)$. 
Then 
\begin{align*}
&(((1,a)^{-1}\odot\hat\Delta^{\mathcal W})\otimes\mathrm{id})\circ ((1,a)^{-1}\odot\hat\Delta^{\mathcal W})
=\mathrm{Ad}^{-1}_{a(e_1,f_1)} \circ (\hat\Delta^{\mathcal W}\otimes\mathrm{id}) \circ
\mathrm{Ad}^{-1}_{a(e_1,f_1)} \circ \hat\Delta^{\mathcal W}
\\ &
=
\mathrm{Ad}^{-1}_{a(e_1,f_1)}\circ\mathrm{Ad}^{-1}_{a(e_1+f_1,g_1)} \circ 
(\hat\Delta^{\mathcal W}\otimes\mathrm{id})  \circ \hat\Delta^{\mathcal W}
=
\mathrm{Ad}^{-1}_{(d'_2a)(e_1,f_1,g_1)}\circ 
(\hat\Delta^{\mathcal W}\otimes\mathrm{id})  \circ \hat\Delta^{\mathcal W}
\end{align*}
using $\hat\Delta^{\mathcal W}(e_1)=e_1+f_1$ and similarly 
$$
(\mathrm{id}\otimes((1,a)^{-1}\odot\hat\Delta^{\mathcal W}))\circ ((1,a)^{-1}\odot\hat\Delta^{\mathcal W})
=\mathrm{Ad}^{-1}_{(d''_2a)(e_1,f_1,g_1)}\circ 
(\mathrm{id}\otimes\hat\Delta^{\mathcal W})  \circ \hat\Delta^{\mathcal W}. 
$$
Since $d_2(a)=d'_2(a)/d''_2(a)$, \eqref{eq:coass:a:Delta:BIS} is equivalent to 
$$
\mathrm{Ad}^{-1}_{(d_2a)(e_1,f_1,g_1)}\circ\hat\Delta^{\mathcal W,(2)}=\hat\Delta^{\mathcal W,(2)}. 
$$
Lem. \ref{lem:comm:Delta:BIS} then implies that $d_2a\in\mathbf k^\times$. Moreover, 
for $a\in\mathbf k[[u,v]]^\times$, the constant term of the series $d_2a$ is 1, which implies   
$d_2a=1$. Lem. \ref{lemma:acyclicity:BIS} then implies that $a$ is in 
$\mathbf k^\times\cdot\theta(1+t^2\mathbf k[[t]])$. 
the image of $\theta$. 
If $(\lambda,f)\in\mathbf k^\times\times (1+t^2\mathbf k[[t]])$ is such that $a=\lambda\cdot \theta(f)$, 
then $(g,a)=(1,\lambda)\cdot \beta(g,f)$, therefore 
$(g,a) \in \mathbf k^\times\cdot \mathrm{im}(\beta)$. 

(b) One obviously has $1\times\mathbf k^\times\subset \mathrm{Stab}_{\tsup{\mathcal G}}(\hat\Delta^{\mathcal W})$, which 
together with Lem. \ref{lem:4:32:1110} implies the inclusion $\beta(\mathrm{Stab}_{{\tilde{\mathcal G}}}(\hat\Delta^{\mathcal W})) 
\cdot (1\times\mathbf k^\times)\subset \mathrm{Stab}_{\tsup{\mathcal G}}(\hat\Delta^{\mathcal W})$. Let us prove the opposite inclusion. 
Let $(g,a)\in\mathrm{Stab}_{\tsup{\mathcal G}}(\hat\Delta^{\mathcal W})$. By (a), there exists 
$(\lambda,f)\in\mathbf k^\times\times (1+t^2\mathbf k[[t]])$ such that 
$(g,a)=(1,\lambda)\cdot \beta(g,f)$. Since $(1,\lambda)\in\mathrm{Stab}_{\tsup{\mathcal G}}(\hat\Delta^{\mathcal W})$, 
$\beta(g,f)=(1,\lambda)^{-1}\cdot(g,a)\in \mathrm{Stab}_{\tsup{\mathcal G}}(\hat\Delta^{\mathcal W})$.
Then $\beta(g,f)\odot\hat\Delta^{\mathcal W}=\hat\Delta^{\mathcal W}$, which implies 
$(g,f)\odot\hat\Delta^{\mathcal W}=\hat\Delta^{\mathcal W}$, therefore 
$(g,f)\in\mathrm{Stab}_{{\tilde{\mathcal G}}}(\hat\Delta^{\mathcal W})$. Therefore 
$\mathrm{Stab}_{\tsup{\mathcal G}}(\hat\Delta^{\mathcal W})\subset\beta(\mathrm{Stab}_{{\tilde{\mathcal G}}}(\hat\Delta^{\mathcal W})) 
\cdot (1\times\mathbf k^\times)$. 
\end{proof}

\subsection{Isomorphism status of $\underline{\mathrm{Stab}}\mathbf X_{\mathcal G,\mathbf E'}\to
\underline{\mathrm{Stab}}\mathbf X_{\mathcal G,\mathbf E''}$}\label{sect:H:I:2402:BIS}

\begin{prop}\label{prop:4:final:1110}
    The groups $\mathrm{Stab}_{\mathcal G}([\hat\Delta^{\mathcal W}])$ and 
    $\mathrm{Stab}_{\mathcal G}([\!\![\hat\Delta^{\mathcal W}]\!\!])$, relative to the respective actions of $\mathcal G$
    on $(1+t^2\mathbf k[[t]])\backslash\mathrm{Cop}_{e_1}(\hat{\mathcal W})$ and 
    $\mathbf k[[u,v]]^\times\backslash\mathrm{Cop}_{e_1}(\hat{\mathcal W})$, are equal. 
\end{prop}

\begin{proof} In order to avoid an ambiguity in notation, we denote in this proof by 
$\mathrm{pr}^{\tilde{\mathcal G}}_{\mathcal G}$ and by 
$\mathrm{pr}^{\tsup{\mathcal G}}_{\mathcal G}$ the projections
$\tilde{\mathcal G}\to\mathcal G$ and $\tsup{\mathcal G}\to\mathcal G$ denoted 
$\mathrm{pr}_{\mathcal G}$ in this text. 

One has the following equality of subsets of $\mathcal G$: 
\begin{align*}
&\mathrm{Stab}_{\mathcal G}([\hat\Delta^{\mathcal W}])
    =\mathrm{pr}^{\tilde{\mathcal G}}_{\mathcal G}(\mathrm{Stab}_{\tilde{\mathcal G}}(\hat\Delta^{\mathcal W})) 
    =\mathrm{pr}^{\tsup{\mathcal G}}_{\mathcal G}\circ\beta(\mathrm{Stab}_{\tilde{\mathcal G}}(\hat\Delta^{\mathcal W}))
\\&   =\mathrm{pr}^{\tsup{\mathcal G}}_{\mathcal G}
 (\beta(\mathrm{Stab}_{\tilde{\mathcal G}}(\hat\Delta^{\mathcal W}))
    \cdot (1\times\mathbf k^\times))
    =\mathrm{pr}^{\tsup{\mathcal G}}_{\mathcal G}(\mathrm{Stab}_{\tsup{\mathcal G}}
(\hat\Delta^{\mathcal W}))
=\mathrm{Stab}_{\mathcal G}([\!\![\hat\Delta^{\mathcal W}]\!\!])
\end{align*} 
where the first and last equalities follow from  Lem. \ref{lem:general:0611}, the second equality follows from the 
commutativity of $$
\xymatrix{{\tilde{\mathcal G}}\ar^\beta[rr]\ar_{\mathrm{pr}^{\tilde{\mathcal G}}_{\mathcal G}}[rd]&&
\tsup{\mathcal G}
\ar^{\mathrm{pr}^{\tsup{\mathcal G}}_{\mathcal G}}[ld]\\&\mathcal G&}$$
the third equality follows from the fact that $\mathrm{pr}_{\mathcal G}^{\tsup{\mathcal G}}$ 
is a group morphism and that $1\times\mathbf k^\times$ is contained in its kernel, and the fourth 
equality follows by applying $\mathrm{pr}_{\mathcal G}^{\tsup{\mathcal G}}$  to the equality of 
Prop. \ref{prop:4:32:1110}(b).  
\end{proof}

\begin{prop}\label{prop:4:34:nov}
  The morphism $\underline{\mathrm{Stab}}\mathbf X_{\mathcal G,\mathbf E'}\to
\underline{\mathrm{Stab}}\mathbf X_{\mathcal G,\mathbf E''}$ is an isomorphism in 
$\mathbf{Gp}_{\mathcal G}$.   
\end{prop}

\begin{proof}
It follows from the existence of the morphism $\mathbf X_{\mathcal G,\mathbf E'}\to
\mathbf X_{\mathcal G,\mathbf E''}$
in Lem \ref{lem:def:p_UV} that the canonical projection 
$$
((1+t^2\mathbf k[[t]])\backslash\mathrm{Cop}_{e_1}(\hat{\mathcal W}),
[\hat\Delta^{\mathcal W}])\to
((\mathbf k[[u,v]]^\times\backslash\mathrm{Cop}_{e_1}(\hat{\mathcal W}),
[\!\![\hat\Delta^{\mathcal W}]\!\!])
$$
of pointed sets is compatible with the actions of $\mathcal G$ on both sides. 
By Lem. \ref{BASIC}(a), this implies that the stabilizer group of the source is 
contained in its counterpart for the target (inclusion of subgroups of $\mathcal G$). 
These groups are respectively 
$\mathrm{Stab}_{\mathcal G}([\hat\Delta^{\mathcal W}])$ and 
    $\mathrm{Stab}_{\mathcal G}([\!\![\hat\Delta^{\mathcal W}]\!\!])$, therefore
    $$\mathrm{Stab}_{\mathcal G}([\hat\Delta^{\mathcal W}])\subset \mathrm{Stab}_{\mathcal G}([\!\![\hat\Delta^{\mathcal W}]\!\!]). 
$$ 
The morphism $\underline{\mathrm{Stab}}\mathbf X_{\mathcal G,\mathbf E'}\to
\underline{\mathrm{Stab}}\mathbf X_{\mathcal G,\mathbf E''}$ 
is the commutative triangle formed by this injection together with the injections of its source and target in 
$\mathcal G$. Prop. \ref{prop:4:final:1110} then implies the statement. 
\end{proof}

\subsection{Isomorphism status of $\underline{\mathrm{Stab}}\mathbf X_{\mathcal G,\mathbf E''}\to
\underline{\mathrm{Stab}}\mathbf X_{\mathcal G,\mathbf H''}$}\label{sect:4:9:5786}

\begin{lem}\label{lem:general:2212}
(a)    If $f : (S,s)\to(S',s')$ is a morphism of pointed sets and if $\Gamma$ is a group, acting on $S$ and $S'$, and such that $f$ is 
    $\Gamma$-equivariant. Then there is an inclusion $\mathrm{Stab}_\Gamma(s)\subset\mathrm{Stab}_\Gamma(s')$ of subgroups of $\Gamma$.  

(b) If moreover $f$ is locally injective (see Def. \ref{def:local:injectivity:1018}), then one has $\mathrm{Stab}_\Gamma(s)=\mathrm{Stab}_\Gamma(s')$. 
\end{lem}

\begin{proof}
(a) is obvious. Let us show (b). If $\gamma\in\mathrm{Stab}_\Gamma(s')$, then for any $f(s)=s'=\gamma\bullet s'=\gamma\bullet f(s)=f(\gamma\bullet s)$, which implies
by the injectivity of $f$ the equality $s=\gamma\bullet s$, therefore $\gamma\in \mathrm{Stab}_\Gamma(s)$. 
\end{proof}

\begin{lem}\label{lem:4:37:jan}
  (a) The set map underlying $i_{\mathbf E,\mathbf H} : \mathbf E\to\mathbf H$ is injective. 

  (b) The set map underlying $i_{\mathbf E'',\mathbf H'''} : \mathbf E''\to\mathbf H''$ is injective. 
\end{lem}

\begin{proof}
 (a) The map $\hat{\mathcal W}^{\hat\otimes2}\stackrel{\mathrm{Ad}_{e_1}^{\otimes 2}}{\to}\hat{\mathcal W}_r^{\hat\otimes2}\hookrightarrow\hat{\mathcal V}^{\hat\otimes2}=\hat V$
is injective since $\mathrm{Ad}_{e_1}$ is an isomorphism. It follows that the map 
 $\mathrm{Hom}_{\mathcal C\operatorname{-alg}}(\hat{\mathcal W},\hat{\mathcal W}^{\hat\otimes2})\to 
 \mathrm{Hom}_{\mathcal C\operatorname{-alg}}(\hat{\mathcal W},\hat V)$ induced by composition with this map is injective as well, and therefore that the same is true 
 of its restriction to $\mathrm{Cop}_{e_1}(\hat{\mathcal W})$. The statement then follows from the fact that this restriction is the said 
set map. 
 
(b) In Lem. \ref{lem28:1001}(b) and (c), one defines actions of the group $\mathbf k[[u,v]]^\times$ on the source and target of the map from (a), and 
by Lem. \ref{lem28:1001}(d) this map is $\mathbf k[[u,v]]^\times$-equivariant; and by Lem. \ref{lem:def:p_UV}, the resulting map between quotient sets is the map underlying 
$i_{\mathbf E'',\mathbf H'''} : \mathbf E''\to\mathbf H''$. 
The result then follows from the following general statement: if $\Gamma$ is a group, if $X,Y$ are $\Gamma$-sets, and if $X\to Y$ is injective and $\Gamma$-equivariant, 
then the induced map $\Gamma\backslash X\to\Gamma\backslash Y$ is injective. 
\end{proof}

\begin{prop}\label{prop:missing:iso}
  The morphism $\underline{\mathrm{Stab}}\mathbf X_{\mathcal G,\mathbf E''}\to
\underline{\mathrm{Stab}}\mathbf X_{\mathcal G,\mathbf H''}$ is an isomorphism in 
$\mathbf{Gp}_{\mathcal G}$.   
\end{prop}

\begin{proof}
This morphism is the commutative triangle formed by the inclusion of two subgroups of $\mathcal G$ and by their inclusions in $\mathcal G$; 
these subgroups being the stabilizers $\mathrm{Stab}_{\mathcal G}(\mathbf k[[u,v]]^\times
    \odot\hat\Delta^{\mathcal W})$ and $\mathrm{Stab}_{\mathcal G}(\mathbf k[[u,v]]^\times
    \odot\Delta^{\mathcal W}_{r,l})$ corresponding to the actions of $\mathcal G$ on 
    $\mathbf k[[u,v]]^\times\backslash\mathrm{Cop}_{e_1}(\hat{\mathcal W})$ and on 
    $\mathbf k[[u,v]]^\times\backslash\mathrm{Hom}_{\mathcal C\operatorname{-alg}}(\hat{\mathcal W},\hat V)$. 
It follows from Lem. \ref{lem:general:2212}(b) and Lem. \ref{lem:4:37:jan}(b) that these subgroups of $\mathcal G$ are equal, which implies the statement. 
\end{proof}

\newpage 
\part{The group inclusion $\mathrm{Stab}_{\mathcal G}(\mathrm{GL}_3\hat V\bullet \rho_{\mathrm{DT}}
)\subset \mathrm{Stab}_{\mathcal G}(\mathbf k[[u,v]]^\times\bullet\Delta^{\mathcal W}_{r,l}
)$}\label{part 2}

\section{Statement and proof of group inclusion $\mathrm{Stab}_{\mathcal G}(\mathrm{GL}_3\hat V\bullet \rho_{\mathrm{DT}}
)\subset 
\mathrm{Stab}_{\mathcal G}(\mathbf k[[u,v]]^\times\bullet\Delta^{\mathcal W}_{r,l}
)$}

The objective of this section is the proof of the inclusion $\mathrm{Stab}_{\mathcal G}
(\mathrm{GL}_3\hat V\bullet \rho_{\mathrm{DT}}
)\subset \mathrm{Stab}_{\mathcal G}(
\mathbf k[[u,v]]^\times\bullet\Delta^{\mathcal W}_{r,l}
)$
(see Thm. \ref{thm:5:31:3103}). 
This inclusion is obtained as the consequence of a diagram of $\mathcal G$-equivariant pointed sets with 
group actions (see Def. \ref{def:quotient:diagram}), 
which is inspired by \cite{DT,T,EF1} and is constructed in several steps. A pointed set diagram is 
constructed in \S\ref{sect:4:2012} (see \eqref{TAG':2112}). After a preliminary \S\ref{lem:4:2:2012}, 
a group diagram is constructed in 
\S\ref{sec 5.3} (see \eqref{TAG:2112:FIRST}), and the compatibility of the pointed sets and group 
diagrams is obtained in \S\ref{sec 5.4}. 
The overall action of $\mathcal G$ on groups is constructed in §\ref{sec 5.5}, and its 
counterpart for sets is constructed in \S\ref{sec 5.6}; the main consequence (Thm. \ref{thm:5:31:3103}) is drawed in 
\S\ref{sect:5:7}.

\subsection{A diagram of pointed sets}\label{sect:4:2012}

For $R$ a $\mathbf k$-module and $k,l\geq1$, define $M_{k,l}R$ as the set of matrices of size $(k,l)$ with coefficients in $R$; this is 
a $\mathbf k$-algebra if $k=l$ and $R$ is a $\mathbf k$-algebra. When $R=\hat V$, $M_{k,l}\hat V$ is a complete graded $\mathbf k$-module
with degree $n$ component $(M_{k,l}\hat V)_n:=M_{k,l}(\hat V_n)$; it is a complete graded $\mathbf k$-algebra if $k=l$. 

\begin{defn}\label{def:hom:m3:jan:2025}
Recall that $\mathrm{Hom}_{\mathcal C\operatorname{-alg}}(\hat{\mathcal V},M_3\hat V)$ 
denotes the set of morphisms $\rho : \hat{\mathcal V}\to M_3\hat V$ of filtered $\mathbf k$-algebras, both sides being equipped with the 
decreasing filtrations (denoted $F^\bullet$) associated with their complete graded structures. 
\end{defn}

It follows from the freeness of $\hat{\mathcal V}$ that the assignment $\rho\mapsto (\rho(e_0),\rho(e_1))$ sets up a bijection between 
$\mathrm{Hom}_{\mathcal C\operatorname{-alg}}(\hat{\mathcal V},M_3\hat V)$ and $M_3(F^1\hat V)^2$.

\begin{defn}\label{def:5:2:1926}
Define
    $$
    \rho_0:=\begin{pmatrix}
    e_0&0&0\\e_1&f_0&-e_1\\0&0&e_0
\end{pmatrix},
\quad
\rho_1:=\begin{pmatrix}
    1\\-1\\0
\end{pmatrix}\begin{pmatrix}
    e_1&-f_1&0
\end{pmatrix}, \quad 
\mathrm{col}_{\mathrm{DT}}:=\begin{pmatrix}
    1\\-1\\0
\end{pmatrix},\quad 
\mathrm{row}_{\mathrm{DT}}:=\begin{pmatrix}
    e_1&-f_1&0
\end{pmatrix};  
$$
then $\rho_0,\rho_1\in M_3\hat V$,  $\mathrm{col}_{\mathrm{DT}}\in M_{3,1}\hat V$, $\mathrm{row}_{\mathrm{DT}}\in M_{1,3}\hat V$,
and $\rho_1=\mathrm{col}_{\mathrm{DT}}\cdot \mathrm{row}_{\mathrm{DT}}$. 
\end{defn}

\begin{defn}\label{def:5:3:paris}
    Define $\rho_{\mathrm{DT}}$ to be the element of $\mathrm{Hom}_{\mathcal C\operatorname{-alg}}(\hat{\mathcal V},M_3\hat V)$
    such that $\rho_{\mathrm{DT}}(e_0)=\rho_0$ and $\rho_{\mathrm{DT}}(e_1)=\rho_1$. 
\end{defn}

\begin{defn}\label{defn: Hom1}
Define the set $\mathrm{Hom}_{\mathcal C\operatorname{-alg}}^1(\hat{\mathcal V},M_3\hat V)$ as the 
subset of $\mathrm{Hom}_{\mathcal C\operatorname{-alg}}(\hat{\mathcal V},M_3\hat V)$ of all the elements $\rho$ such that 
$\rho(e_1)=\rho_1$. 
\end{defn}

One has $\rho_{\mathrm{DT}}\in \mathrm{Hom}_{\mathcal C\operatorname{-alg}}^1(\hat{\mathcal V},M_3\hat V)$. 

\begin{lem}\label{lem:vania}
 (a) If $\rho\in\mathrm{Hom}_{\mathcal C\operatorname{-alg}}^1(\hat{\mathcal V},M_3\hat V)$, there is a unique
 element $\Delta_\rho\in \mathrm{Hom}_{\mathcal C\operatorname{-alg}}(\hat{\mathcal W},\hat V)$, such that 
 for any $n\geq1$, $\Delta_\rho(e_1e_0^{n-1})=\mathrm{row}_{\mathrm{DT}}\cdot \rho(e_0)^{n-1}\cdot \mathrm{col}_{\mathrm{DT}}$. 
 The assignment $\rho\mapsto\Delta_\rho$ defines a map $\mathrm{Hom}_{\mathcal C\operatorname{-alg}}^1(\hat{\mathcal V},M_3\hat V)\to 
    \mathrm{Hom}_{\mathcal C\operatorname{-alg}}(\hat{\mathcal W},\hat V)$. 

    (b) The map $\rho\mapsto\Delta_\rho$ takes $\rho_{\mathrm{DT}}$ to $\Delta^{\mathcal W}_{r,l}$.    
\end{lem}

\begin{proof}
(a) follows from the fact that $\hat{\mathcal W}$ is freely generated by the family $(e_0^{n-1}e_1)_{n\geq1}$. 
(b) For $n\geq1$, one has 
\begin{align*}
& \Delta_{\rho_{\mathrm{DT}}}(e_0^{n-1}e_1)=\begin{pmatrix}
   e_1 & -f_1
\end{pmatrix}\begin{pmatrix}
   e_0 & 0\\ e_1& f_0 
\end{pmatrix}^{n-1}\begin{pmatrix}
1 \\ -1  
\end{pmatrix}
=\begin{pmatrix}
   e_1 & -f_1
\end{pmatrix}\begin{pmatrix}
   e_0^{n-1} & 0\\ \sum_{k=0}^{n-2}f_0^ke_1e_0^{n-2-k}& f_0^{n-1} 
\end{pmatrix}\begin{pmatrix}
1 \\ -1  
\end{pmatrix}
\\ & =e_1e_0^{n-1}+f_0f_1^{n-1}-\sum_{k=0}^{n-2}f_1f_0^ke_1e_0^{n-2-k}=\Delta^{\mathcal W}_{r,l}(e_0^{n-1}e_1). 
\end{align*}   
\end{proof}
    
\begin{cor}
The canonical inclusion and the map $\rho\mapsto\Delta_\rho$ build up a diagram of pointed sets
\begin{equation}\label{TAG':2112}
 (\mathrm{Hom}_{\mathcal C\operatorname{-alg}}(\hat{\mathcal V},M_3\hat V),\rho_{\mathrm{DT}})\leftarrow
(\mathrm{Hom}_{\mathcal C\operatorname{-alg}}^1(\hat{\mathcal V},M_3\hat V),\rho_{\mathrm{DT}})\to
(\mathrm{Hom}_{\mathcal C\operatorname{-alg}}(\hat{\mathcal W},\hat V),\Delta^{\mathcal W}_{r,l}).    
\end{equation}
\end{cor}

\begin{proof}
Follows from Lem. \ref{lem:vania}(b) and the relation 
$\rho_{\mathrm{DT}}\in \mathrm{Hom}_{\mathcal C\operatorname{-alg}}^1(\hat{\mathcal V},M_3\hat V)$. 
\end{proof}

\begin{rem}\label{rem:5:7:3012}
 The algebra morphism $\rho:  \mathcal V\to M_3\hat V$ from $\rho$ in \cite{EF1}, (5.2.5) given by 
 $$
\rho(e_0)=
\begin{pmatrix}
    e_0&0&0\\0&-e_1+f_0 & -e_1\\0&e_0+e_1-f_0 & e_0+e_1
\end{pmatrix}, \quad \rho(e_1)=\rho_1
$$
(see {\it loc. cit.,} (5.2.9) and Lem. 5.6) 
is related to  $\rho_{\mathrm{DT}}$ by  $\rho=\mathrm{Ad}_U\circ \rho_{\mathrm{DT}}$, where 
$U:=\begin{pmatrix}
    1&0&0\\0&1&0\\-1&-1&1
\end{pmatrix}\in\mathrm{GL}_3(\mathbf k)$. Note that $\mathrm{col}_{\mathrm{DT}},\mathrm{row}_{\mathrm{DT}}$ are denoted 
$\mathrm{col}_1,\mathrm{row}_1$ in \cite{EF1}. 
\end{rem}

\subsection{Computation of the algebra $\mathrm C_3(\rho_1)$}\label{lem:4:2:2012}

\begin{lem}\label{lem:TODO:BIS}
(a) The endomorphisms $x\mapsto e_1\cdot x$ and $x\mapsto f_1\cdot x$ of $\hat V$ are injective. 

(b) For $(\mathrm{col},\mathrm{row})\in M_{3\times 1}\hat V\times M_{1\times 3}\hat V$, the equality 
\begin{equation}\label{1652:1905:BIS}
\mathrm{col}\cdot \mathrm{row}_{\mathrm{DT}}=\mathrm{col}_{\mathrm{DT}}\cdot \mathrm{row}
\end{equation}
is equivalent to the existence of $a\in \hat V$ such that $\mathrm{col}=\mathrm{col}_{\mathrm{DT}}\cdot a$ and 
$\mathrm{row}=a\cdot\mathrm{row}_{\mathrm{DT}}$. 
\end{lem}

\begin{proof}
(a) follows from the fact that a topological basis of $\hat V$ is given by the set of all pairs of a word in $e_0,e_1$ and a word in 
$f_0,f_1$, and that its image by either $x\mapsto e_1x$ or $x\mapsto f_1x$ is a subfamily of itself.  

(b) Let $c_i,r_i\in\hat V$ ($i\in[\!\![1,3]\!\!]$) be such that $\mathrm{col}=\begin{pmatrix}c_1\\c_2\\c_3\end{pmatrix}$ 
and $\mathrm{row}=\begin{pmatrix}r_1&r_2&r_3\end{pmatrix}$. 
\eqref{1652:1905:BIS} is equivalent to 
$$
c_1e_1=r_1,\quad -c_1f_1=r_2,\quad   
c_2e_1=-r_1,\quad -c_2f_1=-r_2,\quad c_3e_1=c_3f_1=r_3=0, 
$$
i.e. to 
$$
r_1=c_1e_1,\quad r_2=-c_1f_1,\quad 
e_1\cdot (r_1+r_2)=e_1r_3=f_1\cdot (c_1+c_2)e_1=(c_1+c_2)f_1=c_3e_1=c_3f_1=r_3=0.   
$$
By (a), the last equation is equivalent to $r_3=c_3=c_1+c_2=0$. Set $a:=c_1$, one then obtains $(c_1,c_2,c_3)=(a,-a,0)$ and,
using the first equation, $(r_1,r_2,r_3)=(ae_1,-af_1,0)$, which implies 
$\mathrm{col}=\mathrm{col}_{\mathrm{DT}}\cdot a$ and 
$\mathrm{row}=a\cdot\mathrm{row}_{\mathrm{DT}}$. This proves one of the implications; its converse is obvious. 
\end{proof}

\begin{lem}\label{LEM1:0301:BIS}
The sequence $\hat V\to \hat V^{\oplus2}\to \hat V$ where first map is $x\mapsto (f_1x,e_1x)$ and the 
second map is $(u,v)\mapsto e_1 u-f_1v$, is exact. 
\end{lem}

\begin{proof}
The said sequence is a complex since $e_1$ and $f_1$ commute. Let us show that it is exact. Let $u,v \in \hat V$ be 
such that $f_1u=e_1v$. By the direct sum decomposition 
\begin{equation}\label{decomp:V:0802:BIS}
\hat{\mathcal V}=\mathbf k1\oplus (\bigoplus_{i\in\{0,1\}} e_i\hat{\mathcal V}). 
\end{equation}
Then $u$, resp. $v$ may be uniquely decomposed as $u=e_1u_1+e_0u_0+1 \otimes u_\emptyset$, where $u_0,u_1 \in \hat V$ and 
$u_\emptyset \in \hat{\mathcal V}$, resp.  
$v=f_1v_1+f_0v_0+v_\emptyset \otimes 1$, where $v_0,v_1 \in \hat V$ and $v_\emptyset\in\hat{\mathcal V}$. 
Then $f_1u=e_1v$ implies
\begin{equation}\label{titi:0802:BIS}
e_1f_1(u_1-v_1)+e_0f_1u_0-e_1f_0v_0+1 \otimes e_1u_\emptyset-e_1v_\emptyset \otimes 1=0. 
\end{equation}
The direct sum decomposition \eqref{decomp:V:0802:BIS} gives rise to a direct sum decomposition
$$ 
\hat V=\mathbf k\cdot 1^{\otimes2}\oplus (\bigoplus_{i\in\{0,1\}} e_i\hat{\mathcal V}\otimes 1) 
\oplus (\bigoplus_{i\in\{0,1\}} 1\otimes e_i\hat{\mathcal V}) \oplus (\bigoplus_{(i,j)\in\{0,1\}^2} e_if_j\hat V). 
$$
\eqref{titi:0802:BIS} then implies that all the summands in its left-hand side are zero, which by the integrity of $\hat V$ 
implies $u_1=v_1$, $u_0=v_0=0$, and $u_\emptyset=v_\emptyset=0$, 
therefore $(u,v)$ is the image of $u_1=v_1$ by the first map of the complex. This proves the claimed exactness.   
\end{proof}

\begin{defn}\label{def:barcol:barrow}
Set 
    $$
    \overline{\mathrm{col}}_{\mathrm{DT}}:=\begin{pmatrix}
        1\\-1
    \end{pmatrix}\in M_{2,1}\hat V,\quad \overline{\mathrm{row}}_{\mathrm{DT}}:=\begin{pmatrix}
        e_1 & -f_1
    \end{pmatrix}\in M_{1,2}\hat V. 
    $$
\end{defn}

\begin{defn}\label{def:van:0301:BIS}
(a) Set $\mathrm{Ann}(\mathrm{col}_{\mathrm{DT}},\mathrm{row}_{\mathrm{DT}})
:=\{A\in M_3\hat V|A\cdot\mathrm{col}_{\mathrm{DT}}=0$ and $\mathrm{row}_{\mathrm{DT}}\cdot A=0\}$, where 
$\mathrm{col}_{\mathrm{DT}},\mathrm{row}_{\mathrm{DT}}$ are as in Lem.-Def. \ref{def:5:2:1926}. 

(b) Set $\mathrm{Ann}(\overline{\mathrm{col}}_{\mathrm{DT}},\overline{\mathrm{row}}_{\mathrm{DT}})
:=\{A\in M_2\hat V|A\cdot\overline{\mathrm{col}}_{\mathrm{DT}}=0$ and $\overline{\mathrm{row}}_{\mathrm{DT}}\cdot A=0\}$. 
\end{defn}

\begin{lem}\label{lem:van:0301:BIS}
$\mathrm{Ann}(\overline{\mathrm{col}}_{\mathrm{DT}},\overline{\mathrm{row}}_{\mathrm{DT}})
=\{\begin{pmatrix}f_1\\e_1\end{pmatrix}
\cdot v\cdot \begin{pmatrix}1&1\end{pmatrix}|v\in\hat V\}$.  
\end{lem}

\begin{proof}
It follows from $\begin{pmatrix}e_1&-f_1\end{pmatrix}\cdot \begin{pmatrix}f_1\\e_1\end{pmatrix}=0$ 
and $ \begin{pmatrix}1&1\end{pmatrix}\cdot\begin{pmatrix}1\\-1\end{pmatrix}=0$ that the right-hand side of the announced
equality is contained in $\mathrm{Ann}(\overline{\mathrm{col}}_{\mathrm{DT}},\overline{\mathrm{row}}_{\mathrm{DT}})$. 
Let us prove that the opposite inclusion.  Let 
$A=
\begin{pmatrix}
a_{11}    & a_{12}\\ a_{21}& a_{22}\end{pmatrix}\in \mathrm{Ann}(\overline{\mathrm{col}}_{\mathrm{DT}},\overline{\mathrm{row}}_{\mathrm{DT}})$. 
Since $\begin{pmatrix}e_1&-f_1\end{pmatrix}\cdot A=0$, one has $e_1\cdot a_{1i}=f_1\cdot a_{2i}$ for $i=1,2$. 
By Lem. \ref{LEM1:0301:BIS}, this implies the existence of $x_i\in\hat V$, where $i=1,2$, with 
$a_{1i}=f_1\cdot x_i$  and $a_{2i}=e_1\cdot x_i$ for $i=1,2$. Then 
$A=
\begin{pmatrix}
f_1x_1    & f_1x_2 \\ e_1x_1 & e_1x_2 
\end{pmatrix}$. Since $A\cdot\begin{pmatrix}1\\-1\end{pmatrix}=0$, one has 
$f_1(x_1-x_2)=e_1(x_1-x_2)=0$, which by the injectivity of $x\mapsto e_1x$ 
implies $x_2=x_1$. This equality implies that $A$ has the announced form, with $v:=x_1$. 
\end{proof}

\begin{lem}\label{lem:5:11:2008}
    $\mathrm{Ann}(\mathrm{col}_{\mathrm{DT}},\mathrm{row}_{\mathrm{DT}})=\{\begin{pmatrix}f_1&0\\e_1&0 \\0&1\end{pmatrix}
\cdot m\cdot \begin{pmatrix}1&1&0\\0&0&1\end{pmatrix}|m\in M_2\hat V\}$. 
\end{lem}

\begin{proof}
Let $M\in M_3\hat V$ be decomposed as $M=\begin{pmatrix}
    a&c\\b&d
\end{pmatrix}$, where $a\in M_2\hat V$, $b\in M_{1,2}\hat V$, $c\in M_{2,1}\hat V$, $d\in\hat V$.  
Then $M\in \mathrm{Ann}(\mathrm{col}_{\mathrm{DT}},\mathrm{row}_{\mathrm{DT}})$ is equivalent to the 
conjunction of 
$$
a\in \mathrm{Ann}(\overline{\mathrm{col}}_{\mathrm{DT}},\overline{\mathrm{row}}_{\mathrm{DT}}), \quad 
b\cdot \overline{\mathrm{col}}_{\mathrm{DT}}=0, \quad \overline{\mathrm{row}}_{\mathrm{DT}}\cdot c=0. 
$$
By Lem. \ref{lem:van:0301:BIS}, the first condition is equivalent to the existence of $\alpha\in\hat V$ such that
$a=\begin{pmatrix}f_1\\e_1\end{pmatrix}
\cdot \alpha\cdot \begin{pmatrix}1&1\end{pmatrix}$. The second condition is equivalent to the existence of $\beta\in\hat V$ such that
$b=\begin{pmatrix}\beta&\beta\end{pmatrix}$. By Lem. \ref{LEM1:0301:BIS}, the third condition is equivalent to the existence of $\gamma\in\hat V$ such that
$c=\begin{pmatrix}f_1\gamma\\e_1\gamma\end{pmatrix}$. The condition $M\in \mathrm{Ann}(\mathrm{col}_{\mathrm{DT}},\mathrm{row}_{\mathrm{DT}})$
is therefore equivalent to the existence of a matrix $m=\begin{pmatrix}
    \alpha & \gamma\\\beta&\delta \end{pmatrix}\in M_2\hat V$, such that 
   $$
   M=\begin{pmatrix}
      f_1\alpha & f_1\alpha&f_1\gamma \\ e_1\alpha & e_1\alpha&e_1\gamma\\ \beta& \beta& \delta
   \end{pmatrix}. 
   $$ 
The statement then follows from the equality of the right-hand side with   
$$
\begin{pmatrix}f_1&0\\e_1&0 \\0&1\end{pmatrix}
\cdot m\cdot \begin{pmatrix}1&1&0\\0&0&1\end{pmatrix}. 
$$
\end{proof}

\begin{defn}
For $A$ an algebra and $X\subset A$ a subset, $\mathrm{C}_A(X)$ is the commutant of $X$ in $A$, i.e. the subset of $A$ of 
elements $a$ such that $ax=xa$ for any $x\in X$.  We will write $\mathrm{C}_k(X)$ instead of $\mathrm C_{A}(X)$ if $A=M_k\hat V$, for any
$k\geq1$.  
\end{defn}

\begin{lem}\label{invertible:commutant}
    For $A$ an algebra and $X\subset A$ a subset, $\mathrm{C}_A(X)$ is a subalgebra of $A$. Its group of units is given by 
    $\mathrm{C}_A(X)^\times=\mathrm C_A(X)\cap A^\times$.  
\end{lem}

\begin{proof}
    The first part is obvious. One has obviously $\mathrm{C}_A(X)^\times\subset \mathrm C_A(X)\cap A^\times$. If now $a\in\mathrm C_A(X)\cap A^\times$, then 
    for any $x\in X$, the relation $ax=xa$ implies, after left and right multiplication by $a^{-1}$, the relations $xa^{-1}=a^{-1}x$, therefore 
    $a^{-1}\in \mathrm C_A(X)$. Therefore $a\in \mathrm{C}_A(X)^\times$. 
\end{proof}

\begin{lem}\label{lem:13:5:1108:BIS}
(a) The product $(\phi,m)\bullet(\phi',m'):=(\phi\phi',\phi(e_1,f_1) m'+m\phi'(e_1,f_1)+m\cdot \mathrm{diag}(e_1+f_1,1)\cdot m')$
defines an algebra structure on $\mathbf k[[u,v]]\oplus M_2\hat V$. 

(b) The map $\mathbf k[[u,v]]\oplus M_2\hat V\to M_3\hat V$ given by 
\begin{equation}\label{def:M:2012}
(\phi,m)\mapsto M(\phi,m):=\phi(e_1,f_1)I_3+\begin{pmatrix}f_1&0\\e_1&0 \\0&1\end{pmatrix}
\cdot m\cdot \begin{pmatrix}1&1&0\\0&0&1\end{pmatrix},    
\end{equation}
where $I_3\in M_3\hat V$ is the identity matrix, defines a $\mathbf k$-algebra isomorphism  
\begin{equation}\label{iso:alg:comm:2012}
(\mathbf k[[u,v]]\oplus M_2\hat V,\bullet)\to\mathrm{C}_3(\rho_1). 
\end{equation}

(c) The map $\mathrm C_3(\rho_1)\to \mathbf k[[u,v]]$ taking $M\in\mathrm C_3(\rho_1)$ to the element $\phi$ such that 
there exists $m\in M_2\hat V$ with $M=M(\phi,m)$, defines an algebra morphism. 
\end{lem}

\begin{proof}
(a) is a direct verification. Let us show (b). 
One checks the said map $\mathbf k[[u,v]]\times M_2\hat V\to M_3\hat V$ to be a $\mathbf k$-algebra morphism.
It is given by 
\begin{equation}\label{form:of:commutant:2912}
(\phi,\begin{pmatrix}
   \alpha & \gamma \\ \beta& \delta 
\end{pmatrix})\mapsto \begin{pmatrix}
      f_1\alpha+\phi(e_1,f_1) & f_1\alpha&f_1\gamma \\ e_1\alpha & e_1\alpha+\phi(e_1,f_1)&e_1\gamma\\ \beta& \beta& \delta+\phi(e_1,f_1)
   \end{pmatrix}
\end{equation}
If $(\phi,\begin{pmatrix}
   \alpha & \gamma \\ \beta& \delta 
\end{pmatrix})$ belongs to its kernel, then $\beta=0$, the equalities $e_1\alpha=e_1\gamma=0$ imply $\alpha=\gamma=0$ by Lem. \ref{lem:TODO:BIS}(a); the equality   
$ f_1\alpha+\phi(e_1,f_1)=0$ then implies $\phi=0$; the equality $\delta+\phi(e_1,f_1)=0$ then implies $\delta=0$. Therefore the said map 
$\mathbf k[[u,v]]\times M_2\hat V\to M_3\hat V$ is injective. 

One checks the equalities $\rho_1\cdot \begin{pmatrix}f_1&0\\e_1&0 \\0&1\end{pmatrix}=0$ and 
$\begin{pmatrix}1&1&0\\0&0&1\end{pmatrix}\cdot \rho_1=0$, as well the the commutation of $e_1$ and $f_1$ with $\rho_1$. 
All this implies that the image of the said map $\mathbf k[[u,v]]\times M_2\hat V\to M_3\hat V$  is contained in $\mathrm C_3(\rho_1)$. 

Let us prove that the image of the said map $\mathbf k[[u,v]]\times M_2\hat V\to M_3\hat V$ is equal to $\mathrm C_3(\rho_1)$. 
Let $A\in \mathrm C_3(\rho_1)$. By $\rho_1=\mathrm{col}_{\mathrm{DT}}\cdot\mathrm{row}_{\mathrm{DT}}$
(see Def. \ref{def:5:2:1926}), one has
$$
(A\cdot \mathrm{col}_{\mathrm{DT}})\cdot\mathrm{row}_{\mathrm{DT}}
=\mathrm{col}_{\mathrm{DT}}\cdot(\mathrm{row}_{\mathrm{DT}}\cdot A). 
$$
By Lem. \ref{lem:TODO:BIS}(b), 
the latter equality implies 
\begin{equation}\label{interm:step:0301:BIS}
\exists\alpha\in\hat V \ | \
A\cdot \mathrm{col}_{\mathrm{DT}}=\mathrm{col}_{\mathrm{DT}}\cdot\alpha \quad\mathrm{and}\quad 
\mathrm{row}_{\mathrm{DT}}\cdot A=\alpha\cdot \mathrm{row}_{\mathrm{DT}}. 
\end{equation}
This statement implies
$$
\alpha\cdot (e_1+f_1)=\alpha\cdot\mathrm{row}_{\mathrm{DT}}\cdot\mathrm{col}_{\mathrm{DT}}
=\mathrm{row}_{\mathrm{DT}}\cdot A\cdot\mathrm{col}_{\mathrm{DT}}
=\mathrm{row}_{\mathrm{DT}}\cdot\mathrm{col}_{\mathrm{DT}}\cdot \alpha
=(e_1+f_1)\alpha, 
$$
where the first and last equalities follow from 
$\mathrm{row}_{\mathrm{DT}}\cdot\mathrm{col}_{\mathrm{DT}}=e_1+f_1$, 
and the two middle equalities follows from the equalities of \eqref{interm:step:0301:BIS}. 
The resulting equality $\alpha\cdot (e_1+f_1)=(e_1+f_1)\cdot\alpha$ implies, by Lem. \ref{lemma:comm:0304:BIS}, 
the existence of $\phi\in\mathbf k[[t,u]]$ such that $\alpha=\phi(e_1,f_1)$. 
The statement \eqref{interm:step:0301:BIS} therefore implies 
\begin{equation}\label{interm:step:next:0301:BIS}
\exists\phi\in\mathbf k[[t,u]] \ | \
A\cdot \mathrm{col}_{\mathrm{DT}}=\mathrm{col}_{\mathrm{DT}}\cdot\phi(e_1,f_1)\quad\mathrm{and}\quad 
\mathrm{row}_{\mathrm{DT}}\cdot A=\phi(e_1,f_1)\cdot \mathrm{row}_{\mathrm{DT}}. 
\end{equation}
Since the entries of $\mathrm{col}_{\mathrm{DT}}$ and $\mathrm{row}_{\mathrm{DT}}$ belong to 
the subalgebra $\mathbf k[[e_1,f_1]]\subset\hat V$ and since this subalgebra is commutative, one has the equalities 
$\mathrm{col}_{\mathrm{DT}}\cdot\phi(e_1,f_1)=\phi(e_1,f_1)\cdot\mathrm{col}_{\mathrm{DT}}$
and $\phi(e_1,f_1)\cdot \mathrm{row}_{\mathrm{DT}}=\mathrm{row}_{\mathrm{DT}}\cdot \phi(e_1,f_1)$, 
therefore \eqref{interm:step:next:0301:BIS} implies 
$$ 
\exists\phi\in\mathbf k[[t,u]]\ | \
(A-\phi(e_1,f_1)I_3)\cdot \mathrm{col}_{\mathrm{DT}}=0,
\quad\mathrm{and}\quad 
\mathrm{row}_{\mathrm{DT}}\cdot (A-\phi(e_1,f_1)I_3)=0,  
$$
which by Def. \ref{def:van:0301:BIS}(a) is equivalent to 
$$
\exists\phi\in\mathbf k[[t,u]]\ | \
A-\phi(e_1,f_1)I_3\in \mathrm{Ann}(\mathrm{col}_{\mathrm{DT}},\mathrm{row}_{\mathrm{DT}}), 
$$
which by Lem. \ref{lem:5:11:2008} is equivalent to 
$$
\exists m\in M_2(\hat V) \ | \ 
A=M(\phi,m). 
$$
Therefore $A$ belongs to the image of the said map $\mathbf k[[u,v]]\times M_2\hat V\to M_3\hat V$. This ends the proof of (b). 

 (c) One checks that the map $(\phi,m)\mapsto \phi$ induces an algebra morphism 
 $(\mathbf k[[u,v]]\times M_2\hat V,\bullet)\to \mathbf k[[u,v]]$. The said map is then the composition of this morphism 
 with the inverse of the algebra isomorphism from (b), and is therefore an algebra morphism. 
\end{proof}

\subsection{A diagram of groups $ \mathrm{GL}_3\hat V\leftarrow \mathrm C_3(\rho_1)^\times\to \mathbf k[[u,v]]^\times$}
\label{sec 5.3}

\begin{lem}\label{lem:15:12:BIS}
(a) One has 
$$
\mathrm C_3(\rho_1)^\times=\{M(\phi,m)|\phi\in\mathbf k[[u,v]]^\times,m=\begin{pmatrix}
    \alpha&\gamma\\\beta&\delta
\end{pmatrix}\in M_2\hat V,\phi(0,0)+\epsilon(\delta)\in\mathbf k^\times\}. 
$$

(b) The map $\mathrm C_3(\rho_1)^\times\to\mathbf k[[u,v]]^\times$ given by $M(\phi,m)\mapsto \phi$ 
is a group morphism, which together with the inclusion arising from the equality 
$\mathrm C_3(\rho_1)^\times=\mathrm C_3(\rho_1)\cap \mathrm{GL}_3\hat V$ (see Lem. \ref{invertible:commutant})
gives rise to a diagram of group morphisms
\begin{equation}\label{TAG:2112:FIRST}
\mathrm{GL}_3\hat V\hookleftarrow \mathrm C_3(\rho_1)^\times\to \mathbf k[[u,v]]^\times.     
\end{equation} 
\end{lem}

\begin{proof}
(a) One has 
\begin{align*}    
&\mathrm C_3(\rho_1)^\times=\mathrm C_3(\rho_1)\cap \mathrm{GL}_3\hat V
=\{M(\phi,m)|\phi\in\mathbf k[[u,v]],m\in M_2\hat V,\epsilon(M(\phi,m))\in\mathrm{GL}_3\mathbf k\}
\\ & =\{M(\phi,m)|\phi\in\mathbf k[[u,v]]^\times,m=\begin{pmatrix}
    \alpha&\gamma\\\beta&\delta
\end{pmatrix}\in M_2\hat V,\phi(0,0)+\epsilon(\delta)\in\mathbf k^\times\}
\end{align*}
where the first equality follows from Lem. \ref{invertible:commutant},  the second equality follows from the combination of 
Lem. \ref{lem:13:5:1108:BIS}(b) and 
\begin{equation}\label{charact:gl3}
\mathrm{GL}_3\hat V=\{P\in M_3\hat V|\epsilon(P)\in \mathrm{GL}_3\mathbf k\}, 
\end{equation}
where $\varepsilon : \hat V\to\mathbf k$ is the augmentation map, and the last equality follows from the identity 
$\epsilon(M(\phi,m))=\begin{pmatrix}
   \phi(0,0) &0&0\\0&\phi(0,0)&0\\\epsilon(\beta)&\epsilon(\beta)&\phi(0,0)+\epsilon(\delta)
\end{pmatrix}$. 
 
(b) The algebra morphism from Lem. \ref{lem:13:5:1108:BIS}(c) induces a group morphism $\mathrm C_3(\rho_1)^\times\to\mathbf k[[u,v]]^\times$. 
Since this algebra morphism is a composition $\mathrm C_3(\rho_1)\stackrel{\sim}{\to}\mathbf k[[u,v]]\times M_2\hat V \to\mathbf k[[u,v]]$, 
this group morphism is the map between groups of units induced by this map, which is  the announced map. 
\end{proof}

\subsection{A diagram of pointed sets with group actions}\label{sec 5.4}

\begin{lem}\label{lem:520:2212:FIRST}
(a) The map 
$$
\mathrm{GL}_3\hat V\times\mathrm{Hom}_{\mathcal C\operatorname{-alg}}(\hat{\mathcal V},M_3\hat V)
\to \mathrm{Hom}_{\mathcal C\operatorname{-alg}}(\hat{\mathcal V},M_3\hat V),\quad 
(P,\rho)\mapsto P\bullet \rho:=\mathrm{Ad}_P\circ \rho
$$
defines an action of the group $\mathrm{GL}_3\hat V$ on the set 
$\mathrm{Hom}_{\mathcal C\operatorname{-alg}}(\hat{\mathcal V},M_3\hat V)$. 

(b) The action from (a) restricts to an action of the group  $\mathrm C_3(\rho_1)^\times$ on the set 
$\mathrm{Hom}_{\mathcal C\operatorname{-alg}}^1(\hat{\mathcal V},M_3\hat V)$
(see  Def. \ref{defn: Hom1}).
\end{lem}

\begin{proof}
 (a) If $P\in \mathrm{GL}_3\hat V$, then $\mathrm{Ad}_P$ 
 is an algebra automorphism of $M_3\hat V$, therefore 
 for any $\rho\in \mathrm{Hom}_{\mathcal C\operatorname{-alg}}(\hat{\mathcal V},M_3\hat V)$, one has 
 $\mathrm{Ad}_P\circ \rho\in
 \mathrm{Hom}_{\mathcal C\operatorname{-alg}}(\hat{\mathcal V},M_3\hat V)$. 
 For $P,P'\in \mathrm{GL}_3\hat V$ and 
 $\rho\in\mathrm{Hom}_{\mathcal C\operatorname{-alg}}(\hat{\mathcal V},M_3\hat V)$, one has 
 $P\bullet(P'\bullet\rho)=\mathrm{Ad}_P\circ\mathrm{Ad}_{P'}\circ\rho=\mathrm{Ad}_{PP'}\circ\rho=(PP')\bullet\rho$. 
Therefore the said formula defines an action.  
 
(b) If $P\in \mathrm C_3(\rho_1)^\times$ and $\rho\in \mathrm{Hom}_{\mathcal C\operatorname{-alg}}^1(\hat{\mathcal V},M_3\hat V)$, then 
$$
(P\bullet\rho)(e_1)=\mathrm{Ad}_P\circ \rho(e_1)
=\mathrm{Ad}_P(\rho_1)=\rho_1, 
$$
where the first equality follows from the definition of the action, the second equality follows from 
$\rho\in \mathrm{Hom}_{\mathcal C\operatorname{-alg}}^1(\hat{\mathcal V},M_3\hat V)$, and the 
last equality follows from $P\in\mathrm C_3(\rho_1)^\times$; therefore $P\bullet\rho\in\mathrm{Hom}_{\mathcal C\operatorname{-alg}}^1(\hat{\mathcal V},M_3\hat V)$. 
\end{proof}

\begin{lem}\label{lem:5:21:2212}
    For any $\rho\in\mathrm{Hom}_{\mathcal C\operatorname{-alg}}^1(\hat{\mathcal V},M_3\hat V)$, one has 
    $$
    \forall a\in\hat{\mathcal V},\quad \Delta_\rho(ae_1)=\mathrm{row}_{\mathrm{DT}}\cdot 
    \rho(a)\cdot\mathrm{col}_{\mathrm{DT}}
    $$
    (equality in $\hat V$). 
\end{lem}

\begin{proof}
Define $\tilde\Delta_\rho$ to be the $\mathbf k$-module morphism $\hat{\mathcal W}\to\hat V$ given by $1\mapsto 1$
and $\alpha 1+ae_1\mapsto \alpha 1+
\mathrm{row}_{\mathrm{DT}}\cdot \rho(a)\cdot\mathrm{col}_{\mathrm{DT}}$ for any 
$\alpha\in\mathbf k,a\in\hat{\mathcal V}$. 
For $(\alpha,a),(\alpha',a')\in\mathbf k\times\hat{\mathcal V}$, one has 
\begin{align*}
    & \tilde\Delta_\rho((\alpha 1+ae_1)\cdot (\alpha'1+a'e_1))
=\tilde\Delta_\rho(\alpha\alpha'1+(\alpha a'+a\alpha'+a e_1 a')e_1)
=\alpha\alpha' 1+\mathrm{row}_{\mathrm{DT}}\cdot
\rho(\alpha a'+a\alpha'+a e_1 a')\cdot\mathrm{col}_{\mathrm{DT}}
\\ & 
=\alpha\alpha' 1
+\alpha \cdot \mathrm{row}_{\mathrm{DT}}\cdot\rho(a')\cdot \mathrm{col}_{\mathrm{DT}}
+ \mathrm{row}_{\mathrm{DT}}\cdot\rho(a)\cdot \mathrm{col}_{\mathrm{DT}}\cdot a'
+\mathrm{row}_{\mathrm{DT}}\cdot
\rho(a)\rho(e_1)\rho(a')\cdot\mathrm{col}_{\mathrm{DT}}
\\ & 
=\alpha\alpha' 1
+\alpha \cdot \mathrm{row}_{\mathrm{DT}}\cdot\rho(a')\cdot \mathrm{col}_{\mathrm{DT}}
+ \mathrm{row}_{\mathrm{DT}}\cdot\rho(a)\cdot \mathrm{col}_{\mathrm{DT}}\cdot a'
+\mathrm{row}_{\mathrm{DT}}\cdot
\rho(a)\cdot \mathrm{col}_{\mathrm{DT}} \cdot \mathrm{row}_{\mathrm{DT}}\cdot \rho(a')\cdot\mathrm{col}_{\mathrm{DT}}
\\ & =(\alpha 1+\mathrm{row}_{\mathrm{DT}}\cdot\rho(a)\cdot \mathrm{col}_{\mathrm{DT}})\cdot 
(\alpha' 1+\mathrm{row}_{\mathrm{DT}}\cdot\rho(a')\cdot \mathrm{col}_{\mathrm{DT}})
=\tilde\Delta_\rho(\alpha 1+ae_1)\cdot \tilde\Delta_\rho(\alpha'1+a'e_1),  
\end{align*}
where the first (resp. fifth) equalities follow from equalities in $\hat{\mathcal V}$ (resp. $\hat V$), 
the second and last equalities follow from the definition of $\tilde\Delta_\rho$, and the third (resp. fourth) 
equalities follow from the algebra morphism status of $\rho$ (resp. from $\rho(e_1)=\rho_1$). It follows that 
$\tilde\Delta_\rho$ is an algebra morphism. 

For any $n\geq1$, one has 
$$
\tilde\Delta_\rho(e_0^{n-1}e_1)=\mathrm{row}_{\mathrm{DT}}\cdot \rho(e_0^{n-1})\cdot\mathrm{col}_{\mathrm{DT}}
=\mathrm{row}_{\mathrm{DT}}\cdot \rho(e_0)^{n-1}\cdot\mathrm{col}_{\mathrm{DT}}=\Delta_\rho(e_0^{n-1}e_1)
$$
where the first (resp. second, third) equality follows from the definition of $\tilde\Delta_\rho$
(resp. the algebra morphism status of $\rho$, the definition of $\Delta_\rho$). 

This implies that both $\tilde\Delta_\rho$ and $\Delta_\rho$ are algebra morphisms 
$\hat{\mathcal W}\to \hat V$ which coincide on 
$(e_0^{n-1}e_1)_{n\geq1}$, which is a generating family of $\hat{\mathcal W}$, which implies that they are equal. 
\end{proof}

\begin{lem}\label{lem:5:22:2212:FIRST}
    (a) The map $\mathrm{Hom}_{\mathcal C\operatorname{-alg}}^1(\hat{\mathcal V},M_3\hat V)\to 
    \mathrm{Hom}_{\mathcal C\operatorname{-alg}}(\hat{\mathcal W},\hat V)$, $\rho\mapsto \Delta_\rho$ of sets 
    (see Lem. \ref{lem:vania}(a)) is compatible  with the group morphism $\mathrm C_3(\rho_1)^\times\to\mathbf k[[u,v]]^\times$ 
    (see Lem. \ref{lem:15:12:BIS}(b)) and with the actions from Lem. \ref{lem:520:2212:FIRST}(b) and Def. \ref{lem28:1001}(c). 

    (b) The canonical inclusion $\mathrm{Hom}_{\mathcal C\operatorname{-alg}}^1(\hat{\mathcal V},M_3\hat V)\to
    \mathrm{Hom}_{\mathcal C\operatorname{-alg}}(\hat{\mathcal V},M_3\hat V)$ (see \eqref{TAG':2112}) 
    is compatible with the group inclusion $\mathrm C_3(\rho_1)^\times\to
    \mathrm{GL}_3\hat V$ (see \eqref{TAG:2112:FIRST}) 
    and with the actions of Lem. \ref{lem:520:2212:FIRST}, (a) and (b). 

     (c) The diagram \eqref{TAG':2112} of pointed sets is compatible with the diagram \eqref{TAG:2112:FIRST} of groups. 
\end{lem}

\begin{proof}
(a)  Let $\rho\in \mathrm{Hom}_{\mathcal C\operatorname{-alg}}^1(\hat{\mathcal V},M_3\hat V)$ and 
$P\in \mathrm C_3(\rho_1)^\times$. By Lem. \ref{lem:13:5:1108:BIS}(b), 
there exists $(\phi,m)\in \mathbf k[[u,v]]^\times\times M_2\hat V$
such that $P=M(\phi,m)$, and by Lem. \ref{lem:15:12:BIS}(b), the image of $(P,1)$ under the morphism 
$\mathrm C_3(\rho_1)^\times\rtimes\mathcal G\to\mathbf k[[u,v]]^\times\times\mathcal G$ 
is $(\phi(e_1,f_1),1)$. Then 
\begin{equation}\label{eq:col:2212}
   P\cdot \mathrm{col}_{\mathrm{DT}}=(\phi(e_1,f_1)I_3+\begin{pmatrix}f_1&0\\e_1&0 \\0&1\end{pmatrix}\cdot 
m\cdot \begin{pmatrix}1&1&0\\0&0&1\end{pmatrix})\cdot\mathrm{col}_{\mathrm{DT}}
=\mathrm{col}_{\mathrm{DT}}\cdot \phi(e_1,f_1),  
\end{equation}
and 
\begin{equation}\label{eq:row:2212}
\mathrm{row}_{\mathrm{DT}}\cdot P=\mathrm{row}_{\mathrm{DT}}\cdot 
(\phi(e_1,f_1)I_3+\begin{pmatrix}f_1&0\\e_1&0 \\0&1\end{pmatrix}\cdot m\cdot \begin{pmatrix}1&1&0\\0&0&1\end{pmatrix})
=\phi(e_1,f_1)\cdot \mathrm{row}_{\mathrm{DT}}, 
\end{equation}
where the first equalities follow from $P=M(\phi,m)$ and the second equalities follow from the 
commutation of the entries of $\mathrm{col}_{\mathrm{DT}},\mathrm{row}_{\mathrm{DT}}$
with $e_1,f_1$, and from the equalities $\begin{pmatrix}1&1&0\\0&0&1\end{pmatrix}
\cdot \mathrm{col}_{\mathrm{DT}}=0$ and 
$\mathrm{row}_{\mathrm{DT}}\cdot \begin{pmatrix}f_1&0\\e_1&0 \\0&1\end{pmatrix}=0$. 

For any $n\geq 1$, one then has 
\begin{align*}
& \Delta_{P\bullet\rho}(e_0^{n-1}e_1)
=\Delta_{\mathrm{Ad}_P \circ \rho}(e_0^{n-1}e_1)
=\mathrm{row}_{\mathrm{DT}}\cdot ((\mathrm{Ad}_P\circ \rho)(e_0))^{n-1}\cdot \mathrm{col}_{\mathrm{DT}}
\\ & =\mathrm{row}_{\mathrm{DT}}\cdot (P\cdot \rho(e_0)\cdot P^{-1})^{n-1}\cdot \mathrm{col}_{\mathrm{DT}}
=\mathrm{row}_{\mathrm{DT}}\cdot P\cdot \rho(e_0)^{n-1}\cdot P^{-1}\cdot \mathrm{col}_{\mathrm{DT}}
\\& =\phi(e_1,f_1)\cdot \mathrm{row}_{\mathrm{DT}}\cdot \rho(e_0)^{n-1}\cdot  \mathrm{col}_{\mathrm{DT}}
\cdot \phi(e_1,f_1)^{-1}
\\& =\phi(e_1,f_1)\cdot \Delta_\rho(e_0^{n-1}e_1)
\cdot \phi(e_1,f_1)^{-1}
=(\mathrm{Ad}_{\phi(e_1,f_1)}\circ \Delta_\rho)(e_0^{n-1}e_1)
=(\phi(e_1,f_1)\bullet\Delta_\rho)(e_0^{n-1}e_1). 
\end{align*}
where all the equalities follow from definitions, except for the fifth one, 
which follows from \eqref{eq:col:2212} and 
\eqref{eq:row:2212}. 
This implies the identity 
\begin{equation}\label{id:equiv:a:2212}
 \forall \rho\in\mathrm{Hom}_{\mathcal C\operatorname{-alg}}^1(\hat{\mathcal V},M_3\hat V),
\forall P\in\mathrm C_3(\rho_1)^\times,\quad \Delta_{P\bullet\rho}=\phi(e_1,f_1)\bullet\Delta_\rho.    
\end{equation}
(equality in $\mathrm{Hom}_{\mathcal C\operatorname{-alg}}^1(\hat{\mathcal W},\hat V)$). 
(b) is obvious.  (c) is a direct consequence of (a) and (b).  
\end{proof}

\begin{rem}
 Lem. \ref{lem:vania} could be alternatively proved a follows: one computes (see Rem. \ref{rem:5:7:3012})  $U=M(1,-\begin{pmatrix}
    0\\1
\end{pmatrix}\begin{pmatrix}
    1&0
\end{pmatrix})$, therefore $U\in\mathrm{ker}(\mathrm{C}_3(\rho_1)^\times\to\mathbf k[[u,v]]^\times)$. 
The relations $\rho_{\mathrm{DT}},\tilde\rho_{\mathrm{DT}}\in 
\mathrm{Hom}_{\mathcal C\operatorname{-alg}}^1(\hat{\mathcal V},M_3\hat V)$ and 
$\tilde\rho_{\mathrm{DT}}=\mathrm{Ad}_U\circ \rho_{\mathrm{DT}}$, together with 
Lem. \ref{lem:5:22:2212:FIRST}(c) then imply that the image by $\rho\mapsto\Delta_\rho$ of 
$\rho_{\mathrm{DT}}$ is equal to that of $\tilde\rho_{\mathrm{DT}}$, which is 
$\Delta^{\mathcal W}_{r,l}$ by Lem. 6.2 in \cite{EF1}. 
\end{rem}

\subsection{Overall action of $\mathcal G$: action of $\mathcal G$ on a group diagram}\label{sec 5.5}

The morphisms of algebras from \S\ref{lem:4:2:2012} fit in the following diagram 
\begin{equation}\label{sequence:alg:2012}
 M_3\hat V\leftarrow \mathrm C_3(\rho_1)\stackrel{\sim}{\leftarrow}\mathbf k[[u,v]]\oplus M_2\hat V\to 
 \mathbf k[[u,v]]. 
\end{equation}
In the following lemma, we equip the algebras of this diagram with $\mathcal G$-actions, which are compatible with the 
morphisms. 

\begin{defn}\label{def:5:14:0105} 
    For $g\in\mathcal G$, set $\mathrm{aut}_g^{V}:=(\mathrm{aut}_g^{\mathcal V})^{\otimes 2}$
    ($\mathrm{aut}_g^{\mathcal V}$ being as in \eqref{def:aut:g:V:(1)})). 
\end{defn}

\begin{lem}\label{lem:alg:morphisms}
    (a) The formula $g*P:=\mathrm{aut}_g^V(P)$ 
defines an action of $\mathcal G$ on the algebra $M_3\hat V$.  

(b) This restricts to an action on the subalgebra $\mathrm C_3(\rho_1)\subset M_3\hat V$. 

(c) The formula $g*(\phi,m):=(\phi,g*m)$ where $g*m:=(\mathrm{aut}_g^{\mathcal V,(1)})^{\otimes 2}(m)$
defines an action of $\mathcal G$ on $(\mathbf k[[u,v]]\oplus M_2\hat V,\bullet)$. 

(d) The algebra isomorphism $(\mathbf k[[u,v]]\oplus M_2\hat V,\bullet)\to \mathrm C_3(\rho_1)$, $(\phi,m)\mapsto M(\phi,m)$
is $\mathcal G$-equivariant. 

(e) The projection map $(\mathbf k[[u,v]]\oplus M_2\hat V,\bullet)\to\mathbf k[[u,v]]$ is a $\mathcal G$-equivariant 
algebra morphism, the action on $\mathcal G$ on the target being trivial. 
\end{lem}

\begin{proof}
(a) follows from the fact that the said formula expresses the tensor product of the action of 
the group $\mathcal G$ on $\hat V$ with its trivial action on $M_3(\mathbf k)$. (b) follows from 
the invariance of $\rho_1$ under the action of $\mathcal G$. Let us show (c). The said formula
obviously defines an action of $\mathcal G$ on $\mathbf k[[u,v]]\times M_2\hat V$ by $\mathbf k$-module
isomorphisms. For $(\phi,m),(\phi',m')\in \mathbf k[[u,v]]\times M_2\hat V$, one has 
\begin{align*}
&(g*(\phi,m))\bullet(g*(\phi',m'))=(\phi,g* m)\bullet(\phi',g* m')
\\ & =(\phi\phi',\phi\cdot(g* m')+(g* m)\cdot\phi'+(g* m)\cdot\mathrm{diag}(e_1+f_1,1)\cdot (g* m'))
\\ & 
=\left(\phi\phi',g* (\phi\cdot m')+g* (m\cdot\phi')+g* (m\cdot\mathrm{diag}(e_1+f_1,1)\cdot m')\right)
=g*\left((\phi,m)\bullet(\phi',m')\right). 
\end{align*}
where the equalities follow from the definitions of the action of $\mathcal G$ and of the product on $\mathbf k[[u,v]]\times M_2\hat V$), 
and the third equality follows from the $\mathcal G$-invariance of $\phi$ and $\phi'$, from the compatibility of the action of 
$\mathcal G$ on $M_2\hat V$ with the algebra structure, and from the invariance under the action of $\mathcal G$ of $\mathrm{diag}(e_1+f_1,1)$. 
Therefore $\mathcal G$ acts by algebra automorphisms, implying (c). For $(\phi,m)\in \mathbf k[[u,v]]\times M_2\hat V$ and $g\in\mathcal G$, one has 
\begin{align*}
&  g* M(\phi,m)=g*(\phi(e_1,f_1)I_3+\begin{pmatrix}f_1&0\\e_1&0 \\0&1\end{pmatrix}\cdot m\cdot \begin{pmatrix}1&1&0\\0&0&1\end{pmatrix})
\\ & =\phi(e_1,f_1)I_3+\begin{pmatrix}f_1&0\\e_1&0 \\0&1\end{pmatrix}\cdot (g* m)\cdot \begin{pmatrix}1&1&0\\0&0&1\end{pmatrix}=M(\phi,g* m)=M(g*(\phi,m))  
\end{align*}
where the second equality follows from the $\mathcal G$-invariance of $e_1,f_1$, and from the 
compatiblity of the action of $\mathcal G$ with the product of matrices, and the other equalities
follow from definitions. This implies (d). (e) follows from the explicit expression of the action $\mathcal G$ on 
$\mathbf k[[u,v]]\times M_2\hat V$. 
\end{proof}

\begin{lem}\label{lem:5:18:2012}
 The sequence \eqref{sequence:alg:2012} gives rise to a sequence 
 \begin{equation}\label{sequence:gps:2012}
 \mathrm{GL}_3\hat V\leftarrow \mathrm C_3(\rho_1)^\times\stackrel{\sim}{\leftarrow}\mathbf k[[u,v]]^\times\times M_2\hat V\to \mathbf k[[u,v]]^\times.    
\end{equation}
of groups with actions of $\mathcal G$. 
\end{lem}

\begin{proof}
This follows from the fact that upon taking units, a morphism of algebras with actions of $\mathcal G$ 
gives rise to a morphism of groups
with actions of $\mathcal G$, from the equality $\mathrm{GL}_3\hat V=(M_3\hat V)^\times$, and 
from the identification of 
$\mathbf k[[u,v]]^\times\times M_2\hat V$ with the group of units of $(\mathbf k[[u,v]]\times M_2\hat V,\bullet)$ 
(see Lem. \ref{lem:15:12:BIS}). 
\end{proof}

\subsection{Overall action of $\mathcal G$: action of $\mathcal G$ on a set diagram}\label{sec 5.6}

\begin{lem}\label{lem:520:2212:BIS}
(a) The map 
$$
\mathcal G\times\mathrm{Hom}_{\mathcal C\operatorname{-alg}}(\hat{\mathcal V},M_3\hat V)
\to \mathrm{Hom}_{\mathcal C\operatorname{-alg}}(\hat{\mathcal V},M_3\hat V),$$  
$$
(g,\rho)\mapsto g* \rho:=\mathrm{aut}_g^V\circ \rho\circ 
(\mathrm{aut}_g^{\mathcal V})^{-1} 
$$
defines an action $*$ of the group $\mathcal G$ on the set 
$\mathrm{Hom}_{\mathcal C\operatorname{-alg}}(\hat{\mathcal V},M_3\hat V)$. 

(b) The action from (a) restricts to an action of $\mathcal G$ on the set 
$\mathrm{Hom}_{\mathcal C\operatorname{-alg}}^1(\hat{\mathcal V},M_3\hat V)$
(see  Def. \ref{defn: Hom1}); the inclusion $\mathrm{Hom}_{\mathcal C\operatorname{-alg}}^1(\hat{\mathcal V},M_3\hat V)
\hookrightarrow \mathrm{Hom}_{\mathcal C\operatorname{-alg}}(\hat{\mathcal V},M_3\hat V)$ is then $\mathcal G$-equivariant.
\end{lem}

\begin{proof}
(a)  If $g\in \mathcal G$, then 
 $\mathrm{aut}_g^{\mathcal V}$ (resp. $\mathrm{aut}_g^{V}$) 
 is an algebra automorphism of $\hat{\mathcal V}$ (resp. $M_3\hat V$), therefore 
 for any $\rho\in \mathrm{Hom}_{\mathcal C\operatorname{-alg}}(\hat{\mathcal V},M_3\hat V)$, one has 
 $\mathrm{aut}_g^V\circ \rho\circ (\mathrm{aut}_g^{\mathcal V})^{-1}\in
 \mathrm{Hom}_{\mathcal C\operatorname{-alg}}(\hat{\mathcal V},M_3\hat V)$. 
 For $g,g'\in \mathcal G$ and 
 $\rho\in\mathrm{Hom}_{\mathcal C\operatorname{-alg}}(\hat{\mathcal V},M_3\hat V)$, one has 
$$
g*(g'*\rho)
 =\mathrm{aut}_g^V\circ (\mathrm{aut}_{g'}^V\circ \rho\circ 
(\mathrm{aut}_{g'}^{\mathcal V})^{-1})\circ 
 (\mathrm{aut}_g^{\mathcal V})^{-1}
=\mathrm{aut}_{g\circledast g'}^V\circ \rho\circ 
(\mathrm{aut}_{g\circledast g'}^{\mathcal V})^{-1}    
=(g\circledast g')*\rho, 
$$
which follows from the facts that
 $g\mapsto \mathrm{aut}_g^{\mathcal V}$ (resp. $g\mapsto \mathrm{aut}_g^{V}$) 
 defines a group morphism $\mathcal G\to\mathrm{Aut}_{\mathcal C\operatorname{-alg}}
 (\hat{\mathcal V})$ (resp. $\mathcal G\to\mathrm{Aut}_{\mathcal C\operatorname{-alg}}
 (\hat V)$, $\mathrm{GL}_3\hat V\to\mathrm{Aut}_{\mathcal C\operatorname{-alg}}(M_3\hat V)$). Therefore the 
 said formula defines an action. 

 (b) If $g\in \mathcal G$ and $\rho\in \mathrm{Hom}_{\mathcal C\operatorname{-alg}}^1(\hat{\mathcal V},M_3\hat V)$, then 
$$
(g*\rho)(e_1)=\mathrm{aut}_g^V\circ \rho\circ 
(\mathrm{aut}_g^{\mathcal V})^{-1}(e_1)
=\mathrm{aut}_g^V\circ \rho(e_1)
=\mathrm{aut}_g^V(\rho_1)
=\rho_1, 
$$
where the first equality follows from the definition of the action, the second equality follows from the 
invariance of $e_1$ under the action of $\mathcal G$, the third equality follows from 
$\rho\in \mathrm{Hom}_{\mathcal C\operatorname{-alg}}^1(\hat{\mathcal V},M_3\hat V)$, 
the fourth equality follows from the invariance of $\rho_1$ under the action of $\mathcal G$. 
\end{proof}

\begin{lem}\label{lem:5:22:2212:BIS}
    (a) The map $\mathrm{Hom}_{\mathcal C\operatorname{-alg}}^1(\hat{\mathcal V},M_3\hat V)\to 
    \mathrm{Hom}_{\mathcal C\operatorname{-alg}}(\hat{\mathcal W},\hat V)$, $\rho\mapsto \Delta_\rho$  
    (see Lem. \ref{lem:vania}(a)) is $\mathcal G$-equivariant, the actions of $\mathcal G$ on the source and 
    target being respectively as in Lems. \ref{lem:520:2212:BIS} and \ref{lem:43:0401:TER}(b). 

     (b) The diagram \eqref{TAG':2112} of pointed sets is $\mathcal G$-equivariant. 
\end{lem}

\begin{proof} 
(a) Let $\rho\in \mathrm{Hom}_{\mathcal C\operatorname{-alg}}^1(\hat{\mathcal V},M_3\hat V)$ and 
$g\in\mathcal G$. For any $a\in\hat{\mathcal V}$, one has 
\begin{align*}
    &
    \Delta_{g*\rho}(ae_1)=\mathrm{row}_{\mathrm{DT}}\cdot (g*\rho)(a)\cdot \mathrm{col}_{\mathrm{DT}}
=\mathrm{row}_{\mathrm{DT}}\cdot (\mathrm{aut}_g^V\circ\rho\circ 
(\mathrm{aut}_g^{\mathcal V})^{-1})(a)\cdot \mathrm{col}_{\mathrm{DT}}
    \\ & 
    = \mathrm{aut}_g^V( \mathrm{row}_{\mathrm{DT}}\cdot ( 
\rho  ((\mathrm{aut}_g^{\mathcal V})^{-1}(a))\cdot \mathrm{col}_{\mathrm{DT}} ) 
=\mathrm{aut}_g^V(\Delta_\rho\left((\mathrm{aut}_g^{\mathcal V})^{-1}(a)e_1\right))
\\& =\mathrm{aut}_g^V(\Delta_\rho((\mathrm{aut}_g^{\mathcal W})^{-1}(ae_1)))
=\mathrm{aut}_g^V\circ \Delta_\rho\circ (\mathrm{aut}_g^{\mathcal W})^{-1}(ae_1)
=(g*\Delta_\rho)(ae_1), 
\end{align*}
where the first and fourth equalities follow from Lem. \ref{lem:5:21:2212}, 
the second (resp. last) equality follows from the definition of the action of $\mathcal G$ on 
$\mathrm{Hom}_{\mathcal C\operatorname{-alg}}^1(\hat{\mathcal V},M_3\hat V)$ (resp. 
$\mathrm{Hom}_{\mathcal C\operatorname{-alg}}^1(\hat{\mathcal W},\hat V)$), 
the third equality follows from the invariance of $\mathrm{col}_{\mathrm{DT}}$ and 
$\mathrm{row}_{\mathrm{DT}}$ under the action of 
$\mathcal G$, and the fifth equality follows from the definition of $\mathrm{aut}_g^{\mathcal W}$. 
This implies the wanted identity
$$ 
\forall \rho\in \mathrm{Hom}_{\mathcal C\operatorname{-alg}}^1(\hat{\mathcal V},M_3\hat V),\forall 
g\in\mathcal G, \quad \Delta_{g*\rho}=g*\Delta_\rho. 
$$

(b) is a direct consequence of (a) and Lem. \ref{lem:520:2212:BIS}(b).   
\end{proof}

\subsection{The group inclusion 
$\mathrm{Stab}_{\mathcal G}(\mathrm{GL}_3\hat V\bullet \rho_{\mathrm{DT}}
)\subset \mathrm{Stab}_{\mathcal G}(
\mathbf k[[u,v]]^\times\bullet\Delta^{\mathcal W}_{r,l}
)$}\label{sect:5:7}

\begin{lem}\label{lem:pre:psga}
    (a) The actions $*$ of $\mathcal G$ on the group $\mathrm{GL}_3\hat V$ and the set $\mathrm{Hom}_{\mathcal C\operatorname{-alg}}(\hat{\mathcal V},M_3\hat V)$ (see Lem.
    \ref{lem:5:18:2012} and    Lem. \ref{lem:520:2212:BIS}(a)), 
    and the action $\bullet$ of this group on this set (see Lem. \ref{lem:520:2212:FIRST}(a)), satisfy the compatibility relation from Def. \ref{def:psga}(a).   

    (b) The actions $*$ of $\mathcal G$ on the group $\mathrm C(\rho_1)^\times$ and the set 
    $\mathrm{Hom}^1_{\mathcal C\operatorname{-alg}}(\hat{\mathcal V},M_3\hat V)$ (see Lem.
    \ref{lem:5:18:2012} and    Lem. \ref{lem:520:2212:BIS}(b)), 
    and the action $\bullet$ of this group on this set (see Lem. \ref{lem:520:2212:FIRST}(b)), satisfy the 
    compatibility relation from Def. \ref{def:psga}(a).   
\end{lem}

\begin{proof}
(a) If $g\in\mathcal G$ and $P\in\mathrm{GL}_3\hat V$, and $\rho\in \mathrm{Hom}_{\mathcal C\operatorname{-alg}}(\hat{\mathcal V},M_3\hat V)$, then  
 \begin{equation*}
g*(P\bullet\rho)
 =\mathrm{aut}_g^V\circ (\mathrm{Ad}_{P}\circ \rho)\circ 
 (\mathrm{aut}_g^{\mathcal V})^{-1}   
=\mathrm{Ad}_{g*P}\circ \mathrm{aut}_{g}^V\circ 
\rho\circ 
(\mathrm{aut}_{g}^{\mathcal V})^{-1} =(g*P)\bullet(g*\rho), 
\end{equation*}
which follows from the 
 identity $\mathrm{aut}_g^V\circ\mathrm{Ad}_P=\mathrm{Ad}_{g*P}\circ \mathrm{aut}_g^V$. 
(b) follows from (a) by specialization.  
\end{proof}

\begin{prop}
    (a) The tuples
    $$
    (\mathrm{Hom}_{\mathcal C\operatorname{-alg}}(\hat{\mathcal V},M_3\hat V),\Delta_{r,l}^{\mathcal W},
    \mathrm{GL}_3\hat V,\bullet,*)\quad\text{and}\quad
    (\mathrm{Hom}^1_{\mathcal C\operatorname{-alg}}(\hat{\mathcal V},M_3\hat V),\Delta_{r,l}^{\mathcal W},
    \mathrm{C}_3(\rho_1)^\times,\bullet,*), 
    $$
where the actions $\bullet$ are as  Lem. \ref{lem:520:2212:FIRST} and the actions $*$ are as in Lems.
\ref{lem:5:18:2012} and \ref{lem:520:2212:BIS}, are objects in $\mathcal G\operatorname{-}\mathbf{PSGA}$. 

    (b) The pairs $(\mathrm{Hom}^1_{\mathcal C\operatorname{-alg}}(\hat{\mathcal V},M_3\hat V)\hookrightarrow
    \mathrm{Hom}_{\mathcal C\operatorname{-alg}}(\hat{\mathcal V},M_3\hat V),\mathrm{C}_3(\rho_1)^\times\hookrightarrow\mathrm{GL}_3\hat V)$ 
    and $(\rho\mapsto\Delta_\rho,\mathrm{C}_3(\rho_1)^\times\to \mathbf k[[u,v]]^\times)$, where the maps are as in Lem. 
    \ref{lem:vania}(a) and Lem. \ref{lem:5:18:2012}, build up the following diagram of morphisms 
 \begin{equation}\label{diag:bm:raphael}
   (\mathrm{Hom}_{\mathcal C\operatorname{-alg}}(\hat{\mathcal V},M_3\hat V),\rho_{\mathrm{DT}},
    \mathrm{GL}_3\hat V,\bullet,*)\leftarrow(\mathrm{Hom}^1_{\mathcal C\operatorname{-alg}}(\hat{\mathcal V},M_3\hat V),\rho_{\mathrm{DT}},
    \mathrm{C}_3(\rho_1)^\times,\bullet,*)
    \to(\mathbf H,\mathbf k[[u,v]]^\times,\bullet,*)  
 \end{equation}
    in $\mathcal G\operatorname{-}\mathbf{PSGA}$, where the last object is as in Lem. \ref{lem:2:15:toto}(b). 
\end{prop}

\begin{proof}
(a) follows from Lem. \ref{lem:pre:psga}. (b) follows from Lems. \ref{lem:5:22:2212:BIS}(b) and \ref{lem:5:22:2212:FIRST}(c). 
\end{proof}

\begin{defn}\label{def:quotient:diagram}
Define 
\begin{align*}
&   (\mathrm{GL}_3\hat V\backslash\mathrm{Hom}_{\mathcal C\operatorname{-alg}}(\hat{\mathcal V},M_3\hat V),
\mathrm{GL}_3\hat V\bullet \rho_{\mathrm{DT}}
,*)
   \leftarrow(\mathrm{C}_3(\rho_1)^\times\backslash\mathrm{Hom}^1_{\mathcal C\operatorname{-alg}}(\hat{\mathcal V},M_3\hat V),\mathrm{C}_3(\rho_1)^\times\bullet \rho_{\mathrm{DT}}
   ,*)
 \\ &    \to(\mathbf k[[u,v]]^\times\backslash\mathrm{Hom}_{\mathcal C\operatorname{-alg}}(\hat{\mathcal W},\hat V),
 \mathbf k[[u,v]]^\times\bullet\Delta^{\mathcal W}_{r,l}
 ,*)    
\end{align*}
to be the diagram in $\mathbf{PS}_{\mathcal G}$ obtained by applying the functor $\mathbf q$ to diagram \eqref{diag:bm:raphael}. 
\end{defn}

\begin{notation}\label{notation:F}
    Denote by
    $$
    (A) : (\mathrm{C}_3(\rho_1)^\times\backslash
    \mathrm{Hom}^1_{\mathcal C\operatorname{-alg}}(\hat{\mathcal V},M_3\hat V),
    \mathrm C_3(\rho_1)^\times\bullet \rho_{\mathrm{DT}}
    ,*)\to
    (\mathbf k[[u,v]]^\times\backslash\mathrm{Hom}_{\mathcal C\operatorname{-alg}}(\hat{\mathcal W},\hat V),
   \mathbf k[[u,v]]^\times\bullet\Delta^{\mathcal W}_{r,l}
   ,*) 
    $$
    the morphism from Def. \ref{def:quotient:diagram}. 
\end{notation}

\begin{lem}\label{lem:loc:inj:2122}
The morphism 
\begin{equation}\label{defs:[[rho]]}
(\mathrm C_3(\rho_1)^\times\backslash\mathrm{Hom}_{\mathcal C\operatorname{-alg}}^1(\hat{\mathcal V},M_3\hat V),
\mathrm C_3(\rho_1)^\times\bullet \rho_{\mathrm{DT}}
)\to 
(\mathrm{GL}_3\hat V\backslash\mathrm{Hom}_{\mathcal C\operatorname{-alg}}(\hat{\mathcal V},M_3\hat V),
\mathrm{GL}_3\hat V\bullet \rho_{\mathrm{DT}}
)
\end{equation}
in $\mathbf{PS}$ from Def. \ref{def:quotient:diagram} is locally injective. 
\end{lem}

\begin{proof}
We will prove the injectivity of the map $\mathrm C_3(\rho_1)^\times\backslash\mathrm{Hom}_{\mathcal C\operatorname{-alg}}^1(\hat{\mathcal V},M_3\hat V)
\to \mathrm{GL}_3\hat V\backslash\mathrm{Hom}_{\mathcal C\operatorname{-alg}}(\hat{\mathcal V},M_3\hat V)$, which implies the claimed local injectivity. 
Let $\alpha,\beta$ belong to the source of this map, with equal images in its target. Let $\tilde\alpha,\tilde\beta\in \mathrm{Hom}_{\mathcal C\operatorname{-alg}}^1(\hat{\mathcal V},M_3\hat V)$
be representatives of $\alpha,\beta$, so $\alpha = \mathrm C_3(\rho_1)^\times\bullet \tilde\alpha$ and $\beta = \mathrm C_3(\rho_1)^\times\bullet \tilde\beta$. 
It follows from the equality of the images of $\alpha$ and $\beta$ that for some $P\in\mathrm{GL}_3\hat V$, one has $\tilde\beta=\mathrm{Ad}_P\circ\tilde\alpha$. This implies the middle
equality in 
$$
\rho_1=\tilde\beta(e_1)=\mathrm{Ad}_P\circ\tilde\alpha(e_1)=\mathrm{Ad}_P(\rho_1), 
$$
where the first and last equalities follow from $\tilde\alpha,\tilde\beta\in \mathrm{Hom}_{\mathcal C\operatorname{-alg}}^1(\hat{\mathcal V},M_3\hat V)$. 
The latter equality implies $P\in \mathrm C_3(\rho_1)^\times$, therefore $\alpha=\beta$. 
\end{proof}

\begin{cor}\label{NEWCOR}
    One has 
$$
\mathrm{Stab_{\mathcal G}}(
\mathrm{GL}_3\hat V\bullet \rho_{\mathrm{DT}}
)
=\mathrm{Stab_{\mathcal G}}(\mathrm{C}_3(\rho_1)^\times\bullet \rho_{\mathrm{DT}}
)
$$
(equality of subgroups of $\mathcal G$). 
\end{cor}

\begin{proof}
This follows from Lem. \ref{lem:loc:inj:2122} and Lem. \ref{lem:general:2212}(b).    
\end{proof}

\begin{thm}\label{thm:5:31:3103} (see Thm. \ref{thm:015})
The stabilizer groups of the pointed $\mathcal G$-sets 
$$
(\mathrm{GL}_3\hat V\backslash\mathrm{Hom}_{\mathcal C\operatorname{-alg}}(\hat{\mathcal V},M_3\hat V),
\mathrm{GL}_3\hat V\bullet \rho_{\mathrm{DT}}
)
$$ 
and $(\mathbf k[[u,v]]^\times\backslash\mathrm{Hom}_{\mathcal C\operatorname{-alg}}
(\hat{\mathcal W},\hat V),\mathbf k[[u,v]]^\times\bullet\Delta^{\mathcal W}_{r,l}
)$ 
(see Def. \ref{def:quotient:diagram}) satisfy the inclusion 
$$
\mathrm{Stab_{\mathcal G}}(\mathrm{GL}_3\hat V\bullet \rho_{\mathrm{DT}}
)
\subset
\mathrm{Stab_{\mathcal G}}(\mathbf k[[u,v]]^\times\bullet\Delta^{\mathcal W}_{r,l}
) 
$$
of subgroups of $\mathcal G$. 
\end{thm}

\begin{proof}
According to Lem. \ref{lem:general:2212}(a), the morphisms of pointed $\mathcal G$-sets from Def. \ref{def:quotient:diagram} induce a diagram of inclusions 
$$
\mathrm{Stab}_{\mathcal G}
(\mathrm{GL}_3\hat V\bullet \rho_{\mathrm{DT}}
)\supset\mathrm{Stab}_{\mathcal G}(
\mathrm{C}_3(\rho_1)^\times\bullet \rho_{\mathrm{DT}}
)\subset
\mathrm{Stab}_{\mathcal G}(\mathbf k[[u,v]]^\times\bullet\Delta^{\mathcal W}_{r,l}
)
$$ 
of subgroups of 
$\mathcal G$. The local injectivity proved in Lem. \ref{lem:loc:inj:2122} implies, together with Lem. \ref{lem:general:2212}(b), that the left inclusion is an equality. 
The combination of this with the right inclusion yields the result. 
\end{proof}

\newpage

\part{Equality between $\mathrm{Stab}_{\mathcal G}(\mathrm{GL}_3\hat V\bullet \rho_{\mathrm{DT}}
)$ 
and $\mathrm{Stab}_{\mathcal G}(\mathbf k[[u,v]]^\times\bullet\Delta^{\mathcal W}_{r,l}
)$}
\label{part 3}

The objective of Part \ref{part 3} is to prove that the inclusion of 
$\mathrm{Stab}_{\mathcal G}(\mathrm{GL}_3\hat V\bullet \rho_{\mathrm{DT}}
)$ 
in $\mathrm{Stab}_{\mathcal G}(\mathbf k[[u,v]]^\times\bullet\Delta^{\mathcal W}_{r,l}
)$ obtained in Part \ref{part 2} is in fact an equality 
(Cor. \ref{cor:11:2:17apr}). The idea of this proof is as follows. By Cor. \ref{NEWCOR}, 
$\mathrm{Stab}_{\mathcal G}(\mathrm{GL}_3\hat V\bullet \rho_{\mathrm{DT}}
)$ is equal to 
$\mathrm{Stab}_{\mathcal G}(\mathrm{C}_3(\rho_1)^\times \bullet \rho_{\mathrm{DT}}
)$ 
and by the proof of Thm. \ref{thm:5:31:3103}, the inclusion 
$\mathrm{Stab}_{\mathcal G}(\mathrm{C}_3(\rho_1)^\times \bullet \rho_{\mathrm{DT}}
)\subset 
\mathrm{Stab}_{\mathcal G}(\mathbf k[[u,v]]^\times\bullet\Delta^{\mathcal W}_{r,l}
)$ 
is induced by the morphism 
$$
(A) : (\mathrm{C}_3(\rho_1)^\times\backslash
\mathrm{Hom}^1_{\mathcal C\operatorname{-alg}}(\hat{\mathcal V},M_3\hat V),\mathrm{C}_3(\rho_1)^\times \bullet \rho_{\mathrm{DT}}
,*)\to
(\mathbf k[[u,v]]^\times\backslash\mathrm{Hom}_{\mathcal C\operatorname{-alg}}(\hat{\mathcal W},\hat V),
\mathbf k[[u,v]]^\times\bullet\Delta^{\mathcal W}_{r,l}
,*) 
$$
(see Notation \ref{notation:F}) in $\mathbf{PS}_{\mathcal G}$.  
One relates the source and target of the morphism (A) by the following 
zig-zag of morphisms in $\mathbf{PS}_{\mathcal G}$ consisting of maps labeled (B)-(E) 
\begin{equation}\label{the:big:diagram}
   \xymatrix@C=20pt{
&(\mathrm C_{21}^{(0)}(\rho_1)^\times\backslash\ar_{(C)}[ld]\mathrm{Hom}^{1,((0)),\bullet}_{\mathcal C\operatorname{-alg}}(\hat{\mathcal V},T_{21}\hat V),
\mathrm C_{21}^{(0)}(\rho_1)^\times\bullet\rho_{\mathrm{DT}}
)
\ar^{(D)}[d]
\\
(\mathrm C_{21}(\rho_1)^\times\backslash\mathrm{Hom}^{1,(0)}_{\mathcal C\operatorname{-alg}}(\hat{\mathcal V},T_{21}\hat V),
\mathrm C_{21}(\rho_1)^\times\bullet \rho_{\mathrm{DT}}
)\ar_{(B)}[d]&
(\mathrm C_{2}(\overline\rho_1)^\times\backslash\mathrm{Hom}^{1,(0)}_{\mathcal C\operatorname{-alg}}(\hat{\mathcal V},M_2\hat V),
\mathrm C_{2}(\overline\rho_1)^\times\bullet \overline\rho_{\mathrm{DT}}
\ar^{(E)}[d]\\
(\mathrm C_3(\rho_1)^\times\backslash\mathrm{Hom}^{1}_{\mathcal C\operatorname{-alg}}(\hat{\mathcal V},M_3\hat V),
\mathrm C_3(\rho_1)^\times\bullet \rho_{\mathrm{DT}}
)
&(\mathbf k[[u,v]]^\times\backslash\mathrm{Hom}_{\mathcal C\operatorname{-alg}}(\hat{\mathcal W},\hat V),
\mathbf k[[u,v]]^\times\bullet\Delta^{\mathcal W}_{r,l}
)
}\tag{6.0.0}
\end{equation}
where $\mathrm C_3(\rho_1)^\times\bullet \rho_{\mathrm{DT}},
\mathbf k[[u,v]]^\times\bullet\Delta^{\mathcal W}_{r,l}
$
are as in \eqref{defs:[[rho]]}. 
One furthermore 
proves these morphisms to be locally injective, which together with Lem. \ref{lem:general:2212}(b) implies 
that $\mathrm{Stab}_{\mathcal G}(\mathrm C_3(\rho_1)^\times\bullet \rho_{\mathrm{DT}}
)\subset 
\mathrm{Stab}_{\mathcal G}(\mathbf k[[u,v]]^\times\bullet\Delta^{\mathcal W}_{r,l}
)$ is an equality (see Thm. \ref{thm:eq:Stab:Stab}).

Diagram \eqref{the:big:diagram}  
is established in \S \ref{sec 6}. The injectivity of 
the map (C) is obtained in \S \ref{sec 7}.  The local injectivities of the maps (B), (D) and (E) are then respectively 
proved in \S \S \ref{sec 8}, \ref{sec 9} and \ref{sec 10}. 
The equality between $\mathrm{Stab}_{\mathcal G}(\mathrm C_3(\rho_1)^\times\bullet \rho_{\mathrm{DT}}
)$ 
and $\mathrm{Stab}_{\mathcal G}(\mathbf k[[u,v]]^\times\bullet\Delta^{\mathcal W}_{r,l}
)$ is then derived in \S \ref{sect:11}.

\section{A diagram of pointed sets with action of $\mathcal G$} \label{sec 6}

The purpose of this section is to construct \eqref{the:big:diagram}, which is a diagram  of pointed sets with group 
actions, equipped with an overall action of $\mathcal G$. After preliminary material is treated in 
§\ref{sect:6:1:5768}, we construct the diagram of sets underlying \eqref{the:big:diagram} in 
§\ref{sect:6:2:5768}, then enhance it to a diagram of pointed sets in §\ref{sect:6:3:5786}. 
We treat an algebraic question in §\ref{sect:6:4:5786}, which leads in §\ref{sect:6:5:5786} to the construction 
of a collection of compatible group actions on that diagram. We construct the overall action of $\mathcal G$ 
on the resulting diagram in §\ref{sect:6:6:5786}, thus obtaining \eqref{the:big:diagram}.  
We summarize the situation in §\ref{sect:6:7:5786}.

\subsection{An algebra morphism $T_{21}\hat V\to M_2\hat V$}\label{sect:6:1:5768}

For $p,q\geq 1$ and $R$ a $\mathbf k$-module, define the subset $T_{p,q}R:=\{\begin{pmatrix}
        a&b\\0&c
    \end{pmatrix}|a\in M_pR,b\in M_{p,q}R,c\in M_qR\}\subset M_{p+q}R$.  

\begin{lem}\label{lem:basic:algebra:2412}
(a) For $p,q\geq 1$, the subset $T_{p,q}\hat V$ is a complete graded $\mathbf k$-subalgebra of 
    $M_{p+q}\hat V$ (see \S\ref{sect:4:2012}), with degree $n$ part given by $(T_{p,q}\hat V)_n:=T_{p,q}(\hat V_n)$.

(b) The map $T_{21}\hat V\to M_2\hat V$, $\begin{pmatrix}
        a&b\\0&c
    \end{pmatrix}\mapsto a$, which will be denoted $t\mapsto \overline t$, is an algebra morphism, which induces a group morphism 
    $(T_{21}\hat V)^\times\to \mathrm{GL}_2\hat V$. 

    (c) The map $a\mapsto \mathrm{diag}(a,1)$ is a group morphism $\mathrm{GL}_2\hat V\to(T_{21}\hat V)^\times$, which is a section of 
   the group morphism from (b).  
\end{lem}

\begin{proof}
    Immediate. 
\end{proof}

\subsection{A diagram of sets}\label{sect:6:2:5768}

For $r\in \hat V$, set 
\begin{equation}\label{def:R:r}
   R_r:=\begin{pmatrix}
    0&0&r
\end{pmatrix}\in M_{1,3}\hat V.  
\end{equation} 

\begin{defn}
We will write $\mathrm{C}_{21}(X)$ instead of $\mathrm C_{A}(X)$ if $X\subset A$ and $A=T_{21}\hat V$.  
\end{defn}

\begin{defn}\label{def:6:4:2912}
(a)  Define   $\mathrm{Hom}_{\mathcal C\operatorname{-alg}}(\hat{\mathcal V},T_{21}\hat V)$ as the set of algebra morphisms 
$\hat{\mathcal V}\to T_{21}\hat V$. Using the algebra inclusion $T_{21}\hat V\subset M_3\hat V$, it can be viewed as a subset of 
$\mathrm{Hom}_{\mathcal C\operatorname{-alg}}(\hat{\mathcal V},M_3\hat V)$. 

(b)  Define   $\mathrm{Hom}_{\mathcal C\operatorname{-alg}}^1(\hat{\mathcal V},T_{21}\hat V)$ intersection of 
$\mathrm{Hom}_{\mathcal C\operatorname{-alg}}(\hat{\mathcal V},T_{21}\hat V)$ and 
$\mathrm{Hom}_{\mathcal C\operatorname{-alg}}^1(\hat{\mathcal V},M_3\hat V)$ (see Def. \ref{defn: Hom1}). 

(c) Define  $\mathrm{Hom}^{(0)}_{\mathcal C\operatorname{-alg}}(\hat{\mathcal V},T_{21}\hat V)$ as the subset of $\mathrm{Hom}_{\mathcal C\operatorname{-alg}}(\hat{\mathcal V},T_{21}\hat V)$ 
of all morphisms $\rho : \hat{\mathcal V}\to T_{21}\hat V$
such that $\rho(e_0)$ is $(T_{21}\hat V)^\times$-conjugate to $\rho_0$ (i.e. for some $\alpha\in (T_{21}\hat V)^\times$, one has $\rho(e_0)
=\alpha \rho_0 \alpha^{-1}$). 

(d) Define  $\mathrm{Hom}^{((0))}_{\mathcal C\operatorname{-alg}}(\hat{\mathcal V},T_{21}\hat V)$ as the subset of
$\mathrm{Hom}^{(0)}_{\mathcal C\operatorname{-alg}}(\hat{\mathcal V},T_{21}\hat V)$ of all morphisms $\rho$
such that $\rho(e_0)\in\rho_0+T_{21}F^2\hat V$ (see Def. \ref{def:5:2:1926}). 

(e)  Define   $\mathrm{Hom}^\bullet_{\mathcal C\operatorname{-alg}}(\hat{\mathcal V},T_{21}\hat V)$ as the subset of $\mathrm{Hom}_{\mathcal C\operatorname{-alg}}(\hat{\mathcal V},T_{21}\hat V)$ 
 of all morphisms $\rho : \hat{\mathcal V}\to T_{21}\hat V$
such that for some $(r,C)\in \hat V\times M_{3,1}F^1\hat V$, one has
$$
\mathrm C_{21}(\rho(\hat{\mathcal V}))=\mathbf k1+C\cdot \mathrm C_{\hat V}(e_0)\cdot R_r\text{ and }R_r\cdot C\in e_0+f_\infty+F^2\hat V. 
$$

(f) For $S$ any subset of the set of indices $\{1,(0),((0)),\bullet\}$, set 
$$
\mathrm{Hom}^S_{\mathcal C\operatorname{-alg}}(\hat{\mathcal V},T_{21}\hat V)
:=\cap_{s\in S}\mathrm{Hom}^s_{\mathcal C\operatorname{-alg}}(\hat{\mathcal V},T_{21}\hat V).
$$ 
\end{defn}

\begin{lem}
The images of $\rho_0,\rho_1$ by the morphism $x\mapsto \overline x$ (see Lem. \ref{lem:basic:algebra:2412}(b)) are given by 
    $$
    \overline\rho_0=\begin{pmatrix}
        e_0&0 \\ e_1&f_0
    \end{pmatrix},\quad \overline\rho_1=\begin{pmatrix}
        1\\-1
    \end{pmatrix}\begin{pmatrix}
        e_1 & -f_1
    \end{pmatrix} 
    $$
so that $\overline\rho_1=\overline{\mathrm{col}}_{\mathrm{DT}}\cdot \overline{\mathrm{row}}_{\mathrm{DT}}$. 
\end{lem}

\begin{proof}
    Obvious.
\end{proof}

\begin{defn}
Define  $\mathrm{Hom}^{1,(0)}_{\mathcal C\operatorname{-alg}}(\hat{\mathcal V},M_2\hat V)$ 
as the subset of $\mathrm{Hom}_{\mathcal C\operatorname{-alg}}(\hat{\mathcal V},M_2\hat V)$ of all 
morphisms $\sigma : \hat{\mathcal V}\to M_2\hat V$ of filtered topological $\mathbf k$-algebras
such that $\sigma(e_1)=\overline\rho_1$ and $\sigma(e_0)$ is $\mathrm{GL}_2\hat V$-conjugate to $\overline\rho_0$. 
\end{defn}

\begin{lem}\label{lem:65:7:2412}
 (a)   Composition with the morphism $T_{21}\hat V\to M_2\hat V$, $x\mapsto \overline x$ induces a map $\mathrm{Hom}_{\mathcal C\operatorname{-alg}}(\hat{\mathcal V},T_{21}\hat V)
    \to\mathrm{Hom}_{\mathcal C\operatorname{-alg}}(\hat{\mathcal V},M_2\hat V)$, which induces a map 
    $\mathrm{Hom}^{1,(0)}_{\mathcal C\operatorname{-alg}}(\hat{\mathcal V},T_{21}\hat V)\to\mathrm{Hom}^{1,(0)}_{\mathcal C\operatorname{-alg}}(\hat{\mathcal V},M_2\hat V)$. 

    (b) For any $\sigma\in \mathrm{Hom}^{1,(0)}_{\mathcal C\operatorname{-alg}}(\hat{\mathcal V},M_2\hat V)$, there is an element 
    $\Delta_\sigma\in \mathrm{Hom}_{\mathcal C\operatorname{-alg}}(\hat{\mathcal W},\hat V)$, defined by $\Delta_\sigma(e_0^{n-1}e_1)
    =\overline{\mathrm{row}}_{\mathrm{DT}}\cdot \sigma(e_0)^{n-1}\cdot \overline{\mathrm{col}}_{\mathrm{DT}}$ for any $n\geq 1$; then 
    $\sigma\mapsto\Delta_\sigma$ is a map $\mathrm{Hom}^{1,(0)}_{\mathcal C\operatorname{-alg}}(\hat{\mathcal V},M_2\hat V)\to
    \mathrm{Hom}_{\mathcal C\operatorname{-alg}}(\hat{\mathcal W},\hat V)$. 
\end{lem}

\begin{proof}
(a) For $\rho \in \mathrm{Hom}_{\mathcal C\operatorname{-alg}}(\hat{\mathcal V},T_{21}\hat V)$, let $\overline\rho$ be the composition of $\rho$ with $x\mapsto \overline x$. Then 
$\overline\rho\in\mathrm{Hom}_{\mathcal C\operatorname{-alg}}(\hat{\mathcal V},M_2\hat V)$. Assume that $\rho\in \mathrm{Hom}^{1,(0)}_{\mathcal C\operatorname{-alg}}(\hat{\mathcal V},T_{21}\hat V)$. 
Then $\overline\rho(e_1)=(x\mapsto \overline x)(\rho(e_1))=\overline\rho_1$ since $\rho(e_1)=\rho_1$. Let $\alpha\in (T_{21}\hat V)^\times$ be such that 
$\rho(e_0)=\alpha\rho_0\alpha^{-1}$. Then $\overline\alpha\in\mathrm{GL}_2\hat V$, and $\overline\rho(e_0)=\overline\alpha\overline\rho_0\overline\alpha^{-1}$, therefore 
$\overline\rho\in \mathrm{Hom}^{1,(0)}_{\mathcal C\operatorname{-alg}}(\hat{\mathcal V},M_2\hat V)$. 

(b) follows from the fact that $\hat{\mathcal W}$ is freely generated by the family $(e_0^{n-1}e_1)_{n\geq1}$. 
\end{proof}

\begin{defn}\label{def:6:8:2412}
The map $\mathrm{Hom}^{1,((0)),\bullet}_{\mathcal C\operatorname{-alg}}(\hat{\mathcal V},T_{21}\hat V)\to
\mathrm{Hom}^{1,(0)}_{\mathcal C\operatorname{-alg}}(\hat{\mathcal V},M_2\hat V)$ is defined to be the composition 
$\mathrm{Hom}^{1,((0)),\bullet}_{\mathcal C\operatorname{-alg}}(\hat{\mathcal V},T_{21}\hat V)\hookrightarrow
\mathrm{Hom}^{1,(0)}_{\mathcal C\operatorname{-alg}}(\hat{\mathcal V},T_{21}\hat V)\stackrel{\mathrm{Lem.\ } \ref{lem:65:7:2412}(a)}{\to}
\mathrm{Hom}^{1,(0)}_{\mathcal C\operatorname{-alg}}(\hat{\mathcal V},M_2\hat V)$. 
\end{defn}
The main object of study of Part 3 is the following diagram of sets
\begin{equation}\label{DIAGSET}
   \xymatrix@C=20pt{
\mathrm{Hom}^{1,(0)}_{\mathcal C\operatorname{-alg}}(\hat{\mathcal V},T_{21}\hat V)\ar@{_{(}->}[d]
&\ar@{_{(}->}[l]\mathrm{Hom}^{1,((0)),\bullet}_{\mathcal C\operatorname{-alg}}(\hat{\mathcal V},T_{21}\hat V)\ar^{\mathrm{Def.\ } \ref{def:6:8:2412}}[r]
&\mathrm{Hom}^{1,(0)}_{\mathcal C\operatorname{-alg}}(\hat{\mathcal V},M_2\hat V)\ar_{\mathrm{Lem.\ } \ref{lem:65:7:2412}(b)}[d]
\\
\mathrm{Hom}^{1}_{\mathcal C\operatorname{-alg}}(\hat{\mathcal V},M_3\hat V)&
&
\mathrm{Hom}_{\mathcal C\operatorname{-alg}}(\hat{\mathcal W},\hat V) 
}
\end{equation}

\begin{lem}\label{lem:comm:square:2712}
    The following diagram of sets is commutative
$$
\xymatrix{\mathrm{Hom}^{1,(0)}_{\mathcal C\operatorname{-alg}}(\hat{\mathcal V},T_{21}\hat V)\ar^{\mathrm{Lem.} \ref{lem:65:7:2412}(a)}[r]\ar@{_{(}->}[d]&
 \mathrm{Hom}^{1,(0)}_{\mathcal C\operatorname{-alg}}(\hat{\mathcal V},M_2\hat V)\ar^{\mathrm{Lem.\ } \ref{lem:65:7:2412}(b)}[d]\\
 \mathrm{Hom}_{\mathcal C\operatorname{-alg}}^1(\hat{\mathcal V},M_3\hat V)\ar_{\mathrm{Lem.\ } \ref{lem:vania}(a)}[r]&
 \mathrm{Hom}_{\mathcal C\operatorname{-alg}}(\hat{\mathcal W},\hat V)}
$$
\end{lem}

\begin{proof}
Let $\rho\in\mathrm{Hom}^{1,(0)}_{\mathcal C\operatorname{-alg}}(\hat{\mathcal V},T_{21}\hat V)$. For some $\alpha\in M_{21}\hat V$ and $\beta\in \hat V$, one has 
$\rho(e_0)=\begin{pmatrix}
    \overline{\rho(e_0)}&\alpha \\0&\beta
\end{pmatrix}$. Moreover, one has $\mathrm{row}_{\mathrm{DT}}=\begin{pmatrix}
  \overline{\mathrm{row}}_{\mathrm{DT}}  & 0
\end{pmatrix}$ and $\mathrm{col}_{\mathrm{DT}}=\begin{pmatrix}
  \overline{\mathrm{col}}_{\mathrm{DT}}  \\ 0
\end{pmatrix}$. For any $n\geq 0$, one then has 
\begin{align*}
&\Delta_\rho(e_0^{n-1}e_1)=\mathrm{row}_{\mathrm{DT}}\cdot \rho(e_0)^{n-1}\cdot \mathrm{col}_{\mathrm{DT}}
=\begin{pmatrix}
  \overline{\mathrm{row}}_{\mathrm{DT}}  & 0
\end{pmatrix}\cdot \begin{pmatrix}
    \overline{\rho(e_0)}&\alpha \\0&\beta
\end{pmatrix}^{n-1}\cdot \begin{pmatrix}
  \overline{\mathrm{col}}_{\mathrm{DT}}  \\ 0
\end{pmatrix}
\\& =\overline{\mathrm{row}}_{\mathrm{DT}}\cdot \overline{\rho(e_0)}^{n-1}\cdot \overline{\mathrm{col}}_{\mathrm{DT}}=\Delta_{\overline\rho}(e_0^{n-1}e_1), 
\end{align*}
therefore $\Delta_\rho$
and $\Delta_{\overline\rho}$ are equal. The statement follows from the fact that these elements are the 
respective images of $\rho$ under the maps $\mathrm{Lem.\ } \ref{lem:vania}(a)$ and 
$\mathrm{Lem.\ } \ref{lem:65:7:2412}(b)\circ\mathrm{Lem.\ } \ref{lem:65:7:2412}(a)$. 
\end{proof}

\begin{rem}

 Both diagrams \eqref{DIAGSET} and \eqref{TAG':2112} from Part 2 are subdiagrams of the following diagram
 $$
 \xymatrix@C=40pt{
 &\mathrm{Hom}^{1,((0)),\bullet}_{\mathcal C\operatorname{-alg}}(\hat{\mathcal V},T_{21}\hat V)\ar@{_{(}->}[d]\ar^{\mathrm{Def.\ } \ref{def:6:8:2412}}[rd] & \\
 &\mathrm{Hom}^{1,(0)}_{\mathcal C\operatorname{-alg}}(\hat{\mathcal V},T_{21}\hat V)\ar^{\mathrm{Lem.} \ref{lem:65:7:2412}(a)}[r]\ar@{_{(}->}[d]&
 \mathrm{Hom}^{1,(0)}_{\mathcal C\operatorname{-alg}}(\hat{\mathcal V},M_2\hat V)\ar^{\mathrm{Lem.\ } \ref{lem:65:7:2412}(b)}[d]\\
 \mathrm{Hom}_{\mathcal C\operatorname{-alg}}(\hat{\mathcal V},M_3\hat V)&
 \ar@{_{(}->}[l]\mathrm{Hom}_{\mathcal C\operatorname{-alg}}^1(\hat{\mathcal V},M_3\hat V)\ar_{\mathrm{Lem.\ } \ref{lem:vania}(a)}[r]&
 \mathrm{Hom}_{\mathcal C\operatorname{-alg}}(\hat{\mathcal W},\hat V)
 }
 $$
 where the square is commutative by Lem. \ref{lem:comm:square:2712} and the triangle is commutative by Def. \ref{def:6:8:2412}. 
\end{rem}

\subsection{A diagram of pointed sets}\label{sect:6:3:5786}

\begin{lem}\label{lem:comm:e0}
    (a) One has $\mathrm{C}_V(e_0)=\mathbf k[e_0]\otimes\mathcal V$  (equality of subalgebras of $V$). 

    (b) The subalgebra $\mathrm{C}_{\hat V}(e_0)$ of $\hat V$ is the graded completion of $\mathrm{C}_V(e_0)$, 
    and is equal to $\mathbf k[[e_0]]\hat\otimes\hat{\mathcal V}$. 
\end{lem}

\begin{proof}
(a) For $a,b\geq0$, let us denote by $\mathcal V_{\mathbb Q,a,b}$ the part of $\mathcal V_{\mathbb Q}$ of bidegree $(a,b)$ for 
the bidegree defined by $(e_0,e_1)$. Then $\mathcal V_{\mathbb Q}=\mathbb Q[e_0]\oplus(\oplus_{a\geq 0,b>0}\mathcal V_{\mathbb Q,a,b})$. 
The kernel of the endomorphism $[e_0,-]$ of $\mathcal V_{\mathbb Q}$ is $\mathbb Q[e_0]$, which implies that 
it induces injective linear maps $\mathcal V_{\mathbb Q,a,b}\to\mathcal V_{\mathbb Q,a+1,b}$ for $a\geq 0$ and $b>0$. 
Set then $\mathcal V_{\mathbb Q}^+:=\oplus_{a\geq0,b>0}\mathcal V_{\mathbb Q,a,b}$ and $S:=
\oplus_{a,b\geq0}S_{a,b}$ where $S_{a,b}$ is a complementary vector subspace of the image of the linear map 
$\mathcal V_{\mathbb Q,a-1,b}\to\mathcal V_{\mathbb Q,a,b}$ induced by $[e_0,-]$. Then 
$\mathcal V_{\mathbb Q}=S\oplus\mathrm{im}([e_0,-])$, and $[e_0,-]$ induces an isomorphism
$\mathcal V_{\mathbb Q}^+\to \mathrm{im}([e_0,-])$. The endomorphism $[e_0,-]$ of $V$ is then the 
composition $V=(\mathbb Q[e_0]\otimes_{\mathbb Q}\mathcal V)\oplus (\mathcal V_{\mathbb Q}^+\otimes_{\mathbb Q}\mathcal V)
\to (S\otimes_{\mathbb Q}\mathcal V)\oplus(\mathrm{im}([e_0,-])\otimes_{\mathbb Q}\mathcal V)=V$ induced by 
$0 \oplus ([e_0,-])\otimes id_{\mathcal V}$, whose kernel is $\mathbb Q[e_0]\otimes_{\mathbb Q}\mathcal V=\mathbf k[e_0]\otimes\mathcal V$.

(b) follows from (a) and from the fact that $e_0$ is homogeneous. 
\end{proof}

\begin{lem}\label{lem:technical:0801}
(a) The endomorphism of $\hat V$ given by $x\mapsto x\cdot f_0-e_0\cdot x$ is injective; the same is true of the 
endomorphism $x\mapsto x\cdot e_0-f_0\cdot x$. 

(b) One has the equality $\mathrm{C}_{\hat V}(e_0)\cap \mathrm{C}_{\hat V}(f_0)\cap \mathbf k[[e_1,f_1]]
=\mathbf k1$ (equality of subspaces of $\hat V$). 
\end{lem}

\begin{proof}
(a) Let $\varphi$ be the endomorphism of $V=\mathcal V^{\otimes2}$ given by $x\mapsto x(1\otimes e_0)-(e_0\otimes 1)x$. 
Then $\varphi$ is graded, therefore $\mathrm{Ker}(\varphi)$ is the direct sum of its homogeneous components. Let $n\geq 0$
and $x\in\mathrm{Ker}(\varphi)$ be homogeneous of degree $n$. Decompose $x=\sum_{r=0}^n x_{r,n-r}$ according to the bidegree
induced by the degree of $\mathcal V$. If $x\neq0$, let $r$ be the smallest integer such that $x_{r,n-r}\neq0$. Then the bidegree
$(r+1,n-r)$ component of $0=x(1\otimes e_0)-(e_0\otimes 1)x$ is $-(e_0\otimes 1)x_{r,n-r}$, hence $-(e_0\otimes 1)x_{r,n-r}=0$, which 
by the integrity of $\mathcal V$ implies $x_{r,n-r}=0$, a contradiction. Hence $x=0$. Therefore $\varphi$ is graded. It follows that 
the graded completion of $\varphi$ is injective, which proves the first statement of (a). Its second statement follows from the 
fact that $x\mapsto x\cdot e_0-f_0\cdot x$ is the conjugation of $x\mapsto x\cdot f_0-e_0\cdot x$ by the permutation of factors if 
$\hat V=\hat{\mathcal V}^{\hat\otimes2}$. 

(b) The spaces in the said intersection are the graded completion of graded subspaces of $V$, namely 
$\mathrm{C}_V(e_0)=\mathbf k[e_0]\otimes\mathcal V$, $\mathrm{C}_V(f_0)=\mathcal V\otimes\mathbf k[e_0]$
(see Lem. \ref{lem:comm:e0}(a)),
and $\mathbf k[e_1,f_1]=\mathbf k[e_1]^{\otimes2}$. The intersection of the two first spaces in $\mathbf k[e_0,f_0]
=\mathbf k[e_0]^{\otimes2}$. 
and its intersection with the latter space is $(\mathbf k[e_0]\cap \mathbf k[e_1])^{\otimes2}=\mathbf k$. 
\end{proof}

\begin{defn}\label{defn:C:DT:R:DT}
    Set $R_{\mathrm{DT}}:=\begin{pmatrix}
        0&0&1
    \end{pmatrix}\in M_{1,3}\hat V$ and $C_{\mathrm{DT}}:=\begin{pmatrix}
        f_1\\e_1\\-e_0-f_\infty
    \end{pmatrix}\in M_{3,1}\hat V$.  
\end{defn}

\begin{lem}\label{lem:commutant:BIS}
(a) The map $\mathbf k\oplus \mathrm{C}_{\hat V}(e_0)\to M_3\hat V$ defined by 
$$
(\phi,a)\mapsto \phi\cdot I_3+C_{\mathrm{DT}}\cdot a\cdot R_{\mathrm{DT}} 
$$
defines a $\mathbf k$-module isomorphism 
$\mathbf k\oplus \mathrm{C}_{\hat V}(e_0)\to\mathrm{C}_3(\rho_{\mathrm{DT}}(\hat{\mathcal V}))$. In particular
\begin{equation}\label{eq:commutant:BIS}
\mathrm C_3(\rho_{\mathrm{DT}}(\hat{\mathcal V}))=\mathbf kI_3+C_{\mathrm{DT}}\cdot \mathrm C_{\hat V}(e_0)\cdot R_{\mathrm{DT}}. 
\end{equation}

(b) One has 
$$
\mathrm C_{21}(\rho_{\mathrm{DT}}(\hat{\mathcal V}))=\mathbf kI_3+C_{\mathrm{DT}}\cdot \mathrm C_{\hat V}(e_0)\cdot R_{\mathrm{DT}}. 
$$
\end{lem}

\begin{proof}
(a) One checks the relations
\begin{equation}\label{rels:RDT:CDT:rho0:rho1}
  R_{\mathrm{DT}}\cdot\rho_1=0,\quad \rho_1\cdot C_{\mathrm{DT}}=0,\quad    
  R_{\mathrm{DT}}\cdot\rho_0=e_0\cdot R_{\mathrm{DT}},\quad \rho_0\cdot C_{\mathrm{DT}}=C_{\mathrm{DT}}\cdot e_0. 
\end{equation}
Let $a\in \mathrm C_{\hat V}(e_0)$. Then 
$$
C_{\mathrm{DT}}\cdot a \cdot R_{\mathrm{DT}}\cdot \rho_1=0=\rho_1\cdot C_{\mathrm{DT}}\cdot a \cdot R_{\mathrm{DT}}, 
$$
where the equalities follow from \eqref{rels:RDT:CDT:rho0:rho1}(a) and  \eqref{rels:RDT:CDT:rho0:rho1}(b). Therefore 
$C_{\mathrm{DT}}\cdot a \cdot R_{\mathrm{DT}}\in \mathrm{C}(\rho_1)$. Moreover, 
$$
C_{\mathrm{DT}}\cdot a \cdot R_{\mathrm{DT}}\cdot \rho_0
=C_{\mathrm{DT}}\cdot a \cdot e_0\cdot R_{\mathrm{DT}}
=C_{\mathrm{DT}}\cdot e_0 \cdot a\cdot R_{\mathrm{DT}}
=\rho_0\cdot C_{\mathrm{DT}} \cdot a\cdot R_{\mathrm{DT}}, 
$$
where the equalities follow from \eqref{rels:RDT:CDT:rho0:rho1}(c), the relation $a\in \mathrm C_{\hat V}(e_0)$, 
and  \eqref{rels:RDT:CDT:rho0:rho1}(d). Therefore 
$C_{\mathrm{DT}}\cdot a \cdot R_{\mathrm{DT}}\in \mathrm{C}(\rho_0)$.

It follows that $C_{\mathrm{DT}}\cdot \mathrm C_{\hat V}(e_0) \cdot R_{\mathrm{DT}}$ is contained in the intersection 
$\mathrm{C}(\rho_1)\cap\mathrm{C}(\rho_0)$. Since 
\begin{equation}\label{identification:intersection}
    \mathrm{C}(\rho_1)\cap\mathrm{C}(\rho_0)=\mathrm{C}(\rho_{\mathrm{DT}}(\hat{\mathcal V})), 
\end{equation}
 one derives $C_{\mathrm{DT}}\cdot \mathrm C_{\hat V}(e_0) \cdot R_{\mathrm{DT}}\subset \mathrm{C}(\rho_{\mathrm{DT}}(\hat{\mathcal V}))$, 
therefore 
\begin{equation}\label{one:way:incl:2512}
\mathbf kI_3+C_{\mathrm{DT}}\cdot \mathrm C_{\hat V}(e_0) \cdot R_{\mathrm{DT}}\subset \mathrm{C}(\rho_{\mathrm{DT}}(\hat{\mathcal V})). 
\end{equation}
Let us prove the opposite inclusion. Let $A\in \mathrm{C}(\rho_{\mathrm{DT}}(\hat{\mathcal V}))$. By \eqref{identification:intersection}, one has 
$A\in \mathrm{C}(\rho_1)\cap\mathrm{C}(\rho_0)$. Since $A\in\mathrm{C}(\rho_1)$, and by Lem. \ref{lem:13:5:1108:BIS}(b), there exists $(\phi,m)\in \mathbf k[[u,v]]\times M_2\hat V$ 
such that 
\begin{equation}\label{relation:A:M:phi:m}
    A=M(\phi,m)
\end{equation}
(see \eqref{def:M:2012}); let $a,b,c,d\in\hat V$ are the elements such that $m=\begin{pmatrix}  a&b\\c&d
\end{pmatrix}$. The relation $A\in\mathrm{C}(\rho_0)$ then imposes 
\begin{equation}\label{condition:M:phi:m:rho0}
  M(\phi,m)\cdot \rho_0=\rho_0\cdot M(\phi,m).   
\end{equation}
The $(1,2)$ and $(3,2)$ entries of this relation
yield $f_1af_0=e_0f_1a$ and $cf_0=e_0c$; since left multiplication by $f_1$ is injective in $\hat V$, the former equation implies 
$af_0=e_0a$, so that both $a$ and $c$ belong to the kernel of the map $x\mapsto xf_0-e_0x$, which by Lem. \ref{lem:technical:0801}(a) is 0; hence 
$a=c=0$. It follows that $m=\begin{pmatrix}
    b\\d
\end{pmatrix}\begin{pmatrix}
    0 & 1 
\end{pmatrix}$, therefore 
\begin{equation}\label{interm:form:M}
 M(\phi,m)=\phi I_3+\begin{pmatrix}
    f_1b\\e_1b\\d
\end{pmatrix}\cdot \begin{pmatrix}
    0&0&1
\end{pmatrix}.    
\end{equation}
Plugging this in \eqref{condition:M:phi:m:rho0} yields
\begin{equation}\label{new:eq:comm:2512}
 \phi\begin{pmatrix}
    e_0&0&0\\e_1&f_0&-e_1\\0&0&e_0
\end{pmatrix}+\begin{pmatrix}
    f_1be_0\\e_1be_0\\de_0
\end{pmatrix}\begin{pmatrix}
    0&0&1
\end{pmatrix}
=\begin{pmatrix}
    e_0&0&0\\e_1&f_0&-e_1\\0&0&e_0
\end{pmatrix}\phi+\begin{pmatrix}
    e_0f_1b\\e_1(f_0+f_1)b-e_1d\\e_0d
\end{pmatrix}\begin{pmatrix}
    0&0&1
\end{pmatrix}.   
\end{equation}
The (1,1) and (2,2) entries of this equality imply the commutation of $\phi$ with $e_0$ and $f_0$, which 
together with $\phi\in\mathbf k[[e_1,f_1]]$ and Lem. \ref{lem:technical:0801}(b) implies 
\begin{equation}\label{phi:in:kk}
   \phi\in\mathbf k.  
\end{equation}
Plugging this in \eqref{new:eq:comm:2512} yields the relation 
$$
 \begin{pmatrix}
    f_1be_0\\e_1be_0\\de_0
\end{pmatrix}
=\begin{pmatrix}
    e_0f_1b\\e_1(f_0+f_1)b-e_1d\\e_0d
\end{pmatrix},     
$$
which using the facts that the left multiplications by $e_1$ and by $f_1$ in $\hat V$ are both injective, implies 
$b,d\in \mathrm{C}_{\hat V}(e_0)$ and $d=(f_0+f_1)b-be_0$. Using $b\in \mathrm{C}_{\hat V}(e_0)$, the latter equation implies
 $d=-(f_\infty+e_0)b$, which together with \eqref{interm:form:M}, \eqref{phi:in:kk}, \eqref{relation:A:M:phi:m} and Def. \ref{defn:C:DT:R:DT}, 
 implies $A\in \mathbf kI_3+C_{\mathrm{DT}}\cdot \mathrm C_{\hat V}(e_0)\cdot R_{\mathrm{DT}}$. This shows the opposite inclusion to 
 \eqref{one:way:incl:2512}, and therefore the equality
 $\mathbf kI_3+C_{\mathrm{DT}}\cdot \mathrm C_{\hat V}(e_0) \cdot R_{\mathrm{DT}}\subset \mathrm{C}(\rho_{\mathrm{DT}}(\hat{\mathcal V}))$. 
 It follows that the map $\mathbf k\times \mathrm{C}_{\hat V}(e_0)\to \mathrm{C}(\rho_{\mathrm{DT}}(\hat{\mathcal V}))$, 
  $(\phi,a)\mapsto \phi\cdot I_3+C_{\mathrm{DT}}\cdot a\cdot R_{\mathrm{DT}}$ is well-defined and surjective. 
 Its injectivity follows from the injectivity of the map $\mathbf k\times \mathrm{C}_{\hat V}(e_0)\to M_3\hat V$, 
  $(\phi,a)\mapsto \phi\cdot I_3+C_{\mathrm{DT}}\cdot a\cdot R_{\mathrm{DT}}$, which follows from  
 the fact that it is graded (with $\mathbf k$ of degree 0 and the degree of $\mathrm{C}_{\hat V}(e_0)$
being shifted by one) and that its graded components are injective (obvious in degree 0, follows from Lem. \ref{lem:TODO}(a)
in degree $>0$).

(b) follows from (a), from $\mathrm C_{21}(\rho_{\mathrm{DT}}(\hat{\mathcal V}))=\mathrm C_3(\rho_{\mathrm{DT}}(\hat{\mathcal V}))\cap T_{21}\hat V$, and from the fact that the 
right-hand side of \eqref{eq:commutant:BIS} is contained in $T_{21}\hat V$. 
\end{proof}

\begin{cor}\label{cor:6:15}
(a) The element  $\rho_{\mathrm{DT}}$ belongs to the set  $\mathrm{Hom}^{1,((0)),\bullet}_{\mathcal C\operatorname{-alg}}(\hat{\mathcal V},T_{21}\hat V)$.  

(b) By associating $\rho_{\mathrm{DT}}$ with the first three terms of the diagram \eqref{DIAGSET} (read from left to right) 
and $\overline\rho_{\mathrm{DT}}$ and $\Delta^{\mathcal W}_{r,l}$ with the two last terms, this diagram is upgraded to the diagram
$$
   \xymatrix@C=20pt{
(\mathrm{Hom}^{1,(0)}_{\mathcal C\operatorname{-alg}}(\hat{\mathcal V},T_{21}\hat V),\rho_{\mathrm{DT}})\ar@{_{(}->}[d]
&\ar@{_{(}->}[l]
(\mathrm{Hom}^{1,((0)),\bullet}_{\mathcal C\operatorname{-alg}}(\hat{\mathcal V},T_{21}\hat V),\rho_{\mathrm{DT}})
\ar[r]
&(\mathrm{Hom}^{1,(0)}_{\mathcal C\operatorname{-alg}}(\hat{\mathcal V},M_2\hat V),\overline\rho_{\mathrm{DT}})
\ar[d]
\\
(\mathrm{Hom}^{1}_{\mathcal C\operatorname{-alg}}(\hat{\mathcal V},M_3\hat V),\rho_{\mathrm{DT}})&
&
(\mathrm{Hom}_{\mathcal C\operatorname{-alg}}(\hat{\mathcal W},\hat V),\Delta^{\mathcal W}_{r,l}) 
}
$$
in the category $\mathbf{PS}$ 
of pointed sets. 
\end{cor}

\begin{proof}
(a) The relation $\rho_{\mathrm{DT}}\in \mathrm{Hom}^1_{\mathcal C\operatorname{-alg}}(\hat{\mathcal V},M_3\hat V)$ follows from $\rho_{\mathrm{DT}}(e_1)=\rho_1$. The relation
$\rho_{\mathrm{DT}}\in \mathrm{Hom}_{\mathcal C\operatorname{-alg}}
(\hat{\mathcal V},T_{21}\hat V)$ follows from $\rho_0,\rho_1\in T_{21}\hat V$, 
and the relation $\rho_{\mathrm{DT}}\in \mathrm{Hom}_{\mathcal C\operatorname{-alg}}^{((0))}
(\hat{\mathcal V},T_{21}\hat V)$ follows from $\rho_{\mathrm{DT}}(e_0)=\rho_0$. 
The relation $\rho_{\mathrm{DT}}\in \mathrm{Hom}^\bullet_{\mathcal C\operatorname{-alg}}(\hat{\mathcal V},T_{21}\hat V)$ follows from Lem. \ref{lem:commutant:BIS}(b)
and Def. \ref{def:6:4:2912}(e) (with $r=1$ and $C:=-C_{\mathrm{DT}}$). 
All this implies the statement. 

(b) If follows from (a) that $\rho_{\mathrm{DT}}$ belongs to the third  term of \eqref{DIAGSET}. Its image in the 
first two terms is then $\rho_{\mathrm{DT}}$, which can be viewed as an element of these terms. The image of 
$\rho_{\mathrm{DT}}$ in the fourth term is $\overline\rho_{\mathrm{DT}}$ by Lem. \ref{lem:65:7:2412}(a) and its image 
in the fifth term is $\Delta_{\overline\rho_{\mathrm{DT}}}$, which is equal to $\Delta_{\rho_{\mathrm{DT}}}$ by Lem. 
\ref{lem:comm:square:2712}, which is itself equal to $\Delta_{r,l}^{\mathcal W}$ by Lem. \ref{defn: Hom1}(b). 
\end{proof}

\subsection{Computation of the algebra $\mathrm C_2(\overline\rho_1)$ and of the group $\mathrm C_2(\overline\rho_1)^\times$}\label{sect:6:4:5786}

\begin{lem}\label{lem:TODO}
For $(\mathrm{col},\mathrm{row})\in M_{2,1}\hat V\times M_{1,2}\hat V$, the equality 
\begin{equation}\label{1652:1905:bis}
\mathrm{col}\cdot \overline{\mathrm{row}}_{\mathrm{DT}}
=\overline{\mathrm{col}}_{\mathrm{DT}}\cdot \mathrm{row}
\end{equation}
is equivalent to the existence of $a\in \hat V$ such that $\mathrm{col}=\overline{\mathrm{col}}_{\mathrm{DT}}\cdot a$ and 
$\mathrm{row}=a\cdot  \overline{\mathrm{row}}_{\mathrm{DT}}$. 

\end{lem}

\begin{proof}
Let $c_i,r_i\in\hat V$ ($i\in\{1,2\}$) be such that $\mathrm{col}=\begin{pmatrix}c_1\\c_2\end{pmatrix}$ 
and $\mathrm{row}=\begin{pmatrix}r_1&r_2\end{pmatrix}$. 
\eqref{1652:1905:bis} is equivalent to 
$$
c_1e_1=r_1,\quad -c_1f_1=r_2,\quad   
c_2e_1=-r_1,\quad -c_2f_1=-r_2, 
$$
i.e. to 
$$
r_1=c_1e_1,\quad r_2=-c_1f_1,\quad 
e_1\cdot (r_1+r_2)=f_1\cdot (c_1+c_2)e_1=(c_1+c_2)f_1=0. 
$$
By Lem. \ref{lem:TODO:BIS}(a), the last equation is equivalent to $c_1+c_2=0$. Set $a:=c_1$, one then obtains $(c_1,c_2)=(a,-a)$ and,
using the first equation, $(r_1,r_2)=(ae_1,-af_1)$, which implies 
$\mathrm{col}=\begin{pmatrix}1\\-1\end{pmatrix}\cdot a$ and 
$\mathrm{row}=a\cdot \begin{pmatrix}e_1&-f_1\end{pmatrix}$. This proves one of the implications; its converse is obvious. 
\end{proof}

\begin{lem}\label{lem:6:16:2912}
    (a) When equipped with the product $(\phi,v)\cdot (\phi',v'):=(\phi\phi',\phi(e_1,f_1)v'+v\phi(e_1,f_1)+v(e_1+f_1)v')$, the $\mathbf k$-module
    $\mathbf k[[u,v]]\oplus\hat V$ is a $\mathbf k$-algebra. 

    (b) The map $(\mathbf k[[u,v]]\oplus\hat V,\cdot)\to \mathrm{C}_2(\overline\rho_1)$ given by $(\phi,v)\mapsto \overline M(\phi,v):=\phi(e_1,f_1) I_2+\begin{pmatrix}
f_1\\ e_1 
\end{pmatrix}v\begin{pmatrix}
    1&1
\end{pmatrix}$ is an algebra isomorphism. 

(c) There is an algebra morphism $\mathrm{C}_2(\overline\rho_1)\to \mathbf k[[u,v]]$, $B\mapsto \phi_B$, where for each $B\in \mathrm{C}_2(\overline\rho_1)$, 
one denotes by $\phi_B$ the unique element in $\mathbf k[[u,v]]$ such that there exists $v\in \hat V$ with $B=\overline M(\phi_B,v)$. 
\end{lem}

\begin{proof}
(a) is a direct verification. 
For any $(\phi,v),(\phi',v')\in \mathbf k[[u,v]]\times\hat V$, one has 
\begin{align*}
&\overline M(\phi,v)\overline M(\phi',v')=(\phi(e_1,f_1) I_2+\begin{pmatrix}
f_1\\ e_1 
\end{pmatrix}v\begin{pmatrix}
    1&1
\end{pmatrix})\cdot (\phi'(e_1,f_1) I_2+\begin{pmatrix}
f_1\\ e_1 
\end{pmatrix}v'\begin{pmatrix}
    1&1
\end{pmatrix})
\\&=\phi(e_1,f_1)\phi'(e_1,f_1)I_2+\begin{pmatrix}
f_1\\ e_1 
\end{pmatrix}(v\phi'(e_1,f_1)+\phi(e_1,f_1)v'+v(e_1+f_1)v')\begin{pmatrix}
    1&1
\end{pmatrix}=\overline M((\phi,v)\cdot(\phi',v') , 
\end{align*}
where the second equality follows from the commutations of $\phi(e_1,f_1)$ with $\begin{pmatrix}
f_1\\ e_1 
\end{pmatrix}$ and of $\phi'(e_1,f_1)$ with $\begin{pmatrix}
    1&1
\end{pmatrix}$, as well as $\begin{pmatrix}
    1&1
\end{pmatrix}\begin{pmatrix}
f_1\\ e_1 
\end{pmatrix}=e_1+f_1$. 
Therefore the map $(\mathbf k[[u,v]]\times\hat V,\cdot)\to M_2\hat V$, $(\phi,v)\mapsto \overline M(\phi,v)$ is an algebra morphism. 

For any $(\phi,v)\in \mathbf k[[u,v]]\times\hat V$, one has 
\begin{align*}
 &   \overline M(\phi,v)\cdot\overline\rho_1=(\phi(e_1,f_1)I_2+\begin{pmatrix}
f_1\\ e_1 
\end{pmatrix}v\begin{pmatrix}
    1&1
\end{pmatrix})\cdot \overline{\mathrm{col}}_{\mathrm{DT}} \cdot \overline{\mathrm{row}}_{\mathrm{DT}} 
=\phi(e_1,f_1)  \cdot\overline{\mathrm{col}}_{\mathrm{DT}} \cdot \overline{\mathrm{row}}_{\mathrm{DT}}
\\& = \overline{\mathrm{col}}_{\mathrm{DT}} \cdot \overline{\mathrm{row}}_{\mathrm{DT}}\cdot\phi(e_1,f_1) 
=\overline{\mathrm{col}}_{\mathrm{DT}} \cdot \overline{\mathrm{row}}_{\mathrm{DT}}\cdot (\phi(e_1,f_1)I_2+\begin{pmatrix}
f_1\\ e_1 
\end{pmatrix}v\begin{pmatrix}
    1&1
\end{pmatrix})
=\overline\rho_1\cdot\overline M(\phi,v)
\end{align*}
where the second equality follows from $\begin{pmatrix}
    1&1
\end{pmatrix}\overline{\mathrm{col}}_{\mathrm{DT}}=0$, the third equality follows from the commutation of 
$\phi(e_1,f_1)$ with $\overline{\mathrm{row}}_{\mathrm{DT}}$ and $\overline{\mathrm{col}}_{\mathrm{DT}}$, and
the fourth equality follows from $\overline{\mathrm{row}}_{\mathrm{DT}}
\begin{pmatrix}
f_1\\ e_1 
\end{pmatrix}=0$. 
Therefore the image of the map $(\phi,v)\mapsto \overline M(\phi,v)$ is contained in 
$\mathrm C_2(\overline\rho_1)$. 

Let $B\in \mathrm C_2(\overline\rho_1)$. Then $B\overline\rho_1=\overline\rho_1B$ implies 
$(B\cdot\overline{\mathrm{col}}_{\mathrm{DT}})\cdot \overline{\mathrm{row}}_{\mathrm{DT}}
=\overline{\mathrm{col}}_{\mathrm{DT}}\cdot (\overline{\mathrm{row}}_{\mathrm{DT}}\cdot B)$, where
$B\cdot\overline{\mathrm{col}}_{\mathrm{DT}}\in M_{2,1}\hat V$ and $\overline{\mathrm{row}}_{\mathrm{DT}}\cdot B\in M_{1,2}\hat V$. 
By 
By Lem. \ref{lem:TODO}, this implies the existence of $\phi\in\hat V$ such that 
\begin{equation}\label{zangwill:1108:TER}
B\cdot \overline{\mathrm{col}}_{\mathrm{DT}}=\overline{\mathrm{col}}_{\mathrm{DT}}\cdot \phi
\quad\text{and}\quad
\overline{\mathrm{row}}_{\mathrm{DT}}\cdot B=\phi\cdot \overline{\mathrm{row}}_{\mathrm{DT}}. 
\end{equation}
Then 
$$
\phi\cdot (e_1+f_1)=\phi\cdot
\overline{\mathrm{row}}_{\mathrm{DT}}\cdot \overline{\mathrm{col}}_{\mathrm{DT}}
=\overline{\mathrm{row}}_{\mathrm{DT}}\cdot B\cdot \overline{\mathrm{col}}_{\mathrm{DT}}
=\overline{\mathrm{row}}_{\mathrm{DT}}\cdot \overline{\mathrm{col}}_{\mathrm{DT}}\cdot\phi
=(e_1+f_1)\cdot\phi, 
$$
where the two middle equalities follow from \eqref{zangwill:1108:TER}, 
therefore, by Lem. \ref{lemma:comm:0304:BIS}, $\phi\in\mathbf k[[e_1,f_1]]$.  
Then \eqref{zangwill:1108:TER}, together with the commutation of $\phi$ with the entries of 
$\overline{\mathrm{row}}_{\mathrm{DT}}$ and $\overline{\mathrm{col}}_{\mathrm{DT}}$, implies
the relation $B-\phi I_2\in \mathrm{Ann}(\overline{\mathrm{col}}_{\mathrm{DT}},\overline{\mathrm{row}}_{\mathrm{DT}})$, 
which by Lem. \ref{lem:van:0301:BIS} is equal to  $\{\begin{pmatrix}f_1\\e_1\end{pmatrix}
\cdot v\cdot \begin{pmatrix}1&1\end{pmatrix}|v\in\hat V\}$. Therefore there exists $v\in\hat V$ such that 
$B=\phi(e_1,f_1)I_2+\overline{\mathrm{col}}_{\mathrm{DT}}\cdot v\cdot \overline{\mathrm{row}}_{\mathrm{DT}}=\overline M(\phi,v)$. 
Therefore $B$ is contained in the image of $(\phi,v)\mapsto \overline M(\phi,v)$. 
It follows that the algebra morphism $\mathbf k[[u,v]]\oplus\hat V\to \mathrm C_2(\overline\rho_1)$, $(\phi,v)\mapsto \overline M(\phi,v)$
is surjective. 

The injectivity of the $\mathbf k$-module morphism $\mathbf k[[u,v]]\times\hat V\to M_2\hat V$, $(\phi,v)\mapsto \overline M(\phi,v)$ follows from that of the 
endomorphisms $v\mapsto e_1v$ and $v\mapsto f_1v$ of $\hat V$. All this implies (b). 

One checks that the map $(\mathbf k[[u,v]]\times\hat V,\cdot)\to\mathbf k[[u,v]]$, $(\phi,v)\mapsto \phi$ is an algebra morphism. The map describe in (c) is then 
the composition of this morphism with the inverse to the isomorphism from (a), and is therefore an algebra morphism.  
\end{proof}

\begin{lem}\label{lem:comp:C2:rho1:times}
The map $(\phi,v)\mapsto \overline M(\phi,v)$ from Lem. \ref{lem:6:16:2912}(c) induces a group isomorphism 
$(\mathbf k[[u,v]]^\times\times\hat V,\cdot)\to \mathrm{C}_2(\overline\rho_1)^\times$. 
\end{lem}

\begin{proof}
It follows from Lem. \ref{invertible:commutant} that $\mathrm{C}_2(\overline\rho_1)^\times=\mathrm{C}_2(\overline\rho_1)\cap\mathrm{GL}_2\hat V$. 
For $(\phi,v)\in \mathbf k[[u,v]]\times \hat V$, one has $(\overline M(\phi,v)\in \mathrm{GL}_2\hat V)\iff
(\epsilon(\overline M(\phi,v))\in\mathrm{GL}_2\mathbf k)\iff(\phi\in\mathbf k[[u,v]]^\times)$, where the second equivalence follows from 
the identity $\epsilon(\overline M(\phi,v))=\phi(0,0)I_2$, which together with Lem. \ref{lem:6:16:2912}(c) implies the result. 
\end{proof}

\subsection{Morphisms of sets with group actions}\label{sect:6:5:5786}

\subsubsection{Group morphisms}

\begin{lem}\label{lem:diag:algebras:2912}
   (a)  If $f : A\to B$ is an algebra morphism and $a\in A$, then $f$ induces an algebra morphism 
    $\mathrm C_A(a)\to \mathrm{C}_B(f(a))$. 

    (b) The diagram of algebras $M_3\hat V\supset T_{21}\hat V\to M_2\hat V$, where the second map is $x\mapsto \overline x$, induces a 
    diagram of algebras $\mathrm C_3(\rho_1)\supset \mathrm C_{21}(\rho_1)\to \mathrm C_2(\overline\rho_1)$. 
\end{lem}

\begin{proof}
  (a) is obvious, and (b) is a direct consequence.  
\end{proof}

\begin{lem}\label{lem:fibered:pdts}
(a) If $A\stackrel{f}{\to} C\supset B$ is a diagram of unital $\mathbf k$-algebras, then the fibered product $\mathbf k$-algebra $A\oplus_C B$ is  
equal to the unital $\mathbf k$-subalgebra $f^{-1}(B)$ of $A$. 

(b)  If $G\stackrel{\phi}{\to} K\supset H$ is a diagram of groups, then the fibered product group $G\times_KH$ is equal to the subgroup $\phi^{-1}(H)$ of $G$. 
    
(c) In the situation of (a), the following equality of groups holds 
$(A\oplus_C B)^\times=A^\times\times_{C^\times}B^\times$, the right-hand side being relative to the 
diagram of groups  $A^\times\stackrel{f^\times}{\to} C\supset B^\times$ induced by $A\stackrel{f}{\to} C\supset B$.   
\end{lem}

\begin{proof}
It follows from the fact that the fibered product corresponding to a diagram 
$X\stackrel{\phi}{\to}Z\stackrel{\psi}{\leftarrow}Y$ of groups (resp. unital $\mathbf k$-algebras), which  
is denoted $X\times_ZY$ (resp. $X\oplus_ZY$), is defined as $\{(x,y)\in X\times Y|\phi(x)=\psi(z)\}$. 
\end{proof}

It follows from the fact that $\rho_0\in T_{21}\hat V$ is homogeneous (of degree 1) that $\mathrm C_{21}(\rho_0)$ is a complete graded
$\mathbf k$-subalgebra of $T_{21}\hat V$, i.e. $\mathrm C_{21}(\rho_0)=\hat\oplus_{n\geq0}\mathrm C_{21}(\rho_0)_n$, where  
$\mathrm C_{21}(\rho_0)_n=\mathrm C_{21}(\rho_0)\cap (T_{21}\hat V)_n$ for $n\geq0$ (see Lem. \ref{lem:basic:algebra:2412}(a)), so that  
$\mathrm C_{21}(\rho_0)_0$ is a subalgebra of $(T_{21}\hat V)_0$. 
\begin{defn}\label{def:C0}
Define
$$
\mathrm C^{(0)}_{21}(\rho_1):=\mathrm{C}_{21}(\rho_1)\oplus_{(T_{21}\hat V)_0}\mathrm C_{21}(\rho_0)_0. 
$$
to be the unital $\mathbf k$-subalgebra of $\mathrm{C}_{21}(\rho_1)$ obtained by applying the construction of Lem. \ref{lem:fibered:pdts}(a) to 
the diagram $\mathrm{C}_{21}(\rho_1)\to (T_{21}\hat V)_0\supset\mathrm C_{21}(\rho_0)_0$, where
$\mathrm{C}_{21}(\rho_1)\to (T_{21}\hat V)_0$ is the restriction to $\mathrm{C}_{21}(\rho_1)$
of the projection $T_{21}\hat V=\hat\oplus_{n\geq0}(T_{21}\hat V)_n\to(T_{21}\hat V)_0$.     
\end{defn}

\begin{lem}
The following equality of groups holds
$$
\mathrm C^{(0)}(\rho_1)^\times=\mathrm{C}_{21}(\rho_1)^\times\times_{(T_{21}\hat V)_0^\times}\mathrm C_{21}(\rho_0)_0^\times,  
$$
the right-hand side being relative to the diagram of groups $\mathrm{C}_{21}(\rho_1)^\times\to (T_{21}\hat V)_0^\times\supset\mathrm C_{21}(\rho_0)_0^\times$. 
\end{lem}

\begin{proof}
    This follows from Lem. \ref{lem:fibered:pdts}(c). 
\end{proof}
The combination of the diagram of algebras from Lem. \ref{lem:diag:algebras:2912}(b) with the inclusions $\mathrm C^{(0)}_{21}(\rho_1)\subset 
\mathrm C_{21}(\rho_1)$ (see Def. \ref{def:C0}) and $\mathrm C_3(\rho_1)\subset M_3\hat V$ and 
with the algebra morphism from Lem. \ref{lem:6:16:2912}(c) gives rise to the following diagram of algebras
$$
M_3\hat V\supset\mathrm C_3(\rho_1)\supset \mathrm C_{21}(\rho_1)\supset\mathrm C_{21}^{(0)}(\rho_1)\to \mathrm C_{2}(\overline\rho_1)\to\mathbf k[[u,v]]
$$
which upon taking groups of units (invertible elements) gives rise to the following diagram of groups
\begin{equation}\label{diagram:gp:morphisms}
\mathrm{GL}_3\hat V\supset\mathrm C_3(\rho_1)^\times\supset \mathrm C_{21}(\rho_1)^\times\supset\mathrm C_{21}^{(0)}(\rho_1)^\times\to 
\mathrm C_{2}(\overline\rho_1)^\times\to\mathbf k[[u,v]]^\times. 
\end{equation}

\subsubsection{Actions of groups on sets}

\begin{lem}\label{lem:actions:2912}
    (a) The action of the group $(T_{21}\hat V)^\times$ on $\mathrm{Hom}_{\mathcal C\operatorname{-alg}}(\hat{\mathcal V},T_{21}\hat V)$ given by 
    $(g,\rho)\mapsto g\bullet\rho:=\mathrm{Ad}_g\circ\rho$ restricts to an action of 
    the subgroup $\mathrm C_{21}(\rho_1)^\times$ on the subset $\mathrm{Hom}_{\mathcal C\operatorname{-alg}}^{1,(0)}(\hat{\mathcal V},T_{21}\hat V)$. 

    (b) The action of $\mathrm C_{21}(\rho_1)^\times$ from (a) induces an action of $\mathrm C_{21}^{(0)}(\rho_1)^\times$ on  
    $\mathrm{Hom}_{\mathcal C\operatorname{-alg}}^{1,((0))}(\hat{\mathcal V},T_{21}\hat V)$. 

    (c) The action of $\mathrm C_{21}^{(0)}(\rho_1)^\times$ from (b) restricts to an action on  
    $\mathrm{Hom}_{\mathcal C\operatorname{-alg}}^{1,((0)),\bullet}(\hat{\mathcal V},T_{21}\hat V)$. 
    
  (d) The action of the group $\mathrm{GL}_2\hat V$ on $\mathrm{Hom}_{\mathcal C\operatorname{-alg}}(\hat{\mathcal V},M_2\hat V)$  given by 
    $(g,\rho)\mapsto g\bullet\rho:=\mathrm{Ad}_g\circ\rho$ restricts to an action of 
  $\mathrm C_{2}(\overline\rho_1)^\times$ on $\mathrm{Hom}_{\mathcal C\operatorname{-alg}}^{1,(0)}(\hat{\mathcal V},M_2\hat V)$. 
\end{lem}

\begin{proof}
(a) If $g\in \mathrm C_{21}(\rho_1)^\times$ and $\rho\in \mathrm{Hom}_{\mathcal C\operatorname{-alg}}^{1,(0)}(\hat{\mathcal V},T_{21}\hat V)$, then 
$\mathrm{Ad}_g\circ\rho\in \mathrm{Hom}_{\mathcal C\operatorname{-alg}}(\hat{\mathcal V},T_{21}\hat V)$ is such that 
$\mathrm{Ad}_g\circ\rho(e_1)=\mathrm{Ad}_g(\rho_1)=\rho_1$; moreover, $\mathrm{Ad}_g\circ\rho(e_0)$ is $(T_{21}\hat V)^\times$-conjugate to 
$\rho(e_0)$, which is itself $(T_{21}\hat V)^\times$-conjugate to $\rho_0$, therefore $\mathrm{Ad}_g\circ\rho(e_0)$ is $(T_{21}\hat V)^\times$-conjugate to 
$\rho_0$; therefore $\mathrm{Ad}_g\circ\rho\in \mathrm{Hom}_{\mathcal C\operatorname{-alg}}^{1,(0)}(\hat{\mathcal V},T_{21}\hat V)$. 

(b) If $g\in \mathrm C_{21}^{(0)}(\rho_1)^\times$ and $\rho\in \mathrm{Hom}_{\mathcal C\operatorname{-alg}}^{1,((0))}(\hat{\mathcal V},T_{21}\hat V)$, then 
$\mathrm{Ad}_g\circ\rho\in \mathrm{Hom}_{\mathcal C\operatorname{-alg}}^{1,(0)}(\hat{\mathcal V},T_{21}\hat V)$ by (a); moreover, 
$\mathrm{Ad}_g\circ\rho(e_0)\equiv\mathrm{Ad}_g(\rho_0)\equiv \mathrm{Ad}_{g_0}(\rho_0)=\rho_0$ mod $T_{21}F^1\hat V$, where
the first relation follows from $\rho(e_0)\equiv\rho_0$ mod $T_{21}F^1\hat V$, where $g_0\in (T_{21}\hat V)_0^\times$ is the degree 0 component 
of $g\in (T_{21}\hat V)^\times$, and where the last equality follows from $g_0\in \mathrm C_{21}(\rho_0)^\times$. This implies
$\mathrm{Ad}_g\circ\rho\in \mathrm{Hom}_{\mathcal C\operatorname{-alg}}^{1,((0))}(\hat{\mathcal V},T_{21}\hat V)$. 

(c) Since $\mathrm{Hom}_{\mathcal C\operatorname{-alg}}^{1,((0)),\bullet}(\hat{\mathcal V},T_{21}\hat V)$ is the intersection of 
$\mathrm{Hom}_{\mathcal C\operatorname{-alg}}^{1,((0))}(\hat{\mathcal V},T_{21}\hat V)$ with 
$\mathrm{Hom}_{\mathcal C\operatorname{-alg}}^{\bullet}(\hat{\mathcal V},T_{21}\hat V)$, and in view of (b), it suffices to 
prove that the action of  $(T_{21}\hat V)^\times$ on $\mathrm{Hom}_{\mathcal C\operatorname{-alg}}(\hat{\mathcal V},T_{21}\hat V)$
leaves the subset  $\mathrm{Hom}_{\mathcal C\operatorname{-alg}}^{\bullet}(\hat{\mathcal V},T_{21}\hat V)$ stable. 
Let then $T\in (T_{21}\hat V)^\times$ and $\rho\in \mathrm{Hom}_{\mathcal C\operatorname{-alg}}^{\bullet}(\hat{\mathcal V},T_{21}\hat V)$.
Let $(r,C)\in \hat V\times M_{3,1}F^1\hat V$ be such that 
$\mathrm C_{21}(\rho(\hat{\mathcal V}))=\mathbf k1+C\cdot \mathrm C_{\hat V}(e_0)\cdot R_r$ and $R_r\cdot C\in e_0+f_\infty+F_2\hat V$, and 
$a\in \mathrm{GL}_2\hat V$, $b\in M_{2,1}\hat V$ and $c\in\hat V^\times$ be such that $T=\begin{pmatrix}
    a & b\\0 & c
\end{pmatrix}$. 
Then $R_rg^{-1}=R_{rc^{-1}}$, therefore $(rc^{-1},TC) \in\hat V\times M_{3,1}F^1\hat V$ is such that $R_{rc^{-1}}\cdot TC=R_c T^{-1}TC=R_rC=e_0+f_\infty
+F_2\hat V$ and 
\begin{align*}
    & \mathrm C_{21}(\mathrm{Ad}_T\circ \rho(\hat{\mathcal V}))=\mathrm{Ad}_T(\mathrm C_{21}(\rho(\hat{\mathcal V})))
=T\cdot (\mathbf k1+C\cdot \mathrm C_{\hat V}(e_0)\cdot R_r)\cdot T^{-1}
=\mathbf k1+(TC)\cdot\mathrm C_{\hat V}(e_0) \cdot(R_rT^{-1})\\ & 
=\mathbf k1+(TC)\cdot\mathrm C_{\hat V}(e_0) \cdot R_{rc^{-1}}, 
\end{align*}
which implies $\mathrm{Ad}_T\circ \rho\in \mathrm{Hom}_{\mathcal C\operatorname{-alg}}^{\bullet}(\hat{\mathcal V},T_{21}\hat V)$. 

(d) The proof is similar to that of (a), replacing $T_{21}\hat V$, $\rho_1$, $\rho_0$  by $M_2\hat V$, $\overline\rho_1$, $\overline\rho_0$. 
\end{proof}

\subsubsection{Compatibility of group and set morphisms}

\begin{lem}  \label{lem:6:20}
    (a) The natural group and set injections induce a morphism of sets with group actions 
    $(\mathrm C_{21}(\rho_1)^\times,\mathrm{Hom}^{1,(0)}_{\mathcal C\operatorname{-alg}}(\hat{\mathcal V},T_{21}\hat V))\to (\mathrm C_3(\rho_1)^\times,\mathrm{Hom}^1_{\mathcal C\operatorname{-alg}}(\hat{\mathcal V},M_3\hat V))$, the actions being as 
    in Lem. \ref{lem:actions:2912}(a) and Lem. \ref{lem:520:2212:FIRST}(b). 

(b) The natural group and set morphisms induced by the algebra morphism $T_{21}\hat V\to M_2\hat V$, $x\mapsto \overline x$ induce a 
morphism of sets with group actions 
$$
(\mathrm C_{21}^{(0)}(\rho_1)^\times,\mathrm{Hom}^{1,((0)),\bullet}_{\mathcal C\operatorname{-alg}}(\hat{\mathcal V},T_{21}\hat V))\to
(\mathrm C_2(\overline\rho_1)^\times,\mathrm{Hom}^{1,(0)}_{\mathcal C\operatorname{-alg}}(\hat{\mathcal V},M_2\hat V));
$$
the actions being as in Lem. \ref{lem:actions:2912}(c) and (d). 
    
(c) The pair formed by the group morphism $\mathrm C_2(\overline\rho_1)^\times\to\mathbf k[[u,v]]^\times$ from Lem. \ref{lem:6:16:2912}(c) 
and by the set morphism $\sigma\mapsto\Delta_\sigma$ (see Lem. \ref{lem:65:7:2412}(b)) induces a morphism of sets with group actions 
$(\mathrm C_2(\overline\rho_1)^\times,\mathrm{Hom}^{1,(0)}_{\mathcal C\operatorname{-alg}}(\hat{\mathcal V},M_2\hat V))\to
    (\mathbf k[[u,v]]^\times,\mathrm{Hom}_{\mathcal C\operatorname{-alg}}(\hat{\mathcal W},\hat V))$, 
the actions being as in Lem. \ref{lem:actions:2912}(d) and Lem. \ref{lem28:1001}(c). 
\end{lem}

\begin{proof}
(a) follows from the combination of the fact that the natural group and set injections induces a morphism of sets with group actions
$$
((T_{21}\hat V)^\times,\mathrm{Hom}_{\mathcal C\operatorname{-alg}}(\hat{\mathcal V},T_{21}\hat V))\to 
(\mathrm{GL}_3\hat V,\mathrm{Hom}_{\mathcal C\operatorname{-alg}}(\hat{\mathcal V},M_3\hat V)),
$$
that the group (resp. set) injection maps the subgroup (resp. subset) 
$\mathrm C_{21}(\rho_1)^\times$ to $\mathrm C_3(\rho_1)^\times$ 
(resp. $\mathrm{Hom}^{1,(0)}_{\mathcal C\operatorname{-alg}}(\hat{\mathcal V},T_{21}\hat V)$ to 
$\mathrm{Hom}^1_{\mathcal C\operatorname{-alg}}(\hat{\mathcal V},M_3\hat V)$), and that these subsets are preserved by the actions of these subgroups. 

(b) follows from the fact that the group and set morphisms attached to the algebra morphism $T_{21}\hat V\to M_2\hat V$, $x\mapsto \overline x$
induce a morphism of sets with group actions 
$$
((T_{21}\hat V)^\times,\mathrm{Hom}_{\mathcal C\operatorname{-alg}}(\hat{\mathcal V},T_{21}\hat V))\to
(\mathrm{GL}_2\hat V,\mathrm{Hom}_{\mathcal C\operatorname{-alg}}(\hat{\mathcal V},M_2\hat V)),
$$
that this group (resp. set) morphism maps the 
subgroup $\mathrm C_{21}^{(0)}(\rho_1)^\times$ to the subgroup $\mathrm C_2(\overline\rho_1)^\times$
(resp. the subset $\mathrm{Hom}^{1,((0)),\bullet}_{\mathcal C\operatorname{-alg}}(\hat{\mathcal V},T_{21}\hat V))$ to the subset 
$\mathrm{Hom}^{1,(0)}_{\mathcal C\operatorname{-alg}}(\hat{\mathcal V},M_2\hat V)$), 
and that these subsets are preserved by the actions of these subgroups. 

Let us prove (c). 
Let $\sigma\in \mathrm{Hom}_{\mathcal C\operatorname{-alg}}^1(\hat{\mathcal V},M_3\hat V)$ and $A\in \mathrm C_2(\overline\rho_1)^\times$. 
By Lem. \ref{lem:6:16:2912}, there exists $(\phi,v)\in \mathbf k[[u,v]]^\times\times\hat V$ such that $Q=\overline M(\phi,v)$, and 
the image of $Q$ under the morphism $\mathrm C_2(\rho_1)^\times\to\mathbf k[[u,v]]^\times$ is $\phi(e_1,f_1)$. Then 
\begin{equation}\label{eq:col:2212:BIS}
Q\cdot \overline{\mathrm{col}}_{\mathrm{DT}}=(\phi(e_1,f_1)I_2+\begin{pmatrix}f_1\\e_1\end{pmatrix}\cdot 
v\cdot \begin{pmatrix}1&1\end{pmatrix})\cdot\overline{\mathrm{col}}_{\mathrm{DT}}
=\overline{\mathrm{col}}_{\mathrm{DT}}\cdot \phi(e_1,f_1),  
\end{equation}
and 
\begin{equation}\label{eq:row:2212:BIS}
\overline{\mathrm{row}}_{\mathrm{DT}}\cdot Q=\overline{\mathrm{row}}_{\mathrm{DT}}\cdot 
(\phi(e_1,f_1)I_2+\begin{pmatrix}f_1\\e_1\end{pmatrix}\cdot m\cdot \begin{pmatrix}1&1\end{pmatrix})
=\phi(e_1,f_1)\cdot \overline{\mathrm{row}}_{\mathrm{DT}}, 
\end{equation}
where the first equalities follow from $P=\overline M(\phi,m)$ and the second equalities follow from the 
commutation of the entries of $\overline{\mathrm{col}}_{\mathrm{DT}},\overline{\mathrm{row}}_{\mathrm{DT}}$
with $e_1,f_1$, and from the equalities $\begin{pmatrix}1&1\end{pmatrix}
\cdot \overline{\mathrm{col}}_{\mathrm{DT}}=0$ and 
$\overline{\mathrm{row}}_{\mathrm{DT}}\cdot \begin{pmatrix}f_1\\e_1\end{pmatrix}=0$. 

For any $n\geq 1$, one then has 
\begin{align*}
& \Delta_{Q\bullet\sigma}(e_0^{n-1}e_1)
=\Delta_{\mathrm{Ad}_Q \circ \sigma}(e_0^{n-1}e_1)
=\overline{\mathrm{row}}_{\mathrm{DT}}\cdot ((\mathrm{Ad}_Q\circ \sigma)(e_0))^{n-1}\cdot \overline{\mathrm{col}}_{\mathrm{DT}}
\\ & =\overline{\mathrm{row}}_{\mathrm{DT}}\cdot (Q\cdot \sigma(e_0)\cdot Q^{-1})^{n-1}\cdot \overline{\mathrm{col}}_{\mathrm{DT}}
=\overline{\mathrm{row}}_{\mathrm{DT}}\cdot Q\cdot \sigma(e_0)^{n-1}\cdot Q^{-1}\cdot \overline{\mathrm{col}}_{\mathrm{DT}}
\\& =\phi(e_1,f_1)\cdot \overline{\mathrm{row}}_{\mathrm{DT}}\cdot \sigma(e_0)^{n-1}\cdot  \overline{\mathrm{col}}_{\mathrm{DT}}
\cdot \phi(e_1,f_1)^{-1}
\\& =\phi(e_1,f_1)\cdot \Delta_\sigma(e_0^{n-1}e_1)
\cdot \phi(e_1,f_1)^{-1}
=(\mathrm{Ad}_{\phi(e_1,f_1)}\circ \Delta_\sigma)(e_0^{n-1}e_1)
=(\phi\bullet\Delta_\rho)(e_0^{n-1}e_1). 
\end{align*}
where all the equalities follow from definitions, except for the fifth one, 
which follows from \eqref{eq:col:2212:BIS} and 
\eqref{eq:row:2212:BIS}. 
This implies 
$$
 \forall \sigma\in\mathrm{Hom}_{\mathcal C\operatorname{-alg}}^1(\hat{\mathcal V},M_2\hat V),
\forall Q\in\mathrm C_2(\overline\rho_1)^\times,\quad \Delta_{Q\bullet\sigma}=\phi\bullet\Delta_\sigma. 
$$
(equality in $\mathrm{Hom}_{\mathcal C\operatorname{-alg}}^1(\hat{\mathcal W},\hat V)$), which implies the claim. 
\end{proof}

\begin{lem}
There is a diagram of pointed sets with group actions
\begin{equation}\label{diag:set:gp:acts}
   \xymatrix@C=20pt{
\substack{(\mathrm{Hom}^{1,(0)}_{\mathcal C\operatorname{-alg}}(\hat{\mathcal V},T_{21}\hat V),\\ \mathrm C_{21}(\rho_1)^\times,
\rho_{\mathrm{DT}},\bullet)}
\ar@{_{(}->}[d]
&\ar@{_{(}->}[l]\substack{(\mathrm{Hom}^{1,((0)),\bullet}_{\mathcal C\operatorname{-alg}}
(\hat{\mathcal V},T_{21}\hat V),\\ \mathrm C_{21}^{(0)}(\rho_1)^\times,\rho_{\mathrm{DT}},\bullet)}
\ar[r]
&
\substack{(\mathrm{Hom}^{1,(0)}_{\mathcal C\operatorname{-alg}}(\hat{\mathcal V},M_2\hat V),
\\
\mathrm{C}_2(\overline\rho_1)^\times,\overline\rho_{\mathrm{DT}},\bullet)}
\ar[d]
\\
\substack{(\mathrm{Hom}^{1}_{\mathcal C\operatorname{-alg}}(\hat{\mathcal V},M_3\hat V),\\
\mathrm{C}_3(\rho_1)^\times,\rho_{\mathrm{DT}},\bullet)}&
&
\substack{(\mathrm{Hom}_{\mathcal C\operatorname{-alg}}(\hat{\mathcal W},\hat V),\\
\mathbf k[[u,v]]^\times,\Delta^{\mathcal W}_{r,l},\bullet)} 
}
\end{equation}
where the notation $(X,A,x)\hookrightarrow(Y,B,y)$ for a morphism $(X,A,x)\to(Y,B,y)$ means that 
both the set morphism $X\to Y$ and the group morphism $A\to B$ are injective. 
\end{lem}

\begin{proof}
This follows from the combination of Lem. \ref{lem:6:20} and Lem. \ref{cor:6:15}(b)
\end{proof}


\subsection{Overall action of $\mathcal G$}\label{sect:6:6:5786}

\subsubsection{Action of $\mathcal G$ on sets}

\begin{lem}\label{lem:621:jan2025}
$g*\sigma:=\mathrm{aut}_g^V\circ \sigma\circ(\mathrm{aut}_g^{\mathcal V})^{-1}$
defines an action of $\mathcal G$ on the set $\mathrm{Hom}_{\mathcal C\operatorname{-alg}}(\hat{\mathcal V},M_2\hat V)$. 
\end{lem}

\begin{proof}
This follows from the fact that $g\mapsto \mathrm{aut}_g^{\mathcal V}$, $g\mapsto \mathrm{aut}_g^V$ define actions of
$\mathcal G$ on the algebras $\hat{\mathcal V}$ and $M_2\hat V$. 
\end{proof}

\begin{lem}\label{lem:1920:0307}
For $g\in\mathcal G$, there exists a unique family $(\varphi_k)_{k\geq 0}$ of elements of 
$\hat{\mathcal V}=\mathbf k\langle\langle e_0,e_1\rangle\rangle$ such that 
\begin{equation}\label{def:varphi:k:ja:2024}
e_1 g(e_0,e_1)=\sum_{k\geq0}\varphi_k(e_0,e_1)e_1e_0^k     
\end{equation}
(equality in $\hat{\mathcal V}$). Define elements of $\hat V$ by  
\begin{equation}\label{def:alpha:g:beta:g:3007}
\alpha_g:=\sum_{k\geq0}\varphi_k(e_0,e_1)f_0^k,\quad 
\gamma_g:=-\sum_{k\geq 0}\varphi_k(e_0,e_1)e_1\frac{e_0^k-f_0^k}{e_0-f_0}     
\end{equation}
and
$$
M_g:=\begin{pmatrix}
  g(e_0,e_1)g(f_0,f_1)  &0&0\\
  -g(f_0,f_1)\gamma_g &g(f_0,f_1)\alpha_g&g(f_0,f_1)\gamma_g\\
  0&0&g(e_0,e_1)g(f_0,f_1)
\end{pmatrix}, \quad
\overline M_g:=\begin{pmatrix}
   g(e_0,e_1)g(f_0,f_1)  &0\\
  -g(f_0,f_1)\gamma_g &g(f_0,f_1)\alpha_g
\end{pmatrix}, 
$$
then $M_g\in (T_{21}\hat V)^\times$, $\overline M_g\in\mathrm{GL}_2\hat V$ and 
\begin{equation}\label{act:g:rho0:barrho0}
\mathrm{aut}_g^V(\rho_0)=M_g\cdot\rho_0
\cdot M_g^{-1} ,\quad 
\mathrm{aut}_g^V(\overline\rho_0)=\overline M_g\cdot \overline\rho_0
\cdot \overline M_g^{-1} 
\end{equation}
(equalities in $T_{21}\hat V$ and $M_2\hat V$). 
\end{lem}

\begin{proof}
The map $\oplus_{k\geq 0}\mathcal V\to \mathcal V$, $(\varphi_k)_{k\geq 0}\mapsto \sum_{k\geq 0}\varphi_k\cdot e_1e_0^k$ defines a 
bijection from the source to the part of $\mathcal V$ of positive $e_1$-degree. It follows that the same formula defines a bijection 
between $\prod_{k\geq 0}\hat{\mathcal V}$ and the part of $\hat{\mathcal V}$ of positive $e_1$-degree. The first statement then 
follows from the fact that $e_1 g(e_0,e_1)$ has positive $e_1$-degree. 

One has $M_g\in T_{21}\hat V,\overline M_g\in M_2\hat V$. The part of total degree 0 (in $\hat{\mathcal V}$) 
of $g(e_0,e_1)$ is equal to $1$, which implies that the part of total degree 0 (in $\hat{\mathcal V}$) of $\varphi_0(e_0,e_1)$ is 1. 
This implies that the part of total degree 0 (in $\hat V$) of $\alpha_g$ is 1. It follows that the images of 
$M_g$ in $\overline M_g$ in $T_{21}\mathbf k$ and $M_2\mathbf k$ are $I_3$ and $I_2$, therefore $M_g\in(T_{21}\hat V)^\times$, 
$\overline M_g\in\mathrm{GL}_2\hat V$. 

Then 
\begin{align}\label{rels:alpha:g:gamma:g}
    & -\gamma_g e_0+\alpha_g e_1=-\gamma_g f_0+\gamma_g(f_0-e_0)+\alpha_g e_1
    =-f_0\gamma_g+\sum_{k\geq0}\varphi_k(e_0,e_1)e_1(e_0^k-f_0^k)+\sum_{k\geq0}\varphi_k(e_0,e_1)f_0^ke_1
    \\ & \nonumber=-f_0\gamma_g+\sum_{k\geq 0}\varphi_k(e_0,e_1)e_1e_0^k=-f_0\gamma_g+e_1g(e_0,e_1), 
\end{align}
where the second equality follows from \eqref{def:alpha:g:beta:g:3007}, the third equality follows from the commutation of $e_1$ and $f_0$, 
and the last equality follows from \eqref{def:varphi:k:ja:2024}. 

Then 
\begin{align*}
&
\mathrm{aut}_g^{V}(\rho_0) M_g
\\& = \begin{pmatrix}
    \mathrm{Ad}_{g(e_0,e_1)}(e_0)&0&0\\e_1&\mathrm{Ad}_{g(f_0,f_1)}(f_0)&-e_1\\0&0&\mathrm{Ad}_{g(e_0,e_1)}(e_0)
\end{pmatrix}
\begin{pmatrix}
  g(e_0,e_1)g(f_0,f_1)  &0&0\\
  -g(f_0,f_1)\gamma_g &g(f_0,f_1)\alpha_g&g(f_0,f_1)\gamma_g\\0&0&g(e_0,e_1)g(f_0,f_1)
\end{pmatrix}
\\& 
= \begin{pmatrix}
g(e_0,e_1)e_0g(f_0,f_1)&0&0\\
e_1g(e_0,e_1)g(f_0,f_1)-g(f_0,f_1)f_0\gamma_g&g(f_0,f_1)f_0\alpha_g & -e_1g(e_0,e_1)g(f_0,f_1)+g(f_0,f_1)f_0\gamma_g\\
0&0&g(e_0,e_1)e_0g(f_0,f_1)
\end{pmatrix}
\\& 
= \begin{pmatrix}
g(e_0,e_1)g(f_0,f_1)e_0&0&0\\
g(f_0,f_1)(e_1g(e_0,e_1)-f_0\gamma_g)&g(f_0,f_1)f_0\alpha_g & g(f_0,f_1)(-e_1g(e_0,e_1)+f_0\gamma_g)\\
0&0&g(e_0,e_1)g(f_0,f_1)e_0
\end{pmatrix}
\\& 
= \begin{pmatrix}
    g(e_0,e_1)g(f_0,f_1)e_0&0&0\\
    g(f_0,f_1)(-\gamma_g e_0+\alpha_g e_1)&g(f_0,f_1)\alpha_g f_0&g(f_0,f_1)(-\alpha_g e_1+\gamma_g e_0)
\\0&0&g(e_0,e_1)g(f_0,f_1)e_0
\end{pmatrix}
\\& 
=\begin{pmatrix}
  g(e_0,e_1)g(f_0,f_1)  &0&0\\
  -g(f_0,f_1)\gamma_g &g(f_0,f_1)\alpha_g&g(f_0,f_1)\gamma_g\\0&0&g(e_0,e_1)g(f_0,f_1)
\end{pmatrix}
\begin{pmatrix}
    e_0&0&0\\e_1&f_0&-e_1\\0&0&e_0
\end{pmatrix}=M_g\rho_0, 
\end{align*}
where the third equality follows from the commutation of $g(f_0,f_1)$ with $e_0$ and $e_1g(e_0,e_1)$, 
the fourth equality follows from \eqref{rels:alpha:g:gamma:g}, and the other equalities follow from definitions. 
This implies the first equality in \eqref{act:g:rho0:barrho0}. The second equality follows from the first by 
applying the morphism $T_{21}\hat V\to M_2\hat V$, $x\mapsto \overline x$. 
\end{proof}

\begin{lem}\label{lem:6:21.3112}
   (a) The action of $\mathcal G$ on the set $\mathrm{Hom}_{\mathcal C\operatorname{-alg}}^1(\hat{\mathcal V},M_3\hat V)$ from Lem. \ref{lem:520:2212:BIS}(b) 
    preserves the subsets $\mathrm{Hom}_{\mathcal C\operatorname{-alg}}^{1,(0)}(\hat{\mathcal V},T_{21}\hat V)$ and 
    $\mathrm{Hom}_{\mathcal C\operatorname{-alg}}^{1,((0)),\bullet}(\hat{\mathcal V},T_{21}\hat V)$.

     (b) The action of $\mathcal G$ from Lem. \ref{lem:621:jan2025} preserves the subset 
     $\mathrm{Hom}_{\mathcal C\operatorname{-alg}}^{1,(0)}(\hat{\mathcal V},M_2\hat V)$. 
\end{lem}

\begin{proof}
(a) Since the action $g\mapsto\mathrm{aut}_g^{V}$ of $\mathcal G$ on the algebra $M_3\hat V$ involved in the action of $\mathcal G$ on 
$\mathrm{Hom}_{\mathcal C\operatorname{-alg}}(\hat{\mathcal V},M_3\hat V)$ is entrywise, 
it restricts to an action on the subalgebra $T_{21}\hat V$, which implies that the action of $\mathcal G$ on 
the set $\mathrm{Hom}_{\mathcal C\operatorname{-alg}}(\hat{\mathcal V},M_3\hat V)$ from Lem. \ref{lem:520:2212:BIS}(a) preserves the 
subset 
\begin{equation}\label{subset:1:jan2025}
\mathrm{Hom}_{\mathcal C\operatorname{-alg}}(\hat{\mathcal V},T_{21}\hat V).     
\end{equation}
By Lem. \ref{lem:520:2212:BIS}(b), it also preserves the 
subset 
\begin{equation}\label{subset:2:jan2025}
\mathrm{Hom}_{\mathcal C\operatorname{-alg}}^1(\hat{\mathcal V},M_3\hat V). 
\end{equation}
If $\rho\in \mathrm{Hom}_{\mathcal C\operatorname{-alg}}(\hat{\mathcal V},T_{21}\hat V)$ is such that $\rho(e_0)$ is 
$(T_{21}\hat V)^\times$-conjugate to $\rho_0$, then if $t\in (T_{21}\hat V)^\times$ is such that that $\rho(e_0)=t\rho_0t^{-1}$, one has 
for any $g\in\mathcal G$ the relation
\begin{align}\label{calcul:qui:sera:répété}
    & (g*\rho)(e_0)=\mathrm{aut}_g^V\circ \rho\circ \mathrm{aut}^{\mathcal V}_{g^{\circledast-1}}(e_0)
=\mathrm{aut}_g^V\circ \rho(g^{\circledast-1}(e_0,e_1)e_0g^{\circledast-1}(e_0,e_1)^{-1})
\\ & \nonumber=
\mathrm{aut}_g^V\circ \rho(g^{\circledast-1}(e_0,e_1)) \cdot
\mathrm{aut}_g^V(t\rho_0t^{-1}) \cdot 
\mathrm{aut}_g^V\circ \rho(g^{\circledast-1}(e_0,e_1))^{-1}
\\&  \nonumber=\mathrm{aut}_g^V(\rho(g^{\circledast-1}(e_0,e_1))t)M_g\cdot \rho_0\cdot (\mathrm{aut}_g^V(\rho(g^{\circledast-1}(e_0,e_1))t)M_g)^{-1}
\end{align}
where $g^{\circledast-1}$ is the inverse of $g$ in $\mathcal G$ (equipped with its product $\circledast$), 
and where $\mathrm{aut}_g^V(\rho(g^{\circledast-1}(e_0,e_1))t)M_g
\in (T_{21}\hat V)^\times$; here the fourth equality follows from the first part of \eqref{act:g:rho0:barrho0}. 
It follows that the action of $\mathcal G$ on 
$ \mathrm{Hom}_{\mathcal C\operatorname{-alg}}(\hat{\mathcal V},T_{21}\hat V)$ preserves the subset 
\begin{equation}\label{subset:3:jan2025}
\{\rho\in \mathrm{Hom}_{\mathcal C\operatorname{-alg}}(\hat{\mathcal V},T_{21}\hat V)|\rho(e_0)\text{ is }
(T_{21}\hat V)^\times\text{-conjugate to }\rho_0\}. 
\end{equation}
If $\rho\in \mathrm{Hom}_{\mathcal C\operatorname{-alg}}(\hat{\mathcal V},T_{21}\hat V)$ is such that $\rho(e_0) \equiv \rho_0$ mod 
$T_{21}F^2\hat V$ and $g\in\mathcal G$, then 
\begin{align*}
    & (g*\rho)(e_0)=\mathrm{aut}_g^V\circ \rho\circ \mathrm{aut}^{\mathcal V}_{g^{\circledast-1}}(e_0)
=\mathrm{aut}_g^V\circ \rho(g^{\circledast-1}(e_0,e_1)e_0g^{\circledast-1}(e_0,e_1)^{-1})
    \\& \equiv \mathrm{aut}_g^V(\rho(e_0))\equiv \mathrm{aut}_g^V(\rho_0)\equiv \rho_0\textrm{ mod }T_{21}F^2\hat V, 
\end{align*}
where the third equality follows from $g^{\circledast-1}\equiv 1$ mod $F^1\hat{\mathcal V}$, 
 the fourth equality follows from the assumption on $\rho$, and  
 the fifth equality follows from $\mathrm{aut}_g^V(x) \equiv x$ mod $F^{d+1}\hat V$ for any $d\geq0$ and $x\in F^d\hat V$.  
Therefore the action of $\mathcal G$ on $\mathrm{Hom}_{\mathcal C\operatorname{-alg}}(\hat{\mathcal V},T_{21}\hat V)$ preserves the subset 
\begin{equation}\label{subset:supplementaire}
\{\rho\in \mathrm{Hom}_{\mathcal C\operatorname{-alg}}(\hat{\mathcal V},T_{21}\hat V)|\rho(e_0)\equiv\rho_0\textrm{ mod }T_{21}F^2\hat V\}. 
\end{equation}

Assume $g\in\mathcal G$ and $\rho\in \mathrm{Hom}_{\mathcal C\operatorname{-alg}}(\hat{\mathcal V},T_{21}\hat V)$ is such that for some 
$(r,C)\in \hat V\times M_{3,1}F^1\hat V$ with $R_r\cdot C\in e_0+f_\infty+F_2\hat V$, one has 
$\mathrm C_{21}(\rho(\hat{\mathcal V}))=\mathbf k+C \cdot \mathrm C_{\hat V}(e_0)\cdot R_r$. Set
$$
(\tilde r,\tilde C):=(g(e_0,e_1)^{-1}\mathrm{aut}_g^V(r),\mathrm{aut}_g^V(C)g(e_0,e_1))\in \hat V\times M_{3,1}F^1\hat V.  
$$ 
Then
$$
R_{\tilde r}\tilde C=R_{g(e_0,e_1)^{-1}\mathrm{aut}_g^V(r)}\cdot \mathrm{aut}_g^V(C)g(e_0,e_1)
=g(e_0,e_1)^{-1}\cdot R_{\mathrm{aut}_g^V(r)}\mathrm{aut}_g^V(C)\cdot g(e_0,e_1)
=\mathrm{Ad}_{g(e_0,e_1)^{-1}}\circ \mathrm{aut}_g^V(R_r C).
$$
Since $R_r\cdot C\in e_0+f_\infty+F_2\hat V$ and the automorphism $\mathrm{Ad}_{g(e_0,e_1)^{-1}}\circ \mathrm{aut}_g^V$ of
$\hat V$ is filtered with associated graded equal to the identity, it follows that 
$$
R_{\tilde r}\tilde C\in e_0+f_\infty+F_2\hat V.
$$
Moreover, one has 
\begin{align*}
    & \mathrm C_{21}((g*\rho)(\hat{\mathcal V}))
=\mathrm C_{21}(\mathrm{aut}_g^V(\rho(\hat{\mathcal V})))
=\mathrm{aut}_g^V(\mathrm C_{21}(\rho(\hat{\mathcal V})))
=\mathrm{aut}_g^V(\mathbf k+C \cdot \mathrm C_{\hat V}(e_0)\cdot R_r)\\ & 
=\mathbf k+\mathrm{aut}_g^V(C)\cdot\mathrm{aut}_g^V(\mathrm C_{\hat V}(e_0))
\cdot\mathrm{aut}_g^V(R_r)
=\mathbf k+\mathrm{aut}_g^V(C)g(e_0,e_1)\cdot\mathrm C_{\hat V}(e_0)
\cdot R_{g(e_0,e_1)^{-1}\mathrm{aut}_g^V(r)}
\\&=\mathbf k+\tilde C \cdot \mathrm C_{\hat V}(e_0)\cdot R_{\tilde r}, 
\end{align*}
where the first (resp. second, third, fourth) equality follows from $(g\bullet\rho)(\hat{\mathcal V})=\mathrm{aut}_g^V(\rho(\hat{\mathcal V}))$
(resp. the commutation of centralizers and automorphisms, the assumption on $\rho$, the compatibility of $\mathrm{aut}_g^V$ with products), 
and the fifth equality follows from 
$$
\mathrm{aut}_g^V(\mathrm C_{\hat V}(e_0))=\mathrm C_{\hat V}(\mathrm{aut}_g^V(e_0))
=\mathrm C_{\hat V}(g(e_0,e_1)e_0g(e_0,e_1)^{-1})
=g(e_0,e_1)\cdot\mathrm C_{\hat V}(e_0)\cdot g(e_0,e_1)^{-1}
$$
and $\mathrm{aut}_g^V(R_r)=R_{\mathrm{aut}_g^V(r)}$ and $g(e_0,e_1)^{-1}R_{\mathrm{aut}_g^V(r)}=R_{g(e_0,e_1)^{-1}\mathrm{aut}_g^V(r)}$. 
It follows that the action of $\mathcal G$ on 
$ \mathrm{Hom}_{\mathcal C\operatorname{-alg}}(\hat{\mathcal V},T_{21}\hat V)$ preserves the subset 
\begin{align}\label{subset:4:jan2025}
& \nonumber \{\rho\in \mathrm{Hom}_{\mathcal C\operatorname{-alg}}(\hat{\mathcal V},T_{21}\hat V)|
\text{for some } 
(r,C)\in \hat V\times M_{3,1}F^1\hat V
\text{ with }R_r\cdot C\in e_0+f_\infty+F_2\hat V,
\\& \text{ one has } 
\mathrm C_{21}(\rho(\hat{\mathcal V}))=\mathbf k+C \cdot \mathrm C_{\hat V}(e_0)\cdot R_r
\}.     
\end{align}
The result follows from the fact that $\mathrm{Hom}_{\mathcal C\operatorname{-alg}}^{1,(0)}(\hat{\mathcal V},T_{21}\hat V)$
(resp. $\mathrm{Hom}_{\mathcal C\operatorname{-alg}}^{1,((0)),\bullet}(\hat{\mathcal V},T_{21}\hat V)$) 
is the intersection of the subsets \eqref{subset:1:jan2025}, \eqref{subset:2:jan2025} and \eqref{subset:3:jan2025} of 
$\mathrm{Hom}_{\mathcal C\operatorname{-alg}}(\hat{\mathcal V},M_3\hat V)$ (resp. its intersection with the subsets \eqref{subset:1:jan2025}, \eqref{subset:2:jan2025},  
\eqref{subset:3:jan2025}, \eqref{subset:supplementaire} and \eqref{subset:4:jan2025}).  

(b) Let $g\in\mathcal G$ and $\sigma\in \mathrm{Hom}_{\mathcal C\operatorname{-alg}}(\hat{\mathcal V},M_2\hat V)$. 
If $\sigma$ is such that $\sigma(e_1)=\overline\rho_1$, then 
$$
(g*\sigma)(e_1)=\mathrm{aut}_g^V\circ \sigma\circ(\mathrm{aut}_g^{\mathcal V})^{-1}(e_1)
=\mathrm{aut}_g^V\circ \sigma(e_1)=\mathrm{aut}_g^V(\overline\rho_1)=\overline\rho_1, 
$$
where the equalities follow from $\mathrm{aut}_g^{\mathcal V}(e_1)=e_1$ and $\mathrm{aut}_g^V(\overline\rho_1)=\overline\rho_1$, which 
implies that the action of $\mathcal G$ on $\mathrm{Hom}_{\mathcal C\operatorname{-alg}}(\hat{\mathcal V},M_2\hat V)$ preserves 
\begin{equation}\label{subset:7:jan2025}
\{\sigma\in \mathrm{Hom}_{\mathcal C\operatorname{-alg}}(\hat{\mathcal V},M_2\hat V)|\sigma(e_1)=\overline\rho_1\}. 
\end{equation}
If now $\sigma$ is such that $\sigma(e_0)$ is $\mathrm{GL}_2\hat V$-conjugate to $\overline\rho_0$, then if $t\in \mathrm{GL}_2\hat V$ 
is such that that $\sigma(e_0)=t\overline\rho_0t^{-1}$, the analogue of the computation \eqref{calcul:qui:sera:répété} yields, 
using the second part of \eqref{act:g:rho0:barrho0} instead of its first part, 
$$
(g*\sigma)(e_0)=\mathrm{aut}_g^V\circ \sigma\circ \mathrm{aut}^{\mathcal V}_{g^{\circledast-1}}(e_0)
=\mathrm{aut}_g^V(\sigma(g^{\circledast-1}(e_0,e_1))t)\overline M_g\cdot \overline\rho_0\cdot 
(\mathrm{aut}_g^V(\sigma(g^{\circledast-1}(e_0,e_1))t)\overline M_g)^{-1}, 
$$
where $g^{\circledast-1}$ is as in (a), and where $\mathrm{aut}_g^V(\sigma(g^{\circledast-1}(e_0,e_1))t)\overline M_g
\in \mathrm{GL}_2\hat V$. It follows that the action of $\mathcal G$ on 
$\mathrm{Hom}_{\mathcal C\operatorname{-alg}}(\hat{\mathcal V},T_{21}\hat V)$ preserves the subset 
\begin{equation}\label{subset:8:jan2025}
\{\sigma\in \mathrm{Hom}_{\mathcal C\operatorname{-alg}}(\hat{\mathcal V},M_2\hat V)|\sigma(e_0)\text{ is }
\mathrm{GL}_2\hat V\text{-conjugate to }\overline\rho_0\}. 
\end{equation}
The result follows from the fact that $\mathrm{Hom}_{\mathcal C\operatorname{-alg}}^{1,(0)}(\hat{\mathcal V},M_2\hat V)$
 is the intersection of the subsets \eqref{subset:7:jan2025} and \eqref{subset:8:jan2025} of 
$\mathrm{Hom}_{\mathcal C\operatorname{-alg}}(\hat{\mathcal V},M_2\hat V)$.  
\end{proof}

\begin{lem}
(a) The map $\mathrm{Hom}_{\mathcal C\operatorname{-alg}}^{1,((0)),\bullet}(\hat{\mathcal V},T_{21}\hat V)\to 
\mathrm{Hom}_{\mathcal C\operatorname{-alg}}^{1,(0)}(\hat{\mathcal V},M_2\hat V)$
from Def. \ref{def:6:8:2412} is equivariant with respect to the actions of $\mathcal G$ from Lem. \ref{lem:6:21.3112}(a) and (b).  

(b) The map $\mathrm{Hom}_{\mathcal C\operatorname{-alg}}^{1,(0)}(\hat{\mathcal V},M_2\hat V)\to 
\mathrm{Hom}_{\mathcal C\operatorname{-alg}}(\hat{\mathcal W},\hat V)$ from Lem. \ref{lem:65:7:2412}(b) 
is equivariant with respect to the actions of $\mathcal G$ from Lem. \ref{lem:6:21.3112}(c) and Lem. \ref{lem:43:0401:TER}(b).
\end{lem}

\begin{proof}
(a) Since the action of $\mathcal G$ on both $T_{21}\hat V$ and $M_2\hat V$ is entrywise, the algebra morphism 
$T_{21}\hat V\to M_2\hat V$, $x\mapsto \overline x$ is $\mathcal G$-equivariant. It follows that 
the map $\mathrm{Hom}_{\mathcal C\operatorname{-alg}}(\hat{\mathcal V},T_{21}\hat V)\to 
\mathrm{Hom}_{\mathcal C\operatorname{-alg}}(\hat{\mathcal V},M_2\hat V)$ induced by composition with $x\mapsto \overline x$
is $\mathcal G$-equivariant as well. The statement then follows from the facts that this map induces the said map 
between subsets of its source and target, and that the action of $\mathcal G$ restricts to the said actions on these subsets. 

(b) Let $\sigma\in \mathrm{Hom}_{\mathcal C\operatorname{-alg}}^{1,(0)}(\hat{\mathcal V},M_2\hat V)$ and 
$g\in\mathcal G$. For any $a\in\hat{\mathcal V}$, one has 
\begin{align*}
    &
    \Delta_{g*\sigma}(ae_1)=\overline{\mathrm{row}}_{\mathrm{DT}}\cdot (g*\sigma)(a)\cdot \overline{\mathrm{col}}_{\mathrm{DT}}
=\overline{\mathrm{row}}_{\mathrm{DT}}\cdot (\mathrm{aut}_g^V\circ\sigma\circ 
(\mathrm{aut}_g^{\mathcal V})^{-1})(a)\cdot \overline{\mathrm{col}}_{\mathrm{DT}}
    \\ & 
    = \mathrm{aut}_g^V( \overline{\mathrm{row}}_{\mathrm{DT}}\cdot ( 
\sigma  ((\mathrm{aut}_g^{\mathcal V})^{-1}(a))\cdot \overline{\mathrm{col}}_{\mathrm{DT}} ) 
=\mathrm{aut}_g^V(\Delta_\sigma\left((\mathrm{aut}_g^{\mathcal V})^{-1}(a)e_1\right))
\\& =\mathrm{aut}_g^V(\Delta_\sigma((\mathrm{aut}_g^{\mathcal W})^{-1}(ae_1)))
=\mathrm{aut}_g^V\circ \Delta_\sigma\circ (\mathrm{aut}_g^{\mathcal W})^{-1}(ae_1)
=(g*\Delta_\rho)(ae_1), 
\end{align*}
where the first and fourth equalities follow from Lem. \ref{lem:65:7:2412}(b), 
the second (resp. last) equality follows from the definition of the action of $\mathcal G$ on 
$\mathrm{Hom}_{\mathcal C\operatorname{-alg}}^1(\hat{\mathcal V},M_2\hat V)$ (resp. 
$\mathrm{Hom}_{\mathcal C\operatorname{-alg}}^1(\hat{\mathcal W},\hat V)$), 
the third equality follows from the invariance of $\overline{\mathrm{col}}_{\mathrm{DT}}$ and 
$\overline{\mathrm{row}}_{\mathrm{DT}}$ under the action of 
$\mathcal G$, and the fifth equality follows from the definition of $\mathrm{aut}_g^{\mathcal W}$. 
\end{proof}

\subsubsection{Action of $\mathcal G$ on groups}

\begin{lem}\label{lem:6:25:jan}
    (a) The action of $\mathcal G$ on the group $\mathrm{GL}_3\hat V$ from Lem. \ref{lem:5:18:2012}
    induces an action of $\mathcal G$ on the subgroups $\mathrm C_{21}(\rho_1)^\times$ and $\mathrm C_{21}^{(0)}(\rho_1)^\times$
    (see Definition \ref{def:C0}).

    (b) The map $(g,P)\mapsto g*P:=\mathrm{aut}_g^V(P)$ defines an action of $\mathcal G$ on the group $\mathrm{GL}_2\hat V$. 

    (c) The action of (b) restricts to an action of $\mathcal G$ on the subgroup $\mathrm C_{2}(\overline\rho_1)^\times$. 

    (d) The morphism $\mathrm C_{21}^{(0)}(\rho_1)^\times\to \mathrm C_{2}(\overline\rho_1)^\times$ from \eqref{diagram:gp:morphisms} 
    is equivariant with respect to the actions of $\mathcal G$ defined in (a) and (c). 

    (e) The morphism from Lem. \ref{lem:6:16:2912}(c) induces a group morphism $\mathrm C_{2}(\overline\rho_1)^\times\to\mathbf k[[u,v]]^\times$, which is 
 $\mathcal G$-invariant, the action of  $\mathcal G$ on the source being as in (c). 

    (f) The sequence of group morphisms \eqref{diagram:gp:morphisms} is $\mathcal G$-equivariant, the actions of $\mathcal G$ being given by  
    Lem. \ref{lem:5:18:2012}, (a), (c) and the action being trivial on $\mathbf k[[u,v]]^\times$. 
\end{lem}

\begin{proof}
(a) The action of $\mathcal G$ on $M_3\hat V$ preserves the subalgebra $T_{21}\hat V$ by the proof of Lem. \ref{lem:6:21.3112}(a),
and the subalgebra $\mathrm C_3(\rho_1)$ by Lem. \ref{lem:5:18:2012}, which implies that it preserves their intersection, which is 
$\mathrm C_{21}(\rho_1)$. The augmentation morphism $\hat V\to \hat V_0=\mathbf k$ is $\mathcal G$-invariant, which implies the 
$\mathcal G$-invariance of the algebra morphism $T_{21}\hat V\to (T_{21}\hat V)_0$, and therefore the fact that the action of 
$\mathcal G$ on $T_{21}\hat V$ preserves the preimage under this morphism of any subalgebra of $(T_{21}\hat V)_0$, and in particular the preimage 
of $\mathrm C_{21}(\rho_0)_0$, and therefore preserves also intersection of this preimage with $\mathrm C_{21}(\rho_1)$, which is $\mathrm C_{21}^{(0)}(\rho_1)$. 
It therefore induces an action of $\mathcal G$ on the groups of units of $\mathrm C_{21}(\rho_1)$ and $\mathrm C_{21}^{(0)}(\rho_1)$, which 
are $\mathrm C_{21}(\rho_1)^\times$ and $\mathrm C_{21}^{(0)}(\rho_1)^\times$.
(b) This is the action induced by the entrywise action of $\mathcal G$ on $M_2\hat V$. (c) Since $\overline\rho_1$ is $\mathcal G$-invariant, 
the action of $\mathcal G$ on $M_2\hat V$ restricts to an action on the algebra $\mathrm C_2(\overline\rho_1)$, which then induces an action
on the group $\mathrm C_2(\overline\rho_1)^\times$. (d)  The inclusion morphism 
$\mathrm C_2^{(0)}(\overline\rho_1)\subset \mathrm C_2(\overline\rho_1)$ is $\mathcal G$-equivariant. 
Since the action of $\mathcal G$ on $T_{21}\hat V$ is entrywise, the algebra morphism 
$T_{21}\hat V\to M_2\hat V$, $x\mapsto \overline x$ is $\mathcal G$-equivariant. As this induces a morphism 
$\mathrm C_{21}(\rho_1)\to \mathrm C_2(\overline\rho_1)$ where the source and target as preserved by $\mathcal G$, 
this morphism is $\mathcal G$-equivariant as well. The composed morphism $\mathrm C_{21}^{(0)}(\rho_1)\to \mathrm C_{2}(\overline\rho_1)$ 
is therefore $\mathcal G$-equivariant. The result then follows from the fact that the said morphism is the associated morphism between groups of units. 
(e) The map $(g,(\phi,v))\mapsto g*(\phi,v):=(\phi,\mathrm{aut}_g^{V}(v))$ defines an action of 
$\mathcal G$ on $\mathbf k[[u,v]]\times\hat V$, and it is such that  
for any $g\in\mathcal G$, $\phi\in\mathbf k[[u,v]]$, $v\in\hat V$, one has
\begin{align*}
    &
    \overline M(g*(\phi,v))=\overline M(\phi,\mathrm{aut}_g^{V}(v))
=\phi(e_1,f_1) I_2+\begin{pmatrix}
f_1\\ e_1 
\end{pmatrix}\mathrm{aut}_g^{V}(v)\begin{pmatrix}
    1&1
\end{pmatrix}
\\& =\mathrm{aut}_g^{V}(\phi(e_1,f_1) I_2+\begin{pmatrix}
f_1\\ e_1 
\end{pmatrix}v\begin{pmatrix}
    1&1
\end{pmatrix})
    =\mathrm{aut}_g^{V}(\overline M(\phi,v)),
\end{align*}
which implies that the bijection $\mathbf k[[u,v]]\times\hat V\to \mathrm C_2(\overline\rho_1)$ induced by 
$\overline M$ (see Lem. \ref{lem:6:16:2912}(b)) is $\mathcal G$-equivariant. As the projection 
$\mathbf k[[u,v]]\times\hat V\to \mathbf k[[u,v]]$ is $\mathcal G$-invariant, is follows that the map 
from Lem. \ref{lem:6:16:2912}(c) is $\mathcal G$-invariant. This map is an algebra morphism, therefore 
the underlying group morphism $\mathrm C_{21}(\rho_1)^\times\to \mathrm C_{2}(\overline\rho_1)^\times$ 
is $\mathcal G$-invariant. (f) then follows from (a)-(e). 
\end{proof}

\subsubsection{Compatibility of actions of $\mathcal G$ on groups and on sets}

\begin{lem}\label{lem:G:eq:diag}
(a) When equipped with the actions $*$ of $\mathcal G$ defined in Lems. 
\ref{lem:621:jan2025}, \ref{lem:6:21.3112} and \ref{lem:6:25:jan},
each of the sets with group action $\bullet$ from diagram \eqref{diag:set:gp:acts} satisfies the identity of Def. \ref{def:psga}(a), 
and therefore builds up an object in $\mathcal G\operatorname{-}\mathbf{PSGA}$.

(b) Equipping the objects of this diagram with these actions of $\mathcal G$, \eqref{diag:set:gp:acts} is upgraded to the following diagram in 
$\mathcal G\operatorname{-}\mathbf{PSGA}$
\begin{equation}\label{diag:set:gp:acts:with:*}
   \xymatrix@C=20pt{
\substack{(\mathrm{Hom}^{1,(0)}_{\mathcal C\operatorname{-alg}}(\hat{\mathcal V},T_{21}\hat V),\\ \mathrm C_{21}(\rho_1)^\times,
\rho_{\mathrm{DT}},\bullet,*)}
\ar@{_{(}->}[d]
&\ar@{_{(}->}[l]\substack{(\mathrm{Hom}^{1,((0)),\bullet}_{\mathcal C\operatorname{-alg}}
(\hat{\mathcal V},T_{21}\hat V),\\ \mathrm C_{21}^{(0)}(\rho_1)^\times,\rho_{\mathrm{DT}},\bullet,*)}
\ar[r]
&
\substack{(\mathrm{Hom}^{1,(0)}_{\mathcal C\operatorname{-alg}}(\hat{\mathcal V},M_2\hat V),
\\\mathrm{C}_2(\overline\rho_1)^\times,\overline\rho_{\mathrm{DT}},\bullet,*)}
\ar[d]
\\
\substack{(\mathrm{Hom}^{1}_{\mathcal C\operatorname{-alg}}(\hat{\mathcal V},M_3\hat V),\\
\mathrm{C}_3(\rho_1)^\times,\rho_{\mathrm{DT}},\bullet,*)}&
&
\substack{(\mathrm{Hom}_{\mathcal C\operatorname{-alg}}(\hat{\mathcal W},\hat V),\\
\mathbf k[[u,v]]^\times,\Delta^{\mathcal W}_{r,l},\bullet,*)} 
}
\end{equation}
\end{lem}

\begin{proof}
(a) By Lem. \ref{lem:pre:psga}(b), 
$(\mathrm{Hom}^1_{\mathcal C\operatorname{-alg}}(\hat{\mathcal V},M_3\hat V),\mathrm{C}_3(\rho_1),\bullet,*)$
satisfies the identity of Def. \ref{def:psga}(a). 
The fact that 
$$
(\mathrm{Hom}^{1,(0)}_{\mathcal C\operatorname{-alg}}(\hat{\mathcal V},T_{21}\hat V),\mathrm{C}_{21}(\rho_1),\bullet,*)\textrm{ and }
(\mathrm{Hom}^{1,(0),\bullet}_{\mathcal C\operatorname{-alg}}(\hat{\mathcal V},T_{21}\hat V),\mathrm{C}_{21}(\rho_1),\bullet,*)
$$
satisfy the identity of Def. \ref{def:psga}(a) then follows from the following statement: if $(X,A,\bullet,*)$ satisfies the identity of Def. \ref{def:psga}(a), 
and if $X'\subset X$ is a subset, $A'\subset A$ is a subgroup which are both preserved by $\mathcal G$ and such that the action of $A'$ preserves $X'$, then $(X',A',\bullet)$ 
satisfies the identity of Def. \ref{def:psga}(a); and from the fact that the assumptions  
of this statement are satisfied with $(X,A,\bullet,*)=
(\mathrm{Hom}^1_{\mathcal C\operatorname{-alg}}(\hat{\mathcal V},M_3\hat V),\mathrm{C}_3(\rho_1),\bullet,*)
$ 
and $(X',A',\bullet)$ equal to one of these tuples. 

Replacing $M_3\hat V$ by $M_2\hat V$ in the proof of Lem. \ref{lem:pre:psga}(a), one shows that  
$(\mathrm{Hom}_{\mathcal C\operatorname{-alg}}(\hat{\mathcal V},M_2\hat V),\mathrm{GL}_2\hat V,\bullet,*)$ satisfies the identity of Def. \ref{def:psga}(a). 
The fact the pair formed by this tuple and by 
$$
(\mathrm{Hom}_{\mathcal C\operatorname{-alg}}^{1,(0)}(\hat{\mathcal V},M_2\hat V),\mathrm{C}_2(\overline\rho_1),\bullet,*)
$$ 
satisfy the assumptions of the previous general statement implies that the latter tuple satisfies the identity of Def. \ref{def:psga}(a). 

Finally, the fact that the tuple $(\mathbf k[[u,v]]^\times,\mathrm{Hom}_{\mathcal C\operatorname{-alg}}(\hat{\mathcal W},\hat V),\bullet,*)$
is derived from the tuple $(\mathbf H,\mathbf k[[u,v]]^\times,\bullet,*)$, which has been shown in Cor. \ref{cor:HH} to be an object of 
$\mathcal G\operatorname{-}\mathbf{PSGA}$, therefore it satisfies the identity of Def. \ref{def:psga}(a). 

(b) follows from (a) and from the $\mathcal G$-equivariance of the morphisms in \eqref{diag:set:gp:acts}, which follows from Lem. \ref{lem:6:25:jan}(f). 
\end{proof}

\subsection{A diagram of pointed sets with actions of $\mathcal G$}\label{sect:6:7:5786}

\begin{lem}\label{lem:6:30:3103}
The diagram of pointed sets with $\mathcal G$-action (a diagram in $\mathbf{PS}_{\mathcal G}$) obtained by 
applying the functor $\mathbf q$ to diagram \eqref{diag:set:gp:acts:with:*} from Lem. \ref{lem:G:eq:diag}(b) is 
\eqref{the:big:diagram}. 
\end{lem}

\begin{proof}
    Immediate. 
\end{proof}

We will prove that 
(C) is injective and that (B),(D),(E) are locally injective. 

\section{Injectivity of the map 
 (C)}\label{sec 7}

\begin{lem}\label{lem:C:inj}
    The map (C) (see \eqref{the:big:diagram}) is injective.  
\end{lem}

\begin{proof}
Let $\alpha,\beta\in \mathrm C_{21}^{(0)}(\rho_1)^\times\backslash\mathrm{Hom}^{1,((0)),\bullet}_{\mathcal C\operatorname{-alg}}(\hat{\mathcal V},T_{21}\hat V)$ 
be elements with the same image in 
$$
\mathrm C_{21}(\rho_1)^\times\backslash\mathrm{Hom}^{1,(0)}_{\mathcal C\operatorname{-alg}}(\hat{\mathcal V},T_{21}\hat V).
$$
Let $\rho_\alpha,\rho_\beta\in \mathrm{Hom}^{1,((0)),\bullet}_{\mathcal C\operatorname{-alg}}(\hat{\mathcal V},T_{21}\hat V)$ be representatives of
$\alpha,\beta$. Then for some $g\in \mathrm C_{21}(\rho_1)^\times$; one has $\rho_\beta=\mathrm{Ad}_g\circ\rho_\alpha$. One has therefore
$\rho_\beta(e_0)=\mathrm{Ad}_g\circ\rho_\alpha(e_0)$. Since $\rho_\alpha(e_0)$ and $\rho_\beta(e_0)$ belong to $\rho_0+F^2T_{21}\hat V$, this is an equality in 
$F^1T_{21}\hat V$. Its image in $F^1T_{21}\hat V/F^2T_{21}\hat V=T_{21}\hat V_1$ is the equality $\rho_0=\mathrm{Ad}_{g_0}(\rho_0)$, which implies 
$g_0\in \mathrm C(\rho_0)_0^\times$, therefore $g\in  \mathrm C_{21}^{(0)}(\rho_1)^\times$. It follows that $\alpha=\beta$. 
\end{proof}


\section{Local injectivity of the morphism (B)}\label{sec 8}

This section is devoted to the proof of the local injectivity of the map (B). 
In §\ref{present:section:1}, we compute the commutant of $\rho_1$ in $T_{21}\hat V$, and then its group of 
invertible elements.  In §\ref{sect:8:2:5786}, we compute the kernel of a linear map with target $\hat V$.  This is 
used in §\ref{present:section} for the computation of the commutants over $\overline\rho_0$ in $M_2\hat V$, and of $\rho_0$ 
in $T_{21}\hat V$ and in $M_2\hat V$. In §\ref{sect:8:4:5786}, we compute the groups of invertible elements of these commutants. 
We use these results in §\ref{sect:8:5:5786} to prove a relation between these groups, and we then derive the local injectivity 
in (B) in §\ref{sect:8:6:5786}.  

\subsection{Computation of $\mathrm{C}_{21}(\rho_1)$ and $\mathrm{C}_{21}(\rho_1)^\times$ }\label{present:section:1}

\begin{lem}\label{lem:15:12} 
One has 
$$
\mathrm{C}_{21}(\rho_1)
=\{M(\phi,m)|\phi\in\mathbf k[[u,v]],
m\in T_{11}\hat V\},  
$$
where $T_{11}\hat V=\{\begin{pmatrix}\alpha&\gamma\\0&\delta\end{pmatrix}|\alpha,\gamma,\delta\in\hat V\}\subset M_2\hat V$ and 
the map $\mathbf k[[u,v]]\oplus M_2\hat V\to M_3\hat V$, $(\phi,m)\mapsto M(\phi,m)$
is as in \eqref{def:M:2012}. 
\end{lem}

\begin{proof}
One has 
\begin{align*}
&\mathrm{C}_{21}(\rho_1)=\mathrm{C}_{3}(\rho_1)\cap T_{21}\hat V=\{M(\phi,m)|(\phi,m)\in 
\mathbf k[[u,v]]\times M_2\hat V\quad\text{and}\quad M(\phi,m)\in  T_{21}\hat V\}.
\\& =\{M(\phi,m)|(\phi,m)\in 
\mathbf k[[u,v]]\times T_{11}\hat V\}, 
\end{align*} 
where the second equality follows from the fact that the map $(\phi,m)\mapsto M(\phi,m)$ sets up a bijection 
$\mathbf k[[u,v]]\oplus M_2\hat V\to \mathrm C_3(\rho_1)$ (see Lem. \ref{lem:13:5:1108:BIS}(b)), and the last equality follows from the equivalence  
$(M(\phi,m)\in T_{21}\hat V)\iff(\beta=0)\iff (m\in T_{11}\hat V)$ for 
$(\phi,m)\in \mathbf k[[u,v]]\oplus M_2\hat V$ and $\alpha,\beta,\gamma,\delta\in\hat V$ such that $m=\begin{pmatrix}
    \alpha&\gamma\\\beta&\delta
\end{pmatrix}$, which follows from the explicit form \eqref{form:of:commutant:2912} of the map  $(\phi,m)\mapsto M(\phi,m)$. 
\end{proof}

\begin{lem}\label{lem:c:21:rho:1}
    One has 
    $$
    \mathrm C_{21}(\rho_1)^\times=\{M(\phi,m)|\phi\in\mathbf k[[u,v]]^\times,
m=\begin{pmatrix}\alpha&\gamma\\0&\delta\end{pmatrix}\in T_{11}\hat V \quad \text{and}\quad  \phi+\delta\in \hat V^\times\}, 
    $$
where $(\phi,m)\mapsto M(\phi,m)$ is as in \eqref{def:M:2012} and $\epsilon : \hat V\to\mathbf k$ is the augmentation morphism. 
\end{lem}

\begin{proof}
    By Lem. \ref{invertible:commutant}, one has $\mathrm C_{21}(\rho_1)^\times=\mathrm C_{21}(\rho_1)\cap (T_{21}\hat V)^\times$, which 
by 
\begin{equation}\label{invertibles:t:21}
    (T_{21}\hat V)^\times=T_{21}\hat V\cap\mathrm{GL}_3\hat V
\end{equation}
and $\mathrm C_{21}(\rho_1)\subset T_{21}\hat V$ implies $\mathrm C_{21}(\rho_1)^\times=\mathrm C_{21}(\rho_1)\cap \mathrm{GL}_3\hat V$. 
Lem. \ref{lem:15:12} then implies 
$$
\mathrm C_{21}(\rho_1)^\times=\{M(\phi,m)|\phi\in\mathbf k[[u,v]],
m\in T_{11}\hat V \quad \text{and}\quad  \epsilon(M(\phi,m))\in\mathrm{GL}_3\mathbf k\}.     
    $$
For $m=\begin{pmatrix}\alpha&\gamma\\0&\delta\end{pmatrix}\in T_{11}\hat V$, one computes
$ \epsilon(M(\phi,m))=\mathrm{diag}(\phi(0,0),\phi(0,0),\phi(0,0)+\epsilon(\delta))$, which then implies the result. 
\end{proof}

\subsection{Solution of the equation $f_0 v-v e_0=e_1 \mathbf e-\mathbf f e_1$}\label{sect:8:2:5786}

\begin{defn}
    $\varphi$ is the endomorphism of $V$ defined by $v\mapsto v\cdot e_0-f_0\cdot v$. 
\end{defn}

In the present \S\ref{present:section}, we define: 

$\bullet$ for $p,q\geq0$, $V_{p,q}$ to be the part of $V$ of degree $(p,q)$ with respect to the
degree for which $e_1$ has degree $(1,0)$, $f_1$ has degree $(0,1)$, and $e_0$ and $f_0$ both have degree $(0,0)$. 
Then $V=\oplus_{p,q\geq 0}V_{p,q}$. 

$\bullet$ $m : V\to\mathcal V$ to be the product map $V=\mathcal V^{\otimes 2}\ni a\otimes b\mapsto 
a\cdot b\in\mathcal V$. 

\begin{lem}\label{lem:ressource:0707}
(a) One has $\mathrm{C}_V(e_0)=\oplus_{q\geq0}V_{0,q}$ and $\mathrm{C}_V(f_0)=\oplus_{p\geq0}V_{p,0}$.  

(b) The map $\varphi$ is compatible with the bigrading of $V$, i.e. induces an endomorphism of $V_{p,q}$ for any $p,q\geq0$. 

(c) One has $m\circ\varphi=0$. 

(d) For any $q>0$, the sequence $V_{0,q}\stackrel{\varphi}{\to}V_{0,q}\stackrel{m}{\to}\mathcal V$ is exact. 

(e) One has the identity $m(e_1\cdot x)=e_1\cdot m(x)$ for any $x\in V$. 

(f) One has the identity $\varphi(e_1\cdot x)=e_1\cdot \varphi(x)$ for any $x\in V$. 

(g) The map $\mathrm{C}_V(f_0)\to\mathcal V$, $x\mapsto m(x\cdot e_1)$ is injective. 

(h) $m(V_{0,0})\subset \mathbf k[e_0]$. 

(i) The preimage of $e_1\cdot \mathbf k[e_0]$ under the map $V_{0,0}\to\mathcal V$, 
$x\mapsto m(x\cdot e_1)$ is $\mathbf k[f_0]$. 
\end{lem}

\begin{proof}
The first statement of (a) follows from Lem. \ref{lem:comm:e0}(a); the second statement follows from 
$\mathrm{C}_V(f_0)=\mathcal V\otimes\mathbf k[e_0]$, which follows from it by applying the automorphism 
of exchange of factors of $V=\mathcal V^{\otimes2}$. (b) follows from the fact that the bigrading of $V$ is compatible 
with the algebra structure, and that $e_0$ and $f_0$ both have bidegree $(0,0)$. (c) follows from the equalities 
$m\circ\varphi(a\otimes b)=m(ae_0\otimes b-a\otimes e_0b)=(ae_0)b-a(e_0b)=0$ for any $a,b\in\mathcal V$. 

Let us prove (d). For $n\geq 0$, denote by $\mathcal V_n$ the part of $\mathcal V$ of degree $n$, where $e_1$ has degree $1$
and $e_0$ has degree $0$. The complex $\mathbf k[u,v]\to\mathbf k[u,v]\to\mathbf k[u]$, where the first map is 
$P(u,v)\mapsto P(u,v)\cdot(u-v)$ and the second map is $Q(u,v)\mapsto Q(u,u)$. This complex is known to be exact,
which since $e_1\mathcal V_{q-1}$ is a free $\mathbf k$-module implies the exactness of its tensor product 
with $e_1\mathcal V_{q-1}$, which is
\begin{equation}\label{interm:complex:0807}
\mathbf k[u,v]\otimes e_1\mathcal V_{q-1}\to\mathbf k[u,v]\otimes e_1\mathcal V_{q-1}
\to\mathbf k[u]\otimes e_1\mathcal V_{q-1}. 
\end{equation}
There are isomorphisms 
$\mathbf k[u,v]\otimes e_1\mathcal V_{q-1}\to V_{0,q}$ and $\mathbf k[u]\otimes e_1\mathcal V_{q-1}\to V_{0,q}$, 
respectively given by $u^pv^q\otimes x\mapsto e_0^p\otimes (e_0^qx)$ and $u^p\otimes x\mapsto e_0^px$. These 
isomorphisms take \eqref{interm:complex:0807} to the complex 
$V_{0,q}\stackrel{\varphi}{\to}V_{0,q}\stackrel{m}{\to}\mathcal V_q$, which is therefore exact. 
As $m : V_{0,q}\to\mathcal V$ is the composition $V_{0,q}\stackrel{m}{\to}\mathcal V_q\subset\mathcal V$, 
this implies exactness claimed in (d). 

(e) follows the the associativity of the product in $\mathcal V$. 
(f) follows from the commutation of $f_0$ and $e_1$. 

Let us prove (g). One checks that the $\mathbf k$-module morphisms $\alpha : \mathcal V\otimes \mathbf k[t]\to 
\oplus_{q>0}\mathcal V_q$ induced by $v\otimes t^p\mapsto v\cdot e_1e_0^p$ and $\beta : \oplus_{q>0}\mathcal V_q\to 
\mathcal V\otimes \mathbf k[t]$ induced by $x\mapsto \sum_{p\geq 0} \partial_1\partial_0^p(x)\otimes t^p$ are mutually 
inverse, where $\partial_0,\partial_1$ are the endomorphisms of $\mathcal V$ defined by
$x=\epsilon(x)1+\partial_0(x)e_0+\partial_1(x)e_1$, $\epsilon$ being the augmentation map of $\mathcal V$. It follows that 
$\alpha$ is a linear isomorphism, therefore that the map $\mathcal V\otimes \mathbf k[t]\to\mathcal V$ induced by 
$v\otimes t^p\mapsto v\cdot e_1e_0^p$ is injective. One the other hand, the map 
$\mathcal V\otimes\mathbf k[t]\to\mathrm{C}_V(f_0)$, $v\otimes t^p\mapsto v\otimes e_0^p$
is a $\mathbf k$-module isomorphism. (g) then follows from the fact that the map 
$\mathrm{C}_V(f_0)\to\mathcal V$, $x\mapsto m(x\cdot e_1)$ is the composition of 
the map $\mathcal V\otimes \mathbf k[t]\to\mathcal V$ with the inverse to the isomorphism 
$\mathcal V\otimes\mathbf k[t]\to\mathrm{C}_V(f_0)$. 

(h) follows from $V_{0,0}=\mathbf k[e_0,f_0]=\mathbf k[e_0]^{\otimes2}$.
(i) The composition of the said map with the algebra isomorphism $\mathbf k[u,v]\to V_{0,0}$ given by 
$u\mapsto e_0$, $v\mapsto f_0$ is $u^pv^q\mapsto e_0^pe_1e_0^q$, so the image of $P(u,v)=\sum_{p,q}a_{p,q}u^pv^q$
is $\sum_{p,q}a_{p,q}e_0^pe_1e_0^q$. It belongs to $e_1\cdot \mathbf k[e_0]$ if and only if $a_{p,q}=0$ for any $p>0$, 
which is equivalent to $P(u,v)\in\mathbf k[v]$, i.e. to the statement that the corresponding element of 
$V_{0,0}=\mathbf k[e_0,f_0]$ belongs to $\mathbf k[f_0]$. 
\end{proof}

\begin{lem}\label{pbms:1:2} 
Let $v\in \hat V$, $\mathbf f\in\mathrm{C}_{\hat V}(f_0)$, $\mathbf e\in\mathrm{C}_{\hat V}(e_0)$ 
be such that 
\begin{equation}\label{main:condition} 
f_0 v-v e_0=e_1 \mathbf e-\mathbf f e_1
\end{equation}
(equality in $\hat V$). Then there exists $\Pi\in \mathbf k[[t]]$ and 
$C\in \mathrm{C}_{\hat V}(e_0)$, such that 
$$
(\mathbf f,v,\mathbf e)=(0,e_1 C,(f_0-e_0) C)+(\Pi,0,\Pi).
$$ 
\end{lem}

\begin{proof}
Define $\mathbf k$-modules
$$
A:=\mathbf k[f_0]\oplus \mathrm{C}_V(e_0),\quad 
B:=\mathrm{C}_V(f_0)\oplus V\oplus \mathrm{C}_V(e_0) 
$$
and $\mathbf k$-module morphisms
\begin{equation}\label{eq:alpha:0807}
\alpha : A\to B,\quad
(\Pi,C)\mapsto (0,e_1 C,(f_0-e_0) C)+(\Pi,0,\Pi)
\end{equation}
and 
\begin{equation}\label{eq:beta:0807}
\beta : B\to V,\quad 
(\mathbf f,v,\mathbf e)\mapsto \varphi(v) +e_1\mathbf e-\mathbf f e_1. 
\end{equation}
It follows from Lem. \ref{lem:ressource:0707}(f), from the identity $\varphi(C)=(e_0-f_0)C$ for $C\in \mathrm{C}_V(e_0)$
and from the commutation of $e_1$ and $f_0$ that $\beta\circ\alpha=0$, so that there is a complex
\begin{equation}\label{exact:seq:abv}
A\stackrel{\alpha}{\to}B\stackrel{\beta}{\to}V. 
\end{equation}
Defined bigradings on $A$ and $B$ by 
$$
\forall p,q\geq0,\quad 
A_{p,q}:=\mathbf k[f_0]_{p,q}\oplus \mathrm{C}_V(e_0)_{p-1,q},\quad
B_{p,q}:=\mathrm{C}_V(f_0)_{p-1,q}\oplus V_{p,q}\oplus \mathrm{C}_V(e_0)_{p-1,q}, 
$$
where $\mathbf k[f_0]_{1,0}=\mathbf k[f_0]$ and $\mathbf k[f_0]_{p,q}=0$ if $(p,q)\neq(1,0)$, and the gradings on 
$\mathrm{C}_V(f_0)$ and $\mathrm{C}_V(e_0)$ are those induced by their status of graded submodules of
$V$ (see Lem. \ref{lem:ressource:0707}(a)). The maps $\alpha,\beta$ are compatible with these bigradings, 
therefore the above complex splits up as the direct sum over $p,q\geq 0$ of complexes 
$$
\mathbf C_{p,q}:=(A_{p,q}\stackrel{\alpha}{\to}B_{p,q}\stackrel{\beta}{\to}V_{p,q}). 
$$
\begin{itemize}
    \item If $p=0$, or if $p>1$ and $q>0$, then $\mathbf C_{p,q}=(0\to V_{p,q}\stackrel{\varphi}{\to}V_{p,q})$. It then follows 
from the injectivity of $\varphi$ (see Lem. \ref{lem:technical:0801}(a)) that $\mathbf C_{p,q}$ is exact. 
    \item
If $p=1$ and $q>0$, then by Lem. \ref{lem:ressource:0707}(a) one has 
$\mathbf C_{p,q}=(V_{0,q}\stackrel{\alpha}{\to} V_{1,q}\oplus \mathrm{C}_V(e_0)_{0,q}
\stackrel{\beta}{\to}V_{1,q})$, where $\alpha(C)=(e_1 C,(f_0-e_0) C)$ and 
$\beta(v,\mathbf e)=\varphi(v)+e_1 \mathbf e$. Let $(v,\mathbf e)\in\mathrm{ker}(\beta)$. Then 
$\varphi(v)+e_1 \mathbf e=0$. Applying $m$ and using Lem. \ref{lem:ressource:0707}(c), one obtains $m(e_1 \mathbf e)=0$, 
which by Lem. \ref{lem:ressource:0707}(e) implies $e_1\cdot m(\mathbf e)=0$. The fact that $\mathcal V$ is a domain then 
implies $m(\mathbf e)=0$. Lem. \ref{lem:ressource:0707}(d) then implies the existence of $C\in V_{0,q}$, such that 
$\mathbf e=\varphi(C)$. Then $\varphi(v)=-e_1 \mathbf e=-e_1 \varphi(C)=\varphi(-e_1 C)$, where the last 
equality follows from Lem. \ref{lem:ressource:0707}(f). Lem. \ref{lem:technical:0801}(a) then implies $v=-e_1 C$. 
Moreover, $\mathbf e=\varphi(C)=(e_0-f_0)\cdot C$, where the last equality follows from
$C\in V_{0,q}\subset \mathrm{C}_V(e_0)$. All this implies $(v,\mathbf e)=-\alpha(C)$, therefore
$(v,\mathbf e)\in\mathrm{im}(\alpha)$. This proves the exactness of $\mathbf C_{p,q}$.  
    \item 
If $q=0$ and $p>1$, then $\mathbf C_{p,q}=(0\to \mathrm{C}_V(f_0)_{p-1,0}\oplus 
V_{p,0}\stackrel{\beta}{\to}V_{p,0})$, where  
$\beta(\mathbf f,v)=\varphi(v)-\mathbf f e_1$. Let $(\mathbf f,v)\in\mathrm{ker}(\alpha)$. Then 
$\varphi(v)=\mathbf f e_1$. Applying $m$ and using Lem. \ref{lem:ressource:0707}(c), one obtains 
$m(\mathbf f e_1)=0$. Lem. \ref{lem:ressource:0707}(g) then implies $\mathbf f=0$. It follows that 
$\varphi(v)=0$, which by Lem.  \ref{lem:technical:0801}(a) implies $v=0$. Hence $(\mathbf f,v)=0$. This implies  
$\mathrm{ker}(\alpha)=0$ and therefore the exactness of $\mathbf C_{p,q}$. 
    \item 
If $(p,q)=(1,0)$, then $\mathbf C_{p,q}=(\mathbf k[f_0]\oplus V_{0,0}\stackrel{\alpha}{\to} 
V_{0,0}\oplus V_{1,0}\oplus V_{0,0}
\stackrel{\beta}{\to} V_{1,0})$, where $\alpha,\beta$ are as in \eqref{eq:alpha:0807}, \eqref{eq:beta:0807}. 
Let $(\mathbf f,v,\mathbf e)\in\mathrm{ker}(\beta)$. Then 
\begin{equation}\label{temp:eq:0807}
\varphi(v)=\mathbf f e_1-e_1 \mathbf e.     
\end{equation}
Applying $m$ to this equality and using Lem. \ref{lem:ressource:0707}(c), one obtains 
$m(\mathbf f e_1)=m(e_1 \mathbf e)$. Then $m(e_1 \mathbf e)=e_1\cdot m(\mathbf e)\in e_1\cdot \mathbf k[e_0]$, 
where the equality follows from Lem. \ref{lem:ressource:0707}(e) and the relation follows from Lem. \ref{lem:ressource:0707}(h), 
therefore $m(\mathbf f e_1)\in e_1\cdot \mathbf k[e_0]$. Lem. \ref{lem:ressource:0707}(i) then implies 
$\mathbf f\in\mathbf k[f_0]$. Let $P(u,v)\in\mathbf k[u,v]$ and $Q(v)\in\mathbf k[v]$ be such that 
$\mathbf e=P(e_0,f_0)$ and $\mathbf f=Q(f_0)$. Then 
$m(\mathbf f e_1)=m(Q(f_0) e_1)=e_1\cdot Q(e_0)$, and 
$m(e_1 \mathbf e)=m(e_1\cdot P(e_0,f_0))=e_1\cdot P(e_0,e_0)$. Then 
$m(\mathbf f e_1)=m(e_1 \mathbf e)$ implies $e_1\cdot Q(e_0)=e_1\cdot P(e_0,e_0)$, 
which since $\mathcal V$ is a domain and since the algebra morphism $\mathbf k[v]\to \mathcal V$ defined by 
$v\mapsto e_0$ is injective, implies $Q(v)=P(v,v)$. Let $R(u,v):=(P(u,v)-P(v,v))/(u-v)$. 
Then 
\begin{equation}\label{next:temp:0807}
\varphi(R(e_0,f_0))=(e_0-f_0)\cdot R(e_0,f_0)=P(e_0,f_0)-P(f_0,f_0)=P(e_0,f_0)-Q(f_0),     
\end{equation}
where the first equality follows from the commutation of $e_0$ and $f_0$. The right-hand side of \eqref{temp:eq:0807} is 
\begin{align*}
&\mathbf f e_1-e_1 \mathbf e=Q(f_0) e_1-e_1\cdot P(e_0,f_0)=e_1
(Q(f_0)-P(e_0,f_0))
\\ & =-e_1\cdot \varphi(R(e_0,f_0))=\varphi(-e_1\cdot R(e_0,f_0)), 
\end{align*}
where the second equality follows from the commutation of $e_1$ and $f_0$, the third equality follows from \eqref{next:temp:0807}, 
and the last equality follows from Lem. \ref{lem:ressource:0707}(f). Then  \eqref{temp:eq:0807}
implies $\varphi(v)=\varphi(-e_1\cdot R(e_0,f_0))$, which by Lem. \ref{lem:technical:0801}(a) implies 
$v=-e_1\cdot R(e_0,f_0)$. 
Recall that $\mathbf f=Q(f_0)$ and $\mathbf e=P(e_0,f_0)=(e_0-f_0) R(e_0,f_0)+Q(f_0)$, 
then 
$$
(\mathbf f,v,\mathbf e)=(0,-e_1\cdot R(e_0,f_0),(e_0-f_0)R(e_0,f_0))+(Q(f_0),0,Q(f_0))
=\alpha(Q,-R(e_0,f_0))
$$
therefore $(\mathbf f,v,\mathbf e)\in\mathrm{im}(\alpha)$. Therefore $\mathrm{im}(\alpha)=\mathrm{ker}(\beta)$, which 
implies that $\mathbf C_{p,q}$ is exact. 
\end{itemize}
It follows that for any $p,q\geq0$, the complex $\mathbf C_{p,q}$ is exact. 
Since the sequence \eqref{exact:seq:abv} is the direct sum of these complexes, it is also exact. 

Beside the above bidegree, $V$ (resp. $\mathbf k[t]$) is also graded
for the total degree, for which $e_0,e_1,f_0,f_1$ (resp. $t$) 
all have degree 1, and for $x\in\{e,f\}$, $\mathrm{C}_V(x_0)$ is a graded subalgebra of $V$; for $X$ a $\mathbf k$-module 
which is graded for the total degree and $n\geq0$, denote by $X_n$ the part of $X$ of total degree $n$, so that 
$X=\oplus_{n\geq0}X_n$. 

For $n\geq0$, let $A_n:=\mathbf k[t]_{n-1}\oplus\mathrm{C}_V(e_0)_{n-2}$, $B_n:=\mathrm{C}_V(f_0)_{n-1}\oplus V_{n-1} 
\oplus\mathrm{C}_V(e_0)_{n-1}$. This defines gradings on $A$ and $B$, namely  $A=\oplus_{n\geq 0}A_n$, $B=\oplus_{n\geq 0}B_n$.
By convention $X_n=0$ for $n<0$ and $X$ is equal to $\mathbf k[t]$ or $\mathrm{C}_V(x_0)$ for $x\in\{e,f\}$. 

The maps $\alpha$ and $\beta$ are then homogeneous for these gradings, therefore \eqref{exact:seq:abv} splits up as a direct sum of 
complexes $\oplus_{n\geq 0}\mathbf C_n$, where $\mathbf C_n:=(A_n\stackrel{\alpha_n}{\to}B_n\stackrel{\beta_n}{\to}V_n)$. 
The exactness of \eqref{exact:seq:abv} implies the exactness of the complex $\mathbf C_n$ for each $n\geq0$, which in its turn implies 
the acyclicity of the complete direct sum of these complexes, which is the complex 
\begin{equation}\label{REF:0907}
\hat A\stackrel{\hat\alpha}{\to}\hat B\stackrel{\hat\beta}{\to}\hat V,     
\end{equation}
where $\hat A=\mathbf k[[t]]\oplus \mathrm{C}_{\hat V}(e_0)$ and $\hat B=\mathrm{C}_{\hat V}(f_0)\oplus \hat V\oplus 
\mathrm{C}_{\hat V}(e_0)$ and $\hat\alpha,\hat\beta$ are the degree completions of 
$\alpha,\beta$. 
Lem. \ref{pbms:1:2} follows from the acyclicity of \eqref{REF:0907}.
\end{proof}

\subsection{Computation of $\mathrm{C}_2(\overline\rho_0)$, $\mathrm{C}_3(\rho_0)$ and $\mathrm{C}_{21}(\rho_0)$}\label{present:section}

Define the map 
\begin{equation}\label{map:overlineM:Pi:m}
\mathbf k[[f_0]]\times\mathrm C_{\hat V}(e_0)\ni(\Pi,C)\mapsto \overline X(\Pi,C)
:=\Pi I_2+\begin{pmatrix}
   e_0-f_0\\ e_1 
\end{pmatrix}C\begin{pmatrix}
    1&0
\end{pmatrix}
\in M_2\hat V. 
\end{equation}

\begin{lem}
One has 
\begin{equation}\label{lem:15:8:0908:debut}
\mathrm{C}_2(\overline\rho_0)=\{\overline X(\Pi,C)|\Pi\in\mathbf k[[f_0]],C\in \mathrm C_{\hat V}(e_0)\}.   
\end{equation}
\end{lem}

\begin{proof}
For $m=\begin{pmatrix}
a&b\\c&d\end{pmatrix}\in M_2\hat V$, the (1,2) element of $[m,\overline\rho_0]$ is 
$e_0b-bf_0$. Therefore  $m\in\mathrm{C}_2(\overline\rho_0)$ implies 
$e_0b=bf_0$, which by Lem. \ref{lem:technical:0801}(a) implies $b=0$, so that $m$ is lower-triangular. 
If now $m=\begin{pmatrix}
a&0\\c&d\end{pmatrix}$ is a lower-triangular matrix, $m\in \mathrm{C}_2(\overline\rho_0)$ is equivalent to the system 
\begin{equation}\label{system:comm:sigma0}
a\in \mathrm{C}_{\hat V}(e_0),\quad d\in \mathrm{C}_{\hat V}(f_0),\quad 
f_0c-ce_0=-e_1a+de_1. 
\end{equation}    
By Lem. \ref{pbms:1:2}, this system is equivalent to the existence of $\Pi\in\mathbf k[[f_0]]$
and $C\in \mathrm C_{\hat V}(e_0)$, such that $(a,c,d)=(\Pi,e_1C,\Pi+(e_0-f_0)C)$.  
This implies the result. 
\end{proof}

\begin{defn}
(a) Set $A:=\begin{pmatrix} e_0-f_0&1\\e_1&0\\0&1\end{pmatrix}\in M_{3,2}\hat V$ and 
$B:=\begin{pmatrix}1&0&0\\0&0&1\end{pmatrix}\in M_{2,3}\hat V$.

(b) Define the map 
\begin{equation}\label{map:M:Pi:m}
\mathbf k[[f_0]]\times M_2\mathrm C_{\hat V}(e_0)\ni(\Pi,m)\mapsto X(\Pi,m)\in M_3\hat V
\end{equation}
by 
\begin{equation}\label{def:map:M:Pi:m}
X(\Pi,m):=\Pi I_3+AmB
=\Pi I_3+\begin{pmatrix}(e_0-f_0)a+b&0&(e_0-f_0)c+d\\e_1a&0&e_1c\\b&0&d\end{pmatrix} 
\end{equation} 
if $m=\begin{pmatrix}a&c\\b&d\end{pmatrix}\in M_2\mathrm{C}_{\hat V}(e_0)$. 
\end{defn}

It follows from the injectivity of $x\mapsto e_1x$ that the map $(\Pi,m)\mapsto X(\Pi,m)$ is injective. 

\begin{lem}\label{lem:comp:comm:1108}
One has 
\begin{equation}\label{value:c:m3:rho0}
    \mathrm C_{3}(\rho_0)=\{X(\Pi,m)|\Pi\in\mathbf k[[f_0]],m\in M_2(\mathrm C_{\hat V}(e_0))\}
\end{equation}
and
\begin{equation}\label{value:c:a:rho0}
    \mathrm C_{21}(\rho_0)=\{X(\Pi,m)|\Pi \in\mathbf k[[f_0]],m\in T_{11}\mathrm C_{\hat V}(e_0)\}, 
\end{equation}
where $T_{11}C_{\hat V}(e_0):=\{\begin{pmatrix} a&b\\0&d
    \end{pmatrix}|a,b,d\in C_{\hat V}(e_0)\}$. 
\end{lem}

\begin{proof}
One checks that that the entries of $\rho_0$ commute with $f_0$, which implies that for any $\Pi\in\mathbf k[[f_0]]$, 
$\Pi I_3$ commutes with $\rho_0$. One also checks the equalities $\rho_0A=Ae_0$ and $e_0B=B\rho_0$. Then 
for any $m\in M_2 \mathrm C_{\hat V}(e_0)$, one has 
$$
\rho_0\cdot AmB=Ae_0mB=Ame_0B=AmB\cdot \rho_0, 
$$
where the middle equality follows from $m\in M_2 \mathrm C_{\hat V}(e_0)$, therefore $AmB$ commutes with  
$\rho_0$. All this proves the inclusion 
\begin{equation}\label{incl:one:dir}
\{X(\Pi,m)|\Pi\in\mathbf k[[f_0]],m\in M_2\mathrm C_{\hat V}(e_0)\}\subset \mathrm{C}_3(\rho_0). 
\end{equation}
We now prove the opposite inclusion. For $M\in M_3\hat V$, we will denote by $M_{11}\in M_2\hat V$, 
$M_{12}\in M_{2,1}\hat V$, $M_{21}\in M_{1,2}\hat V$, $M_{22}\in\hat V$ the matrices such that 
\begin{equation}\label{Mij:parts}
M=\begin{pmatrix}M_{11}&M_{12}\\M_{21}&M_{22}\end{pmatrix}
\end{equation} ($M_{ij}$ will be referred to as the $(i,j)$ part of $M$). 
We also set $\kappa:=\begin{pmatrix}
0\\-e_1\end{pmatrix}$, so that $\rho_0=\begin{pmatrix}\overline\rho_0&\kappa\\0&e_0\end{pmatrix}$.   

Let 
\begin{equation}\label{belonging:M}
M\in \mathrm{C}_3(\rho_0). 
\end{equation}
 The (2,1) part of the equality
$\rho_0M=M\rho_0$ gives the equality $e_0M_{21}=M_{21}\overline\rho_0$ (in $M_{1,2}\hat V$). 
Let $b,b'\in \hat V$ be the elements such that $M_{21}=\begin{pmatrix}b&b'\end{pmatrix}$, then 
this equality is equivalent to the system 
$$
e_0b=be_0+b'e_1,\quad e_0b'=b'f_0. 
$$
By Lem. \ref{lem:technical:0801}(a), the second equality implies $b'=0$, and the first equality then implies 
$b\in C_{\hat V}(e_0)$, therefore 
\begin{equation}\label{temp:a:0908}
M_{21}=\begin{pmatrix}b&0\end{pmatrix}
\end{equation}
with $b\in C_{\hat V}(e_0)$. 

One has 
$X(0,\begin{pmatrix}0\\1\end{pmatrix}b\begin{pmatrix}1&0\end{pmatrix})=
\begin{pmatrix}0\\0\\1\end{pmatrix}b\begin{pmatrix}1&0&0\end{pmatrix}$, therefore 
\begin{equation}\label{temp:b:0908}
X(0,\begin{pmatrix}0\\1\end{pmatrix}b\begin{pmatrix}1&0\end{pmatrix})_{21}=
\begin{pmatrix}b&0\end{pmatrix}. 
\end{equation} 

Set then 
\begin{equation}\label{def:M'}
M':=M-X(0,\begin{pmatrix}0\\1\end{pmatrix}b\begin{pmatrix}1&0\end{pmatrix}).
\end{equation} 
It follows from \eqref{temp:a:0908} and \eqref{temp:b:0908} that $M'_{21}=0$, therefore 
$M'=\begin{pmatrix}M'_{11}&M'_{12}\\0&M'_{22}\end{pmatrix}$. 
Since $b\in C_{\hat V}(e_0)$, one has $\begin{pmatrix}0\\1\end{pmatrix}b\begin{pmatrix}1&0\end{pmatrix}\in 
M_2(\mathrm C_{\hat V}(e_0))$, therefore by \eqref{incl:one:dir} one derives  
$M(0,\begin{pmatrix}0\\1\end{pmatrix}b\begin{pmatrix}1&0\end{pmatrix})\in 
\mathrm{C}_3(\rho_0)$, which by \eqref{belonging:M} implies 
\begin{equation}\label{belonging:M'}
M'\in \mathrm{C}_3(\rho_0). 
\end{equation}
This implies the equality $\rho_0M'=M'\rho_0$. Its (1,1) part gives the equality $\overline\rho_0M'_{11}=M'_{11}\overline\rho_0$, 
therefore $M'_{11}\in C_2(\overline\rho_0)$. By \eqref{lem:15:8:0908:debut}, this implies the 
existence of $\Pi\in\mathbf k[[f_0]]$ and $a\in \mathrm C_{\hat V}(e_0)$ such that 
\begin{equation}\label{temp:c:0908}
M'_{11}=\Pi I_2+\begin{pmatrix}   e_0-f_0\\ e_1 
\end{pmatrix}a\begin{pmatrix}    1&0
\end{pmatrix}. 
\end{equation}

Set then 
\begin{equation}\label{def:M''}
M'':=M'-X(\Pi,\begin{pmatrix}1\\0\end{pmatrix}a\begin{pmatrix}1&0\end{pmatrix}).
\end{equation}  
One has  
$X(\Pi,
\begin{pmatrix}1\\0\end{pmatrix}a\begin{pmatrix}1&0\end{pmatrix})=\Pi I_3+\begin{pmatrix}e_0-f_0\\e_1\\0\end{pmatrix}a
\begin{pmatrix}1&0&0\end{pmatrix}$, therefore 
$$
X(\Pi,
\begin{pmatrix}1\\0\end{pmatrix}a\begin{pmatrix}1&0\end{pmatrix})_{11}=
\Pi I_2+\begin{pmatrix}   e_0-f_0\\ e_1 
\end{pmatrix}a\begin{pmatrix}    1&0
\end{pmatrix},\quad 
X(\Pi,\begin{pmatrix}1\\0\end{pmatrix}a\begin{pmatrix}1&0\end{pmatrix})_{21}=0.
$$
which together with \eqref{temp:c:0908} implies $M''_{11}=M''_{21}=0$, therefore 
\begin{equation}\label{form:M''}
M''=\begin{pmatrix}0&M''_{12}\\0&M''_{22}\end{pmatrix}. 
\end{equation}
Since $a\in \mathrm C_{\hat V}(e_0)$, one has $\begin{pmatrix}1\\0\end{pmatrix}a\begin{pmatrix}1&0\end{pmatrix}\in 
M_2\mathrm C_{\hat V}(e_0)$, which together with $\Pi\in \mathbf k[[f_0]]$ and \eqref{incl:one:dir} implies 
$X(\Pi,\begin{pmatrix}1\\0\end{pmatrix}a\begin{pmatrix}1&0\end{pmatrix})\in 
\mathrm{C}_3(\rho_0)$, which together with \eqref{belonging:M'} implies 
$M''\in \mathrm{C}_3(\rho_0)$. 
This implies the equality $\rho_0M''=M''\rho_0$, whose (2,2) and (1,2) parts are respectively
$$
e_0M''_{22}=M''_{22}e_0,\quad \overline\rho_0M''_{12}+\kappa M''_{22}=M''_{12}e_0. 
$$
The first equality is equivalent to 
\begin{equation}\label{belonging:M''22}
M''_{22}\in \mathrm C_{\hat V}(e_0), 
\end{equation}
and if $c,c'\in \hat V$ are such that 
\begin{equation}\label{explicit:M''12}
M''_{12}=\begin{pmatrix}c\\c'\end{pmatrix},
\end{equation} 
the second equality is written as the system
$$
e_0c=ce_0,\quad 
f_0c'-c'e_0=e_1(M''_{22}-c).  
$$
The first equation implies $c\in\mathrm C_{\hat V}(e_0)$, which together with \eqref{belonging:M''22} implies 
$M''_{22}-c\in\mathrm C_{\hat V}(e_0)$. Lem. \ref{pbms:1:2} then implies the existence of $\Pi'\in \mathbf k[[f_0]]$ and 
$C\in \mathrm{C}_{\hat V}(e_0)$, such that $(0,c',M''_{22}-c)=
(\Pi',e_1 C,\Pi'+(f_0-e_0)C)$. One necessarily has $\Pi'=0$, therefore 
\begin{equation}\label{values:c:c'}
c=M''_{22}+(e_0-f_0)C,\quad c'=e_1C. \end{equation} 
Since $(C,M''_{22})\in \mathrm C_{\hat V}(e_0)^2$, one has $\begin{pmatrix}0&C\\0&M''_{22}\end{pmatrix}
\in M_2\mathrm C_{\hat V}(e_0)$. By \eqref{incl:one:dir}, it follows that 
$X(0,\begin{pmatrix}0&C\\0&M''_{22}\end{pmatrix})\in \mathrm{C}_3(\rho_0)$. 
Equations \eqref{explicit:M''12} and \eqref{values:c:c'} then imply 
\begin{equation}\label{eq:M''}
M''=X(0,\begin{pmatrix}0&C\\0&M''_{22}\end{pmatrix}).
\end{equation} 
Then 
$$
M=M'+X(0,\begin{pmatrix}0&0\\b&0\end{pmatrix})
=M''+X(\Pi,\begin{pmatrix}a&0\\b&0\end{pmatrix})
=X(\Pi,\begin{pmatrix}a&C\\b&M''_{22}\end{pmatrix})$$
where the first (resp. second, third) equality follows from \eqref{def:M'} (resp. \eqref{def:M''}, 
\eqref{eq:M''}), and from the linearity if $(\Pi,m)\mapsto X(\Pi,m)$. This implies the opposite inclusion to  
\eqref{incl:one:dir}, and therefore the equality \eqref{value:c:m3:rho0}. 

Since the $(3,1)$ element of any matrix of $T_{21}\hat V$ is 0, and in view of \eqref{def:map:M:Pi:m}
one has $(X(\Pi,m)\in T_{21}\hat V)\implies(b=0)$. Conversely,  
since $b=0$ implies that the  $(3,1)$ and $(3,2)$ elements of 
$X(\Pi,m)$ are 0, one sees that $(b=0)\implies(X(\Pi,m)\in T_{21}\hat V)$. The equivalence $(X(\Pi,m)\in T_{21}\hat V)\iff(b=0)$, 
together with $\mathrm C_{21}(\rho_0)=\mathrm C_3(\rho_0)\cap T_{21}\hat V$ and \eqref{value:c:m3:rho0}, implies 
\eqref{value:c:a:rho0}.
\end{proof}

\subsection{Computation of $\mathrm C_2(\overline\rho_0)^\times$, $\mathrm{C}_3(\rho_0)^\times$ and $\mathrm{C}_{21}(\rho_0)^\times$}\label{sect:8:4:5786}

\begin{lem}\label{lem:8:3:departrennes}
With $(\Pi,C)\mapsto\overline X(\Pi,C)$ as in \eqref{map:overlineM:Pi:m}, one has 
$$
\mathrm C_2(\overline\rho_0)^\times
=\{\overline X(\Pi,C)|\Pi\in
\mathbf k[[f_0]]^\times,C\in \mathrm{C}_{\hat V}(e_0)\}.
$$ 
\end{lem}

\begin{proof}
One has 
\begin{align*}
\mathrm C_2(\overline\rho_0)^\times
=\mathrm C_2(\overline\rho_0)\cap\mathrm{GL}_2(\hat V)
=\{x|\exists (\Pi,C)\in \mathbf k[[f_0]]\times\mathrm{C}_{\hat V}(e_0), 
x=\overline X(\Pi,C)
\mathrm{\ and\ } x\in \mathrm{GL}_2(\hat V)\}. 
\end{align*}
where the first equality follows from Lem. \ref{invertible:commutant}, and
the second equality follows from \eqref{lem:15:8:0908:debut}. The statement then 
follows from the equivalences $(\overline X(\Pi,C)\in \mathrm{GL}_2\hat V)\iff 
(\epsilon(\overline X(\Pi,C))\in \mathrm{GL}_2\mathbf k) \iff
(\epsilon(\Pi)\in\mathbf k^\times)\iff(\Pi\in \mathbf k[[f_0]]^\times)$, where 
the first equivalence follows from \eqref{charact:gl3}, and the second equivalence 
follows from the identity $\epsilon(X(\Pi,C))=\epsilon(\Pi)I_2$.  
\end{proof}

\begin{lem}\label{lem:15:13:1408}
(a) One has 
$$
\mathrm{C}_{3}(\rho_0)^\times
=\{X(\Pi,\begin{pmatrix}a&c\\b&d\end{pmatrix})|
a,b,c,d\in\mathrm C_{\hat V}(e_0), \Pi\in \mathbf k[[f_0]]^\times, 
\quad \text{and} \quad \Pi+b+d\in\hat V^\times\}. 
$$

(b) One has 
$$
\mathrm{C}_{21}(\rho_0)^\times
=\{X(\Pi,\begin{pmatrix}a&c\\0&d\end{pmatrix})|
a,c,d\in\mathrm C_{\hat V}(e_0), \Pi\in \mathbf k[[f_0]]^\times, 
\quad \text{and} \quad \Pi+d\in\hat V^\times\}.
$$
\end{lem}

\begin{proof}
(a) One has 
\begin{align*}
&\mathrm{C}_3(\rho_0)^\times=\mathrm{C}_3(\rho_0)\cap \mathrm{GL}_3\hat V
=\{x\in \mathrm{C}_3(\rho_0)|\epsilon(x)\in \mathrm{GL}_3(\mathbf k)\}
\\&\nonumber =\{X(\Pi,m)|(\Pi,m)\in\mathbf k[[f_0]]\times M_2\mathrm{C}_{\hat V}(e_0),\epsilon(X(\Pi,m))\in \mathrm{GL}_3(\mathbf k)\}.
\end{align*}
where the first equality follows from Lem. \ref{invertible:commutant},  
the second equality follows from \eqref{charact:gl3}, and the third equality follows from 
\eqref{value:c:m3:rho0}, the map $(\Pi,m)\mapsto X(\Pi,m)$ being as in 
\eqref{def:map:M:Pi:m}. For $\Pi\in\mathbf k[[f_0]]$ and $m=\begin{pmatrix}a&c\\b&d\end{pmatrix}\in 
M_2\mathrm{C}_{\hat V}(e_0)$, one has 
$$
\epsilon(X(\Pi,m))=\begin{pmatrix}\epsilon(\Pi) +\epsilon(b)&0&\epsilon(d)\\0&\epsilon(\Pi) &0\\
\epsilon(b)&0&\epsilon(\Pi) +\epsilon(d)\end{pmatrix} 
$$ 
whose determinant is equal to $\epsilon(\Pi)^2\cdot \epsilon(\Pi+b+d)$. 
Then $(\epsilon(X(\Pi,m))\in\mathrm{GL}_3(\mathbf k))\iff (\epsilon(\Pi)^2\cdot \epsilon(\Pi+b+d)\in\mathbf k^\times)
\iff (\epsilon(\Pi)$ and $\epsilon(\Pi+b+d)\in\mathbf k^\times) \iff (\Pi\in\mathbf k[[f_0]]^\times$ and 
$\Pi+b+d\in\hat V^\times)$. This implies the result.

(b) One has 
\begin{align*}
&\mathrm{C}_{21}(\rho_0)^\times=\mathrm{C}_{21}(\rho_0)\cap (T_{21}\hat V)^\times
=\mathrm{C}_{21}(\rho_0)\cap \mathrm{GL}_3\hat V
\\&
=\{X(\Pi,m)|(\Pi,m)\in \mathbf k[[f_0]]\times T_{11}\mathrm C_{\hat V}(e_0)
\text{\ and\ }X(\Pi,m)\in\mathrm{GL}_3\hat V\} 
\\&
=\{X(\Pi,m)|(\Pi,m)\in \mathbf k[[f_0]]\times T_{11}\mathrm C_{\hat V}(e_0)
\text{\ and\ }\epsilon(X(\Pi,m))\in\mathrm{GL}_3(\mathbf k)\} 
\end{align*}
where $(\Pi,m)\mapsto X(\Pi,m)$ is as in \eqref{map:M:Pi:m}, where the first equality follows from 
Lem. \ref{invertible:commutant}, the second equality follows from 
 \eqref{invertibles:t:21} and $\mathrm{C}_{21}(\rho_0)\subset T_{21}\hat V$, the third 
 equality follows from \eqref{value:c:a:rho0}, and the fourth equality follows from \eqref{charact:gl3}. 
 With $m=\begin{pmatrix}a&c\\0&d\end{pmatrix}\in
T_{11}\mathrm C_{\hat V}(e_0)$, one has $\epsilon(X(\Pi,m))=\begin{pmatrix}\epsilon(\Pi)&0&\epsilon(d)\\0&\epsilon(\Pi)&0\\
0&0&\epsilon(\Pi+d)\end{pmatrix}$, whose determinant is $\epsilon(\Pi+d)\epsilon(\Pi)^2$, hence 
$(\epsilon(X(\Pi,m))\in\mathbf k^\times)\iff(\epsilon(\Pi)\in\mathbf k^\times$ and $\epsilon(\Pi+d)\in\mathbf k^\times)
\iff(\Pi\in\mathbf k[[f_0]]^\times$ and $\Pi+d\in\hat V^\times)$, which implies the result. 
\end{proof}

\subsection{The equality $\mathrm C_3(\rho_1)^\times\cdot \mathrm C_3(\rho_0)^\times\cap (T_{21}\hat V)^\times
=\mathrm C_{21}(\rho_1)^\times\cdot \mathrm C_{21}(\rho_0)^\times$}\label{sect:8:5:5786}

\begin{lem}
One has 
\begin{equation}\label{toto:1408}
\mathrm C_3(\rho_1)^\times\cdot \mathrm C_3(\rho_0)^\times\cap (T_{21}\hat V)^\times
=\mathrm C_{21}(\rho_1)^\times\cdot \mathrm C_{21}(\rho_0)^\times
\end{equation}
(equality of subsets of $\mathrm{GL}_3\hat V$). 
 \end{lem}

\begin{proof}
One has $\mathrm{C}_{21}(\rho_i)^\times\subset \mathrm{C}_3(\rho_i)^\times$
($i=0,1$), therefore $\mathrm{C}_{21}(\rho_1)^\times\cdot \mathrm{C}_{21}(\rho_0)^\times\subset 
\mathrm{C}_3(\rho_1)^\times\cdot \mathrm{C}_3(\rho_0)^\times$. Moreover, 
one has $\mathrm{C}_{21}(\rho_i)^\times\subset (T_{21}\hat V)^\times$ ($i=0,1$), and $(T_{21}\hat V)^\times$ 
is a group, therefore 
$\mathrm{C}_{21}(\rho_1)^\times\cdot \mathrm{C}_{21}(\rho_0)^\times\subset (T_{21}\hat V)^\times$. All this 
implies that the right-hand side of \eqref{toto:1408} is contained in its left-hand side. 

Let $x$ belong to  $\mathrm{C}_3(\rho_1)^\times\cdot \mathrm{C}_3(\rho_0)^\times$. 
Then there exist $g_i\in \mathrm{C}_3(\rho_i)^\times$ ($i=0,1$) and $x=g_1\cdot g_0$. 

By Lem. \ref{lem:15:13:1408}(a), there exist $\Pi \in\mathbf k[[f_0]]^\times$ and 
$\begin{pmatrix}a&c\\b&d\end{pmatrix}\in M_2\mathrm C_{\hat V}(e_0)$ with 
$\Pi+b+d\in\hat V^\times$, such that 
$$
g_0=X(\Pi,\begin{pmatrix}a&c\\b&d\end{pmatrix})
$$
and by Lem. \ref{lem:15:12:BIS}(a), there exist $\phi\in\mathbf k[[e_1,f_1]]^\times$ and 
$\begin{pmatrix}\alpha&\gamma\\\beta&\delta\end{pmatrix}\in M_2\hat V$, with $\phi+\delta\in\hat V^\times$,  
such that 
$$
g_1=M(\phi,\begin{pmatrix}\alpha&\gamma\\\beta&\delta\end{pmatrix}).
$$
For $\alpha\in\{0,1,x\}$ and $i,j\in\{1,2\}$, let $M_{ij}^\alpha$ be the $(i,j)$ part of $g_\alpha$ if $\alpha\in\{0,1\}$
and of $x$ if $\alpha=x$ (see \eqref{Mij:parts}). Then 
\begin{equation}\label{M21:x}
M_{21}^x=M_{22}^1M_{21}^0+M_{21}^1M_{11}^0. 
\end{equation}
One computes 
$$
 M_{22}^1=\phi(e_1,f_1)+\delta,\quad
 M_{21}^0=\begin{pmatrix}b&0\end{pmatrix},\quad 
M_{21}^1=\begin{pmatrix}\beta&\beta\end{pmatrix}, \quad 
M_{11}^0=\Pi I_2+\begin{pmatrix}(e_0-f_0)a+b\\e_1a\end{pmatrix}\begin{pmatrix}1&0\end{pmatrix}. 
$$
By \eqref{M21:x}, one obtains 
$$
M_{21}^x=\big((\phi(e_1,f_1)+\delta)b+\beta((e_0+e_1-f_0)a+b)\big)\begin{pmatrix}1&0\end{pmatrix}+\beta\Pi\begin{pmatrix}1&1\end{pmatrix}
$$
Assume also that $x\in (T_{21}\hat V)^\times$. Then $x\in T_{21}\hat V$, therefore $M_{21}^x=0$, therefore 
$$
\beta\Pi =0= (\phi(e_1,f_1)+\delta)b+\beta((e_0+e_1-f_0)a+b). 
$$
Since $\Pi$ is invertible, the first equality implies $\beta=0$, which by Lem. \ref{lem:15:13:1408}(b) implies 
$g_1\in \mathrm C_{21}(\rho_1)^\times$. The second equality then implies 
$(\phi(e_1,f_1)+\delta)b=0$, which since $\phi(e_1,f_1)+\delta\in \hat V^\times$ implies $b=0$, which by Lem. \ref{lem:c:21:rho:1} implies 
$g_0\in \mathrm C_{21}(\rho_0)^\times$. Therefore $x\in \mathrm C_{21}(\rho_1)^\times\cdot 
\mathrm C_{21}(\rho_0)^\times$. All this implies that the right-hand side of \eqref{toto:1408} is contained in 
its left-hand side. 
\end{proof}

\subsection{Local injectivity of the morphism (B)}\label{sect:8:6:5786}

\begin{defn}
    Define $\mathrm{Hom}^{1,(0)}_{\mathcal C\operatorname{-alg}}(\hat{\mathcal V},M_3\hat V)$ to be the subset of 
    $\mathrm{Hom}^1_{\mathcal C\operatorname{-alg}}(\hat{\mathcal V},M_3\hat V)$ of all morphisms $\rho : \hat{\mathcal V}\to 
    M_3\hat V$, such that $\rho(e_0)$ is $\mathrm{GL}_3\hat V$-conjugate to $\rho_0$.  
\end{defn}

\begin{lem}\label{lem:7:13:toto}
    (a) The image of the canonical map $\mathrm{Hom}^{1,(0)}_{\mathcal C\operatorname{-alg}}(\hat{\mathcal V},T_{21}\hat V)\to 
    \mathrm{Hom}^1_{\mathcal C\operatorname{-alg}}(\hat{\mathcal V},M_3\hat V)$ is contained in 
    $\mathrm{Hom}^{1,(0)}_{\mathcal C\operatorname{-alg}}(\hat{\mathcal V},M_3\hat V)$. 

    (b) The action $\bullet$ of $\mathrm{GL}_3\hat V$ on $\mathrm{Hom}^1_{\mathcal C\operatorname{-alg}}(\hat{\mathcal V},M_3\hat V)$
    (see Lem. \ref{lem:520:2212:FIRST}(a)) preserves the subset $\mathrm{Hom}^{1,(0)}_{\mathcal C\operatorname{-alg}}(\hat{\mathcal V},M_3\hat V)$. 
\end{lem}

 \begin{proof}
(a) follows from the definition of $\mathrm{Hom}^{1,(0)}_{\mathcal C\operatorname{-alg}}(\hat{\mathcal V},T_{21}\hat V)$ (see Def. \ref{def:6:4:2912}(f)). 
(b)  If $P\in \mathrm{GL}_3\hat V$ and $\rho\in \mathrm{Hom}^1_{\mathcal C\operatorname{-alg}}(\hat{\mathcal V},M_3\hat V)$, then 
$(P\bullet\rho)(e_0)=P\rho(e_0)P^{-1}=(P\alpha)\rho_0(P\alpha)^{-1}$ 
where $\alpha\in \mathrm{GL}_3\hat V$ is such that $\rho(e_0)=\alpha\rho_0\alpha^{-1}$, therefore $P\bullet\rho\in \mathrm{Hom}^1_{\mathcal C\operatorname{-alg}}(\hat{\mathcal V},M_3\hat V)$. 
\end{proof}

\begin{lem}\label{lem:7:14:toto}
    (a) The map $\mathrm{GL}_3\hat V\to \mathrm{Hom}^{1,(0)}_{\mathcal C\operatorname{-alg}}(\hat{\mathcal V},M_3\hat V)$ given by $g\mapsto \rho_g$, 
    where $\rho_g : \hat{\mathcal V}\to M_3\hat V$ is the algebra morphism defined by $e_1\mapsto \rho_1$, $e_0\mapsto 
    g \rho_0g^{-1}$, 
    induces bijections
    $$
 \mathrm{GL}_3\hat V/\mathrm C_3(\rho_0)^\times\to   \mathrm{Hom}^{1,(0)}_{\mathcal C\operatorname{-alg}}(\hat{\mathcal V},M_3\hat V)  
    $$
    and
    $$
    \mathrm C_3(\rho_1)^\times\backslash\mathrm{GL}_3\hat V/\mathrm C_3(\rho_0)^\times\to 
    \mathrm C_3(\rho_1)^\times\backslash\mathrm{Hom}^{1,(0)}_{\mathcal C\operatorname{-alg}}(\hat{\mathcal V},M_3\hat V), 
    $$ the latter taking the class 
    $\mathrm C_3(\rho_1)^\times\cdot\mathrm C_3(\rho_0)^\times$ of $I_3$ to the class $\mathrm C_3(\rho_1)^\times\bullet 
    \rho_{\mathrm{DT}}
    $ of $\rho_{\mathrm{DT}}$. 

    (b)  The map $g\mapsto \rho_g$ from (a) induces a map 
    $(T_{21}\hat V)^\times\to \mathrm{Hom}^{1,(0)}_{\mathcal C\operatorname{-alg}}(\hat{\mathcal V},T_{21}\hat V)$, which induces bijections
    $$
    (T_{21}\hat V)^\times/\mathrm C_{21}(\rho_0)^\times\to  \mathrm{Hom}^{1,(0)}_{\mathcal C\operatorname{-alg}}(\hat{\mathcal V},T_{21}\hat V)
    $$
    and
    $$
    \mathrm C_{21}(\rho_1)^\times\backslash(T_{21}\hat V)^\times/\mathrm C_{21}(\rho_0)^\times
    \to  \mathrm \mathrm C_{21}(\rho_1)^\times\backslash\mathrm{Hom}^{1,(0)}_{\mathcal C\operatorname{-alg}}(\hat{\mathcal V},T_{21}\hat V),
    $$
    the latter taking the class 
    $\mathrm C_{21}(\rho_1)^\times\cdot\mathrm C_{21}(\rho_0)^\times$ of $I_3$ to the class 
    $\mathrm C_{21}(\rho_1)^\times\bullet \rho_{\mathrm{DT}}
    $ of $\rho_{\mathrm{DT}}$.

    (c) These bijections build up a commutative diagram of pointed sets
    $$
    \xymatrix{(\mathrm C_{21}(\rho_1)^\times\backslash\mathrm{Hom}^{1,(0)}_{\mathcal C\operatorname{-alg}}(\hat{\mathcal V},T_{21}\hat V),\mathrm C_{21}(\rho_1)^\times\bullet\rho_{\mathrm{DT}}
    )\ar[r]\ar_\sim[d]&\ar^\sim[d]
    (\mathrm C_3(\rho_1)^\times\backslash\mathrm{Hom}^{1,(0)}_{\mathcal C\operatorname{-alg}}(\hat{\mathcal V},M_3\hat V),
    \mathrm C_3(\rho_1)^\times\bullet \rho_{\mathrm{DT}}
    )\\
    (\mathrm C_{21}(\rho_1)^\times\backslash(T_{21}\hat V)^\times/\mathrm C_{21}(\rho_0)^\times,\mathrm C_{21}(\rho_1)^\times\cdot\mathrm C_{21}(\rho_0)^\times)\ar[r]&
    (\mathrm C_3(\rho_1)^\times\backslash\mathrm{GL}_3\hat V/\mathrm C_3(\rho_0)^\times,\mathrm C_3(\rho_1)^\times\cdot\mathrm C_3(\rho_0)^\times)}
    $$
    where the top horizontal map is induced by composition with the canonical 
    injection $T_{21}\hat V\hookrightarrow M_3\hat V$ and the bottom horiziontal map is induced by the group inclusion $(T_{21}\hat V)^\times\hookrightarrow\mathrm{GL}_3\hat V$. 
\end{lem}

\begin{proof}
(a)  It follows from the definition of $\mathrm{Hom}^{1,(0)}_{\mathcal C\operatorname{-alg}}(\hat{\mathcal V},M_3\hat V)$ and from the fact that 
$\hat{\mathcal V}$ is freely generated by $e_0$ and $e_1$ that the map $\rho\mapsto \rho(e_0)$ sets up a bijection between 
this set and the $\mathrm{GL}_3\hat V$-conjugation class of $\rho_0$ in $M_3\hat V$. The map from $\mathrm{GL}_3\hat V$ to this conjugation class, taking $g$ to 
$g\rho_0g^{-1}$, induces a bijection from $\mathrm{GL}_3\hat V/\mathrm C_3(\rho_0)^\times$ to this conjugation class.  
The combination of these bijections is a bijection $\mathrm{GL}_3\hat V/\mathrm C_3(\rho_0)^\times\to \mathrm{Hom}^{1,(0)}_{\mathcal C\operatorname{-alg}}(\hat{\mathcal V},M_3\hat V)$, 
which implies the first statement. The composition of this map 
with the natural projection $\mathrm{GL}_3\hat V\to\mathrm{GL}_3\hat V/\mathrm C_3(\rho_0)^\times$ is then the map $g\mapsto \rho_g$. One checks that the map 
$\mathrm{GL}_3\hat V/\mathrm C_3(\rho_0)^\times\to \mathrm{Hom}^{1,(0)}_{\mathcal C\operatorname{-alg}}(\hat{\mathcal V},M_3\hat V)$ is $\mathrm C_3(\rho_1)^\times$-equivariant, 
which upon taking quotients with respect to this action implies the second statement. The last statement follows from the equality $\rho_{\mathrm{DT}}=\rho_{I_3}$. 
(b) The proof is similar to that of (a), replacing $M_3\hat V$, $\mathrm{GL}_3\hat V$, $\mathrm C_3(\rho_0)^\times$, $\mathrm C_3(\rho_1)^\times$
by $T_{21}\hat V$, $(T_{21}\hat V)^\times$, $\mathrm C_{21}(\rho_0)^\times$, $\mathrm C_{21}(\rho_1)^\times$. (c) follows from the commutativity of the diagram 
$$
\xymatrix{(T_{21}\hat V)^\times\ar[r]\ar[d]&\mathrm{GL}_3\hat V\ar[d]\\
\mathrm{Hom}^{1,(0)}_{\mathcal C\operatorname{-alg}}(\hat{\mathcal V},T_{21}\hat V)\ar[r]&\mathrm{Hom}^{1,(0)}_{\mathcal C\operatorname{-alg}}(\hat{\mathcal V},M_3\hat V)}
$$
which follows that the two composed maps to this diagram are given by $g\mapsto\rho_g$. 
\end{proof}

\begin{lem}\label{lem:7:15:toto}
(a) The morphism of pointed sets $(\mathrm C_{21}(\rho_1)^\times\backslash(T_{21}\hat V)^\times/\mathrm C_{21}(\rho_0)^\times,\mathrm C_{21}(\rho_1)^\times\cdot\mathrm C_{21}(\rho_0)^\times)\to
    (\mathrm C_3(\rho_1)^\times\backslash\mathrm{GL}_3\hat V/\mathrm C_3(\rho_0)^\times,\mathrm C_3(\rho_1)^\times\cdot\mathrm C_3(\rho_0)^\times)$ is induced by composition with the canonical 
    injection $T_{21}\hat V\hookrightarrow M_3\hat V$ is locally injective.  

(b) The morphism of pointed sets 
$$
(\mathrm C_{21}(\rho_1)^\times\backslash\mathrm{Hom}^{1,(0)}_{\mathcal C\operatorname{-alg}}(\hat{\mathcal V},T_{21}\hat V),\mathrm C_{21}(\rho_1)^\times\bullet \rho_{\mathrm{DT}}
)\to
(\mathrm C_3(\rho_1)^\times\backslash\mathrm{Hom}^{1,(0)}_{\mathcal C\operatorname{-alg}}(\hat{\mathcal V},M_3\hat V),
\mathrm C_3(\rho_1)^\times\bullet \rho_{\mathrm{DT}}
)
$$ 
induced by composition with the canonical injection $T_{21}\hat V\hookrightarrow M_3\hat V$ is locally injective. 
\end{lem}

\begin{proof}
(a) the preimage of the element $\mathrm C_3(\rho_1)^\times\cdot\mathrm C_3(\rho_0)^\times$ by the said map is the image of the 
preimage of the the same element by the map $(T_{21}\hat V)^\times\to \mathrm C_3(\rho_1)^\times\backslash\mathrm{GL}_3\hat V/\mathrm C_3(\rho_0)^\times$. 
This preimage is the subset $(T_{21}\hat V)^\times\cap \mathrm C_3(\rho_1)^\times\cdot\mathrm C_3(\rho_0)^\times$ of $(T_{21}\hat V)^\times$, 
which by Lem. \ref{toto:1408} is equal to the subset 
$\mathrm C_{21}(\rho_1)^\times\cdot\mathrm C_{21}(\rho_0)^\times$ of $(T_{21}\hat V)^\times$, and whose image is the element 
$\mathrm C_{21}(\rho_1)^\times\cdot\mathrm C_{21}(\rho_0)^\times$ of the source. (b) follows from (a) and from Lem. \ref{lem:7:14:toto}(c).  
\end{proof}

\begin{prop}\label{prop:B:loc:inj}
    The morphism (B) of pointed sets is locally injective. 
\end{prop}

\begin{proof}
It follows from Lem. \ref{lem:7:13:toto} that this morphism admits the decomposition
\begin{align*}
    & (\mathrm C_{21}(\rho_1)^\times\backslash\mathrm{Hom}^{1,(0)}_{\mathcal C\operatorname{-alg}}(\hat{\mathcal V},T_{21}\hat V),\mathrm C_{21}(\rho_1)^\times\bullet \rho_{\mathrm{DT}}
    )\to
(\mathrm C_3(\rho_1)^\times\backslash\mathrm{Hom}^{1,(0)}_{\mathcal C\operatorname{-alg}}(\hat{\mathcal V},M_3\hat V),
\mathrm C_3(\rho_1)^\times\bullet \rho_{\mathrm{DT}}
)
    \\& \to
(\mathrm C_3(\rho_1)^\times\backslash\mathrm{Hom}^1_{\mathcal C\operatorname{-alg}}(\hat{\mathcal V},M_3\hat V),
\mathrm C_3(\rho_1)^\times\bullet \rho_{\mathrm{DT}}
)
\end{align*}
The last morphism is obviously injective, therefore locally injective. The local injectivity of the first morphism follows from Lem. \ref{lem:7:15:toto}(b).   
The statement then follows form the fact that the composition of two locally injective morphisms of pointed sets is locally injective. 
\end{proof}

\newpage

\section{Local injectivity of the morphism (D)}\label{sec 9}

This section is devoted to the proof of the local injectivity of the map 
$$
(D) :  (\mathrm C_{21}^{(0)}(\rho_1)^\times\backslash\mathrm{Hom}^{1,((0)),\bullet}_{\mathcal C\operatorname{-alg}}
(\hat{\mathcal V},T_{21}\hat V),\mathrm C_{21}^{(0)}(\rho_1)^\times\bullet \rho_{\mathrm{DT}}
)\to 
(\mathrm C_2({\overline\rho}_1)^\times\backslash\mathrm{Hom}^{1,(0)}_{\mathcal C\operatorname{-alg}}(\hat{\mathcal V},
M_2\hat V),\mathrm C_2({\overline\rho}_1)^\times\bullet\overline\rho_{\mathrm{DT}}
$$
which is obtained in Prop. \ref{prop:D:loc:inj} (§\ref{sect:9:4}) as the result of the following steps.  

§\ref{sect:9:1} contains the proof of the surjectivity of two group morphisms (Lem. \ref{lem:c:rho1:surj} and 
Lem. \ref{lem:c:rho0:surj}).  In \S\ref{sect:9:2}, based on Lem. \ref{lem:c:rho0:surj}, 
one identifies the map $\mathrm{Hom}^{1,((0))}_{\mathcal C\operatorname{-alg}}
(\hat{\mathcal V},T_{21}\hat V)\to \mathrm{Hom}^{1,(0)}_{\mathcal C\operatorname{-alg}}
(\hat{\mathcal V},M_2\hat V)$ with a coset space morphism, and uses this to compute 
the preimage of $\overline\rho_{\mathrm{DT}}$ by this map (Lem. \ref{lem:9:8:toto}(b)). 
\S\ref{sect:9:3} gives algebraic arguments (from Lem. \ref{lem:prelim:rho} to Lem. \ref{lem:quotients})
which prove the inclusion\footnote{The analogue of this statement without the superscript 
$\bullet$ in the left-hand side is incorrect; this is the main motivation for the introduction of 
$\mathrm{Hom}^{\bullet}_{\mathcal C\operatorname{-alg}}
(\hat{\mathcal V},T_{21}\hat V)$} 
$$
    \{g\in U^{(0)}|g\cdot\rho_{\mathrm{DT}}\in \mathrm{Hom}^{\bullet}_{\mathcal C\operatorname{-alg}}
(\hat{\mathcal V},T_{21}\hat V)\}\subset (U^{(0)}\cap \mathrm{C}_{21}(\rho_1)^\times)\cdot(U^{(0)}\cap \mathrm{C}_{21}
(\rho_0)^\times)
$$
(Lem. \ref{lem:MAIN}) where $U^{(0)}$ is given by Def. \ref{def:U}. 
In \S\ref{sect:9:4}, one proves the equality 
$\mathrm{ker}(\mathrm C_{21}^{(0)}(\rho_1)^\times\to\mathrm C_2({\overline\rho}_1)^\times)
=U^{(0)}\cap \mathrm C_{21}^{(0)}(\rho_1)^\times$ (Lem. \ref{computation:kernel:of:group:morphism}), which 
together with Lem. \ref{lem:9:8:toto}(b) and Lem. \ref{lem:MAIN} implies the local injectivity of the
map 
$$
    (\mathrm{ker}(\mathrm C_{21}^{(0)}(\rho_1)^\times\to\mathrm C_2({\overline\rho}_1)^\times)\backslash\mathrm{Hom}^{1,((0)),\bullet}_{\mathcal C\operatorname{-alg}}
(\hat{\mathcal V},T_{21}\hat V),
\mathrm{ker}(\mathrm C_{21}^{(0)}(\rho_1)^\times\to\mathrm C_2({\overline\rho}_1)^\times)\bullet\rho_{\mathrm{DT}})\to 
(\mathrm{Hom}^{1,(0)}_{\mathcal C\operatorname{-alg}}
(\hat{\mathcal V},M_2\hat V),\overline\rho_{\mathrm{DT}}).  
$$ 
(Cor. \ref{loc:inj:U}). 


The combination of this statement with the surjectivity of the group morphism
$\mathrm{C}_{21}^{(0)}(\rho_1)^\times\to \mathrm C_2(\overline\rho_1)^\times$ (Lem. \ref{lem:c:rho1:surj})
and a general statement on local injectivity (Lem. \ref{lem:loc:inj:succ:quotients}) then 
leads to the local injectivity of $(D)$ (Prop. \ref{prop:D:loc:inj}).

\subsection{Surjectivity of the group morphisms $\mathrm{C}_{21}^{(0)}(\rho_1)^\times\to 
\mathrm C_2(\overline\rho_1)^\times$ and $\mathrm{C}_{21}(\rho_0)^\times\to \mathrm C_2(\overline\rho_0)^\times$}
\label{sect:9:1}

\begin{lem}\label{lem:c:rho1:surj}
    The group morphism $\mathrm{C}_{21}^{(0)}(\rho_1)^\times\to \mathrm C_2(\overline\rho_1)^\times$ associated to the composed algebra morphism 
        $\mathrm{C}_{21}^{(0)}(\rho_1)\subset \mathrm{C}_{21}(\rho_1)\to \mathrm{C}_2(\overline\rho_1)$, where the second morphism is as in
        Lem. \ref{lem:diag:algebras:2912}(b), is surjective. 
\end{lem}

\begin{proof}
By Lem. \ref{lem:c:21:rho:1}, one has $\mathrm{C}_{21}(\rho_1)
=\{M(\phi,m)|\phi\in\mathbf k[[u,v]],m\in T_{11}\hat V\}$, where $T_{11}\hat V$ is as in Lem. \ref{lem:15:12}. 
By \eqref{form:of:commutant:2912}, the map 
$\mathrm{C}_{21}(\rho_1)\to T_{21}\mathbf k$ is given by $M(\phi,\begin{pmatrix}
   \alpha & \gamma \\ 0& \delta 
\end{pmatrix})\mapsto \mathrm{diag}(\phi(0,0),\phi(0,0),\phi(0,0)+\epsilon(\delta))$. 
It follows from \eqref{value:c:a:rho0} that 
\begin{equation}\label{comp:C:21:rho0:0}
    \mathrm C_{21}(\rho_0)_0=\{\begin{pmatrix}
    \Pi&0&\Sigma-\Pi\\0&\Pi&0\\0&0&\Sigma
\end{pmatrix}|\Pi,\Sigma\in\mathbf k\}\subset T_{21}\mathbf k. 
\end{equation}
It follows that $\mathrm{C}^{(0)}_{21}(\rho_1)
=\{M(\phi,\begin{pmatrix}
   \alpha & \gamma \\ 0& \delta 
\end{pmatrix})|\phi\in\mathbf k[[u,v]],\alpha,\gamma,\delta\in\hat V,\epsilon(\delta)=0\}$, and therefore 
$\mathrm{C}^{(0)}_{21}(\rho_1)^\times
=\{M(\phi,\begin{pmatrix}
   \alpha & \gamma \\ 0& \delta 
\end{pmatrix})|\phi\in\mathbf k[[u,v]]^\times,\alpha,\gamma,\delta\in\hat V,\epsilon(\delta)=0\}$.  
On the other hand, by Lem. \ref{lem:comp:C2:rho1:times},  
one has $\mathrm{C}_2(\overline\rho_1)^\times=\{\overline M(\phi,v)|(\phi,v)\in\mathbf k[[u,v]]^\times\times\hat V\}$. 

One checks that the map $T_{21}\hat V\to M_2\hat V$, $x\mapsto \overline x$ is such that 
$M(\phi,m)\mapsto\overline M(\phi,\alpha)$ for any $\phi\in\mathbf k[[u,v]],
m=\begin{pmatrix}\alpha&\gamma\\0&\delta\end{pmatrix}\in T_{11}\hat V$. Then for any 
$(\phi,v)\in\mathbf k[[u,v]]^\times\times\hat V$, $M(\phi,\begin{pmatrix}v&0\\0&0\end{pmatrix})\in \mathrm C_{21}^{(0)}(\rho_1)^\times$
and the map $\mathrm{C}_{21}^{(0)}(\rho_1)^\times\to \mathrm C_2(\overline\rho_1)^\times$ is such that 
$M(\phi,\begin{pmatrix}v&0\\0&0\end{pmatrix})\mapsto\overline M(\phi,v)$, which implies the claimed surjectivity. 
\end{proof}

\begin{lem}\label{lem:c:rho0:surj}
    The group morphism $\mathrm{C}_{21}(\rho_0)^\times\to \mathrm C_2(\overline\rho_0)^\times$ arising from the specialization of Lem. \ref{lem:diag:algebras:2912}(a) to 
    the algebra morphism $T_{21}\hat V\to M_2\hat V$,  $x\mapsto \overline x$, and $a=\rho_0$, is surjective. 
\end{lem}

\begin{proof}
By Lem. \ref{lem:15:13:1408}, one has 
$$
\mathrm{C}_{21}(\rho_0)^\times
=\{X(\Pi,\begin{pmatrix}a&c\\0&d\end{pmatrix})|
a,c,d\in\mathrm C_{\hat V}(e_0), \Pi\in \mathbf k[[f_0]]^\times, 
\quad \text{and} \quad \Pi+d\in\hat V^\times\}.
$$
and by Lem. \ref{lem:8:3:departrennes},  
one has $\mathrm{C}_2(\overline\rho_0)^\times=\{\overline X(\Pi,C)|\Pi\in
\mathbf k[[f_0]]^\times,C\in \mathrm{C}_{\hat V}(e_0)\}$. 

One checks that the map $T_{21}\hat V\to M_2\hat V$, $x\mapsto \overline x$ is such that 
$X(\Pi,m)\mapsto\overline X(\Pi,a)$ for any $\Pi\in\mathbf k[[f_0]],
m=\begin{pmatrix}a&c\\0&d\end{pmatrix}\in T_{11}\mathrm C_{\hat V}(e_0)$. Then for any 
$(\Pi,C)\in\mathbf k[[f_0]]^\times\times\mathrm{C}_{\hat V}(e_0)$, $X(\Pi,\begin{pmatrix}
    C&0\\0&0
\end{pmatrix})\in \mathrm C_{21}(\overline\rho_0)^\times$
and the map $\mathrm{C}_{21}(\rho_0)^\times\to \mathrm C_2(\overline\rho_0)^\times$ is such that 
$X(\Pi,\begin{pmatrix}
    C&0\\0&0
\end{pmatrix})\mapsto\overline X(\Pi,C)$, which implies the claimed surjectivity. 
\end{proof}

\subsection{Identification of $\mathrm{Hom}^{1,((0))}_{\mathcal C\operatorname{-alg}}
(\hat{\mathcal V},T_{21}\hat V)\to \mathrm{Hom}^{1,(0)}_{\mathcal C\operatorname{-alg}}
(\hat{\mathcal V},M_2\hat V)$ with a coset space morphism}\label{sect:9:2}

\begin{defn}
Define $$T_{21}^{(0)}\hat V:=T_{21}\hat V\oplus_{(T_{21}\hat V)_0}\mathrm C_{21}(\rho_0)_0$$ 
to be the unital $\mathbf k$-subalgebra of $T_{21}\hat V$ obtained by applying the construction of Lem. \ref{lem:fibered:pdts}(a) to 
the diagram $T_{21}\hat V\to (T_{21}\hat V)_0\supset\mathrm C_{21}(\rho_0)_0$, where the first map is the projection 
$T_{21}\hat V=\hat\oplus_{n\geq0}(T_{21}\hat V)_n\to(T_{21}\hat V)_0$. 
\end{defn}

\begin{lem}
(a) The following equality of groups holds
$$
(T_{21}^{(0)}\hat V)^\times=(T_{21}\hat V)^\times\times_{(T_{21}\hat V)_0^\times}\mathrm C_{21}(\rho_0)_0^\times,  
$$
the right-hand side being relative to the diagram of groups $(T_{21}\hat V)^\times\to (T_{21}\hat V)_0^\times\supset\mathrm C_{21}(\rho_0)_0^\times$. 

(b) There is a group inclusion $\mathrm C_{21}(\rho_0)^\times\subset(T_{21}^{(0)}\hat V)^\times$. \end{lem}

\begin{proof}
    (a) follows from Lem. \ref{lem:fibered:pdts}(c). (b) The projection $T_{21}\hat V\to(T_{21}\hat V)_0$ maps $\mathrm C_{21}(\rho_0)$
    to $\mathrm C_{21}(\rho_0)_0$, which implies the inclusion of algebras $\mathrm C_{21}(\rho_0)\subset T_{21}^{(0)}\hat V$. This implies the 
    claimed inclusion by taking the groups of units of both sides. 
\end{proof}

\begin{lem}\label{lem:comm:diag:toto}
(a) For $g\in (T_{21}^{(0)}\hat V)^\times$ and $\rho\in \mathrm{Hom}^{1,((0))}_{\mathcal C\operatorname{-alg}}
(\hat{\mathcal V},T_{21}\hat V)$, define $g\cdot\rho$ to be the morphism $\hat{\mathcal V}\to T_{21}\hat V$ such that 
$e_1\mapsto \rho_1$ and $e_0\mapsto \mathrm{Ad}_g(\rho(e_0))$. Then $(g,\rho)\mapsto g\cdot\rho$ 
defines a transitive action of the group $(T_{21}^{(0)}\hat V)^\times$ on the set $\mathrm{Hom}^{1,((0))}_{\mathcal C\operatorname{-alg}}
(\hat{\mathcal V},T_{21}\hat V)$. The stabilizer of the element $\rho_{\mathrm{DT}}$ is $\mathrm C_{21}(\rho_0)^\times$, so that action on this element defines an 
isomorphism 
$$
(T_{21}^{(0)}\hat V)^\times/\mathrm C_{21}(\rho_0)^\times\to \mathrm{Hom}^{1,((0))}_{\mathcal C\operatorname{-alg}}
(\hat{\mathcal V},T_{21}\hat V)
$$
of $(T_{21}^{(0)}\hat V)^\times$-pointed sets, whose inverse takes $\rho$ to the class 
$g\mathrm C_{21}(\rho_0)^\times$, where $g\in (T_{21}^{(0)}\hat V)^\times$ is any element such that 
$\rho(e_0)=\mathrm{Ad}_g(\rho_0)$. 

(b)  For $h\in \mathrm{GL}_2\hat V$ and $\sigma\in \mathrm{Hom}^{1,(0)}_{\mathcal C\operatorname{-alg}}
(\hat{\mathcal V},M_2\hat V)$, define $h\cdot\sigma$ to be the morphism $\hat{\mathcal V}\to M_2\hat V$ such that 
$e_1\mapsto \overline\rho_1$ and $e_0\mapsto \mathrm{Ad}_h(\sigma(e_0))$. Then $(h,\sigma)\mapsto h\cdot\sigma$ 
defines a transitive action of the group $\mathrm{GL}_2\hat V$ on the set $\mathrm{Hom}^{1,(0)}_{\mathcal C\operatorname{-alg}}
(\hat{\mathcal V},M_2\hat V)$. The stabilizer of the element $\overline\rho_{\mathrm{DT}}$ is $\mathrm C_{2}(\overline\rho_0)^\times$, so that action on 
this element defines an isomorphism 
$$
\mathrm{GL}_2\hat V/\mathrm C_{2}(\overline\rho_0)^\times\to \mathrm{Hom}^{1,(0)}_{\mathcal C\operatorname{-alg}}
(\hat{\mathcal V},M_2\hat V)
$$
of $\mathrm{GL}_2\hat V$-pointed sets, whose inverse takes $\sigma$ to the class 
$g\mathrm C_{21}(\overline\rho_0)^\times$, where $g\in \mathrm{GL}_2\hat V$ is any element such that 
$\sigma(e_0)=\mathrm{Ad}_g(\overline\rho_0)$.  

(c) The map $\mathrm{Hom}^{1,((0))}_{\mathcal C\operatorname{-alg}}
(\hat{\mathcal V},T_{21}\hat V)\to \mathrm{Hom}^{1,(0)}_{\mathcal C\operatorname{-alg}}
(\hat{\mathcal V},M_2\hat V)$, the inverses of the isomorphisms from (a),(b) and the morphism of pointed
sets $(T_{21}^{(0)}\hat V)^\times/\mathrm C_{21}(\rho_0)^\times
\to\mathrm{GL}_2\hat V/\mathrm C_{2}(\overline\rho_0)^\times$ induced by the group morphism 
$(T_{21}^{(0)}\hat V)^\times\to \mathrm{GL}_2\hat V$, $x\mapsto \overline x$ fit in the 
following diagram of pointed sets 
$$\xymatrix{\mathrm{Hom}^{1,((0))}_{\mathcal C\operatorname{-alg}}
(\hat{\mathcal V},T_{21}\hat V)\ar[r]\ar[d]&\mathrm{Hom}^{1,(0)}_{\mathcal C\operatorname{-alg}}
(\hat{\mathcal V},M_2\hat V)\ar[d]\\
(T_{21}^{(0)}\hat V)^\times/\mathrm C_{21}(\rho_0)^\times\ar[r]&\mathrm{GL}_2\hat V/\mathrm C_{2}(\overline\rho_0)^\times}$$

(d) One has 
\begin{equation}\label{relation:actions}
\forall g\in \mathrm{C}_{21}^{(0)}(\rho_1)^\times,\forall\rho\in 
\mathrm{Hom}^{1,((0))}_{\mathcal C\operatorname{-alg}}
(\hat{\mathcal V},T_{21}\hat V),\quad 
g\cdot\rho=g\bullet\rho.     
\end{equation}
\end{lem}
 
\begin{proof}
(a) The map $\rho\mapsto \rho(e_0)$ sets up a bijection between $\mathrm{Hom}^{1,((0))}_{\mathcal C\operatorname{-alg}}
(\hat{\mathcal V},T_{21}\hat V)$ and $\{x\in T_{21}\hat V|x\equiv \rho_0$ mod $T_{21}F^2\hat V$ and $x$ is $(T_{21}\hat V)^\times$-conjugate 
to $\rho_0\}$, and $g\mapsto \mathrm{Ad}_g(\rho_0)$ sets up a bijection between $(T_{21}^{(0)}\hat V)^\times/\mathrm C_{21}(\rho_0)^\times$ and the latter set.  

(b) The map $\sigma\mapsto \sigma(e_0)$ sets up a bijection between $\mathrm{Hom}^{1,(0)}_{\mathcal C\operatorname{-alg}}
(\hat{\mathcal V},M_2\hat V)$ and $\{x\in M_2\hat V|x\equiv \overline\rho_0$ mod $M_2F^2\hat V$ and $x$ is 
$\mathrm{GL}_2\hat V$-conjugate 
to $\overline\rho_0\}$, and $h\mapsto \mathrm{Ad}_h(\overline\rho_0)$ sets up a bijection between 
$\mathrm{GL}_2\hat V/\mathrm C_{2}(\overline\rho_0)^\times$ and the latter set.  

(c) follows from the bijections used in the proofs of (a),(b) and from the identity 
$\overline{\mathrm{Ad}_g(\rho_0)}=\mathrm{Ad}_{\overline g}(\overline\rho_0)$ for 
$g\in (T_{21}^{(0)}\hat V)^\times$. 

(d) One has $g\bullet\rho=\mathrm{Ad}_g\circ\rho$ and $g\cdot \rho$ is given by $e_1\mapsto \rho_1$, $e_0\mapsto\mathrm{Ad}_g(\rho(e_0))$, therefore 
$(g\bullet\rho)(e_0)=\mathrm{Ad}_g\circ\rho(e_0)=(g\cdot\rho)(e_0)$ where the 
equalities follow from the definitions, and 
$(g\bullet\rho)(e_1)=\mathrm{Ad}_g\circ\rho(e_1)=\mathrm{Ad}_g(\rho_1)=\rho_1=(g\cdot\rho)(e_1)$,
where the first and last equalities follow from the definitions, the second equality follows from 
$\rho\in \mathrm{Hom}^{1,((0))}_{\mathcal C\operatorname{-alg}}
(\hat{\mathcal V},T_{21}\hat V)$, and the third equality follows from 
$g\in \mathrm{C}_{21}^{(0)}(\rho_1)^\times$. 
\end{proof}

\begin{defn}\label{def:U}
    Define $U:=\mathrm{ker}(T_{21}\hat V^\times\to\mathrm{GL}_2\hat V, t\mapsto \overline t)$ and $U^{(0)}:=U\cap T_{21}^{(0)}\hat V^\times$. 
\end{defn}

\begin{lem}\label{lem:9:8:toto}
(a) Let $\phi : G\to H$ be a group morphism and $G_0\subset G$, $H_0\subset H$ be subgroups, such that $\phi(G_0)=H_0$. Then $\phi$ induces a map 
coset spaces $G/G_0\to H/H_0$, and the preimage by this map of 
the coset $H_0$ is $\{kG_0|k\in\mathrm{ker}\phi\}$. 

(b) The preimage of $\overline\rho_{\mathrm{DT}}$ by the map $\mathrm{Hom}^{1,((0))}_{\mathcal C\operatorname{-alg}}
(\hat{\mathcal V},T_{21}\hat V)\to \mathrm{Hom}^{1,(0)}_{\mathcal C\operatorname{-alg}}
(\hat{\mathcal V},M_2\hat V)$ is equal to $U^{(0)}\cdot\rho_{\mathrm{DT}}$. 
\end{lem}

\begin{proof} (a) Let $\alpha$ belong to the fiber of $H_0$ and $g\in G$ be a representative of $\alpha$. Then 
$\phi(g)\in H_0$. Since $\phi : G_0\to H_0$ is surjective, there exists 
$g_0\in G_0$ such that $\phi(g_0)=\phi(g)$, therefore $\phi(gg_0^{-1})=1$, therefore 
$gg_0^{-1}\in \mathrm{ker}\phi$, therefore $\alpha\in\{kG_0|k\in\mathrm{ker}\phi\}$. Conversely, 
for $k\in\mathrm{ker}\phi$, the image of $kG_0$ is $\phi(k)H_0=H_0$. This proves the claim. 

(a) and Lem. \ref{lem:c:rho0:surj}, together with the equality $\mathrm{ker}((T_{21}^{(0)}\hat V)^\times
\to\mathrm{GL}_2\hat V)=(T_{21}^{(0)}\hat V)^\times\cap U=U^{(0)}$, 
imply that the preimage of the element $\mathrm C_{2}(\overline\rho_0)^\times$ by the map 
$(T_{21}^{(0)}\hat V)^\times/\mathrm C_{21}(\rho_0)^\times\to\mathrm{GL}_2\hat V/\mathrm C_{2}(\overline\rho_0)^\times$
is equal to the image of the map $U^{(0)}\to (T_{21}^{(0)}\hat V)^\times/\mathrm C_{21}(\rho_0)^\times$ induced by action on 
$\mathrm C_{21}(\rho_0)^\times$. Lem. \ref{lem:comm:diag:toto}(c), together with the $(T_{21}^{(0)}\hat V)^\times$-equivariance of the 
bijection $(T_{21}^{(0)}\hat V)^\times/\mathrm C_{21}(\rho_0)^\times\to \mathrm{Hom}^{1,((0))}_{\mathcal C\operatorname{-alg}}
(\hat{\mathcal V},T_{21}\hat V)$, and the fact that this bijection takes $\mathrm C_{21}(\rho_0)^\times$ to $\rho_{\mathrm{DT}}$, 
then implies the result. \end{proof}

\subsection{Inclusion of a fiber in a product of groups
}
\label{sect:9:3}

\begin{lem}\label{lem:prelim:rho}
    If $\rho\in \mathrm{Hom}^{1,\bullet}_{\mathcal C\operatorname{-alg}}(\hat{\mathcal V},T_{21}\hat V)$ (see Def. \ref{def:6:4:2912}(f)) and 
    $(r,C)\in \hat V\times M_{3,1}F^1\hat V$ are such that 
$$
\mathrm C_{21}(\rho(\hat{\mathcal V}))=\mathbf k1+C\cdot \mathrm C_{\hat V}(e_0)\cdot R_r\text{ and }R_r\cdot C\in e_0+f_\infty+F^2\hat V  
$$
where $R_r$ is as in \eqref{def:R:r}, then: 

(a) $r\in \hat V^\times$, 

(b) there exists $(s,t)\in\hat V\times F^1\hat V$, such that $C=\begin{pmatrix}
    f_1s\\e_1s\\t
\end{pmatrix}$

(c) $rt\in \mathrm C_{\hat V}(e_0)\cap (e_0+f_\infty+F^2\hat V)$. 
\end{lem}

\begin{proof}
Let $u,v,t\in F^1\hat V$ be such that $C=\begin{pmatrix}
    u\\v\\t
\end{pmatrix}$.  The equality $R_r\cdot C=rt$ implies $rt\in e_0+f_\infty+F^2\hat V$. It follows that 
$r_0t_1=e_0+f_\infty$, where $r_0\in\mathbf k$ and $t_1\in \mathbf ke_0\oplus\mathbf ke_1\oplus \mathbf kf_0\oplus\mathbf kf_1$ are the degree 0 and 1 components 
of $r$ and $t$. This relation implies $r_0(t_1|e_0)=1$, where $(t_1|e_0)$ is the coordinate of $t_1$ in the basis $(e_0,e_1,f_0,f_1)$, therefore 
$r_0\in\mathbf k^\times$, therefore $r\in\hat V^\times$. 

Then $CR_r\in \mathrm C_{21}(\rho(\hat{\mathcal V}))\subset \mathrm{C}_{21}(\rho_1)$, where the 
last inclusion follows from $\rho(e_1)=\rho_1$. By Lem. \ref{lem:15:12}, this implies the existence of $(\phi,m)\in \mathbf k[[u,v]]\times T_{11}\hat V$, 
such that $CR_r=M(\phi,m)$. Let $\alpha,\gamma,\delta\in\hat V$ be such that $m=\begin{pmatrix}
    \alpha&\gamma\\0&\delta
\end{pmatrix}$, then $CR_r=M(\phi,m)$ gives 
$$
\begin{pmatrix}
    0&0&ur \\0&0&vr\\0&0&tr
\end{pmatrix}
= \begin{pmatrix}
      f_1\alpha+\phi(e_1,f_1) & f_1\alpha&f_1\gamma \\ e_1\alpha & e_1\alpha+\phi(e_1,f_1)&e_1\gamma\\ 0& 0& \delta+\phi(e_1,f_1)
   \end{pmatrix}
$$
which implies $(e_1u-f_1v)r=0$. Since $r\in\hat V^\times$, this implies $e_1u=f_1v$, which by Lem. \ref{LEM1:0301:BIS} implies the existence of $s\in \hat V$ such that 
$(u,v)=(f_1s,e_1s)$, which implies the statement on the form of $C$. 

Since $rt=R_r\cdot C$, one has $rt\in e_0+f_\infty+F^2\hat V$. 
Let us prove that  $rt\in \mathrm C_{\hat V}(e_0)$. Since $\mathrm C_{21}(\rho(\hat{\mathcal V}))$
is an algebra and $CR_r\in \mathrm C_{21}(\rho(\hat{\mathcal V}))$, one has $(CR_r)^2\in \mathrm C_{21}(\rho(\hat{\mathcal V}))$, therefore
there exists $(\lambda,a)\in \mathbf k\times \mathrm C_{\hat V}(e_0)$ such that $(CR_r)^2=\lambda I_3+CaR_r$, i.e. 
$\lambda I_3+C(a-rt)R_r=0$. The (1,1) entry of this relation implies $\lambda=0$ and its (3,3) entry then implies 
$t(a-rt)r=0$, which by the invertibility of $r$, nonvanishing of $t$, and integrity of $\hat V$, implies $rt=a$, therefore
$rt\in  \mathrm C_{\hat V}(e_0)$. 
\end{proof}

\begin{lem}\label{lem:9:10:1345}
    (a) The map $u : (v_1,v_2,w)\mapsto \begin{pmatrix}
        1&0&v_1\\0&1&v_2\\0&0&w
    \end{pmatrix}$ defines a bijection $\hat V^2\times\hat V^\times\to U$  (see Def. \ref{def:U}).

    (b) One has $U^{(0)}=\{u(v_1,v_2,w)|(v_1,v_2,w)\in\hat V^2\times\hat V^\times, \epsilon(v_2)=0,\epsilon(v_1)=\epsilon(w)-1\}$ 
    (see Def. \ref{def:U}).  

    (c) One has $U^{(0)}\cap\mathrm C_{21}(\rho_1)^\times=\{u(f_1\gamma,e_1\gamma,w)|(\gamma,w)\in \hat V\times\hat V^\times, \epsilon(w)=1\}$.

    (d) One has $U^{(0)}\cap\mathrm C_{21}(\rho_0)^\times=\{u((e_0-f_0)c+w-1,e_1c,w)|
    (c,w)\in\mathrm{C}_{\hat V}(e_0)\times\mathrm{C}_{\hat V}(e_0)^\times\}$.

    (e) One has 
\begin{align}\label{comp:product:subgroups}
&    (U^{(0)}\cap\mathrm C_{21}(\rho_1)^\times)\cdot (U^{(0)}\cap\mathrm C_{21}(\rho_0)^\times)
\\&   \nonumber=\{u(f_1\gamma+(e_0+f_\infty)c+d-1,e_1\gamma,\delta d)|(\gamma,\delta)\in\hat V\times\hat V^\times, (c,d)\in \mathrm{C}_{\hat V}(e_0)
\times\mathrm{C}_{\hat V}(e_0)^\times,\epsilon(\delta)=1\}, 
 \end{align}
   where for $A,B$ subgroups of a group $G$, we set $A\cdot B:=\{ab|a\in A,b\in B\}$.
\end{lem}

\begin{proof}
(a) follows from  Lem. \ref{lem:basic:algebra:2412}. For further use, let us note the identity
\begin{equation}\label{pdt:in:U}
    u(v_1,v_2,w)u(v'_1,v'_2,w')=u(v_1w+v'_1,v_2w+v'_2,ww'). 
\end{equation}
  
It follows from \eqref{comp:C:21:rho0:0} that 
$$\mathrm C_{21}(\rho_0)_0^\times=\{\begin{pmatrix}
    \Pi&0&v\\0&\Pi&0\\0&0&\Sigma
\end{pmatrix}|\Pi,\Sigma\in\mathbf k^\times,v=\Sigma-\Pi\}\subset (T_{21}\mathbf k)^\times,  
$$ 
which implies (b). It follows from Lem. \ref{lem:c:21:rho:1} that 
$$
\mathrm C_{21}(\rho_1)^\times=\{M(\phi,\alpha,\gamma,\delta)|
\phi\in\mathbf k[[u,v]]^\times,\alpha,\gamma,\delta\in \hat V,\phi(e_1,f_1)+\delta\in\hat V^\times\}
\subset (T_{21}\hat V)^\times,  
$$ 
where 
$$
M(\phi,\alpha,\gamma,\delta):=\phi(e_1,f_1)I_3+\begin{pmatrix}
f_1\alpha&f_1\alpha&f_1\gamma\\e_1\alpha&e_1\alpha&e_1\gamma\\0&0&\delta
\end{pmatrix}. 
$$
If $(v_1,v_2,w)\in\hat V^2\times\hat V^\times$ and $(\phi,\alpha,\gamma,\delta)\in\mathbf k[[u,v]]^\times\times\hat V^3$ are such that 
$\epsilon(v_2)=0,\epsilon(v_1)=\epsilon(w)-1$, $\phi(e_1,f_1)+\delta\in\hat V^\times$, then the equality $u(v_1,v_2,w)=M(\phi,\alpha,\gamma,\delta)$
implies $\alpha=0$ (by inspection of the (1,2) entry and using the integrity of $\hat V$), which then implies $\phi=1$ (by inspection of the (1,1) entry, 
$(v_1,v_2)=(f_1\gamma,e_1\gamma)$ (by inspection of entries (1,3) and (2,3)), and 
$0=\epsilon(f_1\gamma)=\epsilon(v_1)=\epsilon(w)-1$ (using $v_1=f_1\gamma$). It follows that 
$(v_1,v_2)\in(f_1,e_1)\hat V$ and $\epsilon(w)=1$. Conversely, if $(v_1,v_2,w)$ satisfies these conditions, with 
$(v_1,v_2)=(f_1\gamma,e_1\gamma)$, then $u(v_1,v_2,w)=M(1,0,\gamma,w-1)$, which belongs to $\mathrm C_{21}(\rho_1)^\times$
since $1+(w-1)\in\hat V^\times$. (c) follows. 

It follows from Lem. \ref{lem:15:13:1408}(b) that
$$
\mathrm C_{21}(\rho_0)^\times=\{X(\Pi,a,c,d)|
\Pi\in\mathbf k[[u]]^\times,a,c,d\in \mathrm C_{\hat V}(e_0),\Pi(f_0)+d\in\hat V^\times\}
\subset (T_{21}\hat V)^\times,  
$$
where 
\begin{equation}\label{def:X}
X(\Pi,a,c,d):=\Pi(f_0)I_3+\begin{pmatrix}
(e_0-f_0)a&0&(e_0-f_0)c+d\\e_1a&0&e_1c\\0&0&d
\end{pmatrix}.     
\end{equation}
If $(v_1,v_2,w)\in\hat V^2\times\hat V^\times$ and $(\Pi,a,c,d)\in\mathbf k[[u]]^\times\times \mathrm C_{\hat V}(e_0)^3$ are such that 
$\epsilon(v_2)=0,\epsilon(v_1)=\epsilon(c)-1$, $\Pi(f_0)+d\in\hat V^\times$, then the equality $u(v_1,v_2,w)=M(\Pi,a,c,d)$ implies 
$a=0$ (by inspection of the (1,2) entry and using the integrity of $\hat V$), $\Pi=1$ (by inspection of entry (2,2)), 
$w=d+1$ (by inspection of the (2,2) entry and using $\Pi=1$), which implies $w\in \mathrm C_{\hat V}(e_0)$, which by 
$w\in\hat V^\times$ implies $w\in \mathrm C_{\hat V}(e_0)^\times$, $v_2=e_1c$ (by inspection of the (2,3) entry), 
and $v_1=(e_0-f_0)c+w-1$ (by inspection of the (1,3) entry and using $w=d+1$).  Therefore $w\in \mathrm C_{\hat V}(e_0)^\times$
and $(v_1,v_2)=((e_0-f_0)c+w-1,e_1c)$, where $c\in \mathrm C_{\hat V}(e_0)$. Conversely, if $w,v_1,v_2,c$ satisfy these conditions, 
then $u(v_1,v_2,w)=X(1,0,c,w-1)$, with $c,w-1\in \mathrm C_{\hat V}(e_0)$ and $1+(w-1)\in\hat V^\times$, therefore 
$u(v_1,v_2,w)\in \mathrm C_{21}(\rho_0)^\times$. All this implies (d). 

For $(\gamma,\delta)\in \hat V\times\hat V^\times$ with $\epsilon(\delta)=1$ and 
$(c,d)\in\mathrm{C}_{\hat V}(e_0)\times\mathrm{C}_{\hat V}(e_0)^\times$, one has 
\begin{align}\label{useful:for:eq:of:gps}
&\nonumber u(f_1\gamma,e_1\gamma,\delta)u((e_0-f_0)c+d-1,e_1c,d)=u(f_1\gamma d+(e_0-f_0)c+d-1,e_1\gamma d+e_1c,\delta d)
\\&\nonumber =u(f_1(\gamma d+c)+(e_0-f_0-f_1)c+d-1,e_1(\gamma d+c),\delta d)
\\&=u(f_1(\gamma d+c)+(e_0+f_\infty)c+d-1,e_1(\gamma d+c),\delta d)
\end{align}
where the first equality follows from \eqref{pdt:in:U}. Since $(\gamma d+c,\delta)\in\hat V\times\hat V^\times$, 
$(c,d)\in\mathrm{C}_{\hat V}(e_0)\times\mathrm{C}_{\hat V}(e_0)^\times$ and $\epsilon(\delta)=1$,  
(c) and (d) then imply the inclusion of the left-hand side of \eqref{comp:product:subgroups} in its right-hand side. 
Conversely, \eqref{useful:for:eq:of:gps} implies that for $(\gamma,\delta)\in\hat V\times\hat V^\times, (c,d)\in \mathrm{C}_{\hat V}(e_0)
\times\mathrm{C}_{\hat V}(e_0)^\times,\epsilon(\delta)=1$, there holds  
$$
u(f_1\gamma+(e_0+f_\infty)c+d-1,e_1\gamma,\delta d)=u(f_1(\gamma-c)\delta^{-1},e_1(\gamma-c)\delta^{-1},\delta)u((e_0-f_0)c+d-1,e_1c,d), 
$$
where $((\gamma-c)\delta^{-1},\delta)\in\hat V\times\hat V^\times$ and $\epsilon(\delta)-1$, and 
$(c,d)\in\mathrm{C}_{\hat V}(e_0)\times\mathrm{C}_{\hat V}(e_0)^\times$, which 
by using (b) and (c) proves the opposite inclusion. This proves (e). 
\end{proof}

\begin{lem}\label{lem:33to31}
If $(C',r')\in M_{3,1}\hat V\times\hat V$ is such that $C'\mathrm{C}_{\hat V}(e_0)R_{r'}\subset \mathrm C_{21}(\rho_0)$ (inclusion of subsets of $M_3\hat V$), 
then 
$$
C'\mathrm{C}_{\hat V}(e_0)r'\subset 
\begin{pmatrix}
    e_0-f_0\\e_1\\0
\end{pmatrix}\mathrm{C}_{\hat V}(e_0)+\begin{pmatrix}
    1\\0\\1
\end{pmatrix}\mathrm{C}_{\hat V}(e_0)
$$ 
(inclusion of subsets of $M_{3,1}\hat V$). 
\end{lem}

\begin{proof}
By \eqref{value:c:a:rho0}, the hypothesis means that for any $x\in \mathrm{C}_{\hat V}(e_0)$, there exists $(\Pi_x,a_x,c_x,d_x)\in 
\mathbf k[[u]]\times \mathrm{C}_{\hat V}(e_0)^3$ such that 
$C'\mathrm{C}_{\hat V}(e_0)R_{r'}=X(\Pi_x,a_x,c_x,d_x)$, where the right-hand side is given by 
\eqref{def:X}. The (2,2) entry of this equality implies $\Pi_x=0$, and its (2,1) entry together with the integrity of 
$\hat V$ implies $a_x=0$, therefore $C'\mathrm{C}_{\hat V}(e_0)R_{r'}=X(0,0,c_x,d_x)$, which implies
\begin{equation}\label{temp:33to31}
C'\mathrm{C}_{\hat V}(e_0)R_{r'}\begin{pmatrix}
    0\\0\\1
\end{pmatrix}=X(0,0,c_x,d_x)\begin{pmatrix}
    0\\0\\1
\end{pmatrix}.     
\end{equation}
One checks the equalities
$$
R_{r'}\begin{pmatrix}
    0\\0\\1
\end{pmatrix}=r',\quad X(0,0,c_x,d_x)\begin{pmatrix}
    0\\0\\1
\end{pmatrix}=\begin{pmatrix}
    e_0-f_0\\e_1\\0
\end{pmatrix}c_x+\begin{pmatrix}
    1\\0\\1
\end{pmatrix}d_x,  
$$
therefore \eqref{temp:33to31} gives 
$$
C'\mathrm{C}_{\hat V}(e_0)r'=\begin{pmatrix}
    e_0-f_0\\e_1\\0
\end{pmatrix}c_x+\begin{pmatrix}
    1\\0\\1
\end{pmatrix}d_x,  
$$
which implies the result. 
\end{proof}

\begin{lem}\label{lem240506:1638}
(a)   $\mathrm{C}_V(e_0+f_\infty)=\mathbf k[e_0,f_\infty]$ (equality of subalgebras of $V$). 

(b) For $x\in V$, one has $[x,e_0+f_\infty]\in\mathrm{C}_{\hat V}(e_0)$ iff  $x\in\mathrm{C}_V(e_0)$.
\end{lem}

\begin{proof}
(a)  
Let us denote by $\mathcal V^+$ the sum of the parts of 
$\mathcal V$ of positive degree with respect to the grading $\mathrm{deg}(e_0)=0$, $\mathrm{deg}(e_1)=1$. 
Then $\mathcal V=\mathbf k[e_0]\oplus \mathcal V^+$. The $\mathbf k$-module $\mathcal V^+$ is further graded by the 
 grading $\mathrm{deg}(e_1)=0$, $\mathrm{deg}(e_0)=1$, let $\mathcal V^+=\oplus_{k>0}\mathcal V^+\{k\}$ be the 
 corresponding decomposition. The endomorphism $x\mapsto [e_0,x]$ of $\mathcal V$ has kernel equal to $\mathbf k[e_0]$
 and induces injections $\mathcal V^+\{k\}\to\mathcal V^+\{k+1\}$ for any $k>0$. 

Let now $(N,g)$ be a pair of a $\mathbf k$-module $N$ and of an endomorphism $g$ of $N$. 
Then $\mathcal V\otimes N$ is decomposed as 
$(\mathbf k[e_0]\otimes N) \oplus (\mathcal V\{1\}\otimes N)\oplus (\mathcal V\{2\}\otimes N)\oplus\cdots$, 
and the endomorphism $[e_0,-]\otimes id_N+id_{\mathcal V}\otimes g$ of $\mathcal V\otimes N$ 
has filtration degree 1 with respect to the increasing 
filtration associated with this grading, the associated graded morphism being the direct sum of the zero map on 
$\mathbf k[e_0]\otimes N$ and of the maps $[e_0,-]\otimes id_N : \mathcal V\{k\}\otimes N
\to \mathcal V\{k+1\}\otimes N$ over $k>0$, whose kernel is $\mathbf k[e_0]\otimes N$, which implies  
$\mathrm{ker}([e_0,-]\otimes id_N+id_{\mathcal V}\otimes g)\subset\mathrm k[e_0]\otimes N$. 
It follows that $\mathrm{ker}([e_0,-]\otimes id_N+id_{\mathcal V}\otimes g)=\mathbf k[e_0]\otimes \mathrm{ker}(g)$.  

Applying this statement to $(N,g):=(\mathcal V,[e_\infty,-])$, one sees that the kernel of $[e_0+f_\infty,-]$ 
is equal to $\mathbf k[e_0] \otimes \mathbf k[e_\infty]=\mathbf k[e_0,f_\infty]$.

(b) One has 
$\{x\in V|[x,e_0+f_\infty]\in\mathrm{C}_V(e_0)\}=\{x\in V|[e_0,[x,e_0+f_\infty]]=0\}
=\{x\in V|[e_0+f_\infty,[e_0,x]]=0\}$, where the second 
equality follows from the commutation of $e_0$ and $e_0+f_\infty$.  
This implies 
\begin{equation}\label{toto:1305}
\mathrm{C}_V(e_0) \subset \{x \in V|[x,e_0+f_\infty] \in \mathrm{C}_V(e_0)\}.
\end{equation}
Let us show the opposite inclusion. One has $\{x\in V|[x,e_0+f_\infty]\in\mathrm{C}_V(e_0)\}
=\{x\in V|[e_0+f_\infty,[e_0,x]]=0\}=\{x\in V|[e_0,x]\in\mathbf k[e_0,f_\infty]\}$, where the last equality follows from (a). 
The direct sum decomposition $\mathcal V=\mathbf k[e_0] \oplus \mathcal V^+$ induces a direct sum decomposition 
$V=(\mathbf k[e_0] \otimes \mathcal V) \oplus (\mathcal V^+ \otimes \mathcal V)$; let $x \in 
\{x\in V|[e_0,x]\in\mathbf k[e_0,f_\infty]\}$ let $x=x_e+x_+$ be the corresponding decomposition of $x$. Then 
$[e_0,x]=[e_0,x_+] \in \mathcal V^+ \otimes \mathcal V$ as the endomorphism $[e_0,-]$ of $\mathcal V$ is 0 when restricted to 
$\mathbf k[e_0]$ and preserves $\mathcal V^+$; since the intersection of $\mathcal V^+ \otimes \mathcal V$  with 
$\mathbf k[e_0,f_\infty]$ is 0, this implies $[e_0,x]=0$, therefore $x \in \mathrm{C}_V(e_0)$. Combining this with 
\eqref{toto:1305}, one gets
$$\mathrm{C}_V(e_0) = \{x \in V|[x,e_0+f_\infty] \in \mathrm{C}_V(e_0)\}.$$
\end{proof}

\begin{lem}\label{lem:9:12:conj:in:CVe0}
Let $\kappa\in \mathrm C_{\hat V}(e_0)\cap F^2\hat V$ and 
$\alpha\in\hat V^\times$ be such that 
\begin{equation}\label{eq:to:solve:IN:V}
 \alpha (e_0+f_\infty+\kappa)\alpha^{-1}\subset 
\mathrm C_{\hat V}(e_0),    
\end{equation}
 then 
$\alpha\in \mathrm C_{\hat V}(e_0)^\times$. 
\end{lem}

\begin{proof}
Let $(x_i)_{i>0}$ be free noncommutative variables, with $\mathrm{deg}x_i=i$. Define $(P_i(x_1,\ldots,x_i))_{i>0}$ be the family of 
polynomials in these variables, such that $\mathrm{deg}(P_i(x_1,\ldots,x_i))=i$ and $1+\sum_{i>0}P_i(x_1,\ldots,x_i)=(1+\sum_{i>0}x_i)^{-1}$. 
Then there exists a family of polynomials $(Q_i(x_1,\ldots,x_{i-1}))_{i>0}$ such that $P_i(x_1,\ldots,x_i)=-x_i+Q_i(x_1,\ldots,x_{i-1})$. 

Let $\alpha=\sum_{i\geq0}\alpha_i$, $\kappa=\sum_{i\geq 0}\kappa_i$ be the degree expansions of $\alpha,\kappa$ 
($e_s,f_s$ being of degree 1 for $s=0,1$). One has $\kappa_0=\kappa_1=0$, $\kappa_i\in \mathrm C_{\hat V}(e_0)$ for 
$i\geq2$ and one may assume $\alpha_0=1$. Let us prove inductively 
on $d>0$ that $\alpha_d\in\mathrm C_{\hat V}(e_0)$. The degree 2 component
of \eqref{eq:to:solve:IN:V} yields $[\alpha_1,e_0+f_\infty]+\kappa_2\in \mathrm C_{\hat V}(e_0)$, which together with $\kappa_2\in \mathrm C_{\hat V}(e_0)$
and Lem. \ref{lem240506:1638}(b) implies $\alpha_1\in \mathrm C_{\hat V}(e_0)$. Let $d>0$ and assume that $\alpha_1,\ldots,\alpha_d\in 
\mathrm C_{\hat V}(e_0)$. Let $\alpha^{-1}=\sum_{i\geq0}(\alpha^{-1})_i$ be the degree expansion of $\alpha^{-1}$.
For $i=1,\ldots,d$, one has $(\alpha^{-1})_i=P_i(\alpha_1,\ldots,\alpha_i)$, therefore 
$(\alpha^{-1})_1,\ldots,(\alpha^{-1})_d\in \mathrm C_{\hat V}(e_0)$. Moreover, 
$(\alpha^{-1})_{d+1}=-\alpha_{d+1}+Q_{d+1}(\alpha_1,\ldots,\alpha_d)$. 
Then the degree $d+2$ part of \eqref{eq:to:solve:IN:V} yields 
$$
[\alpha_{d+1},e_0+f_\infty]-(e_0+f_\infty)Q_{d+1}(\alpha_1,\ldots,\alpha_d)
+\sum_{j\geq 2, i+j+k=d+2}\alpha_i \kappa_j(\alpha^{-1})_k\in \mathrm C_{\hat V}(e_0),
$$
Since the indices $i,k$ in $\sum_{j\geq 2, i+j+k=d+2}\alpha_i \kappa_j(\alpha^{-1})_k\in \mathrm C_{\hat V}(e_0)$ are $\leq d$, 
this expression belongs to $\mathrm C_{\hat V}(e_0)$, therefore $[\alpha_{d+1},e_0+f_\infty]\in \mathrm C_{\hat V}(e_0)$, which by 
Lem. \ref{lem240506:1638}(b) implies $\alpha_{d+1}\in \mathrm C_{\hat V}(e_0)$, proving the induction step.
\end{proof}

\begin{lem}\label{lem:9:13}
(a) For any $v\in\hat V$, the relation $v\cdot (e_0+f_\infty)\in e_1\hat V$ implies $v\in e_1\hat V$. 

(b) Let $\alpha\in e_0+f_\infty+F^2\hat V$. Then
$$
\forall v\in \hat V,\quad (v\cdot \alpha\in e_1\hat V)\implies (v\in e_1\hat V). 
$$

\end{lem}

\begin{proof}
(a) Right multiplication by $e_0+f_\infty$ induces an endomorphism of the $\mathbf k$-module $V$ which preserves the 
$\mathbf k$-submodule $e_1V$, and therefore an endomorphism of the quotient $V/e_1V$, which we denote by $\phi$. This quotient is 
isomorphic to $(\mathcal V/e_1\mathcal V)\otimes \mathcal V$, and by the direct sum decomposition 
$\mathcal V=\mathbf k\oplus e_0\mathcal V\oplus e_1\mathcal V$, also to the tensor product 
$(\mathbf k\oplus e_0\mathcal V)\otimes\mathcal V$. Since  $\mathbf k\oplus e_0\mathcal V$ is a subalgebra of $\mathcal V$, 
this tensor product has an algebra structure, and $\phi$ is conjugated with the isomorphism $V/e_1V\simeq 
(\mathbf k\oplus e_0\mathcal V)\otimes\mathcal V$ to the endomorphism given by right multiplication by 
$e_0\otimes1+1\otimes e_\infty$. The algebra $\mathbf k\oplus e_0\mathcal V$ is graded by the $e_0$-degree 
for which $\mathrm{deg}(e_0)=1$, $\mathrm{deg}(e_1)=0$. The tensor product of this grading with the 
trivial grading of $\mathcal V$ (for which the degree 0 part is $\mathcal V$ itself) is a grading on 
$(\mathbf k\oplus e_0\mathcal V)\otimes\mathcal V$. 

This induces an increasing $\mathbf k$-module filtration $\mathbf k \otimes\mathcal V=F_0\subset F_1\subset\cdots$ of  
$(\mathbf k\oplus e_0\mathcal V)\otimes\mathcal V$ which is total (i.e. this $\mathbf k$-module is equal to $\cup_i F_i$), 
and the associated graded of right multiplication by $e_0\otimes1+1\otimes e_\infty$ is filtered of degree 1 (i.e., takes 
$F_i$ to $F_{i+1}$). The associated graded endomorphism is then the endomorphism 
$(x\mapsto xe_0)\otimes id_{\mathcal V}$ 
of  $(\mathbf k\oplus e_0\mathcal V)\otimes\mathcal V$ 
given by right multiplication by $e_0\otimes 1$, which is equal to the tensor product of the endomorphism of 
$\mathbf k\oplus e_0\mathcal V$ given by right multiplication by $e_0$ with the identity endomorphism of $\mathcal V$. 
It follows from the integrity of $\mathcal V$ that right multiplication of $e_0$ is an injective endomorphism of 
$\mathbf k\oplus e_0\mathcal V$. It follows that the endomorphism 
$(x\mapsto xe_0)\otimes id_{\mathcal V}$  is injective as well, therefore that the endomorphism of 
$\mathbf k \otimes\mathcal V=F_0\subset F_1\subset\cdots$ of  
$(\mathbf k\oplus e_0\mathcal V)\otimes\mathcal V$ induced by right multiplication by $e_0\otimes1+1\otimes e_\infty$ 
is injective, and therefore that so is the endomorphism $\phi$ of $V/e_1V$.

The usual grading of $V$ (for which $e_0,f_0,e_1,f_1$ have degree 1) induces a grading of $V/e_1V$, for which the endomorphism 
$\phi$ has degree 1. Its graded completion $\hat\phi$ is injective since $\phi$ is so, and is the endomorphism of $\hat V/e_1\hat V$
induced by right multiplication by $e_0+f_\infty$. (a) follows from the injectivity of this endomorphism. 

Let $\alpha$ be as in (b) and assume $v\in\hat V$ and $v\alpha\in e_1\hat V$. 
Let $\alpha=\sum_{i\geq1}\alpha_i$, $v=\sum_{i\geq0}v_i$ be the degree expansions of $\alpha,v$, 
with $\alpha_1=e_0+f_\infty$.  Let us show by induction on $d\geq0$ that $v_d\in e_1\hat V$. The degree 1 part of
$v\alpha\in e_1\hat V$ is the relation $v_0(e_0+f_\infty)\in e_1\hat V$, which implies $v_0=0$ as $v_0\in\mathbf k$. 
Let $d\geq0$ and assume $v_0,\ldots,v_d\in e_1\hat V$. The degree $d+2$ component of $v\alpha\in e_1\hat V$ implies
$v_{d+1}(e_0+f_\infty)+\sum_{0}^{d}v_i\alpha_{d+2-i}\in e_1\hat V$, which by the assumptions on $v_0,\ldots,v_d$
implies $v_{d+1}(e_0+f_\infty)\in e_1\hat V$. (a) then implies $v_{d+1}\in e_1\hat V$, thus proving the induction step. (b) follows. 
\end{proof}

\begin{lem}\label{lem:van:coh}
(a) The diagram of $\mathbf k$-modules 
\begin{equation}\label{complexe:6mars}
V\oplus \mathrm C_V(e_0)\oplus \mathrm C_V(e_0)\to 
V\oplus \mathrm C_V(e_0)\oplus V\to V,     
\end{equation}
where the first map is $(v,h,k) \mapsto 
((e_0+f_\infty)v-h,h\cdot (e_0+f_\infty)-(e_0+f_\infty)k,v\cdot (e_0+f_\infty)-k)$
and the second map is
$(\alpha,\tau,s)\mapsto\alpha\cdot(e_0+f_\infty)-(e_0+f_\infty)s+\tau$, is an acyclic complex. 

(b) Diagram \eqref{complexe:6mars} decomposes as the direct sum over $n\geq 0$ of the diagrams 
$$
V_{n-1}\oplus \mathrm C_V(e_0)_n\oplus \mathrm C_V(e_0)_n\to 
V_n\oplus \mathrm C_V(e_0)_{n+1}\oplus V_n\to V_{n+1},  
$$
where the indices refer to the total degree (for which $e_0,e_1,f_0,f_1$ all have degree $1$), which are all acyclic. 
\end{lem}

\begin{proof}
(a) Diagram \eqref{complexe:6mars} is obviously a complex. It is graded for the $e_1$-degree, which is the algebra degree on $V$ for which $e_1$ 
has degree $1$ and
$e_0,f_0,f_1$ have degree $0$, and for which $\mathrm C_V(e_0)$ lies in degree 0. For $d\geq0$, denote by $V\{d\}\subset V$ 
the part of $e_1$-degree $d$. 
Then $V=\oplus_{d\geq0}V\{d\}$ and $V\{0\}=\mathrm C_V(e_0)$. The parts of \eqref{complexe:6mars} of $e_1$-degree $d$ are respectively 
complexes 
\begin{equation}\label{complex:deg:O}
 \mathrm C_V(e_0)^{\oplus3}\to\mathrm C_V(e_0)^{\oplus3}\to\mathrm C_V(e_0)   
\end{equation}
if $d=0$ and 
\begin{equation}\label{complex:deg:>0}
V\{d\}\to V\{d\}\oplus V\{d\}\to V\{d\}   
\end{equation}
if $d>0$. 

The acyclic complex from Lem. \ref{lem:exact:seq:0505} is graded for the $e_1$-degree, and for any $d>0$, its part of $e_1$-degree
$d$ is isomorphic to \eqref{complex:deg:>0}. It follows that \eqref{complex:deg:>0} is acyclic. 

It follows from the equality $\mathrm C_V(e_0)=\mathbf k[e_0]\otimes\mathcal V$ (see Lem. \ref{lem:comm:e0}(a))
that the $\mathbf k$-module $\mathrm C_V(e_0)$ is graded by the $e_0$-degree; this induces an 
increasing filtration, defined by $F_n\mathrm C_V(e_0):=
\mathbf k[e_0]_{\leq n}\otimes\mathcal V$ for $n\geq 0$ (the index $\leq n$ meaning polynomials of degree $\leq n$), which is total.
This induces a filtration on the complex \eqref{complex:deg:O}, whose 
degree $n$ step is given by the subcomplex 
$$
F_{n-1}\mathrm C_V(e_0)\oplus F_n\mathrm C_V(e_0)\oplus F_n\mathrm C_V(e_0) \to
F_n\mathrm C_V(e_0)\oplus F_{n+1}\mathrm C_V(e_0)\oplus F_n\mathrm C_V(e_0)\to F_{n+1}\mathrm C_V(e_0)
$$
for $n\in\mathbb Z$. The associated graded complex is a complex
\begin{equation}\label{ass:graded:complex}
\mathbf k[e_0]\otimes(\mathcal V\oplus \mathcal V\oplus \mathcal V) \to
\mathbf k[e_0]\otimes(\mathcal V\oplus \mathcal V\oplus \mathcal V)\to \mathbf k[e_0]\otimes\mathcal V,     
\end{equation}
whose degree $n$ part is the complex
$$
\mathbf ke_0^{n-1}\otimes\mathcal V\oplus \mathbf ke_0^n\otimes\mathcal V\oplus \mathbf ke_0^n\otimes\mathcal V \to
\mathbf ke_0^n\otimes\mathcal V\oplus \mathbf ke_0^{n+1}\otimes\mathcal V\oplus \mathbf ke_0^n\otimes\mathcal V\to \mathbf ke_0^{n+1}\otimes\mathcal V
$$
where the maps are given by $(e_0^{n-1}\otimes v,e_0^n\otimes h,e_0^n\otimes k)\mapsto (e_0^n\otimes(v-h),e_0^{n+1}\otimes(h-k),e_0^n\otimes(v-k))$
and $(e_0^n\otimes\alpha,e_0^{n+1}\otimes\tau,e_0^n\otimes s)\mapsto e_0^{n+1}\otimes(\alpha+\tau-s)$. 
For $n\geq1$, this is isomorphic to the complex
$$
\mathcal V\oplus \mathcal V\oplus \mathcal V \to
\mathcal V\oplus \mathcal V\oplus\mathcal V\to \mathcal V
$$
where the maps are given by $(v,h,k)\mapsto (v-h,h-k,v-k)$
and $(\alpha,\tau,s)\mapsto \alpha+\tau-s$, which is acyclic, since $\alpha+\tau-s=0$
implies that $(\alpha,\tau,s)$ is the image of $(v,h,k):=(0,-\alpha,-s)$. For $n=0$, 
this is isomorphic to the complex
$$
\mathcal V\oplus \mathcal V \to
\mathcal V\oplus \mathcal V\oplus\mathcal V\to \mathcal V
$$
where the maps are given by $(h,k)\mapsto (-h,h-k,-k)$
and $(\alpha,\tau,s)\mapsto \alpha+\tau-s$, which is acyclic, since $\alpha+\tau-s=0$
implies that $(\alpha,\tau,s)$ is the image of $(h,k):=(-\alpha,-s)$. For $n=-1$, this is 
isomorphic to the complex $0\to \mathcal V\stackrel{id}{\to}\mathcal V$, which is acyclic. 
For $n\leq-2$, this is the complex $0\to 0\to 0$, which is acyclic. 

All this implies that \eqref{ass:graded:complex}, which is the associated graded complex of \eqref{complex:deg:O}, is acyclic. 
This implies that \eqref{complex:deg:O} is acyclic. 

(b) The first statement follows from the fact that $e_0+f_\infty$ is homogeneous of degree $1$ for the total degree. The second statement 
then follows from this and from (a). 
\end{proof}

\begin{lem}\label{lemma:in:hatV}
Let $\alpha\in\hat V^\times$, $\tau\in(e_0+f_\infty+F^2\hat V)\cap \mathrm{C}_{\hat V}(e_0)$ and $s\in \hat V$
satisfy the relation
$$
\alpha\tau=(e_0+f_\infty)s. 
$$
Then there exist $\tilde\alpha,\tilde\tau\in \mathrm{C}_{\hat V}(e_0)$ and $\tilde\gamma\in \mathrm{C}_{\hat V}(e_0)^\times$
with $\epsilon(\tilde\gamma)=\epsilon(\tilde\tau)=1$, 
such that 
$$
\alpha=\epsilon(\alpha)(1+(e_0+f_\infty)\cdot\tilde\alpha)\tilde\gamma^{-1},\quad \tau=\tilde\gamma\cdot(e_0+f_\infty)\cdot\tilde\tau,\quad 
s=\epsilon(\alpha)(1+\tilde\alpha\cdot(e_0+f_\infty))\tilde\tau; 
$$
in particular, $s\in \mathrm{C}_{\hat V}(e_0)$. 
\end{lem}

\begin{proof}
Let $X$ be the set of triples $(\alpha,\tau,s)$ as in the hypothesis of this statement. 
Then $\alpha\tau\in\epsilon(\alpha)(e_0+f_\infty)+F^2\hat V$. On the other 
hand, $(e_0+f_\infty)s\in\epsilon(s)(e_0+f_\infty)+F^2\hat V$. The equality $\alpha\tau=(e_0+f_\infty)s$ then implies 
$\epsilon(s)=\epsilon(\alpha)$, therefore $s$ is invertible. It follows that 
$$
X=\{(\alpha,\tau,s)\in 
\hat V^\times\times\mathrm C_{\hat V}(e_0)\times\hat V^\times|\alpha\tau=(e_0+f_\infty)s\}
$$
and that 
\begin{equation}\label{on:X}
    \forall (\alpha,\tau,s)\in X,\quad \epsilon(s)=\epsilon(\alpha)\in\mathbf k^\times. 
\end{equation}
Let $(\mathrm{C}_{\hat V}(e_0),\cdot_{e_0+f_\infty})$ be the unitless algebra defined by $v\cdot_{e_0+f_\infty}v':=v(e_0+f_\infty)v'$
and $\mathbf k\oplus (\mathrm{C}_{\hat V}(e_0),\cdot_{e_0+f_\infty})$ be the corresponding algebra with unity. Its group of invertible elements 
$(\mathbf k\oplus (\mathrm{C}_{\hat V}(e_0),\cdot_{e_0+f_\infty}))^\times$ is $\mathbf k^\times\times\mathrm{C}_{\hat V}(e_0)$ equipped with the product 
$(\lambda,v)\cdot(\lambda',v'):=(\lambda\lambda',\lambda v'+v\lambda'+v\cdot (e_0+f_\infty)v')$.  There are algebra morphisms 
$\mathrm{C}_{\hat V}(e_0)\leftarrow\mathbf k\oplus (\mathrm{C}_{\hat V}(e_0),\cdot_{e_0+f_\infty})\to \mathrm{C}_{\hat V}(e_0)$
given by $\lambda+(e_0+f_\infty)v\mapsfrom(\lambda,v)\mapsto\lambda+v\cdot (e_0+f_\infty)$, inducing group morphisms
$\mathrm{C}_{\hat V}(e_0)^\times\leftarrow(\mathbf k\oplus (\mathrm{C}_{\hat V}(e_0),\cdot_{e_0+f_\infty}))^\times\to \mathrm{C}_{\hat V}(e_0)^\times$. 

Let $\mathrm{C}_{\hat V}(e_0)^\times_{1}$ be the kernel of the group morphism $\epsilon : \mathrm{C}_{\hat V}(e_0)^\times\to \mathbf k^\times$. The product 
group 
$G:=(\mathbf k\oplus (\mathrm{C}_{\hat V}(e_0),\cdot_{e_0+f_\infty}))^\times\times (\mathrm{C}_{\hat V}(e_0)^\times_{1})^2$ acts on $X$ by 
$$
((\lambda,v),h,k)\bullet(\alpha,\tau,s):=
((\lambda+(e_0+f_\infty)v)\alpha h^{-1},h\tau k^{-1},(\lambda+v\cdot (e_0+f_\infty))s k^{-1}).
$$
For $n\geq1$, let us denote by $X_n\subset X$ the subset of triples $(\alpha,\tau,s)$ such that 
$\alpha,s\in 1+F^n\hat V$, $\tau\in e_0+f_\infty+F^{n+1}\hat V$ and by $G_n\subset G$
the subgroup of tuples $((\lambda,v),h,k)$ defined by $\lambda=1$, $v\in F^{n-1}\hat V$, $h,k\in 1+F^n\hat V$. 
Set also $X_0:=X$, $G_0:=G$. 

Let us prove the inclusion $X_n\subset G_n\cdot X_{n+1}$ for any $n\geq0$. 
It follows from \eqref{on:X} that for any $(\alpha,\tau,s)\in X_0=X$, 
$((\epsilon(\alpha)^{-1},0),1,1)\bullet(\alpha,\tau,s)\in X_1$, which proves $X_0\subset G_0\cdot X_1$. 
If now $n>0$ and $(\alpha,\tau,s)\in X_n$, then denoting by $x=\sum_{n\geq0}x_n$ the degree expansion of an 
element of $\hat V$, one has $\alpha\in 1+\alpha_n+F^{n+1}\hat V$, $s\in 1+s_n+F^{n+1}\hat V$, 
$\tau\in e_0+f_\infty+\tau_{n+1}+F^{n+2}\hat V$. The relation $\alpha\tau=(e_0+f_\infty)s$ then implies 
$\tau_{n+1}+\alpha_n\cdot (e_0+f_\infty)-(e_0+f_\infty)s_n=0$. Lem. \ref{lem:van:coh} then implies the existence of 
a triple $(v_{n-1},h_n,k_n)\in \hat V_{n-1}\times \mathrm C_{\hat V}(e_0)_n\times \mathrm C_{\hat V}(e_0)_n$ such that 
$$
(\alpha_n,\tau_{n+1},s_n)=((e_0+f_\infty)v_{n-1}-h_n,h_n\cdot(e_0+f_\infty)-(e_0+f_\infty)k_n,v_{n-1}\cdot(e_0+f_\infty)-k_n). 
$$
Then $g:=((1,v_{n-1}),h_n,k_n)\in G_n$, and 
\begin{align*}
    & g^{-1}\bullet x\in(1+\alpha_n-(e_0+f_\infty)v_{n-1}+h_n+F^{n+1}\hat V, 
    \\&
e_0+f_\infty+\tau_{n+1}-h_n\cdot(e_0+f_\infty)+(e_0+f_\infty)k_n+F^{n+2}\hat V,
1+s_n-v_{n-1}\cdot(e_0+f_\infty)+k_n+F^{n+1}\hat V),     
\end{align*}
therefore $g^{-1}\bullet x\in X_{n+1}$.   
The inclusion $X_n\subset G_n\cdot X_{n+1}$ follows. 

Let now $x\in G$ and define inductively $(x_n)_{n\geq0}$, $(g_n)_{n\geq0}$ by 
$x_n\in X_n$, $g_n\in G_n$, $x_0:=x$ and $x_{n+1}=g_n^{-1}\bullet x_n$. Then the 
product $g:=g_0g_1\cdots$ converges in $G$, and $g^{-1}\bullet x$ belongs to 
$\cap_{n\geq0}X_n$, which is the element $(1,e_0+f_\infty,1)$. Therefore 
$x=g\bullet(1,e_0+f_\infty,1)$. Denoting by $\lambda\in\mathbf k^\times$, $\tilde v\in 
\mathrm C_{\hat V}(e_0)$, $\tilde\gamma,\tilde\tau\in \mathrm{C}_{\hat V}(e_0)^\times_{1}$
the elements such that $g=((\lambda,\tilde v),\tilde\gamma,\tilde\tau^{-1})$, one then obtains 
$$
\alpha=\lambda(1+(e_0+f_\infty)\cdot\tilde\alpha)\tilde\gamma^{-1},\quad \tau=\lambda\cdot(e_0+f_\infty)\cdot\tilde\tau,\quad 
s=\epsilon(\alpha)(1+\tilde\alpha\cdot(e_0+f_\infty))\tilde\tau,  
$$
which implies the announced formulas since $\epsilon(\tilde\gamma)=\epsilon(\tilde\tau)=1$. 
\end{proof}

\begin{lem}\label{lem:quotients} (a)  If $a\in \hat V$ is such that $a\cdot  (e_0+f_\infty)\in \mathrm C_{\hat V}(e_0)$, 
then $a\in \mathrm C_{\hat V}(e_0)$.

(b) If $a\in \hat V$ and $\alpha\in e_0+f_\infty+(F^2\hat V\cap\mathrm C_{\hat V}(e_0))$ are such that $a\alpha\in \mathrm C_{\hat V}(e_0)$, 
then $a\in \mathrm C_{\hat V}(e_0)$. 
\end{lem}

\begin{proof}
(a) Let $\mathcal V=\oplus_{d\geq0}\mathcal V\{d\}$ be the decomposition of $\mathcal V$ for the $e_0$-grading (for which 
$e_0,e_1$ have degrees $1,0$). This induces a $\mathbf k$-module grading $V=\mathcal V^{\otimes2}=\oplus_{d\geq0}\mathcal V\{d\}\otimes\mathcal V$
on $V$; the corresponding increasing filtration is given by $F_dV=\oplus_{i=0}^d \mathcal V\{i\}\otimes\mathcal V$. 
The endomorphism $\phi : x\mapsto x\cdot  (e_0+f_\infty)$ of $V$ has degree 1 for this filtration, i.e. $\phi(F_dV)\subset F_{d+1}V$ for any 
$d\geq0$. The associated graded endomorphism is the endomorphism of $V$ given by $(x\mapsto xe_0)\otimes id_{\mathcal V}$. Since 
$\mathcal V$ is a domain, its endomorphism $x\mapsto xe_0$ is injective. It follows that $\phi$ is injective. 

Let now $\mathcal V=\oplus_{d\geq0}\mathcal V\{\{d\}\}$ be the decomposition of $\mathcal V$ for the $e_1$-grading (for which 
$e_0,e_1$ have degrees $0,1$). It induces a decomposition $V=\oplus_{d\geq0}\mathcal V\{\{d\}\}\otimes\mathcal V$. 
The endomorphism $\phi$ is compatible with this decomposition, therefore induces an endomorphism $\phi\{\{d\}\}$ of 
$\mathcal V\{\{d\}\}\otimes\mathcal V$ for any $d\geq0$, which is injective since $\phi$ is. 

The decomposition of $V$ for the total degree (for which $e_0,e_1,f_0,f_1$ all have degree 1) is denoted $V=\oplus_d V_d$.
It is compatible with the decomposition $V=\oplus_{d\geq0}\mathcal V\{\{d\}\}\otimes\mathcal V$ and with the endomorphism $\phi$. 
It follows that the completion $\hat V$ admits a decomposition $\hat V=\hat\oplus_{d\geq0}(\mathcal V\{\{d\}\}\otimes\mathcal V)^\wedge$, 
that the completion $\hat\phi$ of $\phi$ restricts to endomorphisms of each $(\mathcal V\{\{d\}\}\otimes\mathcal V)^\wedge$ for $d\geq0$, 
and that this restriction is injective if $d>0$. Together with the equality $(\mathcal V\{\{0\}\}\otimes\mathcal V)^\wedge=\mathrm C_{\hat V}(e_0)$
, this implies the statement. 

(b) Let $a,\alpha$ be as in the hypothesis of (b). Let $a=\sum_{d\geq 0}$, $\alpha=\sum_{d\geq2}\alpha_d$ be the decompositions of 
$a,\alpha$ for the total degree of $\hat V$. Let us prove by induction of $d\geq 0$ that $a_d\in\mathrm{C}_{\hat V}(e_0)$. 
One has $a_0\in\hat V_0=\mathbf k\subset\mathrm{C}_{\hat V}(e_0)$. Assume $d\geq0$ and $a_0,\ldots,a_d\in \mathrm{C}_{\hat V}(e_0)$. 
The degree $d+2$ part of the relation $a\alpha\in \mathrm C_{\hat V}(e_0)$ gives $a_{d+1}\cdot (e_0+f_\infty)+\sum_{i=0}^{d}a_i\alpha_{d+2-i}
\in \mathrm C_{\hat V}(e_0)$, which implies  $a_{d+1}\cdot (e_0+f_\infty)\in \mathrm C_{\hat V}(e_0)$, which by (a) implies  
$a_{d+1}\in \mathrm C_{\hat V}(e_0)$, which proves the induction step. 
\end{proof}

\begin{lem}\label{lem:MAIN}
    One has the inclusion 
    $$
    \{g\in U^{(0)}|g\cdot\rho_{\mathrm{DT}}\in \mathrm{Hom}^{\bullet}_{\mathcal C\operatorname{-alg}}
(\hat{\mathcal V},T_{21}\hat V)\}\subset (U^{(0)}\cap \mathrm{C}_{21}(\rho_1)^\times)\cdot(U^{(0)}\cap \mathrm{C}_{21}(\rho_0)^\times). 
    $$
\end{lem}

\begin{proof} 
Let $g\in U^{(0)}$ be such that $g\cdot\rho_{\mathrm{DT}}\in \mathrm{Hom}^{\bullet}_{\mathcal C\operatorname{-alg}}
(\hat{\mathcal V},T_{21}\hat V)$. Since $g\in U^{(0)}$ and by Lem. \ref{lem:9:10:1345}(b), there exist $(v_1,v_2,w)\in \hat V^2\times\hat V^\times$ 
with $\epsilon(v_2)=0$ and 
\begin{equation}\label{epsilon:v1+1}
    \epsilon(v_1)=\epsilon(w)-1
\end{equation}
such that $g=u(v_1,v_2,w)$. Since 
$g\cdot\rho_{\mathrm{DT}}\in \mathrm{Hom}^{\bullet}_{\mathcal C\operatorname{-alg}}
(\hat{\mathcal V},T_{21}\hat V)$, for some $(r,C)\in \hat V\times M_{3,1}F^1\hat V$, one has
\begin{equation}\label{temporary:proof}
\mathrm C_{21}((g\cdot\rho_{\mathrm{DT}})(\hat{\mathcal V}))
=\mathbf k1+C\cdot \mathrm C_{\hat V}(e_0)\cdot R_r\text{ and }R_r\cdot C\in e_0+f_\infty+F^2\hat V. 
    \end{equation}
By Lem. \ref{lem:comm:diag:toto}(a), $g\cdot\rho_{\mathrm{DT}}\in \mathrm{Hom}^{1,\bullet}_{\mathcal C\operatorname{-alg}}
(\hat{\mathcal V},T_{21}\hat V)$, which together with Lem. \ref{lem:prelim:rho} 
implies that $r\in \hat V^\times$, and that for some $(s,t)\in\hat V\times F^1\hat V$, one has 
\begin{equation}\label{rappel:resultats:C:rt}
    C=\begin{pmatrix}
    f_1s\\e_1s\\t
\end{pmatrix}\quad \text{and}\quad rt\in \mathrm C_{\hat V}(e_0)\cap (e_0+f_\infty+F^2\hat V). 
\end{equation}
\eqref{temporary:proof} implies the equality in 
$$
(g^{-1}C)\cdot \mathrm C_{\hat V}(e_0)\cdot (R_rg)
\subset\mathbf k1+(g^{-1}C)\cdot \mathrm C_{\hat V}(e_0)\cdot (R_rg)=
\mathrm C_{21}((g^{-1}\bullet (g\cdot\rho_{\mathrm{DT}}))(\hat{\mathcal V}))\subset 
\mathrm C_{21}(\rho_0), 
$$
where the equality follows from $\mathrm C_{21}((g^{-1}\bullet (g\cdot\rho_{\mathrm{DT}}))(\hat{\mathcal V}))
=\mathrm C_{21}(\mathrm{Ad}_{g^{-1}}\circ(g\cdot\rho_{\mathrm{DT}})(\hat{\mathcal V}))
=\mathrm{Ad}_{g^{-1}}(\mathrm C_{21}((g\cdot\rho_{\mathrm{DT}})(\hat{\mathcal V})))$ and from 
\eqref{temporary:proof}, 
the last inclusion follows from the fact that $g^{-1}\bullet (g\cdot\rho_{\mathrm{DT}})$ is the morphism 
$\hat{\mathcal V}\to T_{21}\hat V$ given by $e_0\mapsto \rho_0$, $e_1\mapsto \mathrm{Ad}_{g^{-1}}(\rho_1)$, so that 
$(g^{-1}\bullet (g\cdot\rho_{\mathrm{DT}}))(\hat{\mathcal V})$ is the subalgebra of $T_{21}\hat V$ generated by 
$\rho_0$ and $\mathrm{Ad}_{g^{-1}}(\rho_1)$. 
Combining the resulting inclusion with the equality $R_rg=R_{rw}$, one obtains the inclusion 
$$
(g^{-1}C)\cdot \mathrm C_{\hat V}(e_0)\cdot R_{rw}\subset \mathrm C_{21}(\rho_0). 
$$
which by Lem. \ref{lem:33to31} implies 
\begin{equation}\label{31:inclusion:in:MAIN}
(g^{-1}C)\mathrm{C}_{\hat V}(e_0)rw\subset 
\begin{pmatrix}
    e_0-f_0\\e_1\\0
\end{pmatrix}\mathrm{C}_{\hat V}(e_0)+\begin{pmatrix}
    1\\0\\1
\end{pmatrix}\mathrm{C}_{\hat V}(e_0)   
\end{equation}
(inclusion of subsets of $M_{3,1}\hat V$). Taking the image of this inclusion by the map 
$M_{3,1}\hat V\to \hat V$, $X\mapsto \begin{pmatrix}
    0&0&1
\end{pmatrix}X$, and using $\begin{pmatrix}
    0&0&1
\end{pmatrix}(g^{-1}C)=w^{-1}t$, one obtains
$$
w^{-1}t\cdot\mathrm{C}_{\hat V}(e_0)\cdot rw\subset\mathrm{C}_{\hat V}(e_0). 
$$
which implies $w^{-1}trw\in \mathrm{C}_{\hat V}(e_0)$, which since $r\in\hat V^\times$ can be expressed as follows
$$
(rw)^{-1}rt(rw)\in \mathrm{C}_{\hat V}(e_0). 
$$
The combination of this relation, of the relation $rt\in \mathrm C_{\hat V}(e_0)\cap (e_0+f_\infty+F^2\hat V)$
(see \eqref{rappel:resultats:C:rt}) and of Lem. \ref{lem:9:12:conj:in:CVe0} then implies
\begin{equation}\label{rel:rw}
rw\in \mathrm{C}_{\hat V}(e_0)^\times. 
\end{equation}
The combination of this relation with the second part of \eqref{rappel:resultats:C:rt} implies 
\begin{equation}\label{belonging:w-1t}
    w^{-1}t\in \mathrm{C}_{\hat V}(e_0). 
\end{equation}
Relation \eqref{rel:rw} also implies 
$(rw)^{-1}\in \mathrm{C}_{\hat V}(e_0)$, which implies $g^{-1}C=(g^{-1}C)(rw)^{-1}rw\in (g^{-1}C)\mathrm{C}_{\hat V}(e_0)rw$, 
which together with \eqref{31:inclusion:in:MAIN} yields  
\begin{equation}\label{àvoir}
g^{-1}C\in
\begin{pmatrix}
    e_0-f_0\\e_1\\0
\end{pmatrix}\mathrm{C}_{\hat V}(e_0)+\begin{pmatrix}
    1\\0\\1
\end{pmatrix}\mathrm{C}_{\hat V}(e_0).     
\end{equation}
The equality $g^{-1}C=\begin{pmatrix}
f_1s-v_1w^{-1}t    \\ e_1s-v_2w^{-1}t \\w^{-1}t 
\end{pmatrix}$ then implies the existence of $c,d\in \mathrm{C}_{\hat V}(e_0)$, such that 
\begin{equation}\label{some:step}
\begin{pmatrix}
 f_1s-v_1w^{-1}t    \\ e_1s-v_2w^{-1}t \\w^{-1}t
\end{pmatrix}=\begin{pmatrix}
    e_0-f_0\\e_1\\0
\end{pmatrix}c+\begin{pmatrix}
    1\\0\\1
\end{pmatrix}d.    
\end{equation}
The third component of this equality implies $d=w^{-1}t$, which when plugged in 
the two first components of \eqref{some:step} yields
\begin{equation}\label{next:step}
\begin{pmatrix}
 f_1s-v_1w^{-1}t    \\ e_1s-v_2w^{-1}t 
\end{pmatrix}=\begin{pmatrix}
    e_0-f_0\\e_1
\end{pmatrix}c+\begin{pmatrix}
    1\\0
\end{pmatrix}w^{-1}t.
\end{equation}
The second component of \eqref{next:step} implies $v_2w^{-1}t=e_1\cdot (s-c)$, therefore 
$$
v_2w^{-1}t\in e_1\hat V, 
$$
while the invertibility of $r,w$ and the second part of \eqref{rappel:resultats:C:rt} imply 
\begin{equation}\label{conds:w-1t}
    \epsilon(rw)\in\mathbf k^\times, \quad w^{-1}t=\epsilon(rw)^{-1}(e_0+f_\infty)+F^2\hat V, 
\end{equation}
One has therefore 
$$
\epsilon(rw)v_2w^{-1}t \in e_1\hat V, \quad \epsilon(rw)w^{-1}t\in (e_0+f_\infty)+F^2\hat V, 
$$
which by Lem. \ref{lem:9:13}(b) (applied with $v:=v_2,\alpha:=\epsilon(rw)w^{-1}t$) implies $v_2\in e_1\hat V$, so that there exists 
$\gamma\in\hat V$ such that 
\begin{equation}\label{eq:v2}
v_2=e_1\gamma. 
\end{equation}
Combining this equality with the equality $v_2w^{-1}t=e_1\cdot (s-c)$
arising from the second component of \eqref{next:step}, and using the injectivity of the 
endomorphism $x\mapsto e_1x$ of $\hat V$, one obtains 
\begin{equation}\label{eq:c:s}
c=s-\gamma w^{-1}t.     
\end{equation}
 Combining this equality with the first component of \eqref{next:step} then yields
$$
 f_1s-v_1w^{-1}t =(e_0-f_0)(s-\gamma w^{-1}t)+w^{-1}t
$$
which implies 
\begin{equation}\label{main:final:equation}
    \big( -(v_1+1)+(e_0-f_0)\gamma \big) w^{-1}t =(e_0+f_\infty)s
\end{equation}
It follows from \eqref{epsilon:v1+1} that $\epsilon(v_1+1)=\epsilon(w)\in\mathbf k^\times$, and since 
$\epsilon((e_0-f_0)\gamma)=0$, this implies 
\begin{equation}\label{-v1-1-smth:invertible}
 -(v_1+1)+(e_0-f_0)\gamma\in\hat V^\times.    
\end{equation}
The combination of \eqref{belonging:w-1t}, 
\eqref{conds:w-1t}, \eqref{main:final:equation} and \eqref{-v1-1-smth:invertible} 
implies that one may apply Lem. \ref{lemma:in:hatV} with $\alpha,\tau$ from this statement
respectively equal to $ \epsilon(rw)^{-1}(-(v_1+1)+(e_0-f_0)\gamma) $ and $\epsilon(rw)w^{-1}\tau$, 
and obtain 
\begin{equation}\label{eq:s:in:comm}
s\in \mathrm C_{\hat V}(e_0)    
\end{equation}
and the existence of $\tilde\alpha,\tilde\tau\in \mathrm{C}_{\hat V}(e_0)$ and $\tilde\gamma\in \mathrm{C}_{\hat V}(e_0)^\times$
with $\epsilon(\tilde\gamma)=1$ such that 
$$
-(v_1+1)+(e_0-f_0)\gamma =\epsilon(-(v_1+1))(1+(e_0+f_\infty)\tilde\alpha)\tilde\gamma^{-1}
$$
which by $\epsilon(v_1+1)=\epsilon(w)$ implies
$$
v_1=(e_0-f_0)\gamma +\epsilon(w)(1+(e_0+f_\infty)\tilde\alpha)\tilde\gamma^{-1}-1 
$$
therefore
\begin{equation}\label{eq:v1}
v_1=f_1\gamma+(e_0+f_\infty)(\gamma +\epsilon(w)\tilde\alpha\tilde\gamma^{-1})
+\epsilon(w)\tilde\gamma^{-1}-1     
\end{equation}
The combination of \eqref{eq:c:s}, $c\in \mathrm C_{\hat V}(e_0)$ and \eqref{eq:s:in:comm} yields
$\gamma w^{-1}t\in \mathrm C_{\hat V}(e_0)$, which together with \eqref{belonging:w-1t} and  \eqref{conds:w-1t}, 
and by Lem. \ref{lem:quotients}(b), implies 
\begin{equation}\label{belonging:gamma!}
\gamma\in \mathrm C_{\hat V}(e_0). 
\end{equation}
One also has 
\begin{equation}\label{eq:w}
w=\epsilon(w)^{-1}w\tilde\gamma\cdot \epsilon(w)\tilde\gamma^{-1}
\end{equation}
Set $\boldsymbol\gamma:=\gamma$, $\boldsymbol{\delta}:=\epsilon(w)^{-1}w\tilde\gamma$, $\mathbf c:=\gamma +\epsilon(w)\tilde\alpha\tilde\gamma^{-1}$, 
$\mathbf d:= \epsilon(w)\tilde\gamma^{-1}$. Then 
\begin{equation}\label{crucial:belonging}
(\boldsymbol\gamma,\boldsymbol\delta)\in\hat V\times\hat V^\times, \quad \epsilon(\boldsymbol\delta)=1, \quad 
(\mathbf c,\mathbf d)\in \mathrm{C}_{\hat V}(e_0)\times\mathrm{C}_{\hat V}(e_0)^\times,
\end{equation}
where all the relations follow from definitions, the equality also follows from 
 $\epsilon(\tilde\gamma)=1$, and the relation 
$\mathbf c\in \mathrm{C}_{\hat V}(e_0)$ also follows from \eqref{belonging:gamma!}. 
Then \eqref{eq:v1}, \eqref{eq:v2}, \eqref{eq:w} imply the equality 
$$
(v_1,v_2,w)=(f_1\boldsymbol\gamma+(e_0+f_\infty)\mathbf c+\mathbf d-1,e_1\boldsymbol\gamma,\boldsymbol\delta \mathbf d)
$$
which together with \eqref{crucial:belonging}, \eqref{comp:product:subgroups} and the equality $g=u(v_1,v_2,w)$
implies $g\in (U^{(0)}\cap\mathrm C_{21}(\rho_1)^\times)\cdot (U^{(0)}\cap\mathrm C_{21}(\rho_0)^\times)$. The statement 
follows. 
\end{proof}

\subsection{Local injectivity of the morphism (D)}\label{sect:9:4}

\begin{lem}\label{computation:kernel:of:group:morphism}
    One has $\mathrm{ker}(\mathrm C_{21}^{(0)}(\rho_1)^\times\to\mathrm C_2(\overline\rho_1)^\times)
=\mathrm C_{21}^{(0)}(\rho_1)^\times\cap U^{(0)}$ (equality of subgroups of $\mathrm C_{21}^{(0)}(\rho_1)^\times$). 
\end{lem}

\begin{proof}
One has 
$$
\mathrm{ker}(\mathrm C_{21}^{(0)}(\rho_1)^\times\to\mathrm C_2(\overline\rho_1)^\times)
=\mathrm C_{21}^{(0)}(\rho_1)^\times\cap \mathrm{ker}((T_{21}\hat V)^\times\to\mathrm{GL}_2\hat V)=
\mathrm C_{21}^{(0)}(\rho_1)^\times\cap U=\mathrm C_{21}^{(0)}(\rho_1)^\times\cap U^{(0)}, 
$$
where the last equality follows from 
$\mathrm C_{21}^{(0)}(\rho_1)^\times=\mathrm C_{21}(\rho_1)^\times\cap U^{(0)}\subset U^{(0)}$. 
\end{proof}

\begin{cor}\label{loc:inj:U}
    The map $\mathrm{Hom}^{1,((0)),\bullet}_{\mathcal C\operatorname{-alg}}
(\hat{\mathcal V},T_{21}\hat V)\to \mathrm{Hom}^{1,(0)}_{\mathcal C\operatorname{-alg}}
(\hat{\mathcal V},M_2\hat V)$ is invariant under the action on the source of the subgroup 
$\mathrm{ker}(\mathrm C_{21}^{(0)}(\rho_1)^\times\to\mathrm C_2(\overline\rho_1)^\times)$ of $\mathrm C_{21}^{(0)}(\rho_1)^\times$, 
    and the resulting morphism of pointed sets  
\begin{equation}\label{16oct2025}
    \mathrm{ker}(\mathrm C_{21}^{(0)}(\rho_1)^\times\to\mathrm C_2(\overline\rho_1)^\times)\backslash\mathrm{Hom}^{1,((0)),\bullet}_{\mathcal C\operatorname{-alg}}
(\hat{\mathcal V},T_{21}\hat V)\to 
\mathrm{Hom}^{1,(0)}_{\mathcal C\operatorname{-alg}}
(\hat{\mathcal V},M_2\hat V)
\end{equation}
(the distinguished elements being the class of $\rho_{\mathrm{DT}}$ in the source and $\overline\rho_{\mathrm{DT}}$ in 
the target) is locally injective. 
\end{cor}

\begin{proof} 
The sequence of maps 
$$
\mathrm{Hom}^{1,((0)),\bullet}_{\mathcal C\operatorname{-alg}}
(\hat{\mathcal V},T_{21}\hat V)\hookrightarrow
\mathrm{Hom}^{1,((0))}_{\mathcal C\operatorname{-alg}}
(\hat{\mathcal V},T_{21}\hat V)\to \mathrm{Hom}^{1,(0)}_{\mathcal C\operatorname{-alg}}
(\hat{\mathcal V},M_2\hat V)
$$
is compatible with the sequence of group morphisms 
$\mathrm C_{21}^{(0)}(\rho_1)^\times=\mathrm C_{21}^{(0)}(\rho_1)^\times\to\mathrm C_2(\overline\rho_1)^\times$, 
where the second morphism is induced by the morphism
$(T_{21}\hat V)^\times\to\mathrm{GL}_2\hat V$, $x\mapsto\overline x$, and therefore 
with the sequence of group morphisms 
$\mathrm{ker}(\mathrm C_{21}^{(0)}(\rho_1)^\times\to\mathrm C_2(\overline\rho_1)^\times)=
\mathrm{ker}(\mathrm C_{21}^{(0)}(\rho_1)^\times\to\mathrm C_2(\overline\rho_1)^\times)\to1$, which 
proves the first statement. 
This induces the sequence of maps
\begin{align*}    
&\mathrm{ker}(\mathrm C_{21}^{(0)}(\rho_1)^\times\to\mathrm C_2(\overline\rho_1)^\times)
\backslash\mathrm{Hom}^{1,((0)),\bullet}_{\mathcal C\operatorname{-alg}}
(\hat{\mathcal V},T_{21}\hat V)\to 
\mathrm{ker}(\mathrm C_{21}^{(0)}(\rho_1)^\times\to\mathrm C_2(\overline\rho_1)^\times)
\backslash\mathrm{Hom}^{1,((0))}_{\mathcal C\operatorname{-alg}}
(\hat{\mathcal V},T_{21}\hat V)
\\& \to 
\mathrm{Hom}^{1,(0)}_{\mathcal C\operatorname{-alg}}
(\hat{\mathcal V},M_2\hat V)
\end{align*} 
between quotient spaces, whose composition is \eqref{16oct2025}.
It follows from Lem. \ref{computation:kernel:of:group:morphism}
that this can be identified with a map
$$
 (U^{(0)}\cap \mathrm C_{21}^{(0)}(\rho_1)^\times)\backslash\mathrm{Hom}^{1,((0)),\bullet}_{\mathcal C\operatorname{-alg}}
(\hat{\mathcal V},T_{21}\hat V)\to 
\mathrm{Hom}^{1,(0)}_{\mathcal C\operatorname{-alg}}
(\hat{\mathcal V},M_2\hat V)
$$
such that $(U^{(0)}\cap \mathrm C_{21}^{(0)}(\rho_1)^\times)\cdot\rho_{\mathrm{DT}}\mapsto \overline\rho_{\mathrm{DT}}$; 
let us prove the local injectivity of this map of pointed sets. 


Let $\alpha\in(U^{(0)}\cap \mathrm C_{21}^{(0)}(\rho_1)^\times)\backslash\mathrm{Hom}^{1,((0)),\bullet}_{\mathcal C\operatorname{-alg}}
(\hat{\mathcal V},T_{21}\hat V)$ belong to the preimage of $\overline\rho_{\mathrm{DT}}$. Let 
$\rho$ be a representative of $\alpha$, then 
$\overline\rho=\overline\rho_{\mathrm{DT}}$, which by Lem. \ref{lem:9:8:toto}(b) implies $\rho\in U^{(0)}\cdot\rho_{\mathrm{DT}}$. Let $u\in U^{(0)}$
be such that $\rho=u\cdot\rho_{\mathrm{DT}}$, then Lem. \ref{lem:MAIN} implies the existence of $u_1\in U^{(0)}\cap \mathrm{C}_{21}(\rho_1)^\times$ 
and $u_0\in U^{(0)}\cap \mathrm{C}_{21}(\rho_0)^\times$, such that 
$u=u_1u_0$. One has $u_0\cdot\rho_{\mathrm{DT}}=\rho_{\mathrm{DT}}$ by Lem. \ref{lem:comm:diag:toto}(a), 
therefore $\rho=u\cdot \rho_{\mathrm{DT}}=u_1\cdot\rho_{\mathrm{DT}}=u_1\bullet\rho_{\mathrm{DT}}$. 

Then 
$$
\rho=u\cdot \rho_{\mathrm{DT}}=u_1\cdot\rho_{\mathrm{DT}}=u_1\bullet\rho_{\mathrm{DT}}
\in (U^{(0)}\cap \mathrm{C}_{21}(\rho_0)^\times)\bullet \rho_{\mathrm{DT}}, 
$$
where the second equality follows from $u_0\cdot\rho_{\mathrm{DT}}=\rho_{\mathrm{DT}}$ 
(see Lem. \ref{lem:comm:diag:toto}(a)) and the last equality follows from 
$u_1\in U^{(0)}\cap \mathrm{C}_{21}(\rho_0)^\times=U\cap (T_{21}^{(0)}\hat V)^\times\cap\mathrm{C}_{21}(\rho_1)^\times
=U\cap \mathrm{C}_{21}^{(0)}(\rho_1)^\times\subset \mathrm{C}_{21}^{(0)}(\rho_1)^\times$ 
and \eqref{relation:actions}. This implies 
$\alpha=(U^{(0)}\cap \mathrm{C}_{21}(\rho_0)^\times)\bullet \rho_{\mathrm{DT}}$, and therefore the claimed
local injectivity. 
\end{proof}

\begin{lem}\label{lem:loc:inj:succ:quotients}
    Let $(G,X,x_0)\to(H,Y,y_0)$ be a morphism of pointed sets with group actions. It gives rise to 
    morphism of pointed sets $(K\backslash X,Kx_0)\to (Y,y_0)$ and $(G\backslash X,Gx_0)\to (H\backslash Y,Hy_0)$, 
    where $K\subset G$ be the kernel of the group morphism $G\to H$. 
    
    If the group morphism $G\to H$ is surjective, and if the morphism of pointed sets $(K\backslash X,Kx_0)\to (Y,y_0)$ is locally 
    injective, and  then so is the morphism of pointed sets $(G\backslash X,Gx_0)\to (H\backslash Y,Hy_0)$. 
\end{lem}

\begin{proof}
Denote by $f : X\to Y$, $\phi : G\to H$ the set and group morphisms underlying the morphism $(G,X,x_0)\to(H,Y,y_0)$. 
Let $\alpha\in G\backslash X$ belong to the preimage of $Hy_0$ by the map $G\backslash X\to H\backslash Y$. If $x\in X$ is a representative of $\alpha$, then 
$f(x)\in Hy_0$. Let $h\in H$ be such that $f(x)=hy_0$. By the surjectivity of $\phi$, there exists $g\in G$ such that $\phi(g)=h$. Then 
$f(g^{-1}x)=h^{-1}f(x)=y_0$. This implies that the class $Kg^{-1}x$ of $g^{-1}x$ in $K\backslash X$ belongs to the preimage of $y_0$ by the map 
$K\backslash X\to Y$. The local injectivity of the morphism of pointed sets $(K\backslash X,Kx_0)\to (Y,y_0)$ then implies $Kg^{-1}x=Kx_0$. Therefore there exists 
$k\in K$ such that $g^{-1}x=kx_0$. Then $x=gkx_0\in Gx_0$, which implies $\alpha=Gx_0$. This implies the claimed local injectivity. 
\end{proof}

\begin{prop}\label{prop:D:loc:inj}
    The morphism (D) of pointed sets given by 
    $$
    (\mathrm C_{21}^{(0)}(\rho_1)^\times\backslash\mathrm{Hom}^{1,((0)),\bullet}_{\mathcal C\operatorname{-alg}}
(\hat{\mathcal V},T_{21}\hat V),\mathrm C_{21}^{(0)}(\rho_1)^\times\bullet \rho_{\mathrm{DT}}
)\to 
(\mathrm C_2({\overline\rho}_1)^\times\backslash\mathrm{Hom}^{1,(0)}_{\mathcal C\operatorname{-alg}}(\hat{\mathcal V},M_2\hat V),\mathrm C_2({\overline\rho}_1)^\times\bullet \overline\rho_{\mathrm{DT}}
)
$$ 
is locally injective. 
\end{prop}

\begin{proof}
This follows by applying Lem. \ref{lem:loc:inj:succ:quotients} to the morphism of pointed sets with group actions 
$(\mathrm C_{21}(\rho_1)^\times,\mathrm{Hom}^{1,(0),\bullet}_{\mathcal C\operatorname{-alg}}
(\hat{\mathcal V},T_{21}\hat V),\rho_{\mathrm{DT}})\to 
(\mathrm C_2({\overline\rho}_1)^\times,\mathrm{Hom}^{1,(0)}_{\mathcal C\operatorname{-alg}}(\hat{\mathcal V},M_2\hat V),\overline\rho_{\mathrm{DT}})$ 
defined in Lem. \ref{lem:6:20}(b), using Lem. \ref{lem:c:rho1:surj} 
and Cor. \ref{loc:inj:U}. 
\end{proof}


\section{Local injectivity of the morphism (E)}\label{sec 10}

This section is devoted to the proof of the local injectivity of the map (E). This map is expressed as the composition 
of various maps (E1)-(E5) as follows 
\begin{equation}\label{diag:E}
\xymatrix{
\mathrm C_2(\overline\rho_1)^\times\backslash\mathrm{GL}_2\hat V/\mathrm C_2(\overline\rho_0)^\times\ar_{(E1)}[dd]
\ar^{\!\!\!\!(E2)}[r]&
\mathbf k[[u,v]]^\times\backslash(M_{12}F^1\hat V\times M_{21}\hat V)/\mathrm C_2(\overline\rho_0)^\times\ar^{\ \ \ \ \ \ (E3)}[r]&
\mathbf k[[u,v]]^\times\backslash F^1\hat{\mathbf C}\ar^{(E4)}[d]\\
&&
\mathbf k[[u,v]]^\times\backslash\prod_{n\geq1}F^n\hat V\ar^{(E5)}[d]\\
\mathrm C_2(\overline\rho_1)^\times\backslash\mathrm{Hom}_{\mathcal C\operatorname{-alg}}^{1,(0)}(\hat{\mathcal V},M_2\hat V)
\ar_{(E)}[rr]&&
\mathbf k[[u,v]]^\times\backslash\mathrm{Hom}_{\mathcal C\operatorname{-alg}}(\hat{\mathcal W},\hat V)}
\end{equation}
The map (E1) is introduced and proved to be bijective in §\ref{sect:10:1}. The map (E2) is introduced in §\ref{sect:10:2} 
and proved there to be injective, using in particular the computation of the commutant 
of $\overline\rho_1$ in $M_2\hat V$ from §\ref{sec 6}. The construction of the map (E3), based on a $(V,V)$-bimodule $\mathbf C$, 
is done in  §\ref{sect:10:3}. The map (E4), which involves a map from $\mathbf C$ to a space of maps $\{n|n\geq1\}\to\hat V$, 
is defined in §\ref{sect:10:4}. The map (E5), which identifies the latter space of maps with a set of algebra morphisms 
$\hat{\mathcal W}\to\hat V$, is defined and shown to be bijective in §\ref{sect:10:5}. §\ref{sect:10:6} is devoted to a study 
of $\mathbf C$, which then enables us to prove the local injectivity of (E4) in §\ref{sect:10:7}. The local injectivity of 
(E3) is then proved in §\ref{sect:10:8}. The whole material of the section is used in §\ref{sect:10:8} to derive the local 
injectivity of (E).

\subsection{The map (E1)}\label{sect:10:1}

\begin{lem}\label{lem:basic:E1}
(a) For $P\in \mathrm{GL}_2\hat V$ and $\sigma\in \mathrm{Hom}_{\mathcal C\operatorname{-alg}}^{1,(0)}(\hat{\mathcal V},M_2\hat V)$, 
define $P\cdot\sigma : \hat{\mathcal V}\to M_2\hat V$ to be the topological $\mathbf k$-algebra morphism such that $e_1\mapsto \overline\rho_1=\sigma(e_1)$ 
and $e_0\mapsto \mathrm{Ad}_P(\sigma(e_0))$. Then $(P,\sigma)\mapsto P\cdot\sigma$ defines an action of the group $\mathrm{GL}_2\hat V$
on the set $\mathrm{Hom}_{\mathcal C\operatorname{-alg}}^{1,(0)}(\hat{\mathcal V},M_2\hat V)$. 

(b) The action from (a) is transitive, the isotopy group of the element $\overline\rho_{\mathrm{DT}}$ is the subgroup 
$\mathrm C_2(\overline\rho_0)^\times$, so that the map $P\mapsto P\cdot\overline\rho_{\mathrm{DT}}$ induces a bijection 
$\mathrm{GL}_2\hat V/\mathrm C_2(\overline\rho_0)^\times\to\mathrm{Hom}_{\mathcal C\operatorname{-alg}}^{1,(0)}(\hat{\mathcal V},M_2\hat V)$. 

(c) The bijection from (b) is equivariant under the action of $\mathrm C_2(\overline\rho_1)^\times$, the action on the source being by 
left multiplication and its action on the target being $(P,\sigma)\mapsto P\bullet\sigma$ (see Lem. \ref{lem:actions:2912}). 
\end{lem}

\begin{proof}
(a) The map from $\mathrm{Hom}_{\mathcal C\operatorname{-alg}}^{1,(0)}(\hat{\mathcal V},M_2\hat V)$ to the $\mathrm{GL}_2\hat V$-conjugacy 
class of $\overline\rho_0$ in $M_2\hat V$ given by $\sigma\mapsto \sigma(e_0)$ is a bijection. The said map $(P,\sigma)\mapsto P\cdot\sigma$ corresponds
under this bijection to the conjugation action of $\mathrm{GL}_2\hat V$ on this conjugacy class; it follows that this is an action. 

(b) The first statement follows from (a) and from the transitivity of the action of $\mathrm{GL}_2\hat V$ on the 
$\mathrm{GL}_2\hat V$-conjugacy class of $\overline\rho_0$ in $M_2\hat V$. The second statement follows from the fact that the bijection from (a) 
takes $\overline\rho_{\mathrm{DT}}$ to $\overline\rho_0$ and from the fact that the isotropy group of $\overline\rho_0$ for the conjugation action is 
$\mathrm C_2(\overline\rho_0)^\times$. 

(c) For $c\in \mathrm C_2(\overline\rho_1)^\times$ and $P\in \mathrm{GL}_2\hat V$, one has 
$c\bullet(P\cdot\overline\rho_{\mathrm{DT}})=\mathrm{Ad}_c \circ (P\cdot\overline\rho_{\mathrm{DT}})=
\mathrm{Ad}_c \circ (e_1\mapsto\overline\rho_1,e_0\mapsto \mathrm{Ad}_P(\overline\rho_0))
=(e_1\mapsto \mathrm{Ad}_c(\overline\rho_1),e_0\mapsto \mathrm{Ad}_c(\mathrm{Ad}_P(\overline\rho_0)))
=(e_1\mapsto \overline\rho_1,e_0\mapsto \mathrm{Ad}_{cP}(\overline\rho_0))
=cP\cdot\overline\rho_{\mathrm{DT}}$, where $(e_0\mapsto \alpha,e_1\mapsto\beta)$ denote the algebra morphism 
$\hat{\mathcal V}\to M_2\hat V$ such that $e_0\mapsto \alpha,e_1\mapsto\beta$ and where the third equality follows from 
$c\in \mathrm C_2(\overline\rho_1)^\times$. This implies the $\mathrm C_2(\overline\rho_1)^\times$-equivariance of the map 
$\mathrm{GL}_2\hat V\to\mathrm{Hom}_{\mathcal C\operatorname{-alg}}^{1,(0)}(\hat{\mathcal V},M_2\hat V)$, 
$P\mapsto P\cdot \overline\rho_{\mathrm{DT}}$, which implies the result. 
\end{proof}

\begin{defn}\label{defE1}
    The bijection 
    $$
    \mathrm C_2(\overline\rho_1)^\times\backslash\mathrm{GL}_2\hat V/\mathrm C_2(\overline\rho_0)^\times\to
    \mathrm C_2(\overline\rho_1)^\times\backslash\mathrm{Hom}_{\mathcal C\operatorname{-alg}}^{1,(0)}(\hat{\mathcal V},M_2\hat V)
    $$
    between coset spaces arising from Lem. \ref{lem:basic:E1}(c) is denoted (E1). 
\end{defn}

\subsection{The map (E2)}\label{sect:10:2}

\begin{lem}\label{lem:10:3:7mars}
(a) The assignment
$$
(\phi,(\mathrm{row},\mathrm{col}))\mapsto \phi\bullet(\mathrm{row},\mathrm{col}):=(\phi(e_1,f_1)\cdot \mathrm{row},\mathrm{col}\cdot \phi^{-1}(e_1,f_1))
$$
defines an action of the group $\mathbf k[[u,v]]^\times$ on the set $M_{12}F^1\hat V\times M_{21}\hat V$. 

(b) The map 
$$
\mathrm{GL}_2\hat V\to M_{12}F^1\hat V\times M_{21}\hat V,\quad 
g\mapsto (\overline{\mathrm{row}}_{\mathrm{DT}}\cdot g,g^{-1}\cdot \overline{\mathrm{col}}_{\mathrm{DT}})
$$
builds up, together with the group morphism $\mathrm C_2(\overline\rho_1)^\times\to \mathbf k[[u,v]]^\times$ induced by Lem. 
\ref{lem:6:16:2912}(c), a morphism of sets with group actions
$$
(\mathrm C_2(\overline\rho_1)^\times,\mathrm{GL}_2\hat V,\cdot)\to (\mathbf k[[u,v]]^\times,
M_{12}F^1\hat V\times M_{21}\hat V,\bullet),  
$$
where the group action in the source is induced by the left multiplication of 
$\mathrm C_2(\overline\rho_1)^\times$ on $\mathrm{GL}_2\hat V$, 
and the group action in the target is induced by (a). 

(c) The map between coset spaces 
$$
\mathrm C_2(\overline\rho_1)^\times\backslash\mathrm{GL}_2\hat V\to\mathbf k[[u,v]]^\times\backslash(M_{12}F^1\hat V\times M_{21}\hat V)
$$
induced by the morphism of sets with group actions from (b) is injective. 
\end{lem}

\begin{proof}
(a) follows from the equality $\psi\bullet(\phi\bullet(\mathrm{row},\mathrm{col}))=\psi\bullet
(\phi(e_1,f_1)\cdot \mathrm{row},\mathrm{col}\cdot \phi^{-1}(e_1,f_1))
=(\psi(e_1,f_1)\cdot(\phi(e_1,f_1)\cdot \mathrm{row}),(\mathrm{col}\cdot \phi^{-1}(e_1,f_1))\cdot\psi^{-1}(e_1,f_1))
=(\psi\phi(e_1,f_1)\cdot \mathrm{row},\mathrm{col}\cdot (\psi\phi)^{-1}(e_1,f_1))
=\psi\phi\bullet(\mathrm{row},\mathrm{col})$ for any $\phi,\psi\in\mathbf k[[u,v]]^\times$.     

(b) Let $c\in \mathrm C_2(\overline\rho_1)^\times$, then by Lem. \ref{lem:comp:C2:rho1:times} there exists a pair $(\phi,v)\in 
\mathbf k[[u,v]]^\times\times\hat V$ such that $c=\overline M(\phi,v)$, and the image of $c$ by 
$\mathrm C_2(\overline\rho_1)^\times\to \mathbf k[[u,v]]^\times$ is $c$. Then 
\begin{equation}\label{identité:8mars:a}
\overline{\mathrm{row}}_{\mathrm{DT}}\cdot c=\overline{\mathrm{row}}_{\mathrm{DT}}\cdot \overline M(\phi,v)
=\begin{pmatrix}
e_1&-f_1
\end{pmatrix}\cdot(\phi(e_1,f_1)I_2+\begin{pmatrix}
    f_1\\e_1
\end{pmatrix}v\begin{pmatrix}
    1&1
\end{pmatrix})=\phi(e_1,f_1)\begin{pmatrix}
e_1&-f_1
\end{pmatrix}=\phi(e_1,f_1)\overline{\mathrm{row}}_{\mathrm{DT}}.     
\end{equation}
and 
$$
c\cdot\overline{\mathrm{col}}_{\mathrm{DT}}=\overline M(\phi,v)\cdot\overline{\mathrm{col}}_{\mathrm{DT}} 
=
(\phi(e_1,f_1)I_2+\begin{pmatrix}
    f_1\\e_1
\end{pmatrix}v\begin{pmatrix}
    1&1
\end{pmatrix})\cdot \begin{pmatrix}
1\\-1
\end{pmatrix}=\begin{pmatrix}
1\\-1
\end{pmatrix}\phi(e_1,f_1)=\overline{\mathrm{col}}_{\mathrm{DT}}\phi(e_1,f_1),  
$$
where the equalities follow from Lem. \ref{lem:6:16:2912}(b) and Def. \ref{def:barcol:barrow}; 
since the morphism $\mathrm C_2(\overline\rho_1)^\times\to \mathbf k[[u,v]]^\times$ is such that 
$c^{-1}\mapsto \phi^{-1}$, one also has 
\begin{equation}\label{identité:8mars:b}
c^{-1}\cdot\overline{\mathrm{col}}_{\mathrm{DT}}=\overline M(\phi,v)\cdot\overline{\mathrm{col}}_{\mathrm{DT}} 
=\overline{\mathrm{col}}_{\mathrm{DT}}\phi^{-1}(e_1,f_1). 
\end{equation} 
Then for any $P\in\mathrm{GL}_2\hat V$ and $c\in \mathrm C_2(\overline\rho_1)^\times$, one has 
$$
(\overline{\mathrm{row}}_{\mathrm{DT}}\cdot cP,(cP)^{-1}\cdot \overline{\mathrm{col}}_{\mathrm{DT}})
=
(\phi(e_1,f_1)\overline{\mathrm{row}}_{\mathrm{DT}}\cdot P,P^{-1}\cdot \overline{\mathrm{col}}_{\mathrm{DT}}
(\phi^{-1}(e_1,f_1))
=\phi(e_1,f_1)\bullet 
(\overline{\mathrm{row}}_{\mathrm{DT}}\cdot P,P^{-1}\cdot \overline{\mathrm{col}}_{\mathrm{DT}})
$$
where the equalities follow from \eqref{identité:8mars:a} and \eqref{identité:8mars:b}, which implies the statement.  

(c) Let $P,Q\in\mathrm{GL}_2\hat V$ be such that their images by the map from (b) are related by the action of 
$\mathbf k[[u,v]]^\times$. 
Let then $\phi\in\mathbf k[[u,v]]^\times$ be such that 
$$
(\overline{\mathrm{row}}_{\mathrm{DT}}\cdot Q,Q^{-1}\cdot \overline{\mathrm{col}}_{\mathrm{DT}})
=\phi\bullet(\overline{\mathrm{row}}_{\mathrm{DT}}\cdot P,P^{-1}\cdot \overline{\mathrm{col}}_{\mathrm{DT}}),  
$$
i.e. $(\overline{\mathrm{row}}_{\mathrm{DT}}\cdot Q,Q^{-1}\cdot \overline{\mathrm{col}}_{\mathrm{DT}})
=(\phi(e_1,f_1)\overline{\mathrm{row}}_{\mathrm{DT}}\cdot P,P^{-1}\cdot \overline{\mathrm{col}}_{\mathrm{DT}}\phi(e_1,f_1)^{-1})$. 
Then if $x:=\phi(e_1,f_1)^{-1}PQ^{-1}\in\mathrm{GL}_2\hat V$, one obtains
$$
\overline{\mathrm{row}}_{\mathrm{DT}}\cdot x=\overline{\mathrm{row}}_{\mathrm{DT}},\quad 
x\cdot \overline{\mathrm{col}}_{\mathrm{DT}}=\overline{\mathrm{col}}_{\mathrm{DT}}, 
$$
therefore $x-I_2\in\mathrm{Ann}(\overline{\mathrm{row}}_{\mathrm{DT}},\overline{\mathrm{col}}_{\mathrm{DT}})$, which by Lem. \ref{lem:van:0301:BIS}
implies the existence of $v\in\hat V$ such that 
$$
x-I_2=\begin{pmatrix}
    f_1\\e_1
\end{pmatrix}v\begin{pmatrix}
    1&1
\end{pmatrix}. 
$$
It follows that $x=\overline M(1,v)$ (see Lem. \ref{lem:6:16:2912}(b)), which implies the second equality in 
$PQ^{-1}=\phi(e_1,f_1)x=\phi(e_1,f_1)\overline M(1,v)
=\overline M(\phi,\phi(e_1,f_1)v)\in \mathrm C_2(\overline\rho_1)^\times$, where the first equality follows from the definition of $x$. 
It follows that the classes of $P,Q$ in $\mathrm C_2(\overline\rho_1)^\times\backslash\mathrm{GL}_2\hat V$ are equal, which proves the statement. 
\end{proof}

\begin{lem}\label{lem:10:4:7mars}
(a) The assignment
$$
((\mathrm{row},\mathrm{col}),c)\mapsto (\mathrm{row},\mathrm{col})\bullet c:=(\mathrm{row}\cdot c,c^{-1}\cdot \mathrm{col})
$$
defines a right action of the group $\mathrm C_2(\overline\rho_0)^\times$ on the set $M_{12}F^1\hat V\times M_{21}\hat V$, which 
commutes with the left action of $\mathbf k[[u,v]]^\times$ from Lem. \ref{lem:10:3:7mars}(a), 
and therefore induces a right action of $\mathrm C_2(\overline\rho_0)^\times$ on the coset space 
$\mathbf k[[u,v]]^\times\backslash(M_{12}F^1\hat V\times M_{21}\hat V)$. 

(b) The map from Lem. \ref{lem:10:3:7mars}(c) is $\mathrm C_2(\overline\rho_0)^\times$-equivariant, the action of 
$\mathrm C_2(\overline\rho_0)^\times$ on the source being by right multiplication and its action on the target being as in (a). 

(c) The map between coset spaces 
$$
\mathrm C_2(\overline\rho_1)^\times\backslash\mathrm{GL}_2\hat V/\mathrm C_2(\overline\rho_0)^\times
\to\mathbf k[[u,v]]^\times\backslash(M_{12}F^1\hat V\times M_{21}\hat V)/\mathrm C_2(\overline\rho_0)^\times
$$
induced by (b) is injective. 
\end{lem}

\begin{proof}
(a) is immediate. (b) follows from the  $\mathrm C_2(\overline\rho_0)^\times$-equivariance of the map from Lem. \ref{lem:10:3:7mars}(b), which follows from 
$$
(\overline{\mathrm{row}}_{\mathrm{DT}}\cdot (Pc),(Pc)^{-1}\cdot \overline{\mathrm{col}}_{\mathrm{DT}})
=
((\overline{\mathrm{row}}_{\mathrm{DT}}\cdot P)\cdot c,c^{-1}\cdot (P^{-1}\cdot \overline{\mathrm{col}}_{\mathrm{DT}}))
$$
for any $c\in\mathrm C_2(\overline\rho_0)^\times,P\in \mathrm{GL}_2\hat V$. 
(c) follows from (b) and Lem. \ref{lem:10:3:7mars}(c).   
\end{proof}

The map from Lem. \ref{lem:10:4:7mars}(c) will be denoted (E2). 

\subsection{The $(V,V)$-bimodule $\mathbf C$ and the map (E3)}\label{sect:10:3}

The following lemma gathers some basic facts on bimodules and bimodule morphisms. 

\begin{lem}\label{lem:10:5}
\begin{itemize}
    \item[(a)] If $A,B$ are $\mathbf k$-algebras and $M,M'$ are $(A,B)$-bimodules, then the set of $(A,B)$-bimodule morphisms $M\to M'$
is a $\mathbf k$-module. 
    \item[(b)] If $A,B,C$ are $\mathbf k$-algebras, $M$ is an $(A,B)$-bimodule and $N$ is a $(B,C)$-bimodule, the $(A,C)$-bimodule 
$M\otimes_B N$ is defined as the cokernel of the $(A,C)$-bimodule morphism $M\otimes B\otimes N\to M\otimes N$, 
$m\otimes b\otimes n\mapsto mb\otimes n-m\otimes bn$. 
In this situation: 
    \begin{itemize}
        \item[(b1)] an algebra morphism $\mathbf u: B'\to B$ induces (by pull-back) an $(A,B')$-bimodule structure on $M$ 
        and a $(B',C)$-bimodule 
    structure on $N$, and therefore an $(A,C)$-bimodule structure $M\otimes_{B'}N$, as well as an $(A,C)$-bimodule morphism 
    $\varphi_{\mathbf u} : M\otimes_{B'}N\to M\otimes_B N$ taking the class of $m\otimes n$ in $M\otimes_{B'}N$ to the class of $m\otimes n$ in 
    $M\otimes_BN$ for $m\in M,n\in N$;  
        \item[(b2)] any element $z$ of the center of $B$ 
        induces an $(A,B)$-bimodule endomorphism $\mu_z$ of $M\otimes_B N$, induced by 
$m\otimes n\mapsto mz\otimes n$.
    \end{itemize}
\end{itemize}    
\end{lem}

\begin{proof}
    Immediate. 
\end{proof}

Let $X$ be one of the subalgebras $\mathrm{C}_V(e_0),\mathrm{C}_V(f_0)$ or $\mathbf k[e_0,f_0]$ of $V$. Then $X$ is a graded 
subalgebra, $V$ being equipped with the total degree (for which $e_0,e_1,f_0,f_1$ all have degree 1).   
Therefore $V\otimes_X V$ is a graded $(V,V)$-bimodule, where the first (resp. second) factor $V$ is viewed as a $(V,X)$-bimodule
(resp. $(X,V)$-bimodule). 

For $M:=\oplus_{n\in\mathbb Z}M_n$ a $\mathbb Z$-graded $\mathbf k$-module, let $M[1]$ be the graded module with 
$M[1]_n:=M_{n+1}$.

\begin{defn}\label{def:CC:3003}
(a) 
Define the graded $(V,V)$-bimodules
\begin{equation}\label{def:SS:TT}
 \mathbf S:=V\otimes_{\mathbf k[e_0,f_0]}V,\quad 
\mathbf T^{\mathbf e}:=V\otimes_{\mathrm{C}_V(e_0)}V,\quad \mathbf T^{\mathbf f}:=V\otimes_{\mathrm{C}_V(f_0)}V, \quad
 \mathbf T:=\mathbf T^{\mathbf e}\oplus \mathbf T^{\mathbf f} \oplus \mathbf S[1]. 
\end{equation}

(b) Let $\mathbf S\to\mathbf T$ be the $(V,V)$-bimodule morphism given by the direct sum 
$(-\varphi_{\mathbf e})\oplus \varphi_{\mathbf f}\oplus \mu_{e_0-f_0}$, where $\varphi_{\mathbf e},\varphi_{\mathbf f}$
are the $(V,V)$-bimodule morphisms from $V\otimes_{\mathbf k[e_0,f_0]}V$ to 
$V\otimes_{\mathrm{C}_V(e_0)}V$ and $V\otimes_{\mathrm{C}_V(f_0)}V$
corresponding (see Lem. \ref{lem:10:5}(b1) above) to the inclusions of algebras $\mathbf e,\mathbf f$ of $\mathbf k[e_0,f_0]$
in $\mathrm{C}_V(e_0)$ and $\mathrm{C}_V(f_0)$, 
where $-\varphi_{\mathbf e}$ is the opposite of $\varphi_{\mathbf e}$ (see Lem. \ref{lem:10:5}(a)), and 
where $\mu_{e_0-f_0}$ is the $(V,V)$-bimodule endomorphism of $V\otimes_{\mathbf k[e_0,f_0]}V$ 
induced by the central element 
$e_0-f_0\in\mathbf k[e_0,f_0]$ (see Lem. \ref{lem:10:5}(b2)); the map $\mathbf S\to\mathbf T$ is therefore such that 
$$
[v\otimes w]_{\mathbf k[e_0,f_0]}\mapsto -[v\otimes w]_{\mathrm C_V(e_0)}+[v\otimes w]_{\mathrm C_V(f_0)}
+[v\otimes (e_0-f_0)w]_{\mathbf k[e_0,f_0]} 
$$
for $v,w\in V$, where $x\mapsto [x]_X$ is the projection map $(\hat V\hat\otimes \hat V)\to (V\otimes_XV)^\wedge$. 
where we identify an element of each of the summands of 
$\mathbf T=\mathbf T^{\mathbf e}\oplus \mathbf T^{\mathbf f} \oplus \mathbf S[1]$
with its image in  $\mathbf T$. 
 
(c) Define the $(V,V)$-bimodule $\mathbf C:=\mathrm{coker}(\mathbf S\to\mathbf T)$, the morphism $\mathbf S\to\mathbf T$ being as in (b). 
\end{defn}

For $M:=\oplus_{n\geq0}M_n$ a $\mathbb Z_+$-graded $\mathbf k$-module, recall that $F^iM=\oplus_{j\geq i}M_j$
for any $i\in\mathbb Z$, $\hat M=\hat\oplus_{n\geq0}M_n$, and $F^i\hat M=\hat\oplus_{j\geq i}M_j$. 

\begin{lem}\label{lem:10:7:1603}
(a) The morphism $\mathbf S\to\mathbf T$ is homogeneous of degree $0$, therefore $\mathbf C$ is a graded $(V,V)$-bimodule. 
Its $\mathbf k$-submodule $F^1\mathbf C$ is a graded sub-$(V,V)$-bimodule. The morphism 
$\mathbf S\to\mathbf T$ induces a morphism $F^1\mathbf S\to F^1\mathbf T$ of $(V,V)$-bimodules, and 
$F^1\mathbf C=\mathrm{coker}(F^1\mathbf S\to F^1\mathbf T)$. 

(b) $\hat{\mathbf C}$ is a $(\hat V,\hat V)$-bimodule, equal to $\mathrm{coker}(\hat{\mathbf S}\to\hat{\mathbf T})$, and 
$F^1\hat{\mathbf C}$ is a sub-$(\hat V,\hat V)$-bimodule, equal to $\mathrm{coker}(F^1\hat{\mathbf S}\to F^1\hat{\mathbf T})$, and

(c) The $\mathbf k$-module $\hat{\mathbf C}$ is equipped with an action of the group $\mathbf k[[u,v]]^\times$ by 
$\phi\bullet q:=\phi(e_1,f_1)\cdot q\cdot \phi(e_1,f_1)^{-1}$, which preserves the $\mathbf k$-submodule $F^1\hat{\mathbf C}$.     
\end{lem}

\begin{proof}
The fact that $\mathbf S\to\mathbf T$ is graded follows from the fact that $e_0-f_0$
is homogeneous of degree 1. This implies the first part of (a). The rest (a) and (b) are direct consequences. 

The $\mathbf k$-module $\hat{\mathbf C}$ is equipped with an action of the group $\hat V^\times$ by 
$a\bullet q:=a\cdot g\cdot a^{-1}$, which preserves $F^1\hat{\mathbf C}$; the structure from (c) is its pull-back by the group 
morphism $\mathbf k[[u,v]]^\times\to\hat V^\times$, $\phi(u,v)\mapsto\phi(e_1,f_1)$. 
\end{proof}

Denote by $t\mapsto[t]_{\mathbf C}$ the projection $\hat{\mathbf T}\to\hat{\mathbf C}$. 

\begin{lem}\label{lem:invce:eqvce:cc}
    (a)The assignment
$$
\boldsymbol{\kappa} : (\begin{pmatrix}
    \alpha & \beta
\end{pmatrix},\begin{pmatrix}
    a \\ b
\end{pmatrix})\mapsto \Big[[\alpha\otimes a]_{\mathrm C_V(e_0)}+[\beta\otimes b]_{\mathrm C_V(f_0)}+
[\beta e_1\otimes a]_{\mathbf k[e_0,f_0]}\Big]_{\mathbf C}.  
$$
defines a map 
\begin{equation}\label{the:map:E3}
\boldsymbol{\kappa} : M_{12}F^1\hat V\times M_{21}\hat V\to F^1\hat{\mathbf C}. 
\end{equation}

(b) The map $\boldsymbol{\kappa}$ is $\mathbf k[[u,v]]^\times$-equivariant and $\mathrm C(\overline\rho_0)^\times$-invariant. 
\end{lem}

\begin{proof}
(a) For $\alpha,\beta\in F^1\hat V$ and $a,b\in \hat V$, one has $\alpha\otimes a,\beta\otimes b\in F^1(\hat V\hat\otimes \hat V)$, 
and $\beta e_1\otimes a\in F^2(\hat V\hat\otimes \hat V)=F^1(\hat V\hat\otimes \hat V)[1]$, therefore
$[\alpha\otimes a]_{\mathrm C_V(e_0)}+[\beta\otimes b]_{\mathrm C_V(f_0)}+
[\beta e_1\otimes a]_{\mathbf k[e_0,f_0]}\in F^1\hat{\mathbf T}$, which implies
$ \Big[[\alpha\otimes a]_{\mathrm C_V(e_0)}+[\beta\otimes b]_{\mathrm C_V(f_0)}+
[\beta e_1\otimes a]_{\mathbf k[e_0,f_0]}\Big]_{\mathbf C}\in F^1\hat{\mathbf C}$.  

(b) For $\phi\in\mathbf k[[u,v]]^\times$ and $(\begin{pmatrix}
    \alpha & \beta
\end{pmatrix},\begin{pmatrix}
    a \\ b
\end{pmatrix})\in M_{12}F^1\hat V\times M_{21}\hat V$, one has 
\begin{align*}
&\boldsymbol{\kappa}(\phi\bullet(\begin{pmatrix}
    \alpha & \beta
\end{pmatrix},\begin{pmatrix}
    a \\ b
\end{pmatrix}))=\boldsymbol{\kappa}(\begin{pmatrix}
    \phi(e_1,f_1)\alpha & \phi(e_1,f_1)\beta
\end{pmatrix},\begin{pmatrix}
    a\phi(e_1,f_1)^{-1} \\ b\phi(e_1,f_1)^{-1}
\end{pmatrix})
\\& =\Big[[\phi(e_1,f_1)\alpha\otimes a\phi(e_1,f_1)^{-1}]_{\mathrm C_V(e_0)}+[\phi(e_1,f_1)\beta\otimes b\phi(e_1,f_1)^{-1}]_{\mathrm C_V(f_0)}
\\& +
[\phi(e_1,f_1)\beta e_1\otimes a\phi(e_1,f_1)^{-1}]_{\mathbf k[e_0,f_0]}\Big]_{\mathbf C}
\\& =\phi\bullet
\Big[[\alpha\otimes a]_{\mathrm C_V(e_0)}+[\beta\otimes b]_{\mathrm C_V(f_0)}+
[\beta e_1\otimes a]_{\mathbf k[e_0,f_0]}\Big]_{\mathbf C}
=\phi\bullet \boldsymbol{\kappa}(\begin{pmatrix}
    \alpha & \beta
\end{pmatrix},\begin{pmatrix}
    a \\ b
\end{pmatrix}), 
\end{align*}
which implies the claimed $\mathbf k[[u,v]]^\times$-equivariance. 

Let $c\in \mathrm C(\overline\rho_0)^\times$ and 
$(\begin{pmatrix}
    \alpha & \beta
\end{pmatrix},\begin{pmatrix}
    a \\ b
\end{pmatrix})\in M_{12}F^1\hat V\times M_{21}\hat V$. By Lem. \ref{lem:8:3:departrennes}, there exists 
$\Pi\in\mathbf k[[f_0]]^\times,C\in \mathrm{C}_{\hat V}(e_0)$ such that $c=\overline X(\Pi,C)$. 
Then 
\begin{equation}\label{comp:actions}
\begin{pmatrix}
    \alpha & \beta
\end{pmatrix}\cdot c=\begin{pmatrix}
    \alpha(\Pi+(e_0-f_0)C)+\beta e_1C & \beta\Pi
\end{pmatrix}    \quad \mathrm{and}\quad 
c\cdot \begin{pmatrix}
    a \\ b
\end{pmatrix}=\begin{pmatrix}
    (\Pi+(e_0-f_0)C)a \\ e_1Ca+\Pi b
\end{pmatrix}. 
\end{equation}
Then 
\begin{align}\label{big:equality}
&\nonumber\boldsymbol{\kappa}((\begin{pmatrix}
    \alpha & \beta
\end{pmatrix}\cdot c,\begin{pmatrix}
    a \\ b
\end{pmatrix}))
\\&\nonumber =\Big[[(\alpha(\Pi+(e_0-f_0)C)+\beta e_1C)\otimes a]_{\mathrm C_V(e_0)}+[\beta\Pi\otimes b]_{\mathrm C_V(f_0)}+
[\beta\Pi e_1\otimes a]_{\mathbf k[e_0,f_0]}\Big]_{\mathbf C}
\\&\nonumber =\Big[[\alpha(\Pi+(e_0-f_0)C)\otimes a]_{\mathrm C_V(e_0)}
+[\beta e_1C\otimes a]_{\mathrm C_V(e_0)}+[\beta\Pi\otimes b]_{\mathrm C_V(f_0)}+
[\beta\Pi e_1\otimes a]_{\mathbf k[e_0,f_0]}\Big]_{\mathbf C}
\\&\nonumber =\Big[[\alpha\otimes (\Pi+(e_0-f_0)C)a]_{\mathrm C_V(e_0)}
+[\beta\otimes \Pi b]_{\mathrm C_V(f_0)}+
[\beta e_1\Pi\otimes a]_{\mathbf k[e_0,f_0]}+[\beta e_1\otimes Ca]_{\mathrm C_V(e_0)}\Big]_{\mathbf C}
\\&\nonumber=\Big[[\alpha\otimes (\Pi+(e_0-f_0)C)a]_{\mathrm C_V(e_0)}
+[\beta\otimes \Pi b]_{\mathrm C_V(f_0)}+
[\beta e_1\otimes \Pi a]_{\mathbf k[e_0,f_0]}+[\beta e_1\otimes Ca]_{\mathrm C_V(f_0)}
\\&\nonumber+[\beta e_1\otimes (e_0-f_0)Ca]_{\mathbf k[e_0,f_0]}\Big]_{\mathbf C}
\\&\nonumber=\Big[[\alpha\otimes (\Pi+(e_0-f_0)C)a]_{\mathrm C_V(e_0)}+[\beta\otimes e_1Ca]_{\mathrm C_V(f_0)}
+[\beta\otimes \Pi b]_{\mathrm C_V(f_0)}
+
[\beta e_1\otimes \Pi a]_{\mathbf k[e_0,f_0]}
\\&\nonumber+[\beta e_1\otimes (e_0-f_0)Ca]_{\mathbf k[e_0,f_0]}\Big]_{\mathbf C}
\\&\nonumber=\Big[[\alpha\otimes (\Pi+(e_0-f_0)C)a]_{\mathrm C_V(e_0)}+[\beta\otimes (e_1Ca+\Pi b)]_{\mathrm C_V(f_0)}+
[\beta e_1\otimes (\Pi+(e_0-f_0)C)a]_{\mathbf k[e_0,f_0]}\Big]_{\mathbf C}
\\&=\boldsymbol{\kappa}((\begin{pmatrix}
    \alpha & \beta
\end{pmatrix},c\cdot \begin{pmatrix}
    a \\ b
\end{pmatrix}))    
\end{align}
where the the first and seventh equalities follow from definitions and \eqref{comp:actions}, the second and sixth equalities are linear expansions, the 
third equality follows from $C\in \mathrm C_{\hat V}(e_0)$, $\Pi+(e_0-f_0)C\in \mathrm C_{\hat V}(e_0)$, 
$\Pi\in \mathrm C_{\hat V}(f_0)$ and the commutation of $e_1$ and $\Pi$, the fourth equality follows from 
$\Pi\in \mathrm C_{\hat V}(f_0)$ and the identity 
$[[u\otimes (e_0-f_0)v]_{\mathbf k[e_0,f_0]}+[u\otimes v]_{\mathrm C_V(f_0)}-[u\otimes v]_{\mathrm C_V(e_0)}]_{\mathbf C}=0$, the fifth 
equality follows from $e_1\in \mathrm C_{\hat V}(f_0)$. Then for any $\mathrm{row}\in M_{12}F^1\hat V$, $\mathrm{col}\in M_{21}\hat V$ 
and $c\in \mathrm C(\overline\rho_0)^\times$, one has  
$$
\boldsymbol{\kappa}((\mathrm{row},\mathrm{col})\bullet c)
=\boldsymbol{\kappa}(\mathrm{row}\cdot c,c^{-1}\cdot\mathrm{col})=\boldsymbol{\kappa}(\mathrm{row},c\cdot c^{-1}\cdot\mathrm{col})
=\boldsymbol{\kappa}(\mathrm{row},\mathrm{col}). 
$$
where the second equality follows from \eqref{big:equality}, which proves the claimed $\mathrm C(\overline\rho_0)^\times$-invariance.  
\end{proof}

\begin{defn}
The map 
$$
\mathbf k[[u,v]]^\times\backslash(M_{12}F^1\hat V\times M_{21}\hat V)/\mathrm C(\overline\rho_0)^\times \to 
\mathbf k[[u,v]]^\times\backslash F^1\hat{\mathbf C} 
$$
induced by the map $\boldsymbol{\kappa}$ and its invariance and equivariance properties (see Lem. \ref{lem:invce:eqvce:cc}) is denoted (E3).
\end{defn}

\subsection{The map (E4)}\label{sect:10:4}

For any $\alpha\in\mathbb Z$, the $\mathbf k$-module $\prod_{n\geq1}F^{n+\alpha}\hat V$ is equal to the set of maps 
$\delta : \mathbb Z_{\geq1}\to\hat V$ such that for any $n\geq 1$, $\delta(n)\in F^{n+\alpha}\hat V$. As each $F^{n+\alpha}\hat V$ is a 
$(\hat V,\hat V)$-bimodule, this $\mathbf k$-module 
is equipped with the product bimodule structure. It is equipped with a decreasing filtration of bimodules given by 
$F^d(\prod_{n\geq1}F^{n+\alpha}\hat V):=\prod_{n\geq1}F^{n+d+\alpha}\hat V$ for $d\geq0$. 

One has $F^\alpha\hat{\mathbf T}:=F^\alpha(V\otimes_{\mathrm{C}_V(e_0)}V)^\wedge\oplus F^\alpha(V\otimes_{\mathrm{C}_V(f_0)}V)^\wedge \oplus F^{\alpha+1}(V\otimes_{\mathbf k[e_0,f_0]}V)^\wedge$ for any $\alpha\in\mathbb Z$, therefore $F^{-1}\hat{\mathbf T}
=\hat{\mathbf T}$. 

\begin{lem}\label{lem:10:11:BIS}
(a) The map $V\otimes V\to \prod_{n\geq1}V$, $v\otimes w\mapsto \delta_{v\otimes w}^e$, where $\delta_{v\otimes w}^e$ is defined by 
$n\mapsto ve_0^{n-1}w$ for any $n\geq 1$, induces a  filtered $(V,V)$-bimodule morphism 
$V\otimes_{\mathrm C_V(e_0)} V\to \prod_{n\geq1}F^{n-1}V$, which leads to a $(\hat V,\hat V)$-bimodule morphism
$F^1(V\otimes_{\mathrm C_V(e_0)} V)^\wedge\to \prod_{n\geq1}F^n\hat V$. 

(b) The map $V\otimes V\to \prod_{n\geq1}V$, $v\otimes w\mapsto \delta_{v\otimes w}^f$, where $\delta_{v\otimes w}^f$ is defined 
by $n\mapsto vf_0^{n-1}w$ for any $n\geq 1$, induces a filtered $(V,V)$-bimodule morphism 
$V\otimes_{\mathrm C_V(f_0)} V\to \prod_{n\geq1}F^{n-1}V$, which leads to a $(\hat V,\hat V)$-bimodule morphism
$F^1(V\otimes_{\mathrm C_V(f_0)} V)^\wedge\to \prod_{n\geq1}F^n\hat V$. 

(c) The map $V\otimes V\to \prod_{n\geq1}V$, $v\otimes w\mapsto \delta_{v\otimes w}^{ef}$, where $\delta_{v\otimes w}^{ef}$ is defined by 
$n\mapsto v{e_0^{n-1}-f_0^{n-1}\over e_0-f_0}w$ for any $n\geq 1$, induces a filtered $(V,V)$-bimodule morphism 
$V\otimes_{\mathbf k[e_0,f_0]} V\to \prod_{n\geq1}F^{n-2}V$, which leads to a $(\hat V,\hat V)$-bimodule morphism
$F^2(V\otimes_{\mathrm C_V(e_f)} V)^\wedge\to \prod_{n\geq1}F^n\hat V$. 

(d) The $(V,V)$-bimodule morphism ${\mathbf T}\to \prod_{n\geq1}F^{n-2}V$
defined as the sum of the $(V,V)$-bimodule morphisms from (a), (b),(c) factors through a
$(V,V)$-bimodule morphism 
\begin{equation}\label{morph:C}
\mathbf C\to\prod_{n\geq1}F^{n-2}V     
\end{equation}
which is filtered, and therefore induces a $(\hat V,\hat V)$-bimodule morphism 
$$
F^1\hat{\mathbf C}\to\prod_{n\geq1}F^n\hat V. 
$$

(e) The map $F^1\hat{\mathbf C}\to \prod_{n\geq1}F^n\hat V$
from (d) is $\mathbf k[[u,v]]^\times$-equivariant, the action on the target space being defined by $(\phi\bullet\delta)(n):=
\phi(e_1,f_1)\delta(n)\phi(e_1,f_1)^{-1}$ for $\phi\in\mathbf k[[u,v]]^\times$, $\delta\in\prod_{n\geq1}F^n\hat V$, $n\geq1$,  
and therefore induces a map 
$$
\mathbf k[[u,v]]^\times\backslash F^1\hat{\mathbf C}\to \mathbf k[[u,v]]^\times\backslash\prod_{n\geq1}F^n\hat V. 
$$
which will be denoted (E4). 
\end{lem}

\begin{proof}
(a) For any $v,w\in V$ and $c\in\mathrm C_V(e_0)$, one has $vce_0^{n-1}w=ve_0^{n-1c}w$ for any $n\geq1$, hence 
$\delta^e_{vc\otimes w}=\delta^e_{v\otimes cw}$. It follows that $v\otimes w\mapsto \delta^e_{v\otimes w}$ induces a map 
$V\otimes_{\mathrm C_V(e_0)}V\to \prod_{n\geq1}V$. One has $v\cdot \delta_{v'\otimes w'}^e\cdot w=\delta_{vv'\otimes w'w}^e$, which implies that this map is a $(V,V)$-bimodule morphism.
For any $d\geq 0$, this map takes the degree $d$ part of the source to 
$\prod_{n\geq1}V_{d+n-1}$. The direct product over $d\geq0$ of the maps $(V\otimes_{\mathrm C_V(e_0)}V)_d\to\prod_{n\geq1}V_{d+n-1}$
is a map $(V\otimes_{\mathrm C_V(e_0)}V)^\wedge\to\prod_{d\geq0}\prod_{n\geq1}V_{d+n-1}=\prod_{n\geq1}F^{n-1}\hat V$, which 
takes $F^d(V\otimes_{\mathrm C_V(e_0)}V)^\wedge$ to $\prod_{d'\geq d}\prod_{n\geq1}V_{d'+n-1}
=\prod_{n\geq1}F^{n+d-1}\hat V$ and is a $(\hat V,\hat V)$-bimodule morphism as it is the completion of a 
$(V,V)$-bimodule morphism. 

The proof of (b), (c) is similar to that of (a). 

(d)  Let $(u,v)$ belong to $V\times V$, then the image of $[u\otimes v]_{\mathbf k[e_0,f_0]}$ by 
the map $\mathbf S\to \mathbf T$ is 
$-[u\otimes v]_{\mathrm C_V(e_0)}+[u\otimes v]_{\mathrm C_V(f_0)}+[u(e_0-f_0)\otimes v]_{\mathbf k[e_0,f_0]}$. 
The image of this element in $\prod_{n\geq1}F^{n-2}V$ is 
$\delta^e_{u\otimes v}-\delta^f_{u\otimes v}+\delta^{ef}_{u(e_0-f_0)\otimes v}$, which is given by 
$$
n\mapsto -ue_0^{n}v+uf_0^n v+u(e_0-f_0){e_0^n -f_0^n\over e_0-f_0}v=0
$$
for $n\geq1$, so that the said image is zero. This implies the vanishing of the composed map $\mathbf S\to\mathbf T\to 
\prod_{n\geq1}F^{n-2}V$, which implies the first statement. The morphism \eqref{morph:C} takes 
$\mathbf C_d=\mathrm{im}((V\otimes_{\mathrm C_V(e_0)} V)_d\oplus(V\otimes_{\mathrm C_V(f_0)} V)_d
\oplus (V\otimes_{\mathbf k[e_0,f_0]} V)_{d+1}\to\mathbf C)$ to $\prod_{n\geq 1}V_{n+d-1}$, hence 
induces a morphism from $\hat{\mathbf C}=\prod_{d\geq0}\mathbf C_d$ to $\prod_{d\geq0}\prod_{n\geq 1}V_{n+d-1}
=\prod_{n\geq1}F^{n-1}\hat V$, which restricts to a morphism from $F^1\hat{\mathbf C}=\prod_{d\geq1}\mathbf C_d$ to 
$\prod_{d\geq1}\prod_{n\geq 1}V_{n+d-1}=\prod_{n\geq1}F^{n}\hat V$. These morphisms are $(\hat V,\hat V)$-bimodule morphisms 
as they are completions of $(V,V)$-bimodule morphisms. 

(e) A $(\hat V,\hat V)$-bimodule $M$ is equipped with a $(\mathbf k[[u,v]],\mathbf k[[u,v]])$-bimodule structure, obtained
by pull-back by the algebra morphism $\mathbf k[[u,v]]\to\hat V$, $\phi\mapsto \phi(e_1,f_1)$, and therefore with a 
 $\mathbf k[[u,v]]^\times$-module structure obtained by $\phi\bullet m:=\phi\cdot m\cdot\phi^{-1}$ for any $m\in M$ and $\phi\in 
 \mathbf k[[u,v]]^\times$. Applying this to the $(\hat V,\hat V)$-bimodule structure of $\prod_{n\geq1}F^n\hat V$, 
 one obtains the said $\mathbf k[[u,v]]^\times$-module structure on $\prod_{n\geq1}F^n\hat V$. The statement then follows from the 
$(\hat V,\hat V)$-bimodule morphism status of the map from (d). 
\end{proof}

\subsection{The bijection (E5)}\label{sect:10:5}

\begin{lem}\label{lem:10:12:1603}
For $\delta\in \prod_{n\geq1}F^n\hat V$, there is an element $\Delta_\delta\in\prod_{n\geq1}F^n\hat V\to 
\mathrm{Hom}_{\mathcal C\operatorname{-alg}}(\hat{\mathcal W},\hat V)$ uniquely determined by the relations
$\Delta_\delta(e_0^{n-1}e_1):=\delta(n)$ for any $n\geq1$. The map $\delta\mapsto\Delta_\delta$ sets up a bijection 
$\prod_{n\geq1}F^n\hat V\to \mathrm{Hom}_{\mathcal C\operatorname{-alg}}(\hat{\mathcal W},\hat V)$
which is $\mathbf k[[u,v]]^\times$-equivariant, the actions of $\mathbf k[[u,v]]^\times$
on the source and on the target being as in Lem. \ref{lem:10:11:BIS}(e) and in Lem. \ref{lem28:1001}(c). 
The resulting bijection 
$$
\mathbf k[[u,v]]^\times\backslash\prod_{n\geq1}F^n\hat V\to 
\mathbf k[[u,v]]^\times\backslash\mathrm{Hom}_{\mathcal C\operatorname{-alg}}(\hat{\mathcal W},\hat V)
$$ 
is denoted (E5). 
\end{lem}

\begin{proof}
  It follows from the fact that $\hat{\mathcal W}$ is topologically generated by the elements $e_0^{n-1}e_1$, $n\geq0$ 
with $e_0^{n-1}e_1$ of degree $n$ that the map $\mathrm{Hom}_{\mathcal C\operatorname{-alg}}(\hat{\mathcal W},\hat V)\to\prod_{n\geq1}
\hat V$ given by $\Delta\mapsto (n\mapsto \Delta(e_0^{n-1}e_1))$ defines a bijection 
 $\mathrm{Hom}_{\mathcal C\operatorname{-alg}}(\hat{\mathcal W},\hat V)\to\prod_{n\geq1}
F^n\hat V$, denoted $\Delta\mapsto\delta_\Delta$; the map $\delta\mapsto\Delta_\delta$ is then its inverse. 
For $\phi\in\mathbf k[[u,v]]^\times$
and $\Delta\in\mathrm{Hom}_{\mathcal C\operatorname{-alg}}(\hat{\mathcal W},\hat V)$, one has 
$\delta_{\phi\bullet\Delta}(n)=(\phi\bullet\Delta)(e_0^{n-1}e_1)=\phi(e_1,f_1)\Delta(e_0^{n-1}e_1)\phi(e_1,f_1)^{-1}
=\phi(e_1,f_1)\delta_\Delta(n)\phi(e_1,f_1)^{-1}=(\phi\bullet\delta_\Delta)(n)$ for any $n\geq1$, therefore 
$\delta_{\phi\bullet\Delta}=\phi\bullet\delta_\Delta$, which 
implies the $\mathbf k[[u,v]]^\times$-equivariance of the bijection 
$\mathrm{Hom}_{\mathcal C\operatorname{-alg}}(\hat{\mathcal W},\hat V)\to\prod_{n\geq1}F^n\hat V$, 
from which the claimed $\mathbf k[[u,v]]^\times$-equivariance follows. 
\end{proof}

\subsection{A $\mathbf k$-module isomorphic to $\mathbf C$}\label{sect:10:6}

\begin{defn}
    Define $\mathcal W_l$ and $\mathcal W_r$ to be the $\mathbf k$-subalgebras of $\mathcal V$ respectively given by 
    $$
    \mathcal W_l:=\mathbf k\oplus\mathcal Ve_1,\quad \mathcal W_r:=\mathbf k\oplus e_1\mathcal V. 
    $$
\end{defn}
Then  $\mathcal W_l$ and $\mathcal W_r$ are graded subalgebras of $\mathcal V$, and $\hat{\mathcal W}$ (see Def. 
\ref{def:mrs:1608:BIS:BIS}(b)) and 
$\hat{\mathcal W}_r$ (see Lem. \ref{lem:toto:0930}(a)) are their graded completions. 

\begin{lem}\label{lem:complements}
    Let $S,T,\Sigma$ be graded $\mathbf k$-modules and let $\tau : S\to T$, $\sigma : S\to \Sigma$ be homogeneous $\mathbf k$-module 
    morphisms, therefore $\mathrm{im}(\tau\oplus\sigma)$ is a graded submodule of 
    $T[\mathrm{deg}\tau]\oplus \Sigma[\mathrm{deg}\sigma]$. If $\sigma$ is injective and if $\Sigma'$ is a graded complement of 
    $\mathrm{im}(\sigma)$ in $\Sigma$, then a graded complement of $\mathrm{im}(\tau\oplus\sigma : S\to T\oplus\Sigma)$ in 
    $T[\mathrm{deg}\tau]\oplus \Sigma[\mathrm{deg}\sigma]$ is $T[\mathrm{deg}\tau]\oplus\Sigma'[\mathrm{deg}\sigma]$. 
\end{lem}

\begin{proof}
Both $\mathrm{im}(\tau\oplus\sigma)$ and $T[\mathrm{deg}\tau]\oplus\Sigma'[\mathrm{deg}\sigma]$  
are graded $\mathbf k$-submodules of $T[\mathrm{deg}\tau]\oplus \Sigma[\mathrm{deg}\sigma]$, so that 
it suffices to prove that a complement of 
$\mathrm{im}(\tau\oplus\sigma : S\to T\oplus\Sigma)$ in $T\oplus \Sigma$ is 
    $T\oplus\Sigma'$ (in the category of $\mathbf k$-modules).

Let $x\in \mathrm{im}(\tau\oplus\sigma)\cap (T\oplus\Sigma')$; let $t\in T,s\in S$ be such that $x=t\oplus s$. 
Since $s\in T\oplus\Sigma'$, one has $s\in \Sigma'$. Let $s_0\in S$ be such that $x=(\tau\oplus\sigma)(s_0)$. 
Then $s=\sigma(s_0)$, therefore $s\in\mathrm{im}(\sigma)$. Since $\Sigma'\cap \mathrm{im}(\sigma)=0$, it follows that 
$s=0$. Since $\sigma$ is injective, 
it follows that $s_0=0$; it follows that $x=0$. All this proves that $\mathrm{im}(\tau\oplus\sigma)\cap (T\oplus\Sigma')=0$. 

Let now $x\in T\oplus \Sigma$, and let $t\in T,s\in S$ be such that $x=t\oplus s$. Since $\Sigma=\mathrm{im}(\sigma)+\Sigma'$,
there exists $(s_0,s')\in S\times \Sigma'$ such that $s=\sigma(s)+s'$. Let $t':=t-\tau(s_0)\in T$, 
then $x=(\tau\oplus\sigma)(s_0)+(t'\oplus s')$, where $(\tau\oplus\sigma)(s_0)\in \mathrm{im}(\tau\oplus\sigma)$ and 
$t'\oplus s'\in T\oplus\Sigma'$. It follows that $T\oplus\Sigma=\mathrm{im}(\tau\oplus\sigma)+(T\oplus\Sigma')$. 
\end{proof}

\begin{lem}\label{lem:1014:1803}
Define 
$$
\boldsymbol{\Sigma}:=\mathrm{im}((\mathcal W_r\otimes\mathcal V)\otimes (\mathcal W_l\otimes\mathcal V)
\hookrightarrow V\otimes V\to V\otimes_{\mathbf k[e_0,f_0]}V)\subset V\otimes_{\mathbf k[e_0,f_0]}V=\mathbf S.
$$
Then $\mathbf{T^e}\oplus\mathbf{T^f}\oplus\boldsymbol{\Sigma}[1]$ is a graded complement of $\mathrm{im}(\mathbf S\to\mathbf T)$ 
in $\mathbf T$. 
\end{lem}

\begin{proof}
Recall that $-\varphi_{\mathbf e}$, $\varphi_{\mathbf f}$, $\mu_{e_0-f_0}$ are morphisms from $\mathbf S$ to 
$\mathbf{T^e},\mathbf{T^f},\mathbf S$ homogeneous of degrees 0,0,1, 
and that the morphism $\mathbf S\to\mathbf{T^e}\oplus\mathbf{T^f}\oplus \mathbf S$ is 
$(-\varphi_{\mathbf e})\oplus \varphi_{\mathbf f}\oplus \mu_{e_0-f_0}$. 

Recall that the product induces right and left $\mathbf k[e_0]$-module isomorphisms
$$
    \mathcal W_r\otimes\mathbf k[e_0]\simeq \mathcal V,\quad \mathbf k[e_0]\otimes\mathcal W_l\simeq\mathcal V,  
$$
which induce a $\mathbf k$-module isomorphism
\begin{equation}\label{iso:V:V}
\mathcal V\otimes_{\mathbf k[e_0]}\mathcal V \simeq 
(\mathcal W_r\otimes\mathbf k[e_0])\otimes_{\mathbf k[e_0]}(\mathbf k[e_0]\otimes\mathcal W_l)
=\mathcal W_r\otimes \mathbf k[e_0]\otimes\mathcal W_l.  
\end{equation}

The sequence of  isomorphisms 
\begin{align*}
& \mathbf S=V\otimes_{\mathbf k[e_0,f_0]}V \simeq
(\mathcal V\otimes\mathcal V)\otimes_{\mathbf k[e_0]\otimes\mathbf k[e_0]}(\mathcal V\otimes\mathcal V)
\stackrel{\sim}{\to} (\mathcal V\otimes_{\mathbf k[e_0]}\mathcal V)\otimes (\mathcal V\otimes_{\mathbf k[e_0]}\mathcal V)
\\& \simeq(\mathcal W_r\otimes\mathbf k[e_0]\otimes\mathcal W_l)\otimes(\mathcal W_r\otimes\mathbf k[e_0]\otimes\mathcal W_l)
\simeq (\mathcal W_r\otimes\mathcal W_r)\otimes \mathbf k[e_0,f_0]\otimes (\mathcal W_l\otimes\mathcal W_l), 
\end{align*}
where the first and last isomorphisms follow from $\mathbf k[e_0]\otimes\mathbf k[e_0]\simeq\mathbf k[e_0,f_0]$, 
the second isomorphism follows from the isomorphism 
\begin{equation}\label{iso:swap}
    (M\otimes M')\otimes_{A\otimes A'}(N\otimes N')\simeq 
(M\otimes_A M')\otimes (N\otimes_B N'), 
\end{equation}
the third isomorphism follows from \eqref{iso:V:V}, 
sets up a graded $\mathbf k$-module isomorphism 
\begin{equation}\label{iso:S:explicit}
\mathbf S\simeq (\mathcal W_r\otimes\mathcal W_r)\otimes \mathbf k[e_0,f_0]\otimes (\mathcal W_l\otimes\mathcal W_l). 
\end{equation}
The morphism $\mu_{e_0-f_0}$ is intertwined by this isomorphism with the morphism 
$$
id\otimes ((e_0-f_0)\cdot-)\otimes id : 
(\mathcal W_r\otimes\mathcal W_r)\otimes \mathbf k[e_0,f_0]\otimes (\mathcal W_l\otimes\mathcal W_l)\to(\mathcal W_r\otimes
\mathcal W_r)\otimes \mathbf k[e_0,f_0]\otimes (\mathcal W_l\otimes\mathcal W_l)[1], 
$$
where $(e_0-f_0)\cdot-$ is the morphism $\mathbf k[e_0,f_0]\to\mathbf k[e_0,f_0][1]$ of multiplication by $e_0-f_0$.
Since $\mathbf k[e_0,f_0]$ is a domain, $(e_0-f_0)\cdot-$ is injective, therefore 
\begin{equation}\label{mu:e0f0:inj}
    \mu_{e_0-f_0}\text{ is injective}. 
\end{equation}
The fact that $\mu_{e_0-f_0}$ is intertwined with $id\otimes ((e_0-f_0)\cdot-)\otimes id$ also implies that 
\begin{equation}\label{explicit:im:mu:e0f0}
\text{the image of }\mathrm{im}(\mu_{e_0-f_0})\text{ by \eqref{iso:S:explicit} is }  
(\mathcal W_r\otimes\mathcal W_r)\otimes (e_0-f_0)\mathbf k[e_0,f_0]\otimes (\mathcal W_l\otimes\mathcal W_l).   
\end{equation}

The map $(\mathcal W_r\otimes\mathcal V)\otimes (\mathcal W_l\otimes\mathcal V)\to V\otimes_{\mathbf k[e_0,f_0]}V)$
in the definition of $\boldsymbol{\Sigma}$ admits a factorization 
$$
(\mathcal W_r\otimes\mathcal V)\otimes (\mathcal W_l\otimes\mathcal V)\twoheadrightarrow 
(\mathcal W_r\otimes\mathcal V)\otimes_{\mathbf k\otimes \mathbf k[e_0]} (\mathcal W_l\otimes\mathcal V)\to V\otimes_{\mathbf k[e_0,f_0]}V)
$$
where the first map is surjective, therefore 
\begin{equation}\label{new:eq:Sigma}
\boldsymbol{\Sigma}=\mathrm{im}((\mathcal W_r\otimes\mathcal V)\otimes_{\mathbf k\otimes \mathbf k[e_0]} 
(\mathcal W_l\otimes\mathcal V)\to V\otimes_{\mathbf k[e_0,f_0]}V). 
\end{equation}
The sequence of isomorphisms 
\begin{align*}
    & (\mathcal W_r\otimes\mathcal V)\otimes_{\mathbf k\otimes \mathbf k[e_0]} (\mathcal W_l\otimes\mathcal V)
\simeq (\mathcal W_r\otimes\mathcal W_l)\otimes(\mathcal V\otimes_{\mathbf k[e_0]}\mathcal V)
\simeq (\mathcal W_r\otimes\mathcal W_l)\otimes (\mathcal W_r\otimes\mathbf k[e_0]\otimes\mathcal W_l)
    \\& \simeq (\mathcal W_r\otimes\mathcal W_r)\otimes\mathbf k[f_0]\otimes(\mathcal W_l\otimes\mathcal W_l)
\end{align*}
where the first isomorphism follows from \eqref{iso:swap} and the second isomorphism follows from 
\eqref{iso:V:V}, gives rise to an isomorphism 
\begin{equation}\label{iso:source}
(\mathcal W_r\otimes\mathcal V)\otimes_{\mathbf k\otimes \mathbf k[e_0]} (\mathcal W_l\otimes\mathcal V)
\simeq (\mathcal W_r\otimes\mathcal W_r)\otimes\mathbf k[f_0]\otimes(\mathcal W_l\otimes\mathcal W_l)    
\end{equation}
which fits in a commutative diagram 
$$
\xymatrix{(\mathcal W_r\otimes\mathcal V)\otimes_{\mathbf k\otimes \mathbf k[e_0]} (\mathcal W_l\otimes\mathcal V)\ar[r]
\ar^\sim_{\eqref{iso:source}}[d]&
V\otimes_{\mathbf k[e_0,f_0]}V\ar^{\eqref{iso:S:explicit}}_{\sim}[d]\\
(\mathcal W_r\otimes\mathcal W_r)\otimes\mathbf k[f_0]\otimes(\mathcal W_l\otimes\mathcal W_l)\ar[r]&
(\mathcal W_r\otimes\mathcal W_r)\otimes \mathbf k[e_0,f_0]\otimes (\mathcal W_l\otimes\mathcal W_l)}
$$
where the top map is the map from \eqref{new:eq:Sigma} and the bottom map is the canonical inclusion. 
It follows that 
\begin{equation}\label{explicit:Sigma}
\text{the image of }\boldsymbol{\Sigma}\text{ by \eqref{iso:S:explicit} is }  
(\mathcal W_r\otimes\mathcal W_r)\otimes \mathbf k[f_0]\otimes (\mathcal W_l\otimes\mathcal W_l).   
\end{equation}
The sequence $(e_0-f_0)\mathbf k[e_0,f_0]\hookrightarrow\mathbf k[e_0,f_0]\twoheadrightarrow \mathbf k[f_0]$ is 
exact, where the second map is the algebra morphism such that $e_0\mapsto f_0$, $f_0\mapsto f_0$, and admits as a 
splitting the algebra morphism $\mathbf k[e_0]\to \mathbf k[e_0,f_0]$ given by $f_0\mapsto f_0$. It follows that 
 $\mathbf k[f_0]$ is a complement of $(e_0-f_0)\mathbf k[e_0,f_0]$ in $\mathbf k[e_0,f_0]$. 
Together with \eqref{explicit:im:mu:e0f0} and \eqref{explicit:Sigma}, this implies: 
\begin{equation}\label{statement:graded:complement}
   \boldsymbol{\Sigma}\text{ is a graded complement of }\mathrm{im}(\mu_{e_0-f_0})\text{ in }\mathbf S.  
\end{equation}
This and \eqref{mu:e0f0:inj} enables one to apply Lem. \ref{lem:complements}, where $S,T,\Sigma,\Sigma',\tau,\sigma$ are 
respectively taken to be equal to $\mathbf S,\mathbf{T^e}\oplus\mathbf{T^f},\mathbf S[1],\boldsymbol{\Sigma}[1]$, 
$(-\varphi_{\mathbf e})\oplus\varphi_{\mathbf f}$, $\mu_{e_0-f_0}$, from where the result follows.  
\end{proof}

\begin{defn}\label{def:10:13:1703}
(a)     Define the graded $\mathbf k$-modules
$$
\underline{\mathbf C}^{\mathbf e}:=(\mathcal W_r\otimes\mathbf k[e_0]\otimes\mathcal W_l)\otimes\mathcal V,\quad 
\underline{\mathbf C}^{\mathbf f}:=\mathcal V\otimes(\mathcal W_r\otimes\mathbf k[e_0]\otimes\mathcal W_l),\quad 
\underline{\mathbf C}^{\mathbf{ef}}:=(\mathcal W_r\otimes\mathcal W_r)\otimes\mathbf k[f_0]\otimes(\mathcal W_l\otimes\mathcal W_l)[1], 
$$
\begin{equation}\label{def:undeline:CC}
\underline{\mathbf C}:=\underline{\mathbf C}^{\mathbf e}\oplus \underline{\mathbf C}^{\mathbf f}
\oplus \underline{\mathbf C}^{\mathbf{ef}}  
\end{equation}
where in all the factors, $e_0,e_1,f_0,f_1$ are of degree 1. 

(b) $(\mathcal W_r\otimes\mathbf k[e_0]\otimes\mathcal W_l)\otimes\mathcal V\to V\otimes_{\mathrm C_V(e_0)}V$ is the graded 
$\mathbf k$-module morphism induced by the assignment $(w\otimes e_0^a\otimes w') \otimes v\mapsto 
[(we_0^a\otimes v)\otimes(w'\otimes 1)]_{\mathrm C_V(e_0)}$. 

(c)  $\mathcal V\otimes(\mathcal W_r\otimes\mathbf k[e_0]\otimes\mathcal W_l)\to V\otimes_{\mathrm C_V(f_0)}V$ is the graded 
$\mathbf k$-module morphism induced by the assignment $v\otimes (w\otimes e_0^a\otimes w')\mapsto [(v\otimes we_0^a)\otimes
(1\otimes w')]_{\mathrm C_V(f_0)}$. 

(d) $(\mathcal W_r\otimes\mathcal W_r)\otimes\mathbf k[f_0]\otimes(\mathcal W_l\otimes\mathcal W_l)\to 
(V\otimes_{\mathbf k[e_0,f_0]}V)$ is the graded $\mathbf k$-module morphism induced by the assignment
$(w\otimes w')\otimes f_0^a\otimes (w''\otimes w''') \mapsto [(w\otimes w')\otimes
(w''e_0^a\otimes w''')]_{\mathbf k[e_0,f_0]}$

(e) $\underline{\mathbf C}\to \mathbf T$ is the graded $\mathbf k$-module morphism  
given by the direct sum of the $\mathbf k$-linear maps from 
(b), (c), (d).  
\end{defn}

\begin{lem}\label{lem:iso:underlineC:C}
The map $\underline{\mathbf C}\to\mathbf C$, defined as the composition of the map of Def. \ref{def:10:13:1703}(e) with the projection 
$\mathbf T\to\mathbf C$, is an isomorphism of graded $\mathbf k$-modules. 
\end{lem}

\begin{proof}  
The linear maps from Def. \ref{def:10:13:1703}(b), (c), (d) are respectively maps $\underline{\mathbf C}^{\mathbf e}\to\mathbf T^{\mathbf e}$, 
$\underline{\mathbf C}^{\mathbf f}\to\mathbf T^{\mathbf f}$, $\underline{\mathbf C}^{\mathbf{ef}}\to\mathbf S$. 

The composition of the morphism $\underline{\mathbf C}^{\mathbf{ef}}\to\mathbf S$ with the isomorphism \eqref{iso:S:explicit}
is the canonical inclusion 
$$
(\mathcal W_r\otimes\mathcal W_r)\otimes \mathbf k[f_0]\otimes (\mathcal W_l\otimes\mathcal W_l)
\hookrightarrow
(\mathcal W_r\otimes\mathcal W_r)\otimes \mathbf k[e_0,f_0]\otimes (\mathcal W_l\otimes\mathcal W_l). 
$$ 
It then follows from \eqref{explicit:Sigma} that 
that 
\begin{equation}\label{iso:stat:0}
    \text{the image of }\underline{\mathbf C}^{\mathbf{ef}}\to\mathbf S\text{ is contained in $\boldsymbol{\Sigma}$,}
\end{equation}
and that
\begin{equation}\label{iso:stat:1}
    \text{the map }\underline{\mathbf C}^{\mathbf{ef}}\to\boldsymbol{\Sigma}\text{ is an isomorphism of $\mathbf k$-modules.}
\end{equation}
It follows from \eqref{iso:stat:0} that the composed map $\underline{\mathbf C}\to\mathbf T\to\mathbf C$ fits in a diagram 
\begin{equation}\label{double:triangle}
\xymatrix{&\mathbf T^{\mathbf{e}}\oplus\mathbf T^{\mathbf{f}}\oplus\boldsymbol{\Sigma}\ar[d]\ar^{(b)}[rd]&\\
\underline{\mathbf C}\ar[r]\ar^{(a)}[ru]&\mathbf T^{\mathbf e}\oplus\mathbf T^{\mathbf{f}}\oplus\mathbf S\ar[r]&\mathbf C}
\end{equation}

The following 
\begin{align*}
& \mathbf T^{\mathbf{e}}=V\otimes_{\mathrm C_V(e_0)}V \simeq
(\mathcal V\otimes\mathcal V)\otimes_{\mathbf k[e_0]\otimes \mathcal V}(\mathcal V\otimes\mathcal V)
\stackrel{\sim}{\to} (\mathcal V\otimes_{\mathbf k[e_0]}\mathcal V)
\otimes 
(\mathcal V\otimes_{\mathcal V}\mathcal V)
\\& \stackrel{\sim}{\to} (\mathcal W_r\otimes\mathbf k[e_0]\otimes\mathcal W_l)\otimes\mathcal V=\underline{\mathbf C}^{\mathbf{e}}
\end{align*}
is a sequence of $\mathbf k$-module isomorphisms, where the first isomorphism follows from Lem. \ref{lem:comm:e0}(a), 
the second isomorphism follows from \eqref{iso:swap}, the third isomorphism follows from \eqref{iso:V:V}. The map 
$\underline{\mathbf C}^{\mathbf{e}}\to\mathbf T^{\mathbf{e}}$ is such that $\underline{\mathbf C}^{\mathbf{e}}\to\mathbf T^{\mathbf{e}}
\to \underline{\mathbf C}^{\mathbf{e}}$ is the 
identity, therefore 
\begin{equation}\label{iso:stat:2}
    \text{the map }\underline{\mathbf C}^{\mathbf{e}}\to\mathbf T^{\mathbf{e}}\text{ is an isomorphism of $\mathbf k$-modules.}
\end{equation}

One similarly derives from the sequence of isomorphisms 
\begin{align*}
& \mathbf T^{\mathbf{f}}=V\otimes_{\mathrm C_V(f_0)}V \simeq
(\mathcal V\otimes\mathcal V)\otimes_{\mathcal V\otimes\mathbf k[e_0]}(\mathcal V\otimes\mathcal V)
\stackrel{\sim}{\to} (\mathcal V\otimes_{\mathcal V}\mathcal V)\otimes (\mathcal V\otimes_{\mathbf k[e_0]}\mathcal V)
\\& \stackrel{\sim}{\to} \mathcal V\otimes(\mathcal W_r\otimes\mathbf k[e_0]\otimes\mathcal W_l)=\underline{\mathbf C}^{\mathbf{f}}
\end{align*}
that 
\begin{equation}\label{iso:stat:3}
    \text{the map }\underline{\mathbf C}^{\mathbf{f}}\to\mathbf T^{\mathbf{f}}\text{ is an isomorphism of $\mathbf k$-modules.}
\end{equation}

Since the map (a) from \eqref{double:triangle} is the direct sum of the maps from \eqref{iso:stat:1}, \eqref{iso:stat:2} 
and \eqref{iso:stat:3}, the conjunction of these statements implies that this map is a $\mathbf k$-module isomorphism. 
The map (b) in \eqref{double:triangle} is a $\mathbf k$-module isomorphism by Lem. \ref{lem:1014:1803}. The commutativity 
of \eqref{double:triangle} then implies that the composed morphism $\underline{\mathbf C}\to \mathbf T\to\mathbf C$
is a $\mathbf k$-module isomorphism. The result then follows from the fact that this morphism is also graded. 
\end{proof}

\subsection{Local injectivity of (E4)}\label{sect:10:7}

There is a bijection between the set $\{e_0,e_1\}^*$ of words in $e_0,e_1$ and 
the set $\{(k,\underline a)|k\geq0,\underline a\in \mathbb Z_{\geq0}^{k+1}\}$, whose inverse takes 
the pair $(k,(a_0,\ldots,a_k))$ to the word $e_0^{a_0}e_1\cdots e_1 e_0^{a_k}$ (which is $e_0^{a_0}$ if 
$k=0$). For $w\in\{e_0,e_1\}^*$ a word, let $\mathrm{ht}(w):=\mathrm{max}(a_0,\ldots,a_k)$, where $(k,(a_0,\ldots,a_k))$
is the sequence corresponding to $w$.

\begin{defn}
Fix $N\geq 0$. 

(a) For any $s=0,1$, define $\{e_0,e_1\}^*_{s,N}$ to be the set of 
words such that the corresponding pair $(k,(a_0,\ldots,a_k))$ satisfies $|\{i|a_{i}\geq N\}|=s$. 

(b) Set $\{e_0,e_1\}^*_{\leq1,N}:=\{e_0,e_1\}^*_{0,N}\cup\{e_0,e_1\}^*_{1,N}$. 
\end{defn}

One checks that $\{e_0,e_1\}^*_{0,N}$ and $\{e_0,e_1\}^*_{1,N}$ are disjoint.

\begin{lem}\label{lem:en:cinq:points}
(a) $\{e_0,e_1\}^*_{0,N}$ is the set of words such that $\mathrm{ht}(w)<N$. 

(b) Let $w_r\in \{1\}\cup \{e_0,e_1\}^*e_1$, $w_l\in \{1\}\cup e_1\{e_0,e_1\}^*$ and $\alpha\in\mathbb Z_{\geq0}$. 
Then 
\begin{equation}\label{10.7.1:2103}
 w_r e_0^\alpha w_l\in \{e_0,e_1\}^*_{\leq1,\max(\mathrm{ht}(w_r),\mathrm{ht}(w_l))+1}.    
\end{equation}
More precisely, 
\begin{equation}\label{10.7.2:2103}
\text{$w_r e_0^\alpha w_l \in \{e_0,e_1\}^*_{0,\max(\mathrm{ht}(w_r),\mathrm{ht}(w_l))+1}$ iff 
$\alpha\leq \max(\mathrm{ht}(w_r),\mathrm{ht}(w_l))$} 
\end{equation}
and
\begin{equation}\label{10.7.3:2103}
\text{$w_r e_0^\alpha w_l\in \{e_0,e_1\}^*_{1,\max(\mathrm{ht}(w_r),\mathrm{ht}(w_l))+1}$ iff 
$\alpha>\max(\mathrm{ht}(w_r),\mathrm{ht}(w_l))$. }    
\end{equation}

(c) For $N\geq0$, the assignment 
\begin{equation}\label{assignement:kappa}
\{e_0,e_1\}^*_{1,N}\ni w=e_0^{a_0}e_1\cdots e_1e_0^{a_k}\mapsto (e_0^{a_0}e_1\cdots 
e_0^{a_{\alpha-1}}e_1,e_0^{a_\alpha-N},e_1e_0^{a_{\alpha+1}}\cdots e_1e_0^{a_{k}}),    
\end{equation}
where $\alpha$ is the unique element in $\{0,\ldots,k\}$ such that $a_\alpha\geq N$, 
with the convention that $e_0^{a_0}e_1\cdots 
e_0^{a_{\alpha-1}}e_1$ (resp. $e_1e_0^{a_{\alpha+1}}\cdots e_1e_0^{a_{k}}$) is 
$1$ if $\alpha=0$ (resp. $\alpha=k$), defines a map 
\begin{equation}\label{kappa:N:sets}
\kappa_N : \{e_0,e_1\}^*_{1,N}\to(\{1\}\cup \{e_0,e_1\}^*e_1)\times\{e_0\}^*\times(\{1\}\cup e_1\{e_0,e_1\}^*). 
\end{equation}

(d) One has 
\begin{equation}\label{increasing:words}
\forall N\leq M,\quad \{e_0,e_1\}^*_{\leq1,N}\subset\{e_0,e_1\}^*_{\leq1,M}. 
\end{equation}

(e) Let $w_r\in \{1\}\cup \{e_0,e_1\}^*e_1$, $w_l\in \{1\}\cup e_1\{e_0,e_1\}^*$ and 
$s\geq0$. Then for $n\geq\max(\mathrm{ht}(w_r),\mathrm{ht}(w_l))+1$, one has 
$$
w_re_0^{s+n}w_l\in\{e_0,e_1\}^*_{1,n}
$$
and
$$
\kappa_n(w_re_0^{s+n}w_l)=(w_r,s,w_l). 
$$
\end{lem}

\begin{proof}
(a) follows from the definitions.  

(b) There exist pairs $(p,(a_1,\ldots,a_p))$ and $(q,(b_1,\ldots,b_q))$ with $p,q\geq0$ such that 
the sequences associated with $w_r$ and $w_l$ are respectively $(a_1,\ldots,a_p,0)$ and $(0,b_1,\ldots,b_q)$; 
moreover, $\mathrm{ht}(w_r)=\max(a_1,\ldots,a_p)$ and $\mathrm{ht}(w_l)=\max(b_1,\ldots,b_q)$. 
The sequence associated to $w_re_0^\alpha w_l$ is then $(a_1,\ldots,a_p,\alpha,b_1,\ldots,b_q)$, which implies 
both \eqref{10.7.2:2103} and \eqref{10.7.3:2103}. \eqref{10.7.1:2103} then follows from the conjunction of 
these two results. 

(c) follows from the fact that $e_0^{a_0}e_1\cdots 
e_0^{a_{\alpha-1}}e_1$ (resp. $e_1e_0^{a_{\alpha+1}}\cdots e_1e_0^{a_{k}}$) necessarily belongs to 
the subset $\{1\}\cup\{e_0,e_1\}^*e_1$ (resp. $\{1\}\cup e_1\{e_0,e_1\}^*$) of $\{e_0,e_1\}^*$. 

(d) follows from the fact that $\{e_0,e_1\}^*_{\leq1,N}$ is the set of words such that the corresponding pair $(k,(a_0,\ldots,a_k))$ 
satisfies $|\{i|a_{i}\geq N\}|\leq1$. 

(e) Let $(p,(a_1,\ldots,a_p))$ and $(q,(b_1,\ldots,b_q))$ be as in (b). The sequence associated with $w_re_0^{s+n}w_l$
is then $(a_1,\ldots,a_p,s+n,b_1,\ldots,b_q)$. The first statement follows from 
$$
s+n\geq n>\max(\mathrm{ht}(w_r),\mathrm{ht}(w_l))=\max(a_1,\ldots,a_p,b_1,\ldots,b_q)
$$
and the second statement follows from the definition of $\kappa_n$. \end{proof}

Recall that $\mathcal V=\mathbf k\{e_0,e_1\}^*$, 
where $\mathbf kS$ is the free $\mathbf k$-module generated by a set $S$; for any $N\geq0$, set then  
$$
\mathcal V_{0,N}:=\mathbf k\{e_0,e_1\}^*_{0,N},\quad \mathcal V_{\leq1,N}:=\mathbf k\{e_0,e_1\}^*_{\leq1,N}=\mathbf k\{e_0,e_1\}^*_{0,N}\oplus \mathbf k\{e_0,e_1\}^*_{1,N}\subset \mathcal V. 
$$
Since $\mathcal W_r$ (resp. $\mathcal W_l$) is the free $\mathbf k$-module generated by 
$\{1\}\cup\{e_0,e_1\}^*e_1$ (resp. $\{1\}\cup e_1\{e_0,e_1\}^*$), the linearization of 
\eqref{kappa:N:sets} is a $\mathbf k$-module morphism 
\begin{equation}\label{kappa:partial}
\mathbf k\{e_0,e_1\}^*_{1,N}\to\mathcal W_r\otimes \mathbf k[e_0]\otimes \mathcal W_l.    
\end{equation}
 
\begin{defn} Let $N\geq 0$. 

(a) The $\mathbf k$-module morphism
$$
\kappa_N : \mathcal V_{\leq1,N}\to \mathcal W_r\otimes \mathbf k[e_0]\otimes \mathcal W_l   
$$
is the direct sum of the morphism \eqref{kappa:partial} with the zero morphism with source $\mathbf k\{e_0,e_1\}^*_{0,N}$.  

(b) The $\mathbf k$-module morphism
$$
\overline\kappa_N : \mathcal V_{\leq1,N}\to \mathcal W_r\otimes \mathcal W_l   
$$
is the composition of $\kappa_N$ with $id\otimes \epsilon\otimes id$, where $\epsilon : \mathbf k[e_0]\to\mathbf k$
is the $\mathbf k$-algebra morphism induced by $e_0\mapsto 0$. 
\end{defn}
 
\begin{lem}\label{lem:10:20:2103}
(a) For $N\leq M$, one has $\mathcal V_{\leq1,N}\subset\mathcal V_{\leq1,M}$, therefore 
the maps $\kappa_M : \mathcal V_{\leq1,N}\to \mathcal W_r\otimes \mathbf k[e_0]\otimes\mathcal W_l$, 
$\overline\kappa_M : \mathcal V_{\leq1,N}\to \mathcal W_r\otimes \mathcal W_l$ are well-defined for $M\geq N$. 
The restrictions of these maps to $\mathcal V_{0,N}\subset \mathcal V_{\leq1,N}\subset\mathcal V_{\leq1,M}$ are zero. 

(b) For any $c\in \underline{\mathbf C}^{\mathbf{e}}=(\mathcal W_l\otimes\mathbf k[e_0]\otimes\mathcal W_r)\otimes\mathcal V$, 
    there exists $N^e_c,M^e_c\geq 0$ such that for any $n\geq 1$, 
    $\delta^e_c(n)\in \mathcal V_{\leq1,N^e_c}\otimes \mathcal V_{0,M^e_c}$. For any $n\geq N^e_c$, 
    $(\kappa_{n}\otimes id)(\delta^e_c(n+1))$ is well-defined and  
    $$
    (\kappa_{n}\otimes id)(\delta^e_c(n+1))=c
    $$ 
    (equality in $(\mathcal W_l\otimes\mathbf k[e_0]\otimes\mathcal W_r)\otimes\mathcal V$). 

    (c) For any $c\in \underline{\mathbf C}^{\mathbf{f}}=\mathcal V\otimes(\mathcal W_l\otimes\mathbf k[e_0]\otimes\mathcal W_r)$, there exists $N^f_c,M^f_c\geq 0$ such that for any $n\geq 1$, 
    $\delta^f_c(n)\in \mathcal V_{0,N^f_c}\otimes \mathcal V_{\leq1,M^f_c}$. For any $n\geq M^f_c$, 
    $(id\otimes\kappa_{n})(\delta^f_c(n+1))$ is well-defined and  
    $$
    (id\otimes\kappa_{n})(\delta^f_c(n+1))=c
    $$ 
    (equality in $\mathcal V\otimes(\mathcal W_l\otimes\mathbf k[e_0]\otimes\mathcal W_r)$). 

    (d)  For any $c\in \underline{\mathbf C}^{\mathbf{ef}}\simeq(\mathcal W_r\otimes\mathbf k[e_0]\otimes\mathcal W_l)\otimes (\mathcal W_r\otimes\mathcal W_l)$, there exists $N^{ef}_c,M^{ef}_c\geq 0$ such that for any $n\geq 1$, 
    $\delta^{ef}_c(n)\in \mathcal V_{\leq1,N^{ef}_c}\otimes \mathcal V_{\leq1,M^{ef}_c}$, and for any $n\geq\max(N^{ef}_c,M^{ef}_c)$, 
    $(\kappa_{n}\otimes \overline\kappa_{n})(\delta^{ef}_c(2n+2))$ is well-defined and one has  
    $$
    (\kappa_{n}\otimes \overline\kappa_{n})(\delta^{ef}_c(2n+2))=c
    $$ 
    (equality in $(\mathcal W_r\otimes\mathbf k[e_0]\otimes\mathcal W_l)\otimes (\mathcal W_r\otimes\mathcal W_l)$). 
\end{lem}

\begin{proof}
The first part of (a) follows from \eqref{increasing:words} and its second part follows from 
$\{e_0,e_1\}^*_{0,N}\subset \{e_0,e_1\}^*_{0,M}$, which follows from Lem. \ref{lem:en:cinq:points}(a). 

(b)  Set $\mathcal B_e:=(\{1\}\cup\{e_0,e_1\}^*e_1)\times\mathbb Z_{\geq0}\times(\{1\}\cup e_1\{e_0,e_1\}^*)\times\{e_0,e_1\}^*$
and for $(w_r,s,w_l,v)\in\mathcal B_e$, set $c(w_r,s,w_l,v):=(w_r\otimes e_0^s\otimes w_l)\otimes v\in \underline{\mathbf C}^{\mathbf{e}}$. 
Then a basis of $\underline{\mathbf C}^{\mathbf{e}}$ is the set of elements $c(w_r,s,w_l,v)$, where 
$(w_r,s,w_l,v)$  belongs to $\mathcal B_e$. 

Fix $(w_r,s,w_l,v)\in \mathcal B_e$. For $n\geq1$, 
$\delta_{c(w_r,s,w_l,v)}^e(n)=w_r e_0^{n-1+s} w_l\otimes v$. Then $v\in \{e_0,e_1\}^*_{0,\mathrm{ht}(v)+1}\subset 
\mathcal V_{0,\mathrm{ht}(v)+1}$ by Lem. \ref{lem:en:cinq:points}(a) and $w_r e_0^{n-1+s} w_l\in  
\{e_0,e_1\}^*_{\leq1,\max(\mathrm{ht}(w_l),\mathrm{ht}(w_r))+1}\subset 
\mathcal V_{\leq1,\max(\mathrm{ht}(w_l),\mathrm{ht}(w_r))+1}$ by Lem. \ref{lem:en:cinq:points}(b), therefore 
\begin{equation}\label{inclusion:delta}
 \forall n\geq1,\quad \delta_{c(w_r,s,w_l,v)}^e(n)\in \mathcal V_{\leq1,\max(\mathrm{ht}(w_l),\mathrm{ht}(w_r))+1}\otimes \mathcal V_{0,\mathrm{ht}(v)+1} 
\end{equation}
therefore
$$
\forall n\geq0,\quad \delta_{c(w_r,s,w_l,v)}^e(n+1)\in \mathcal V_{\leq1,\max(\mathrm{ht}(w_l),\mathrm{ht}(w_r))+1}\otimes \mathcal V_{0,\mathrm{ht}(v)+1}.  
$$
Together with the first part of (a), this relation implies that for $n\geq\max(\mathrm{ht}(w_r),\mathrm{ht}(w_l))+1$, one has 
$$\delta_{c(w_r,s,w_l,v)}^e(n+1)\in\mathcal V_{\leq1,n}\otimes\mathcal V_{0,\mathrm{ht}(v)+1},
$$
so that 
\begin{equation}\label{interm:incl:1}
    \text{$\forall n\geq\max(\mathrm{ht}(w_r),\mathrm{ht}(w_l))+1$, 
    $(\kappa_n\otimes id)(\delta_{c(w_r,s,w_l,v)}^e(n+1))$ is well-defined. }
\end{equation}
For $n\geq0$, $\delta_{c(w_r,s,w_l,v)}^e(n+1)=w_r e_0^{n+s} w_l\otimes v$.
This implies the first equality in  
\begin{equation}\label{interm:incl:2}
\forall n\geq\max(\mathrm{ht}(w_r),\mathrm{ht}(w_l))+1, \quad 
(\kappa_n\otimes id)(\delta_{c(w_r,s,w_l,v)}^e(n+1))=w_r e_0^s w_l\otimes v=c(w_r,s,w_l,v), 
\end{equation}
where the second equality follows from definitions. \eqref{interm:incl:1} and \eqref{interm:incl:2} then imply 
\begin{align}\label{equality:delta:firstcase}
    &\text{for all $n\geq\max(\mathrm{ht}(w_r),\mathrm{ht}(w_l))+1$, } 
    \\&\nonumber 
    \text{$(\kappa_n\otimes id)(\delta_{c(w_r,s,w_l,v)}^e(n+1))$
    is well-defined and equal to $c(w_r,s,w_l,v)$.} 
\end{align}
Let now $c\in \underline{\mathbf C}^{\mathbf{e}}$ and $\mathrm{supp}(c)\subset \mathcal B_e$ be the set of tuples 
$(w_r,s,w_l,v)$ such that the coefficient of $c(w_r,s,w_l,v)$ in the decomposition of $c$ is nonzero. Set 
$$
N_c^e:=1+\max(\{\mathrm{ht}(w_r)|(w_r,s,w_l,v)\in\mathrm{supp}(c)\}\cup\{\mathrm{ht}(w_l)|(w_r,s,w_l,v)\in\mathrm{supp}(c)\}), 
$$
$$
M_c^e:=1+\max\{\mathrm{ht}(v)|(w_r,s,w_l,v)\in\mathrm{supp}(c)\},  
$$
then \eqref{inclusion:delta} implies $\delta^e_c(n)\in \mathcal V_{\leq1,N^e_c}\otimes \mathcal V_{0,M^e_c}$, and  
\eqref{equality:delta:firstcase} implies that for any $n\geq N_e^c$, 
$(\kappa_n\otimes id)(\delta^e_c(n+1))$ is well-defined and is equal to $c$. 
This proves (b). 

(c) can be derived from (b) by applying 
the automorphism of exchange of factors in $V=\mathcal V\otimes\mathcal V$.

(d) Let $\mathcal B_{ef}:=(\{1\}\cup\{e_0,e_1\}^*e_1)^2\times\mathbb Z_{\geq0}\times (\{1\}\cup e_1\{e_0,e_1\}^*)^2$
and for $(w_r,w'_r,s,w_l,w'_l)\in \mathcal B_{ef}$, set 
$c(w_r,w'_r,s,w_l,w'_l):=(w_r\otimes w'_r)\otimes e_0^s(w_l\otimes w'_l)
\in \underline{\mathbf C}^{\mathbf{ef}}$. 
Then a basis of $\underline{\mathbf C}^{\mathbf{ef}}$ is the set of elements $c(w_r,w'_r,s,w_l,w'_l)$, where 
$(w_r,w'_r,s,w_l,w'_l)$  belongs to $\mathcal B_{ef}$. 

For $(w_r,w'_r,s,w_l,w'_l)\in\mathcal B_{ef}$. For any $n\geq1$, 
$$
\delta^{ef}_{c(w_r,w'_r,s,w_l,w'_l)}(n)=\sum_{\alpha=0}^{n-2} w_re_0^{s+\alpha}w_l\otimes w'_re_0^{n-2-\alpha}w'_l
$$
(this sum is 0 when $n=1$). By Lem. \ref{lem:en:cinq:points}(b), one has for any $\alpha\in\{0,\ldots,n-2\}$
the relations $w_re_0^{s+\alpha}w_l\in\{e_0,e_1\}^*_{\leq1,\max(\mathrm{ht}(w_r),\mathrm{ht}(w_l))+1}$ and 
$w'_re_0^{n-2-\alpha}w'_l\in\{e_0,e_1\}^*_{\leq1,\max(\mathrm{ht}(w'_r),\mathrm{ht}(w'_l))+1}$, which by the 
relation $\{e_0,e_1\}^*_{\leq1,\alpha}\subset \mathcal V_{\leq1,\alpha}$ imply
\begin{equation}\label{inclusion:ef}
\forall n\geq1,\quad \delta^{ef}_{c(w_r,w'_r,s,w_l,w'_l)}(n)\in \mathcal V_{\leq1,\max(\mathrm{ht}(w_r),\mathrm{ht}(w_l))+1}
\otimes \mathcal V_{\leq1,\max(\mathrm{ht}(w'_r),\mathrm{ht}(w'_l))+1}, 
\end{equation}
therefore 
$$
\forall n\geq 0,  \delta^{ef}_{c(w_r,w'_r,s,w_l,w'_l)}(2n+2)\in \mathcal V_{\leq1,\max(\mathrm{ht}(w_r),\mathrm{ht}(w_l))+1}
\otimes \mathcal V_{\leq1,\max(\mathrm{ht}(w'_r),\mathrm{ht}(w'_l))+1}. 
$$
Together with the first part of (a), this relation implies that
$$
\forall n\geq\max(\mathrm{ht}(w_r),\mathrm{ht}(w_l),\mathrm{ht}(w'_r),\mathrm{ht}(w'_l))+1,\quad 
 \delta^{ef}_{c(w_r,w'_r,s,w_l,w'_l)}(2n+2)\in \mathcal V_{\leq1,n}
\otimes \mathcal V_{\leq1,n}
$$
therefore 
\begin{equation}\label{interm:11}
 \text{$\forall n\geq\max(\mathrm{ht}(w_r),\mathrm{ht}(w_l),\mathrm{ht}(w'_r),\mathrm{ht}(w'_l))+1$, 
$(\kappa_n\otimes\overline\kappa_n)(\delta^{ef}_{c(w_r,w'_r,s,w_l,w'_l)}(2n+2))$ is well-defined.}   
\end{equation}
For $n\geq0$, 
$$
\delta^{ef}_{c(w_r,w'_r,s,w_l,w'_l)}(2n+2)=\sum_{\alpha=0}^{2n} w_re_0^{s+\alpha}w_l\otimes w'_re_0^{2n-\alpha}w'_l. 
$$
For $n\geq\max(\mathrm{ht}(w'_r),\mathrm{ht}(w'_l))+1$ and 
$\alpha\in\{0,\ldots,2n\}$, one has $\overline\kappa_n(w'_re_0^{2n-\alpha}w'_l)=\delta_{\alpha,n}w'_r\otimes w'_l$. 
On the other hand, for $n\geq\max(\mathrm{ht}(w_r),\mathrm{ht}(w_l))+1$, one has $\kappa_n(w_re_0^{s+n}w_l)=w_r\otimes e_0^s\otimes w_l$. 
All this implies the first equality in  
\begin{align}\label{interm:12}
    & \forall n\geq\max(\mathrm{ht}(w_r),\mathrm{ht}(w_l),\mathrm{ht}(w'_r),\mathrm{ht}(w'_l))+1,
    \\& \nonumber (\kappa_n\otimes\overline\kappa_n)(\delta^{ef}_{c(w_r,w'_r,s,w_l,w'_l)}(2n+2))=(w_r\otimes e_0^s\otimes w_l)\otimes
(w'_r\otimes w'_l)=c(w_r,w'_r,s,w_l,w'_l). 
\end{align}
\eqref{interm:11} and \eqref{interm:12} imply 
\begin{align}\label{interm:13}
  &\text{$\forall n\geq\max(\mathrm{ht}(w_r),\mathrm{ht}(w_l),\mathrm{ht}(w'_r),\mathrm{ht}(w'_l))+1$,}
  \\& \nonumber \text{ 
$(\kappa_n\otimes\overline\kappa_n)(\delta^{ef}_{c(w_r,w'_r,s,w_l,w'_l)}(2n+2))$ is well-defined
 and equal to $c(w_r,w'_r,s,w_l,w'_l)$.}   
\end{align}

Let now $c\in \underline{\mathbf C}^{\mathbf{ef}}$ and $\mathrm{supp}(c)\subset \mathcal B_{ef}$ be the set of tuples 
$(w_r,w'_r,s,w_l,w'_l)$ such that the coefficient of $c(w_r,w'_r,s,w_l,w'_l)$ in the decomposition of $c$ is nonzero. Set 
$$
N_c^{ef}:=1+\max(\{\mathrm{ht}(w_r)|(w_r,w'_r,s,w_l,w'_l)\in\mathrm{supp}(c)\}\cup
\{\mathrm{ht}(w_l)|(w_r,w'_r,s,w_l,w'_l)\in\mathrm{supp}(c)\}), 
$$
$$
M_c^{ef}:=1+\max(\{\mathrm{ht}(w'_r)|(w_r,w'_r,s,w_l,w'_l)\in\mathrm{supp}(c)\}\cup
\{\mathrm{ht}(w'_l)|(w_r,w'_r,s,w_l,w'_l)\in\mathrm{supp}(c)\}),  
$$
then \eqref{inclusion:ef} implies that for any $n\geq0$, 
$\delta^{ef}_c(n)\in \mathcal V_{\leq1,N^{ef}_c}\otimes \mathcal V_{0,M^{ef}_c}$, and 
\eqref{interm:13} implies that for any 
$$
n\geq 1+\max\{\max(\mathrm{ht}(w_r),\mathrm{ht}(w_l),\mathrm{ht}(w'_r),\mathrm{ht}(w'_l))|
(w_r,w'_r,s,w_l,w'_l)\in\mathrm{supp}(c)\}=\max(N_c^{ef},M_c^{ef}), 
$$
$(\kappa_n\otimes\overline\kappa_n)(\delta^{ef}_c(2n+2))$ is well-defined and equal to $c$. 
\end{proof}

\begin{lem}\label{lem:10:21}
    (a) Let $c=c_e\oplus c_f\oplus c_{ef}\in \underline{\mathbf C}^{\mathbf{e}}\oplus \underline{\mathbf C}^{\mathbf{f}}\oplus
    \underline{\mathbf C}^{\mathbf{ef}}=\underline{\mathbf C}$. There exist $N_c,M_c\geq 0$ such that for 
    $n\geq 1$, $\delta_c(n)\in \mathcal V_{\leq1,N_c}\otimes \mathcal V_{\leq1,M_c}$. 
    For $n\geq \max(N_c,M_c)$, $(\kappa_n\otimes\overline\kappa_n)(\delta_c(2n+2))$
    is well-defined and equal to $c_{ef}$. 

    (b) Let $c=c_e\oplus c_f\in \underline{\mathbf C}^{\mathbf{e}}\oplus \underline{\mathbf C}^{\mathbf{f}}\subset\underline{\mathbf C}$. 
     There exist $\tilde N_c,\tilde M_c\geq 0$ such that for 
    $n\geq 1$, $\delta_c(n)\in \mathcal V_{\leq1,\tilde N_c}\otimes \mathcal V_{\leq1,\tilde M_c}$. 
    For $n\geq \tilde M_c$, $(id\otimes\kappa_n)(\delta_c(n+1))$
    is well-defined and equal to $c_{f}$. 
\end{lem}

\begin{proof}
(a) Let $m\geq1$. By Lem. \ref{lem:10:20:2103}(b), $\delta_{c_e}^e(m)\in \mathcal V_{\leq 1,N^e_{c_e}}\otimes 
\mathcal V_{0,M^e_{c_e}}$, $\delta_{c_f}^f(m)\in \mathcal V_{0,N^f_{c_f}}\otimes 
\mathcal V_{\leq 1,M^f_{c_f}}$ and $\delta_{c_{ef}}^{ef}(m)\in\mathcal V_{\leq 1,N^{ef}_{c_{ef}}}\otimes 
\mathcal V_{\leq 1,M^{ef}_{c_{ef}}}$. Set $N_c:=\max(N^e_{c_e},N^f_{c_f},N^{ef}_{c_{ef}})$ and 
$M_c:=\max(M^e_{c_e},M^f_{c_f},M^{ef}_{c_{ef}})$, then each of the elements $\delta_{c_e}^e(m)$, 
$\delta_{c_f}^f(m)$ and $\delta_{c_{ef}}^{ef}(m)$ belongs to $\mathcal V_{\leq1,N_c}\otimes 
\mathcal V_{\leq1,M_c}$, therefore so does their sum $\delta_c(n)$. 

Then follows from the first part of Lem. \ref{lem:10:20:2103}(a) 
that for any $m\geq1$ and any $n\geq \max(N_c,M_c)$, all the terms in the following equality are well-defined and
the equality holds
$$
(\kappa_n\otimes\overline\kappa_n)(\delta_c(m))
=(\kappa_n\otimes\overline\kappa_n)(\delta^e_{c_e}(m))
+(\kappa_n\otimes\overline\kappa_n)(\delta^f_{c_f}(m))
+(\kappa_n\otimes\overline\kappa_n)(\delta^{ef}_{c_{ef}}(m)). 
$$
By the relation $n\geq N^e_{c_e}$ and the first part of Lem. \ref{lem:10:20:2103}(a), 
$(\kappa_n\otimes id)(\delta^e_{c_e}(m))=0$ hence 
$(\kappa_n\otimes\overline\kappa_n)(\delta^e_{c_e}(m))=0$. 
Similarly, the relation $n\geq N^f_{c_f}$ implies $(id\otimes\overline\kappa_n)(\delta^f_{c_f}(m))=0$.
Then $m\geq 1$ and $n\geq\max(N_c,M_c)$ implies 
$$
(\kappa_n\otimes\overline\kappa_n)(\delta_c(m))
=(\kappa_n\otimes\overline\kappa_n)(\delta^{ef}_{c_{ef}}(m)). 
$$
Lem. \ref{lem:10:20:2103}(d) then implies that if $n\geq \max(N_c,M_c,N^{ef}_{c_{ef}},M^{ef}_{c_{ef}})=\max(N_c,M_c)$, then 
$$
(\kappa_n\otimes\overline\kappa_n)(\delta_c(2n+2))
=c_{ef}. 
$$

(b) Let $m\geq1$. By Lem. \ref{lem:10:20:2103}(b), $\delta_{c_e}^e(m)\in \mathcal V_{\leq 1,N^e_{c_e}}\otimes 
\mathcal V_{0,M^e_{c_e}}$ and $\delta_{c_f}^f(m)\in \mathcal V_{0,N^f_{c_f}}\otimes 
\mathcal V_{\leq 1,M^f_{c_f}}$. Set $\tilde N_c:=\max(N^e_{c_e},N^f_{c_f})$ and 
$\tilde M_c:=\max(M^e_{c_e},M^f_{c_f})$, then both $\delta_{c_e}^e(m)$ and 
$\delta_{c_f}^f(m)$ belong to $\mathcal V_{\leq1,\tilde N_c}\otimes 
\mathcal V_{\leq1,\tilde M_c}$, therefore so does their sum $\delta_c(m)$. This proves the first statement. 

It then follows from the first part of Lem. \ref{lem:10:20:2103}(a) that for any $n\geq \tilde M_c$,
$(id\otimes\kappa_n)(\delta_c(m))$, $(id\otimes\kappa_n)(\delta_{c_e}^e(m))$ and 
$(id\otimes\kappa_n)(\delta_{c_f}^f(m))$ are well-defined, and 
$$
(id\otimes\kappa_n)(\delta_c(m))=(id\otimes\kappa_n)(\delta_{c_e}^e(m))
+(id\otimes\kappa_n)(\delta_{c_f}^f(m)). 
$$
Moreover, $\delta_{c_e}^e(m)\in \mathcal V_{\leq 1,N^e_{c_e}}\otimes 
\mathcal V_{0,M^e_{c_e}}\subset \mathcal V_{\leq1,\tilde N_c}\otimes 
\mathcal V_{0,\tilde M_c}$ therefore by the second part of Lem. \ref{lem:10:20:2103}(a), 
$(id\otimes\kappa_n)(\delta_{c_e}^e(m))=0$
therefore 
$$
(id\otimes\kappa_n)(\delta_c(m))=(id\otimes\kappa_n)(\delta_{c_f}^f(m)). 
$$
It follows that 
$$
\forall n\geq  \tilde M_c\quad 
(id\otimes\kappa_n)(\delta_c(n+1))=(id\otimes\kappa_n)(\delta_{c_f}^f(n+1)). 
$$
Since $\tilde M_c\geq M^f_{c_f}$, the right-hand side is equal to $c_f$, which implies the second statement. 
\end{proof}

\begin{lem}\label{lem:E4:inj}
(a) The maps 
$$
\mathbf C\to\prod_{n\geq1}F^{n-2}V\text{ and }F^1\hat{\mathbf C}\to\prod_{n\geq1}F^n\hat V
$$
from Lem. \ref{lem:10:11:BIS}(d) are injective. 

(b) The map (E4) is injective.     
\end{lem}

\begin{proof}
 (a) Let us prove the injectivity of the composed map $\underline{\mathbf C}\to \mathbf C\to\prod_{n\geq1}F^{n-2}V$, 
 where the map $\underline{\mathbf C}\to \mathbf C$ is as in Lem. \ref{lem:iso:underlineC:C}. 
Let $c\in \underline{\mathbf C}$ belong to the kernel of this map, so for any $n\geq 1$, $\delta_c(n)=0$. 
Let $c_e\in\underline{\mathbf C}^{\mathbf{e}}$, 
$c_f\in\underline{\mathbf C}^{\mathbf{f}}$, $c_{ef}\in\underline{\mathbf C}^{\mathbf{ef}}$ be such that $c=c_e\oplus c_f\oplus c_{ef}$. 
Let $N_c,M_c$ be as in Lem. \ref{lem:10:21}(a), then this result says that 
for $n\geq\max(N_c,M_c)$, $c_{ef}=(\kappa_n\otimes\overline\kappa_n)(\delta_c(2n+2))$. Therefore $c_{ef}=0$. 
Therefore $c=c_e\oplus c_f\in \underline{\mathbf C}^{\mathbf{e}}\oplus \underline{\mathbf C}^{\mathbf{f}}\subset \underline{\mathbf C}$. 
Let $\tilde N_c,\tilde M_c$ be as in Lem. \ref{lem:10:21}(b), then this result says that 
for $n\geq\max(\tilde M_c)$, $c_f=(id\otimes\kappa_n)(\delta_c(n+1))$. Therefore $c_f=0$, so $c=c_e\in \underline{\mathbf C}^{\mathbf{e}}
\subset \underline{\mathbf C}$, therefore for any $n\geq 1$, $\delta^e_{c_e}(n)=0$. Then Lem. \ref{lem:10:20:2103}(b) 
says that for any $n\geq N_{c_e}^e$, $(\kappa_n\otimes id)(\delta^e_{c_e}(n+1))=c_e$. It follows that $c_e=0$, therefore 
$c=0$. 
It follows that $\underline{\mathbf C}\to \mathbf C\to\prod_{n\geq1}F^{n-2}V$ is injective.
Lem. \ref{lem:iso:underlineC:C} then implies the injectivity of $\mathbf C\to\prod_{n\geq1}F^{n-2}V$. 

By the proof of Lem. \ref{lem:10:11:BIS}(d), this map admits a factorization 
$$
\mathbf C=\oplus_{d\geq-1}\mathbf C_d\to \oplus_{d\geq-1}(\textstyle{\prod}_{n\geq 1}V_{d+n-1})\hookrightarrow\textstyle{\prod}_{n\geq1}
(\oplus_{d\geq-1}V_{d+n-1})=\textstyle{\prod}_{n\geq 1}F^{n-2}V,
$$
where the middle map is the direct some over $d\geq-1$
of linear maps $\mathbf C_d\to \prod_{n\geq 1}V_{d+n-1}$. The injectivity of $\mathbf C\to\prod_{n\geq1}F^{n-2}V$
then implies the injectivity of the map $\mathbf C_d\to \prod_{n\geq 1}V_{d+n-1}$ for each $d\geq-1$. One derives
the injectivity of the direct product of these maps, which is 
$\hat{\mathbf C}=\prod_{d\geq-1}\mathbf C_d\to \prod_{d\geq-1} \prod_{n\geq 1}V_{d+n-1}=\prod_{n\geq 1}F^{n-2}\hat V$, 
as well as of $F^1\hat{\mathbf C}=\prod_{d\geq1}\mathbf C_d\to \prod_{d\geq1} \prod_{n\geq 1}V_{d+n-1}
=\prod_{n\geq 1}F^n\hat V$. 
 
(b) Follows from (a) and for the fact that if $G$ is a group, $X,Y$ are sets with actions of $G$
and $X\to Y$ is a $G$-equivariant map, then the injectivity of $X\to Y$ implies the injectivity of the 
induced map $G\backslash X\to G\backslash Y$. 
 \end{proof}

\subsection{Local injectivity of (E3)}\label{sect:10:8}

\begin{defn}
    $\mathbf c_{\mathrm{DT}}$ is the element of $F^1\hat{\mathbf C}$ defined by 
    $$
    \mathbf c_{\mathrm{DT}}:= \Big[[e_1\otimes 1]_{\mathrm C_V(e_0)}+[f_1\otimes 1]_{\mathrm C_V(f_0)}
-[e_1f_1\otimes 1]_{\mathbf k[e_0,f_0]}\Big]_{\mathbf C}\in F^1\hat{\mathbf C}
    $$
\end{defn}

Recall 
 $(\overline{\mathrm{row}}_{\mathrm{DT}},\overline{\mathrm{col}}_{\mathrm{DT}})=(\begin{pmatrix}
    e_1 & -f_1
\end{pmatrix},\begin{pmatrix}
    1\\-1
\end{pmatrix})$. Then  $\boldsymbol{\kappa}$ (see \eqref{the:map:E3}) is such that 
$$
\boldsymbol{\kappa}(\overline{\mathrm{row}}_{\mathrm{DT}},\overline{\mathrm{col}}_{\mathrm{DT}})
=\mathbf c_{\mathrm{DT}}. 
$$ 

It follows from Lem. \ref{lem:invce:eqvce:cc}(b) that $\boldsymbol{\kappa}$ 
admits a factorization through a $\mathbf k[[u,v]]^\times$-equivariant map 
\begin{equation}\label{map:MM/C:Q}
\overline{\boldsymbol{\kappa}} : (M_{12}F^1\hat V\times M_{21}\hat V)/\mathrm C(\overline\rho_0)^\times\to F^1\hat{\mathbf C}
\end{equation}
such that 
\begin{equation}\label{underlinec:effect:on:marked:points}
(\overline{\mathrm{row}}_{\mathrm{DT}},\overline{\mathrm{col}}_{\mathrm{DT}})\bullet 
\mathrm C_2(\overline\rho_0)^\times \mapsto \mathbf c_{\mathrm{DT}}\in F^1\hat{\mathbf C}.
\end{equation}

\begin{defn}
    (a) Set $\mathbf X:=M_{12}F^1\hat V\times M_{21}\hat V$. 

    (b) For $n\geq 0$, set $F^n\mathbf X:=(\overline{\mathrm{row}}_{\mathrm{DT}}
    +M_{12}F^{n+1}\hat V)\times(\overline{\mathrm{col}}_{\mathrm{DT}}+M_{21}F^n\hat V)$. 

    (c) Set $F^n(\mathrm C_2(\overline\rho_0)^\times):=\mathrm C_2(\overline\rho_0)^\times$ and for 
    $n\geq1$, set $F^n(\mathrm C_2(\overline\rho_0)^\times):=I_2+(\mathrm C_2(\overline\rho_0)\cap M_2F^n\hat V)$. 
\end{defn}

\begin{lem}
(a) There is a decreasing filtration of sets $\mathbf X=F^0\mathbf X\supset F^1\mathbf X\supset\cdots$. 

(b) There is a decreasing filtration of groups $\mathrm C_2(\overline\rho_0)^\times = F^0\mathrm C_2(\overline\rho_0)^\times \supset 
F^1\mathrm C_2(\overline\rho_0)^\times\supset \cdots$. 

(c) For each $n\geq0$, the right action of $\mathrm C_2(\overline\rho_0)^\times$ on $\mathbf X$ restricts to a right action of 
$F^n(\mathrm C_2(\overline\rho_0)^\times)$ on $F^n\mathbf X$. 

(d) For each $n\geq0$, the map $\boldsymbol{\kappa}$ induces a map $\boldsymbol{\kappa} : F^n\mathbf X\to 
\mathbf c_{\mathrm{DT}}+F^{n+1}\hat{\mathbf C}$. 
\end{lem}

\begin{proof}
(a) is obvious. (b) follows from the fact that for any $n\geq1$, $\mathrm C_2(\overline\rho_0)\cap M_2F^n\hat V$ is a 
subalgebra without unit of $\mathrm C_2(\overline\rho_0)$, and from the convergence of the series $\sum_{i\geq0}x^i$
for $x\in \mathrm C_2(\overline\rho_0)\cap M_2F^n\hat V$. (c) follows from the relations 
$\overline{\mathrm{row}}_{\mathrm{DT}}\cdot M_2F^n\hat V\subset M_{12}F^{n+1}\hat V$, 
$M_2F^n\hat V\cdot \overline{\mathrm{col}}_{\mathrm{DT}}\subset M_{12}F^n\hat V$, 
$M_{12}F^{n+1}\hat V\cdot M_2F^n\hat V\subset M_{12}F^{n+1}\hat V$ and
$M_2F^n\hat V\cdot M_{12}F^n\hat V\subset M_{12}F^n\hat V$ for $n\geq0$. 
(d) for $n=0$, this follows from $\mathbf c_{\mathrm{DT}}\in F^1\hat{\mathbf C}$ which implies 
$\mathbf c_{\mathrm{DT}}+F^{1}\hat{\mathbf C}=F^{1}\hat{\mathbf C}$. For $n\geq0$, this follows from the equality
\begin{align*}
    & \boldsymbol{\kappa}(\overline{\mathrm{row}}_{\mathrm{DT}}+\begin{pmatrix}
    \alpha & \beta
\end{pmatrix},\overline{\mathrm{col}}_{\mathrm{DT}}+\begin{pmatrix}
    a \\ b
\end{pmatrix})
    =\mathbf c_{\mathrm{DT}}+\Big[
[e_1\otimes a+\alpha\otimes 1+\alpha\otimes a]_{\mathrm C_V(e_0)}
\\& +[-f_1\otimes b-\beta\otimes 1+\beta\otimes b]_{\mathrm C_V(e_0)}
+[-f_1e_1\otimes a+\beta e_1\otimes 1+\beta e_1\otimes a]_{\mathbf k[e_0,f_0]}\Big]_{\mathbf C}
\end{align*}
and the fact that for $\alpha,\beta\in F^{n+1}\hat V$ and $a,b\in F^n\hat V$, the arguments of 
$[-]_{V\otimes_{\mathrm C_V(e_0)}V}$, $[-]_{V\otimes_{\mathrm C_V(f_0)}V}$ and
$[-]_{\mathbf k[e_0,f_0]}$ in the right-hand side respectively belong to 
$F^{n+1}(\hat V\hat\otimes\hat V)$, $F^{n+1}(\hat V\hat\otimes\hat V)$ and $F^{n+2}(\hat V\hat\otimes\hat V)$. 
\end{proof}

\begin{lem}\label{golda}
For $C\in V$, the relation $[1\otimes C-C\otimes 1]_{\mathrm C_V(e_0)}=0$
(equality in $V\otimes_{\mathrm C_V(e_0)}V$) implies $C\in \mathrm C_V(e_0)$.     
\end{lem}

\begin{proof}
Let $C\in V$. Let $\{e_0,e_1\}^*_1$ be the set of words in $e_0,e_1$ whose initial and 
final letters (which may coincide) are both $e_1$. Let $(C_\alpha)_{\alpha\geq0}$, $(C_{\alpha,w,\beta})_{\alpha,\beta\geq0,
w\in \{e_0,e_1\}^*_1}$ be the elements of $\mathcal V$ such that 
$$
C=\sum_{\alpha_\geq0}e_0^\alpha\otimes C_\alpha
+\sum_{\alpha,\beta\geq0,w\in \{e_0,e_1\}^*_1}e_0^\alpha we_0^\beta\otimes C_{\alpha,w,\beta},  
$$
then the image of $[1\otimes C-C\otimes 1]_{\mathrm C_V(e_0)}$ in 
$(\mathcal W_r\otimes\mathbf k[e_0]\otimes\mathcal W_l)\otimes\mathcal V$ is 
$$
\sum_{\alpha,\beta\geq0,w\in \{e_0,e_1\}^*_1}(1\otimes e_0^\alpha \otimes we_0^\beta-e_0^\alpha w\otimes 
e_0^\beta\otimes1)\otimes C_{\alpha,w,\beta}
$$
whose image with the tensor product of the counit $\mathcal W_r\to\mathbf k$, of the product map 
$\mathbf k[e_0]\otimes\mathcal W_l\to\mathcal V$ and the identity of $\mathcal V$, is 
$\sum_{\alpha,\beta\geq0,w\in \{e_0,e_1\}^*_1}e_0^\alpha we_0^\beta\otimes C_{\alpha,w,\beta}$. 
The relation $[1\otimes C-C\otimes 1]_{\mathrm C_V(e_0)}=0$ therefore implies the vanishing of this element, 
therefore $C=\sum_{\alpha_\geq0}e_0^\alpha\otimes C_\alpha$, hence $C\in\mathrm C_V(e_0)$.    
\end{proof}

\begin{lem}\label{lem:10:27:3003}
(a)   If $\tilde\alpha,\tilde a,\tilde b\in V$ are such that  
  $$
  \Big[[e_1\otimes \tilde a+\tilde\alpha\otimes 1]_{\mathrm C_V(e_0)}
+[-f_1\otimes \tilde b]_{\mathrm C_V(f_0)}
+[-f_1e_1\otimes \tilde a]_{\mathbf k[e_0,f_0]}\Big]_{\mathbf C}=0
  $$
  (relation in $\mathbf C$), then there exists $C\in \mathrm C_V(e_0)$ such that 
  $$
  (\tilde\alpha,\tilde a,\tilde b)=(-e_1(e_0+f_\infty)C,(e_0-f_0)C,e_1 C). 
  $$

  (b) If moreover $n\geq 1$ and $\tilde\alpha\in V_{n+1}$ and $\tilde a,\tilde b\in V_n$, then $C$ has degree $n-1$. 
\end{lem}

\begin{proof}
(a)    Let $\tilde \alpha,\tilde a,\tilde b\in V$. Let $(\tilde a_{kl})_{k,l\geq 0}$ be the elements of 
$\mathcal W_r^{\otimes 2}$ such that $\tilde a=\sum_{k,l\geq 0}e_0^kf_0^l\tilde a_{kl}$. Set also 
$$
C(\tilde a):=\sum_{k,l\geq0}{e_0^k-f_0^k\over e_0-f_0}f_0^l\tilde a_{kl}\in V, 
$$
so that 
\begin{equation}\label{eq:DE:tildea}
\tilde a=(e_0-f_0)C(\tilde a)+\sum_{k,l\geq 0}f_0^{k+l}\tilde a_{k,l}. 
\end{equation}
 Then the image of $-f_1e_1\otimes \tilde a$ in 
$V\otimes_{\mathbf k[e_0,f_0]}V\simeq \mathcal W_r^{\otimes 2}\otimes\mathbf k[e_0,f_0]
\otimes \mathcal W_l^{\otimes 2}$ is 
\begin{align*}
& -f_1e_1\otimes (\sum_{k,l\geq 0}e_0^kf_0^l\tilde a_{kl})
\simeq 
-\sum_{k,l\geq 0}f_1e_1\otimes e_0^kf_0^l\otimes \tilde a_{kl}
\\& =
-\sum_{k,l\geq 0}f_1e_1\otimes f_0^{k+l}\otimes \tilde a_{kl}
-\sum_{k,l\geq 0}f_1e_1\otimes (e_0-f_0){e_0^k-f_0^k\over e_0-f_0}f_0^l\otimes \tilde a_{kl},  
\end{align*}
therefore the image of 
\begin{align*}
& \Big[[e_1\otimes \tilde a+\tilde\alpha\otimes 1]_{\mathrm C_V(e_0)}
+[-f_1\otimes \tilde b]_{\mathrm C_V(f_0)}
+[-f_1e_1\otimes \tilde a]_{\mathbf k[e_0,f_0]}\Big]_{\mathbf C}\in 
\mathbf C\simeq \underline{\mathbf C}
\\&
=((\mathcal W_r\otimes \mathbf k[e_0]\otimes\mathcal W_l)\otimes\mathcal V)
\oplus(\mathcal V\otimes (\mathcal W_r\otimes \mathbf k[e_0]\otimes\mathcal W_l))
\oplus (\mathcal W_r^{\otimes2}\otimes\mathbf k[f_0]\otimes \mathcal W_l^{\otimes2}),
\end{align*}
is such that its third component is
$$
-\sum_{k,l\geq 0}f_1e_1\otimes f_0^{k+l}\otimes \tilde a_{kl}\in \mathcal W_r^{\otimes2}\otimes\mathbf k[f_0]\otimes \mathcal W_l^{\otimes2}. 
$$

If now $(\begin{pmatrix}
    \tilde\alpha & 0
\end{pmatrix},\begin{pmatrix}
    \tilde a \\ \tilde b
\end{pmatrix})$ belongs to $\mathrm{ker}(M_{12}V_{n+1}\times M_{21}V_n\to\mathbf C_{n+1})$, this element is 0, which 
together with \eqref{eq:DE:tildea} gives 
\begin{equation}\label{expr:tilde:a}
 \tilde a=(e_0-f_0)C(\tilde a).    
\end{equation}
The assumption on $(\begin{pmatrix}
    \tilde\alpha & 0
\end{pmatrix},\begin{pmatrix}
    \tilde a \\ \tilde b
\end{pmatrix})$ then implies the first equality in 
\begin{align*}
&0=\Big[[e_1\otimes \tilde a+\tilde\alpha\otimes 1]_{\mathrm C_V(e_0)}
+[-f_1\otimes \tilde b]_{\mathrm C_V(f_0)}
+[-f_1e_1\otimes \tilde a]_{\mathbf k[e_0,f_0]}\Big]_{\mathbf C}
\\& =\Big[[e_1\otimes (e_0-f_0)C(\tilde a)+\tilde\alpha\otimes 1]_{\mathrm C_V(e_0)}
+[-f_1\otimes \tilde b]_{\mathrm C_V(f_0)}
+[-f_1e_1\otimes (e_0-f_0)C(\tilde a)]_{\mathbf k[e_0,f_0]}\Big]_{\mathbf C} 
\\&=\Big[[e_1\otimes (e_0-f_0)C(\tilde a)+\tilde\alpha\otimes 1-f_1e_1\otimes C(\tilde a)]_{\mathrm C_V(e_0)}
+[-f_1\otimes \tilde b+f_1e_1\otimes C(\tilde a)]_{\mathrm C_V(f_0)}
\Big]_{\mathbf C} 
\\&=\Big[[e_1(e_0+f_\infty)\otimes C(\tilde a)+\tilde\alpha\otimes 1]_{\mathrm C_V(e_0)}
+[f_1\otimes (e_1C(\tilde a)-\tilde b)]_{\mathrm C_V(f_0)}
\Big]_{\mathbf C} 
\end{align*}
where the second equality follows from \eqref{expr:tilde:a}, the third equality follows from 
Def. \ref{def:CC:3003}, the fourth equality follows from the commutativity of $f_1$ with $e_1$ and from 
$e_0-f_0\in \mathrm C_V(e_0)$ and from $e_1\in\mathrm C_V(f_0)$. The composed map 
$\mathbf T^{\mathbf e}\oplus\mathbf T^{\mathbf e}\hookrightarrow \mathbf T\to \mathbf C\simeq\underline{\mathbf C}$
as injective as it coincides with the injection in the two first summands of $\underline{\mathbf C}$. 
It follows that 
\begin{equation}\label{vissotsky}
[e_1(e_0+f_\infty)\otimes C(\tilde a)+\tilde\alpha\otimes 1]_{\mathrm C_V(e_0)}=0,\quad 
[f_1\otimes (e_1C(\tilde a)-\tilde b)]=0.     
\end{equation}
(equalities in $(V\otimes V)_{\mathrm C_V(e_0)}$ and $(V\otimes V)_{\mathrm C_V(f_0)}$ respectively). 
The product map $V\otimes V\to V$ induces linear maps $V\otimes_X V\to V$ for $X$ equal to $\mathrm C_V(e_0)$ or $\mathrm C_V(f_0)$. 
The images of these equalities by these maps yield 
\begin{equation}\label{chertok}
    \tilde\alpha=-e_1(e_0+f_\infty)C(\tilde a),\quad f_1\cdot (e_1C(\tilde a)-\tilde b)=0. 
\end{equation}
(equalities in $V$). Combining with the latter equality the injectivity of the endomorphism $x\mapsto f_1x$ of $V$ then yields 
\begin{equation}\label{expression:tilde:b}
\tilde b=e_1 C(\tilde a),     
\end{equation}
and combining the former equality with the first equality of \eqref{vissotsky} gives 
\begin{equation}\label{okujava}
0=[e_1(e_0+f_\infty)\otimes C(\tilde a)-e_1(e_0+f_\infty)C(\tilde a)\otimes 1]_{\mathrm C_V(e_0)}\in (V\otimes V)_{\mathrm C_V(e_0)}.     
\end{equation}
The isomorphism 
\begin{equation}\label{akhmatova}
(V\otimes V)_{\mathrm C_V(e_0)}\simeq (\mathcal W_r\otimes\mathbf k[e_0]\otimes \mathcal W_l)\otimes \mathcal V    
\end{equation}
intertwines the endomorphism of $(V\otimes V)_{\mathrm C_V(e_0)}$ given by left multiplication by $e_1$ with the tensor product 
of left multiplication by $e_1$ in $\mathcal W_r$ by the identity of $\mathbf k[e_0]\otimes \mathcal W_l\otimes \mathcal V$, 
which is injective, therefore  left multiplication by $e_1$ on $(V\otimes V)_{\mathrm C_V(e_0)}$ is injective; 
and the isomorphism $(V\otimes V)_{\mathrm C_V(e_0)}\simeq (\mathcal V\otimes \mathcal W_l)\otimes \mathcal V$
obtained by composing \eqref{akhmatova} with the tensor product of the inverse of the isomorphism 
$\mathcal W_r\otimes\mathbf k[e_0]\to\mathcal V$ with $id_{\mathcal W_r\otimes\mathcal V}$ intertwines
the endomorphism of $(V\otimes V)_{\mathrm C_V(e_0)}$ given by left multiplication by $e_0+f_\infty$
with the tensor product of the identity of $\mathcal W_l$ with right mutplication by $e_0\otimes1+1\otimes e_\infty$
in $\mathcal V\otimes\mathcal V$, which is injective, therefore 
the endomorphism of $(V\otimes V)_{\mathrm C_V(e_0)}$ given by left multiplication by $e_0+f_\infty$
is injective. \eqref{okujava} therefore implies 
$$
[1\otimes C(\tilde a)-C(\tilde a)\otimes 1]_{\mathrm C_V(e_0)}=0
$$
(equality in $V\otimes_{\mathrm C_V(e_0)}V$). Lem. \ref{golda} then implies $C(\tilde a)\in \mathrm C_V(e_0)$.
The claim then follows from this together with \eqref{expr:tilde:a}, the first equality in \eqref{chertok} and 
\eqref{expression:tilde:b}. 

(b) follows from (a) and from the fact that the endomorphism $C\mapsto e_1C$ of $V$ is injective and of degree 1. 
\end{proof}

\begin{lem}\label{lem:10:28:3003}
    For any $n\geq 1$, the sequence 
\begin{equation}\label{complex:E:todo}
    \mathbf kf_0^n\oplus \mathrm C_V(e_0)_{n-1}\to M_{12}V_{n+1}\times M_{21}V_n
        \stackrel{\partial}{\to} \mathbf C_{n+1}, 
\end{equation}
where the first map is 
$$
(\Pi,C)\mapsto (-\overline{\mathrm{row}}_{\mathrm{DT}}\cdot \overline X(\Pi,C),\overline X(\Pi,C)\cdot 
\overline{\mathrm{col}}_{\mathrm{DT}})=(\begin{pmatrix}
    -e_1\Pi-e_1(e_0+f_\infty)C & f_1\Pi
\end{pmatrix},\begin{pmatrix}
    \Pi+(e_0-f_0)C\\-\Pi+e_1C
\end{pmatrix})
$$
(see \eqref{map:overlineM:Pi:m})
and the second map is 
$$
(\begin{pmatrix}
    \alpha & \beta
\end{pmatrix},\begin{pmatrix}
    a \\ b
\end{pmatrix})\mapsto 
\Big[
[e_1\otimes a+\alpha\otimes 1]_{\mathrm C_V(e_0)}
+[-f_1\otimes b-\beta\otimes 1]_{\mathrm C_V(f_0)}
+[-f_1e_1\otimes a+\beta e_1\otimes 1]_{\mathbf k[e_0,f_0]}\Big]_{\mathbf C}
$$
is an acyclic complex. 
\end{lem}

\begin{proof}
Let $(\Pi,C)\in \mathbf kf_0^n\oplus \mathrm C_V(e_0)_{n-1}$ and let $(\begin{pmatrix}
    \alpha & \beta
\end{pmatrix},\begin{pmatrix}
    a \\ b
\end{pmatrix})$ be the image of $(\Pi,C)$ by the first map. Then 
\begin{align*}
& [e_1\otimes a+\alpha\otimes 1]_{\mathrm C_V(e_0)}
=[e_1\otimes (\Pi+(e_0-f_0)C)+(-e_1\Pi-e_1(e_0+f_\infty)C)\otimes 1]_{\mathrm C_V(e_0)}
\\&=[e_1\otimes (e_0-f_0)C+(-e_1(e_0+f_\infty)C)\otimes 1]_{\mathrm C_V(e_0)}=
[e_1(e_0-f_0)\otimes C-e_1(e_0+f_\infty)\otimes C]_{\mathrm C_V(e_0)}
\\&=[e_1f_1\otimes C]_{\mathrm C_V(e_0)}
\end{align*}
where the second equality follows from $\Pi\in \mathrm C_V(e_0)$, the third equality follows from $C,e_0-f_0\in \mathrm C_V(e_0)$; moreover
$$
[-f_1\otimes b-\beta\otimes 1]_{\mathrm C_V(f_0)}
=[-f_1\otimes (-\Pi+e_1C)-f_1\Pi\otimes 1]_{\mathrm C_V(f_0)}
=-[e_1f_1\otimes C]_{\mathrm C_V(f_0)}
$$
where the second equality follows from $e_1,\Pi\in \mathrm C_V(f_0)$; and 
$$
[-f_1e_1\otimes a+\beta e_1\otimes 1]_{\mathbf k[e_0,f_0]}
=[-f_1e_1\otimes (\Pi+(e_0-f_0)C)+f_1\Pi e_1\otimes 1]_{\mathbf k[e_0,f_0]}
=[-e_1f_1\otimes (e_0-f_0)C]_{\mathbf k[e_0,f_0]}
$$
where the second equality follows from the commutation of $\Pi$ with $e_1$ and from $\Pi\in \mathbf k[e_0,f_0]$. 
Then 
\begin{align*}
&\Big[[e_1\otimes a+\alpha\otimes 1]_{\mathrm C_V(e_0)}
+[-f_1\otimes b-\beta\otimes 1]_{\mathrm C_V(f_0)}
+[-f_1e_1\otimes a+\beta e_1\otimes 1]_{\mathbf k[e_0,f_0]}\Big]_{\mathbf C}
\\&
=\Big[[e_1f_1\otimes C]_{\mathrm C_V(e_0)}-[e_1f_1\otimes C]_{\mathrm C_V(f_0)}+[-e_1f_1\otimes (e_0-f_0)C]_{\mathbf k[e_0,f_0]}\Big]_{\mathbf C}=0
\end{align*}
by Def. \ref{def:CC:3003}. This proves that the said sequence of maps is a complex. 

Let $(\begin{pmatrix}
    \alpha & \beta
\end{pmatrix},\begin{pmatrix}
    a \\ b
\end{pmatrix})\in M_{12}V_{n+1}\times M_{21}V_n$. 
Let $(a_{kl})_{k,l\geq 0}$ be the elements of $\mathcal W_r^{\otimes 2}$ such that $a=\sum_{k,l\geq 0}e_0^kf_0^la_{kl}$
and $(\beta_k)_{k\geq 0}$ be the elements of $\mathcal V\otimes \mathcal W_r$ such that $\beta=\sum_{k\geq 0}\beta_kf_0^k$.  
The image of $-f_1e_1\otimes a+\beta e_1\otimes 1$ in 
$V\otimes_{\mathbf k[e_0,f_0]}V\simeq \mathcal W_r^{\otimes 2}\otimes\mathbf k[e_0,f_0]
\otimes \mathcal W_l^{\otimes 2}$ is 
\begin{align*}
& -f_1e_1\otimes (\sum_{k,l\geq 0}e_0^kf_0^la_{kl})
+(\sum_{k\geq 0}\beta_kf_0^k) e_1\otimes 1\simeq 
-\sum_{k,l\geq 0}f_1e_1\otimes e_0^kf_0^l\otimes a_{kl}
+\sum_{k\geq 0}\beta_ke_1\otimes f_0^k\otimes 1
\\& =
\sum_{k\geq 0}\beta_ke_1\otimes f_0^k\otimes 1
-\sum_{k,l\geq 0}f_1e_1\otimes f_0^{k+l}\otimes a_{kl}
-\sum_{k,l\geq 0}f_1e_1\otimes (e_0-f_0){e_0^k-f_0^k\over e_0-f_0}f_0^l\otimes a_{kl},  
\end{align*}
therefore the image of 
\begin{align*}
&
(\begin{pmatrix}
    \alpha & \beta
\end{pmatrix},\begin{pmatrix}
    a \\ b
\end{pmatrix})\in M_{12}V_{n+1}\times M_{21}V_n
\to \mathbf C_{n+1} \subset \mathbf C\simeq \underline{\mathbf C}
\\&
=((\mathcal W_r\otimes \mathbf k[e_0]\otimes\mathcal W_l)\otimes\mathcal V)
\oplus(\mathcal V\otimes (\mathcal W_r\otimes \mathbf k[e_0]\otimes\mathcal W_l))
\oplus (\mathcal W_r^{\otimes2}\otimes\mathbf k[f_0]\otimes \mathcal W_l^{\otimes2})
\end{align*}
is such that its third component is
$$
\sum_{k\geq 0}\beta_ke_1\otimes f_0^k\otimes 1
-\sum_{k,l\geq 0}f_1e_1\otimes f_0^{k+l}\otimes a_{kl}\in \mathcal W_r^{\otimes2}\otimes\mathbf k[f_0]\otimes \mathcal W_l^{\otimes2}. 
$$
Assume now $(\begin{pmatrix}
    \alpha & \beta
\end{pmatrix},\begin{pmatrix}
    a \\ b
\end{pmatrix})\in \mathrm{ker}(M_{12}V_{n+1}\times M_{21}V_n\to\mathbf C_{n+1})$. This relation implies 
the collection of relations $\beta_ke_1\otimes 1=f_1e_1\otimes(\sum_{s+t=k}a_{st})$ for any $k\geq 0$, which 
implies $\beta_k\in\mathbf k f_1$, therefore $\beta\in f_1\mathbf k[f_0]$. The relation 
$\beta\in V_{n+1}$ then implies $\beta\in \mathbf kf_1f_0^n$, therefore $\beta\in f_1\Pi$ for some
$\Pi\in\mathbf kf_0^n$. 

Set 
\begin{equation}\label{def:tilde:alpha:etc}
\tilde\alpha:=\alpha+e_1\Pi,\quad\tilde a:=a-\Pi,\quad \tilde b:=b+\Pi.
\end{equation}  
It follows from the fact that $(\begin{pmatrix}-e_1\Pi & f_1\Pi\end{pmatrix},\begin{pmatrix}
    \Pi\\-\Pi
\end{pmatrix})\in\mathrm{ker}(M_{12}V_{n+1}\times M_{21}V_n\to \mathbf C_{n+1})$ that 
$$
(\begin{pmatrix}
    \tilde\alpha & 0
\end{pmatrix},\begin{pmatrix}
   \tilde a\\\tilde b 
\end{pmatrix})\in \mathrm{ker}(M_{12}V_{n+1}\times M_{21}V_n\to \mathbf C_{n+1}), 
$$
therefore 
$$
\Big[[e_1\otimes \tilde a+\tilde\alpha\otimes 1]_{\mathrm C_V(e_0)}
+[-f_1\otimes \tilde b]_{\mathrm C_V(f_0)}
+[-f_1e_1\otimes \tilde a]_{\mathbf k[e_0,f_0]}\Big]_{\mathbf C}=0. 
$$
Lem. \ref{lem:10:27:3003} then implies the existence of $C\in\mathrm C_V(e_0)_{n-1}$, such that 
$(\tilde\alpha,\tilde a,\tilde b)=(-e_1(e_0+f_\infty)C,(e_0-f_0)C,e_1 C)$. Together with \eqref{def:tilde:alpha:etc}, 
this implies that $(\begin{pmatrix}
    \alpha & \beta
\end{pmatrix},\begin{pmatrix}
    a \\ b
\end{pmatrix})$ belongs to the image of the first map of \eqref{complex:E:todo}, thus proving the claimed 
acyclicity. 
\end{proof}

\begin{lem}\label{lem:induct:kappa}
    Let $n\geq 0$. Then $x\in F^n\mathbf X$ and $\boldsymbol{\kappa}(x)\in \mathbf c_{\mathrm{DT}}+F^{n+2}\hat{\mathbf C}$
    implies the existence of $g\in F^n(\mathrm C_2(\overline\rho_0)^\times)$ such that $x\bullet g\in F^{n+1}\mathbf X$. 
\end{lem}

\begin{proof}
Let us prove that statement for $n=0$. Let $x\in \mathbf X$ be such that $\boldsymbol{\kappa}(x)\in \mathbf c_{\mathrm{DT}}+F^2\hat{\mathbf C}$. 
Let $\alpha,\beta\in F^1\hat V$ and $a,b\in \hat V$ be such that $x=(\begin{pmatrix}
    \alpha & \beta
\end{pmatrix},\begin{pmatrix}
    a\\ b
\end{pmatrix})$. Let $\alpha=\sum_{d\geq1}\alpha_d$,  $\beta=\sum_{d\geq1}\beta_d$,
$a=\sum_{d\geq0}a_d$,  $b=\sum_{d\geq0}b_d$ be the degree decompositions of $\alpha,\beta,a,b$. 
The assumption implies 
$$
 \Big[[\alpha_1\otimes a_0]_{\mathrm C_V(e_0)}+[\beta_1\otimes b_0]_{\mathrm C_V(f_0)}+
[\beta_1 e_1\otimes a_0]_{\mathbf k[e_0,f_0]}\Big]_{\mathbf C}
=\mathbf{c}_{\mathrm{DT}}  
$$
therefore
$$
 \Big[[\alpha_1\otimes a_0]_{\mathrm C_V(e_0)}+[\beta_1\otimes b_0]_{\mathrm C_V(f_0)}+
[\beta_1 e_1\otimes a_0]_{\mathbf k[e_0,f_0]}\Big]_{\mathbf C}
=\Big[[e_1\otimes 1]_{\mathrm C_V(e_0)}+[f_1\otimes 1]_{\mathrm C_V(f_0)}
-[e_1f_1\otimes 1]_{\mathbf k[e_0,f_0]}\Big]_{\mathbf C} 
$$
therefore (as $a_0,b_0\in\mathbf k$)
$$
 \Big[[(\alpha_1a_0-e_1)\otimes 1]_{\mathrm C_V(e_0)}+[(\beta_1 b_0-f_1)\otimes1]_{\mathrm C_V(f_0)}
 +[(\beta_1a_0+f_1) e_1\otimes 1]_{\mathbf k[e_0,f_0]}\Big]_{\mathbf C}
=0
$$
i.e., setting $\alpha_1=s_0e_0+s_1e_1+t_0f_0+t_1f_1$, $\beta_1=u_0e_0+u_1e_1+v_0f_0+v_1f_1$ with $s_i,t_i,u_i,v_i\in\mathbf k$,  
\begin{align*}
    & \Big[[(a_0s_0e_0+(a_0s_1-1)e_1+a_0t_0f_0+a_0t_1f_1)\otimes 1]_{\mathrm C_V(e_0)}
 +[(b_0u_0e_0+b_0u_1e_1+b_0v_0f_0+(b_0v_1-1)f_1)\otimes1]_{\mathrm C_V(f_0)}
    \\&
 +[(a_0u_0e_0+a_0u_1e_1+a_0v_0f_0+(a_0v_1+1)f_1) e_1\otimes 1]_{\mathbf k[e_0,f_0]}\Big]_{\mathbf C}
=0  .  
\end{align*}
 The image of the left-hand side by the isomorphism
 $\mathbf C\to\underline{\mathbf C}=((\mathcal W_r\otimes\mathbf k[e_0]\otimes\mathcal W_r)\otimes\mathcal V)\oplus 
 (\mathcal V\otimes (\mathcal W_r\otimes\mathbf k[e_0]\otimes\mathcal W_r))\oplus (\mathcal W_r^{\otimes2}\otimes 
 \mathbf k[f_0]\otimes 
 \mathcal W_l^{\otimes2})$ is 
\begin{align*}
    &\Big( a_0s_0(1\otimes e_0\otimes 1)\otimes 1+(a_0s_1-1)(e_1\otimes1\otimes1)\otimes1+(1\otimes1\otimes1)\otimes (a_0t_0e_0+a_0t_1e_1)\Big)
    \\&\oplus\Big((b_0u_0e_0+b_0u_1e_1)\otimes(1\otimes1\otimes1)+b_0v_0\cdot 1\otimes (1\otimes e_0\otimes 1)
    +(b_0v_1-1)1\otimes (e_1\otimes 1\otimes 1)\Big)
    \\&\oplus\Big(((a_0u_0e_0+a_0u_1e_1)e_1\otimes1)\otimes 1\otimes 1^{\otimes2}
    +a_0v_0\cdot 1^{\otimes 2}\otimes f_0\otimes1^{\otimes 2}+(a_0v_1+1)(e_1\otimes e_1)\otimes 1\otimes 1^{\otimes2}\Big).  
\end{align*}
Since the vectors arising in this expression form a free family, their coefficients are zero, 
therefore $1=a_0s_1=b_0v_1=-a_0v_1$, therefore $a_0=-b_0\in\mathbf k^\times$, therefore
$s_0=t_0=t_1=0=v_0=u_0=u_1$. Therefore 
$$
\begin{pmatrix}
    \alpha & \beta
\end{pmatrix}\in \begin{pmatrix}
    e_1 & -f_1
\end{pmatrix}a_0^{-1}+M_{12}F^2\hat V,\quad \begin{pmatrix}
    a\\ b
\end{pmatrix}\in a_0\begin{pmatrix}
    1\\ -1
\end{pmatrix}+M_{21}F^1\hat V, 
$$
which implies the claim with $g=a_0I_2\in \mathrm C_2(\overline\rho_0)^\times$. 

Assume now $n>0$. Let $x=(\mathrm{row},\mathrm{col})\in F^n\mathbf X$. 
Then the degree expansions of $\mathrm{row},\mathrm{col}$ are 
$\mathrm{row}=\overline{\mathrm{row}}_{\mathrm{DT}}+\sum_{d\geq n+1}\mathrm{row}_d$ and 
$\mathrm{col}=\overline{\mathrm{col}}_{\mathrm{DT}}+\sum_{d\geq n}\mathrm{col}_d$, and the degree expansion of
$\boldsymbol{\kappa}(x)$ is $\boldsymbol{\kappa}(x)=\mathbf c_{\mathrm{DT}}+\partial(\mathrm{row}_{n+1},\mathrm{col}_n)
+F^{n+2}\hat{\mathbf C}$, where the two first summands have degrees $1,n+1$ and $\partial$ is as in \eqref{complex:E:todo}. 
If one further assumes $\boldsymbol{\kappa}(x)\in \mathbf c_{\mathrm{DT}}+F^{n+2}\hat{\mathbf C}$, this implies 
$\partial(\mathrm{row}_{n+1},\mathrm{col}_n)=0$. By Lem. \ref{lem:10:28:3003}, this implies the existence of 
$(\Pi,C)\in\mathbf kf_0^n\oplus \mathrm C_V(e_0)_{n-1}$, such that 
$(\mathrm{row}_{n+1},\mathrm{col}_n)= (-\overline{\mathrm{row}}_{\mathrm{DT}}\cdot \overline X(\Pi,C),\overline X(\Pi,C)\cdot 
\overline{\mathrm{col}}_{\mathrm{DT}})$. Then $\overline X(\Pi,C)\in \mathrm C_2(\overline\rho_0)_n$, therefore 
$g:=I_2+\overline X(\Pi,C)\in F^n(\mathrm C_2(\overline\rho_0)^\times)$. 
Then $x\bullet g=(\mathrm{row}\cdot g,g^{-1}\cdot \mathrm{col})\in (\overline{\mathrm{row}}_{\mathrm{DT}}+M_{12}F^{n+2}\hat V)
\times(\overline{\mathrm{col}}_{\mathrm{DT}}+M_{21}F^{n+1}\hat V)=F^{n+1}\mathbf X$. 
\end{proof}

 \begin{lem}\label{lem:loc:inj:overlinec}
     The morphism of pointed sets induced by $\overline{\boldsymbol{\kappa}}$ (see \eqref{map:MM/C:Q}) and \eqref{underlinec:effect:on:marked:points}
     is locally injective. 
 \end{lem}

\begin{proof}
Let $\alpha\in \mathbf X/\mathrm C_2(\overline\rho_0)^\times$ be such that $\overline{\boldsymbol{\kappa}}(\alpha)
=\mathbf{c}_{\mathrm{DT}}$. Let $x\in \mathbf X$ be a representative of $\alpha$, then $\boldsymbol{\kappa}(x)=\mathbf{c}_{\mathrm{DT}}$. 
We construct inductively on $n\geq 0$ a sequence $(g_n)_{n\geq0}$, with $g_n\in F^n\mathrm C_2(\overline\rho_0)^\times$
for any $n$, such that the sequence $(x_n)_{n\geq0}$ defined by $x_0:=x$ and $x_{n+1}:=x_n\bullet g_n$ for any $n\geq0$ is such that 
$x_n\in F^n\mathbf X$ for any $n$. 
Indeed, the existence of $g_0$ follows from Lem. \ref{lem:induct:kappa} for $n=0$, and for any $n\geq0$, 
the existence of $g_{n+1}$ given $g_0,\ldots,g_n$ follows from Lem. \ref{lem:induct:kappa}
and from $\boldsymbol{\kappa}(x_n)=\boldsymbol{\kappa}(x\bullet (g_0\cdots g_n))=\boldsymbol{\kappa}(x)=\mathbf c_{\mathrm{DT}}$, 
which follows from the right $\mathrm C_2(\overline\rho_0)^\times$-invariance of $\boldsymbol{\kappa}$.

Then the sequence $(\gamma_n)_{n\geq 0}$ defined by $\gamma_n:=g_0\cdots g_n$ has a limit $\gamma\in \mathrm C_2(\overline\rho_0)^\times$, 
and for any $n\geq0$, $x\bullet\gamma=x_n\bullet(g_ng_{n+1}\cdots)$ where $g_ng_{n+1}\cdots\in F^n(\mathrm C_2(\overline\rho_0)^\times)$, 
therefore $x\in F^n\mathbf X$, therefore $x\bullet\gamma\in\cap_{n\geq0}F^n\mathbf X$, hence 
$x\bullet\gamma=(\overline{\mathrm{row}}_{\mathrm{DT}},\overline{\mathrm{col}}_{\mathrm{DT}})$, therefore 
$x\in (\overline{\mathrm{row}}_{\mathrm{DT}},\overline{\mathrm{col}}_{\mathrm{DT}})\bullet \mathrm C_2(\overline\rho_0)^\times$, 
therefore $\alpha=(\overline{\mathrm{row}}_{\mathrm{DT}},\overline{\mathrm{col}}_{\mathrm{DT}})\bullet \mathrm C_2(\overline\rho_0)^\times$,
which implies the claim. 
\end{proof}

\begin{lem}\label{quot:loc:inj:generalite}
    If $G$ is a group and $(X,x_0)\to (Y,y_0)$ is a $G$-equivariant map of pointed sets which is locally injective, then 
the map $(G\backslash X,Gx_0)\to(G\backslash Y,Gy_0)$ is locally injective.
\end{lem} 

\begin{proof}
Let $\alpha\in G\backslash X$ belong to the preimage of $Gy_0$. Choose $x\in \alpha$. The image of $x\in X\to Y$ belongs to 
$Gy_0$, therefore there exists $g\in G$, such that this image is $gy_0$. It follows that the image of $g^{-1}x\in X\to Y$ is 
$y_0$ which by the local injectivity of $X\to Y$ implies $g^{-1}x=x_0$. Therefore $x=gx_0$. Then $\alpha=Gx=Gx_0$. It follows that 
the preimage of $Gy_0$ by $G\backslash X\to G\backslash Y$ is $\{Gx_0\}$, therefore that $(G\backslash X,Gx_0)\to(G\backslash Y,Gy_0)$ 
is locally injective.  
\end{proof}

\begin{lem}\label{lem:E3:loc:inv}
    The map (E3) is locally injective. 
\end{lem}

\begin{proof}
This follows from Lem. \ref{quot:loc:inj:generalite}, applied to the map $\overline{\boldsymbol{\kappa}}$ and to the action of the group 
$\mathbf k[[u,v]]^\times$, using Lem. \ref{lem:loc:inj:overlinec}. 
\end{proof}

\subsection{Local injectivity of the morphism (E)}\label{sect:10:9}

\begin{lem}
(a)  The diagram \eqref{diag:E} is commutative.  

(b) The map (E) is equal to the composition $(E5)\circ (E4)\circ(E3)\circ(E2)\circ(E1)^{-1}$.
\end{lem}

\begin{proof}
(a)  Let $\alpha\in \mathrm{GL}_2\hat V/\mathrm C_2(\overline\rho_0)^\times$ and $P\in \mathrm{GL}_2\hat V$ be a representative of 
$\alpha$. 

The image of $\alpha$ by the map from Lem. \ref{lem:comm:diag:toto}(b) is $\sigma_P\in 
\mathrm{Hom}_{\mathcal C\operatorname{-alg}}^{1,(0)}(\hat{\mathcal V},M_2\hat V)$ such that 
$e_1\mapsto \overline\rho_1$, $e_0\mapsto \mathrm{Ad}_P(\overline\rho_0)$. The image of the latter element
by the map from Lem. \ref{lem:65:7:2412}(b) is the element $\Delta_{\sigma_P}\in 
\mathrm{Hom}_{\mathcal C\operatorname{-alg}}(\hat{\mathcal W},\hat V)$ such that 
$\Delta_{\sigma_P}(e_0^{n-1}e_1)=\overline{\mathrm{row}}_{\mathrm{DT}}\cdot \sigma_P(e_0)^{n-1}
\cdot \overline{\mathrm{col}}_{\mathrm{DT}}$ for any $n\geq1$; $\sigma_P(e_0)=\mathrm{Ad}_P(\overline\rho_0)$ then 
implies
\begin{equation}\label{eq:Delta:rho:g}
\forall n\geq 1, \quad 
\Delta_{\sigma_P}(e_0^{n-1}e_1)=\overline{\mathrm{row}}_{\mathrm{DT}}\cdot P\cdot \overline\rho_0^{n-1}
\cdot P^{-1}\cdot  \overline{\mathrm{col}}_{\mathrm{DT}}.     
\end{equation}
The image of $\alpha$ by the map $\mathrm{GL}_2\hat V/\mathrm C_2(\overline\rho_0)^\times\to 
(M_{12}F^1\hat V\times M_{21}\hat V)/\mathrm C_2(\overline\rho_0)^\times$ induced by 
Lems. \ref{lem:10:3:7mars}(b) and \ref{lem:10:4:7mars}(b) is the class in 
$(M_{12}F^1\hat V\times M_{21}\hat V)/\mathrm C_2(\overline\rho_0)^\times$ of the pair
$(\overline{\mathrm{row}}_{\mathrm{DT}}\cdot P,P^{-1}\cdot \overline{\mathrm{col}}_{\mathrm{DT}})$.
The image of the latter element by the map from Lem. \ref{lem:invce:eqvce:cc} is the element 
$\Big[[\alpha\otimes a]_{\mathrm C_V(e_0)}+[\beta\otimes b]_{\mathrm C_V(f_0)}+
[\beta e_1\otimes a]_{\mathbf k[e_0,f_0]}\Big]_{\mathbf C}\in \hat{\mathbf C}$, where $\alpha,\beta\in F^1\hat V$, $a,b\in\hat V$ 
are such that 
\begin{equation}\label{conventions:abalphabeta}
 (\overline{\mathrm{row}}_{\mathrm{DT}}\cdot P,P^{-1}\cdot \overline{\mathrm{col}}_{\mathrm{DT}})
=(\begin{pmatrix}
    \alpha&\beta
\end{pmatrix},\begin{pmatrix}
    a\\b
\end{pmatrix}).    
\end{equation}
The image of the latter element by the map from Lem. \ref{lem:10:11:BIS}(d)
is the element $\delta:=\delta^e_{\alpha\otimes a}+\delta^f_{\beta\otimes b}+\delta^{ef}_{\beta e_1\otimes a}
\in\prod_{n\geq1}F^n\hat V$, given by 
$$
\forall n\geq1,\quad \delta(n)=\alpha e_0^{n-1}a+\beta f_0^{n-1}b+\beta e_1{e_0^{n-1}-f_0^{n-1}\over e_0-f_0}a, 
$$
and the image of the latter element by the map from Lem. \ref{lem:10:12:1603} is $\Delta_\delta$. 
For any $n\geq1$, one has 
\begin{align*}
    & \Delta_\delta(e_0^{n-1}e_1)=\delta(n)=\alpha e_0^{n-1}a+\beta f_0^{n-1}b+\beta e_1{e_0^{n-1}-f_0^{n-1}\over e_0-f_0}a
=\begin{pmatrix}
    \alpha & \beta
\end{pmatrix}
\overline\rho_0^{n-1}
\begin{pmatrix}
    a\\b 
\end{pmatrix}
    \\& =\overline{\mathrm{row}}_{\mathrm{DT}}\cdot P\cdot \overline\rho_0^{n-1}
\cdot P^{-1}\cdot  \overline{\mathrm{col}}_{\mathrm{DT}}=\Delta_{\sigma_g}(e_0^{n-1}e_1)
\end{align*}
where the second equality follows from \eqref{eq:Delta:rho:g}, the third equality follows from 
$$
\overline\rho_0^{n-1}=\begin{pmatrix}
    e_0^{n-1}&0\\e_1{e_0^{n-1}-f_0^{n-1}\over e_0-f_0}&f_0^{n-1}
\end{pmatrix},
$$
and the fourth equality follows from \eqref{conventions:abalphabeta}. 
Since the family $(e_0^{n-1}e_1)_{n\geq1}$ generates $\hat{\mathcal W}$, the equality $\Delta_{\sigma_P}=\Delta_\delta$ follows. 

This implies that the diagram of set maps  
\begin{equation}\label{diag:pre:E}
\xymatrix{
\mathrm{GL}_2\hat V/\mathrm C_2(\overline\rho_0)^\times\ar_{\text{Lem. \ref{lem:comm:diag:toto}(b)}}[dd]
\ar^{\!\!\!\!\!\!\!\!\!\!\!\!\substack{\text{Lems. \ref{lem:10:3:7mars}(b)}\\ \text{and \ref{lem:10:4:7mars}(b)}}}[r]&
(M_{12}F^1\hat V\times M_{21}\hat V)/\mathrm C_2(\overline\rho_0)^\times\ar^{\ \ \ \ \ \ \ \ \ \ 
\text{Lem. \ref{lem:invce:eqvce:cc}}}[r]&
F^1\hat{\mathbf C}\ar^{\text{Lem. \ref{lem:10:11:BIS}(d)}}[d]\\
&&
\prod_{n\geq1}F^n\hat V\ar^{\text{Lem. \ref{lem:10:12:1603}}}[d]\\
\mathrm{Hom}_{\mathcal C\operatorname{-alg}}^{1,(0)}(\hat{\mathcal V},M_2\hat V)
\ar_{\text{Lem. \ref{lem:65:7:2412}(b)}}[rr]&&
\mathrm{Hom}_{\mathcal C\operatorname{-alg}}(\hat{\mathcal W},\hat V)}
\end{equation}
is commutative. 
It follows from Lem. \ref{lem:10:3:7mars}(a), Lem. \ref{lem:10:7:1603}(c), Lem. \ref{lem:10:11:BIS}(e), 
Lem. \ref{lem:actions:2912}(d) and Lem. \ref{lem28:1001}(c)  that each of the sets of \eqref{diag:pre:E} is 
equipped with an action of the corresponding 
group of the following commutative diagram of group morphisms
\begin{equation}\label{diag:gp:pre:E}
\xymatrix{
\mathrm C_2(\overline\rho_1)^\times\ar_{id}[dd]
\ar_{\text{Lem. \ref{lem:6:16:2912}(c)}}[r]&
\mathbf k[[u,v]]^\times\ar^{id}[r]&
\mathbf k[[u,v]]^\times\ar^{id}[d]\\
&&
\mathbf k[[u,v]]^\times\ar^{id}[d]\\
\mathrm C_2(\overline\rho_1)^\times
\ar_{\text{Lem. \ref{lem:6:16:2912}(c)}}[rr]&&
\mathbf k[[u,v]]^\times}
\end{equation}
and from Lem. \ref{lem:basic:E1}(c), Lem. \ref{lem:10:3:7mars}(b) and Lem. \ref{lem:10:4:7mars}(b), 
Lem. \ref{lem:invce:eqvce:cc}, Lem. \ref{lem:10:11:BIS}(e), Lem. \ref{lem:10:12:1603}, and Lem. \ref{lem:6:20},  
that each of the maps of \eqref{diag:pre:E} is  compatible
with the actions and with the corresponding group morphism of \eqref{diag:gp:pre:E}. 
The statement follows from the fact that \eqref{diag:E} is the corresponding diagram between 
orbit spaces. 

(b) follows from (a) and from the bijectivity of (E1), which follows from   
Def. \ref{defE1} and Lem.  \ref{lem:basic:E1}(c)
\end{proof}

\begin{prop}\label{prop:E:loc:inj}
The morphism of pointed sets 
$$
(E) : \mathrm C_2(\overline\rho_1)^\times\backslash\mathrm{Hom}_{\mathcal C\operatorname{-alg}}^{1,(0)}(\hat{\mathcal V},M_2\hat V)
\to \mathbf k[[u,v]]^\times\backslash\mathrm{Hom}_{\mathcal C\operatorname{-alg}}(\hat{\mathcal W},\hat V)
$$
(see \eqref{the:big:diagram}) is locally injective. 
\end{prop}

\begin{proof}
The diagram \eqref{diag:E} is upgraded to a diagram of morphisms of pointed sets by the adjunction of the following diagram of
elements 
$$
\xymatrix{
\mathrm C_2(\overline\rho_0)^\times\cdot I_2\cdot \mathrm C_2(\overline\rho_0)^\times\ar[r]\ar[dd] &
\mathbf k[[u,v]]^\times\bullet(\overline{\mathrm{row}}_{\mathrm{DT}},\overline{\mathrm{col}}_{\mathrm{DT}})\bullet 
\mathrm C_2(\overline\rho_0)^\times\ar[r] &
\mathbf k[[u,v]]^\times\bullet\mathbf c_{\mathrm{DT}}\ar[d] \\ 
&&\mathbf k[[u,v]]^\times\bullet(n\mapsto \Delta^{\mathcal W}_{r,l}(e_0^{n-1}e_1))\ar[d]\\
\mathrm C_2(\overline\rho_0)^\times\bullet \overline\rho_{\mathrm{DT}}\ar[rr]& & \mathbf k[[u,v]]^\times\bullet \Delta^{\mathcal W}_{r,l}}
$$
In this diagram, (E1) is bijective, (E2) is injective by Lem. \ref{lem:10:4:7mars}(c), 
(E3) is locally injective is Lem. \ref{lem:E3:loc:inv}, (E4) is injective by Lem. \ref{lem:E4:inj}, 
(E5) is a bijection by Lem. \ref{lem:10:12:1603}. The statement follows. 
\end{proof}

\section{Equality between $\mathrm{Stab}_{\mathcal G}(\mathrm{GL}_3\hat V\bullet \rho_{\mathrm{DT}}
)$ 
and $\mathsf{DMR}_0(\mathbf k)$}\label{sect:11}

\begin{thm}\label{thm:eq:Stab:Stab}
    The group inclusion from Thm. \ref{thm:5:31:3103} is an equality, therefore one has the equality 
    $$
\mathrm{Stab}_{\mathcal G}(\mathrm{GL}_3\hat V\bullet \rho_{\mathrm{DT}}
)
=\mathrm{Stab}_{\mathcal G}(\mathbf k[[u,v]]^\times\bullet\Delta^{\mathcal W}_{r,l}
)
    $$
of subgroups of $\mathcal G$. 
\end{thm}

\begin{proof}
By Lem. \ref{lem:6:30:3103}, \eqref{the:big:diagram} is a diagram of pointed sets with actions of $\mathcal G$, where the morphisms 
are denoted (B)-(E). It follows from Prop. \ref{prop:B:loc:inj} (resp. Prop. \ref{prop:D:loc:inj}, Prop. \ref{prop:E:loc:inj})
that (B) (resp. (D), (E)) is locally injective, and from Lem. \ref{lem:C:inj} that (C) 
is injective. By the first statement of Lem. \ref{lem:3010:1129}, this implies the equalities 
\begin{align*}
&
\mathrm{Stab}_{\mathcal G}(\mathrm{C}_3(\rho_1)^\times\bullet \rho_{\mathrm{DT}}
)
=\mathrm{Stab}_{\mathcal G}(\mathrm C_{21}(\rho_1)^\times\bullet \rho_{\mathrm{DT}})
=\mathrm{Stab}_{\mathcal G}((\mathrm C_{21}^{(0)}(\rho_1)^\times\bullet \rho_{\mathrm{DT}})
\\& 
=\mathrm{Stab}_{\mathcal G}(\mathrm C_2(\overline \rho_1)^\times\bullet \overline\rho_{\mathrm{DT}}
)
=\mathrm{Stab}_{\mathcal G}(\mathbf k[[u,v]]^\times\bullet\Delta^{\mathcal W}_{r,l}
),     
\end{align*}
where the intermediate terms are the stabilizer groups of the successive pointed $\mathcal G$-sets from 
\eqref{the:big:diagram}. The statement follows from the combination of this with the equality 
$\mathrm{Stab}_{\mathcal G}(\mathrm{C}_3(\rho_1)^\times\bullet \rho_{\mathrm{DT}}
)
=\mathrm{Stab}_{\mathcal G}(\mathrm{GL}_3\hat V\bullet \rho_{\mathrm{DT}}
)$ (see Cor. \ref{NEWCOR}). 
\end{proof}

\begin{cor}\label{cor:11:2:17apr} (see Thm. \ref{thm:016})
    The subgroups  $\mathrm{Stab}_{\mathcal G}(\mathrm{GL}_3\hat V\bullet \rho_{\mathrm{DT}}
    )$ 
and $\mathsf{DMR}_0(\mathbf k)$ of $\mathcal G$ are equal. 
\end{cor}

\begin{proof}
    This follows by combining Thm. \ref{thm:eq:Stab:Stab} and Thm. \ref{thm:main}. 
\end{proof}

\newpage

\part{Relationship of $\mathrm{Stab}_{\mathcal G}(\mathrm{GL}_3\hat V\bullet 
\rho_{\mathrm{DT}})$ with inertia-preserving automorphisms}\label{part 4}

The objective of Part \ref{part 4} is twofold: on the one hand, we prove the inclusion of 
$\mathrm{Stab}_{\mathcal G}(\mathrm{GL}_3\hat V\bullet\rho_{\mathrm{DT}})$ in $\mathcal G_{\mathrm{inert}}$ 
(\S \ref{sect:12}, Thm. \ref{thm:13:22:0205}) ; on the other hand, we prove the stability of 
$\mathrm{Stab}_{\mathcal G}(\mathrm{GL}_3\hat V\bullet\rho_{\mathrm{DT}})$ under the involution $\Theta$ 
of $\mathcal G_{\mathrm{inert}}$ (\S\ref{sect:13}, Thm. \ref{thm:13:36:17apr}).

\section{Inclusion of $\mathrm{Stab}_{\mathcal G}(\mathrm{GL}_3\hat V\bullet\rho_{\mathrm{DT}})$ in 
$\mathcal G_{\mathrm{inert}}$}\label{sect:12}

In this section, we prove the inclusion of $\mathrm{Stab}_{\mathcal G}(\mathrm{GL}_3\hat V \bullet \rho_{\mathrm{DT}})$ 
in $\mathcal G_{\mathrm{inert}}$ (Thm. \ref{thm:13:22:0205}). 
We first prove the statements of the Introduction on $\mathcal G_{\mathrm{inert}}$ and its 
automorphism $\Theta$ (\S\ref{sect:25jan1356}). We formulate the definition of the former group in 
terms of a subset $\mathrm{Tor}_g$ of $\mathrm{GL}_3\hat V$ (§\ref{sect:12:1}); the main result 
is obtained in §\ref{sect:12:2} as the consequence of the following steps: in Lem. \ref{lem:13:9:1301}, we 
characterize the "generalized eigen(co-)vectors" corresponding to a character $\chi$ (Def. \ref{def:chi:3001}) in terms 
of row and column matrices $R_{\mathrm{DT}}$ and $C_{\mathrm{DT}}$, which will be shown to be 
the building blocks of the centralizer algebra $C_3(\rho_{\mathrm{DT}}(\hat{\mathcal V}))$ (see Lem. \ref{lem:13:8}). 
In Lem. \ref{relations:ROW:1301} (resp. Lem. \ref{rels:COL:1301}), we derive behavior of the $\mathrm{Tor}_g$ with 
respect to $R_{\mathrm{DT}}$ (resp. $C_{\mathrm{DT}}$) under the assumption 
$g\in\mathrm{Stab}_{\mathcal G}(\mathrm{GL}_3\hat V \bullet \rho_{\mathrm{DT}})$, and in Lem. \ref{lem:12:8},
we derive from this and from $R_{\mathrm{DT}}\cdot C_{\mathrm{DT}}=f_0+f_1-e_0$ an identity relating $g$ with $f_0+f_1-e_0$. 
Together with Prop. \ref{prop:Z} (to be proved in Appendix \ref{appC}), this implies Thm. \ref{thm:13:22:0205}.

\subsection{The group $\mathcal G_{\mathrm{inert}}$ and its involution $\Theta$}
\label{sect:25jan1356}

\begin{lem}\label{lem02:2506} (see Lem. \ref{lem02:2506:wo:proof})
    (a) If $g\in \mathcal G$ is such that there exists $h\in \mathcal G$ such that 
$\mathrm{Ad}_g(e_0)+e_1+\mathrm{Ad}_h(e_\infty)=0$ (equality in $\mathfrak{lie}_{\{0,1\}}^\wedge$), then 
$h$ is unique; it will be denoted $h_g$. 

(b) The subset $\mathcal G_{\mathrm{inert}}\subset  \mathcal G$ of all elements 
$g$ as in (a) is a subgroup of 
$(\mathcal G,\circledast)$.

(c) There is a unique automorphism $s_{(0,\infty)}$ of 
$\mathfrak{lie}_{\{0,1\}}^\wedge$, such that 
$e_1\mapsto e_1$ and $e_0\leftrightarrow e_\infty$; it is an involution. 

(d) The map $\Theta : g\mapsto s_{(0,\infty)}(h_g)$ is an involutive automorphism of $(\mathcal G_{\mathrm{inert}},\circledast)$.
\end{lem}





\begin{proof}
Let $(g,h)$ be as in (a) and let $h'\in\mathcal G$ be such that $\mathrm{Ad}_g(e_0)+e_1+\mathrm{Ad}_{h'}(e_\infty)=0$. Then 
$\mathrm{Ad}_{h'}(e_\infty)=\mathrm{Ad}_h(e_\infty)$, which implies the existence of $\upsilon\in\mathbf k$ such that 
$h'=h\cdot\mathrm{exp}(\upsilon\cdot e_\infty)$. This implies the middle equality in $0=(h'|e_0)=(h|e_0)+\upsilon=\upsilon$, where 
the first (resp. last) equality follows from $h\in\mathcal G$ (resp. $h'\in\mathcal G$). Therefore $h'=h$; this proves (a). 

Let $g,g' \in \mathcal G_{\mathrm{inert}}$. It follows from $\mathrm{Ad}_g(e_0)+e_1+\mathrm{Ad}_{h_g}(e_\infty)=0$ that there is 
a unique automorphism of $\mathfrak{lie}_{\{0,1\}}^\wedge$, such that $e_0\mapsto \mathrm{Ad}_g(e_0)$, 
$e_1\mapsto e_1$, and $e_\infty\mapsto\mathrm{Ad}_{h_g}(e_\infty)$. Applying this automorphism to the equality 
$\mathrm{Ad}_{g'}(e_0)+e_1+\mathrm{Ad}_{h_{g'}}(e_\infty)=0$, one gets 
$$
\mathrm{Ad}_{g\circledast g'}(e_0)+e_1+\mathrm{Ad}_{h_{g'}(\mathrm{Ad}_g(e_0),e_1)\cdot h_g(e_0,e_1)}(e_\infty)=0, 
$$
which implies $g\circledast g'\in \mathcal G_{\mathrm{inert}}$, with 
\begin{equation}\label{toto:3006}
h_{g \circledast g'}=h_{g'}(g(e_0,e_1) \cdot e_0 \cdot  
g(e_0,e_1)^{-1},e_1)\cdot h_g(e_0,e_1).     
\end{equation}

Let $g\in\mathcal G$. Its inverse for $\circledast$ is $g^{\circledast-1}:=(\mathrm{aut}_g^{\mathcal V})^{-1}(g^{-1})$, where 
$\mathrm{aut}_g^{\mathcal V}$ is the 
automorphism of $\mathfrak{lie}_{\{0,1\}}^\wedge$ given by $e_1 \mapsto e_1$, $e_0 \mapsto g(e_0,e_1)\cdot e_0\cdot g(e_0,e_1)^{-1}$;  
If in addition $g\in\mathcal G_{\mathrm{inert}}$, there exist uniquely defined $g',h'\in\mathcal G$ such that 
the automorphism $(\mathrm{aut}_g^{\mathcal V})^{-1}$ is given by $e_0\mapsto \mathrm{Ad}_{g'}(e_0)$, 
$e_1\mapsto e_1$, $e_\infty\mapsto \mathrm{Ad}_{h'}(e_\infty)$; in particular 
$\mathrm{Ad}_{g'}(e_0)+e_1+\mathrm{Ad}_{h'}(e_\infty)=0$. Then $\mathrm{aut}_g^{\mathcal V} \circ 
(\mathrm{aut}_g^{\mathcal V})^{-1}=id$ implies 
$\mathrm{aut}_g^{\mathcal V}(g')\cdot g=\mathrm{aut}_g^{\mathcal V}(h')\cdot h_g=1$. 
It follows that $g'=g^{\circledast-1}$, and $h'=(\mathrm{aut}_g^{\mathcal V})^{-1}(h_g)$.  
Then $\mathrm{Ad}_{g^{\circledast-1}}(e_0)+e_1+\mathrm{Ad}_{(\mathrm{aut}_g^{\mathcal V})^{-1}(h_g)}(e_\infty)=0$, 
which implies $g^{\circledast-1}\in \mathcal G_{\mathrm{inert}}$, and $h_{g^{\circledast-1}}=
(\mathrm{aut}_g^{\mathcal V})^{-1}(h_g)$. All this implies (b). (c) is obvious. 

Let $g\in\mathcal G_{\mathrm{inert}}$. Applying $s_{(0,\infty)}$ to the equality 
$\mathrm{Ad}_g(e_0)+e_1+\mathrm{Ad}_{h_g}(e_\infty)=0$, one obtains 
$\mathrm{Ad}_{s_{(0,\infty)}(h_g)}(e_0)+e_1+\mathrm{Ad}_{s_{(0,\infty)}(g)}(e_\infty)=0$, which 
implies that $\Theta(g):=s_{(0,\infty)}(h_g)\in\mathcal G_{\mathrm{inert}}$ and that 
$h_{\Theta(g)}=s_{(0,\infty)}(g)$. Then $\Theta^2(g)=s_{(0,\infty)}(h_{\Theta(g)})=s_{(0,\infty)}(s_{(0,\infty)}(g))=g$, which 
implies that $\Theta$ is involutive. Finally for $g,g'\in\mathcal G_{\mathrm{inert}}$, 
\begin{align*}
& \Theta(g)\circledast\Theta(g')=s_{(0,\infty)}(h_g)\circledast s_{(0,\infty)}(h_{g'})=
a(s_{(0,\infty)}(h_{g'}))\cdot s_{(0,\infty)}(h_g)
\\ & =
h_{g'}(g(e_\infty,e_1) \cdot e_\infty \cdot  
g(e_\infty,e_1)^{-1},e_1)\cdot h(e_\infty,e_1)
=s_{(0,\infty)}(h_{g\circledast g'})=\Theta(g\circledast g'),    
\end{align*}
where $a$ is the automorphism of 
$\mathfrak{lie}_{\{0,1\}}^\wedge$ such that $e_1\mapsto e_1$, $e_0\mapsto 
s_{(0,\infty)}(h_g)\cdot e_0\cdot s_{(0,\infty)}(h_g)^{-1}$, where the first and last equalities follows from the definition of 
$\Theta$, the second equality follows from the definition of $\circledast$, the third equality follows from the fact that $a$ is 
such that $e_\infty\mapsto s_{(0,\infty)}(g)\cdot e_\infty\cdot s_{(0,\infty)}(g)^{-1}$, which follows from applying 
$s_{(0,\infty)}$ to $\mathrm{Ad}_g(e_0)+e_1+\mathrm{Ad}_{h_g}(e_\infty)=0$, 
the fourth equality follows from \eqref{toto:3006}; 
this implies that $\Theta$ is a group automorphism. All this proves (d). 
\end{proof}

\subsection{The $\mathbf k$-module $\mathrm{Tor}_g$ and the set $\mathrm{Tor}_g^\times$}\label{sect:12:1}

\begin{defn}\label{def:12:1:sbg}
For $g\in\mathcal G$, define: 

(a) $\mathrm{Tor}_g\subset M_3\hat V$ as the subset of elements $P$ such that (following 
the notation of Lem. \ref{lem:520:2212:BIS}) 
$$
\forall x\in\hat{\mathcal V},\quad (g*\rho_{\mathrm{DT}})(x)\cdot P=P\cdot \rho_{\mathrm{DT}}(x).   
$$ 

(b) $\mathrm{Tor}_g^\times$ as the subset of $\mathrm{GL}_3\hat V$ of all elements $P$ such that 
   $g*\rho_{\mathrm{DT}}=\mathrm{Ad}_P\circ\rho_{\mathrm{DT}}$ (equality in 
   $\mathrm{Hom}_{\mathcal C\operatorname{-alg}}(\hat{\mathcal V},M_3\hat V)$). 
\end{defn}

One checks that: (a) $\mathrm{Tor}_g^\times=\mathrm{Tor}_g\cap\mathrm{GL}_3\hat V$;  
(b) $\mathrm{Tor}_g$ is a right $\mathrm C_3(\rho_{\mathrm{DT}}(\hat{\mathcal V}))$-submodule of $M_3\hat V$;  
(c) $\mathrm{Tor}_g^\times$ is a right $\mathrm C_3(\rho_{\mathrm{DT}}(\hat{\mathcal V}))^\times$-subset of $\mathrm{GL}_3\hat V$. 

\begin{lem}\label{lem:last:1706:BIS}
For $g\in\mathcal G$, the following statements are equivalent: 

(a) $g\in \mathrm{Stab}_{\mathcal G}(\mathrm{GL}_3\hat V\bullet \rho_{\mathrm{DT}}
)$; 

(b) $\mathrm{Tor}_g^\times\neq\emptyset$. 
\end{lem}

\begin{proof} (a) is equivalent to $g*(\mathrm{GL}_3\hat V\bullet \rho_{\mathrm{DT}})
=\mathrm{GL}_3\hat V\bullet \rho_{\mathrm{DT}}$ (equality in $\mathrm{GL}_3\hat V\backslash 
\mathrm{Hom}_{\mathcal C\operatorname{-alg}}(\hat{\mathcal V},M_3\hat V)$, where the action of 
$\mathcal G$ is induced by Lem. \ref{lem:pre:psga}(a)). The left-hand side is equal to 
$\mathrm{GL}_3\hat V\bullet(g*\rho_{\mathrm{DT}})$, therefore this equality is equivalent to 
$\mathrm{GL}_3\hat V\bullet(g*\rho_{\mathrm{DT}})=\mathrm{GL}_3\hat V\bullet \rho_{\mathrm{DT}}$. This is 
again equivalent to $g*\rho_{\mathrm{DT}}\in \mathrm{GL}_3\hat V\bullet \rho_{\mathrm{DT}}$, which is 
equivalent to the existence of $P\in \mathrm{GL}_3\hat V$ such that 
$g*\rho_{\mathrm{DT}}=\mathrm{Ad}_P\circ\rho_{\mathrm{DT}}$, i.e. to (b).
This proves the equivalence of (a) and (b). 
\end{proof}

\subsection{Inclusion of $\mathrm{Stab}_{\mathcal G}(\mathrm{GL}_3\hat V\bullet\rho_{\mathrm{DT}})$ in $\mathcal G_{\mathrm{inert}}$}\label{sect:12:2}

\begin{defn}\label{def:chi:3001}
Let $\chi : \hat{\mathcal V}\to \hat V=\hat{\mathcal V}^{\hat\otimes2}$ be the algebra morphism defined by 
$e_0\mapsto e_0=e_0\otimes1$, $e_1\mapsto 0$. 
\end{defn}

Recall from Def. \ref{defn:C:DT:R:DT} the elements $C_{\mathrm{DT}}\in M_{13}\hat V$, $R_{\mathrm{DT}}\in M_{31}\hat V$.

\begin{lem}\label{lem:13:9:1301}
(a) The map $a\mapsto C_{\mathrm{DT}}\cdot a$ defines a bijection $\mathrm{C}_{\hat V}(e_0)\to \{v\in M_{31}\hat V|
\forall x\in\hat{\mathcal V}, \rho_{\mathrm{DT}}(x)v=v\chi(x)\}$. 

(b) The map $a\mapsto a\cdot R_{\mathrm{DT}}$ defines a bijection $\mathrm{C}_{\hat V}(e_0)\to \{\xi\in M_{13}\hat V|
\forall x\in\hat{\mathcal V},\xi\rho_{\mathrm{DT}}(x)=\chi(x)\xi\}$. 
\end{lem}

\begin{proof}
(a)  Set $\mathcal C:=\{v\in M_{31}\hat V|
\forall x\in\hat{\mathcal V}, \rho_{\mathrm{DT}}(x)v=v\chi(x)\}$. It follows from the fact that 
the image of $\chi$ is contained in $\mathbf k[[e_0]]\subset \hat V$ that
\begin{equation}\label{stability:mathcalC}
    \text{$\mathcal C$ is stable by 
right multiplication by $\mathrm C_{\hat V}(e_0)$. } 
\end{equation}

One has 
\begin{equation}\label{criterion:CDT}
\mathcal C
= \{v\in M_{31}\hat V|
 \rho_1\cdot v=0, \quad \rho_0\cdot v=ve_0\}
 = \{v\in M_{31}\hat V|
 \mathrm{row}_{\mathrm{DT}}\cdot v=0, \quad \rho_0\cdot v=ve_0\}. 
\end{equation}
where the first equality follows from the fact that $\hat{\mathcal V}$ is generated by $e_0,e_1$, 
Def. \ref{def:5:3:paris} and the fact that $\chi$ is an algebra morphism, and the second equality follows from 
$\rho_1=\mathrm{col}_{\mathrm{DT}}\cdot\mathrm{row}_{\mathrm{DT}}$ and the fact that one of the entries of
$\mathrm{col}_{\mathrm{DT}}$ is equal to 1. 

One checks that $\mathrm{row}_{\mathrm{DT}}\cdot C_{\mathrm{DT}}=0$ and 
$\rho_0\cdot C_{\mathrm{DT}}=C_{\mathrm{DT}}e_0$, which by 
\eqref{criterion:CDT} implies $C_{\mathrm{DT}}\in\mathcal C$. 
\eqref{stability:mathcalC} then implies 
\begin{equation}\label{mathcalC:contains:rightmodule}
 C_{\mathrm{DT}}\cdot \mathrm C_{\hat V}(e_0)\subset\mathcal C.   
\end{equation}  
Let us prove the opposite inclusion. 
Let $v=\begin{pmatrix} a\\b\\ c\end{pmatrix}\in \mathcal C$, so by \eqref{criterion:CDT} one has 
\begin{equation}\label{conditions:v}
    \mathrm{row}_{\mathrm{DT}}\cdot v=0, \quad 
\rho_0\cdot v=ve_0.
\end{equation}
The first condition is equivalent 
to the equality $e_1\cdot a=f_1\cdot b$, which by Lem. \ref{LEM1:0301:BIS} is equivalent to
the existence of $d\in\hat V$ such that 
\begin{equation}\label{eqs:a:b:d}
    \text{$a=f_1d$ and $b=e_1d$.} 
\end{equation}
The second condition is equivalent to the conjunction of the following equalities 
\begin{equation}\label{eqs:rho0:paris}
e_0a=ae_0,\quad e_0c=ce_0,\quad e_1(a-c)+f_0b-be_0=0. 
\end{equation}
Plugging the first equality of \eqref{eqs:a:b:d} in the first of these equalities and using the commutation of 
$e_0$ and $f_1$ implies the relation $f_1\cdot (e_0d-de_0)=0$, which by the injectivity of the endomorphism 
$x\mapsto f_1x$ of $\hat V$ implies $e_0d=de_0$ therefore 
\begin{equation}\label{d:in:CVe0}
    d\in \mathrm C_{\hat V}(e_0). 
\end{equation}
Plugging \eqref{eqs:a:b:d} in the third equality of \eqref{eqs:rho0:paris} gives 
$e_1(f_1d-c)+f_0e_1d-e_1de_0=0$, which using the commutativity of $f_0$ and $e_1$ and 
the injectivity of the endomorphism $x\mapsto e_1x$ of $\hat V$ implies
$f_1d-c+f_0d-de_0=0$. Using \eqref{d:in:CVe0}, this equality implies 
$c=-(e_0+f_\infty)d$. The combination of this equality with \eqref{eqs:a:b:d} gives 
$v=C_{\mathrm{DT}}\cdot d$, which by \eqref{d:in:CVe0}
implies $v\in C_{\mathrm{DT}}\cdot \mathrm C_{\hat V}(e_0)$. It follows that 
$\mathcal C\subset C_{\mathrm{DT}}\cdot \mathrm C_{\hat V}(e_0)$, which together  
with \eqref{mathcalC:contains:rightmodule} implies $\mathcal C=C_{\mathrm{DT}}\cdot \mathrm C_{\hat V}(e_0)$. 
This implies the surjectivity of the map $\mathrm C_{\hat V}(e_0)\to \mathcal C$, 
$a\mapsto C_{\mathrm{DT}}\cdot a$; its injectivity follows from that of the endomorphism of 
$\hat V$ given by $v\mapsto e_1v$. 

(b) Set $\mathcal R:=\{\xi\in M_{13}\hat V|
\forall x\in\hat{\mathcal V}, \xi\rho_{\mathrm{DT}}(x)=\chi(x)\xi\}$. It follows from the fact that 
the image of $\chi$ is contained in $\mathbf k[[e_0]]\subset \hat V$ that
\begin{equation}\label{stability:mathcalR}
    \text{$\mathcal R$ is stable by 
left multiplication by $\mathrm C_{\hat V}(e_0)$. } 
\end{equation}

One has 
\begin{equation}\label{criterion:RDT}
\mathcal C
= \{\xi\in M_{13}\hat V|
 \xi\cdot\rho_1=0, \quad \xi\cdot\rho_0=e_0\xi\}
 = \{\xi\in M_{13}\hat V|
 \xi \cdot \mathrm{row}_{\mathrm{DT}}=0, \quad \xi\cdot \rho_0=e_0\xi\}. 
\end{equation}
where the first equality follows from the fact that $\hat{\mathcal V}$ is generated by $e_0,e_1$, 
Def. \ref{def:5:3:paris} and the fact that $\chi$ is an algebra morphism, and the second equality follows from 
$\rho_1=\mathrm{col}_{\mathrm{DT}}\cdot\mathrm{row}_{\mathrm{DT}}$, the fact that one of the entries of
$\mathrm{col}_{\mathrm{DT}}$ is equal to $e_1$ and the injectivity of the endomorphism of $\hat V$ given by 
$x\mapsto e_1 x$. 

One checks that $R_{\mathrm{DT}}\cdot \mathrm{col}_{\mathrm{DT}}=0$ and 
$R_{\mathrm{DT}}\cdot \rho_0=e_0 R_{\mathrm{DT}}$, which by 
\eqref{criterion:RDT} implies $R_{\mathrm{DT}}\in\mathcal R$. 
\eqref{stability:mathcalR} then implies 
\begin{equation}\label{mathcalR:contains:rightmodule}
\mathrm C_{\hat V}(e_0)\cdot R_{\mathrm{DT}}\subset\mathcal R.   
\end{equation}  
Let us prove the opposite inclusion. Let $\xi=\begin{pmatrix}
    \alpha & \beta & \gamma
\end{pmatrix}\in \mathcal R$, so by \eqref{criterion:RDT} one has 
$$
 \xi \cdot \mathrm{row}_{\mathrm{DT}}=0, \quad \xi\cdot \rho_0=e_0\xi. 
$$
The first of these equalities implies 
\begin{equation}\label{alpha=beta}
    \alpha=\beta, 
\end{equation} 
while the second one implies the conjunction of the following equalities
\begin{equation}\label{arr:sbg}
\alpha e_0+\beta e_1=e_0\alpha,\quad \alpha f_0=e_0\alpha,\quad \gamma e_0=e_0\gamma. 
\end{equation}
By Lem. \ref{lem:technical:0801}, the second equality implies 
\begin{equation}\label{alpha=0}
    \alpha=0, 
\end{equation}
 \eqref{alpha=beta} implies 
 \begin{equation}\label{beta=0}
    \beta=0.  
\end{equation} 
The third equation of \eqref{arr:sbg} implies 
\begin{equation}\label{where:is:gamma}
    \gamma\in \mathrm C_{\hat V}(e_0).
\end{equation} 
Then \eqref{alpha=0} and \eqref{beta=0} imply $\xi=\gamma R_{\mathrm{DT}}$, which by 
\eqref{where:is:gamma} implies $\xi\in \mathrm C_{\hat V}(e_0)\cdot R_{\mathrm{DT}}$. 
Hence $\mathcal C\subset \mathrm C_{\hat V}(e_0)\cdot R_{\mathrm{DT}}$, with together with \eqref{mathcalR:contains:rightmodule}
implies $\mathcal C\subset \mathrm C_{\hat V}(e_0)\cdot R_{\mathrm{DT}}$. This implies the 
surjectivity of the map $\mathrm C_{\hat V}(e_0)\to \mathcal R$, 
$a\mapsto a\cdot C_{\mathrm{DT}}$; its injectivity follows from the fact that one of the entries of 
$C_{\mathrm{DT}}$ is 1. 
\end{proof}

\begin{lem}\label{LABEL:3001}
Let $\chi$ be as in Def. \ref{def:chi:3001}. Then for any $g\in\mathcal G$, one has 
$\chi \circ \mathrm{aut}_g^{\mathcal V}=\chi$ (equality of algebra endomorphisms of 
$\hat{\mathcal V}$). 
\end{lem}

\begin{proof}
One has $\chi \circ \mathrm{aut}_g^{\mathcal V}(e_1)=\chi(e_1)$ $(=0)$ 
as $\mathrm{aut}_g^{\mathcal V}(e_1)=e_1$ and $\chi\circ \mathrm{aut}_g^{\mathcal V} 
(e_0)=\chi(g(e_0,e_1)\cdot e_0 \cdot g(e_0,e_1)^{-1})=g(e_0,0)\cdot e_0\cdot g(e_0,0)^{-1}
=e_0=\chi(e_0)$, where the third equation follows from 
the commutativity of the subalgebra of $\hat{\mathcal V}$ generated by $e_0$. 
\end{proof}

\begin{lem}\label{relations:ROW:1301}
Let $g\in\mathrm{Stab}_{\mathcal G}(\mathrm{GL}_3\hat V\bullet\rho_{\mathrm{DT}})$. 

(a) There exists a linear map $\mathrm{Tor}_g\to\mathrm{C}_{\hat V}(e_0)$, $P\mapsto\kappa_g^P$, such that 
$$
\forall P\in \mathrm{Tor}_g,\quad 
g(e_0,e_1)^{-1}(\mathrm{aut}_g^{\mathcal V})^{\otimes2}(R_{\mathrm{DT}}) \cdot P
=\kappa_g^P\cdot R_{\mathrm{DT}}.
$$

(b) The map $P\mapsto\kappa_g^P$ restricts to a map  $\mathrm{Tor}_g^\times\to\mathrm{C}_{\hat V}(e_0)^\times$ and satisfies the identity
$$
\forall P\in \mathrm{Tor}_g,\forall(\phi,a)\in\mathbf k\times\mathrm C_{\hat V}(e_0),\quad 
\kappa_g^{P\cdot P_0(\phi,a)}=\kappa_g^P\cdot (\phi-(e_0+f_\infty)a), 
$$
where $(\phi,a)\mapsto P_0(\phi,a):=\phi I_3+C_{\mathrm{DT}}\cdot a \cdot R_{\mathrm{DT}}$ is the bijection 
$\mathbf k\times \mathrm{C}_{\hat V}(e_0)\to\mathrm{C}_3(\rho_{\mathrm{DT}}(e_0))$ from Lem. \ref{lem:commutant:BIS}(a). 
\end{lem}

\begin{proof}
(a) Let $P\in\mathrm{Tor}_g$. Then by Def. \ref{def:12:1:sbg}(a), one has  
\begin{equation}\label{identité:2304:0304}
\forall a\in\hat{\mathcal V},\quad P\cdot \rho_{\mathrm{DT}} ( \mathrm{aut}_g^{\mathcal V}(a))=
(\mathrm{aut}_g^{\mathcal V})^{\otimes2} (\rho_{\mathrm{DT}}(a))\cdot P,     
\end{equation}
which implies the first equality in 
\begin{align}\label{toto:3001}
&\forall a\in\hat{\mathcal V},\quad 
(\mathrm{aut}_g^{\mathcal V})^{\otimes2}(R_{\mathrm{DT}}) \cdot P\cdot \rho_{\mathrm{DT}} ( \mathrm{aut}_g^{\mathcal V}(a))
=(\mathrm{aut}_g^{\mathcal V})^{\otimes2}(R_{\mathrm{DT}})
\cdot (\mathrm{aut}_g^{\mathcal V})^{\otimes2} (\rho_{\mathrm{DT}}(a))\cdot P
\\ & \nonumber 
=(\mathrm{aut}_g^{\mathcal V})^{\otimes2}(R_{\mathrm{DT}}
\cdot \rho_{\mathrm{DT}}(a))\cdot P
=(\mathrm{aut}_g^{\mathcal V})^{\otimes2}(\chi(a)R_{\mathrm{DT}})\cdot P
\\&\nonumber=g(e_0,e_1)\chi(a)g(e_0,e_1)^{-1}(\mathrm{aut}_g^{\mathcal V})^{\otimes2}(R_{\mathrm{DT}})\cdot P
=g(e_0,e_1)\chi(\mathrm{aut}_g^{\mathcal V}(a))g(e_0,e_1)^{-1}(\mathrm{aut}_g^{\mathcal V})^{\otimes2}(R_{\mathrm{DT}})\cdot P,  
\end{align}
where the third equality follows from Lem. \ref{lem:13:9:1301}(b), the fourth equality follows from the facts that the 
restriction of $(\mathrm{aut}_g^{\mathcal V})^{\otimes2}$ to $\mathbf k[[e_0]]$ coincides with conjugation by $g(e_0,e_1)$, and 
that the image of $\chi$ is contained in $\mathbf k[[e_0]]$, and the fifth equality follows from Lem. \ref{LABEL:3001}. 
Right multiplying by $g(e_0,e_1)^{-1}$ and replacing $a$ by its preimage by $(\mathrm{aut}_g^{\mathcal V})^{\otimes2}$ in 
the resulting identity, one obtains
$$
\forall a\in\hat{\mathcal V},\quad 
g(e_0,e_1)^{-1}(\mathrm{aut}_g^{\mathcal V})^{\otimes2}(R_{\mathrm{DT}}) \cdot P\cdot \rho_{\mathrm{DT}} ( a)
=\chi(a)g(e_0,e_1)^{-1}(\mathrm{aut}_g^{\mathcal V})^{\otimes2}(R_{\mathrm{DT}})\cdot P,  
$$
therefore $g(e_0,e_1)^{-1}(\mathrm{aut}_g^{\mathcal V})^{\otimes2}(R_{\mathrm{DT}}) \cdot P$ belongs to 
$\{\xi\in M_{13}\hat V|\forall x\in\hat{\mathcal V}, \xi\rho_{\mathrm{DT}}(x)=\chi(x)\xi\}$. 
Lem. \ref{lem:13:9:1301}(b) then implies the existence of $\kappa_g^P\in\mathrm{C}_{\hat V}(e_0)$ such that  
\begin{equation}\label{eq:2022:0304}
g(e_0,e_1)^{-1}(\mathrm{aut}_g^{\mathcal V})^{\otimes2}(R_{\mathrm{DT}}) \cdot P
=\kappa_g^P\cdot R_{\mathrm{DT}}.     
\end{equation}
The map $P\mapsto \kappa_g^P$ is $\mathbf k$-linear since the map $P\mapsto 
g(e_0,e_1)^{-1}(\mathrm{aut}_g^{\mathcal V})^{\otimes2}(R_{\mathrm{DT}}) \cdot P$ is linear and since the map from 
Lem. \ref{lem:13:9:1301}(b) is an isomorphism of $\mathbf k$-modules.

(b) Assume that $P\in\mathrm{Tor}_g^\times$. Projecting \eqref{eq:2022:0304} in degree 0 and denoting with an index 0 
the degree 0 parts of elements of graded $\mathbf k$-modules, one obtains $R_{\mathrm{DT}}\cdot P_0=(\kappa_g^P)_0\cdot R_{\mathrm{DT}}$, therefore 
$(\kappa_g^P)_0$ is an eigenvalue of $P_0^t$. 
Since $P_0\in\mathrm{GL}_3(\mathbf k)$, its eigenvalues belong to $\mathbf k^\times$, which implies 
$(\kappa_g^P)_0\in\mathbf k^\times$, therefore $\kappa_g^P\in\hat V^\times$, which then implies 
$\kappa_g^P\in\mathrm{C}_{\hat V}(e_0)^\times$. 

If $(\phi,a)\in\mathbf k\times\mathrm C_{\hat V}(e_0)$, and $P\in\mathrm{Tor}_g$, then $P\cdot P_0(\phi,a)\in \mathrm{Tor}_g$. 
Then 
\begin{align*}
&
\kappa_g^{P\cdot P_0(\phi,a)}\cdot R_{\mathrm{DT}}
=g(e_0,e_1)^{-1}(\mathrm{aut}_g^{\mathcal V})^{\otimes2}(R_{\mathrm{DT}}) \cdot P\cdot P_0(\phi,a)
=\kappa_g^{P}\cdot R_{\mathrm{DT}}\cdot P_0(\phi,a)
\\&
=\kappa_g^{P}\cdot R_{\mathrm{DT}}\cdot(\phi I_3+C_{\mathrm{DT}}\cdot a \cdot R_{\mathrm{DT}})
=\kappa_g^{P}\cdot (\phi-(e_0+f_\infty)a) R_{\mathrm{DT}}
\end{align*}
where the two first equalities follow from \eqref{eq:2022:0304} applied to $P$ and $P\cdot P_0(\phi,a)$. 
Then statement then follows from the fact that one of the entries of $R_{\mathrm{DT}}$ is 1. 
\end{proof}

Denote by $g\mapsto g^{\circledast-1}$ the operation of taking the inverse in the group $(\mathcal G,\circledast)$.  
Then for $g\in\mathcal G$, one has 
$(\mathrm{aut}_g^{\mathcal V})^{-1}(e_0)=\mathrm{Ad}_{g^{\circledast-1}}(e_0)
=g^{\circledast-1}(e_0,e_1)\cdot e_0\cdot g^{\circledast-1}(e_0,e_1)^{-1}$.
Notice that $\mathrm{Ad}_{(g^{\circledast-1} \otimes 1)^{-1}} \circ ((\mathrm{aut}_g^{\mathcal V})^{-1})^{\otimes2}$
is an automorphism of $\mathrm{C}_{\hat V}(e_0)$.

\begin{lem}\label{rels:COL:1301}
Let $g\in\mathrm{Stab}_{\mathcal G}(\mathrm{GL}_3\hat V\bullet\rho_{\mathrm{DT}})$. 

(a) There exists a linear map $\mathrm{Tor}_g\to\mathrm{C}_{\hat V}(e_0)$, $P\mapsto\upsilon_g^P$, such that 
$$
\forall P\in\mathrm{Tor}_g,\quad
((\mathrm{aut}_g^{\mathcal V})^{-1})^{\otimes2}(P\cdot C_{\mathrm{DT}}) g^{\circledast-1}(e_0,e_1)
= C_{\mathrm{DT}}\cdot\upsilon_g^P. 
$$

(b) The map $P\mapsto\upsilon_g^P$ restricts to a map $\mathrm{Tor}_g^\times\to\mathrm{C}_{\hat V}(e_0)^\times$, 
and satisfies the identity
$$
\forall P\in \mathrm{Tor}_g,\forall(\phi,a)\in\mathbf k\times\mathrm C_{\hat V}(e_0),\quad 
\upsilon^{P\cdot P_0(\phi,a)}=\upsilon_g^P\cdot \mathrm{Ad}_{(g^{\circledast-1} \otimes 1)^{-1}} \circ ((\mathrm{aut}_g^{\mathcal V})^{-1})^{\otimes2}(\phi-a(e_0+f_\infty))
$$
where $(\phi,a)\mapsto P_0(\phi,a)$ is as in Lem. \ref{relations:ROW:1301}(b). 
\end{lem}

\begin{proof}
(a) Let $P\in\mathrm{Tor}_g$. 
One has 
$$
\forall a\in\hat{\mathcal V},\quad P\cdot C_{\mathrm{DT}}\cdot \chi(a)=
P\cdot C_{\mathrm{DT}}\cdot \chi(\mathrm{aut}_g^{\mathcal V}(a))=
P\cdot \rho_{\mathrm{DT}} ( \mathrm{aut}_g^{\mathcal V}(a))\cdot C_{\mathrm{DT}}=
(\mathrm{aut}_g^{\mathcal V})^{\otimes2} (\rho_{\mathrm{DT}}(a))\cdot P\cdot C_{\mathrm{DT}},     
$$
where the first equality follows from Lem. \ref{LABEL:3001}, the second equality follows from Lem. \ref{lem:13:9:1301}(a), 
and the third equality follows from \eqref{identité:2304:0304}. Applying the inverse of $(\mathrm{aut}_g^{\mathcal V})^{\otimes2}$
to the resulting equality given the first equality in 
\begin{align*}
&\forall a\in\hat{\mathcal V},\quad
 \rho_{\mathrm{DT}}(a)\cdot ((\mathrm{aut}_g^{\mathcal V})^{-1})^{\otimes2}(P\cdot C_{\mathrm{DT}})
= ((\mathrm{aut}_g^{\mathcal V})^{-1})^{\otimes2}(P\cdot C_{\mathrm{DT}})\cdot ((\mathrm{aut}_g^{\mathcal V})^{-1})^{\otimes2}(\chi(a))
\\&
= ((\mathrm{aut}_g^{\mathcal V})^{-1})^{\otimes2}(P\cdot C_{\mathrm{DT}})\cdot 
g^{\circledast-1}(e_0,e_1)\chi(a) g^{\circledast-1}(e_0,e_1)^{-1} 
\end{align*}
where the second equality follows from the facts that the restriction of $((\mathrm{aut}_g^{\mathcal V})^{-1})^{\otimes2}$
to $\mathbf k[[e_0]]\subset \hat V$ coincides with the conjugation by $g^{\circledast-1}(e_0,e_1)$, and that the image of 
$\chi$ is contained in $\mathbf k[[e_0]]$. One derives the identity 
$$
\forall a\in\hat{\mathcal V},\quad
\rho_{\mathrm{DT}}(a)\cdot ((\mathrm{aut}_g^{\mathcal V})^{-1})^{\otimes2}(P\cdot C_{\mathrm{DT}})
g^{\circledast-1}(e_0,e_1)
= ((\mathrm{aut}_g^{\mathcal V})^{-1})^{\otimes2}(P\cdot C_{\mathrm{DT}})\cdot 
g^{\circledast-1}(e_0,e_1)\chi(a)  
$$
therefore $((\mathrm{aut}_g^{\mathcal V})^{-1})^{\otimes2}(P\cdot C_{\mathrm{DT}}) g^{\circledast-1}(e_0,e_1)$ belongs to 
$\{v\in M_{31}\hat V|\forall x\in\hat{\mathcal V}, \rho_{\mathrm{DT}}(x)v=v\chi(x)\}$. 
Lem. \ref{lem:13:9:1301}(a) then implies the existence of $\upsilon_g^P\in\mathrm{C}_{\hat V}(e_0)$ such that  
\begin{equation}\label{eq:2357:0304}
((\mathrm{aut}_g^{\mathcal V})^{-1})^{\otimes2}(P\cdot C_{\mathrm{DT}}) g^{\circledast-1}(e_0,e_1)
= C_{\mathrm{DT}}\cdot\upsilon_g^P.     
\end{equation}
The map $P\mapsto \upsilon_g^P$ is $\mathbf k$-linear since the map $P\mapsto 
((\mathrm{aut}_g^{\mathcal V})^{-1})^{\otimes2}(P\cdot C_{\mathrm{DT}}) g^{\circledast-1}(e_0,e_1)$ is linear and since the map from 
Lem. \ref{lem:13:9:1301}(a) is an isomorphism of $\mathbf k$-modules.

(b) Assume that $P\in\mathrm{Tor}_g^\times$. The projection in degree 1 of \eqref{eq:2357:0304} implies 
$P_0\cdot C_{\mathrm{DT}}=C_{\mathrm{DT}}\cdot\upsilon_{g,0}^P$, where the indices $0$ denote the degree 
0 parts, which then implies 
$P_0\cdot \begin{pmatrix}
    1\\0\\0
\end{pmatrix}=\begin{pmatrix}
    1\\0\\0
\end{pmatrix}\cdot\upsilon_{g,0}^P$ and $P_0\cdot \begin{pmatrix}
    0\\1\\0
\end{pmatrix}=\begin{pmatrix}
    0\\1\\0
\end{pmatrix}\cdot\upsilon_{g,0}^P$. Since $P_0\in\mathrm{GL}_3(\mathbf k)$, this implies 
$\upsilon_{g,0}^P\in\mathbf{k}^\times$, therefore $\upsilon_g^P\in\hat V^\times$, therefore 
$\upsilon_g^P\in\mathrm{C}_{\hat V}(e_0)^\times$. 

Let $P\in\mathrm{Tor}_g$ and $(\phi,a)\in\mathbf k\times\mathrm C_{\hat V}(e_0)$. Then 
$P\cdot P_0(\phi,a)\in \mathrm{Tor}_g$, therefore 
\begin{align*}
& C_{\mathrm{DT}}\cdot\upsilon_g^{P\cdot P_0(\phi,a)}
=((\mathrm{aut}_g^{\mathcal V})^{-1})^{\otimes2}(P\cdot P_0(\phi,a)\cdot C_{\mathrm{DT}}) g^{\circledast-1}(e_0,e_1)
\\& =((\mathrm{aut}_g^{\mathcal V})^{-1})^{\otimes2}(P\cdot C_{\mathrm{DT}}(\phi-a(e_0+f_\infty))) g^{\circledast-1}(e_0,e_1)
\\& = ((\mathrm{aut}_g^{\mathcal V})^{-1})^{\otimes2}(P\cdot C_{\mathrm{DT}})g^{\circledast-1}(e_0,e_1)
\cdot g^{\circledast-1}(e_0,e_1)^{-1}
((\mathrm{aut}_g^{\mathcal V})^{-1})^{\otimes2}(\phi-a(e_0+f_\infty)) g^{\circledast-1}(e_0,e_1)
\\& = C_{\mathrm{DT}}\cdot\upsilon_g^{P}
\cdot g^{\circledast-1}(e_0,e_1)^{-1}
((\mathrm{aut}_g^{\mathcal V})^{-1})^{\otimes2}(\phi-a(e_0+f_\infty)) g^{\circledast-1}(e_0,e_1)
\\& = C_{\mathrm{DT}}\cdot\upsilon_g^{P}
\cdot \mathrm{Ad}_{(g^{\circledast-1} \otimes 1)^{-1}} \circ ((\mathrm{aut}_g^{\mathcal V})^{-1})^{\otimes2}(\phi-a(e_0+f_\infty))
\end{align*}
which implies the last claim, since one of the components of $C_{\mathrm{DT}}$ is $e_1$ and since $x\mapsto e_1x$ is an injective 
endomorphism of $\hat V$.   
\end{proof}

\begin{lem}\label{lem:12:8}
For any $g\in\mathrm{Stab}_{\mathcal G}(\mathrm{GL}_3\hat V\bullet\rho_{\mathrm{DT}})$, one has 
\begin{equation}\label{id:V2:stab:0703}
(\kappa_g^P)^{-1}g(e_0,e_1)^{-1}\cdot (\mathrm{aut}_g^{\mathcal V})^{\otimes2}(f_0+f_1-e_0)\cdot
(\mathrm{aut}_g^{\mathcal V})^{\otimes2}(\upsilon_g^P\cdot g^{\circledast-1}(e_0,e_1)^{-1})
=f_0+f_1-e_0. 
\end{equation}
\end{lem}

\begin{proof}
Let $g\in\mathcal G$. For any $P\in \mathrm{Tor}_g$, one has 
\begin{equation}\label{equation:P:first}
g(e_0,e_1)^{-1}(\mathrm{aut}_g^{\mathcal V})^{\otimes2}(R_{\mathrm{DT}}) \cdot P
=\kappa_g^P\cdot R_{\mathrm{DT}}.
\end{equation}
by Lem. \ref{relations:ROW:1301}(a) and 
\begin{equation}\label{equation:P:second}
 (\mathrm{aut}_g^{\mathcal V})^{\otimes2}(C_{\mathrm{DT}}\cdot\upsilon_g^P\cdot g^{\circledast-1}(e_0,e_1)^{-1})
=P\cdot C_{\mathrm{DT}}. 
\end{equation}
by Lem. \ref{rels:COL:1301}(a). If now $g\in\mathrm{Stab}_{\mathcal G}(\mathrm{GL}_3\hat V\bullet\rho_{\mathrm{DT}})$, 
then there exists $P\in \mathrm{Tor}_g^\times$, which satisfies both \eqref{equation:P:first} and  
$$
P^{-1}\cdot (\mathrm{aut}_g^{\mathcal V})^{\otimes2}(C_{\mathrm{DT}}\cdot\upsilon_g^P\cdot g^{\circledast-1}(e_0,e_1)^{-1})
=C_{\mathrm{DT}} 
$$
which follows from \eqref{equation:P:second} and from the invertibility of $P$. 
Combined with the latter equation, \eqref{equation:P:first} implies
$$
g(e_0,e_1)^{-1}(\mathrm{aut}_g^{\mathcal V})^{\otimes2}(R_{\mathrm{DT}}
C_{\mathrm{DT}}\cdot\upsilon_g^P\cdot g^{\circledast-1}(e_0,e_1)^{-1})
=\kappa_g^P\cdot R_{\mathrm{DT}}
C_{\mathrm{DT}}. 
$$
The statement then follows from $R_{\mathrm{DT}}
C_{\mathrm{DT}}=f_0+f_1-e_0$ and the from invertibility of $\kappa_g^P$. 
\end{proof}
 
\begin{lem}\label{TOTO}
For any $g\in\mathrm{Stab}_{\mathcal G}(\mathrm{GL}_3\hat V \bullet \rho_{\mathrm{DT}})$, 
there exist $a,b \in \hat{\mathcal V}^\times$  such that $a\cdot (e_0+e_1)=
\mathrm{aut}_g^{\mathcal V}(e_0+e_1)\cdot b$.
\end{lem}

\begin{proof}
This follows by applying the algebra morphism $\epsilon \otimes\mathrm{id}: \hat{\mathcal V}^{\hat\otimes2}\to\hat{\mathcal V}$ 
to \eqref{id:V2:stab:0703}, and by using the fact that both 
$(\kappa_g^P)^{-1}g(e_0,e_1)^{-1}$ and $
(\mathrm{aut}_g^{\mathcal V})^{\otimes2}(\upsilon_g^P\cdot g^{\circledast-1}(e_0,e_1)^{-1})$ belong to $(\hat{\mathcal V}^{\hat\otimes2})^\times$.  
\end{proof}

\begin{thm}\label{thm:13:22:0205}
There holds the following inclusion 
\begin{equation}\label{stab:in:intert}
\mathrm{Stab}_{\mathcal G}(\mathrm{GL}_3\hat V\bullet\rho_{\mathrm{DT}})\subset \mathcal G_{\mathrm{inert}}
\end{equation}
of subgroups of $\mathcal G$; it fits 
in a diagram of group inclusions
    $$
\xymatrix{
\mathcal G_{\mathrm{inert}}^\Theta\ar@{^{(}->}[r]
& \mathcal G_{\mathrm{inert}}\ar@{^{(}->}[r]&
(\mathcal G,\circledast)\\
\mathsf{GRT}_1(\mathbf k)^{\mathrm{op}}\ar@{^{(}->}[r]\ar@{^{(}->}[u]&
\mathrm{Stab}_{\mathcal G}(\mathrm{GL}_3\hat V\bullet\rho_{\mathrm{DT}})\ar@{^{(}->}[r]\ar@{^{(}->}[u]&
\mathsf{DMR}_0(\mathbf k)\ar@{^{(}->}[u]}
$$
\end{thm}

\begin{proof}
The inclusion $\mathsf{GRT}_1(\mathbf k)^{\mathrm{op}}\subset \mathcal G_{\mathrm{inert}}^\Theta$
follows from the fact that any $g\in\mathsf{GRT}_1(\mathbf k)^{\mathrm{op}}$ satisfies the identity 
$\mathrm{Ad}_{g(e_0,e_1)}(e_0)+e_1+\mathrm{Ad}_{g(e_\infty,e_1)}(e_\infty)=0$ (see \cite{Dr}, (5.14) and (5.12)). 

Let us prove the inclusion $\mathrm{Stab}_{\mathcal G}(\mathrm{GL}_3\hat V\bullet\rho_{\mathrm{DT}})\subset 
\mathcal G_{\mathrm{inert}}$. Let $g\in\mathrm{Stab}_{\mathcal G}(\mathrm{GL}_3\hat V\bullet\rho_{\mathrm{DT}})$. 
Combining Lem. \ref{TOTO} and Prop. \ref{prop:Z} with $x:=e_0+e_1$, $y:=e_1$, $z:=\mathrm{aut}_g^{\mathcal V}(e_0+e_1)-(e_0+e_1)$,
one obtains the existence of $h\in \mathrm{exp}(\mathfrak{lie}_{\{0,1\}}^\wedge)$ such that 
$\mathrm{aut}_g^{\mathcal V}(e_0+e_1)=h\cdot (e_0+e_1)\cdot h^{-1}$. Since 
$\mathrm{aut}_g^{\mathcal V}(e_0+e_1)=g\cdot e_0\cdot g^{-1}+e_1$ and since 
$g\in\mathcal G$, this implies that the degree 2 part of $h\cdot (e_0+e_1)\cdot h^{-1}$ vanishes. 
Since this degree 2 part is equal to $[h_1,e_0+e_1]$, where $h_1$ is the degree 1 part of $h$, it follows
that $h_1$ is proportional to $e_0+e_1$; let $\upsilon\in\mathbf k$ be such that 
$h_1=\upsilon\cdot (e_0+e_1)$, then $h\cdot\mathrm{exp}(-\upsilon\cdot(e_0+e_1))$ belongs to $\mathcal G$
and is such that $\mathrm{Ad}_g(e_0)+e_1+\mathrm{Ad}_{h\cdot\mathrm{exp}(-\upsilon\cdot(e_0+e_1))}(e_\infty)=0$, therefore 
$g\in\mathcal G_{\mathrm{inert}}$. 
\end{proof}

\section{Stability of $\mathrm{Stab}_{\mathcal G}(\mathrm{GL}_3\hat V\bullet\rho_{\mathrm{DT}})$ 
under the involution $\Theta$ of $\mathcal G_{\mathrm{inert}}$}

Recall that the group $\mathrm{Stab}_{\mathcal G}(\mathrm{GL}_3\hat V\bullet\rho_{\mathrm{DT}})$ corresponds to 
an action of $\mathcal G$ on the coset space $\mathrm{GL}_3\hat V\backslash \mathrm{Hom}_{\mathcal C\operatorname{-alg}}
(\hat{\mathcal V},M_3\hat V)$. The purpose of this section is to prove, based on the inclusion of 
this group in $\mathcal G_{\mathrm{inert}}$ 
(see Thm. \ref{thm:13:22:0205}), its stability under the involution $\Theta$ of $\mathcal G_{\mathrm{inert}}$ 
(see Lem. \ref{lem02:2506}(d)); this is obtained in Thm. \ref{thm:13:36:17apr} (§\ref{sect:13:9}) as the result 
of the following steps: 

(a) the inclusion of Thm. \ref{thm:13:22:0205} implies  
$\mathrm{Stab}_{\mathcal G}(\mathrm{GL}_3\hat V\bullet\rho_{\mathrm{DT}})=\mathrm{Stab}_{\mathcal G_{\mathrm{inert}}}(\mathrm{GL}_3\hat V\bullet\rho_{\mathrm{DT}})$; 

(b) a coset space $\mathrm{GL}_3\hat V\backslash 
\mathrm{Hom}_{\mathcal C\operatorname{-alg}}(\mathcal V[z]^\wedge,M_3\hat V)$
can be constructed and equipped with an action of $\mathcal G_{\mathrm{inert}}$ and a
$\mathcal G_{\mathrm{inert}}$-equivariant map to $\mathrm{GL}_3\hat V\backslash 
\mathrm{Hom}_{\mathcal C\operatorname{-alg}}(\hat{\mathcal V},M_3\hat V)$;

(c) a lift $\mathrm{GL}_3\hat V\bullet\tilde\rho_{\mathrm{DT}}$ of the element 
$\mathrm{GL}_3\hat V\bullet\rho_{\mathrm{DT}}$ of the latter space can be constructed, leading together with (b) to a 
morphism of $\mathcal G_{\mathrm{inert}}$-pointed spaces
\begin{equation}\label{diag:12oct}
(\mathrm{GL}_3\hat V\backslash 
\mathrm{Hom}_{\mathcal C\operatorname{-alg}}(\hat{\mathcal V},M_3\hat V),\mathrm{GL}_3\hat V\bullet\rho_{\mathrm{DT}})\to(\mathrm{GL}_3\hat V\backslash 
\mathrm{Hom}_{\mathcal C\operatorname{-alg}}({\mathcal V}[z]^\wedge,M_3\hat V),\mathrm{GL}_3\hat V\bullet\tilde\rho_{\mathrm{DT}}); 
\end{equation}

(d) the inclusion $\mathrm{Stab}_{\mathcal G_{\mathrm{inert}}}(\mathrm{GL}_3\hat V\bullet\tilde\rho_{\mathrm{DT}})\subset 
\mathrm{Stab}_{\mathcal G_{\mathrm{inert}}}(\mathrm{GL}_3\hat V\bullet\rho_{\mathrm{DT}})$ of subgroups of 
$\mathcal G_{\mathrm{inert}}$ induced by \eqref{diag:12oct} is an equality (\S\ref{sect:inj:H}, Cor. \ref{cor:13:29}); 

(e) the involution $\Theta$ of $\mathcal G_{\mathrm{inert}}$ gives rise to a semidirect product 
$\mathcal G_{\mathrm{inert}}\rtimes(\mathbb Z/2\mathbb Z)$, which acts on 
$\mathrm{GL}_3\hat V\backslash \mathrm{Hom}_{\mathcal C\operatorname{-alg}}(\mathcal V[z]^\wedge,M_3\hat V)$
by extending the action of $\mathcal G_{\mathrm{inert}}$ from (b); 

(f) the element $\overline 1\in\mathbb Z/2\mathbb Z\subset 
\mathcal G_{\mathrm{inert}}\rtimes(\mathbb Z/2\mathbb Z)$ belongs to 
$\mathrm{Stab}_{\mathcal G_{\mathrm{inert}}\rtimes(\mathbb Z/2\mathbb Z)}(\mathrm{GL}_3
\hat V\bullet\tilde\rho_{\mathrm{DT}})$;

(g) (f) then implies the $\Theta$-stability of 
$\mathrm{Stab}_{\mathcal G_{\mathrm{inert}}}(\mathrm{GL}_3\hat V\bullet\tilde\rho_{\mathrm{DT}})$, which by (d)
implies that $\mathrm{Stab}_{\mathcal G_{\mathrm{inert}}}(\mathrm{GL}_3\hat V\bullet\rho_{\mathrm{DT}})$
is $\Theta$-stable, which by (a) implies that $\mathrm{Stab}_{\mathcal G}(\mathrm{GL}_3\hat V\bullet\rho_{\mathrm{DT}})$
is $\Theta$-stable. 

Steps (b) and (c) are carried out in §§\ref{sect:13}-\ref{sect:13:4}; more precisely, 
in \S\ref{sect:13}, we construct the action of $\mathcal G_{\mathrm{inert}}$ on 
$\mathcal V[z]^\wedge$ and $M_3\hat V$, which induces its action on 
$\mathrm{Hom}_{\mathcal C\operatorname{-alg}}(\mathcal V[z]^\wedge,M_3\hat V)$, compatible with its
action on $\mathrm{Hom}_{\mathcal C\operatorname{-alg}}(\hat{\mathcal V},M_3\hat V)$; in \S\ref{sect:13:2}, we construct 
the element $\tilde\rho_{\mathrm{DT}}$ of this set lifting $\rho_{\mathrm{DT}}$; in \S\ref{sect:13:3}, we 
construct the action of $\mathrm{GL}_3\hat V$
on it, and in \S\ref{sect:13:4}, the overall action of $\mathcal G_{\mathrm{inert}}$. 

Step (d) cannot be obtained from \eqref{diag:12oct} since this map is 
not locally injective; to establish it, one relates the source and target of \eqref{diag:12oct}
by a zig-zag of morphisms of $\mathcal G_{\mathrm{inert}}$-pointed spaces
\begin{align*}
    & (\mathrm{GL}_3\hat V\backslash 
\mathrm{Hom}_{\mathcal C\operatorname{-alg}}(\hat{\mathcal V},M_3\hat V),\mathrm{GL}_3\hat V\bullet\rho_{\mathrm{DT}}) 
\\ & \leftarrow((\mathrm{GL}_3\hat V\times \mathbf k[[e_0,f_\infty]]^\times)\backslash\mathcal S,
(\mathrm{GL}_3\hat V\times \mathbf k[[e_0,f_\infty]]^\times)\bullet(\rho_{\mathrm{DT}},R_{\mathrm{DT}},C_{\mathrm{DT}}))
    \\& \stackrel{(F)}{\to}(\mathrm{GL}_3\hat V\backslash 
\mathrm{Hom}_{\mathcal C\operatorname{-alg}}({\mathcal V}[z]^\wedge,M_3\hat V),
\mathrm{GL}_3\hat V\bullet\tilde\rho_{\mathrm{DT}}), 
\end{align*}
(see \eqref{diag:13:4:11}) and shows $(F)$ to be locally injective. The construction of this diagram is done in \S\S  
\ref{sect:13:2}-\ref{sect:13:4}; the proof of local injectivity if $(F)$ is done in \S\ref{sect:inj:H}
based on algebraic results proved in §\ref{sect:13:5}.  

Steps (e),(f),(g) are then respectively carried out in \S\S\ref{sect:amrzpswga},\ref{sect:13:8},\ref{sect:13:9}. 

\subsection{Actions of the groups $\mathcal G_{\mathrm{inert}}$ and 
$\mathcal G_{\mathrm{inert}}\rtimes(\mathbb Z/2\mathbb Z)$}\label{sect:13}

\begin{defn}
(a) $\mathcal V[z]$ is the free polynomial algebra over $\mathcal V$ in the indeterminate $z$; it is graded, 
the generators $e_0,e_1$ and $z$ being of degree 1. 

(b) $\mathcal V[z]^\wedge$ is the degree completion of the graded $\mathbf k$-algebra $\mathcal V[z]$. 
\end{defn}
There is a natural injection
\begin{equation}\label{inclusion:V}
i_{\mathcal V,\mathcal V[z]} : \hat{\mathcal V}\to \mathcal V[z]^\wedge. 
\end{equation}

\begin{lem}\label{lem:149:0105}
(a) For any $g\in\mathcal G$, there is a unique topological $\mathbf k$-algebras automorphism $\mathrm{aut}_g^{\mathcal V[z]}$ 
of $\mathcal V[z]^\wedge$, such that \eqref{inclusion:V} intertwines $\mathrm{aut}_g^{\mathcal V[z]}$ and 
$\mathrm{aut}_g^{\mathcal V}$ and that $\mathrm{aut}_g^{\mathcal V[z]}(z)=z$. 

(b) The map $g\mapsto\mathrm{aut}_g^{\mathcal V[z]}$ defines an group morphism $\mathcal G\to
\mathrm{Aut}_{\mathcal C\operatorname{-alg}}(\mathcal V[z]^\wedge)$. 
\end{lem}

\begin{proof}
 The topological algebra $\mathcal V[z]^\wedge$ is isomorphic 
 to the completed tensor product $\hat{\mathcal V}\hat\otimes 
 \mathbf k[[z]]$. For any $g \in \mathcal G$, 
 $\mathrm{aut}_g^{\mathcal V[z]}$ is equal to 
 $\mathrm{aut}_g^{\mathcal V} \hat\otimes 
 id_{\mathbf k[[z]]}$. It follows that 
 $g\mapsto\mathrm{aut}_g^{\mathcal V[z]}$ is the map 
 corresponding to the tensor product of the action of 
 $\mathcal G$ on $\hat{\mathcal V}$ by 
 $g\mapsto \mathrm{aut}_g^{\mathcal V}$ with the trivial 
 action of $\mathcal G$ on $\mathbf k[[z]]$, which implies 
 both (a) and (b).    
\end{proof}

\begin{defn}\label{s:O:infty:V[z]}
$s_{(0,\infty)}^{\mathcal V[z]}$ the automorphism of $\mathcal V[z]^\wedge$ defined by 
$z\mapsto z$, $e_0\mapsto e_\infty-z$, $e_\infty\mapsto e_0+z$. 
\end{defn}
One checks that this automorphism is involutive, and it such that $e_1\mapsto e_1$. 

\begin{lem}\label{lem:action:1342:0105}
(a) There is a unique group morphism 
$$
\mathcal G_{\mathrm{inert}}\rtimes(\mathbb Z/2\mathbb Z) \to\mathrm{Aut}_{\mathcal C\operatorname{-alg}}(\hat{\mathcal V}), 
\quad g\mapsto \mathrm{aut}_g^{\mathcal V}, 
$$
where the source is as in Lem. \ref{lem04:2506}(c), extending 
the restriction of the morphism from Lem. \ref{lem03:2506} to $\mathcal G_{\mathrm{inert}}$ and 
such that $(\mathbb Z/2\mathbb Z)\ni\overline 1\mapsto s_{(0,\infty)}$ (see \ref{lem02:2506}(d)). 
One has 
\begin{equation}\label{useful:duality:1553}
    \forall g\in \mathcal G_{\mathrm{inert}},\quad \mathrm{aut}_g^{\mathcal V}(e_\infty)=h_g\cdot e_\infty\cdot h_g^{-1}, 
\end{equation}
the map $g\mapsto h_g$ being as in Lem. \ref{lem02:2506}(a). 

(b) There is a unique group morphism 
$$
\mathcal G_{\mathrm{inert}}\rtimes(\mathbb Z/2\mathbb Z)\to \mathrm{Aut}_{\mathcal C\operatorname{-alg}}(\hat V), \quad
g\mapsto\mathrm{aut}_g^{V},  
$$
extending the restriction of the morphism $\mathcal G\ni g\mapsto \mathrm{aut}_g^V=(\mathrm{aut}_g^{\mathcal V})^{\otimes2}$ 
from Def. \ref{def:5:14:0105} to $\mathcal G_{\mathrm{inert}}$ and such that 
$\mathbb Z/2\mathbb Z\ni\overline 1\mapsto \mathrm{aut}_{\overline 1}^V=\mathrm{sw}\circ 
s_{(0,\infty)}^{\otimes2}$, where $\mathrm{sw}$ is the automorphism of $\hat V=\hat{\mathcal V}^{\hat\otimes2}$ 
given by the exchange of factors. 

(c) There is a unique group morphism 
$$
\mathcal G_{\mathrm{inert}}\rtimes(\mathbb Z/2\mathbb Z)\to 
\mathrm{Aut}_{\mathcal C\operatorname{-alg}}(\mathcal V[z]^\wedge), \quad
g\mapsto\mathrm{aut}_g^{\mathcal V[z]} 
$$
extending the morphism from Lem. \ref{lem:149:0105}(b) and such that $\mathrm{aut}_{\overline 1}^{\mathcal V[z]}:=
s_{(0,\infty)}^{\mathcal V[z]}$.
\end{lem}

\begin{proof}
(a) The composition of the group morphism from Lem. \ref{lem04:2506}(d) with the inclusion of its target in its ambient group 
is a group morphism $\mathcal G_{\mathrm{inert}}\rtimes(\mathbb Z/2\mathbb Z) \to
\mathrm{Aut}_{\mathcal C\operatorname{-alg}}(\hat{\mathcal V})$, which 
has the announced properties. One also checks that a morphism with the said properties is necessarily unique. 
For $g\in \mathcal G_{\mathrm{inert}}$, one has $\mathrm{aut}_g^{\mathcal V}(e_\infty)=-\mathrm{aut}_g^{\mathcal V}(e_1)
-\mathrm{aut}_g^{\mathcal V}(e_\infty)=-e_1-ge_0g^{-1}=h_g\cdot e_\infty\cdot h_g^{-1}$, where the last equality follows from 
Lem. \ref{lem02:2506}(a). 

(b) We will use the following general fact: 
\begin{align}\label{general:fact:morphisms}
&    \text{If $G,K$ are groups, $\alpha,\beta : G\to
H$ are group morphisms, such that for any $g,g'\in G$, }
\\&\nonumber
\text{the elements 
$\alpha(g)$ and $\beta(g')$ commute, then $G\to H$,  $g\mapsto 
\alpha(g)\beta(g)$ is an group morphism.}
\end{align}
 
Let $\chi : \mathcal G_{\mathrm{inert}}\rtimes(\mathbb Z/2\mathbb Z) \to\mathbb Z/2\mathbb Z$ be the canonical morphism. 
Since $\chi$ is a group morphism and $\mathrm{sw}^2=id$, the map $ \mathcal G_{\mathrm{inert}}\rtimes(\mathbb Z/2\mathbb Z)
\to\mathrm{Aut}_{\mathcal C\operatorname{-alg}}(\hat V)$, $g\mapsto \mathrm{sw}^{\chi(g)}$ 
is a group morphism. It follows from (a) that the same is true of the map $g\mapsto (\mathrm{aut}_g^{\mathcal V})^{\otimes 2}$. 
Moreover, for any $g \in \mathcal G$, $(\mathrm{aut}_g^{\mathcal V})^{\otimes 2}$ commutes with 
$\mathrm{sw}$. It then follows from \eqref{general:fact:morphisms} that the map $\mathcal G_{\mathrm{inert}}\rtimes(\mathbb Z/2\mathbb Z)
\to\mathrm{Aut}_{\mathcal C\operatorname{-alg}}(\hat V)$ defined by $g\mapsto 
\mathrm{sw}^{\chi(g)}\circ (\mathrm{aut}_g^{\mathcal V})^{\otimes2}$ is a group morphism. The statement then follows from the facts that 
a group morphism with source $\mathcal G_{\mathrm{inert}}\rtimes(\mathbb Z/2\mathbb Z)$ is uniquely determined by its restrictions to 
$\mathcal G_{\mathrm{inert}}$ and $\mathbb Z/2\mathbb Z$, and from $\mathrm{aut}_{\bar 1}^{\mathcal V}=s_{(0,\infty)}$ (see (a)).  

(c)  Let us prove 
\begin{equation}\label{toprove:0406:b}
s_{(0,\infty)}^{\mathcal V[z]} \circ \mathrm{aut}_g^{\mathcal V[z]}
\circ s_{(0,\infty)}^{\mathcal V[z]}
=\mathrm{aut}_{s_{(0,\infty)}(h_g)}^{\mathcal V[z]}
\end{equation}
for any $g\in\mathcal G_{\mathrm{inert}}$. 

Denote by $s_{(0,\infty)}$ the extension of the automorphism $s_{(0,\infty)}$ of $\hat{\mathcal V}$ 
to an algebra automorphism of $\mathcal V[z]^\wedge$ by $z\mapsto z$. Let also $\tau$ be the algebra 
automorphism of $\mathcal V[z]^\wedge$ given by $e_0\mapsto e_0-z$, $e_\infty\mapsto e_\infty+z$, 
$e_1\mapsto e_1$, $z\mapsto z$. Then 
\begin{equation}\label{s:s:tau}
    s_{(0,\infty)}^{\mathcal V[z]}=s_{(0,\infty)}\circ\tau.  
\end{equation}
Since $z$ is central, the degree $\geq 2$ part of the complete Lie subalgebra of $\mathcal V[z]^\wedge$ generated 
by $e_0,e_1$ is pointwise fixed $\tau$, therefore so is the image of this subset by the exponential map, which 
coincides with that of $\mathcal G\hookrightarrow\hat{\mathcal V}\hookrightarrow\mathcal V[z]^\wedge$. Together with 
\eqref{s:s:tau}, this implies 
\begin{equation}\label{eq:of:two:autos:0406}
\forall g\in\mathcal G,\quad s_{(0,\infty)}^{\mathcal V[z]}(g)=s_{(0,\infty)}(g). 
\end{equation}
Let us now prove \eqref{toprove:0406:b}. Let us compare the images of $e_0$ by both sides: 
\begin{align*}
&
s_{(0,\infty)}^{\mathcal V[z]} \circ \mathrm{aut}_g^{\mathcal V[z]}
\circ s_{(0,\infty)}^{\mathcal V[z]}(e_0)
=s_{(0,\infty)}^{\mathcal V[z]} \circ \mathrm{aut}_g^{\mathcal V[z]}(e_\infty-z)
=s_{(0,\infty)}^{\mathcal V[z]}  (h_g\cdot e_\infty\cdot h_g^{-1}-z)
\\ & \nonumber =s_{(0,\infty)}^{\mathcal V[z]}  (h_g\cdot (e_\infty-z)\cdot h_g^{-1})
=s_{(0,\infty)}^{\mathcal V[z]} (h_g) \cdot e_0\cdot 
s_{(0,\infty)}^{\mathcal V[z]} (h_g)^{-1} 
=s_{(0,\infty)}(h_g) \cdot e_0\cdot 
s_{(0,\infty)}(h_g)^{-1} 
\\&\nonumber =\mathrm{aut}_{s_{(0,\infty)}(h_g)}^{\mathcal V}(e_0)
=\mathrm{aut}_{s_{(0,\infty)}(h_g)}^{\mathcal V[z]}(e_0)\end{align*}
where the first and fourth (resp. second) equality follows from the definition of $s_{(0,\infty)}^{\mathcal V[z]}$
(resp. \eqref{useful:duality:1553}), the 
third equality follows from the centrality of $z$, the fifth equality follows from 
$h_g\in\hat{\mathcal V}^\times$, 
the sixth (resp. last) equality follows from 
\eqref{def:aut:g:V:(1)} (resp., the definition of $g\mapsto \mathrm{aut}_g^{\mathcal V[z]}$).  
Therefore the images of $e_0$ by both sides of \eqref{toprove:0406:b} are equal. 
These two sides are such that $e_1\mapsto e_1$ and $z\mapsto z$, and are both automorphisms of 
$\mathcal V[z]^\wedge$, which implies \eqref{toprove:0406:b}.

For $g\in \mathcal G_{\mathrm{inert}}$, 
one has $\overline 1\circledast g\circledast \overline 1=s_{(0,\infty)}(h_g)$ (equality in 
$\mathcal G_{\mathrm{inert}}\rtimes(\mathbb Z/2\mathbb Z)$), therefore \eqref{toprove:0406:b} implies 
\begin{equation}\label{toprove:0406:a}
s_{(0,\infty)}^{\mathcal V[z]}\circ \mathrm{aut}_g^{\mathcal V[z]}\circ 
s_{(0,\infty)}^{\mathcal V[z]}
=\mathrm{aut}_{\overline 1\circledast g\circledast \overline 1}^{\mathcal V[z]}
\end{equation}
for any $g\in\mathcal G_{\mathrm{inert}}$.

Since $s_{(0,\infty)}^{\mathcal V[z]}$ is involutive, 
it defines a group morphism 
$$
\mathbb Z/2\mathbb Z\to \mathrm{Aut}_{\mathcal C\operatorname{-alg}}(\mathcal V[z]^\wedge),\quad 
\overline 1\mapsto s_{(0,\infty)}^{\mathcal V[z]}. 
$$
Restriction of the morphism from Lem. \ref{lem:149:0105}(b) to the subgroup $\mathcal G_{\mathrm{inert}}\subset \mathcal G$
also defines a group morphism 
$$
\mathcal G_{\mathrm{inert}}\to \mathrm{Aut}_{\mathcal C\operatorname{-alg}}(\mathcal V[z]^\wedge), \quad 
g\mapsto \mathrm{aut}_g^{\mathcal V[z]}. 
$$
The claim follows from the conjunction of these facts with \eqref{toprove:0406:a}.  
\end{proof}

\subsection{A diagram of pointed sets}\label{sect:13:2}

\begin{defn}
    Recall that $\mathrm{Hom}_{\mathcal C\operatorname{-alg}}(\mathcal V[z]^\wedge,M_3\hat V)$ is the the set of morphisms 
    of filtered $\mathbf k$-algebras, both sides being equipped with the decreasing filtrations
associated with their complete graded structures.
\end{defn}

\begin{lem}\label{lem:elt}
There is a unique morphism of filtered $\mathbf k$-algebras $\tilde\rho_{\mathrm{DT}} : \mathcal V[z]^\wedge\to M_3\hat V$, 
whose restriction to $\hat{\mathcal V}$ is $\rho_{\mathrm{DT}}$ (see Def. \ref{def:5:3:paris}) and such that 
$\tilde\rho_{\mathrm{DT}}(z)=C_{\mathrm{DT}}\cdot R_{\mathrm{DT}}$ (see Def. \ref{defn:C:DT:R:DT}). Then 
 $\tilde\rho_{\mathrm{DT}}\in\mathrm{Hom}_{\mathcal C\operatorname{-alg}}(\mathcal V[z]^\wedge,M_3\hat V)$. 
\end{lem}

\begin{proof}
Since $\mathcal V[z]$ is a polynomial extension of $\mathcal V$ and since $-C_{\mathrm{DT}}\cdot R_{\mathrm{DT}}$ commutes with the image of 
$\rho_{\mathrm{DT}} : \mathcal V\to M_3V$ (see \eqref{eq:commutant:BIS}), there is a unique algebra morphism $\mathcal V[z]\to M_3V$ extending 
$\rho_{\mathrm{DT}}$ and such that $z\mapsto -C_{\mathrm{DT}}\cdot R_{\mathrm{DT}}$. Since this morphism is graded, it extends to a morphism 
between graded completions, which is determined by the said properties. 
\end{proof}

\begin{defn}
    $\mathcal S$ is the set of tuples $(\rho,R,C)$ such that $\rho\in \mathrm{Hom}_{\mathcal C\operatorname{-alg}}(\hat{\mathcal V},M_3\hat V)$
    (see Def. \ref{def:hom:m3:jan:2025}) and $(R,C)\in M_{13}\hat V\times M_{31}F^1\hat V$ are such that 
    $$
    \mathrm C_3(\rho(\hat{\mathcal V}))=\mathbf kI_3+C\cdot \mathrm C_{\hat V}(e_0)\cdot R, \quad RC=-(e_0+f_\infty). 
    $$
\end{defn}

\begin{lem}\label{lem:13:8}
The triple $(\rho_{\mathrm{DT}},R_{\mathrm{DT}},C_{\mathrm{DT}})$ belongs to $\mathcal S$. 
\end{lem}

\begin{proof}
This follows from \eqref{eq:commutant:BIS} and from Def. \ref{defn:C:DT:R:DT}, which implies in particular 
$C_{\mathrm{DT}}R_{\mathrm{DT}}=-(e_0+f_\infty)$. 
\end{proof}

\begin{lem}\label{morphisms:pointed:sets:7apr}
(a) The assignment $\tilde\rho\mapsto \tilde\rho\circ i_{\mathcal V,\mathcal V[z]}=\tilde\rho|_{\hat{\mathcal V}}$ 
(where $i_{\mathcal V,\mathcal V[z]}$ is as in \eqref{inclusion:V})
induces a morphism of pointed sets (i.e. in the category $\mathbf{PS}$) 
\begin{equation}\label{pre:F}
    (\mathrm{Hom}_{\mathcal C\operatorname{-alg}}(\mathcal V[z]^\wedge,M_3\hat V),\tilde\rho_{\mathrm{DT}})
    \to(\mathrm{Hom}_{\mathcal C\operatorname{-alg}}(\hat{\mathcal V},M_3\hat V),\rho_{\mathrm{DT}}). 
\end{equation}

    (b) For any $(\rho,R,C)\in \mathcal S$ there exists a unique 
    $\tilde\rho\in \mathrm{Hom}_{\mathcal C\operatorname{-alg}}(\mathcal V[z]^\wedge,M_3\hat V)$ whose restriction to 
    $\hat{\mathcal V}$ is $\rho$ and such that $\tilde\rho(z)=CR$. 
    The assignment $(\rho,R,C)\mapsto \tilde \rho$ defines a morphism of pointed sets 
\begin{equation}\label{pre:G}
    (\mathcal S,(\rho_{\mathrm{DT}},R_{\mathrm{DT}},C_{\mathrm{DT}}))\to 
    (\mathrm{Hom}_{\mathcal C\operatorname{-alg}}(\mathcal V[z]^\wedge,M_3\hat V),\tilde\rho_{\mathrm{DT}}). 
\end{equation}

    (c) The composition of the morphisms (b) and (a) is the morphism of pointed sets 
\begin{equation}\label{pre:H}
    (\mathcal S,(\rho_{\mathrm{DT}},R_{\mathrm{DT}},C_{\mathrm{DT}}))\to
    (\mathrm{Hom}_{\mathcal C\operatorname{-alg}}(\hat{\mathcal V},M_3\hat V),\rho_{\mathrm{DT}})
\end{equation}
    given by $(\rho,R,C)\mapsto \rho$. 
\end{lem}

\begin{proof}
 (a) follows from the fact that $i_{\mathcal V,\mathcal V[z]}$ is a morphism in $\mathcal C\operatorname{-alg}$
 the definition of $\tilde\rho_{\mathrm{DT}}$. 

 (b) If $(\rho,R,C) \in \mathcal S$, then $CR \in M_3F^1\hat V$ and $CR$ belongs to the commutant 
 $\mathrm C_3(\rho(\hat{\mathcal V}))$, which implies first statement. The second 
 statement follows from the definition of $\tilde\rho_{\mathrm{DT}}$.

 (c) follows from the fact that if $(\rho,R,C) \in \mathcal S$, then the morphism 
 $\tilde\rho$ associated to it in (b) is such that its restriction to $\hat{\mathcal V}$ 
 is $\rho$. 
\end{proof}

The morphisms of pointed sets from Lem. \ref{morphisms:pointed:sets:7apr} fit in a diagram
\begin{equation}\label{diag:ptt:sets:7avril}
\xymatrix{&(\mathcal S,(\rho_{\mathrm{DT}},R_{\mathrm{DT}},C_{\mathrm{DT}}))\ar[ld]\ar[rd] & \\ 
(\mathrm{Hom}_{\mathcal C\operatorname{-alg}}(\mathcal V[z]^\wedge,M_3\hat V),\tilde\rho_{\mathrm{DT}})
\ar[rr]& & 
(\mathrm{Hom}_{\mathcal C\operatorname{-alg}}(\hat{\mathcal V},M_3\hat V),\rho_{\mathrm{DT}})}
\end{equation}

\subsection{A diagram of pointed sets with group actions}\label{sect:13:3}

Since the elements $e_0,f_\infty\in V$ commute, the subalgebra of $V$ generated by them, denoted $\mathbf k[e_0,f_\infty]$, is 
commutative; its closure in $\hat V$ will be denoted $\mathbf k[[e_0,f_\infty]]$. 

\begin{lem}\label{lem:actions:0704}
    (a) The map $(P,\tilde\rho)\mapsto P\bullet\tilde\rho:=\mathrm{Ad}_P\circ \tilde\rho$ defines an action of 
    the group $\mathrm{GL}_3\hat V$ on the set 
    $\mathrm{Hom}_{\mathcal C\operatorname{-alg}}(\mathcal V[z]^\wedge,M_3\hat V)$. 

    (b) The map 
    $$
    ((P,\varphi),(\rho,R,C))\mapsto (P,\varphi)\bullet(\rho,R,C):=(\mathrm{Ad}_P\circ\rho,\varphi RP^{-1},PC\varphi^{-1})
    $$
    defines an action of the group $\mathrm{GL}_3\hat V\times\mathbf k[[e_0,f_\infty]]^\times$ on the set 
    $\mathcal S$. 
\end{lem}

\begin{proof}
 The proof of (a) is similar to that of Lem. \ref{lem:actions:0704}(a). Let $(\rho,R,C)\in \mathcal S$ and 
 $(P,\varphi)\in \mathrm{GL}_3\hat V\times\mathbf k[[e_0,f_\infty]]^\times$. Then $\mathrm{Ad}_P\circ\rho\in 
 \mathrm{Hom}_{\mathcal C\operatorname{-alg}}(\hat{\mathcal V},M_3\hat V)$ as $\mathrm{Ad}_P$ is an 
 algebra automorphism of $M_3\hat V$, and by one has clearly $\varphi RP^{-1}\in M_{13}\hat V$, $PC\varphi^{-1}\in M_{31}F^1\hat V$. 
 Moreover, 
\begin{align*}
    &\mathrm C_3(\mathrm{Ad}_P\circ\rho(\hat{\mathcal V}))=\mathrm{Ad}_P(\mathrm C_3(\rho(\hat{\mathcal V}))) 
=P(\mathbf kI_3+C\mathrm{C}_{\hat V}(e_0)R)P^{-1}=\mathbf kI_3+PC\cdot \mathrm{C}_{\hat V}(e_0)\cdot RP^{-1}
    \\&=\mathbf kI_3+PC\varphi^{-1}\cdot \mathrm{C}_{\hat V}(e_0)\cdot \varphi RP^{-1}
\end{align*}
 and $\varphi RP^{-1}\cdot PC\varphi^{-1}=\varphi RC\varphi^{-1}=-\varphi (e_0+f_\infty)\varphi^{-1}=-(e_0+f_\infty)$. 
 It follows that $(P,\varphi)\bullet(\rho,R,C)\in \mathcal S$. If now $(P',\varphi')\in 
 \mathrm{GL}_3\hat V\times\mathbf k[[e_0,f_\infty]]^\times$, then 
\begin{align*}
& (P',\varphi')\bullet((P,\varphi)\bullet(\rho,R,C))= (P',\varphi')\bullet(\mathrm{Ad}_P\circ\rho,\varphi RP^{-1},PC\varphi^{-1})
\\ &=(\mathrm{Ad}_{P'}\circ\mathrm{Ad}_P\circ\rho,\varphi'\varphi RP^{-1}(P')^{-1},P'PC\varphi^{-1}(\varphi')^{-1})
=(\mathrm{Ad}_{P'P}\circ\rho,\varphi'\varphi R(P'P)^{-1},P'PC(\varphi'\varphi)^{-1})
\\& = (P'P,\varphi'\varphi)\bullet(\rho,R,C)=((P',\varphi')\cdot(P,\varphi))\bullet(\rho,R,C)
\end{align*}

 All this implies (b). 
\end{proof}

\begin{lem}\label{lem:13:10:7avril}
    (a) The group morphism $\mathrm{GL}_3\hat V\times\mathbf k[[e_0,f_\infty]]^\times\to\mathrm{GL}_3\hat V$ defined by projection on the first factor is compatible with the actions of its source on $\mathcal S$ (see Lem. \ref{lem:actions:0704}(b)) and of its target on 
    $\mathrm{Hom}_{\mathcal C\operatorname{-alg}}(\mathcal V[z]^\wedge,M_3\hat V)$ (see Lem. \ref{lem:actions:0704}(a)) 
    and with the morphism of pointed sets \eqref{pre:F}, so that \eqref{pre:F} gives rise to a morphism of pointed sets with group actions 
    $$
    (\mathcal S,(\rho_{\mathrm{DT}},R_{\mathrm{DT}},C_{\mathrm{DT}}),
    \mathrm{GL}_3\hat V\times\mathbf k[[e_0,f_\infty]]^\times,\bullet)
    \to(\mathrm{Hom}_{\mathcal C\operatorname{-alg}}(\mathcal V[z]^\wedge,M_3\hat V),\tilde\rho_{\mathrm{DT}},\mathrm{GL}_3\hat V,\bullet).
    $$

    (b) The morphism of pointed sets \eqref{pre:G} is equivariant with respect to the actions of $\mathrm{GL}_3\hat V$ on its source 
    (see Lem. \ref{lem:actions:0704}(a)) and on its target (see Lem. \ref{lem:520:2212:FIRST}(a)),
    so that 
    \eqref{pre:G} gives rise to a morphism of pointed sets with group actions 
    $$
(\mathrm{Hom}_{\mathcal C\operatorname{-alg}}(\mathcal V[z]^\wedge,M_3\hat V),\tilde\rho_{\mathrm{DT}},\mathrm{GL}_3\hat V,\bullet)
\to (\mathrm{Hom}_{\mathcal C\operatorname{-alg}}(\hat{\mathcal V},M_3\hat V),\rho_{\mathrm{DT}},\mathrm{GL}_3\hat V,\bullet). 
    $$
\end{lem}

\begin{proof}
Denote by $(\rho,R,C)\mapsto \tilde\rho_{(\rho,R,C)}$ the map \eqref{pre:G}. 
Let $(\rho,R,C)\in \mathcal S$ and $(P,\varphi)\in\mathrm{GL}_3\hat V\times\mathbf k[[e_0,f_\infty]]^\times$. 
Then $\tilde\rho_{(\mathrm{Ad}_P\circ\rho,\varphi RP^{-1},
PC\varphi^{-1})}|_{\hat{\mathcal V}}=\mathrm{Ad}_P\circ\rho=(\mathrm{Ad}_P\circ\tilde\rho)|_{\hat{\mathcal V}}$
and 
$$
\tilde\rho_{(\mathrm{Ad}_P\circ\rho,\varphi RP^{-1},
PC\varphi^{-1})}(z)=PC\varphi^{-1}\cdot \varphi RP^{-1}=PCRP^{-1}=\mathrm{Ad}_P\circ\tilde\rho_{(\rho,R,C)}(z),
$$
therefore 
$\tilde\rho_{(\mathrm{Ad}_P\circ\rho,\varphi RP^{-1},
PC\varphi^{-1})}=\mathrm{Ad}_P\circ\tilde\rho_{(\rho,R,C)}$, which is the middle equality in 
$$
\tilde\rho_{(P,\varphi)\bullet(\rho,R,C)}=\tilde\rho_{(\mathrm{Ad}_P\circ\rho,\varphi RP^{-1},
PC\varphi^{-1})}=\mathrm{Ad}_P\circ\tilde\rho_{(\rho,R,C)}=P\bullet \tilde\rho_{(\rho,R,C)}.
$$
One derives $\tilde\rho_{(P,\varphi)\bullet(\rho,R,C)}=P\bullet \tilde\rho_{(\rho,R,C)}$, which implies (a). 
(b) follows from $(\mathrm{Ad}_P\circ \tilde\rho)\circ i_{\mathcal V,\mathcal V[z]}=\mathrm{Ad}_P\circ (\tilde\rho\circ 
i_{\mathcal V,\mathcal V[z]})$ for any $\tilde\rho\in \mathrm{Hom}_{\mathcal C\operatorname{-alg}}(\mathcal V[z]^\wedge,M_3\hat V)$
and $P\in\mathrm{GL}_3\hat V$. 
\end{proof}

By Lem. \ref{lem:13:10:7avril}, the diagram of pointed sets \eqref{diag:ptt:sets:7avril} is upgraded to a diagram of pointed sets with group actions, 
where the underlying group diagram is 
\begin{equation}\label{diag:psga:final:10apr}
\xymatrix{&
\substack{(\mathcal S,(\rho_{\mathrm{DT}},R_{\mathrm{DT}},C_{\mathrm{DT}}), \\ \mathrm{GL}_3\hat V\times\mathbf k[[e_0,f_\infty]]^\times,\bullet)}
\ar[ld]\ar[rd] & \\ 
\substack{(\mathrm{Hom}_{\mathcal C\operatorname{-alg}}(\mathcal V[z]^\wedge,M_3\hat V),\\
\tilde\rho_{\mathrm{DT}},\mathrm{GL}_3\hat V,\bullet)}
\ar[rr]& & 
\substack{(\mathrm{Hom}_{\mathcal C\operatorname{-alg}}(\hat{\mathcal V},M_3\hat V),\\
\rho_{\mathrm{DT}},\mathrm{GL}_3\hat V,\bullet)}}    
\end{equation}

where the diagonal maps are projections and the horizontal map is the identity. 

\subsection{Overall action of $\mathcal G_{\mathrm{inert}}$: a diagram of $\mathcal G_{\mathrm{inert}}$-pointed sets 
with group actions}\label{sect:13:4}

Replacing $\mathcal G$ by $\mathcal G_{\mathrm{inert}}$ in Def. \ref{def:psga}, one defines the notions of 
$\mathcal G_{\mathrm{inert}}$-equivariant pointed sets with a group action ($\mathcal G_{\mathrm{inert}}$-PSGA) and 
$\mathcal G_{\mathrm{inert}}$-PSGA morphisms; proves that these notions build up a category $\mathcal G_{\mathrm{inert}}$-$\mathbf{PSGA}$;  
defines the category $\mathbf{PS}_{\mathcal G_{\mathrm{inert}}}$ of pointed sets with an action of $\mathcal G_{\mathrm{inert}}$. Doing the same replacement in Lem. \ref{lem:foncteur:q}, one defines a functor
$$
\mathbf{q}_{\mathrm{inert}} : \mathbf{PSGA}_{\mathcal G_{\mathrm{inert}}}\to \mathbf{PS}_{\mathcal G_{\mathrm{inert}}}, 
$$
which takes an object $(X,x_0,A,\bullet,*)$, with $(X,x_0)$ a pointed set, $A$ a group, $\bullet$ an action of $A$ on $X$, $*$ 
a pair of compatible actions of $\mathcal G_{\mathrm{inert}}$ on $A$ and $X$, to the triple $(A\backslash X,A\bullet x_0,*)$, 
where $(A\backslash X,A\bullet x_0)$ is the pointed set  obtained by factorization by the action of $A$ and $*$ is the action of 
$\mathcal G_{\mathrm{inert}}$ on $A\backslash X$ induced by its action on $X$. 

\begin{lem}\label{lem:13:11:8april}
(a)  The restriction to $\mathcal G_{\mathrm{inert}}$ of the action of $\mathcal G$ on the set 
$\mathrm{Hom}_{\mathcal C\operatorname{-alg}}(\hat{\mathcal V},M_3\hat V)$ defined 
in Lem. \ref{lem:520:2212:BIS}(a) defines an action of $\mathcal G_{\mathrm{inert}}$ on 
$\mathrm{Hom}_{\mathcal C\operatorname{-alg}}(\hat{\mathcal V},M_3\hat V)$, denoted $(g,\rho)\mapsto g*\rho$. 

(b)  The restriction to $\mathcal G_{\mathrm{inert}}$ of the action of $\mathcal G$ on the group $\mathrm{GL}_3\hat V$ defined 
in Lem. \ref{lem:5:18:2012} defines an action of $\mathcal G_{\mathrm{inert}}$ on $\mathrm{GL}_3\hat V$, denoted $(g,P)\mapsto g*P$. 
 
(c) Together with the action $\bullet$ from Lem. \ref{lem:520:2212:FIRST}(a), the actions from (a) and (b) 
satisfy the compatibility condition from Def. \ref{def:psga}(a). 

(d) The tuple $(\mathrm{Hom}_{\mathcal C\operatorname{-alg}}(\hat{\mathcal V},M_3\hat V),\rho_{\mathrm{DT}},
\mathrm{GL}_3\hat V,\bullet,*)$ is a  $\mathcal G_{\mathrm{inert}}$-pointed set with group action, equal to the restriction 
from $\mathcal G$ to $\mathcal G_{\mathrm{inert}}$ of the diagram with the same notation from \eqref{diag:set:gp:acts:with:*}. 
\end{lem}

\begin{proof}
(a), (b) and (c) respectively follow from Lem. \ref{lem:520:2212:BIS}(a), Lem. \ref{lem:520:2212:FIRST}(a) and Lem. \ref{lem:pre:psga}(a), 
and (d) is a direct consequence.
\end{proof}

\begin{lem}\label{lem:13:12:10apr}
 (a)    The assignment $(g,\tilde\rho)\mapsto g*\tilde\rho:=\mathrm{aut}_g^V\circ \tilde\rho\circ (\mathrm{aut}_g^{\mathcal V[z]})^{-1}$ 
defines an action of  $\mathcal G_{\mathrm{inert}}$ on the set 
$\mathrm{Hom}_{\mathcal C\operatorname{-alg}}(\mathcal V[z]^\wedge,M_3\hat V)$. 

 (b) The action of $\mathcal G_{\mathrm{inert}}$ from (a) and its action on $\mathrm{GL}_3\hat V$ from Lem. \ref{lem:13:11:8april}(b)
 satisfy, together with the action $\bullet$ from Lem. \ref{lem:actions:0704}(a), the  compatibility condition from Def. \ref{def:psga}(a). 
 
(c) The tuple $(\mathrm{Hom}_{\mathcal C\operatorname{-alg}}(\mathcal V[z]^\wedge,M_3\hat V),\tilde\rho_{\mathrm{DT}},
\mathrm{GL}_3\hat V,\bullet,*)$ is a  $\mathcal G_{\mathrm{inert}}$-pointed set with group action. 
\end{lem}

\begin{proof}
 (a) follows from Lem. \ref{lem:action:1342:0105}. Let us prove (b). For $g\in\mathcal G_{\mathrm{inert}}$, 
 $P\in\mathrm{GL}_3\hat V$ and $\tilde\rho\in\mathrm{Hom}_{\mathcal C\operatorname{-alg}}(\mathcal V[z]^\wedge,M_3\hat V)$, 
 one has 
 $$
g*(P\bullet\tilde\rho)
 =\mathrm{aut}_g^V\circ (\mathrm{Ad}_{P}\circ \tilde\rho)\circ 
 (\mathrm{aut}_g^{\mathcal V[z]})^{-1}   
=\mathrm{Ad}_{g*P}\circ \mathrm{aut}_{g}^V\circ 
\tilde\rho\circ 
(\mathrm{aut}_{g}^{\mathcal V[z]})^{-1} =(g*P)\bullet(g*\tilde\rho), 
 $$
 using the equality $\mathrm{aut}_g^V\circ\mathrm{Ad}_P=\mathrm{Ad}_{g*P}\circ \mathrm{aut}_g^V$ used in 
 the proof of Lem. \ref{lem:520:2212:FIRST}(a). (c) is a direct consequence of (a) and (b). 
\end{proof}

\begin{lem}\label{lemm:commutants}
Let $A$ be an algebra, $\alpha$ be an automorphism of $A$, and $B\subset A$ be a subalgebra. Then
$\mathrm{C}_A(\alpha(B))=\alpha(\mathrm{C}_A(B))$. 
\end{lem}

\begin{proof}
One has 
\begin{align*}
& \mathrm{C}_A(\alpha(B))=\cap_{x\in \alpha(B)}\mathrm{ker}(\mathrm{ad}_x)=
 \cap_{b\in B}\mathrm{ker}(\mathrm{ad}_{\alpha(b)})
 =\cap_{b\in B}\mathrm{ker}(\alpha\circ \mathrm{ad}_{b}\circ\alpha^{-1})
\\ &  
=\cap_{b\in B}\alpha(\mathrm{ker}(\mathrm{ad}_{b}))
=\alpha(\cap_{b\in B}\mathrm{ker}(\mathrm{ad}_{b}))
=\alpha(\mathrm{C}_A(B)), 
\end{align*}
 where for $a\in A$, $\mathrm{ad}_a : A\to A$ is $x\mapsto [a,x]$. The first and 
 last equalities follow from the definition of $\mathrm{C}_A(\alpha(B))$ and $\mathrm{C}_A(B)$, 
 the second equality is a substitution, the third equality follows from 
 $\mathrm{ad}_{\alpha(b)}=\alpha\circ \mathrm{ad}_{b}\circ\alpha^{-1}$, 
the fourth equality follows from $\mathrm{ker}(\beta\circ\phi\circ\alpha^{-1})=\alpha(\mathrm{ker}(\phi))$
for $\phi$ a module morphism and $\alpha,\beta$ automorphisms of its source and target, and the fifth equality follows
from the compatibility of module automorphisms with intersections. 
\end{proof}

\begin{lem}\label{lem:13:14:10apr}
(a) The assignment 
\begin{equation}\label{def:act:mathcalGinert:SET}
(g,(\rho,R,C))\mapsto g*(\rho,R,C):=(g*\rho,(g\otimes h_g)^{-1}\mathrm{aut}_g^V(R),\mathrm{aut}_g^V(C)(g\otimes h_g)), 
\end{equation}
 where $g*\rho$ is as in Lem. \ref{lem:13:11:8april}(a) defines an action of the group $\mathcal G_{\mathrm{inert}}$ on the set $\mathcal S$. 

(b) The product of the action of $\mathcal G_{\mathrm{inert}}$ on $\mathrm{GL}_3\hat V$ given by restriction 
of the action of $\mathcal G$ from Lem. \ref{lem:5:18:2012} with the trivial action on $\mathbf k[[e_0,f_\infty]]^\times$ 
is an action of 
$\mathcal G_{\mathrm{inert}}$ on the group $\mathrm{GL}_3\hat V\times\mathbf k[[e_0,f_\infty]]^\times$ given by 
$(g,(P,\varphi))\mapsto g*(P,\varphi)=(g*P,\varphi)$.  

(c)  The action of $\mathcal G_{\mathrm{inert}}$ on the set $\mathcal S$ from (a) and its action on the group 
$\mathrm{GL}_3\hat V\times \mathbf k[[e_0,f_\infty]]^\times$ from (b)
 satisfy, together with the action $\bullet$ from Lem. \ref{lem:actions:0704}(b), the  compatibility condition from Def. \ref{def:psga}(a). 
 
(d) The tuple $(\mathcal S,(\rho_{\mathrm{DT}},R_{\mathrm{DT}},C_{\mathrm{DT}}),
\mathrm{GL}_3\hat V\times\mathbf k[[e_0,f_\infty]]^\times,\bullet,*)$ is a  
$\mathcal G_{\mathrm{inert}}$-pointed set with group action. 
\end{lem}

\begin{proof}
Let us prove (a). 

Let $g\in \mathcal G_{\mathrm{inert}}$ and $(\rho,R,C)\in \mathcal S$ and let us show that $g*(\rho,R,C)\in \mathcal S$. 
By  Lem. \ref{lem:13:11:8april}(a), $g*\rho\in \mathrm{Hom}_{\mathcal C\operatorname{-alg}}(\hat{\mathcal V},M_3\hat V)$; one also 
checks $\mathrm{aut}_g^V(R)(g\otimes g)\in M_{13}\hat V$ and $(g\otimes g)^{-1}\mathrm{aut}_g^V(C)\in M_{31}F^1\hat V$, so 
$g*(\rho,R,C)\in \mathrm{Hom}_{\mathcal C\operatorname{-alg}}(\hat{\mathcal V},M_3\hat V)\times M_{13}\hat V\times M_{31}F^1\hat V$.  

One has 
\begin{align*}
& \mathrm{aut}_g^V(\mathrm{C}_{\hat V}(e_0))=\mathrm{C}_{\hat V}(\mathrm{aut}_g^V(e_0))
=\mathrm{C}_{\hat V}(g(e_0,e_1)e_0g(e_0,e_1)^{-1})=g(e_0,e_1)\mathrm C_{\hat V}(e_0)g(e_0,e_1)^{-1}
\\&
=(g\otimes h_g)\mathrm C_{\hat V}(e_0)(g\otimes h_g)^{-1},
\end{align*}
where the first (resp. third) equality follows from the automorphism status
of $\mathrm{aut}_g^V$ (resp. conjugation by $g(e_0,e_1)$) and Lem. \ref{lemm:commutants}, the second equality follows from 
$\mathrm{aut}_g^V(e_0)=g(e_0,e_1)e_0g(e_0,e_1)^{-1}$, and the last equality follows from $g(e_0,e_1)=g\otimes 1$ and 
$\mathrm C_{\hat V}(e_0)=\mathbf k[[e_0]]\hat\otimes\hat{\mathcal V}$, therefore 
\begin{equation}\label{commutant:autom}
\mathrm{aut}_g^V(\mathrm C_V(e_0))=(g\otimes h_g)\mathrm C_{\hat V}(e_0)(g\otimes h_g)^{-1}. 
\end{equation}
One has 
\begin{align}\label{in:SET:a}
    & \mathrm C_{\hat V}((g*\rho)(\hat{\mathcal V}))
    =\mathrm C_{\hat V}(\mathrm{aut}_g^V\circ\rho\circ (\mathrm{aut}_g^{\mathcal V})^{-1}(\hat{\mathcal V}))
=\mathrm C_{\hat V}(\mathrm{aut}_g^V(\rho(\hat{\mathcal V})))
=\mathrm{aut}_g^V(\mathrm C_{\hat V}(\rho(\hat{\mathcal V})))
=\mathrm{aut}_g^V(\mathbf k I_3+C\mathrm C_V(e_0)R)
    \\& \nonumber
=\mathbf k I_3+\mathrm{aut}_g^V(C)\mathrm{aut}_g^V(\mathrm C_V(e_0))\mathrm{aut}_g^V(R)
=\mathbf k I_3+\mathrm{aut}_g^V(C)(g\otimes h_g)\mathrm C_V(e_0)(g\otimes h_g)^{-1}\mathrm{aut}_g^V(R),
\end{align}
where the first equality follows from the definition of $g*\rho$, the second (resp. fifth) 
equality follows from the automorphism status of $\mathrm{aut}_g^{\mathcal V}$ (resp. $\mathrm{aut}_g^{\mathcal V}$), 
the third equality follows from Lem. \ref{lemm:commutants}, the fourth equality follows from 
$(\rho,R,C)\in \mathcal S$, and the last equality follows from \eqref{commutant:autom}. 

Moreover, since $(\rho,R,C)\in \mathcal S$, 
\begin{equation}\label{eq:RC}
RC=-(e_0+f_\infty).     
\end{equation} 
Then
\begin{align}\label{in:SET:b}
& (g\otimes h_g)^{-1}\mathrm{aut}_g^V(R) \cdot \mathrm{aut}_g^V(C)(g\otimes h_g)
=(g\otimes h_g)^{-1}\mathrm{aut}_g^V(R)\mathrm{aut}_g^V(C)(g\otimes h_g)
\\&\nonumber=(g\otimes h_g)^{-1}\mathrm{aut}_g^V(RC)(g\otimes h_g)=(g\otimes h_g)^{-1}\mathrm{aut}_g^V(-(e_0+f_\infty)^{-1})(g\otimes h_g)
\\& \nonumber =- (e_0+f_\infty)
\end{align}
where the second equality follows from the automorphism status of $\mathrm{aut}_g^V$, the third equality follows from 
\eqref{eq:RC}, and the fourth equality follows from the combination of \eqref{def:aut:g:V:(1)}
and \eqref{useful:duality:1553}. \eqref{in:SET:a} and \eqref{in:SET:b} imply 
\begin{equation}\label{a:first:part}
    \forall (g,(\rho,R,C))\in\mathcal G_{\mathrm{inert}}\times \mathcal S, \quad g*(\rho,R,C)\in \mathcal S. 
\end{equation}
 
Let now $g,g'\in \mathcal G_{\mathrm{inert}}$ and $(\rho,R,C)\in \mathcal S$. Then 
\begin{align}\label{a:second:part}
& \nonumber g*(g'*(\rho,R,C))=g*(g'*\rho,(g'\otimes h_{g'})^{-1}\mathrm{aut}_{g'}^V(R),\mathrm{aut}_{g'}^V(C)(g'\otimes h_{g'}))
\\& \nonumber =(g*(g'*\rho),(g\otimes h_{g})^{-1}\mathrm{aut}_{g}^V((g'\otimes h_{g'})^{-1}\mathrm{aut}_{g'}^V(R)),
\mathrm{aut}_{g}^V(\mathrm{aut}_{g'}^V(C)(g'\otimes h_{g'}))(g\otimes h_{g}))
\\& \nonumber=((g\circledast g')*\rho,
(g\otimes h_{g})^{-1}\mathrm{aut}_{g}^V((g'\otimes h_{g'})^{-1})\mathrm{aut}_{g}^V\circ\mathrm{aut}_{g'}^V(R)),
\mathrm{aut}_{g}^V\circ\mathrm{aut}_{g'}^V(C)\mathrm{aut}_{g}^V(g'\otimes h_{g'})(g\otimes h_{g}))
\\&=\nonumber((g\circledast g')*\rho,
(\mathrm{aut}_{g}^{\mathcal V}(g')g\otimes \mathrm{aut}_{g}^{\mathcal V}(h_{g'})h_{g}))^{-1}\mathrm{aut}_{g\circledast g'}^V(R),
\mathrm{aut}_{g\circledast g'}^V(C)(\mathrm{aut}_{g}^{\mathcal V}(g')g\otimes \mathrm{aut}_{g}^{\mathcal V}(h_{g'})h_{g}))
\\&=\nonumber((g\circledast g')*\rho,
((g\circledast g')\otimes h_{g\circledast g'})^{-1}\mathrm{aut}_{g\circledast g'}^V(R),
\mathrm{aut}_{g\circledast g'}^V(C)((g\circledast g')\otimes h_{g\circledast g'}))
\\& =(g\circledast g')*(\rho,R,C),      
\end{align}
where the two first and last equalities follows from \eqref{def:act:mathcalGinert:SET}, the third equality follows from 
Lem. \ref{lem:520:2212:BIS}(a) and the algebra automorphism status of $\mathrm{aut}_{g}^V$, the fourth equality follows from 
Lem. \ref{lem:2:8:jan25}(c), the fifth equality follows from the identities
$$
g\circledast g'=\mathrm{aut}_g^{\mathcal V}(g')g, \quad h_{g\circledast g'}=\mathrm{aut}_{g}^{\mathcal V}(h_{g'})h_g
$$
of which the former follows from \eqref{FORMULA:1912} and \eqref{def:aut:g:V:(1)} and the latter follows from \eqref{toto:3006}. 
The identities \eqref{a:first:part} and \eqref{a:second:part} imply (a). (b) is obvious. 

(c) Let $g\in \mathcal G_{\mathrm{inert}}$, $(P,\varphi)\in\mathrm{GL}_3\hat V\times\mathbf k[[e_0,f_\infty]]^\times$
and $(\rho,R,C)\in \mathcal S$. Then  
\begin{align*}
 & g*((P,\varphi)\bullet(\rho,R,C))
 =g*(\mathrm{Ad}_P\circ\rho,\varphi RP^{-1},PC\varphi^{-1})
 \\&\nonumber=(g*(\mathrm{Ad}_P\circ\rho),(g\otimes h_g)^{-1}\mathrm{aut}_g^V(\varphi RP^{-1}),\mathrm{aut}_g^V(PC\varphi^{-1})(g\otimes h_g))
\\&\nonumber=(\mathrm{Ad}_{\mathrm{aut}_g^V(P)}\circ(g*\rho),(g\otimes h_g)^{-1}
\mathrm{aut}_g^V(\varphi)\mathrm{aut}_g^V(R)\mathrm{aut}_g^V(P)^{-1},
\mathrm{aut}_g^V(P)\mathrm{aut}_g^V(C)\mathrm{aut}_g^V(\varphi)^{-1}(g\otimes h_g))
\\&\nonumber=(\mathrm{Ad}_{\mathrm{aut}_g^V(P)}\circ(g*\rho),
\varphi(g\otimes h_g)^{-1}\mathrm{aut}_g^V(R)\mathrm{aut}_g^V(P)^{-1},
\mathrm{aut}_g^V(P)\mathrm{aut}_g^V(C)(g\otimes h_g)\varphi^{-1})
\\& 
=(\mathrm{aut}_g^V(P),\varphi)\bullet(g*\rho,(g\otimes h_g)^{-1}\mathrm{aut}_g^V(R),\mathrm{aut}_g^V(C)(g\otimes h_g))
\\ & =(g*(P,\varphi))\bullet(g*(\rho,C,R)),  
\end{align*}
where the first and fifth equalities follows from the definition in Lem. \ref{lem:actions:0704}(b), the second equality follows 
from \eqref{def:act:mathcalGinert:SET}, the third equality follows from the identity 
$\mathrm{aut}_g^V\circ \mathrm{Ad}_P=\mathrm{Ad}_{\mathrm{aut}_g^V(P)} \circ \mathrm{aut}_g^V$
and the algebra automorphism status of $\mathrm{aut}_g^V$, 
the fourth equality follows from the fact that the restriction of $\mathrm{aut}_g^V$ to 
$\mathbf k[[e_0,f_\infty]]$ coincides with $\mathrm{Ad}_{g\otimes h_g}$, which follows from 
$\mathrm{aut}_g^{\mathcal V}(e_0)=\mathrm{Ad}_g(e_0)$ and $\mathrm{aut}_g^{\mathcal V}(e_\infty)=\mathrm{Ad}_{h_g}(e_\infty)$, 
and the last equality follows from (b) and \eqref{def:act:mathcalGinert:SET}.  

(d) follows from (a)-(c). 
\end{proof}

\begin{lem}\label{lem:13:15:10apr}
(a)  The group morphism $\mathrm{GL}_3\hat V\times\mathbf k[[e_0,f_\infty]]^\times\to\mathrm{GL}_3\hat V$ defined by projection on the first 
factor is equivariant with respect the actions of $\mathcal G_{\mathrm{inert}}$ on its source as in Lem. \ref{lem:13:14:10apr}(b) and on its 
target by restriction of the action of $\mathcal G$ from Lem. \ref{lem:5:18:2012}. 

(b) The map $(\mathcal S,(\rho_{\mathrm{DT}},R_{\mathrm{DT}},C_{\mathrm{DT}}))\to 
    (\mathrm{Hom}_{\mathcal C\operatorname{-alg}}(\mathcal V[z]^\wedge,M_3\hat V),\tilde\rho_{\mathrm{DT}})$ from \eqref{pre:G} 
    is equivariant with respect the actions of $\mathcal G_{\mathrm{inert}}$ on its source by Lem. \ref{lem:13:14:10apr}(a) 
    and on its target by Lem. \ref{lem:13:12:10apr}(a). 

(c) The map $(\mathrm{Hom}_{\mathcal C\operatorname{-alg}}(\mathcal V[z]^\wedge,M_3\hat V),\tilde\rho_{\mathrm{DT}})
    \to(\mathrm{Hom}_{\mathcal C\operatorname{-alg}}(\hat{\mathcal V},M_3\hat V),\rho_{\mathrm{DT}})$ from \eqref{pre:F} 
    is equivariant with respect the actions of $\mathcal G_{\mathrm{inert}}$ on its source by Lem. \ref{lem:13:12:10apr}(a) 
    and on its target by restriction of the action of $\mathcal G$ defined in Lem. \ref{lem:520:2212:BIS}(a). 

(d) The morphisms of pointed sets with group actions from Lem. \ref{lem:13:10:7avril}(a) and (b) induce morphisms 
$$
(\mathcal S,(\rho_{\mathrm{DT}},R_{\mathrm{DT}},C_{\mathrm{DT}}),\mathrm{GL}_3\hat V\times\mathbf k[[e_0,f_\infty]]^\times,\bullet,*)
\to(\mathrm{Hom}_{\mathcal C\operatorname{-alg}}(\mathcal V[z]^\wedge,M_3\hat V),\tilde\rho_{\mathrm{DT}},\mathrm{GL}_3\hat V,\bullet,*).
$$
and
$$
(\mathrm{Hom}_{\mathcal C\operatorname{-alg}}(\mathcal V[z]^\wedge,M_3\hat V),\tilde\rho_{\mathrm{DT}},\mathrm{GL}_3\hat V,\bullet,*)
\to (\mathrm{Hom}_{\mathcal C\operatorname{-alg}}(\hat{\mathcal V},M_3\hat V),\rho_{\mathrm{DT}},\mathrm{GL}_3\hat V,\bullet,*). 
$$
of $\mathcal G_{\mathrm{inert}}$-pointed sets with group actions, leading to the following commutative 
diagram of $\mathcal G_{\mathrm{inert}}$-pointed sets with group actions upgrading the 
diagram of pointed sets with group actions \eqref{diag:psga:final:10apr}:  
\begin{equation}\label{diag:mathcalGintert:psga}
\xymatrix{&
\substack{(\mathcal S,(\rho_{\mathrm{DT}},R_{\mathrm{DT}},C_{\mathrm{DT}}), \\ \mathrm{GL}_3\hat V\times\mathbf k[[e_0,f_\infty]]^\times,\bullet,*)}
\ar[ld]\ar[rd] & \\ 
\substack{(\mathrm{Hom}_{\mathcal C\operatorname{-alg}}(\mathcal V[z]^\wedge,M_3\hat V),\\
\tilde\rho_{\mathrm{DT}},\mathrm{GL}_3\hat V,\bullet,*)}
\ar[rr]& & 
\substack{(\mathrm{Hom}_{\mathcal C\operatorname{-alg}}(\hat{\mathcal V},M_3\hat V),\\
\rho_{\mathrm{DT}},\mathrm{GL}_3\hat V,\bullet,*)}}    
\end{equation}
\end{lem}

\begin{proof}
(a) is obvious. 

(b) Let $g\in \mathcal G_{\mathrm{inert}}$ and $(\rho,R,C)\in \mathcal S$. Then 
\begin{equation}\label{eq:pessah:1}
\forall v\in\hat{\mathcal V},\quad  g*\tilde\rho_{(\rho,R,C)}(v)
=\mathrm{aut}_g^V \circ \tilde\rho_{(\rho,R,C)}\circ (\mathrm{aut}_g^{\mathcal V[z]})^{-1}(v)
=\mathrm{aut}_g^V \circ \rho\circ (\mathrm{aut}_g^{\mathcal V})^{-1}(v)
=(g*\rho)(v),    
\end{equation}
where the first equality follows from the definition of $g*\tilde\rho_{(\rho,R,C)}$, the second equality follows from 
the definition of $\tilde\rho_{(\rho,R,C)}$, from the stability of $\hat{\mathcal V}\subset \mathcal V[z]^\wedge$
by $(\mathrm{aut}_g^{\mathcal V[z]})^{-1}$ and the equality of its restriction to this subspace with 
$(\mathrm{aut}_g^{\mathcal V})^{-1}$, and the last equality follows from the definition of $g*\rho$. 
One also has 
\begin{align}\label{eq:pessah:2}
& \nonumber g*\tilde\rho_{(\rho,R,C)}(z)
=\mathrm{aut}_g^V \circ \tilde\rho_{(\rho,R,C)}\circ (\mathrm{aut}_g^{\mathcal V[z]})^{-1}(z)
=\mathrm{aut}_g^V \circ \tilde\rho_{(\rho,R,C)}(z)
=\mathrm{aut}_g^V (CR)
    \\&
=\mathrm{aut}_g^V (C)(g\otimes h_g)\cdot (g\otimes h_g)^{-1}\mathrm{aut}_g^V(R), 
\end{align}
where 
where the last equality follows from the automorphism status of $\mathrm{aut}_g^V$ and all the other equalities follow from 
definitions. Since $(g*\rho,(g\otimes h_g)^{-1}\mathrm{aut}_g^V(R),\mathrm{aut}_g^V(C)(g\otimes h_g))$ belongs to $\mathcal S$
and is equal to $g*(\rho,R,C)$, \eqref{eq:pessah:1} and \eqref{eq:pessah:2} imply 
$g*\tilde\rho_{(\rho,R,C)}=\tilde\rho_{g*(\rho,R,C)}$. (b) follows.

(c) Let $g\in \mathcal G_{\mathrm{inert}}$ and $\tilde\rho\in \mathrm{Hom}_{\mathcal C\operatorname{-alg}}(\mathcal V[z]^\wedge,M_3\hat V)$.
Then $g*\tilde\rho|_{\hat{\mathcal V}}=\mathrm{aut}_g^V\circ \tilde\rho|_{\hat{\mathcal V}}\circ (\mathrm{aut}_g^{\mathcal V})^{-1}
=\mathrm{aut}_g^V\circ \tilde\rho\circ i_{\mathcal V,\mathcal V[z]}\circ (\mathrm{aut}_g^{\mathcal V})^{-1}
=\mathrm{aut}_g^V\circ \tilde\rho\circ (\mathrm{aut}_g^{\mathcal V[z]})^{-1}\circ i_{\mathcal V,\mathcal V[z]}
=(g*\tilde\rho)\circ i_{\mathcal V,\mathcal V[z]}=(g*\tilde\rho)|_{\hat{\mathcal V}}$, where the third equality follows from 
$i_{\mathcal V,\mathcal V[z]}\circ (\mathrm{aut}_g^{\mathcal V})^{-1}=(\mathrm{aut}_g^{\mathcal V[z]})^{-1}\circ i_{\mathcal V,\mathcal V[z]}$
and all the other equalities follows from definitions, which implies (c). 

(d) follows from (a)-(c) and (morphisms of pointed sets with group actions). 
\end{proof}

\begin{defn}\label{def:13:16:11apr}
The composition of the morphisms of $\mathcal G_{\mathrm{inert}}$-pointed sets 
\begin{align*}
&  ((\mathrm{GL}_3\hat V\times\mathbf k[[e_0,f_\infty]]^\times)\backslash \mathcal S,
(\mathrm{GL}_3\hat V\times\mathbf k[[e_0,f_\infty]]^\times)\bullet(\rho_{\mathrm{DT}},R_{\mathrm{DT}},C_{\mathrm{DT}}),*)
\\& \to(\mathrm{GL}_3\hat V\backslash\mathrm{Hom}_{\mathcal C\operatorname{-alg}}(\mathcal V[z]^\wedge,M_3\hat V),
\mathrm{GL}_3\hat V\bullet\tilde\rho_{\mathrm{DT}},*)
\end{align*}
and
$$ 
(\mathrm{GL}_3\hat V\backslash\mathrm{Hom}_{\mathcal C\operatorname{-alg}}(\mathcal V[z]^\wedge,M_3\hat V),
\mathrm{GL}_3\hat V\bullet\tilde\rho_{\mathrm{DT}},*)
\to (\mathrm{GL}_3\hat V\backslash\mathrm{Hom}_{\mathcal C\operatorname{-alg}}(\hat{\mathcal V},M_3\hat V),
\mathrm{GL}_3\hat V\bullet\rho_{\mathrm{DT}},*)    
$$
obtained by applying the functor $\mathbf q_{\mathrm{inert}}$ to the morphisms of $\mathcal G_{\mathrm{inert}}$-pointed sets 
with group actions from Lem. \ref{lem:13:15:10apr}(d) is denoted 
\begin{align*}
& (F) : 
((\mathrm{GL}_3\hat V\times\mathbf k[[e_0,f_\infty]]^\times)\backslash \mathcal S,
(\mathrm{GL}_3\hat V\times\mathbf k[[e_0,f_\infty]]^\times)\bullet(\rho_{\mathrm{DT}},R_{\mathrm{DT}},C_{\mathrm{DT}}),*)
\\& 
\to (\mathrm{GL}_3\hat V\backslash\mathrm{Hom}_{\mathcal C\operatorname{-alg}}(\hat{\mathcal V},M_3\hat V),
\mathrm{GL}_3\hat V\bullet\rho_{\mathrm{DT}},*).  
\end{align*}
\end{defn}

\begin{lem}\label{lem:13:17:11apr}
The subgroups of $\mathcal G_{\mathrm{inert}}$ obtained as the stabilizers of the $\mathcal G_{\mathrm{inert}}$-pointed sets
$((\mathrm{GL}_3\hat V\times\mathbf k[[e_0,f_\infty]]^\times)\backslash \mathcal S,
(\mathrm{GL}_3\hat V\times\mathbf k[[e_0,f_\infty]]^\times)\bullet(\rho_{\mathrm{DT}},R_{\mathrm{DT}},C_{\mathrm{DT}}),*)$, 
$$(\mathrm{GL}_3\hat V\backslash\mathrm{Hom}_{\mathcal C\operatorname{-alg}}(\mathcal V[z]^\wedge,M_3\hat V),
\mathrm{GL}_3\hat V\bullet\tilde\rho_{\mathrm{DT}},*)$$ and 
$(\mathrm{GL}_3\hat V\backslash\mathrm{Hom}_{\mathcal C\operatorname{-alg}}(\hat{\mathcal V},M_3\hat V),
\mathrm{GL}_3\hat V\bullet\rho_{\mathrm{DT}},,*)$ from Def. \ref{def:13:16:11apr} satisfy the inclusions
$$
\mathrm{Stab}_{\mathcal G_{\mathrm{inert}}}((\mathrm{GL}_3\hat V\times\mathbf k[[e_0,f_\infty]]^\times)\bullet
(\rho_{\mathrm{DT}},R_{\mathrm{DT}},C_{\mathrm{DT}}))\subset 
\mathrm{Stab}_{\mathcal G_{\mathrm{inert}}}(\mathrm{GL}_3\hat V\bullet\tilde\rho_{\mathrm{DT}})
\subset
\mathrm{Stab}_{\mathcal G_{\mathrm{inert}}}(\mathrm{GL}_3\hat V\bullet\rho_{\mathrm{DT}}) . 
$$
\end{lem}

\begin{proof}
The morphisms from Def. \ref{def:13:16:11apr} build up the following commutative triangle of 
$\mathcal G_{\mathrm{inert}}$-pointed sets
\begin{equation}\label{diag:13:4:11}
\xymatrix{&
\substack{((\mathrm{GL}_3\hat V\times\mathbf k[[e_0,f_\infty]]^\times)\backslash \mathcal S, \\
(\mathrm{GL}_3\hat V\times\mathbf k[[e_0,f_\infty]]^\times)
\bullet(\rho_{\mathrm{DT}},R_{\mathrm{DT}},C_{\mathrm{DT}}),*)}
\ar_{}[ld]\ar^{(F)}[rd] & \\ 
\substack{(\mathrm{GL}_3\hat V\backslash\mathrm{Hom}_{\mathcal C\operatorname{-alg}}(\mathcal V[z]^\wedge,M_3\hat V),\\
\mathrm{GL}_3\hat V\bullet\tilde\rho_{\mathrm{DT}},*)}
\ar_{}[rr]& & 
\substack{(\mathrm{GL}_3\hat V\backslash\mathrm{Hom}_{\mathcal C\operatorname{-alg}}(\hat{\mathcal V},M_3\hat V),\\
\mathrm{GL}_3\hat V\bullet\rho_{\mathrm{DT}},*)}}    
\end{equation}
which upon applying the stabilizer groups functor gives rise to the diagram of subgroups of 
$\mathcal G_{\mathrm{inert}}$: 
$$
\xymatrix{& \scriptstyle{\mathrm{Stab}_{\mathcal G_{\mathrm{inert}}}((\mathrm{GL}_3\hat V\times\mathbf k[[e_0,f_\infty]]^\times)\bullet
(\rho_{\mathrm{DT}},R_{\mathrm{DT}},C_{\mathrm{DT}}))}
\ar@{_{(}->}[ld]\ar@{^{(}->}[rd]& \\ 
\scriptstyle{\mathrm{Stab}_{\mathcal G_{\mathrm{inert}}}(\mathrm{GL}_3\hat V\bullet\tilde\rho_{\mathrm{DT}})}
\ar@{^{(}->}[rr]& &
\scriptstyle{\mathrm{Stab}_{\mathcal G_{\mathrm{inert}}}(\mathrm{GL}_3\hat V\bullet\rho_{\mathrm{DT}})} }
$$  
leading to the claimed inclusions. 
\end{proof}

\subsection{Material for \S\ref{sect:inj:H}: algebraic results}\label{sect:13:5}

\begin{lem}\label{lem:A:0505}
If $u\in V$ is such that $[u,e_0]
\in \mathrm{C}_V(e_0)$, then $u\in \mathrm{C}_V(e_0)$.     
\end{lem}

\begin{proof}
Equip $V$ with the algebra grading $V=\oplus_{i\geq0}V_{(i)}$ such that $e_1$ has degree 1 and 
$e_0,f_0,f_1$ have degree 0. Then 
\begin{equation}\label{toto:0505}
 V_{(0)}=\mathbf k[e_0]\otimes\mathcal V=\mathrm{C}_V(e_0),    
\end{equation}
where the first equality follows from the definition of the grading on $V$, and the second equality follows from 
Lem. \ref{lem:comm:e0}(a). Let now $u\in V$ be such that  $[u,e_0]\subset \mathrm{C}_V(e_0)$. 
Let $u=\sum_{i\geq0}u_{(i)}$ be the decomposition of $u$. 
By \eqref{toto:0505}, the assumption on $u$ implies $\sum_{i\geq0}[u_{(i)},e_0]\in V_{(0)}$. 
Since $[u_{(i)},e_0]\in V_{(i)}$ for any $i$, and since 
$V_{(i)}\cap V_{(0)}=0$ for $i>0$, this implies $[u_{(i)},e_0]=0$ for any $i>0$; therefore 
for any $i>0$, $u_{(i)}\in \mathrm{C}_V(e_0)$, which by \eqref{toto:0505} implies $u_{(i)}\in V_{(0)}$;
since $V_{(i)}\cap V_{(0)}=0$, one has $u_{(i)}=0$, therefore $u=u_{(0)}\in V_{(0)}$, which by 
\eqref{toto:0505} implies the statement. 
\end{proof}

\begin{lem}\label{lem:B:0505}
(a) If $u\in \hat V^\times$ is such that $u \mathrm{C}_{\hat V}(e_0)u^{-1}
=\mathrm{C}_{\hat V}(e_0)$, then $u\in \mathrm{C}_{\hat V}(e_0)^\times$.     

(b) Let $u,v\in \hat V^\times$ be such that $u\mathrm C_{\hat V}(e_0)v=\mathrm C_{\hat V}(e_0)$, then $u,v\in \mathrm C_{\hat V}(e_0)^\times$. 
\end{lem}

\begin{proof}
(a) There exists a family of polynomials $(P_n(x_1,\ldots,z))_{n\geq 1}$ in the free noncommutative variables $x_1,\ldots,z$, such that for each $n\geq 1$, $P_n(x_1,\ldots,z)$ has 
total degree $n$, where $\deg(x_i)=i$ and $\deg(z)=1$, and degree one with respect to $z$, and 
$$
(1+\sum_{i\geq 0}x_i)\cdot z\cdot (1+\sum_{i\geq 0}x_i)^{-1}=\sum_{i\geq 1}P_i(x_1,\ldots,z).  
$$
Then $P_1(x_1,\ldots,z)=z$, and for any $i>0$, there exists a polynomial $Q_{i+1}(x_1,\ldots,x_{i-1},z)$
such that  $P_{i+1}(x_1,\ldots,z)=[x_i,z]+Q_{i+1}(x_1,\ldots,x_{i-1},z)$. 

Denote by $\hat V=\hat \oplus_{i\geq0}V_i$ the decomposition of $\hat V$ for the total degree (for which $e_i,f_i$ ($i=0,1$) have degree 1). 
It follows from Lem. \ref{lem:comm:e0} that $\mathrm{C}_{\hat V}(e_0)$ is graded with respect to this grading, and that
$\mathrm{C}_{\hat V}(e_0)=\hat\oplus_{i\geq0}\mathrm{C}_V(e_0)_i$, where $\mathrm{C}_V(e_0)_i$ is the total degree $i$ part of 
$\mathrm{C}_V(e_0)$. 

Let $u$ be as in the assumption of (a) and let $u=\sum_{i\geq0}u_i$ be the its decomposition, so $u_i\in V_i$ for any $i$. 
Let us show inductively on $n\geq 0$ that 
$u_n\in\mathrm{C}_V(e_0)_i$. One has $u_0\in \mathbf k^\times\subset\mathbf k=\mathrm{C}_V(e_0)_0$; dividing 
$u$ by $u_0$, we henceforth assume that $u_0=1$. 
Assume that $n>0$ and that $u_i\in \mathrm{C}_V(e_0)_i$ for $i<n$. Then since $e_0\in \mathrm{C}_{\hat V}(e_0)$
and by the assumption on $u$, one has 
\begin{equation}\label{cond:lambda}
ue_0 u^{-1}\in\mathrm{C}_{\hat V}(e_0).     
\end{equation}
The degree $n+1$ component of $u e_0 u^{-1}$ is equal to 
$[u_n,e_0]+Q_{n+1}(u_1,\ldots,u_{n-1},e_0)$, therefore 
\eqref{cond:lambda} implies $[u_n,e_0]+Q_{n+1}(u_1,\ldots,u_{n-1},e_0)\in \mathrm{C}_{V}(e_0)$, 
which by the induction assumption and since $e_0\in \mathrm{C}_{V}(e_0)$ implies  
$[u_n,e_0]\in \mathrm{C}_{V}(e_0)$, which by Lem. \ref{lem:A:0505}
implies $u_n\in \mathrm{C}_{V}(e_0)$, therefore $u_n\in \mathrm{C}_{V}(e_0)_n$. It follows that 
$u\in \mathrm{C}_{\hat V}(e_0)^\times$, as claimed. 

(b) Let $u,v$ be as in the assumptions of (b). As $1\in\mathrm C_{\hat V}(e_0)$, one has 
\begin{equation}\label{uv:invertible}
uv\in \mathrm C_{\hat V}(e_0)^\times.     
\end{equation} 
Then $u\mathrm C_{\hat V}(e_0) u^{-1}=u\mathrm C_{\hat V}(e_0) v(uv)^{-1}=\mathrm C_{\hat V}(e_0) (uv)^{-1}
=\mathrm C_{\hat V}(e_0)$, where the last equality follows from \eqref{uv:invertible}. 
The resulting equality $u\mathrm C_{\hat V}(e_0) u^{-1}=\mathrm C_{\hat V}(e_0)$ implies  $u\in \mathrm{C}_{\hat V}(e_0)^\times$
by (a). Then $v=u^{-1}(uv)\in\mathrm{C}_{\hat V}(e_0)^\times$, where the equality follows from 
$u\in \mathrm{C}_{\hat V}(e_0)^\times$ and the statement "$\in$" follows from the combination of $u\in \mathrm{C}_{\hat V}(e_0)^\times$
and  \eqref{uv:invertible}. (b) follows. 
\end{proof}

\begin{lem}\label{lem:comm:alg}
Let $n\geq1$ and $(a_1,\ldots,a_n),(b_1,\ldots,b_n)$ be non-colinear vectors in $\mathbb Q^n$. 
Let $x_1,\ldots,x_n$ be free commutative variables, $a:=\sum_i a_ix_i$, $b:=\sum_i b_ix_i$, so 
$a,b\in \mathbf k[x_1,\ldots,x_n]$. Then the sequence 
$$
\mathbf k[x_1,\ldots,x_n]\to\mathbf k[x_1,\ldots,x_n]^{\oplus2}\to\mathbf k[x_1,\ldots,x_n], 
$$
where the first map is $P\mapsto (a\cdot P,b\cdot P)$ and the second map is $(A,B)\mapsto b\cdot A-a\cdot B$, is exact. 
\end{lem}

\begin{proof}
Using the action of a suitable element of $\mathrm{GL}_n(\mathbb Q)$, one may assume $a=x_1$, $b=x_2$, 
in which case the statement follows from an argument on the coefficients the involved polynomials.
\end{proof}

\begin{lem}\label{lem:exact:seq:0505}
The sequence of $\mathbf k$-module morphisms 
\begin{equation}\label{exact:seq:0705}
\mathbf k[[e_0,f_\infty]]\oplus \hat V\to \hat V^{\oplus 2}\to \hat V, 
\end{equation}
where the first map is $(\varphi,\gamma)\mapsto (\varphi,\varphi)+(\gamma(e_0+f_\infty),(e_0+f_\infty)\gamma)$
and the second map is $(u,v)\mapsto (e_0+f_\infty)u-v(e_0+f_\infty)$, is exact. 
\end{lem}

\begin{proof}
For $\Sigma\subset\mathcal V$ a subset, let ${\mathtt M}(\Sigma)$ be the submonoid of $\mathcal V$ generated by $\Sigma$ 
(we denote the unit by $\emptyset$). Then ${\mathtt M}(e_0,e_1)\simeq\{e_0,e_1\}^*$ is a $\mathbf k$-module basis of $\mathcal V$. 
Let ${}_{1}{\mathtt M}^0_{1}$ be its submonoid consisting 
of the unit $\emptyset$ and of the elements starting with and ending in $e_1$; then there are set inclusions 
 ${}_{1}{\mathtt M}^0_{1}\subset{\mathtt M}(e_0,e_1)\subset\mathcal V$. 

There is a bijection ${\mathtt M}(e_0)\sqcup({\mathtt M}(e_0)\times({}_{1}{\mathtt M}^0_{1}\smallsetminus\{\emptyset\})
\times{\mathtt M}(e_0))\to {\mathtt M}(e_0,e_1)$, where ${\mathtt M}(e_0)\to {\mathtt M}(e_0,e_1)$ is the canonical
injection and ${\mathtt M}(e_0)\times({}_{1}{\mathtt M}^0_{1}\smallsetminus\{\emptyset\})
\times{\mathtt M}(e_0)\to {\mathtt M}(e_0,e_1)$ is induced by the product. This implies that  
the $\mathbf k[e_0]$-bimodule $\mathcal{V}$ is decomposed as the direct sum
$$
\mathcal V=\oplus_{w\in{} {}_{1}{\mathtt M}^0_{1}} \mathcal V_0(w),\quad \mathrm{where}\quad 
\mathcal V_0(w):=\mathrm{im}(
\textbf{k}
[e_0]^{\otimes 2}\to \mathcal V,a\otimes b\mapsto a\cdot w\cdot b)
$$
is the $\mathbf k[e_0]$-subbimodule of $\mathcal V$ generated by $w$.  
When $w\neq\emptyset$ (resp. $w=\emptyset$), the map $\mathbf k[e_0]^{\otimes 2}\to \mathcal{V}_0(w)$ induced by 
$a\otimes b\mapsto a\cdot w\cdot b$ (resp. the map $\mathbf k[e_0]\to \mathcal V_0(\emptyset)$ induced by the inclusion 
$\mathbf k[e_0]\subset\mathcal V$) induces an isomorphism of $\mathbf k[e_0]$-bimodules, the 
$\mathbf k[e_0]$-bimodule (i.e. $\mathbf k[e_0^{(l)},e_0^{(r)}]=\mathbf k[e_0]^{\otimes 2}$-module) structure on 
$\mathbf k[e_0^{(l)},e_0^{(r)}]=\mathbf k[e_0]^{\otimes 2}$ (resp. on $\mathbf k[e_0]$) being the regular one
(resp. induced by the product $\mathbf k[e_0]^{\otimes 2}\to\mathbf k[e_0]$, i.e. the morphism 
$\mathbf k[e_0^{(l)},e_0^{(r)}]\to \mathbf k[e_0]$, $e_0^{(l)},e_0^{(r)}\mapsto e_0$). 

Applying to the above situation the algebra automorphism $(e_1\mapsto e_1,e_0\leftrightarrow e_\infty)$ of $\mathcal V$
(recall $e_\infty:=-e_0-e_1$), one defines a submonoid ${\mathtt M}(e_\infty,e_1)$ of $\mathcal V$,
its submonoid ${}_{1}{\mathtt M}^\infty_{1}$, and the set inclusions ${}_{1}{\mathtt M}_{1}^\infty\subset{\mathtt M}(e_\infty,e_1)
\subset \mathcal{V}$; we then have the $\mathbf k[e_\infty]$-bimodule decomposition
$$
\mathcal V=\oplus_{w\in{} {}_{1}{\mathtt M}^\infty_{1}} \mathcal{V}_\infty(w),\quad \mathrm{where}\quad 
\mathcal{V}_\infty(w):=\mathrm{im}(
\textbf{k}
[e_\infty]^{\otimes 2}\to \mathcal{V},a\otimes b\mapsto a\cdot w\cdot b), 
$$
and the $\mathbf k[e_\infty]$-bimodule isomorphisms $\mathbf k[e_\infty]^{\otimes2}\to\mathcal V_\infty(w)$, $a\otimes b\mapsto 
a\cdot w\cdot b$ for $w\neq\emptyset$, and $\mathbf k[e_\infty]\to\mathcal V_\infty(\emptyset)$ induced by the inclusion 
$\mathbf k[e_\infty]\subset\mathcal V$. 

Recall the notation $e_i:=e_i\otimes1$, $f_i:=1\otimes e_i$ in the tensor square $V=\mathcal V^{\otimes2}$ for
$i\in\{0,1,\infty\}$. The tensor product of the above bimodule decompositions gives rise to a
$\mathbf k[e_0,f_\infty]$-bimodule decomposition 
\begin{equation}\label{direct;sum:V:0705}
V=\underset{(w,z)\in{} {}_{1}{\mathtt M}^0_{1}\times {}_{1}{\mathtt M}^\infty_{1}
}{\oplus}{V({w,z})}, \quad \mathrm{where}\quad 
{V}(w,z):=\mathcal V_0(w)\otimes\mathcal V_\infty(z).
\end{equation}
Then $V(w,z)$ is the $\mathbf k [e_0,f_\infty]$-subbimodule of $V$ generated by $w\otimes z$. 
The above bimodule isomorphisms induce  $\mathbf k [e_0,f_\infty]$-bimodule (i.e. $\mathbf k [e_0^{(l)},f_\infty^{(l)},e_0^{(r)},f_\infty^{(r)}]=\mathbf k [e_0,f_\infty]^{\otimes2}$-module) isomorphisms
\begin{equation}\label{case:240508:1833}
    {V}(w,z)\simeq
\begin{cases}
 \mathbf k[e_0,f_\infty] & \text{ if } w=z=\emptyset, \\ 
  \mathbf k[e_0,f_\infty^{(l)},f_\infty^{(r)}] &  \text{ if }  w=\emptyset, z\neq \emptyset,\\ 
  \mathbf k[e_0^{(l)},e_0^{(r)},f_\infty] & \text{ if }  w\neq\emptyset, z=\emptyset, \\ 
  \mathbf k[e_0^{(l)},e_0^{(r)},f_\infty^{(l)},f_\infty^{(r)}] & \text{ if }  w\neq\emptyset, z\neq\emptyset. \\ 
\end{cases}
\end{equation}
the $\mathbf k [e_0^{(l)},f_\infty^{(l)},e_0^{(r)},f_\infty^{(r)}]$-module structures on the algebras in the right-hand sides being 
induced by the algebra morphisms from  $\mathbf k [e_0^{(l)},f_\infty^{(l)},e_0^{(r)},f_\infty^{(r)}]$ to these algebras 
taking $(e_0^{(x)},f_\infty^{(x)})$ ($x\in\{r,l\}$) to $(e_0,f_\infty)$ in the first case, $(e_0,f_\infty^{(x)})$ in the second 
case, $(e_0^{(x)},f_\infty)$ in the third case,  and the identity algebra morphism in the last case. 

If $M$ is a $\mathbf k[e_0,f_\infty]$-bimodule, let $\alpha_M : M\to M^{\oplus 2}$ and 
$\beta_M : M^{\oplus 2}\to M$ be the maps $m\mapsto (m\cdot(e_0+f_\infty),(e_0+f_\infty)\cdot m)$
and $(m,n)\mapsto (e_0+f_\infty)\cdot m-n\cdot(e_0+f_\infty)$. Let also 
$\mathrm{diag} : \mathbf k[e_0,f_\infty]\to\mathbf k[e_0,f_\infty]^{\oplus2}$, $\varphi\mapsto 
(\varphi,\varphi)$ be the diagonal map. Then $\mathrm{diag}$ as well as $\alpha_M,\beta_M$ are 
$\mathbf k[e_0,f_\infty]$-bimodule morphisms. Moreover, the assignments $M\mapsto\alpha_M,\beta_M$ are 
compatible with direct sums, so 
\begin{equation}\label{direct:sums}
\alpha_{M\oplus N}=\alpha_M\oplus\alpha_N, \quad \beta_{M\oplus N}=\beta_M\oplus\beta_N
\end{equation}
for $M,N$ two $\mathbf k[e_0,f_\infty]$-bimodules. 

The formulas from the statement of the lemma induce a sequence of $\mathbf k$-module morphisms
\begin{equation}\label{sequence:mod}
\mathbf k[e_0,f_\infty]\oplus V\to V^{\oplus 2}\to V,
\end{equation}
which coincides (upon passing from $\mathbf k[e_0,f_\infty]$-bimodule structures to $\mathbf k$-module structures) 
with the following sequence of $\mathbf k[e_0,f_\infty]$-bimodule morphisms
\begin{equation}\label{sequence:bimod}
\mathbf k[e_0,f_\infty]\oplus V\stackrel{(i^{\oplus2}\circ \mathrm{diag})\oplus\alpha_V}{\to} V^{\oplus 2}\stackrel{\beta_V}{\to} V,   
\end{equation}
where $i : \mathbf k[e_0,f_\infty]\to V$ is the canonical injection. 
 
It follows from \eqref{direct:sums}, from the direct sum decomposition \eqref{direct;sum:V:0705}, and from the fact that 
the image of $i$ is contained in $V(\emptyset,\emptyset)$,  
that  \eqref{sequence:bimod} decomposes as the direct sum over $(w,z)\in 
{}_{1}{\mathtt M}^0_{1}\times {}_{1}{\mathtt M}^\infty_{1}\setminus\{(\emptyset,\emptyset)\}$
of the sequence of $\mathbf k[e_0,f_\infty]$-bimodule morphisms  
\begin{equation}\label{sequence:z:w}
V(w,z)\stackrel{\alpha_{V(w,z)}}{\to}V(w,z)^{\oplus 2}
\stackrel{\beta_{V(w,z)}}{\to}V(w,z)    
\end{equation}
and of the sequence of $\mathbf k[e_0,f_\infty]$-bimodule morphisms
\begin{equation}\label{sequence:empty:empty}
V(\emptyset,\emptyset)\oplus\mathbf k[e_0,f_\infty]
\stackrel{\alpha_{V(\emptyset,\emptyset)}\oplus\mathrm{diag}}{\to}V(\emptyset,\emptyset)^{\oplus 2}
\stackrel{\beta_{V(\emptyset,\emptyset)}}{\to}V(\emptyset,\emptyset)    
\end{equation}
corresponding to $(w,z)=(\emptyset,\emptyset)$. 

Depending on the values of $(w,z)\neq(\emptyset,\emptyset)$, the isomorphism 
\eqref{case:240508:1833} sets up an isomorphism between the sequence of $\mathbf k[e_0,f_\infty]$-bimodule morphisms
\eqref{sequence:z:w} and the following sequences of $\mathbf k[e_0^{(l)},e_0^{(r)},f_\infty^{(l)},f_\infty^{(r)}]$-module morphisms: 

$\bullet$ if $w\neq\emptyset$ and $z\neq\emptyset$, the sequence 
$$
\mathbf k[e_0^{(l)},e_0^{(r)},f_\infty^{(l)},f_\infty^{(r)}]
\to \mathbf k[e_0^{(l)},e_0^{(r)},f_\infty^{(l)},f_\infty^{(r)}]^{\oplus 2}
\to \mathbf k[e_0^{(l)},e_0^{(r)},f_\infty^{(l)},f_\infty^{(r)}], 
$$
where the first map is $P\mapsto ((e_0^{(r)}+f_\infty^{(r)})\cdot P,(e_0^{(l)}+f_\infty^{(l)})\cdot P)$ 
and the second map is $(A,B)\mapsto (e_0^{(l)}+f_\infty^{(l)})\cdot A-(e_0^{(r)}+f_\infty^{(r)})\cdot B$. 
This is an exact sequence by Lem. \ref{lem:comm:alg}. 

$\bullet$ if $w=\emptyset$ and $z\neq\emptyset$, the sequence 
$$
\mathbf k[e_0,f_\infty^{(l)},f_\infty^{(r)}]
\to \mathbf k[e_0,f_\infty^{(l)},f_\infty^{(r)}]^{\oplus 2}
\to \mathbf k[e_0,f_\infty^{(l)},f_\infty^{(r)}], 
$$
where the first map is $P\mapsto ((e_0+f_\infty^{(r)})\cdot P,(e_0+f_\infty^{(l)})\cdot P)$ 
and the second map is $(A,B)\mapsto (e_0+f_\infty^{(l)})\cdot A-(e_0+f_\infty^{(r)})\cdot B$. 
This is an exact sequence by Lem. \ref{lem:comm:alg}. 

$\bullet$ if $w\neq\emptyset$ and $z=\emptyset$, the sequence 
$$
\mathbf k[e_0^{(l)},e_0^{(r)},f_\infty]
\to \mathbf k[e_0^{(l)},e_0^{(r)},f_\infty]^{\oplus 2}
\to \mathbf k[e_0^{(l)},e_0^{(r)},f_\infty], 
$$
where the first map is $P\mapsto ((e_0^{(r)}+f_\infty)\cdot P,(e_0^{(l)}+f_\infty)\cdot P)$ 
and the second map is $(A,B)\mapsto (e_0^{(l)}+f_\infty)\cdot A-(e_0^{(r)}+f_\infty)\cdot B$. 
This is an exact sequence by Lem. \ref{lem:comm:alg}. 

When $(w,z)=(\emptyset,\emptyset)$, the isomorphism \eqref{case:240508:1833} sets up an isomorphism between the sequence 
of $\mathbf k[e_0,f_\infty]$-bimodule morphisms \eqref{sequence:empty:empty} and the following sequences of 
$\mathbf k[e_0^{(l)},e_0^{(r)},f_\infty^{(l)},f_\infty^{(r)}]$-module morphisms: 
$$
\mathbf k[e_0,f_\infty]^{\oplus 2}
\to \mathbf k[e_0,f_\infty]^{\oplus 2}
\to \mathbf k[e_0,f_\infty], 
$$
where the first map is $(P,Q)\mapsto (P+(e_0+f_\infty)\cdot Q,P+(e_0+f_\infty)\cdot Q)$ 
and the second map is $(A,B)\mapsto (e_0+f_\infty)\cdot (A-B)$. 
This is an exact sequence since $e_0+f_\infty$ is not a zero divisor in $\mathbf k[e_0,f_\infty]$. 

It follows that the sequences of $\mathbf k[e_0,f_\infty]$-bimodule morphisms  
\eqref{sequence:z:w} and \eqref{sequence:empty:empty} are all exact. This implies the exactness of  
the sequence of $\mathbf k[e_0,f_\infty]$-bimodule morphisms \eqref{sequence:bimod}, and therefore of  
the sequence of $\mathbf k$-module morphisms \eqref{sequence:mod}. The latter sequence of $\mathbf k$-module morphisms
is graded, and \eqref{exact:seq:0705} is its graded completion. The exactness of \eqref{sequence:mod} then 
implies that of \eqref{exact:seq:0705}. 
\end{proof}

\begin{lem}\label{lem:shriek}
If $v,u\in\mathrm C_{\hat V}(e_0)^\times$ satisfy the equality $v(e_0+f_\infty)u=e_0+f_\infty$, then there exists 
$\delta\in\hat V$ and $\varphi\in\mathbf k[[e_0,f_\infty]]^\times$, such that 
\begin{equation}\label{equalities:u:v}
v=\varphi (1-(e_0+f_\infty)\delta )^{-1},\quad 
u=(1-\delta (e_0+f_\infty))\varphi^{-1}. 
\end{equation}
\end{lem}

\begin{proof}
Since $v$ is invertible, one has $v^{-1} (e_0+f_\infty)=(e_0+f_\infty)u$. 
By Lem. \ref{lem:exact:seq:0505}, this implies the existence of a pair $(\tilde\varphi,\gamma)\in 
\mathbf k[[e_0,f_\infty]]\times\hat V$, such that 
\begin{equation}\label{beta:alpha:explicites}
v^{-1}=\tilde\varphi+(e_0+f_\infty)\gamma,\quad 
u=\tilde\varphi+\gamma(e_0+f_\infty). 
\end{equation}
Then $\epsilon(\tilde\varphi)=\epsilon(u)\in\mathbf k^\times$, where the equality follows from the
the second equality in \eqref{beta:alpha:explicites} and the statement "$\in$" follows from $u\in\hat V^\times$. 
This implies $\tilde\varphi\in \mathbf k[[e_0,f_\infty]]^\times$; set $\varphi:=\tilde\varphi^{-1}$. Set 
$$
\delta:=-\gamma\tilde\varphi^{-1}\in\hat V, 
$$
then 
\begin{equation}\label{1952:0505}
v^{-1}=(1- (e_0+f_\infty)\delta)\tilde\varphi=(1- (e_0+f_\infty)\delta)\varphi^{-1}, 
\end{equation}
and 
$$
u=\tilde\varphi-\delta\tilde\varphi(e_0+f_\infty)=\tilde\varphi-\delta(e_0+f_\infty)\tilde\varphi=(1-\delta(e_0+f_\infty))\tilde\varphi
=(1-\delta(e_0+f_\infty))\varphi^{-1},   
$$
where the second equality follows from the fact that $\tilde\varphi$ commutes with $e_0+f_\infty$
(as $\mathbf k[[e_0,f_\infty]]$ is commutative).  
This equality implies the second part of the statement, while \eqref{1952:0505} implies its first part. 
\end{proof}

\begin{lem}\label{lemme:jaune}
Let $F : \hat V\to\hat V$ be the map $\delta\mapsto (1-\delta (e_0+f_\infty))
(1-(e_0+f_\infty)\delta)^{-1}$. 
If $\delta\in \hat V$ is such that $F(\delta)\in \mathrm{C}_{\hat V}(e_0)$, then $\delta\in\mathrm{C}_{\hat V}(e_0)$. 
\end{lem}

\begin{proof}
Let $x_0,x_1,\ldots,z$ be free noncommutative variables, with $\mathrm{deg}(x_i)=i$ and $\mathrm{deg}(z)=1$. 
The let $(\underline F_n)_{n \geq0}$ be the family homogeneous polynomials in these variables 
defined from the expansion  
$$
(1-(\sum_{i\geq0}x_i) z) (1-z (\sum_{i\geq0}x_i))^{-1}=\sum_{n\geq0}\underline F_n(x_0,x_1,\ldots,z),\quad 
\mathrm{deg}(\underline F_i)=i. 
$$ 
(equality of noncommutative formal series). Then $\underline F_0=1$, $\underline F_1=[z,x_0]$; for any $n\geq 0$, 
$\underline F_{n+1}$ has the form $[x_n,z]+\underline G_{n+1}(x_0,\ldots,x_{n-1},z)$, where $\underline G_{n+1}$ is a polynomial 
in the variables $x_0,\ldots,x_{n-1},z$. 

For $\delta\in\hat V$ expanded as $\sum_{i\geq 0}\delta_i$ with respect to the total degree, and 
$F(\delta)=\sum_{i\geq 0}F_i(\delta)$ the total degree expansion of $F(\delta)$ in $\hat V$, 
one has $F_0(\delta)=1$ and for $n\geq 0$, 
\begin{equation}\label{toto:0905}
F_{n+1}(\delta)=\underline F_{n+1}(\delta_0,\ldots,e_0+f_\infty)
=[\delta_n,e_0+f_\infty]+\underline G_{n+1}(\delta_0,\ldots,\delta_{n-1},e_0+f_\infty)    
\end{equation}
Assume now that $\delta\in\hat V$ is such that $F(\delta)\in\mathrm{C}_{\hat V}(e_0+f_\infty)$. 
Since $e_0+f_\infty$ is homogeneous, $\mathrm{C}_{\hat V}(e_0+f_\infty)$ is complete graded, therefore  
\begin{equation}\label{titit:0905}
\forall n\geq 0,\quad  F_n(\delta) \in \mathrm{C}_V(e_0+f_\infty).    
\end{equation}
Let us prove by induction on $n \geq 0$ that $\delta_n \in \mathrm{C}_V(e_0+f_\infty)$. Since 
$\delta_0\in\mathbf k$, one has $\delta_0\in \mathrm{C}_V(e_0+f_\infty)$. Assume that $n\geq 0$ 
and $\delta_0,\ldots,\delta_n \in \mathrm{C}_V(e_0+f_\infty)$. Then $\mathrm{C}_V(e_0+f_\infty) \ni 
F_{n+2}(\delta)=[e_0+f_\infty,\delta_{n+1}]+\underline G_{n+1}(\delta_0,\ldots,\delta_n,e_0+f_\infty)$ where $\ni$ follows 
from \eqref{titit:0905} and the equality follows from \eqref{toto:0905}. Since $\delta_0,\ldots,\delta_n,e_0+f_\infty$ 
belong to $\mathrm{C}_V(e_0+f_\infty)$, $\underline G_{n+1}(\delta_0,\ldots,\delta_n,e_0+f_\infty)$ belongs to the same algebra, 
which implies $[e_0+f_\infty,\delta_{n+1}]\in \mathrm{C}_V(e_0+f_\infty)$.  By Lem. \ref{lem240506:1638}, this implies $\delta_{n+1} 
\in \mathrm{C}_V(e_0+f_\infty)$.
\end{proof}

\begin{prop}\label{prop:4:pages}
If $u,v\in\mathrm C_{\hat V}(e_0)^\times$ satisfy the equality $v(e_0+f_\infty)u=e_0+f_\infty$, then there exists 
$\delta\in\mathrm C_{\hat V}(e_0)$ and $\varphi\in\mathbf k[[e_0,f_\infty]]^\times$, such that 
\begin{equation}\label{equalities:u:v:BIS}
u=(1-\delta (e_0+f_\infty))\varphi^{-1},\quad 
v=\varphi (1-(e_0+f_\infty)\delta )^{-1}.     
\end{equation}    
\end{prop}

\begin{proof}
It follows from Lem. \ref{lem:shriek} that there exist $\varphi\in\mathbf k[[e_0,f_\infty]]^\times$ and $\delta\in\hat V$ such that 
\eqref{equalities:u:v} hold. Then $\mathrm C_{\hat V}(e_0)^\times\ni uv=(1-\delta (e_0+f_\infty))(1-(e_0+f_\infty)\delta )^{-1}$, 
where the equality follows from \eqref{equalities:u:v}. By Lem. \ref{lemme:jaune}, this relation implies
$\delta\in\mathrm C_{\hat V}(e_0)$, which implies the statement. 
\end{proof}

\subsection{Local injectivity of the map (F)}\label{sect:inj:H}

\begin{lem}\label{lem:13:21:15apr}
(a)    Let $(R,C)\in M_{13}F^1\hat V\times M_{31}\hat V$ be such that 
    $$\mathbf kI_3+C\mathrm C_{\hat V}(e_0)R
=\mathbf kI_3+C_{\mathrm{DT}}\mathrm C_{\hat V}(e_0)R_{\mathrm{DT}}$$
(equality if subsets of $M_3\hat V$). Then there exist $u,v\in\hat V^\times$, such that
$C=C_{\mathrm{DT}}u$, $R=v R_{\mathrm{DT}}$. 

(b) In the situation of (a), $u$ and $v$ belong to $\mathrm C_{\hat V}(e_0)^\times$. 
\end{lem}

\begin{proof}
(a) 
If $(R,C)\in M_{13}F^1\hat V\times M_{31}\hat V$, then 
\begin{equation}\label{F1:commutant}    
(\mathbf kI_3+C\mathrm C_{\hat V}(e_0)R)\cap M_3F^1\hat V=C\mathrm C_{\hat V}(e_0)R;
\end{equation} 
the inclusion of the right-hand side in the left-hand side is obvious, and if 
$(\lambda,c)\in\mathbf k\times \mathrm C_{\hat V}(e_0)$ is 
such that $\lambda I_3+CcR\in M_3F^1\hat V$, then since 
$C\mathrm C_{\hat V}(e_0)R\subset M_3F^1\hat V$ one gets $\lambda I_3\in M_3 F^1\hat V$, which implies $\lambda=0$; this 
proves the inclusion of the left-hand side in the right-hand side. 

Assume now $\mathbf kI_3+C\mathrm C_{\hat V}(e_0)R
=\mathbf kI_3+C_{\mathrm{DT}}\mathrm C_{\hat V}(e_0)R_{\mathrm{DT}}$. 
By \eqref{F1:commutant}, the intersection of this equality with $M_3F^1\hat V$ yields 
$$
C\mathrm C_{\hat V}(e_0)R=C_{\mathrm{DT}}\mathrm C_{\hat V}(e_0)R_{\mathrm{DT}}. 
$$
Let $\alpha,\beta,\gamma\in F^1\hat V$ and $s,t,u\in\hat V$ be such that $C=\begin{pmatrix}
    \alpha\\\beta\\\gamma
\end{pmatrix}$ and $R=\begin{pmatrix}
    s&t&u
\end{pmatrix}$, then the latter equality is written as
\begin{equation}\label{matrix:equality}
  \begin{pmatrix}
    \alpha\\\beta\\\gamma
\end{pmatrix}\mathrm C_{\hat V}(e_0)\begin{pmatrix}
    s&t&u
\end{pmatrix}=\begin{pmatrix}
    f_1\\e_1\\-(e_0+f_\infty)
\end{pmatrix}\mathrm C_{\hat V}(e_0)\begin{pmatrix}
    0&0&1
\end{pmatrix}  
\end{equation}
(equality of $\mathbf k$-submodules of $M_3F^1\hat V$), whose image by the projection $M_3F^1\hat V\to F^1\hat V$ 
corresponding to the (1,3) entry yields 
\begin{equation}\label{alpha:c:u=f1:c}
\alpha \mathrm C_{\hat V}(e_0)u=f_1\mathrm C_{\hat V}(e_0).
\end{equation}
(equality of $\mathbf k$-submodules of $F^1\hat V$). Let $\alpha=\sum_{i\geq1}\alpha_i$ and $u=\sum_{i\geq0}u_i$ be the degree
decompositions of $\alpha$ and $u$, then the image of the latter equality by the projection $F^1\hat V\to V_1$ is 
$$
\alpha_1\mathbf k u_0=f_1\mathbf k. 
$$
(equality of $\mathbf k$-submodules of $V_1$). 
Let $\alpha_1=\sum_{x\in\{e_0,e_1,f_0,f_1\}}\alpha_{x}x$ be the decomposition of 
$\alpha_1$ is the basis $\{e_0,e_1,f_0,f_1\}$ of the free $\mathbf k$-module $V_1=\oplus_{x\in \{e_0,e_1,f_0,f_1\}}\mathbf kx$; then 
the projection of the latter equality on $\mathbf kf_1$ gives $\alpha_{f_1}\mathbf ku_0=\mathbf k$, which implies that 
 $\alpha_{f_1}\in \mathbf k^\times$ and $u_0\in \mathbf k^\times$. The latter relation implies 
\begin{equation}\label{c:invertible}
 u\in\hat V^\times.   
\end{equation}
The projection of the same equality on $\mathbf kx$ for 
$x\in \{e_0,e_1,f_0\}$ gives $\alpha_{x}\mathbf ku_0=0$, which since $u_0\in\mathbf k^\times$ implies $\alpha_{x}=0$. 
All this implies
\begin{equation}\label{equality:alpha} 
\alpha_1=\alpha_{f_1}f_1,\quad \alpha_{f_1}\in\mathbf k^{\times}.  
\end{equation}

The image of \eqref{matrix:equality} by the projection $M_3F^1\hat V\to F^1\hat V$ corresponding to the (2,3) entry yields
$$
\alpha \mathrm C_{\hat V}(e_0)t=0, 
$$
therefore $\alpha t=0$. Assume $t\neq0$ and let $d\geq0$ be the smallest integer such that the degree $d$ part $t_d$ of $b$ is nonzero. 
Then $\alpha_1t_d=0$, which by $\alpha_{f_1}\in\mathbf k^{\times}$, hence by \eqref{equality:alpha}, $\alpha_{f_1}f_1t_d=0$,
and the injectivity of the map $V_d\to V_{d+1}$, 
$x\mapsto f_1x$ implies $t_d=0$, a contradiction; therefore 
\begin{equation}\label{b:zero}
 t=0.   
\end{equation}
Applying the projection $M_3F^1\hat V\to F^1\hat V$ corresponding to the (1,3) similarly yields 
\begin{equation}\label{a:zero}
s=0. 
\end{equation}
The projection of \eqref{matrix:equality} corresponding to $M_3\hat V\to M_{31}\hat V$ associated with the last column
yields 
$$ 
  \begin{pmatrix}
    \alpha\\\beta\\\gamma
\end{pmatrix}\mathrm C_{\hat V}(e_0)u=\begin{pmatrix}
    f_1\\e_1\\-(e_0+f_\infty)
\end{pmatrix}\mathrm C_{\hat V}(e_0) 
$$
(equality of $\mathbf k$-submodules of $M_{31}F^1\hat V$). 
It follows 
$$ 
  \begin{pmatrix}
    \alpha\\\beta\\\gamma
\end{pmatrix}\mathrm C_{\hat V}(e_0)u\hat V=\begin{pmatrix}
    f_1\\e_1\\-(e_0+f_\infty)
\end{pmatrix}\mathrm C_{\hat V}(e_0)\hat V.  
$$
Since $u\in\hat V^\times$, one has $u\hat V=\hat V$; moreover, $\mathrm C_{\hat V}(e_0)\hat V=\hat V$, therefore 
$$ 
  \begin{pmatrix}
    \alpha\\\beta\\\gamma
\end{pmatrix}\hat V=\begin{pmatrix}
    f_1\\e_1\\-(e_0+f_\infty)
\end{pmatrix}\hat V,   
$$
which implies the existence of $v,w\in\hat V$ such that 
\begin{equation}\label{equation:15apr}
\begin{pmatrix}
    \alpha\\\beta\\\gamma
\end{pmatrix}=\begin{pmatrix}
    f_1\\e_1\\-(e_0+f_\infty)
\end{pmatrix}v,\quad 
\begin{pmatrix}
    f_1\\e_1\\-(e_0+f_\infty)
\end{pmatrix}= \begin{pmatrix}
    \alpha\\\beta\\\gamma
\end{pmatrix} w ; 
\end{equation}
then $\begin{pmatrix}
    f_1\\e_1\\-(e_0+f_\infty)
\end{pmatrix}(1-vw)=0$, which since $x\mapsto f_1x$ is injective implies $vw=1$, therefore $v\in\hat V^\times$. 
The statement follows from this combined with the first equality in \eqref{equation:15apr}, and 
from \eqref{c:invertible}, \eqref{b:zero} and \eqref{a:zero}.  

(b) Combining the first entry of the first relation of \eqref{equation:15apr} with \eqref{alpha:c:u=f1:c}, one obtains 
the equality $f_1v \mathrm C_{\hat V}(e_0)u=f_1\mathrm C_{\hat V}(e_0)$ of subsets of $\hat V$. The injectivity of the endomorphism 
$x\mapsto f_1x$ of $\hat V$ then implies $v \mathrm C_{\hat V}(e_0)u=\mathrm C_{\hat V}(e_0)$, which by Lem. \ref{lem:B:0505}(b) implies 
$u,v\in\mathrm C_{\hat V}(e_0)^\times$, as claimed. 
\end{proof}

\begin{lem}\label{lem:erev:erev:pessah}
The preimage of $\rho_{\mathrm{DT}}$ by the map $\mathcal S\to\mathrm{Hom}_{\mathcal C\operatorname{-alg}}(\hat{\mathcal V},M_3\hat V)$, 
$(\rho,R,C)\mapsto \rho$ is $(\mathrm C_3(\rho_{\mathrm{DT}}(\hat{\mathcal V}))^\times\times\mathbf k[[e_0,f_\infty]]^\times)\bullet(\rho_{\mathrm{DT}},R_{\mathrm{DT}},C_{\mathrm{DT}})$ (recall the algebra inclusion 
$\mathrm C_3(\rho_{\mathrm{DT}}(\hat{\mathcal V}))\subset M_3\hat V$, from which one derives the group inclusion 
$\mathrm C_3(\rho_{\mathrm{DT}}(\hat{\mathcal V}))^\times\subset\mathrm{GL}_3\hat V$). 
\end{lem}

\begin{proof}
 Let $(\rho,R,C)$ belong to the said preimage. Then $\rho=\rho_{\mathrm{DT}}$, and 
\begin{equation}\label{conditions:preimage:15apr}
\mathrm C_3(\rho(\hat{\mathcal V}))=\mathbf kI_3+C\mathrm C_{\hat V}(e_0)R,\quad RC=-(e_0+f_\infty). 
\end{equation}
The first part of \eqref{conditions:preimage:15apr} implies $\mathbf kI_3+C\mathrm C_{\hat V}(e_0)R
=\mathbf kI_3+C_{\mathrm{DT}}\mathrm C_{\hat V}(e_0)R_{\mathrm{DT}}$, which by Lem. \ref{lem:13:21:15apr}(b) implies 
the existence of $u,v\in\mathrm C_{\hat V}(e_0)^\times$, such that 
\begin{equation}\label{form:of:R:C}
    C=C_{\mathrm{DT}}u,\quad R=v R_{\mathrm{DT}}. 
\end{equation}
Then 
$$
e_0+f_\infty=-RC=-v R_{\mathrm{DT}}C_{\mathrm{DT}}u=v(e_0+f_\infty)u
$$
where the first (resp. second) equality follows from second part of \eqref{conditions:preimage:15apr} 
(resp. \eqref{form:of:R:C}). 
By Prop. \ref{prop:4:pages}, the resulting equality $v(e_0+f_\infty)u=e_0+f_\infty$, together with $u,v\in\mathrm C_{\hat V}(e_0)^\times$, 
implies the existence of $\delta\in\mathrm C_{\hat V}(e_0)$ and $\varphi\in\mathbf k[[e_0,f_\infty]]^\times$, such that 
$$
u=(1-\delta (e_0+f_\infty))\varphi^{-1},\quad 
v=\varphi (1-(e_0+f_\infty)\delta )^{-1}. 
$$
Set then 
$$
P:=I_3+C_{\mathrm{DT}}\delta R_{\mathrm{DT}}\in M_3\hat V. 
$$
Then the image of $P$ in $M_3\mathbf k$ is $I_3$, which implies that $M$ is invertible; since $\delta\in \mathrm C_{\hat V}(e_0)$, 
it belongs to $\mathrm C_3(\rho_{\mathrm{DT}})$, hence $P\in \mathrm C_3(\rho_{\mathrm{DT}})^\times$, which implies
$$
\mathrm{Ad}_P\circ\rho_{\mathrm{DT}}=\rho_{\mathrm{DT}}. 
$$
Then 
$$
C=C_{\mathrm{DT}}u=C_{\mathrm{DT}}(1+\delta R_{\mathrm{DT}} C_{\mathrm{DT}})\varphi^{-1}=(1+C_{\mathrm{DT}}\delta 
R_{\mathrm{DT}} )C_{\mathrm{DT}}\varphi^{-1}=PC_{\mathrm{DT}}\varphi^{-1}, 
$$
$$
R=vR_{\mathrm{DT}}=\varphi (1+R_{\mathrm{DT}} C_{\mathrm{DT}}\delta )^{-1}R_{\mathrm{DT}}=\varphi R_{\mathrm{DT}}(1+C_{\mathrm{DT}}\delta 
R_{\mathrm{DT}} )^{-1}=\varphi R_{\mathrm{DT}}P^{-1}. 
$$
The three last identities imply $(\rho,R,C)=(P,\varphi)\bullet(\rho_{\mathrm{DT}},R_{\mathrm{DT}},C_{\mathrm{DT}})$, the action 
being as in Lem. \ref{lem:actions:0704}(b), 
which implies $(\rho,R,C)\in (\mathrm C_3(\rho_{\mathrm{DT}}(\hat{\mathcal V}))^\times
\times\mathbf k[[e_0,f_\infty]]^\times) \bullet(\rho_{\mathrm{DT}},R_{\mathrm{DT}},C_{\mathrm{DT}})$. 
All this implies the inclusion 
$$
\{(\rho,R,C)\in \mathcal S|\rho=\rho_{\mathrm{DT}}\}\subset (\mathrm C_3(\rho_{\mathrm{DT}}(\hat{\mathcal V}))^\times
\times\mathbf k[[e_0,f_\infty]]^\times) \bullet(\rho_{\mathrm{DT}},R_{\mathrm{DT}},C_{\mathrm{DT}}). 
$$
The opposite inclusion follows from the fact that the image of $(\mathrm C_3(\rho_{\mathrm{DT}}(\hat{\mathcal V}))^\times
\times\mathbf k[[e_0,f_\infty]]^\times) \bullet(\rho_{\mathrm{DT}},R_{\mathrm{DT}},C_{\mathrm{DT}})\subset \mathcal S\to 
\mathrm{Hom}_{\mathcal C\operatorname{-alg}}(\hat{\mathcal V},M_3\hat V)$ is 
$\mathrm C_3(\rho_{\mathrm{DT}}(\hat{\mathcal V}))^\times\bullet \rho_{\mathrm{DT}}$, which is equal to 
$\{\rho_{\mathrm{DT}}\}$. All this proves the claim. 
\end{proof}

\begin{prop}\label{prop:inj:H}
The morphism of pointed sets 
\begin{align*}
&((\mathrm{GL}_3\hat V\times\mathbf k[[e_0,f_\infty]]^\times)\backslash \mathcal S,
(\mathrm{GL}_3\hat V\times\mathbf k[[e_0,f_\infty]]^\times)\bullet(\rho_{\mathrm{DT}},R_{\mathrm{DT}},C_{\mathrm{DT}}))
\\&\to 
(\mathrm{GL}_3\hat V\backslash\mathrm{Hom}_{\mathcal C\operatorname{-alg}}(\hat{\mathcal V},M_3\hat V),
\mathrm{GL}_3\hat V\bullet\rho_{\mathrm{DT}})    
\end{align*}
underlying the morphism (F) from Def. \ref{def:13:16:11apr} is locally injective. 
\end{prop}

\begin{proof}
Let $\alpha$ belong to the preimage of $\mathrm{GL}_3\hat V\bullet\rho_{\mathrm{DT}}$ by this map and let $(\rho,R,C)$
belong to $\alpha$. Then the image if $\rho$ by  $\mathcal S\to\mathrm{Hom}_{\mathcal C\operatorname{-alg}}(\hat{\mathcal V},M_3\hat V)$
belongs to $\mathrm{GL}_3\hat V\bullet\rho_{\mathrm{DT}}$, therefore there exists $P\in\mathrm{GL}_3\hat V$, such that 
$\rho=P\bullet \rho_{\mathrm{DT}}$. Since $(P,1)\in\mathrm{GL}_3\hat V\times \mathbf k[[e_0,f_\infty]]^\times$ 
is a lift of $P$, the element $(P,1)^{-1}\bullet(\rho,R,C)$ of $\mathcal S$ belongs to the fiber of $\rho_{\mathrm{DT}}$ by 
$\mathcal S\to\mathrm{Hom}_{\mathcal C\operatorname{-alg}}(\hat{\mathcal V},M_3\hat V)$. By Lem. \ref{lem:erev:erev:pessah}, there exist
$(c,s)\in \mathrm C_3(\rho_{\mathrm{DT}}(\hat{\mathcal V}))^\times\times\mathbf k[[e_0,f_\infty]]^\times$, such that 
$(P,1)^{-1}\bullet(\rho,R,C)=(c,s)\bullet(\rho_{\mathrm{DT}},R_{\mathrm{DT}},C_{\mathrm{DT}})$. Therefore
$(\rho,R,C)=(Pc,s)\bullet (\rho_{\mathrm{DT}},R_{\mathrm{DT}},C_{\mathrm{DT}})$, therefore 
$\alpha=(\mathrm{GL}_3\hat V\times\mathbf k[[e_0,f_\infty]]^\times)\bullet(\rho_{\mathrm{DT}},R_{\mathrm{DT}},C_{\mathrm{DT}})$. 
The claim follows. 
\end{proof}

\begin{cor}\label{cor:13:29}
The subgroups $\mathrm{Stab}_{\mathcal G_{\mathrm{inert}}}(\mathrm{GL}_3\hat V\bullet\tilde\rho_{\mathrm{DT}})$ and 
$\mathrm{Stab}_{\mathcal G_{\mathrm{inert}}}(\mathrm{GL}_3\hat V\bullet\rho_{\mathrm{DT}}). $ of $\mathcal G_{\mathrm{inert}}$ are equal, i.e. 
$$
\mathrm{Stab}_{\mathcal G_{\mathrm{inert}}}(\mathrm{GL}_3\hat V\bullet\tilde\rho_{\mathrm{DT}})
=\mathrm{Stab}_{\mathcal G_{\mathrm{inert}}}(\mathrm{GL}_3\hat V\bullet\rho_{\mathrm{DT}}). 
$$
\end{cor}

\begin{proof}
The combination of Prop. \ref{prop:inj:H} and Lem. \ref{lem:general:2212}(b) implies the equality 
$\mathrm{Stab}_{\mathcal G_{\mathrm{inert}}}((\mathrm{GL}_3\hat V\times\mathbf k[[e_0,f_\infty]]^\times)\bullet
(\rho_{\mathrm{DT}},R_{\mathrm{DT}},C_{\mathrm{DT}}))=\mathrm{Stab}_{\mathcal G_{\mathrm{inert}}}
(\mathrm{GL}_3\hat V\bullet\rho_{\mathrm{DT}})$ of subgroups of $\mathcal G_{\mathrm{inert}}$. 
The result follows from the combination of this equality with Lem. \ref{lem:13:17:11apr}. 
\end{proof}

\subsection{A $\mathcal G_{\mathrm{inert}}\rtimes(\mathbb Z/2\mathbb Z)$-pointed set with group action}\label{sect:amrzpswga}

Replacing $\mathcal G$ by $\mathcal G_{\mathrm{inert}}\rtimes(\mathbb Z/2\mathbb Z)$ in Def. \ref{def:psga}, one defines
the notion of a $\mathcal G_{\mathrm{inert}}\rtimes(\mathbb Z/2\mathbb Z)$-pointed set with group action, 
which is a tuple $(X,x_0,A,\bullet,\tilde*)$, where $(X,x_0,A,\bullet)$ is a pointed set with group action and 
$\tilde*$ is an action of $\mathcal G_{\mathrm{inert}}\rtimes(\mathbb Z/2\mathbb Z)$ both on the group $A$ and on the set $X$, 
which are compatible in the sense of Def. \ref{def:psga}(a). Such a structure gives rise to a 
 $\mathcal G_{\mathrm{inert}}\rtimes(\mathbb Z/2\mathbb Z)$-pointed set $(A\backslash X,A\bullet x_0,\tilde*)$. 
 By restricting the actions of the group $\mathcal G_{\mathrm{inert}}\rtimes(\mathbb Z/2\mathbb Z)$ to the subgroup
 $\mathcal G_{\mathrm{inert}}$, 
 a $\mathcal G_{\mathrm{inert}}\rtimes(\mathbb Z/2\mathbb Z)$-pointed set with group action $(X,x_0,A,\bullet,\tilde*)$ gives rise to a 
 $\mathcal G_{\mathrm{inert}}$-pointed set with group action $(X,x_0,A,\bullet,*)$. The restriction of  
 $\mathcal G_{\mathrm{inert}}\rtimes(\mathbb Z/2\mathbb Z)$ to the subgroup $\mathcal G_{\mathrm{inert}}$
 also takes the $\mathcal G_{\mathrm{inert}}\rtimes(\mathbb Z/2\mathbb Z)$-pointed set $(A\backslash X,A\bullet x_0,\tilde*)$ to the 
 $\mathcal G_{\mathrm{inert}}$-pointed set $(A\backslash X,A\bullet x_0,*)$. 

 \begin{lem}\label{lem:extension:sd:product}
(a) There is an action $(g,P)\mapsto g\tilde* P$ of $\mathcal G_{\mathrm{inert}}\rtimes(\mathbb Z/2\mathbb Z)$ on the group 
$\mathrm{GL}_3\hat V$, uniquely defined by the conditions that it extends the action of $\mathcal G_{\mathrm{inert}}$ on the same group
from Lem. \ref{lem:13:11:8april}(b), and that $\overline 1\tilde* P=s_{(0,\infty)}(P)$ (see Lem. \ref{lem02:2506}(d)). 

(b) There is an action $(g,\tilde\rho)\mapsto g\tilde* \tilde\rho$ of $\mathcal G_{\mathrm{inert}}\rtimes(\mathbb Z/2\mathbb Z)$ on the set 
$\mathrm{Hom}_{\mathcal C\operatorname{-alg}}(\mathcal V[z]^\wedge,M_3\hat V)$, uniquely defined by the conditions that it extends the action 
of $\mathcal G_{\mathrm{inert}}$ on the same set from Lem. \ref{lem:13:12:10apr}(a), and that $\overline 1\tilde* \tilde\rho=
\mathrm{sw}\circ s_{(0,\infty)}^{\otimes2}\circ\tilde\rho\circ s_{(0,\infty)}^{\mathcal V[z]}$, where $\mathrm{sw}$ is as in 
Lem. \ref{lem:action:1342:0105}(b) and $s_{(0,\infty)}^{\mathcal V[z]}$ is as in Def. \ref{s:O:infty:V[z]}. 

(c) The tuple $(\mathrm{Hom}_{\mathcal C\operatorname{-alg}}(\mathcal V[z]^\wedge,M_3\hat V),\tilde\rho_{\mathrm{DT}},
\mathrm{GL}_3\hat V,\bullet,\tilde*)$ is a $\mathcal G_{\mathrm{inert}}\rtimes(\mathbb Z/2\mathbb Z)$-pointed set with group action, extending the 
$\mathcal G_{\mathrm{inert}}$-pointed set with group action from Lem. \ref{lem:13:12:10apr}(c). 

(d) The tuple $(\mathrm{GL}_3\hat V\backslash\mathrm{Hom}_{\mathcal C\operatorname{-alg}}(\mathcal V[z]^\wedge,M_3\hat V),
\mathrm{GL}_3\hat V\bullet\tilde\rho_{\mathrm{DT}},\tilde*)$ obtained by the tuple from (c) 
is a $\mathcal G_{\mathrm{inert}}\rtimes(\mathbb Z/2\mathbb Z)$-pointed set, extending the 
$\mathcal G_{\mathrm{inert}}$-pointed set $(\mathrm{GL}_3\hat V\backslash\mathrm{Hom}_{\mathcal C\operatorname{-alg}}(\mathcal V[z]^\wedge,
M_3\hat V),\mathrm{GL}_3\hat V\bullet\tilde\rho_{\mathrm{DT}},*)$ from Def. \ref{def:13:16:11apr}. 
 \end{lem}

\begin{proof}
(a) Since $\mathcal G_{\mathrm{inert}}\rtimes(\mathbb Z/2\mathbb Z)$ is generated by $\mathcal G_{\mathrm{inert}}$ and 
$\overline 1$, there is at most one action satisfying the conditions from (a). Let us now prove the existence of 
such an action. It follows from Lem. \ref{lem:action:1342:0105}(a) that the assignment $(g,P)\mapsto g\tilde*P:=\mathrm{aut}_g^{\mathcal V}(P)$
defines an action of $\mathcal G_{\mathrm{inert}}\rtimes(\mathbb Z/2\mathbb Z)$ on $\mathrm{GL}_3\hat V$. This action 
extends the action of $\mathcal G_{\mathrm{inert}}$ on $\mathrm{GL}_3\hat V$ from Lem. \ref{lem:13:11:8april}(b)
because of the extension property of the morphism $g\mapsto \mathrm{aut}_g^{\mathcal V}$ from Lem. \ref{lem:action:1342:0105}(a)
mentioned in this statement, and the action of $\overline 1$ is as indicated by the definition of $g\mapsto \mathrm{aut}_g^{\mathcal V}$ in 
Lem. \ref{lem:action:1342:0105}(a). 

(b) Since $\mathcal G_{\mathrm{inert}}\rtimes(\mathbb Z/2\mathbb Z)$ is generated by $\mathcal G_{\mathrm{inert}}$ and 
$\overline 1$, there is at most one action satisfying the conditions from (b). Let us now prove the existence of 
such an action. It follows from Lem. \ref{lem:action:1342:0105}(b,c) that the assignment 
$(g,\tilde\rho)\mapsto g\tilde*\tilde\rho:=\mathrm{aut}_g^V\circ\tilde\rho\circ(\mathrm{aut}_g^{\mathcal V[z]})^{-1}$
defines an action of $\mathcal G_{\mathrm{inert}}\rtimes(\mathbb Z/2\mathbb Z)$ on 
$\mathrm{Hom}_{\mathcal C\operatorname{-alg}}(\mathcal V[z]^\wedge,M_3\hat V)$. This action extends the action of 
$\mathcal G_{\mathrm{inert}}$ on the same set from Lem. \ref{lem:13:12:10apr}(a) by definition of the latter action. 
The action of $\overline 1$ is then as indicated, using the involutivity of $s_{(0,\infty)}^{\mathcal V[z]}$. 

(c) follows from (a) and (b), and (d) follows from (c) and from the compatibility 
between quotients en restrictions mentioned at the beginning of \S\ref{sect:amrzpswga}.  
\end{proof}

\subsection{The relation $\overline 1\in\mathrm{Stab}_{\mathcal G_{\mathrm{inert}}\rtimes(\mathbb Z/2\mathbb Z)}
(\mathrm{GL}_3\hat V\bullet\tilde\rho_{\mathrm{DT}})$}\label{sect:13:8}

\begin{lem}\label{lem:before}
Set $P_0:=\begin{pmatrix} 0&1&0\\    1&0&0\\ 0&0&1\end{pmatrix}
\in \mathrm{GL}_3\hat V$. 
One has (see Defs. \ref{def:5:2:1926}, \ref{defn:C:DT:R:DT} and Lem. \ref{lem:action:1342:0105}(b))
$$ 
\mathrm{aut}_{\overline 1}^V(\mathrm{col}_{\mathrm{DT}})=-P_0\cdot \mathrm{col}_{\mathrm{DT}},\quad 
\mathrm{aut}_{\overline 1}^V(\mathrm{row}_{\mathrm{DT}})=-\mathrm{row}_{\mathrm{DT}}\cdot P_0^{-1},  
$$
$$
\mathrm{aut}_{\overline 1}^V(C_{\mathrm{DT}})=P_0\cdot C_{\mathrm{DT}},\quad 
\mathrm{aut}_{\overline 1}^V(R_{\mathrm{DT}})=R_{\mathrm{DT}}\cdot P_0^{-1}. 
$$

\end{lem}

\begin{proof}
The automorphism $\mathrm{aut}_{\overline 1}^V$ of $\hat V$ is involutive and such that 
$e_1\leftrightarrow f_1$, $e_0\leftrightarrow f_\infty$ and  $f_0\leftrightarrow e_\infty$. Then 
$$
\mathrm{aut}_{\overline 1}^V(\mathrm{col}_{\mathrm{DT}})
=\mathrm{aut}_{\overline 1}^V(\begin{pmatrix}
    1\\-1\\0
\end{pmatrix})
=\begin{pmatrix}
    1\\-1\\0
\end{pmatrix}
=P_0\cdot \begin{pmatrix}
    -1\\1\\0
\end{pmatrix}=-P_0\cdot \mathrm{col}_{\mathrm{DT}}, 
$$
$$
\mathrm{aut}_{\overline 1}^V(\mathrm{row}_{\mathrm{DT}})
=\mathrm{aut}_{\overline 1}^V(\begin{pmatrix}
    e_1&-f_1&0
\end{pmatrix})
=\begin{pmatrix}
    f_1&-e_1&0
\end{pmatrix}
=\begin{pmatrix}
    -e_1&f_1&0
\end{pmatrix}\cdot P_0^{-1}=-\mathrm{row}_{\mathrm{DT}}\cdot P_0^{-1}.  
$$
$$
\mathrm{aut}_{\overline 1}^V(C_{\mathrm{DT}})=\mathrm{aut}_{\overline 1}^V(\begin{pmatrix}
 f_1   \\ e_1\\ -(e_0+f_\infty)
\end{pmatrix})=\begin{pmatrix}
 e_1   \\ f_1\\  -(e_0+f_\infty)
\end{pmatrix}=P_0\cdot \begin{pmatrix}
 f_1   \\ e_1\\  -(e_0+f_\infty)
\end{pmatrix}=P_0\cdot C_{\mathrm{DT}}, 
$$
$$
\mathrm{aut}_{\overline 1}^V(R_{\mathrm{DT}})
=\mathrm{aut}_{\overline 1}^V(\begin{pmatrix}
    0&0&1
\end{pmatrix})=
\begin{pmatrix}
    0&0&1
\end{pmatrix}=\begin{pmatrix}
    0&0&1
\end{pmatrix}\cdot P_0^{-1}=R_{\mathrm{DT}}\cdot P_0^{-1}. 
$$
\end{proof}

\begin{lem}\label{lem:s:belongs}
With $P_0$ as in Lem. \ref{lem:before}, one has 
$$
\overline 1*\tilde\rho_{\mathrm{DT}}=P_0\bullet\tilde\rho_{\mathrm{DT}}
$$
(equality in $\mathrm{Hom}_{\mathcal C\operatorname{-alg}}(\mathcal V[z]^\wedge,M_3\hat V)$), therefore 
$$
\overline 1\in 
\mathrm{Stab}_{\mathcal G_{\mathrm{inert}}\rtimes(\mathbb Z/2\mathbb Z)}
(\mathrm{GL}_3\hat V\bullet\tilde\rho_{\mathrm{DT}}). 
$$
\end{lem}

\begin{proof} By Lem. \ref{lem:action:1342:0105}, one has 
$\overline 1*\tilde\rho_{\mathrm{DT}}=\mathrm{aut}_{\overline 1}^V\circ\tilde\rho_{\mathrm{DT}}\circ
(\mathrm{aut}_{\overline 1}^{\mathcal V[z]})^{-1}=\mathrm{aut}_{\overline 1}^V\circ\tilde\rho_{\mathrm{DT}}\circ
s_{(0,\infty)}^{\mathcal V[z]}$,  therefore 
$$
\overline 1*\tilde\rho_{\mathrm{DT}} :  z\mapsto \mathrm{aut}_{\overline 1}^V(C_{\mathrm{DT}}R_{\mathrm{DT}}),\quad  e_1\mapsto \mathrm{aut}_{\overline 1}^V(\mathrm{col}_{\mathrm{DT}}\mathrm{row}_{\mathrm{DT}}),\quad e_0\mapsto 
\mathrm{aut}_{\overline 1}^V(-\rho_0-\mathrm{col}_{\mathrm{DT}}\mathrm{row}_{\mathrm{DT}}-C_{\mathrm{DT}}R_{\mathrm{DT}}). 
$$
while $P_0\bullet\tilde\rho_{\mathrm{DT}}=\mathrm{Ad}_{P_0}\circ \tilde\rho_{\mathrm{DT}}$ is given by 
$$
P_0\bullet\tilde\rho_{\mathrm{DT}} :  z\mapsto P_0C_{\mathrm{DT}}R_{\mathrm{DT}}P_0^{-1},\quad  e_1\mapsto P_0\mathrm{col}_{\mathrm{DT}}\mathrm{row}_{\mathrm{DT}}P_0^{-1},\quad e_0\mapsto 
P_0\rho_0P_0^{-1}. 
$$
It follows from Lem. \ref{lem:before} that the images of $z$ and $e_1$ by $\overline 1*\tilde\rho_{\mathrm{DT}}$ and 
$P_0\bullet\tilde\rho_{\mathrm{DT}}$ are equal. Moreover, 
$$
\rho_0+\mathrm{col}_{\mathrm{DT}}\mathrm{row}_{\mathrm{DT}}+C_{\mathrm{DT}}R_{\mathrm{DT}}
=\begin{pmatrix}
    e_0&0&0\\e_1&f_0&-e_1\\0&0&e_0
\end{pmatrix}+\begin{pmatrix}
    e_1&-f_1&0\\-e_1&f_1&0\\0&0&0
\end{pmatrix}
+\begin{pmatrix}
    0&0&f_1\\0&0&e_1\\0&0&-(e_0+f_\infty)
\end{pmatrix}=\begin{pmatrix}
   -e_\infty &-f_1&f_1\\0&-f_\infty&0\\0&0&-f_\infty
\end{pmatrix}
$$
therefore
$$
\mathrm{aut}_{\overline 1}^V(-\rho_0-\mathrm{col}_{\mathrm{DT}}\mathrm{row}_{\mathrm{DT}}-C_{\mathrm{DT}}R_{\mathrm{DT}})
=\begin{pmatrix}
   f_0 &e_1&-e_1\\0&e_0&0\\0&0&e_0
\end{pmatrix}=P_0\begin{pmatrix}
    e_0&0&0\\e_1&f_0&-e_1\\0&0&e_0
\end{pmatrix}P_0^{-1}=P_0\rho_0 P_0^{-1}
$$
which implies that the images of $e_0$ by $\overline 1*\tilde\rho_{\mathrm{DT}}$ and 
$P_0\bullet\tilde\rho_{\mathrm{DT}}$ are equal. The result then follows from the fact that 
$z,e_1,e_0$ generate the algebra $\mathcal V[z]^\wedge$. 
\end{proof}

\subsection{Stability of $\mathrm{Stab}_{\mathcal G}(\mathrm{GL}_3\hat V\bullet\rho_{\mathrm{DT}})$ 
under $\Theta$}\label{sect:13:9}

\begin{lem}\label{lem:aux}
Let $G_0$ be a group and $\theta$ be an involutive automorphism of $G_0$, and 
$G:=G_0\rtimes_\theta(\mathbb Z/2\mathbb Z)$ be the induced semidirect product group. 
Let $\overline 1\in G$ be the image of $\overline 1\in\mathbb Z/2\mathbb Z\hookrightarrow G$. Then 
$\overline 1\cdot g_0\cdot \overline 1=\theta(g_0)$ for any $g_0\in G_0$. 

If $H\subset G$ is a subgroup with $\overline 1\in H$, then $H_0:=H\cap G_0$ is stable under the automorphism 
$\theta$ of $G_0$. 
\end{lem}

\begin{proof}
The two subgroups of $G$ given by $G_0$ (by construction of the semidirect product) and $H$ (since $\overline 1\in H$)
are stable under the adjoint action $x\mapsto \overline 1\cdot x\cdot \overline 1^{-1}$, therefore so is their intersection 
$G_0\cap H$; the statement then follows from the coincidence of the restriction of 
$x\mapsto \overline 1\cdot x\cdot \overline 1^{-1}$ to $G_0$ with $\theta$. 
\end{proof}

\begin{lem}\label{lem:13:22:11apr}
    One has $\mathrm{Stab}_{\mathcal G_{\mathrm{inert}}}(\mathrm{GL}_3\hat V\bullet\tilde\rho_{\mathrm{DT}})
    =\mathrm{Stab}_{\mathcal G_{\mathrm{inert}}\rtimes(\mathbb Z/2\mathbb Z)}(\mathrm{GL}_3\hat V\bullet\tilde\rho_{\mathrm{DT}})
    \cap \mathcal G_{\mathrm{inert}}$
    (equality of subgroups of $\mathcal G_{\mathrm{inert}}\rtimes(\mathbb Z/2\mathbb Z)$). 
\end{lem}

\begin{proof}
This follows from the restriction statement of Lem. \ref{lem:extension:sd:product}(d).  
\end{proof}

\begin{lem}\label{lem:13:23:11apr}
    One has $\mathrm{Stab}_{\mathcal G}(\mathrm{GL}_3\hat V\bullet\rho_{\mathrm{DT}})
    =\mathrm{Stab}_{\mathcal G_{\mathrm{inert}}}(\mathrm{GL}_3\hat V\bullet\rho_{\mathrm{DT}})$
    (equality of subgroups of $\mathcal G$). 
\end{lem}

\begin{proof}
It follows from the restriction statement of Lem. \ref{lem:13:11:8april}(d) that 
$$
\mathrm{Stab}_{\mathcal G_{\mathrm{inert}}}(\mathrm{GL}_3\hat V\bullet\rho_{\mathrm{DT}})
=\mathrm{Stab}_{\mathcal G}(\mathrm{GL}_3\hat V\bullet\rho_{\mathrm{DT}})
\cap \mathcal G_{\mathrm{inert}}. 
$$
The result then follows from the inclusion \eqref{stab:in:intert} (see Thm. \ref{thm:13:22:0205}).  
\end{proof}

\begin{thm}\label{thm:13:36:17apr}
$\mathrm{Stab}_{\mathcal G}(\mathrm{GL}_3\hat V\bullet\rho_{\mathrm{DT}})$, which is a subgroup of 
$\mathcal G_{\mathrm{inert}}$ (see Thm. \ref{thm:13:22:0205}), is stable under the involution
$\Theta$ of this group (see Lem. \ref{lem02:2506}(d)).  
\end{thm}

\begin{proof}
Let us set $G_0:=\mathcal G_{\mathrm{inert}}$, let $\theta:=\Theta$. 
 The corresponding semidirect product (in the notation of Lem. \ref{lem:aux}) is 
$G:=\mathcal G_{\mathrm{inert}}\rtimes(\mathbb Z/2\mathbb Z)$.  
Let $H\subset G$ be the subgroup 
$$
\mathrm{Stab}_{\mathcal G_{\mathrm{inert}}\rtimes(\mathbb Z/2\mathbb Z)}
(\mathrm{GL}_3\hat V\bullet\tilde\rho_{\mathrm{DT}}). 
$$
By Lem. \ref{lem:s:belongs}, $\overline 1\in H$. 
Lem. \ref{lem:aux} then implies that $H_0:=G_0\cap H$ is stable under $\theta$. 

One has 
\begin{align*}
\mathrm{Stab}_{\mathcal G_{\mathrm{inert}}}(\mathrm{GL}_3\hat V\bullet\tilde\rho_{\mathrm{DT}})
=
\mathrm{Stab}_{\mathcal G_{\mathrm{inert}}\rtimes(\mathbb Z/2\mathbb Z)}
(\mathrm{GL}_3\hat V\bullet\tilde\rho_{\mathrm{DT}})
\cap \mathcal G_{\mathrm{inert}}
=H\cap G_0=H_0
\end{align*}
where the first equality (equality of subgroups of $\mathcal G_{\mathrm{inert}}\rtimes(\mathbb Z/2\mathbb Z)=G$) 
follows from Lem. \ref{lem:13:22:11apr}, and the next equalities 
follows from the already done identifications. It follows that 
$\mathrm{Stab}_{\mathcal G_{\mathrm{inert}}}(\mathrm{GL}_3\hat V\bullet\tilde\rho_{\mathrm{DT}})$
is stable under $\Theta$. By Cor. \ref{cor:13:29} and 
Lem. \ref{lem:13:23:11apr}, this group is equal to 
$\mathrm{Stab}_{\mathcal G}(\mathrm{GL}_3\hat V\bullet\rho_{\mathrm{DT}})$, which implies the statement. 
\end{proof}

\newpage

\part{Relationship of the double shuffle bitorsor with inertia-preserving bitorsors}\label{part 5}

The results established so far are concerned with the subgroup 
$\mathsf{DMR}_0(\mathbf k)\subset (\mathcal G,\circledast)$, 
and say in particular that this group is contained in the subgroup 
$(\mathcal G_{\mathrm{inert}},\circledast)$, and is invariant under 
the involution $\Theta$ of this subgroup (Cor. \ref{main:cor}(a) and (b)). 
However, the semidirect product $\mathsf{DMR}_0(\mathbf k)\rtimes\mathbf k^\times$ 
is known (\cite{EF3}) to be part of a richer structure, namely the bitorsor 
$$
(\mathsf{DMR}_0(\mathbf k)\rtimes\mathbf k^\times,\sqcup_{\mu\in\mathbf k^\times}
\mathsf{DMR}_\mu(\mathbf k),\mathsf{DMR}^{\mathrm B}(\mathbf k)),
$$ 
where $\sqcup_{\mu\in\mathbf k^\times}
\mathsf{DMR}_\mu(\mathbf k)$ is called the ``double shuffle scheme'' and 
$\mathsf{DMR}^{\mathrm B}(\mathbf k)$ the ``Betti version of the double shuffle group''.
The purpose of this part is to formulate and prove for these objects the analogues 
of the results obtained in Cor. \ref{main:cor}. This is done in \S\ref{last:sectioon} for the 
``double shuffle scheme'' and in \S\ref{section:LAST:SECTION} for the ``Betti 
version of the double shuffle group''. 

\section{Relationship of double shuffle schemes with inertia}\label{last:sectioon}

In this section, we fix $\mu\in\mathbf k$; then the pair 
$(\mathsf{DMR}_0(\mathbf k),\mathsf{DMR}_\mu(\mathbf k))$ is a left torsor. 
The purpose of this section is to formulate and prove the analogues of the results 
of Cor. \ref{main:cor} for $\mathsf{DMR}_\mu(\mathbf k)$. 
In \S\ref{sect:14:1}, we study the relationship of $(\mathcal G,\circledast)$ and 
$(\mathcal G_{\mathrm{inert}},\circledast)$ with tangential 
and inertia-preserving automorphisms. In \S\ref{subsect:14:1:2912}, 
we show using these results that $\mathcal G^\mu_{\mathrm{inert}}\subset \mathcal G$ is a left
torsor under the action of the subgroup $\mathcal G_{\mathrm{inert}}\subset 
(\mathcal G,\circledast)$. In \S\ref{sect:14:2:2912}, we extend the involution $\Theta$ of 
$\mathcal G_{\mathrm{inert}}$ to an involution $(\Theta,\Theta^\mu)$ 
of the torsor $(\mathcal G_{\mathrm{inert}},\mathcal G_{\mathrm{inert}}^\mu)$.
In \S\ref{sect:14:3:2912}, we prove the inclusion of the set $\mathsf M_\mu(\mathbf k)$ 
of associators with parameter $\mu$ in $\mathcal G_{\mathrm{inert}}^\mu$. 
In \S\ref{sect:rodsswi:2912}, we combine this result and Cor. \ref{main:cor}(a) 
to prove the inclusion $\mathsf M_\mu(\mathbf k)\subset \mathcal G^\mu_{\mathrm{inert}}$.

\subsection{Relationship of $(\mathcal G,\circledast)$ and $(\mathcal G_{\mathrm{inert}},\circledast)$ 
with tangential 
and inertia-preserving automorphisms}\label{sect:14:1}

\subsubsection{Tangential and inertia-preserving automorphisms}\label{rogagwotaipa} 
Denote by $\mathcal C(\mathfrak{lie}_{\{0,1\}}^\wedge)$ the quotient of $\mathfrak{lie}_{\{0,1\}}^\wedge$ by the conjugation action of 
$\mathrm{exp}(\mathfrak{lie}_{\{0,1\}}^\wedge)$. An automorphism $\alpha$ of $\mathfrak{lie}_{\{0,1\}}^\wedge$ induces a permutation of 
$\mathcal C(\mathfrak{lie}_{\{0,1\}}^\wedge)$, which will be denoted $\mathcal C(\alpha)$. Moreover, for $\alpha,\beta$ two automorphisms of 
 $\mathfrak{lie}_{\{0,1\}}^\wedge$, one has $\mathcal C(\alpha\circ\beta)=\mathcal C(\alpha)\circ\mathcal C(\beta)$.

\begin{defn}
(a) (see \cite{AT})  
$\mathrm{TAut}_{\{0,1\}}$ is the group of automorphisms $\alpha$ of $\mathfrak{lie}_{\{0,1\}}^\wedge$ 
  such that $\mathcal C(\alpha)(\mathcal C(e_0))=\mathcal C(e_0)$ and 
  $\mathcal C(\alpha)(\mathcal C(e_1))=\mathcal C(e_1)$. 
\\
(b) $\mathrm{IAut}_{\{0,1\}}$ is the subgroup of $\mathrm{TAut}_{\{0,1\}}$ of automorphisms $\alpha$ 
such that $\mathcal C(\alpha)(\mathcal C(e_\infty))=\mathcal C(e_\infty)$. 
\end{defn}
\begin{lem}\label{lem:section:S:2506}
The inner automorphism of $\mathrm{Aut}(\mathfrak{lie}_{\{0,1\}}^\wedge)$ induced by 
conjugation by $s_{(0,\infty)}$ induces an involution of the subgroup $\mathrm{IAut}_{\{0,1\}}$. 
\end{lem}
\begin{proof}
This follows from the fact that $\mathcal C(s_{(0,\infty)})$ leaves $\mathcal C(e_1)$ fixed and
permutes $\mathcal C(e_0)$ and $\mathcal C(e_\infty)$. 
\end{proof}

\subsubsection{Relationship of $(\mathcal G,\circledast)$ and $\mathcal G_{\mathrm{inert}}$ with 
 tangential and inertia-preserving automorphisms}

\begin{defn}
For $g\in \mathrm{exp}(\mathfrak{lie}_{\{0,1\}}^\wedge)$, one denotes
by $\mathrm{aut}_g^{\mathcal V}$ the 
automorphism of $\mathfrak{lie}_{\{0,1\}}^\wedge$ given by \eqref{def:aut:g:V:(1)}. 
\end{defn}

\begin{lem}\label{lem03:2506}
(a) The map $g\mapsto \mathrm{aut}_g^{\mathcal V}$ induces a group morphism 
$(\mathrm{exp}(\mathfrak{lie}_{\{0,1\}}^\wedge),\circledast)\to \mathrm{TAut}_{\{0,1\}}$. 
\\
(b) The subset of $\mathrm{TAut}_{\{0,1\}}$ of all automorphisms $\alpha$ such that 
$\alpha(e_1)=e_1$ and $\alpha(e_0)\equiv e_0$ mod degree $\geq 3$ is a subgroup 
$\mathrm{TAut}_{\{0,1\}}^0$, and the morphism from (a) induces a group isomorphism 
$(\mathcal G,\circledast)\stackrel{\sim}{\to} \mathrm{TAut}_{\{0,1\}}^0$. 
\end{lem}

\begin{proof}
(a) The fact that $g\mapsto \mathrm{aut}_g^{\mathcal V}$ is a group morphism 
$(\mathrm{exp}(\mathfrak{lie}_{\{0,1\}}^\wedge),\circledast)\to 
\mathrm{Aut}(\mathfrak{lie}_{\{0,1\}}^\wedge)$
follows from \cite{EF2}, (1.6.3); its image is obviously contained in $\mathrm{TAut}_{\{0,1\}}$.  
\\
(b) For $k\geq1$, let $F^k\mathfrak{lie}_{\{0,1\}}^\wedge$ be the ideal of 
$\mathfrak{lie}_{\{0,1\}}^\wedge$ which is the completed direct sum of all components of 
degree $\geq k$. Then $\mathrm{TAut}_{\{0,1\}}^0$ is the intersection of the stabilizer 
subgroup of $e_1$ and of the kernel of the morphism 
$\mathrm{Aut}(\mathfrak{lie}_{\{0,1\}}^\wedge)\to\mathrm{Aut}(\mathfrak{lie}_{\{0,1\}}^\wedge
/F^3\mathfrak{lie}_{\{0,1\}}^\wedge)$, which implies that it is a subgroup of  
$\mathrm{TAut}_{\{0,1\}}$.
\\
If $g\in \mathcal G$, then one checks that $\mathrm{aut}_g^{\mathcal V}\in \mathrm{TAut}_{\{0,1\}}^0$;
therefore $g\mapsto \mathrm{aut}_g^{\mathcal V}$ induces a group morphism 
$(\mathcal G,\circledast)\to \mathrm{TAut}_{\{0,1\}}^0$. 
The kernel of the morphism $g\mapsto \mathrm{aut}_g^{\mathcal V}$
is $\{\mathrm{exp}(\lambda e_0)|\lambda\in\mathbf k\}$, whose intersection with $\mathcal G$
is 1; this implies the injectivity of the morphism $(\mathcal G,\circledast)\to 
\mathrm{TAut}_{\{0,1\}}^0$. Let now $\alpha\in\mathrm{TAut}_{\{0,1\}}^0$, then there exists 
$\tilde g\in \mathrm{exp}(\mathfrak{lie}_{\{0,1\}}^\wedge)$ such that $\alpha(e_0)
=\tilde g\cdot e_0\cdot \tilde g^{-1}$. It follows from $\alpha(e_0)\equiv e_0$ modulo 
degree $\geq3$ that the expansion of $\mathrm{log}\tilde g$ is $ue_0+$ degree $\geq 2$, where 
$u\in\mathbf k$. Set then $g:=\tilde g\cdot\mathrm{exp}(-ue_0)$.  Then $g\in \mathcal G$
and $\alpha=\mathrm{aut}_g^{\mathcal V}$, which implies the surjectivity of 
$(\mathcal G,\circledast)\to \mathrm{TAut}_{\{0,1\}}^0$. 
\end{proof}

\begin{defn}
    $\mathrm{IAut}_{\{0,1\}}^0$ is the intersection $\mathrm{IAut}_{\{0,1\}}\cap 
    \mathrm{TAut}_{\{0,1\}}^0$ (intersection of subgroups of 
    $\mathrm{Aut}(\mathfrak{lie}_{\{0,1\}}^\wedge)$);  $\mathrm{IAut}_{\{0,1\}}^0$ is therefore the group of 
    automorphisms $\alpha$ of $\mathfrak{lie}_{\{0,1\}}^\wedge$, such that $\alpha(e_1)=e_1$, $\alpha(e_0)\equiv e_0$ 
    mod $F^3\mathfrak{lie}_{\{0,1\}}^\wedge$, and $\mathcal C(\alpha)$ leaves $\mathcal C(e_0)$ and 
    $\mathcal C(e_\infty)$ fixed. 
\end{defn}

\begin{lem}\label{lem04:2506}
(a) The assignment $g\mapsto\mathrm{aut}_g^{\mathcal V}$ induces a group 
isomorphism $\mathcal G_{\mathrm{inert}}\stackrel{\sim}{\to} \mathrm{IAut}_{\{0,1\}}^0$, 
such that for $g\in \mathcal G_{\mathrm{inert}}$
one has $\mathrm{aut}_g^{\mathcal V} : e_\infty\mapsto\mathrm{Ad}_{h_g}(e_\infty)$. 

(b) Conjugation by $s_{(0,\infty)}$ in $\mathrm{Aut}(\mathfrak{lie}_{\{0,1\}}^\wedge)$ induces an involutive automorphism $\mathrm{Ad}_{s_{(0,\infty)}}$ of $\mathrm{IAut}_{\{0,1\}}^0$. 

(c) The group isomorphism from (a) intertwines the involutive automorphisms $\Theta$ (see Lem. \ref{lem02:2506}(d)) of the source, 
and $\mathrm{Ad}_{s_{(0,\infty)}}$ of the target. 

(d) There is a unique group morphism 
$$
\mathcal G_{\mathrm{inert}}\rtimes (\mathbb Z/2\mathbb Z)\to \mathrm{IAut}_{\{0,1\}}\cdot\langle s_{(0,\infty)}\rangle, 
$$
(the target being equipped from the group structure arising from Lem. \ref{lem:section:S:2506}
) which extends the group morphism from 
(c) and is such that $(\mathbb Z/2\mathbb Z)\ni \overline 1\mapsto s_{(0,\infty)}$. This morphism is injective.  
\end{lem} 

\begin{proof} 
For $g\in \mathcal G_{\mathrm{inert}}$, one has 
$\mathrm{aut}_g^{\mathcal V}(e_\infty)=-\mathrm{aut}_g^{\mathcal V}(e_1)
-\mathrm{aut}_g^{\mathcal V}(e_0)=-e_1-ge_0g^{-1}=h_g\cdot e_\infty\cdot h_g^{-1}$, where the last equality follows from Lem. \ref{lem02:2506}(a). Therefore $\mathrm{aut}_g^{\mathcal V}\in 
\mathrm{IAut}_{\{0,1\}}$. It follows that $g\mapsto\mathrm{aut}_g^{\mathcal V}$ induces a group 
isomorphism $\mathcal G_{\mathrm{inert}}\to \mathrm{IAut}_{\{0,1\}}^0$. This morphism 
compatible with the isomorphism from Lem. \ref{lem03:2506}(b), therefore it 
is injective. Let us prove its surjectivity. Let $\alpha\in \mathrm{IAut}_{\{0,1\}}^0$.
By Lem. \ref{lem03:2506}(b), there exists $g\in\mathcal G$ such that 
$\alpha=\mathrm{aut}_g^{\mathcal V}$. Since $\alpha\in \mathrm{IAut}_{\{0,1\}}$, there exists 
$h\in \mathrm{exp}(\mathfrak{lie}_{\{0,1\}}^\wedge)$ such that 
$\alpha(e_\infty)=\mathrm{Ad}_h(e_\infty)$. Then $\mathrm{Ad}_h(e_\infty)=-\alpha(e_0)-\alpha(e_1)
=-\mathrm{Ad}_g(e_0)-e_1$, which since $g\in\mathcal G$ is equal to $e_\infty$ mod degree $\geq3$. 
This implies the existence of $\nu\in\mathbf k$ such that $\mathrm{log}h\equiv\nu e_\infty$
mod degree $\geq2$. Then $\tilde h:=h\mathrm{exp}(-\nu e_\infty)$ belongs to $\mathcal G$
and is such that $\mathrm{Ad}_g(e_0)+e_1+\mathrm{Ad}_{\tilde h}(e_\infty)=0$, which 
implies $g\in\mathcal G_{\mathrm{inert}}$. This ends to proof of (a). 

If $u\in\mathrm{IAut}_{\{0,1\}}$, 
then for some $a_0,a_1,a_\infty\in \mathrm{exp}(\mathfrak{lie}_{\{0,1\}}^\wedge)$, one has $u : e_0\mapsto \mathrm{Ad}_{a_0}(e_0)$, 
$e_1\mapsto \mathrm{Ad}_{a_1}(e_1)$, $e_\infty\mapsto \mathrm{Ad}_{a_\infty}(e_\infty)$. Then $s_{(0,\infty)}\circ u\circ s_{(0,\infty)}$
is given by $e_0\mapsto \mathrm{Ad}_{s_{(0,\infty)}(a_\infty)}(e_0)$, $e_1\mapsto \mathrm{Ad}_{s_{(0,\infty)}(a_1)}(e_1)$, 
$e_\infty\mapsto \mathrm{Ad}_{s_{(0,\infty)}(a_0)}(e_\infty)$, which implies 
$s_{(0,\infty)}\circ u\circ s_{(0,\infty)}\in \mathrm{IAut}_{\{0,1\}}$. This implies that 
conjugation by $s_{(0,\infty)}$
induces a group automorphism of $\mathrm{IAut}_{\{0,1\}}$; since the conjugation 
preserves $F^3\mathfrak{lie}_{\{0,1\}}^\wedge$, it restricts to an automorphism of
$\mathrm{IAut}_{\{0,1\}}^0$; since this conjugation is an involution of 
$\mathrm{TAut}_{\{0,1\}}$, this automorphism is involutive as well. This proves (b). 
Let $g\in\mathcal G_{\mathrm{inert}}$, then 
$\mathrm{aut}_{\Theta(g)}^{\mathcal V}(e_1)=s_{(0,\infty)}\circ\mathrm{aut}_g^{\mathcal V}\circ s_{(0,\infty)}(e_1)$ 
since the automorphisms $\mathrm{aut}_g^{\mathcal V}$, $\mathrm{aut}_{\Theta(g)}^{\mathcal V}$ and $s_{(0,\infty)}$
all leave $e_1$ fixed; and
\begin{align*}
&\mathrm{aut}_{\Theta(g)}^{\mathcal V}(e_0)=\mathrm{Ad}_{\Theta(g)}(e_1)
=\mathrm{Ad}_{s_{(0,\infty)}(h_g)}(e_0)
=s_{(0,\infty)}(\mathrm{Ad}_{h_g}(e_\infty))
=s_{(0,\infty)}\circ\mathrm{aut}_g^{\mathcal V}(e_\infty)
\\ & =s_{(0,\infty)}\circ\mathrm{aut}_g^{\mathcal V}\circ s_{(0,\infty)}(e_0)
\end{align*}
using in particular the fact that $\mathrm{aut}_g^{\mathcal V}$ is such that $e_\infty\mapsto 
\mathrm{Ad}_{h_g}(e_\infty)$. 
All this implies 
\begin{equation}\label{id:Theta:aut}
\mathrm{aut}_{\Theta(g)}^{\mathcal V}=s_{(0,\infty)}\circ\mathrm{aut}_g^{\mathcal V}\circ s_{(0,\infty)}, 
\end{equation}
proving (c). 
\\
The first statement of (d) follows from (c). It follows from (a) that the intersection of the kernel of 
$\mathcal G_{\mathrm{inert}}\rtimes_\Theta (\mathbb Z/2\mathbb Z)\to \mathrm{IAut}_{\{0,1\}}\cdot\langle s_{(0,\infty)}\rangle$
with $\mathcal G_{\mathrm{inert}}$ is trivial. The abelianization of the image in $\mathrm{Aut}(\mathfrak{lie}_{\{0,1\}}^\wedge)$ 
of any element of  $\mathcal G_{\mathrm{inert}}\cdot\overline 1$ is the automorphism $\overline e_1\mapsto \overline e_1$, 
$\overline e_0\mapsto \overline e_\infty$ of $\mathbf k\overline e_0\oplus\mathbf ke_1$, which implies that the intersection of 
the kernel of 
$\mathcal G_{\mathrm{inert}}\rtimes_\Theta (\mathbb Z/2\mathbb Z)\to \mathrm{IAut}_{\{0,1\}}\cdot\langle s_{(0,\infty)}\rangle$ with 
$\mathcal G_{\mathrm{inert}}\cdot\overline 1$ is empty; this proves the final statement. 
\end{proof}

\subsection{Torsor structure of $\mathcal G^\mu_{\mathrm{inert}}$}
\label{subsect:14:1:2912}

\subsubsection{The torsor $\mathcal G_{\mathrm{inert}}^\mu$ over 
$(\mathcal G_{\mathrm{inert}},\circledast)$}
 
Recall the notation $x*_\mu y=\mu^{-1}\mathrm{log}(e^{\mu x}e^{\mu y})$ for 
$\mu\in\mathbf k^\times$. 

\begin{lem}\label{rules:*} Let $\mu\in\mathbf k$. 

    (a) One has $(x*_\mu y)*_\mu z=x*_\mu(y*_\mu z)$ (identity in 
    $ (\mathfrak{lie}_{x,y,z}\otimes\mathbf k)^\wedge$) and $(-x)*_\mu x=0$ (identity in 
    $(\mathfrak{lie}_{x}\otimes\mathbf k)^\wedge=\mathbf kx$). 

    (b) One has $x*_\mu y*_\mu(-x)=\mathrm{Ad}_{e^{\mu x}}y$ (identity in 
    $ (\mathfrak{lie}_{x,y}\otimes\mathbf k)^\wedge$). 
\end{lem}

\begin{proof}
(a) Both sides of the first equality are checked to be equal to 
$\sum_{k\geq1}\mu^{k-1}\mathrm{cbh}_k(x,y,z)$, where $\sum_{k\geq1}\mathrm{cbh}_k(x,y,z)$
is the degree decomposition of $\mathrm{log}(e^xe^ye^z)$. The second equality is obvious. 

(b) One has $x*_\mu y*_\mu(-x)=\sum_{k\geq1}\mu^{k-1}\mathrm{cbh}_k(x,y,-x)$; on the other 
hand, $\sum_{k\geq1}\mathrm{cbh}_k(x,y,-x)=\sum_{k\geq 0}(\mathrm{ad}x)^k(y)/k!$, therefore 
$\sum_{k\geq1}\mu^{k-1}\mathrm{cbh}_k(x,y,-x)
=\sum_{k\geq 0}\mu^k(\mathrm{ad}x)^k(y)/k!=e^{\mu\mathrm{ad}x}(y)$, whch implies the result. 
\end{proof}

Note that $\mathcal G_{\mathrm{inert}}^0=\mathcal G_{\mathrm{inert}}$.  

\begin{lem}\label{def:Theta:mu:w:proofs} \Add{(see Lem. \ref{def:Theta:mu})}
Let $\mu\in\mathbf k$. 

    (a) If $g\in \mathcal G_{\mathrm{inert}}^\mu$, then there exists a unique 
    $h\in\mathrm{exp}(\mathfrak{lie}_{\{0,1\}}^\wedge)$ such that 
 $\mathrm{Ad}_ge_0*_\mu \mathrm{Ad}_he_\infty
 =e_0+e_\infty
 $ and $\mathrm{log}h\equiv (\mu/2)e_1$ mod $F^2\mathfrak{lie}_{\{0,1\}}^\wedge$; it 
    will be denoted $h_g$.

    (b) If $g\in \mathcal G_{\mathrm{inert}}^\mu$, then 
    $e^{-\mu e_1/2}(s_{(0,\infty)}(h_g))\in \mathcal G_{\mathrm{inert}}^\mu$. 

    (c) The map $\Theta^\mu : g\mapsto e^{-\mu e_1/2}(s_{(0,\infty)}(h_g))$ is an involution of the set $\mathcal G_{\mathrm{inert}}^\mu$. 
\end{lem}

\begin{proof}
  (a) 
Let $g\in\mathcal G_{\mathrm{inert}}^\mu$. Let us prove the existence of $h$. By assumption, 
there exists $\tilde h\in \mathrm{exp}(\mathfrak{lie}_{\{0,1\}}^\wedge)$, such that 
$\mathrm{Ad}_ge_0*_\mu \mathrm{Ad}_{\tilde h}e_\infty=e_0+e_\infty$. Let 
$h_1,h_\infty\in\mathbf k$ be such that $\mathrm{log}h\equiv h_1e_1+h_\infty e_\infty$
mod $F^2\mathfrak{lie}_{\{0,1\}}^\wedge$. Recall that 
$g\in \mathrm{exp}(F^2\mathfrak{lie}_{\{0,1\}}^\wedge)$. Then one has the expansion 
$\mathrm{Ad}_ge_0*_\mu \mathrm{Ad}_{\tilde h}e_\infty
\equiv e_0+e_\infty+(\mu/2)[e_0,e_\infty]+[h_1e_1+h_\infty e_\infty,e_\infty]$ mod 
$F^3\mathfrak{lie}_{\{0,1\}}^\wedge$. 
The equality 
$\mathrm{Ad}_ge_0*_\mu \mathrm{Ad}_{\tilde h}e_\infty=e_0+e_\infty$
then implies $h_1=\mu/2$. Then $h:=\tilde h e^{-h_\infty e_\infty}$ is such that 
$\mathrm{log}h\equiv h_1e_1=(\mu/2)e_1$ mod $F^2\mathfrak{lie}_{\{0,1\}}^\wedge$
and $\mathrm{Ad}_{h}e_\infty=\mathrm{Ad}_{\tilde h}e_\infty$, which implies 
$\mathrm{Ad}_ge_0*_\mu \mathrm{Ad}_{h}e_\infty=e_0+e_\infty$.

Let us now prove the uniqueness of the said $h$. If $h,h'\in\mathcal G$ satisfy the said condition,
  then applying $(-\mathrm{Ad}_ge_0)*_\mu-$ to the resulting equality 
  $\mathrm{Ad}_ge_0*_\mu \mathrm{Ad}_he_\infty=\mathrm{Ad}_ge_0*_\mu \mathrm{Ad}_{h'}e_\infty$ 
  and using Lem. \ref{rules:*}(a) yields
  $\mathrm{Ad}_{h}(e_\infty)=\mathrm{Ad}_{h'}(e_\infty)$, which implies 
  then $h'=h\mathrm{exp}(\alpha e_\infty)$ for some $\alpha\in\mathbf k$, which together with 
  the degree 1 conditions on $h,h'$ implies $h'=h$.

  (b) Let $g\in \mathcal G_{\mathrm{inert}}^\mu$. Then 
  $\mathrm{Ad}_ge_0 *_\mu \mathrm{Ad}_{h_g}e_\infty=e_0+e_\infty$, therefore 
$\mathrm{Ad}_ge_0 *_\mu \mathrm{Ad}_{h_g}e_\infty*_\mu e_1=0$, so 
$\mathrm{Ad}_{h_g}e_\infty*_\mu e_1*_\mu \mathrm{Ad}_ge_0  =0$, which 
upon applying $s_{(0,\infty)}$ gives  
$\mathrm{Ad}_{s_{(0,\infty)}(h_g)}e_0*_\mu e_1*_\mu \mathrm{Ad}_{s_{(0,\infty)}(g)}e_\infty  =0$. 
Applying $(-e_1/2)*-$ and $-*(-e_1/2)$, one then obtains  
$$
((-e_1/2)*_\mu \mathrm{Ad}_{s_{(0,\infty)}(h_g)}e_0*_\mu(e_1/2))
*_\mu((e_1/2)*_\mu \mathrm{Ad}_{s_{(0,\infty)}(g)}e_\infty*_\mu(-e_1/2))=-e_1
$$
therefore by Lem. \ref{rules:*}(b)
$$
\mathrm{Ad}_{\mathrm{exp}(-(\mu/2)e_1)}(\mathrm{Ad}_{s_{(0,\infty)}(h_g)}e_0)
*_\mu \mathrm{Ad}_{\mathrm{exp}((\mu/2)e_1)}(
\mathrm{Ad}_{s_{(0,\infty)}(g)}e_\infty
)=-e_1, 
$$
i.e.
$$
\mathrm{Ad}_{\mathrm{exp}(-(\mu/2)e_1)s_{(0,\infty)}(h_g)}e_0
*_\mu 
\mathrm{Ad}_{\mathrm{exp}((\mu/2)e_1)s_{(0,\infty)}(g)}e_\infty=-e_1.  
$$
It follows from (a) that $\mathrm{exp}(-(\mu/2)e_1)s_{(0,\infty)}(h_g)\in\mathcal G$, which 
together with this equality implies $\mathrm{exp}(-(\mu/2)e_1)s_{(0,\infty)}(h_g)\in\mathcal G^\mu$. 
Moreover, $\mathrm{log}(\mathrm{exp}((\mu/2)e_1)s_{(0,\infty)}(g))\equiv (\mu/2)e_1$ mod
$F^2\mathfrak{lie}_{\{0,1\}}^\wedge$, which implies 
$h_{e^{-\mu e_1/2}s_{(0,\infty)}(h_g)}=e^{\mu e_1/2}s_{(0,\infty)}(g)$. 

(c) For $g\in \mathcal G_{\mathrm{inert}}^\mu$, one then has 
\begin{align*}
&(\Theta^\mu)^2(g)=\Theta^\mu(e^{-\mu e_1/2}s_{(0,\infty)}(h_g))=
e^{-\mu e_1/2}s_{(0,\infty)}(h_{e^{-\mu e_1/2}s_{(0,\infty)}(h_g)})
\\&=e^{-\mu e_1/2}s_{(0,\infty)}(e^{\mu e_1/2}s_{(0,\infty)}(g))=g.
\end{align*}
\end{proof}

\subsubsection{Relationship of $\mathcal G_{\mathrm{inert}}^\mu$ with inertia-preserving isomorphisms}

\begin{defn}\label{def:14:1}

(a) Let $\mathrm{TAut}^\mu_{\{0,1\}}$ be the set of automorphisms $\beta$ of
$\mathfrak{lie}_{\{0,1\}}^\wedge$, such that  $\beta(e_0*_\mu
e_\infty)=e_0+e_\infty$,
the map $\mathcal C(\beta):\mathcal C(e_0)\mapsto \mathcal C(e_0)$, and
$\beta(e_0)\equiv e_0$ mod $F^3\mathfrak{lie}_{\{0,1\}}^\wedge$
($\mathcal C(\beta)$
being as in \S\ref{rogagwotaipa}).

(b) Let $\mathrm{IAut}^\mu_{\{0,1\}}$ be the subset of  $\mathrm{TAut}^\mu_{\{0,1\}}$
of automorphisms $\beta$ which satisfy $\mathcal C(\beta):\mathcal C(e_\infty)\mapsto 
\mathcal C(e_\infty)$. 
\end{defn}

\begin{defn}\label{def:aut:mu}
For $g\in\mathcal G$, let $\mathrm{aut}_g^{\mathcal V,\mu}$ be the automorphism of 
$\mathfrak{lie}_{\{0,1\}}^\wedge$ such that $e_0\mapsto \mathrm{Ad}_g(e_0)$, 
$e_0*_\mu e_\infty\mapsto e_0+e_\infty$.    
\end{defn}

\begin{lem}\label{lem:14:2}
 The map $g\mapsto \mathrm{aut}_g^{\mathcal V,\mu}$ 
sets up a bijection $\mathcal G\to \mathrm{TAut}^\mu_{\{0,1\}}$, 
which restricts to a bijection $\mathcal G_{\mathrm{inert}}^\mu\to\mathrm{IAut}^\mu_{\{0,1\}}$
(with the notation of Def. \ref{def:G:mu:inert}). 
\end{lem}

\begin{proof} Let us first show that $g\mapsto \mathrm{aut}_g^{\mathcal V,\mu}$ 
is a bijection $\mathcal G\to \mathrm{TAut}^\mu_{\{0,1\}}$. If $g\in \mathcal G$, then 
$\mathrm{aut}_g^{\mathcal V,\mu}$ is such that 
$e_0*_\mu e_\infty\mapsto e_0+e_\infty$, and 
$\mathrm{aut}_g^{\mathcal V,\mu}(e_0)\equiv e_0$ mod $F^3\mathfrak{lie}_{\{0,1\}}^\wedge$ 
 as $g\in\mathcal G$. Therefore $\mathrm{aut}_g^{\mathcal V,\mu}\in \mathrm{TAut}^\mu_{\{0,1\}}$. 
 Therefore  $g\mapsto \mathrm{aut}_g^{\mathcal V,\mu}$ 
induces a map $\mathcal G\to \mathrm{TAut}^\mu_{\{0,1\}}$.
Let us show that this map is injective. Let $g,g'\in \mathcal G_{\mathrm{inert}}^\mu$
 such that $\mathrm{aut}_g^{\mathcal V,\mu}=\mathrm{aut}_{g'}^{\mathcal V,\mu}$. Then 
 $\mathrm{Ad}_g(e_0)=\mathrm{Ad}_{g'}(e_0)$. Together with 
 $g,g'\in\mathrm{exp}(\mathfrak{lie}_{\{0,1\}}^\wedge)$, this implies the existence of 
 $\nu\in\mathbf k$, such that $g'=ge^{\nu e_0}$. The relations $g,g'\in\mathcal G$ then
 imply $\nu=0$, hence $g=g'$. Let us prove the surjectivity of the map 
 $\mathcal G\to \mathrm{TAut}^\mu_{\{0,1\}}$. Let $\beta\in \mathrm{TAut}^\mu_{\{0,1\}}$. 
 Then there exists 
 $u\in \mathrm{exp}(\mathfrak{lie}_{\{0,1\}}^\wedge)$ such that 
 $\beta(e_0)=\mathrm{Ad}_u(e_0)$. 
Define $u_0,u_1\in \mathbf k$ by $\mathrm{log}u\equiv u_0e_0+u_1e_1$ mod 
$F^2\mathfrak{lie}_{\{0,1\}}^\wedge$. Then $\beta(e_0) \equiv e_0+u_1[e_1,e_0]$
mod $F^3\mathfrak{lie}_{\{0,1\}}^\wedge$; the relation $\beta(e_0)\equiv e_0$
mod $F^3\mathfrak{lie}_{\{0,1\}}^\wedge$ then implies $u_1=0$, therefore 
$\mathrm{log}u\equiv u_0e_0$ mod $F^2\mathfrak{lie}_{\{0,1\}}^\wedge$. 
Set then $g:=ue^{-u_0e_0}$; one has $g\in\mathcal G$ and $\beta=\mathrm{aut}_g^{\mathcal V,\mu}$.

Let us show that $g\mapsto \mathrm{aut}_g^{\mathcal V,\mu}$
defines a map $\mathcal G_{\mathrm{inert}}^\mu\to\mathrm{IAut}^\mu_{\{0,1\}}$. 
Let $g\in \mathcal G_{\mathrm{inert}}^\mu$. Then 
$\mathrm{aut}_g^{\mathcal V,\mu}$ is such that 
$e_0*_\mu e_\infty\mapsto e_0+e_\infty$, and 
$\mathrm{aut}_g^{\mathcal V,\mu}(e_0)\equiv e_0$ mod $F^3\mathfrak{lie}_{\{0,1\}}^\wedge$ 
 as $g\in\mathcal G$. One has $\mathrm{aut}_g^{\mathcal V,\mu}(e_0)=\mathrm{Ad}_g(e_0)$
 which together with 
  $\mathrm{Ad}_ge_0*_\mu \mathrm{Ad}_{h_g}e_\infty=e_0+e_\infty$
 implies $\mathrm{aut}_g^{\mathcal V,\mu}(e_\infty)=\mathrm{Ad}_{h_g}(e_\infty)$, 
 therefore $\mathcal C(\mathrm{aut}_g^{\mathcal V,\mu})$ leaves 
 $\mathcal C(e_0)$ and $\mathcal C(e_\infty)$ invariant. 

 It remains to prove the surjectivity of $\mathcal G_{\mathrm{inert}}^\mu\to 
 \mathrm{IAut}^\mu_{\{0,1\}}$, $g\mapsto \mathrm{aut}_g^{\mathcal V,\mu}$. 
Let $\beta\in \mathrm{IAut}^\mu_{\{0,1\}}$. By the surjectivity of
the map $\mathcal G\to \mathrm{TAut}^\mu_{\{0,1\}}$, there exists $g\in \mathcal G$
such that $\beta= \mathrm{aut}_g^{\mathcal V,\mu}$, so in particular $\beta(e_0)=\mathrm{Ad}_g(e_0)$. 
Moreover, there exists $v\in \mathrm{exp}(\mathfrak{lie}_{\{0,1\}}^\wedge)$ such that 
$\beta(e_\infty)=\mathrm{Ad}_v(e_\infty)$; then $\beta(e_0*_\mu e_\infty)=e_0+e_\infty$ implies  
$e_0+e_\infty=\mathrm{Ad}_g(e_0)*_\mu\mathrm{Ad}_v(e_\infty)$, which implies 
$g\in \mathcal G_{\mathrm{inert}}^\mu$. 
\end{proof}

Recall that a left torsor is that data of a triple $(G,X,G\times X\to X)$ where $G$ is a group, 
$X$ is a nonempty set, $G\times X\to X$ is a left action of $G$ on $X$, which is transitive 
(i.e. $G\cdot x=X$ some, or equivalently any, $x\in X$) and such that $\mathrm{Stab}_G(x)=1$ for some, 
or equivalently any, $x\in X$. 

\begin{lem}\label{lem:torsor:iaut}
The map $\mathrm{TAut}^0_{\{0,1\}}\times \mathrm{TAut}^\mu_{\{0,1\}} \to 
\mathrm{TAut}^\mu_{\{0,1\}}$, $(\alpha,\beta)\mapsto \alpha\circ\beta$ defines 
a left torsor structure of the set $\mathrm{TAut}^\mu_{\{0,1\}}$ over the group 
$\mathrm{TAut}^0_{\{0,1\}}$, which restricts to a 
left torsor structure $\mathrm{IAut}^0_{\{0,1\}}\times \mathrm{IAut}^\mu_{\{0,1\}} \to 
\mathrm{IAut}^\mu_{\{0,1\}}$ of the set $\mathrm{IAut}^\mu_{\{0,1\}}$ over the group 
$\mathrm{IAut}^0_{\{0,1\}}$. 
\end{lem}

\begin{proof}
Let $\alpha \in \mathrm{TAut}^0_{\{0,1\}}$ and $\beta \in \mathrm{TAut}^\mu_{\{0,1\}}$. 
Then $\alpha \circ \beta(e_0*_\mu e_\infty)=e_0+e_\infty$ 
since $\alpha(e_1)=e_1$ and $\beta(e_0*_\mu e_\infty)=e_0+e_\infty$.
One has  
$\mathcal C(\alpha \circ \beta)(\mathcal C(e_0))=\mathcal C(e_0)$ since 
$\mathcal C(\alpha)$ and $\mathcal C(\beta)$ leave $\mathcal C(e_0)$ fixed. 
Finally $(\alpha \circ \beta)(e_0)\equiv \alpha(e_0) \equiv e_0$ mod 
$F^3\mathfrak{lie}_{\{0,1\}}^\wedge$, where the equivalences follow from $\beta(e_0) \equiv e_0$ 
and $\alpha(e_0) \equiv e_0$ mod $F^3\mathfrak{lie}_{\{0,1\}}^\wedge$ and $x \equiv y$ mod $F^3\mathfrak{lie}_{\{0,1\}}^\wedge$ implies $\alpha(x) \equiv \alpha(y)$ 
mod $F^3\mathfrak{lie}_{\{0,1\}}^\wedge$ for any $x,y \in \mathfrak{lie}_{\{0,1\}}^\wedge$. It follows that $\alpha \circ \beta \in \mathrm{IAut}_{\{0,1\}}^\mu$.

If $\beta,\beta' \in \mathrm{TAut}^\mu_{\{0,1\}}$, then the same arguments imply that $\beta' \circ \beta^{-1}$ belongs to $\mathrm{TAut}^0_{\{0,1\}}$. This implies 
the first torsor statement.

Let $\alpha \in \mathrm{IAut}_{\{0,1\}}^0$ and $\beta \in \mathrm{IAut}_{\{0,1\}}^\mu$. 
Then $\alpha \circ \beta\in \mathrm{TAut}^0_{\{0,1\}}$. Moreover 
$\mathcal C(\alpha \circ \beta)(\mathcal C(e_\infty))=\mathcal C(e_\infty)$ since 
$\mathcal C(\alpha)$ and $\mathcal C(\beta)$ leave $\mathcal C(e_\infty)$ fixed. 
Therefore $\alpha \circ \beta\in\mathrm{IAut}_{\{0,1\}}^\mu$. 

If $\beta,\beta' \in \mathrm{IAut}_{\{0,1\}}^\mu$, then the same arguments imply that $\beta' \circ 
\beta^{-1}$ belongs to $\mathrm{IAut}_{\{0,1\}}^0$, which implies 
the second torsor statement.
\end{proof}

Recall that a morphism of torsors from $(G,X,G\times X\to X)$ to $(H,Y,H\times Y\to Y)$
is the pair of a group morphism $\phi : G\to H$ and a map $f:X\to Y$, which is compatible with the 
actions.

\begin{lem}\label{pull:back:torsors}
If $(\phi,f) : (G,X,G\times X\to X)\to(H,Y,H\times Y\to Y)$ is a morphism of torsors 
and $(H',Y')$ is a subtorsor of $(H,Y,H\times Y\to Y)$
(i.e. $H'\subset H$ is a subgroup, $Y'\subset Y$ is a subset, and  $(H',Y',H'\times Y'\to Y')$), 
then $(\phi^{-1}(H'),f^{-1}(Y'))$ is a subtorsor of $(G,X,G\times X\to X)$. 
\end{lem}

\begin{proof}
    Obvious. 
\end{proof}

\begin{lem}\label{lem:14:6}
$\mathcal G_{\mathrm{inert}}^\mu$ is a left $\mathcal G_{\mathrm{inert}}$-torsor, 
 the action being induced by the left action of $(\mathcal G,\circledast)$ on itself; 
 $(g\mapsto\mathrm{aut}_g^{\mathcal V},g\mapsto\mathrm{aut}_g^{\mathcal V,\mu})$ 
 defines an isomorphism of torsors between $(\mathcal G_{\mathrm{inert}},\mathcal G_{\mathrm{inert}}^\mu)$ 
 and $(\mathrm{IAut}^0_{\{0,1\}},\mathrm{IAut}^\mu_{\{0,1\}})$. 
\end{lem}

\begin{proof}
One has for any $g,h\in\mathcal G$, 
$\mathrm{aut}_g^{\mathcal V} \circ \mathrm{aut}_h^{\mathcal V,\mu}(e_0*_\mu e_\infty)=\mathrm{aut}_g^{\mathcal V}(e_0+e_\infty)=e_0+e_\infty=\mathrm{aut}_{g\circledast h}^{\mathcal V,\mu}(e_0*_\mu e_\infty)$ and 
\begin{align*}
    & \mathrm{aut}_g^{\mathcal V} \circ \mathrm{aut}_h^{\mathcal V,\mu}(e_0)=\mathrm{aut}_g^{\mathcal V}
    (h(e_0,e_1)e_0h(e_0,e_1)^{-1})
    \\& =h(g(e_0,e_1)e_0g(e_0,e_1)^{-1},e_1)g(e_0,e_1)e_0g(e_0,e_1)^{-1}h(g(e_0,e_1)e_0
    g(e_0,e_1)^{-1},e_1)^{-1}=\mathrm{aut}_{g\circledast h}^{\mathcal V,\mu}(e_0), 
\end{align*}
which implies 
\begin{equation}\label{equality:aut:g:g:mu}
\forall g,h\in\mathcal G, \quad \mathrm{aut}_{g\circledast h}^{\mathcal V,\mu}=\mathrm{aut}_g^{\mathcal V}
\circ \mathrm{aut}_h^{\mathcal V,\mu}
\end{equation}
(equality in $\mathrm{TAut}^0_{\{0,1\}}$). 
Together with the first part of Lem. \ref{lem:torsor:iaut}, this implies that 
the pair of group and set morphisms 
$(g\mapsto \mathrm{aut}_g^{\mathcal V},h\mapsto \mathrm{aut}_h^{\mathcal V,\mu})$ 
defines an isomorphism of torsors between the triples $(\mathcal G,\mathcal G,
(g,h)\mapsto g\circledast h)$ and 
$(\mathrm{TAut}^0_{\{0,1\}},\mathrm{TAut}^\mu_{\{0,1\}},(\alpha,\beta)\mapsto \alpha\circ\beta)$.

 By Lems. \ref{lem04:2506}(a) and \ref{lem:14:2}, the preimages in 
 $\mathcal G$ of $\mathrm{IAut}_{\{0,1\}}^0$ and $\mathrm{IAut}_{\{0,1\}}^\mu$ by the 
 group isomorphism and bijection
 $g\mapsto \mathrm{aut}_g^{\mathcal V}$ and $h\mapsto \mathrm{aut}_h^{\mathcal V,\mu}$
 are respectively $\mathcal G_{\mathrm{inert}}$ and $\mathcal G_{\mathrm{inert}}^\mu$. 
 The first statement then follows from Lems. \ref{lem:torsor:iaut} and \ref{pull:back:torsors}; the 
 second statement is a direct consequence. 
\end{proof}

\subsection{Involution of the torsor $(\mathcal G_{\mathrm{inert}},\mathcal G_{\mathrm{inert}}^\mu)$}\label{sect:14:2:2912}

\begin{defn}
    Define $\sigma_\mu$ as the automorphism of 
    $\mathfrak{lie}_{\{0,1\}}^\wedge$ such that 
$$
e_0\mapsto \mathrm{Ad}_{(e^{\mu e_0}e^{\mu e_\infty})^{1/2}}(e_\infty), \quad 
e_\infty\mapsto \mathrm{Ad}_{(e^{\mu e_0}e^{\mu e_\infty})^{-1/2}}(e_0). 
$$    
\end{defn}

\begin{lem}\label{lem:14:7}
(a) $\sigma_\mu$ leaves $e_0*_\mu e_\infty$ fixed. (b) $\sigma_\mu^2=id$. 
\end{lem}

\begin{proof}
  (a) One has $\mathrm{log}(e^{\mu e_0}e^{\mu e_\infty})^{1/2}=(\mu/2)(e_0*_\mu e_\infty)$, 
therefore $\sigma_\mu(e_0)=((1/2)(e_0*_\mu e_\infty))*_\mu e_\infty *_\mu(-(1/2)(e_0*_\mu e_\infty))$ , 
and $\sigma_\mu(e_\infty)=(-(1/2)(e_0*_\mu e_\infty))*_\mu e_0 *_\mu((1/2)(e_0*_\mu e_\infty))$. 
Then 
\begin{align*}
&\sigma_\mu(e_0*_\mu e_\infty)
\\& =((1/2)(e_0*_\mu e_\infty))*_\mu e_\infty *_\mu(-(1/2)(e_0*_\mu e_\infty))
*_\mu (-(1/2)(e_0*_\mu e_\infty))*_\mu e_0 *_\mu((1/2)(e_0*_\mu e_\infty))
\\ & =((1/2)(e_0*_\mu e_\infty))*_\mu e_\infty *_\mu(-(e_0*_\mu e_\infty))*_\mu e_0 *_\mu((1/2)(e_0*_\mu e_\infty))
\\& =((1/2)(e_0*_\mu e_\infty))*_\mu e_\infty *_\mu(-e_\infty)*_\mu(-e_0)*_\mu e_0 *_\mu((1/2)(e_0*_\mu e_\infty))
\\& =((1/2)(e_0*_\mu e_\infty))*_\mu((1/2)(e_0*_\mu e_\infty))=e_0*_\mu e_\infty. 
\end{align*}
 
 (b) Since $(e^{\mu e_0}e^{\mu e_\infty})^{1/2}=e^{(\mu/2)(e_0*_\mu e_\infty)}$, it follows from (a) 
 that $\sigma_\mu$ leaves $(e^{\mu e_0}e^{\mu e_\infty})^{1/2}$ fixed. Then 
 one computes 
 $$
 \sigma_\mu^2(e_0)=\sigma_\mu( \mathrm{Ad}_{(e^{\mu e_0}e^{\mu e_\infty})^{1/2}}(e_\infty))=
 \mathrm{Ad}_{(e^{\mu e_0}e^{\mu e_\infty})^{1/2}}(\sigma_\mu(e_\infty))
 =\mathrm{Ad}_{(e^{\mu e_0}e^{\mu e_\infty})^{1/2}}\circ  \mathrm{Ad}_{(e^{\mu e_0}e^{\mu e_\infty})^{-1/2}}(e_0)
 =e_0
 $$
 and similarly 
  $$
 \sigma_\mu^2(e_\infty)=\sigma_\mu( \mathrm{Ad}_{(e^{\mu e_0}e^{\mu e_\infty})^{-1/2}}(e_0))=
 \mathrm{Ad}_{(e^{\mu e_0}e^{\mu e_\infty})^{-1/2}}(\sigma_\mu(e_0))
 =\mathrm{Ad}_{(e^{\mu e_0}e^{\mu e_\infty})^{-1/2}}\circ  \mathrm{Ad}_{(e^{\mu e_0}e^{\mu e_\infty})^{1/2}}(e_\infty)
 =e_\infty, 
 $$
therefore $ \sigma_\mu^2=id$.   
\end{proof}

\begin{lem}
(a) The assignment  $\ell(s_{(0,\infty)})
r(\sigma_\mu) : \beta\mapsto s_{(0,\infty)}\circ \beta\circ\sigma_\mu$ defines an involution of the set $\mathrm{IAut}^\mu_{\{0,1\}}$. 

(b) The bijection $\mathcal G^\mu_{\mathrm{inert}}\to\mathrm{IAut}^\mu_{\{0,1\}}$, 
$g\mapsto \mathrm{aut}_g^{\mathcal V,\mu}$ intertwines the 
involutions $\Theta^\mu$ (see Lem. \ref{def:Theta:mu}) and 
$\ell(s_{(0,\infty)})r(\sigma_\mu)$ of its source and target, namely 
\begin{equation}\label{id:Theta:s:sigma}
\forall g\in\mathcal G^\mu_{\mathrm{inert}},\quad
s_{(0,\infty)}\circ \mathrm{aut}_g^{\mathcal V,\mu}\circ\sigma_\mu=\mathrm{aut}_{\Theta^\mu(g)}^{\mathcal V,\mu}. 
\end{equation} 
\end{lem}

\begin{proof}
(a) Since $s_{(0,\infty)}$ and $\sigma_\mu$ are involutions in 
$\mathrm{Aut}(\mathfrak{lie}_{\{0,1\}}^\wedge)$, the assignment 
$\beta\mapsto s_{(0,\infty)}\circ \beta\circ\sigma_\mu$ is an 
involutive self-map $\ell(s_{(0,\infty)})
r(\sigma_\mu)$ of $\mathrm{Aut}(\mathfrak{lie}_{\{0,1\}}^\wedge)$.  
Let us now show that $\ell(s_{(0,\infty)})
r(\sigma_\mu)$ maps the subset $\mathrm{IAut}^\mu_{\{0,1\}}$
to itself. 

Let $\beta\in\mathrm{IAut}^\mu_{\{0,1\}}$. Then $\mathcal C(\beta)$ leaves $\mathcal C(e_0)$ 
and $\mathcal C(e_\infty)$ fixed. Since  $\mathcal C(\sigma_\mu)$ and $\mathcal C(s_{(0,\infty)})$ 
both permute these two classes, $\mathcal C(s_{(0,\infty)})\mathcal C(\sigma_\mu)
\mathcal C(\sigma_\mu)$ leaves each of them fixed. Moreover, $\sigma_\mu$ leaves 
$e_0*_\mu e_\infty$ fixed while $s_{(0,\infty)}$ leaves $e_0+e_\infty$ fixed, which together with
$\beta(e_0*_\mu e_\infty)=e_0+e_\infty$ implies 
$s_{(0,\infty)}\beta\sigma_\mu(e_0*_\mu e_\infty)=e_0+e_\infty$. 
Moreover, $e_\infty*_\mu e_0\equiv e_\infty+e_0+(\mu/2)[e_\infty,e_0]$ mod 
$F^3\mathfrak{lie}_{\{0,1\}}^\wedge$ and $\beta(e_0°\equiv e_0$ mod 
$F^3\mathfrak{lie}_{\{0,1\}}^\wedge$ implies $\beta(e_\infty)\equiv e_\infty-(\mu/2)[e_0,e_\infty]$. 
This implies the second equivalence in 
$$
s_{(0,\infty)}\beta\sigma_\mu(e_0)\equiv
s_{(0,\infty)}\beta(e_\infty+(\mu/2)[e_0,e_\infty])\equiv s_{(0,\infty)}(e_\infty)=e_0
\text{ mod }F^3\mathfrak{lie}_{\{0,1\}}^\wedge. 
$$
Therefore $s_{(0,\infty)}\circ \beta\circ\sigma_\mu\in\mathrm{IAut}^\mu_{\{0,1\}}$. 


(b) Let $g\in \mathcal G^\mu_{\mathrm{inert}}$. Then 
$$
s_{(0,\infty)}\circ \mathrm{aut}_g^{\mathcal V,\mu}\circ\sigma_\mu(e_0*_\mu e_\infty)
=s_{(0,\infty)}\circ \mathrm{aut}_g^{\mathcal V,\mu}(e_0*_\mu e_\infty)
=s_{(0,\infty)}(e_0+e_\infty)
=e_0+e_\infty
=\mathrm{aut}_{\Theta^\mu(g)}^{\mathcal V,\mu}(e_0*_\mu e_\infty)
$$
where the first equality follows from Lem. \ref{lem:14:7}(a), and the second and last equalities
follows from Def. \ref{def:aut:mu}. 

Moreover, 
\begin{align*}
&s_{(0,\infty)}\circ \mathrm{aut}_g^{\mathcal V,\mu}\circ\sigma_\mu(e_0)
=s_{(0,\infty)}\circ \mathrm{aut}_g^{\mathcal V,\mu}(\mathrm{Ad}_{(e^{\mu e_0}e^{\mu e_\infty})^{1/2}}(e_\infty))
\\&=s_{(0,\infty)}\circ \mathrm{Ad}_{\mathrm{aut}_g^{\mathcal V,\mu}((e^{\mu e_0}e^{\mu e_\infty})^{1/2})}
(\mathrm{aut}_g^{\mathcal V,\mu}(e_\infty))
=s_{(0,\infty)}\circ \mathrm{Ad}_{e^{(\mu/2)(e_0+e_\infty)}}
(\mathrm{Ad}_{h_g(e_0,e_1)}(e_\infty))
\\&=s_{(0,\infty)}(\mathrm{Ad}_{e^{(\mu/2)e_1}h_g(e_0,e_1)}
(e_\infty))=\mathrm{Ad}_{e^{(\mu/2)e_1}h_g(e_\infty,e_1)}
(e_0)=\mathrm{aut}_{\Theta^\mu(g)}^{\mathcal V,\mu}(e_0). 
\end{align*}
All this implies \eqref{id:Theta:s:sigma}. 
\end{proof}

\begin{lem}\label{lem:14:10}
 Let $\mu\in\mathbf k$. One has $\Theta^\mu(g\circledast g')=\Theta(g)\circledast \Theta^\mu(g')$
for any $g\in \mathcal G_{\mathrm{inert}}$ and $g'\in \mathcal G_{\mathrm{inert}}^\mu$; 
therefore $(\Theta,\Theta^\mu)$ is an involution of the torsor 
$(\mathcal G_{\mathrm{inert}},\mathcal G_{\mathrm{inert}}^{\mu})$. 
\end{lem}

\begin{proof}
One has 
\begin{align*}
&\mathrm{aut}_{\Theta^\mu(g\circledast g')}^{\mathcal V,\mu}
=s_{(0,\infty)}\circ \mathrm{aut}_{g\circledast g'}^{\mathcal V,\mu}\circ\sigma_\mu
=s_{(0,\infty)}\circ \mathrm{aut}_{g}^{\mathcal V}\circ 
\mathrm{aut}_{g'}^{\mathcal V,\mu}\circ\sigma_\mu
\\&=(s_{(0,\infty)}\circ \mathrm{aut}_{g}^{\mathcal V}\circ s_{(0,\infty)}^{-1})
\circ (s_{(0,\infty)}\circ 
\mathrm{aut}_{g'}^{\mathcal V,\mu}\circ\sigma_\mu)
=\mathrm{aut}_{\Theta(g)}^{\mathcal V}\circ\mathrm{aut}_{\Theta^\mu(g')}^{\mathcal V,\mu}
=\mathrm{aut}_{\Theta(g)\circledast\Theta^\mu(g')}^{\mathcal V,\mu}, 
\end{align*}
(equalities in $\mathrm{IAut}^0_{\{0,1\}}$)
where the first equality follows from \eqref{id:Theta:s:sigma}, 
the second and equalities follow from \eqref{equality:aut:g:g:mu}, and the fourth equality 
follows from \eqref{id:Theta:aut} and \eqref{id:Theta:s:sigma}. The result then follows from 
Lem. \ref{lem:14:2}.  
\end{proof}

\begin{rem}
The identity $s_{(0,\infty)}\circ \mathrm{aut}_{g\circledast g'}^{\mathcal V,\mu}\circ\sigma_\mu
=(s_{(0,\infty)}\circ \mathrm{aut}_{g}^{\mathcal V}\circ s_{(0,\infty)}^{-1})
\circ (s_{(0,\infty)}\circ \mathrm{aut}_{g'}^{\mathcal V,\mu}\circ\sigma_\mu)$ used in the above proof
expresses the fact that $(\mathrm{Ad}_{s_{(0,\infty)}},\ell(s_{(0,\infty)})r(\sigma_\mu))$ 
is an involution of the torsor $(\mathrm{IAut}^0_{\{0,1\}},\mathrm{IAut}^\mu_{\{0,1\}})$,  
which maps isomorphically to the involution $(\Theta,\Theta^\mu)$ 
of the torsor $(\mathcal G_{\mathrm{inert}},\mathcal G_{\mathrm{inert}}^{\mu})$ 
by the torsor isomorphism from Lem. \ref{lem:14:6}. 
\end{rem}

\subsection{Associators and $\mathcal G^\mu_{\mathrm{inert}}$}\label{sect:14:3:2912}

Let $\mathbf k$ be a commutative $\mathbb Q$-algebra. 
\begin{defn}
    For $\mu\in\mathbf k$, the set of $\mathbf k$-associators with parameter $\mu$ is the subset 
    $\mathsf M_\mu(\mathbf k)\subset \mathcal G$ defined by the duality, hexagon and pentagon
    conditions 
    $$
    \varphi(e_0,e_1)\varphi(e_1,e_0)=1,\quad e^{\mu e_0/2}\varphi(e_\infty,e_0)e^{\mu e_\infty/2}
    \varphi(e_1,e_\infty)e^{\mu e_1/2}\varphi(e_0,e_1)=1 
    $$
(in $U(\mathfrak{lie}_{\{0,1\}})^\wedge$)
    $$
    \varphi(e_{12},e_{23})\varphi(e_{34},e_{45})\varphi(e_{51},e_{12})\varphi(e_{23},e_{34})
    \varphi(e_{45},e_{51})=0
    $$
    in $U(\mathfrak p_5)^\wedge$, where $\mathfrak p_5$ is the Lie algebra with generators 
    $e_{ij}$ ($i\neq j\in\{1,\ldots,5\}$) and relations $e_{ji}=e_{ij}$ for $i\neq j$, $[e_{ij},e_{kl}]=0$, for distinct $i,j,k,l$, and $\sum_{j|j\neq i}e_{ij}=0$ for any $i$.  
\end{defn}

\begin{thm} \label{thm:dr:fur}
(a) (cf. \cite{Dr}, Prop. 5.3) $\mathsf M_1(\mathbf k)$ is nonempty, and 
$\varphi\mapsto (\mu\bullet\varphi)(e_0,e_1):=\varphi(\mu e_0,\mu e_1)$ defines a map 
$\mathsf M_1(\mathbf k)\to \mathsf M_\mu(\mathbf k)$. 
  
(b)  (cf. \cite{Fu:assandDS,EF2}) One has the inclusion  $\mathsf M_\mu(\mathbf k)\subset 
\mathsf{DMR}_\mu(\mathbf k)$.    
\end{thm}

\begin{lem}\label{lem:14:13}
 Let $\mu\in\mathbf k$. 
 
 (a)  One has $\mu\bullet\mathsf M_1(\mathbf k)
 \subset \mathcal G_{\mathrm{inert}}^\mu$ (inclusion of sets).  

  (b) $\mu\bullet\mathsf M_1(\mathbf k)$ is contained in the subset of 
  $\mathcal G_{\mathrm{inert}}^\mu$ of fixed points of its involution 
  $\Theta^\mu$. 
\end{lem}

\begin{proof} Let $\varphi\in \mathsf M_1(\mathbf k)$.  
 One has 
 $\mathrm{Ad}_{\varphi(e_0,e_1)}(e^{e_0})
 \mathrm{Ad}_{e^{e_1/2}\varphi(e_\infty,e_1)}(e^{e_\infty})=e^{e_0+e_\infty}$ 
 (see e.g. \cite{AET}, \S5.2) i.e. 
  $\mathrm{Ad}_{\varphi(e_0,e_1)}(e_0) *_1
 \mathrm{Ad}_{e^{e_1/2}\varphi(e_\infty,e_1)}(e_\infty)=e_0+e_\infty$. Therefore 
   $$
   \mathrm{Ad}_{(\mu\bullet\varphi)(e_0,e_1)}(e_0) *_\mu
 \mathrm{Ad}_{e^{\mu e_1/2}(\mu\bullet\varphi)(e_\infty,e_1)}(e_\infty)=e_0+e_\infty.
 $$
One has $\mu\bullet\varphi\in\mathcal G$ and 
$\mathrm{log}(e^{\mu e_1/2}(\mu\bullet\varphi)(e_\infty,e_1))\equiv (\mu/2)e_1$ mod
$F^2\mathfrak{lie}_{\{0,1\}}^\wedge$, which together with this equality 
implies $\mu\bullet\varphi\in\mathcal G_{\mathrm{inert}}^\mu$
 and $h_{(\mu\bullet\varphi)(e_0,e_1)}=e^{\mu e_1/2}(\mu\bullet\varphi)(e_\infty,e_1)$. This proves (a). 
One then has $\Theta^\mu((\mu\bullet\varphi)(e_0,e_1))
=e^{-\mu e_1/2}s_{(0,\infty)}(h_{(\mu\bullet\varphi)(e_0,e_1)})
 =e^{-\mu e_1/2}s_{(0,\infty)}(e^{\mu e_1/2}(\mu\bullet\varphi)(e_\infty,e_1))
 =(\mu\bullet\varphi)(e_0,e_1)$, which proves (b). 
\end{proof}

\subsection{Relationship of double shuffle schemes with inertia}\label{sect:rodsswi:2912}

\begin{thm}\label{thm:DMRmu} (see Thm. \ref{thm:0:31:25jan})
  Let $\mu\in\mathbf k$. Then:
  
  (a) the inclusion $\mathsf{DMR}_\mu(\mathbf k)\subset\mathcal G^\mu_{\mathrm{inert}}$ holds
  (inclusion of sets); 
  
  (b) the subset $\mathsf{DMR}_\mu(\mathbf k)$ of $\mathcal G^\mu_{\mathrm{inert}}$ is stable 
  under the involution $\Theta^\mu$ of this set. 
\end{thm}

\begin{proof}
 (a)  Combining Thm. \ref{thm:dr:fur}(a) and (b), one obtains $\mu \bullet \mathsf M_1(\mathbf k) 
 \subset \mathsf{DMR}_\mu(\mathbf k)$. It then follows from Thm. \ref{thm:racinet}(b) that 
 \begin{equation}\label{TOTO:24nov}
     \mathsf{DMR}_\mu(\mathbf k)=\mathsf{DMR}_0(\mathbf k) \circledast (\mu \bullet 
 \mathsf M_1(\mathbf k)). 
 \end{equation} 
 The result then follows from the combination of this equality and the 
 inclusions $\mathsf{DMR}_0(\mathbf k) \subset \mathcal G_{\mathrm{inert}}$ (Thm. \ref{thm:014}(a)), 
 $\mu \bullet \mathsf M_1(\mathbf k) \subset \mathcal G_{\mathrm{inert}}^\mu$ (Lem. \ref{lem:14:13}(a)) 
 and $\mathcal G_{\mathrm{inert}} \circledast \mathcal G_{\mathrm{inert}}^\mu \subset 
 \mathcal G_{\mathrm{inert}}^\mu$ (Lem. \ref{lem:14:6}).
 
(b)  Let $g\in \mathsf{DMR}_\mu(\mathbf k)$. By \eqref{TOTO:24nov}, there exists 
$g'\in \mathsf{DMR}_0(\mathbf k)$ and $g''\in \mu \bullet \mathsf M_1(\mathbf k)$
such that $g=g'\circledast g''$. Then
$$
\Theta^\mu(g)=\Theta^\mu(g'\circledast g'')
=\Theta(g')\circledast\Theta^\mu(g'')=\Theta(g')\circledast g''
\in \mathsf{DMR}_0(\mathbf k)\circledast \mathsf{DMR}_\mu(\mathbf k)=\mathsf{DMR}_\mu(\mathbf k),
$$
where the second equality follows from 
$\mathsf{DMR}_0(\mathbf k)\subset \mathcal G_{\mathrm{inert}}$ (Thm. \ref{thm:014}(a)), 
$\mu \bullet \mathsf M_1(\mathbf k)\subset \mathcal G_{\mathrm{inert}}^\mu$ (Lem. \ref{lem:14:13}) and 
Lem. \ref{lem:14:10}, the third equality follows from Lem. \ref{lem:14:13}(b), the relation 
`$\in$' follows from the combination of the stability of $\mathsf{DMR}_0(\mathbf k)$ 
under $\Theta$ (see Thm. \ref{thm:014}(b))
and the inclusion $\mu \bullet \mathsf M_1(\mathbf k) \subset \mathcal G_{\mathrm{inert}}^\mu$ 
(Lem. \ref{lem:14:13}(a)), and the last equality follows from Thm. \ref{thm:racinet}(b).  
\end{proof}

\section{The Betti double shuffle group and inertia}\label{section:LAST:SECTION}

The purpose of this part is to formulate and prove the analogue of the results of 
Cor. \ref{main:cor} for the ``Betti double shuffle group'' $\mathsf{DMR}^{\mathrm B}(\mathbf k)$. 
This relies on bitorsor results, which are established first. In 
§\ref{sect:15:1:2912}, we recall the bitorsor $(\mathcal G \rtimes \mathbf k^\times,
\mathcal G\times\mathbf k^\times,\mathcal G^{\mathrm B}\rtimes \mathbf k^\times)$ of which
$(\mathsf{DMR}_0(\mathbf k)\rtimes\mathbf k^\times,\sqcup_{\mu\in\mathbf k^\times}\mathsf{DMR}_\mu(\mathbf k),\mathsf{DMR}^{\mathrm B}(\mathbf k))$ is a subbitorsor. 
In §\ref{sect:15:2:2912}, we construct a subgroup $(\mathcal G^{\mathrm B} \rtimes 
\mathbf k^\times)_{\mathrm{inert}}$ of $\mathcal G^{\mathrm B} \rtimes \mathbf k^\times$ and 
its involution $\Theta^{\mathrm B}$. In §\ref{subsect:15:3:rt:2912}, we show that 
$(\mathcal G_{\mathrm{inert}} \rtimes 
\mathbf k^\times,\sqcup_{\mu \in \mathbf k^\times} \mathcal G^\mu_{\mathrm{inert}},
(\mathcal G^{\mathrm B} \rtimes \mathbf k^\times)_{\mathrm{inert}})$ is a subbitorsor of 
$(\mathcal G \rtimes \mathbf k^\times,\sqcup_{\mu \in \mathbf k}^\times \mathcal G,
\mathcal G^{\mathrm B} \rtimes \mathbf k^\times)$, and of which we show that 
$(\Theta \rtimes id,\sqcup_{\mu\in\mathbf k^\times}\Theta^\mu,\Theta^{\mathrm B})$ is 
an involution. In §\ref{sect:15:4:tiXXXaisuXXX:2912}, we combine the results 
obtained in §\ref{last:sectioon} with Lem. 1.13 in \cite{EF3} (see Lem. \ref{lem:ef3})
to obtain the announced results on $\mathsf{DMR}^{\mathrm B}(\mathbf k)$: its inclusion 
in $(\mathcal G^{\mathrm B} \rtimes id)_{\mathrm{inert}}$ and its stability under 
$\Theta^{\mathrm B}$.

\subsection{Reminders from \cite{EF3}}\label{sect:15:1:2912}

For $\Gamma$ a discrete group, set $\Gamma(\mathbf k):=\mathcal G((\mathbf k\Gamma)^\wedge)$, 
where $(\mathbf k\Gamma)^\wedge$ is the topological Hopf algebra obtained from 
$\mathbf k\Gamma$ by completion for the topology of powers of the augmentation ideal, 
and $\mathcal G$ means the set of group-like elements. A group morphism 
$\Gamma\to\Gamma'$ gives rise to a morphism $\Gamma(\mathbf k)\to\Gamma'(\mathbf k)$. 
One has $\mathbb Z^2(\mathbf k)=\mathbf k^2$.

For $\Gamma$ a discrete group, the assignment $\mathbf k\mapsto\Gamma(\mathbf k)$
is a prounipotent $\mathbb Q$-group scheme, with Lie algebra 
$\mathrm{Lie}\Gamma:=\mathcal P((\mathbb Q\Gamma)^\wedge)$, where 
$\mathcal P$ stands for primitive elements. 

\begin{defn}
    Define $\mathcal G^{\mathrm B}:=\mathrm{ker}(F_2(\mathbf k)\to\mathbb Z^2(\mathbf k))$, 
    where the morphism $F_2\to\mathbb Z^2$ is the abelianization morphism of the free group $F_2$ 
    with generators $X_0,X_1$.
\end{defn}

\begin{lem} \label{lem:15:2:14dec}(\cite{EF3}, §2.1.3)
(a) A group structure $\circledast$ is defined on $\mathcal G^{\mathrm B}$ by 
$$
g(X_0,X_1)\circledast h(X_0,X_1)
:=h(g(X_0,X_1)X_0 g(X_0,X_1)^{-1},X_1)  \cdot g(X_0,X_1).
$$
The group $\mathbf k^\times$ acts on $(\mathcal G^{\mathrm B},\circledast)$ by 
$\lambda\bullet g(X_0,X_1):=g(X_0^\lambda,X_1^\lambda)$. The resulting semidirect 
product group $\mathcal G^{\mathrm B}\rtimes \mathbf k^\times$ is the set 
$\mathcal G^{\mathrm B}\times \mathbf k^\times$, equipped with the product 
$$
(g(X_0,X_1),\lambda)\circledast (h(X_0,X_1),\mu)
:=(h(g(X_0,X_1)X_0^\lambda g(X_0,X_1)^{-1},X_1^\lambda)  \cdot g(X_0,X_1),\lambda\mu). 
$$

(b) $(g(X_0,X_1),\lambda)\mapsto [X_0\mapsto g(X_0,X_1)X_0^\lambda
g(X_0,X_1)^{-1},X_1\mapsto X_1^\lambda]$ 
defines an action of the group $\mathcal G^{\mathrm B}\rtimes \mathbf k^\times $ on 
$F_2(\mathbf k)$. 
\end{lem}    

The semidirect product $\mathcal G\rtimes\mathbf k^\times$ is similarly the set 
$\mathcal G\times\mathbf k^\times$, equipped with the product 
$(g(e_0,e_1),\lambda)\circledast(h(e_0,e_1),\mu)
:=(g(e_0,e_1)\circledast h(\lambda e_0,\lambda e_1),\lambda\mu)$, the product 
$\circledast$ being as in \eqref{FORMULA:1912}.

\begin{lem} \label{lem:ef3} (see \cite{EF3}, Lem. 1.13) 
If a bitorsor  $(G,X,H)$ contains subbitorsors 
 $(G',X',H')$ and  $(G'',X'',H'')$ such that $G'\subset G''$ and $X'\subset X''$, then 
 $H'\subset H''$. 
\end{lem}



Define $\hat{\mathcal W}^{\mathrm B}:=\mathbf k+(\mathbf kF_2)^\wedge(X_1-1)\subset 
(\mathbf kF_2)^\wedge$; there is a unique continuous $\mathbf k$-algebra morphism 
$\hat\Delta^{\mathcal W,\mathrm B} : \hat{\mathcal W}^{\mathrm B}\to 
\hat{\mathcal W}^{\mathrm B}\hat\otimes
\hat{\mathcal W}^{\mathrm B}$, such that $X_1^{\pm1}\mapsto X_1^{\pm1}\otimes X_1^{\pm1}$ 
and $X_0^k(1-X_1)\mapsto X_0^k(1-X_1)\otimes1
+1\otimes X_0^k(1-X_1)+\sum_{i=1}^{k-1}X_0^i(1-X_1)\otimes X_0^{k-i}(1-X_1)$ for $k\in\mathbb Z$
(with $\sum_{i=1}^{k-1}f(i)$ being defined as $0$ for $k=1$ and 
as $-f(0)-f(-1)\cdots -f(k)$ for $k\leq0$). 
Let 
$\hat{\mathcal M}^{\mathrm B}:=(\mathbf kF_2)^\wedge/(\mathbf kF_2)^\wedge\cdot (X_0-1)$, 
and denote by $x\mapsto x\cdot 1_{\mathrm B}$ 
the natural projection $(\mathbf kF_2)^\wedge\to \hat{\mathcal M}^{\mathrm B}$. 
Then the map $\hat{\mathcal W}^{\mathrm B}\to \hat{\mathcal M}^{\mathrm B}$, $x\mapsto 
x\cdot 1_{\mathrm B}$ is an isomorphism. Let $\hat\Delta^{\mathcal M,\mathrm B} : 
\hat{\mathcal M}^{\mathrm B}\to 
\hat{\mathcal M}^{\mathrm B}\hat\otimes\hat{\mathcal M}^{\mathrm B}$ be the map such that 
$\hat\Delta^{\mathcal M,\mathrm B}(w\cdot 1_{\mathrm B})
=\hat\Delta^{\mathcal W,\mathrm B}(w)\cdot (1_{\mathrm B}\otimes 1_{\mathrm B})$
for any $w\in \hat{\mathcal W}^{\mathrm B}$
and let $\mathcal G(\hat{\mathcal M}^{\mathrm B})$ be set of group-like elements of 
$(\hat{\mathcal M}^{\mathrm B},\hat\Delta^{\mathcal M,\mathrm B})$.

\subsection{The group $(\mathcal G^{\mathrm B}\rtimes\mathbf k^\times)_{\mathrm{inert}}$ 
and its involution $\Theta^{\mathrm B}$}\label{sect:15:2:2912}

For $g\in F_2(\mathbf k)$, let us denote by $\mathcal C(g)$ the conjugacy class of 
$g$ in $F_2(\mathbf k)$. For $\alpha$ a group automorphism of $F_2(\mathbf k)$, denote 
by $\mathcal C(\alpha)$ the permutation of conjugacy classes induced by $\alpha$. 

\begin{defn}\label{def:TautB:IAutB}
(a) $\mathrm{TAut}^{\mathrm B}_{\{0,1\}}$ is the group of automorphisms $\alpha$ of 
$F_2(\mathbf k)$ such that for some $\lambda\in\mathbf k^\times$ (necessarily unique), 
$\alpha(X_1)=X_1^\lambda$, $\mathcal C(\alpha)(\mathcal C(X_0))=\mathcal C(X_0^\lambda)$ 
and $\alpha(X_0)\equiv X_0^\lambda$ mod $\Gamma^3 F_2(\mathbf k)$. 

(b) $X_\infty\in F_2(\mathbf k)$ is defined by $X_\infty:=(X_1X_0)^{-1}$. 

(c) $\mathrm{IAut}^{\mathrm B}_{\{0,1\}}$ is the subgroup of 
$\mathrm{TAut}^{\mathrm B}_{\{0,1\}}$ of all $\alpha$ such that 
$\mathcal C(\alpha)(\mathcal C(X_\infty))=\mathcal C(X_\infty^\lambda)$, where 
$\lambda$ is as above.
\end{defn}

\begin{lem}\label{iso:gp:betti}
The map from Lem. \ref{lem:15:2:14dec}(b) induces a group isomorphism 
\begin{equation}\label{iso:iota}
   \iota^{\mathrm B} :  (\mathcal G^{\mathrm B}\rtimes\mathbf k^\times,\circledast)\to
    \mathrm{TAut}^{\mathrm B}_{\{0,1\}}. 
\end{equation}
\end{lem}

\begin{proof}
By Lem. \ref{lem:15:2:14dec}(b), this map induces a group morphism 
$(\mathcal G^{\mathrm B}\rtimes\mathbf k^\times,\circledast)\to
\mathrm{Aut}(F_2(\mathbf k))$. Its kernel is the intersection of the 
centralizer of $X_0$ with $\mathcal G^{\mathrm B}$, which is 1, therefore it is injective. 
Its image is contained in $\mathrm{TAut}^{\mathrm B}_{\{0,1\}}$
since $(g\in\mathcal G^{\mathrm B},\lambda\in\mathbf k^\times)
\implies((g,X_0^\lambda)\in\Gamma^3F_2(\mathbf k))$. 
It is also checked to be surjective, 
since $(g\in F_2(\mathbf k),\lambda\in\mathbf k^\times,(g,X_0^\lambda)\in 
\Gamma^3F_2(\mathbf k))\implies(\exists\alpha\in\mathbf k, gX_0^\alpha
\in\mathcal G^{\mathrm B})$. 
\end{proof}


\begin{lem}\label{lem:semi:B} (see Lem. \ref{lem:semi:B:wo:proof})
$((\mathcal G^{\mathrm B}\rtimes\mathbf k^\times)_{\mathrm{inert}},\circledast)$
    is a subgroup of $(\mathcal G^{\mathrm B}\rtimes\mathbf k^\times,\circledast)$, and the 
    isomorphism $\iota^{\mathrm B}$ (see \eqref{iso:iota}) restricts to an isomorphism  
    $\iota^{\mathrm B}_{\mathrm{inert}} : ((\mathcal G^{\mathrm B}\rtimes\mathbf k^\times)_{\mathrm{inert}},\circledast)\to
    \mathrm{IAut}^{\mathrm B}_{\{0,1\}}$.  
\end{lem}

\begin{proof} The preimage of $\mathrm{IAut}^{\mathrm B}_{\{0,1\}}$ by the isomorphism 
$(\mathcal G^{\mathrm B}\rtimes\mathbf k^\times,\circledast)\to 
\mathrm{TAut}^{\mathrm B}_{\{0,1\}}$ is exactly 
$(\mathcal G^{\mathrm B}\rtimes\mathbf k^\times)_{\mathrm{inert}}$, which proves the statement. 
\end{proof}

\begin{lem}\label{lem:15:8}
For any $(g,\lambda)\in (\mathcal G^{\mathrm B}\rtimes
\mathbf k^\times)_{\mathrm{inert}}$, there exists a unique $h\in F_2(\mathbf k)$
such that \eqref{def:dual:g:h} holds and $h\equiv X_1^{(\lambda-1)/2}$ 
mod $\Gamma^2 F_2(\mathbf k)$. It will be denoted $h_{g}$. 
\end{lem}

\begin{proof} The map log sets up a bijection between $F_2(\mathbf k)$ and 
$\mathrm{Lie}F_2(\mathbf k)$, which is the free  
$\mathbf k$-Lie algebra generated by $\xi_0:=\mathrm{log}X_0$ and $\xi_1:=\mathrm{log}X_1$. 
Then the image by log of \eqref{def:dual:g:h} gives 
$\lambda \xi_1+\lambda\xi_0+(1/2)\lambda^2[\xi_1,\xi_0]=
\lambda(\xi_1+\xi_0+(1/2)[\xi_1,\xi_0])+[\mathrm{log}h,\lambda(\xi_1+\xi_0)]$
mod $\Gamma^3\mathrm{Lie}F_2(\mathbf k)$, therefore $[\mathrm{log}h,\xi_1+\xi_0]
=(1/2)(\lambda-1)[\xi_1,\xi_0]$ so $\mathrm{log}h \equiv 
(1/2)(\lambda-1)\xi_1$ mod $\mathbf k(\xi_1+\xi_0)+\Gamma^3\mathrm{Lie}F_2(\mathbf k)$. 
The existence of the said $h$ then follows, after multiplying by the appropriate power of 
$X_1X_0$. The uniqueness follows from the implication $(h,h'\in F_2(\mathbf k)$ 
and $\mathrm{Ad}_h(X_0X_1)=\mathrm{Ad}_{h'}(X_0X_1))\implies(\exists\alpha\in\mathbf k, 
h'=h\cdot (X_0X_1)^\alpha)$. 
\end{proof}

\begin{lem}\label{lem:15:12:new:w:proof} (see Lem. \ref{lem:15:12:new})
(a) There is an involutive automorphism $\sigma$ of $F_2(\mathbf k)$, determined by 
$$
\sigma : X_0\mapsto X_1^{-1/2}X_\infty X_1^{1/2},\quad X_\infty\mapsto X_1^{1/2}X_0X_1^{-1/2},
\quad X_1\mapsto X_1.  
$$

(b) The involution $\mathrm{Ad}_\sigma : 
\alpha\mapsto \sigma\alpha\sigma$ of $\mathrm{Aut}(F_2(\mathbf k))$
restricts to an involution of the group $\mathrm{IAut}^{\mathrm B}_{\{0,1\}}$. 

(c) There is a unique involution $\Theta^{\mathrm B}$ of 
$((\mathcal G^{\mathrm B}\rtimes\mathbf k^\times)_{\mathrm{inert}},\circledast)$
such that 
$$
\forall (g,\lambda)\in (\mathcal G^{\mathrm B}\rtimes\mathbf k^\times)_{\mathrm{inert}}, \quad 
\Theta^{\mathrm B}(g,\lambda)=(X_1^{-\lambda/2} \sigma(h_g)X_1^{1/2},\lambda); 
$$
it is intertwined with the involution $\mathrm{Ad}_\sigma$ (see (b)) under the isomorphism 
$\iota^{\mathrm B}_{\mathrm{inert}}$ 
(see \eqref{iso:iota}). 
\end{lem}

\begin{proof}
    (a) follows from $\sigma(X_\infty)\sigma(X_1)\sigma(X_0)
=X_1^{1/2}X_0X_1^{-1/2}X_1X_1^{-1/2}X_\infty X_1^{1/2}=X_1^{1/2}X_0X_\infty X_1^{1/2}=1$. 

To prove (b) and (c), we first prove: 
\begin{equation}\label{statement:1:for:bc}
    \text{$\Theta^{\mathrm B}$ is a self-map of 
$(\mathcal G^{\mathrm B}\rtimes\mathbf k^\times)_{\mathrm{inert}}$}
\end{equation}
and
\begin{equation}\label{statement:2:for:bc}
    \forall (g,\lambda)\in (\mathcal G^{\mathrm B}\rtimes\mathbf k^\times)_{\mathrm{inert}},
    \quad \iota^{\mathrm B}_{\mathrm{inert}}(\Theta^{\mathrm B}(k,\lambda))=\sigma\circ
\iota^{\mathrm B}_{\mathrm{inert}}(k,\lambda)\circ\sigma
\end{equation}
Let us first prove \eqref{statement:1:for:bc}. Let $(g,\lambda)\in 
(\mathcal G^{\mathrm B}\rtimes\mathbf k^\times)_{\mathrm{inert}}$. 
It follows from Lem. \ref{lem:15:8} that $X_1^{-\lambda/2} \sigma(h_g)X_1^{1/2}
\in\mathcal G^{\mathrm B}$. 
\eqref{def:dual:g:h} is rewritten as 
\begin{equation}\label{def:dual:g:h:bis}
X_1^\lambda\mathrm{Ad}_{g}(X_0^\lambda)\mathrm{Ad}_{h_g}(X_\infty^\lambda)=1.  
\end{equation}
Then
\begin{align*}
    &\mathrm{Ad}_{X_1^{\lambda/2}\sigma(g)X_1^{-1/2}}(X_\infty^\lambda)
X_1^{\lambda}\mathrm{Ad}_{X_1^{-\lambda/2} \sigma(h_g)X_1^{1/2}}(X_0^\lambda)
 =
X_1^{\lambda/2}\mathrm{Ad}_{\sigma(g)X_1^{-1/2}}(X_\infty^\lambda)
\mathrm{Ad}_{\sigma(h_g)X_1^{1/2}}(X_0^\lambda)X_1^{\lambda/2}
\\& 
=\mathrm{Ad}_{X_1^{-\lambda/2}}(
X_1^{\lambda}\mathrm{Ad}_{\sigma(g)X_1^{-1/2}}(X_\infty^\lambda)
\mathrm{Ad}_{\sigma(h_g)X_1^{1/2}}(X_0^\lambda)
) 
=\mathrm{Ad}_{X_1^{-\lambda/2}}\circ\sigma(
X_1^{\lambda}\mathrm{Ad}_{g}(X_0^\lambda)
\mathrm{Ad}_{h_g}(X_\infty^\lambda)
)=1
\end{align*}
where the last equality follows from \eqref{def:dual:g:h:bis}. It follows that 
$(k,\lambda):=(X_1^{-\lambda/2} \sigma(h_g)X_1^{1/2},\lambda)\in 
(\mathcal G^{\mathrm B}\rtimes\mathbf k^\times)_{\mathrm{inert}}$. 
Since $X_1^{\lambda/2}\sigma(g)X_1^{-1/2}\equiv X_1^{(\lambda-1)/2}$ mod 
$\Gamma^2F_2(\mathbf k)$, $h_k=X_1^{\lambda/2}\sigma(g)X_1^{-1/2}$. 
This proves \eqref{statement:1:for:bc}. 

Let us now prove \eqref{statement:2:for:bc}. One has for any $(g,\lambda)\in 
(\mathcal G^{\mathrm B}\rtimes\mathbf k^\times)_{\mathrm{inert}}$ the equality  
\begin{align*}
&\iota^{\mathrm B}_{\mathrm{inert}}(X_1^{-\lambda/2}\sigma(h_g)X_1^{1/2},\lambda)= 
[X_1\mapsto X_1^\lambda,X_0\mapsto X_1^{-\lambda/2}\sigma(h_g)X_1^{1/2}
X_0^\lambda X_1^{-1/2}\sigma(h_g)^{-1}X_1^{\lambda/2},\\ & 
X_\infty\mapsto X_1^{\lambda/2}\sigma(g)X_1^{-1/2}X_\infty^\lambda 
X_1^{1/2}\sigma(g)^{-1}X_1^{-\lambda/2}]
\\ & =\sigma\circ[X_1\mapsto X_1^\lambda,X_0\mapsto gX_0^\lambda g^{-1},
X_\infty\mapsto h_gX_\infty^\lambda h_g^{-1}]\circ\sigma
= \sigma\circ\iota^{\mathrm B}_{\mathrm{inert}}(g,\lambda)\circ\sigma 
\end{align*}
which implies \eqref{statement:2:for:bc}. 

Let us now prove (b). Let $\alpha\in \mathrm{IAut}^{\mathrm B}_{\{0,1\}}$, then 
there exists $(g,\lambda)\in (\mathcal G^{\mathrm B}\rtimes\mathbf k^\times)_{\mathrm{inert}}$
such that $\alpha=\iota^{\mathrm B}(g,\lambda)$. Then by \eqref{statement:1:for:bc}, $\sigma\alpha\sigma=
\iota^{\mathrm B}(\Theta^{\mathrm B}(g,\lambda))$ where $\Theta^{\mathrm B}(g,\lambda)\in (
\mathcal G^{\mathrm B}\rtimes\mathbf k^\times)_{\mathrm{inert}}$ by \eqref{statement:2:for:bc}, 
therefore $\sigma\alpha\sigma\in \mathrm{IAut}^{\mathrm B}_{\{0,1\}}$.  
Since $\alpha\mapsto\sigma\alpha\sigma$ is a set-theoretic involution of 
$\mathrm{Aut}(F_2(\mathbf k))$, it follows that this is a set-theoretic involution of 
$\mathrm{IAut}^{\mathrm B}_{\{0,1\}}$, and therefore an involution of the group 
$\mathrm{IAut}^{\mathrm B}_{\{0,1\}}$. 
 
Let us now prove (c). For $(g,\lambda)\in 
(\mathcal G^{\mathrm B}\rtimes\mathbf k^\times)_{\mathrm{inert}}$, one has 
$$
\iota^{\mathrm B}_{\mathrm{inert}}((\Theta^{\mathrm B})^2(g,\lambda))=\sigma\circ
\iota^{\mathrm B}_{\mathrm{inert}}(\Theta^{\mathrm B}(g,\lambda))\circ\sigma
=\sigma^2\circ
\iota^{\mathrm B}(g,\lambda)\circ\sigma^2=\iota^{\mathrm B}_{\mathrm{inert}}(g,\lambda)
$$
by \eqref{statement:2:for:bc} and since $\sigma^2=id$. The fact that 
$\iota^{\mathrm B}_{\mathrm{inert}}$ is a bijection 
then implies $(\Theta^{\mathrm B})^2(g,\lambda)=(g,\lambda)$, hence $(\Theta^{\mathrm B})^2=id$. 
For $(g,\lambda),(g',\lambda')\in 
(\mathcal G^{\mathrm B}\rtimes\mathbf k^\times)_{\mathrm{inert}}$, one has 
\begin{align*}
& \iota^{\mathrm B}_{\mathrm{inert}}(\Theta^{\mathrm B}(g,\lambda)
\circledast\Theta^{\mathrm B}(g',\lambda'))
=\iota^{\mathrm B}_{\mathrm{inert}}(\Theta^{\mathrm B}(g,\lambda))
\iota^{\mathrm B}_{\mathrm{inert}}(\Theta^{\mathrm B}(g',\lambda'))
=\sigma\iota^{\mathrm B}_{\mathrm{inert}}(g,\lambda)\sigma\sigma
\iota^{\mathrm B}_{\mathrm{inert}}(g',\lambda')\sigma
\\&=\sigma\iota^{\mathrm B}_{\mathrm{inert}}(g,\lambda)
\iota^{\mathrm B}_{\mathrm{inert}}(g',\lambda')\sigma
=\sigma\iota^{\mathrm B}_{\mathrm{inert}}((g,\lambda)\circledast(g',\lambda'))\sigma
=\iota^{\mathrm B}_{\mathrm{inert}}(\Theta^{\mathrm B}((g,\lambda)\circledast(g',\lambda')))
\end{align*}
where the first and fourth equalities follow from Lem. \ref{iso:gp:betti}, the 
second and fifth equalities follow from \eqref{statement:2:for:bc} and 
the third equality follows from (a). Lem. \ref{iso:gp:betti} then implies 
$\Theta^{\mathrm B}(g,\lambda)\circledast\Theta^{\mathrm B}(g',\lambda')=
\Theta^{\mathrm B}((g,\lambda)\circledast(g',\lambda'))$; this equality takes place in 
$\mathcal G^{\mathrm B}\rtimes\mathbf k^\times$, therefore also in 
$(\mathcal G^{\mathrm B}\rtimes\mathbf k^\times)_{\mathrm{inert}}$ since both sides belong
to this subgroup. (c) follows.  
\end{proof}

\subsection{Right torsors}\label{subsect:15:3:rt:2912}

\begin{defn}
Let $\mu\in\mathbf k^\times$. 

(a) $\mathrm{TIso}^\mu_{\{0,1\}}$ is the set of isomorphisms $\gamma : F_2(\mathbf k)\to 
 \mathrm{exp}(\mathfrak{lie}_{\{0,1\}}^\wedge)$ such that 
$\gamma(X_1)=e^{\mu e_1}$, $\mathcal C(\gamma)(\mathcal C(X_0))=\mathcal C(e^{\mu e_0})$
and $\gamma(X_0)\equiv e^{\mu e_0}$ mod $\Gamma^3\mathrm{exp}(\mathfrak{lie}_{\{0,1\}}^\wedge)$. 

(b) $\mathrm{IIso}^\mu_{\{0,1\}}$ is the subset of $\mathrm{TIso}^\mu_{\{0,1\}}$ 
of isomorphisms $\gamma$ such that 
$\mathcal C(\gamma)(\mathcal C(X_\infty))=\mathcal C(e^{\mu e_\infty})$. 
\end{defn}

\begin{defn}
Set
$$
\mathrm{TIso}_{\{0,1\}}:=\{(\gamma,\mu)|\mu\in\mathbf k^\times,\gamma\in 
\mathrm{TIso}^\mu_{\{0,1\}}\}=\sqcup_{\mu\in\mathbf k^\times}\mathrm{TIso}^\mu_{\{0,1\}},
$$
$$
\mathrm{IIso}_{\{0,1\}}:=\{(\gamma,\mu)|\mu\in\mathbf k^\times,\gamma\in 
\mathrm{IIso}^\mu_{\{0,1\}}\}=\sqcup_{\mu\in\mathbf k^\times}\mathrm{IIso}^\mu_{\{0,1\}}.  
$$
Then $\mathrm{IIso}_{\{0,1\}}\subset \mathrm{TIso}_{\{0,1\}}$. 
\end{defn}

\begin{lem}\label{lem:1515marignan}
(a)  The group $\mathrm{TAut}^{\mathrm B}_{\{0,1\}}$ (see Def. \ref{def:TautB:IAutB})  
 acts from the right on $\mathrm{TIso}_{\{0,1\}}$
by $\mathrm{TIso}_{\{0,1\}}\times \mathrm{TAut}^{\mathrm B}_{\{0,1\}}
\ni ((\gamma,\mu),(\alpha,\lambda))\to (\gamma\circ\alpha,\lambda\mu)\in 
\mathrm{TIso}_{\{0,1\}}^{\lambda\mu}$, and 
$(\mathrm{TIso}_{\{0,1\}},\mathrm{TAut}^{\mathrm B}_{\{0,1\}})$ is a right torsor. 

(b)  The action from (a) restricts to a action of $\mathrm{IAut}^{\mathrm B}_{\{0,1\}}$
from the right on $\mathrm{IIso}_{\{0,1\}}$, and $(\mathrm{IIso}_{\{0,1\}},
\mathrm{IAut}^{\mathrm B}_{\{0,1\}})$ is a subtorsor of the right torsor from (a). 

(c) The assignment $\ell(s_{(0,\infty)})r(\sigma) : 
(\gamma,\mu)\mapsto (s_{(0,\infty)}\gamma\sigma,\mu)$ is an 
set-theoretic involution of $\mathrm{IIso}_{\{0,1\}}$, and the pair 
$(\ell(s_{(0,\infty)})r(\sigma),\mathrm{Ad}_\sigma)$, where $\mathrm{Ad}_\sigma$ is as in 
Lem. \ref{lem:15:12:new}(b), is an involution of the right torsor $(\mathrm{IIso}_{\{0,1\}},
\mathrm{IAut}^{\mathrm B}_{\{0,1\}})$. 
\end{lem}

\begin{proof}
(a) and (b) are immediate. (c) follows from the equality
$\ell(s_{(0,\infty)})r(\sigma)(\gamma\alpha)=s_{(0,\infty)}\gamma\alpha\sigma
=s_{(0,\infty)}\gamma\sigma\sigma\alpha\sigma
=\ell(s_{(0,\infty)})r(\sigma)(\gamma)\mathrm{Ad}_\sigma(\alpha)$
for any $\mu\in\mathbf k^\times,\gamma\in 
\mathrm{TIso}^\mu_{\{0,1\}}$ and $(\gamma,\lambda)\in 
\mathrm{TAut}^{\mathrm B}_{\{0,1\}}$. 
\end{proof}

\begin{defn}
$\mathrm{iso}_\mu : F_2(\mathbf k)\to 
 \mathrm{exp}(\mathfrak{lie}_{\{0,1\}}^\wedge)$ is the isomorphism induced by 
 $X_0\mapsto e^{\mu e_0}$, $X_\infty\mapsto e^{\mu e_\infty}$ (recall that 
 $X_\infty X_1X_0=1$ and $e_0+e_1+e_\infty=0$). 
\end{defn} 

\begin{lem}
The map $\iota:\mathcal G\times\mathbf k^\times\to \mathrm{TIso}_{\{0,1\}}$ given by 
$(\phi,\mu)\mapsto (\mathrm{aut}^{\mathcal V,\mu}_\phi\circ\mathrm{iso}_\mu,\mu)$ 
is a bijection. The pair $(\iota,\iota^{\mathrm B})$ is an isomorphism between the right torsors 
$(\mathcal G\times\mathbf k^\times,\mathcal G^{\mathrm B}\rtimes\mathbf k^\times)$ (equipped 
with the structure from Lem. \ref{lem:GGG:bitorsor}) and $(\mathrm{TIso}_{\{0,1\}},
\mathrm{TAut}^{\mathrm B}_{\{0,1\}})$, where $\iota^{\mathrm B}$ is as in \eqref{iso:iota}. 
\end{lem}

\begin{proof}
For $\mu\in\mathbf k^\times$, the map $\beta\mapsto \beta\circ \mathrm{iso}_\mu$ sets up 
a bijection $\mathrm{TAut}_{\{0,1\}}^\mu\to\mathrm{TIso}_{\{0,1\}}^\mu$. 
The composition of this bijection with the bijection from Lem. \ref{lem:14:2}
is a bijection $\mathcal G\to \mathrm{TIso}_{\{0,1\}}^\mu$. As $\iota$ is the disjoint union 
of these maps over $\mu\in\mathbf k^\times$, it is a bijection. 

Let $(\phi,\mu)\in\mathcal G\times\mathbf k^\times$ and $(g,\lambda)\in\mathcal G^{\mathrm B}
\times\mathbf k^\times$. Then 
\begin{align*}
&\iota(\phi,\mu)\iota^{\mathrm B}(g,\lambda)
=[e_0\mapsto \phi e_0 \phi^{-1},e_0*_\mu e_\infty\mapsto e_0+e_\infty][X_0\mapsto e^{\mu e_0},X_\infty\mapsto e^{\mu e_\infty}][X_0\mapsto g X_0^\lambda g^{-1},X_1\mapsto X_1^\lambda]
\\ & =[X_0\mapsto\phi e^{\mu e_0} \phi^{-1},X_1\mapsto e^{\mu e_1}]
[X_0\mapsto g X_0^\lambda g^{-1},X_1\mapsto X_1^\lambda]
\\&
=[X_0\mapsto g(\phi e^{\mu e_0} \phi^{-1},e^{\mu e_1})\phi e^{\mu\lambda e_0}\phi^{-1}
g(\phi e^{\mu e_0} \phi^{-1},e^{\mu e_1})^{-1},X_1\mapsto e^{\mu\lambda e_1}]\\& 
=[e_0\mapsto g(\phi e^{\mu e_0} \phi^{-1},e^{\mu e_1})\phi\cdot e_0\cdot 
\phi^{-1}g(\phi e^{\mu e_0} \phi^{-1},e^{\mu e_1})^{-1},e_0*_{\mu\lambda} e_\infty
\mapsto e_0+e_\infty][X_0\mapsto e^{\mu\lambda e_0},X_\infty\mapsto e^{\mu\lambda e_\infty}]
\\& =\iota(g(\phi e^{\mu e_1}\phi^{-1},e^{\mu e_1})\phi,\mu\lambda)
=\iota((\phi,\mu)\circledast(g,\lambda)). 
\end{align*}
The right torsor morphism property of $(\iota,\iota^{\mathrm B})$ follows.  
\end{proof}

\begin{lem}\label{righttorsor:1510}
(a) $(\sqcup_{\mu\in\mathbf k^\times}\mathcal G_{\mathrm{inert}}^\mu,
(\mathcal G^{\mathrm B}\rtimes\mathbf k^\times)_{\mathrm{inert}})$ 
is a right subtorsor of $(\mathcal G\times\mathbf k^\times,\mathcal G^{\mathrm B}\rtimes
\mathbf k^\times)$ (see notation in Def. \ref{def:G:mu:inert}). 

(b) The pair $(\iota,\iota^{\mathrm B})$ restricts to an isomorphism 
$(\iota_{\mathrm{inert}},\iota_{\mathrm{inert}}^{\mathrm B})$ 
between the right torsors 
$$(\sqcup_{\mu\in\mathbf k^\times}\mathcal G^\mu_{\mathrm{inert}}, 
    (\mathcal G^{\mathrm B}\rtimes\mathbf k^\times)_{\mathrm{inert}})
    $$ 
    and $(\mathrm{IIso}_{\{0,1\}},
    \mathrm{IAut}^{\mathrm B}_{\{0,1\}})$, where $\iota_{\mathrm{inert}}^{\mathrm B}$ is 
    as in Lem. \ref{lem:semi:B}. 
   
\end{lem}

\begin{proof}
For $\mu\in\mathbf k^\times$, the bijection $\mathrm{TAut}_{\{0,1\}}^\mu\to
\mathrm{TIso}_{\{0,1\}}^\mu$ induced by the map $\beta\mapsto \beta\circ \mathrm{iso}_\mu$ 
restricts to a bijection $\mathrm{IAut}_{\{0,1\}}^\mu\to\mathrm{IIso}_{\{0,1\}}^\mu$
 (see the notation in Def. \ref{def:14:1}); combining this with Lem. \ref{lem:14:2}, 
 one sees that the composed bijection $\mathcal G\to \mathrm{TAut}_{\{0,1\}}^\mu\to
\mathrm{TIso}_{\{0,1\}}^\mu$, where the first map is $g\mapsto \mathrm{aut}_g^{\mathcal V,\mu}$, 
restricts to a bijection $\mathcal G^\mu_{\mathrm{inert}}\to\mathrm{IIso}_{\{0,1\}}^\mu$. 
It follows that $\iota$ restricts to a bijection $\iota_{\mathrm{inert}} : 
\sqcup_{\mu\in\mathbf k^\times}\mathcal G_{\mathrm{inert}}^\mu\to \mathrm{TIso}_{\{0,1\}}$. 
(a) and (b) then follow from Lem. \ref{lem:semi:B} and Lem. \ref{lem:1515marignan}. 
\end{proof}

\begin{lem}\label{lem:intertw:Theta:Ad}
The bijection $\iota_{\mathrm{inert}}$ intertwines the automorphisms
$\ell(s_{(0,\infty)})r(\sigma)$ and $\Theta:=\sqcup_{\mu\in\mathbf k^\times}\Theta^\mu$ 
(see Lem. \ref{def:Theta:mu}) of its source and target.  
\end{lem}

\begin{proof}
Let $\mu\in\mathbf k^\times$ and $\phi\in\mathcal G^\mu_{\mathrm{inert}}$. Then 
 $\phi\in\mathcal G$, $h_\phi\in\mathrm{exp}(\mathfrak{lie}_{\{0,1\}}^\wedge)$, 
 $\mathrm{log}h_\phi\equiv(\mu/2)e_1$ mod $F^2\mathfrak{lie}_{\{0,1\}}^\wedge$,  
$\mathrm{Ad}_\phi(e_0)*_\mu \mathrm{Ad}_{h_\phi}(e_\infty)=e_0+e_\infty$, and 
$\Theta(\phi,\mu)=(\Theta^\mu(\phi),\mu)$, where 
$\Theta^\mu(\phi)=e^{-\mu \mathrm{ad}e_1/2}(s_{(0,\infty)}(h_\phi)))$. 

\begin{align*}
& s_{(0,\infty)}\iota_{\mathrm{inert}}(\phi,\mu)\sigma
\\&=[e_0\leftrightarrow e_\infty,e_1\mapsto e_1][e_0\mapsto \phi e_0\phi^{-1},
e_\infty\mapsto h_\phi e_\infty h_\phi^{-1},e_\infty*_\mu e_0\mapsto e_\infty+e_0]
\\ & 
[X_0\mapsto e^{\mu e_0},X_\infty\mapsto e^{\mu e_\infty},X_1\mapsto 
e^{-\mu\cdot e_\infty*_\mu e_0}][X_0\mapsto X_1^{-1/2}X_\infty X_1^{1/2},X_\infty\mapsto X_1^{1/2}X_0 X_1^{-1/2},
X_1\mapsto X_1]
\\&=[e_0\mapsto \mathrm{Ad}_{s_{(0,\infty)}(\phi)}(e_\infty),
e_\infty\mapsto \mathrm{Ad}_{s_{(0,\infty)}(h_\phi)}(e_0),e_\infty*_\mu e_0\mapsto -e_1]
\\ & [X_1\mapsto e^{-\mu \cdot e_\infty*_\mu e_0},X_0\mapsto 
\mathrm{Ad}_{e^{(\mu/2) e_\infty*_\mu e_0}}(e^{\mu e_\infty}),X_\infty\mapsto 
\mathrm{Ad}_{e^{-(\mu/2) e_\infty*_\mu e_0}}(e^{\mu e_0})]
\\ & = [X_1\mapsto e^{\mu e_1},X_0\mapsto \mathrm{Ad}_{e^{-(\mu/2)e_1}s_{(0,\infty)}(h_\phi)}
(e^{\mu e_0}),X_\infty\mapsto \mathrm{Ad}_{e^{(\mu/2)e_1}s_{(0,\infty)}(\phi)}(e^{\mu e_\infty})]
\\&=[e_0\mapsto e^{-\mu e_1/2}s_{(0,\infty)}(h_\phi)e^{\mu e_1/2}\cdot e_0\cdot 
e^{-\mu e_1/2}s_{(0,\infty)}(h_\phi)^{-1}e^{\mu e_1/2},e_\infty*_\mu e_0\mapsto e_\infty+e_0]
\\ & 
[X_0\mapsto e^{\mu e_0},X_\infty\mapsto e^{\mu e_\infty},X_1\mapsto e^{-\mu(e_\infty*_\mu e_0)}]
=\iota_{\mathrm{inert}}(e^{-\mu e_1/2}s_{(0,\infty)}(h_\phi)e^{\mu e_1/2},\mu)
=\iota_{\mathrm{inert}}(\Theta^\mu(\phi),\mu)
\\ & =\iota_{\mathrm{inert}}(\Theta(\phi,\mu)). 
\end{align*}
\end{proof}

\begin{lem}\label{lem:compduality}
 The right action of $(\mathcal G^{\mathrm B}\rtimes \mathbf k^\times)_{\mathrm{inert}}$ 
 on $\sqcup_{\mu\in\mathbf k^\times}\mathcal G_{\mathrm{inert}}^\mu$ from 
 Lem. \ref{righttorsor:1510}(a)
 exhibits the following compatibility with the group involution $\Theta^{\mathrm B}$ 
 (Lem. \ref{lem:15:12:new}(c)) and 
 the involution $\Theta=\sqcup_{\mu\in\mathbf k^\times}\Theta^\mu$ (Lem. \ref{def:Theta:mu}): 
 $$
 \forall \mu\in\mathbf k^\times,\forall \phi\in\mathcal G^\mu_{\mathrm{inert}},
 \forall(g,\lambda)\in
 (\mathcal G^{\mathrm B}\rtimes\mathbf k^\times)_{\mathrm{inert}},\quad 
\Theta((\phi,\mu)\bullet(g,\lambda))=\Theta(\phi,\mu)\bullet\Theta^{\mathrm B}(g,\lambda). 
 $$
\end{lem}

\begin{proof}
For $\mu,\phi,g,\lambda$ as above, one has 
\begin{align*}
&\iota_{\mathrm{inert}}(\Theta((\phi,\mu)\bullet(g,\lambda)))
=s_{(0,\infty)}\iota_{\mathrm{inert}}((\phi,\mu)\bullet(g,\lambda))\sigma
=s_{(0,\infty)}\iota_{\mathrm{inert}}(\phi,\mu)\iota_{\mathrm{inert}}^{\mathrm B}(g,\lambda)\sigma
\\& =s_{(0,\infty)}\iota_{\mathrm{inert}}(\phi,\mu)\sigma\sigma\iota_{\mathrm{inert}}^{\mathrm B}(g,\lambda)\sigma
=\iota_{\mathrm{inert}}(\Theta(\phi,\mu))\iota_{\mathrm{inert}}^{\mathrm B}(\Theta^{\mathrm B}(g,\lambda))
=\iota_{\mathrm{inert}}(\Theta(\phi,\mu)\bullet\Theta^{\mathrm B}(g,\lambda)), 
\end{align*}
where the first equality follows from Lem. \ref{lem:intertw:Theta:Ad}, 
the second and fifth equalities follow from Lem. \ref{righttorsor:1510}(b), 
the third equality follows from $\sigma^2=id$, 
the fourth equality follows from Lem. \ref{lem:15:12:new}(c) and Lem. \ref{lem:intertw:Theta:Ad}. 
\end{proof}

\subsection{The inclusion $\mathsf{DMR}^{\mathrm B}(\mathbf k)\subset 
(\mathcal G^{\mathrm B}\rtimes\mathbf k^\times)_{\mathrm{inert}}$ and its
stability under $\Theta^{\mathrm B}$}\label{sect:15:4:tiXXXaisuXXX:2912}

\begin{lem}\label{lem:15:14}
(a) (see Lem. \ref{lem:actions:ktimes:b}(a)) 
The subgroup $\mathcal G_{\mathrm{inert}}\subset\mathcal G$ is stable under the action of 
$\mathbf k^\times$, giving rise to a group inclusion $\mathcal G_{\mathrm{inert}}
\rtimes\mathbf k^\times\subset \mathcal G\rtimes\mathbf k^\times$. 

(b) The left action of $\mathcal G\rtimes\mathbf k^\times$ on 
$\mathcal G\times\mathbf k^\times=\sqcup_{\mu\in\mathbf k^\times}\mathcal G$ from Lem. \ref{lem:GGG:bitorsor} 
restricts to an action of 
$\mathcal G_{\mathrm{inert}}\rtimes\mathbf k^\times$ on 
$\sqcup_{\mu\in\mathbf k^\times}\mathcal G^\mu_{\mathrm{inert}}$. 

(c) (see Lem. \ref{lem:G:mu}) The left and right actions of $\mathcal G_{\mathrm{inert}}
\rtimes\mathbf k^\times$ and $(\mathcal G^{\mathrm B}\rtimes
\mathbf k^\times)_{\mathrm{inert}}$ on 
$\sqcup_{\mu\in\mathbf k^\times}\mathcal G^\mu_{\mathrm{inert}}$ give rise to a bitorsor 
$(\mathcal G_{\mathrm{inert}}\rtimes\mathbf k^\times,
\sqcup_{\mu\in\mathbf k^\times}\mathcal G_{\mathrm{inert}}^\mu,(\mathcal G^{\mathrm B}\rtimes \mathbf k^\times)_{\mathrm{inert}})$, which is a subbitorsor of 
$(\mathcal G\rtimes\mathbf k^\times,
\sqcup_{\mu\in\mathbf k^\times}\mathcal G,\mathcal G^{\mathrm B}\rtimes \mathbf k^\times)$ (see Lem. \ref{lem:GGG:bitorsor}).  
\end{lem}

\begin{proof}
(a) Let $g\in\mathcal G_{\mathrm{inert}}$, then $h_g\in \mathcal G_{\mathrm{inert}}$
and $\mathrm{Ad}_g(e_0)+e_1+\mathrm{Ad}_{h_g}(e_\infty)=0$. Acting by 
$\lambda\in\mathbf k^\times$ on the left-hand side and dividing by $\lambda$, one obtains 
$\mathrm{Ad}_{\lambda\bullet g}(e_0)+e_1+\mathrm{Ad}_{\lambda\bullet h_g}(e_\infty)=0$, 
where $\lambda\bullet h_g\in\mathcal G_{\mathrm{inert}}$. Therefore 
$\lambda\bullet g\in\mathcal G_{\mathrm{inert}}$. 

(b) By Lem. \ref{lem:14:6}, the action of $\mathcal G$ on 
itself restricts to an action of $\mathcal G_{\mathrm{inert}}$  on 
$\mathcal G^\mu_{\mathrm{inert}}$. Moreover, for each 
$\mu\in\mathbf k^\times$, the action of $\lambda\in\mathbf k^\times$, 
which induces a bijection $\mathcal G^\mu\to\mathcal G^{\lambda\mu}$, restricts to a bijection 
$\mathcal G^\mu_{\mathrm{inert}}\to\mathcal G^{\lambda\mu}_{\mathrm{inert}}$: indeed, if
$\phi\in\mathcal G^\mu_{\mathrm{inert}}$, then 
 $\phi\in\mathcal G$, $h_\phi\in\mathrm{exp}(\mathfrak{lie}_{\{0,1\}}^\wedge)$, and 
$\mathrm{Ad}_\phi(e_0)*_\mu \mathrm{Ad}_{h_\phi}(e_\infty)=e_0+e_\infty$. As the action of 
$\mathbf k^\times$ is by Lie algebra automorphisms, one derives 
$\lambda\bullet(\mathrm{Ad}_\phi(e_0))*_\mu 
\lambda\bullet(\mathrm{Ad}_{h_\phi}(e_\infty))=\lambda(e_0+e_\infty)$, therefore
$(\lambda\mathrm{Ad}_{\lambda\bullet\phi}(e_0))*_\mu
(\lambda\mathrm{Ad}_{\lambda\bullet h_\phi}(e_\infty))=\lambda(e_0+e_\infty)$, 
which implies $\mathrm{Ad}_{\lambda\bullet\phi}(e_0)*_{\lambda\mu}
\mathrm{Ad}_{\lambda\bullet h_\phi}(e_\infty)=e_0+e_\infty$ using the identity
$(\lambda a)*_\mu(\lambda b)=\lambda\cdot(a*_{\lambda\mu}b)$ for any 
 $a,b\in \mathfrak{lie}_{\{0,1\}}^\wedge$.  
This implies the announced restriction of actions. 

(c) The said actions commute as they are restrictions of commuting actions of 
$\mathcal G\rtimes\mathbf k^\times$ and $\mathcal G^{\mathrm B}\rtimes
\mathbf k^\times$ on $\sqcup_{\mu\in\mathbf k^\times}\mathcal G^\mu$. 
These actions also have trivial stabilizers as they are restrictions of actions 
with trivial stabilizer. The action of $\mathcal G_{\mathrm{inert}}$ on 
$\mathcal G^1_{\mathrm{inert}}$ is transitive by Lem. \ref{lem:14:6}, and 
$\lambda\cdot \mathcal G^1_{\mathrm{inert}}=\mathcal G^\lambda_{\mathrm{inert}}$ 
for any $\lambda\in\mathbf k^\times$, therefore the action of 
$\mathcal G_{\mathrm{inert}}\rtimes\mathbf k^\times$ on 
$\sqcup_{\mu\in\mathbf k^\times}\mathcal G^\mu_{\mathrm{inert}}$ is transitive. 
This implies that $(\mathcal G\rtimes\mathbf k^\times,
\sqcup_{\mu\in\mathbf k^\times}\mathcal G^\mu)$ is a left torsor. The result then follows 
from the combination of this an Lem. \ref{righttorsor:1510}(a). 
\end{proof}

\begin{thm}\label{thm:betti:a} (see Thm. \ref{thm:0:33})
There holds the inclusion $\mathsf{DMR}^{\mathrm B}(\mathbf k)\subset 
(\mathcal G^{\mathrm B}\rtimes\mathbf k^\times)_{\mathrm{inert}}$ of subgroups of 
$\mathcal G^{\mathrm B}\rtimes\mathbf k^\times$
(see Lem. \ref{lem:semi:B}). 
\end{thm}

\begin{proof}
  It follows from Lem. \ref{lem:15:14}(c) that  
  $(\mathcal G_{\mathrm{inert}}\rtimes\mathbf k^\times,\sqcup_{\mu\in\mathbf k^\times}
  \mathcal G^\mu_{\mathrm{inert}},(\mathcal G^{\mathrm B}\rtimes
  \mathbf k^\times)_{\mathrm{inert}})$ 
  is a sub-bitorsor of $(\mathcal G\rtimes\mathbf k^\times,\mathcal G\times\mathbf k^\times,
\mathcal G^{\mathrm B}\rtimes\mathbf k^\times)$. It follows from Lem. \ref{lem:15:5}(b)
that the same is true of $(\mathsf{DMR}_0(\mathbf k)\rtimes\mathbf k^\times,
\sqcup_{\mu\in\mathbf k^\times}\mathsf{DMR}_\mu(\mathbf k),\mathsf{DMR}^{\mathrm B}(\mathbf k))$.
Cor. \ref{main:cor}(a) implies the inclusion $\mathsf{DMR}_0(\mathbf k)\rtimes\mathbf k^\times
\subset \mathcal G_{\mathrm{inert}}\rtimes\mathbf k^\times$ and 
Thm. \ref{thm:DMRmu}(a) implies the inclusion 
$\sqcup_{\mu\in\mathbf k^\times}\mathsf{DMR}_\mu(\mathbf k)\subset
\sqcup_{\mu\in\mathbf k^\times}\mathcal G^\mu_{\mathrm{inert}}$. 
Lem. \ref{lem:ef3} then implies the result. 
\end{proof}

\begin{cor}\label{lem:15:16}
    The subbitorsor $(\mathsf{DMR}_0(\mathbf k)\rtimes\mathbf k^\times,
\sqcup_{\mu\in\mathbf k^\times}\mathsf{DMR}_\mu(\mathbf k),
  \mathsf{DMR}^{\mathrm B}(\mathbf k))$ of 
  $$
  (\mathcal G\rtimes\mathbf k^\times,
  \sqcup_{\mu\in\mathbf k^\times} \mathcal G,
\mathcal G^{\mathrm B}\rtimes\mathbf k^\times)
$$ is in fact a subbitorsor of 
 $(\mathcal G_{\mathrm{inert}}\rtimes\mathbf k^\times,\sqcup_{\mu\in\mathbf k^\times}
 \mathcal G_{\mathrm{inert}}^\mu,
(\mathcal G^{\mathrm B}\rtimes\mathbf k^\times)_{\mathrm{inert}})$. 
\end{cor}

\begin{proof}
The inclusion $\mathsf{DMR}_0(\mathbf k)\subset \mathcal G_{\mathrm{inert}}$ 
(see Cor. \ref{main:cor}(a)) implies the inclusion $\mathsf{DMR}_0(\mathbf k)
\rtimes\mathbf k^\times\subset \mathcal G_{\mathrm{inert}}\rtimes\mathbf k^\times$. 
The inclusion $\mathsf{DMR}_\mu(\mathbf k)\subset \mathcal G^\mu_{\mathrm{inert}}$
for any $\mu\in\mathbf k^\times$ (see Thm. \ref{thm:DMRmu}(a)) implies 
$\sqcup_{\mu\in\mathbf k^\times}\mathsf{DMR}_\mu(\mathbf k)\subset 
\sqcup_{\mu\in\mathbf k^\times}\mathcal G^\mu_{\mathrm{inert}}$. 
Finally, the inclusion $\mathsf{DMR}^{\mathrm B}(\mathbf k) \subset (\mathcal G^{\mathrm B} 
\rtimes \mathbf k^\times)_{\mathrm{inert}}$ follows from Thm. \ref{thm:betti:a}. 
The result follows. 
\end{proof}

\begin{lem}\label{lem:15:15}
(a)  (see Lem. \ref{lem:actions:ktimes:b}(b)) 
The action of $\Theta$ on $\mathcal G_{\mathrm{inert}}$ commutes with the action of 
$\mathbf k^\times$, resulting in an involution $\Theta\rtimes id$ of $\mathcal G_{\mathrm{inert}}\rtimes
\mathbf k^\times$. 
 
(b) \Add{(see Lem. \ref{lem:invol:bitorsor})} 
The triple $(\Theta\rtimes id,\sqcup_{\mu\in\mathbf k^\times}\Theta^\mu,\Theta^{\mathrm B})$ 
(see (a), Lem. \ref{def:Theta:mu}, Lem. \ref{lem:15:12:new}) is an involution of the bitorsor  
$(\mathcal G_{\mathrm{inert}}\rtimes
\mathbf k^\times,\sqcup_{\mu\in\mathbf k^\times}\mathcal G_{\mathrm{inert}}^\mu,
(\mathcal G^{\mathrm B}\rtimes \mathbf k^\times)_{\mathrm{inert}})$. 
\end{lem}

\begin{proof}
    (b) Let $\lambda\in\mathbf k^\times$ and $g\in \mathcal G_{\mathrm{inert}}$. 
It follows from Lem. \ref{lem:15:14}(a) that 
$\lambda\bullet g\in\mathcal G_{\mathrm{inert}}$, and that 
$h_{\lambda\bullet g}=\lambda\bullet h_g$. This is the second equality in  
$\Theta(\lambda\bullet g)=s_{(0,\infty)}(h_{\lambda\bullet g})
=s_{(0,\infty)}(\lambda\bullet h_g)=\lambda\bullet s_{(0,\infty)}(h_g)
=\lambda\bullet \Theta(g)$, where the first and last equalities follow from the definition of 
$\Theta$, and the third equality follows from the commutation of the action of $\mathbf k^\times$
with $s_{(0,\infty)}$. 
   
(c) Let $\lambda\in\mathbf k^\times$, $\mu\in\mathbf k^\times$, 
$\phi\in\mathcal G^\mu_{\mathrm{inert}}$. Then 
$$
\Theta^{\lambda\mu}(\lambda\bullet\phi)
=s_{(0,\infty)}(h_{\lambda\bullet\phi})
=s_{(0,\infty)}(\lambda\bullet h_{\phi})
=\lambda\bullet s_{(0,\infty)}(h_{\phi})
=\lambda\bullet\Theta^\mu(\phi), 
$$
where the first and last equalities follow from the definitions of 
$\Theta^\mu,\Theta^{\lambda\mu}$, the second equality follows from the proof of 
Lem. \ref{lem:15:14}, and  the third equality follows from the commutation of the action 
of $\mathbf k^\times$ with $s_{(0,\infty)}$. The combination of this identity with  
Lem. \ref{lem:14:10} implies that the pair
 $(\Theta\rtimes id,\sqcup_{\mu\in\mathbf k^\times}\Theta^\mu)$ 
is an involution of the torsor $(\mathcal G_{\mathrm{inert}}\rtimes
\mathbf k^\times,\sqcup_{\mu\in\mathbf k^\times}\mathcal G_{\mathrm{inert}}^\mu)$. 
The combination of Lem. \ref{lem:15:12:new} and Lem. \ref{lem:compduality} 
implies that the pair 
$(\sqcup_{\mu\in\mathbf k^\times}\Theta^\mu,\Theta^{\mathrm B})$ 
is an involution of the right torsor  
$(\sqcup_{\mu\in\mathbf k^\times}\mathcal G_{\mathrm{inert}}^\mu,
(\mathcal G^{\mathrm B}\rtimes \mathbf k^\times)_{\mathrm{inert}})$.
The statement follows. 
\end{proof}

\begin{thm}\label{thm:betti:b} (see Thm. \ref{thm:0:33})
The subgroup
$\mathsf{DMR}^{\mathrm B}(\mathbf k)$ of 
$(\mathcal G^{\mathrm B}\rtimes\mathbf k^\times)_{\mathrm{inert}}$
is invariant under the action of the involution $\Theta^{\mathrm B}$ of the latter group. 
\end{thm}

\begin{proof}
It follows from Cor. \ref{lem:15:16} and Lem. \ref{lem:15:15}(b) that 
 $(\mathcal G_{\mathrm{inert}}\rtimes\mathbf k^\times,\sqcup_{\mu\in\mathbf k^\times}
 \mathcal G_{\mathrm{inert}}^\mu,
(\mathcal G^{\mathrm B}\rtimes\mathbf k^\times)_{\mathrm{inert}})$ contains
two subbitorsors
$$
(\mathsf{DMR}_0(\mathbf k)\rtimes\mathbf k^\times,
\sqcup_{\mu\in\mathbf k^\times}\mathsf{DMR}_\mu(\mathbf k),
  \mathsf{DMR}^{\mathrm B}(\mathbf k))
$$
and
$$
((\Theta\rtimes id)(\mathsf{DMR}_0(\mathbf k)\rtimes\mathbf k^\times),
(\sqcup_{\mu\in\mathbf k^\times}\Theta^\mu)
(\sqcup_{\mu\in\mathbf k^\times}\mathsf{DMR}_\mu(\mathbf k)),
  \Theta^{\mathrm B}(\mathsf{DMR}^{\mathrm B}(\mathbf k))). 
$$
It follows from Cor. \ref{main:cor}(b) that the two subgroups 
$\mathsf{DMR}_0(\mathbf k)\rtimes\mathbf k^\times$ and 
$(\Theta\rtimes id)(\mathsf{DMR}_0(\mathbf k)\rtimes\mathbf k^\times)$ of 
$\mathcal G_{\mathrm{inert}}\rtimes\mathbf k^\times$ as equal, 
and from Thm. \ref{thm:DMRmu}(b) that the two subsets 
$\sqcup_{\mu\in\mathbf k^\times}\mathsf{DMR}_\mu(\mathbf k)$ and 
$(\sqcup_{\mu\in\mathbf k^\times}\Theta^\mu)
(\sqcup_{\mu\in\mathbf k^\times}\mathsf{DMR}_\mu(\mathbf k))$ of 
$\sqcup_{\mu\in\mathbf k^\times}
 \mathcal G_{\mathrm{inert}}^\mu$ are equal. 
Lem. \ref{lem:ef3} then implies that the two subgroups 
$\mathsf{DMR}^{\mathrm B}(\mathbf k)$ and 
$\Theta^{\mathrm B}(\mathsf{DMR}^{\mathrm B}(\mathbf k))$ of 
$(\mathcal G^{\mathrm B}\rtimes\mathbf k^\times)_{\mathrm{inert}}$ are equal. 
\end{proof}

\newpage

\part{Relationship with the Kashiwara-Vergne Lie algebra}\label{part:kv}

\Add{This part consists in only one section (\S\ref{sect:kv}), which is 
devoted to the proof of the relation between the Lie algebras $\mathfrak{ds}$ and $\mathfrak{krv}$.}

\section{Relationship with the Kashiwara-Vergne Lie algebra}\label{sect:kv}

\Add{In \S\ref{sect:16:1}, we prove the results relating  $\mathfrak{ds}$ and $\mathfrak{krv}$, 
announced in Lem. \ref{main:cor:c} and Thm. \ref{thm:schneps:nu}. In \S\ref{sect:comments:involutions}, 
we study the relation between the resulting sequence of 
inclusions of Lie algebras $\mathfrak{grt}_1\subset \mathfrak{dmr}_0\subset\mathfrak{krv}$ 
and the action of $\mathfrak S_3$ on the Lie algebra
$\underline{\mathfrak{ider}}_{\{0,1,\infty\}}/\underline{\mathfrak{inn}}_{\{0,1,\infty\}}$,
in which this inclusion takes place.}

\subsection{Relationship with the formalisms of \cite{FK,Sch1}}\label{sect:16:1}

Recall from \S\ref{sect:announcement:krv} the definition of
the Lie algebra $\mathbb L$, its subspace $\mathfrak{ds}$, its 
Lie algebras of derivations $\mathfrak{krv}\subset\mathfrak{sder}\subset\mathfrak{der}$, and the map $\nu : 
\mathbb L\to \mathfrak{der}$. 

\begin{defn} (see \cite{Sch1}) 
For any $u,v\in\mathbb L$ such that $[x,u]+[y,v]=0$, denote by $D_{v,u}\in\mathfrak{sder}$ 
the derivation such that $x\mapsto [x,u]$, $y\mapsto [y,v]$.  
\end{defn}


\begin{defn}\label{def:LL:inert}
    Set $\mathbb L_{\mathrm{inert}}:=\nu^{-1}(\mathfrak{sder})$, then $\mathbb L_{\mathrm{inert}}\subset \mathbb L$.  
\end{defn}

\begin{lem}\label{lem:sder}
Let $\alpha,\beta$ be the automorphisms of $\mathbb L$ given by 
$$
\alpha : x\mapsto x, y\mapsto -y\text{ and }
\beta : x\mapsto -x-y, y\mapsto y. 
$$

(a) One has 
\begin{equation}\label{eq:sder}
\mathbb L_{\mathrm{inert}}=\{\tilde f\in\mathbb L|\exists v\in \mathbb L, [x,\alpha(\tilde f)]+[-x-y,v]=0\}.      
\end{equation}

(b) The map $\nu : \mathbb L_{\mathrm{inert}}\to \mathfrak{sder}$ is given by 
$\tilde f\mapsto D_{\beta\alpha(\tilde f),\beta\alpha(\tilde f)-\beta(v)}$ 
for any $\tilde f\in \nu^{-1}(\mathfrak{sder})$, where $v$ is as in \eqref{eq:sder}. 
\end{lem}

\begin{proof}
(a) Let $\Phi,\Psi : \mathbb L^2\to\mathbb L$ be the maps respectively defined by 
$$
\Phi(f,v):=[x,f]+[-x-y,v],\quad \Psi(F,G):=[y,F]+[x,G]. 
$$
The map 
\begin{equation}\label{eq:Upsilon:03012025}
 \Upsilon : (F,G)\mapsto (f,v):=(\beta(F),\beta(F-G))   
\end{equation}
is a linear automorphism of $\mathbb L^2$, whose 
inverse is $(f,v)\mapsto (\beta(f),\beta(f-v))$. 
Then 
$$
(-\beta)\circ\Psi=\Phi\circ\Upsilon,  
$$
which implies that 
\begin{equation}\label{upsilon:bijection}
\text{$\Upsilon$ induces a bijection $\mathrm{ker}(\Psi)\to\mathrm{ker}(\Phi)$.} 
\end{equation}

Fix $\tilde f\in\mathbb L$. The elements $f,F\in\mathbb L$ defined as in (2.2) of \cite{FK} are given by 
the first component of \eqref{eq:Upsilon:03012025} and 
$f=\alpha(\tilde f)$.
Then
\begin{align*}
&
(\tilde f\in\nu^{-1}(\mathfrak{sder}))\iff(d^{\mathrm{FK}}_F\in\mathfrak{sder})\iff(
\exists G\in \mathbb L|[y,F]+[x,G]=0)
\iff(\exists G\in \mathbb L|(F,G)\in\mathrm{ker}(\Psi))
\\&
\iff(\exists G\in \mathbb L|\Upsilon(F,G)\in\mathrm{ker}(\Phi))
\iff(\exists G\in \mathbb L|(\beta(F),\beta(F-G))\in\mathrm{ker}(\Phi))
\\&
\iff(\exists v\in \mathbb L|(f,v)\in\mathrm{ker}(\Phi))
\iff(\exists v\in \mathbb L|[x,f]+[-x-y,v]=0)
\iff(\exists v\in \mathbb L|[x,\alpha(\tilde f)]+[-x-y,v]=0),
\end{align*}
where the first (resp. second, third, fifth, seventh) equivalence follows from the definition of $\nu$ (resp. $\mathfrak{sder}$, 
$\Psi$, $\Upsilon$, $\Phi$), the fourth (resp. sixth, last) equivalence follows from \eqref{upsilon:bijection} (resp. 
\eqref{eq:Upsilon:03012025}, $f=\alpha(\tilde f)$).  

(b) For $\tilde f\in  \mathbb L_{\mathrm{inert}}$, $\nu(\tilde f)=d_F^{\mathrm{FK}}
=d_{\beta\alpha(\tilde f)}^{\mathrm{FK}}=(y\mapsto[y,\beta\alpha(\tilde f)],x+y\mapsto 0)$. 
The relation $[x,\alpha(\tilde f)]+[-x-y,v]=0$ implies $[-x-y,\beta\alpha(\tilde f)]+[x,\beta(v)]=0$
hence $-[y,\beta\alpha(\tilde f)]=[x,\beta\alpha(\tilde f)-\beta(v)]$. It follows that 
$\nu(\tilde f)=(y\mapsto[y,\beta\alpha(\tilde f)],x\mapsto [x,\beta\alpha(\tilde f)-\beta(v)])
=D_{\beta\alpha(\tilde f),\beta\alpha(\tilde f)-\beta(v)}$. 
\end{proof}

\begin{lem}\label{lem08}
  The subspace $\mathbb L_{\mathrm{inert}}$ of $\mathbb L$ is graded, denote by $\mathbb L_{\mathrm{inert},>1}$
  the sum of its components of degree $>1$. Then $i$ induces an isomorphism $\mathbb L_{\mathrm{inert},>1}\to\mathfrak G_{\mathrm{inert}}$, therefore
  $$
\mathfrak G_{\mathrm{inert}}  =i(\mathbb L_{\mathrm{inert},>1}).
  $$ 
\end{lem}

\begin{proof}
The first statement is obvious. Let  $j : \mathbb L\to\mathfrak{lie}_{\{0,1\}}$ be the Lie algebra isomorphism 
induced by $x\mapsto e_0$, $y\mapsto e_1$, then  $j=i\circ\alpha$. It follows from \eqref{eq:sder} that 
$$
\mathbb L_{\mathrm{inert},>1}
=\{\tilde f\in\mathbb L_{>1}|\exists v\in \mathbb L_{>1}, [x,\alpha(\tilde f)]+[-x-y,v]=0\}
$$
where the index "$>1$" has the same meaning as above. Therefore 
$$
\alpha(\mathbb L_{\mathrm{inert},>1})
=\{f\in\mathbb L_{>1}|\exists v\in \mathbb L_{>1}, [x,f]+[-x-y,v]=0\},  
$$
which implies, as $j$ is a graded Lie algebra isomorphism inducing an isomorphism $\mathbb L_{>1}\to\mathfrak G$
and such that $x\mapsto e_0$, $-x-y\mapsto e_\infty$, 
$$
j(\alpha(\mathbb L_{\mathrm{inert},>1}))
=\{u\in\mathfrak G|\exists v'\in \mathfrak G, [e_0,u]+[e_\infty,v']=0\},  
$$
which implies the claim in view of $j=i\circ\alpha$ and Lem. \ref{lem06:1109}(c). 
\end{proof}









\Add{\begin{lem}
One has the inclusion
\begin{equation}\label{main:inc}
\mathfrak{ds}\subset\mathbb L_{\mathrm{inert}}^\wedge.   
\end{equation}
of subspaces of $\mathbb L^\wedge$.
\end{lem}
\begin{proof}
    \eqref{main:inc} follows
from Cor. \ref{main:cor}(c), using Lem. \ref{lem08} and Def. \ref{def:ds}. 
\end{proof}
}

\Add{\begin{lem}\label{main:cor:c:w:proof} (see Lem. \ref{main:cor:c})
If $\tilde f \in\mathfrak{ds}$, then the mould $M:=\mathrm{ma}_{\tilde f}$ 
defined in \cite{FK} satisfies the Ecalle senary relation (see (1) in \cite{FK}). 
\end{lem}
\begin{proof}
If $\tilde f\in\mathfrak{ds}$, then by Cor. \ref{main:cor}(c) 
one has $\tilde f\in\mathbb L_{\mathrm{inert}}$. By Lem. 
\ref{lem:sder}(b), 
it follows that $\nu(\tilde f) \in \mathfrak{sder}$. By the proof of Lem. \ref{lem:sder}(b), 
$\nu(\tilde f)=d^{\mathrm{FK}}_{\beta\alpha(\tilde f)}=d^{\mathrm{FK}}_F$. 
Therefore $d^{\mathrm{FK}}_F \in \mathfrak{sder}$. 
Prop. 1.2 in \cite{FK} then implies the conclusion. 
\end{proof}}

\Add{\begin{thm} \label{thm:schneps:nu:w:proof} (see Thm. \ref{thm:schneps:nu})
(see \cite{Sch1,Sch2}) The map $\nu : \mathbb L\to\mathfrak{der}$ induces an injection of Lie algebras $\mathfrak{ds}\hookrightarrow\mathfrak{krv}$. 
\end{thm}
\begin{proof}
Thm. 1.1 in \cite{Sch1} is based on a particular case of the statement (3.64) in \cite{Ec2}, which is 
equivalent to Lem. \ref{main:cor:c}; this statement is therefore valid. 
It says that some linear map is an injection of graded Lie algebras
$\mathfrak{ds}\hookrightarrow\mathfrak{krv}$. This linear map 
is $\mathfrak{ds}\ni\tilde f\mapsto D_{F,s(F^x)}=(x\mapsto[x,s(F^x)],y\mapsto[y,F])
\in\mathfrak{krv}$, where $f(x,y):=\tilde f(x,-y)$, $F(x,y):=f(-x-y,y)$, $F^x,F^y\in\mathbb Q\langle\langle x,y\rangle\rangle$ are
such that $F=xF^x+yF^y$, and $s$ is given by (1.4) of {\it loc. cit.} Since $\mathfrak{krv}\subset\mathfrak{sder}$, 
$D_{F,s(F^x)}(x+y)=0$, therefore $D_{F,s(F^x)}=d_F^{\mathrm{FK}}$. It follows that the injection from \cite{Sch1}
is induced by $\nu$.
\end{proof}}

\subsection{Relationship with involutions}\label{sect:comments:involutions}

For $S$ a set, denote by $\mathfrak{lie}_S$ the free $\mathbf k$-Lie algebra over $(e_s)_{s\in S}$, 
graded by the condition that all the generators are of degree 1. Define $\underline{\mathfrak{lie}}_S
:=\mathfrak{lie}_S/(\sum_{s\in S}e_s)$.  
Let $\underline{\mathfrak{ider}}_S$ be the set of derivations $D$ of $\underline{\mathfrak{lie}}_S$ such that 
there exists $(a_s)_{s\in S}$ such that $D(e_s)=[a_s,e_s]$ for any $s\in S$. Let $\underline{\mathfrak{inn}}_S$
be the set of inner derivations of  $\underline{\mathfrak{lie}}_S$, then $\underline{\mathfrak{inn}}_S$
is an ideal of $\underline{\mathfrak{ider}}_S$, which gives rise to the quotient Lie algebra
$$
\underline{\mathfrak{ider}}_S/\underline{\mathfrak{inn}}_S
$$
which is equipped with an action of $\mathfrak S_S$.

For $T$ a set, define $\mathfrak{ider}_T$ to be the set of derivations $D$ of $\mathfrak{lie}_T$ such that 
there exists $(b_t)_{t\in T\cup\infty}$ such that $D(e_t)=[e_t,b_t]$ for any $t\in  T$
and $D(\sum_{t\in T}e_t)=[\sum_{t\in T}e_t,b_\infty]$.
Let $\mathfrak{sder}_T\subset \mathfrak{ider}_T$ be the subset defined by the condition $D(\sum_{t\in T}e_t)=0$. 
Then the set $\mathfrak{inn}_T$ of inner derivations of $\mathfrak{lie}_T$ 
is an ideal of $\mathfrak{ider}_T$, and the map 
$\mathfrak{sder}_T\to\mathfrak{ider}_T/\mathfrak{inn}_T$ is an isomorphism in degree $>1$. 

Then $S$ a set and $s_0\in S$, the natural map 
$\mathfrak{lie}_{S\smallsetminus s_0}\to \underline{\mathfrak{lie}}_S$ is a Lie algebra isomorphism, which 
gives rise to an isomorphism $\mathfrak{ider}_{S\smallsetminus s_0}/
\mathfrak{inn}_{S\smallsetminus s_0}\to\underline{\mathfrak{ider}}_S/\underline{\mathfrak{inn}}_S$. 
Combining this with the morphism  
$\mathfrak{sder}_{S\smallsetminus s_0}\to\mathfrak{ider}_{S\smallsetminus s_0}/\mathfrak{inn}_{S\smallsetminus s_0}$, 
one obtains the Lie algebra morphism  
$$
\mathfrak{sder}_{S\smallsetminus s_0}\to 
\underline{\mathfrak{ider}}_S/\underline{\mathfrak{inn}}_S
$$
which is an isomorphism in degree $>1$. For $S:=\{0,1,\infty\}$, this leads to the diagram of Lie algebra morphisms
$$
\mathfrak{sder}_{\{0,\infty\}}\to\underline{\mathfrak{ider}}_{\{0,1,\infty\}}/\underline{\mathfrak{inn}}_{\{0,1,\infty\}}\leftarrow 
\mathfrak{sder}_{\{1,\infty\}}, 
$$
which are isomorphisms in degree $>1$.

\Add{
\begin{lem}\label{lem04:2506:ex:d}
(a) The Lie algebra of the group functor $\mathbf k\mapsto\mathrm{IAut}_{\{0,1\}}^0$ 
is the degree completion of $\mathfrak{ider}_{\{0,1\}}^0:=\{D\in 
\mathfrak{ider}_{\{0,1\}}|D(e_1)=0$ and $D(e_0)$ has degree $\geq 3\}$. 
\\
(b) The Lie algebra isomorphism associated to the isomorphism from Lem. \ref{lem04:2506}(a)
is $\mathfrak G_{\mathrm{inert}}\to \mathfrak{ider}_{\{0,1\}}^0$ given by $a\mapsto D_a$, where 
$D_a$ is given by $e_0\mapsto [a,e_0]$ and $e_1\mapsto 0$. 
\\
(c) One has the following inclusion $\mathfrak{ider}_{\{0,1\}}^0\subset\mathfrak{sder}_{\{0,\infty\}}$
of the Lie algebra of derivations $\mathrm{Der}(\mathfrak{lie}_{\{0,1\}})$ of $\mathfrak{lie}_{\{0,1\}}$; 
it is graded and is an isomorphism in degrees $>1$. 
\\
(d) The Lie algebra morphism $\mathfrak G_{\mathrm{inert}}\to\mathfrak{sder}_{\{0,\infty\}}$
defined as the composition of the morphisms from (b) and (c) is graded and
is an isomorphism in degrees $>1$. 
\end{lem}
\begin{proof}
 (a) follows from the definition of $\mathbf k\mapsto\mathrm{IAut}_{\{0,1\}}^0$.  
 (b) follows from the definition of the isomorphism from Lem. \ref{lem04:2506}(a). 
 (c) follows from definition of $\mathfrak{ider}^0_{\{0,1\}}$ and from 
 the equality of $\mathfrak{sder}_{\{0,\infty\}}$
with the set of derivations of $\mathfrak{lie}_{\{0,1\}}$ of the form 
$e_0\mapsto [a,e_0]$, $e_\infty\mapsto[b,e_\infty]$ for some $a,b\in \mathfrak{lie}_{\{0,1\}}$, 
and such that $e_1=-e_0-e_\infty\mapsto 0$. (d) follows from (b) and (c). 
\end{proof}
}

\begin{lem}
(a) The Lie algebra isomorphism $i(\beta\alpha)^{-1} : \mathbb L\to \underline{\mathfrak{lie}}_{\{0,1,\infty\}}$ 
induces an isomorphism $D\mapsto i(\beta\alpha)^{-1}D\beta\alpha i^{-1}$ 
between the Lie algebra of derivations of its two sides, which 
restricts to an isomorphism $\mathfrak{sder}\to\mathfrak{sder}_{\{1,\infty\}}$. 

(b) The quadrilaterals of the following diagram of Lie algebra morphisms are commutative
(the indicated maps are isomorphisms in degree $>1$)
\begin{equation}\label{diag:1er:mai}
\xymatrix{
\mathfrak{sder}_{\{0,\infty\}}\ar^{\!\!\!\!\!\!\!\!\sim}
[r]&
\underline{\mathfrak{ider}}_{\{0,1,\infty\}}/\underline{\mathfrak{inn}}_{\{0,1,\infty\}}&
\ar_\sim[l]\mathfrak{sder}_{\{1,\infty\}}\\
\mathfrak G_{\mathrm{inert}}\ar^{\mathrm{Lem. \ref{lem04:2506:ex:d}}}_\sim[u]&
\ar_{i\sim}^{\mathrm{Lem. \ref{lem08}}}[l]\mathbb L_{\mathrm{inert},>1}
\ar^{\nu}_{\mathrm{Def. \ref{def:LL:inert}}}[r]&\mathfrak{sder}\ar_{D\mapsto i(\beta\alpha)^{-1}D\beta\alpha i^{-1}}^\sim[u]\\
\mathfrak{dmr}_0\ar@{^{(}->}[u]&\ar_{i\sim}[l]\ar@{^{(}->}[u]\mathfrak{ds}
\ar^{\nu}[r]&\ar@{^{(}->}[u]\mathfrak{krv}
}
\end{equation}
\end{lem}

\begin{proof}
(a) follows from the definitions of $\mathfrak{sder}$ and $\mathfrak{sder}_{\{1,\infty\}}$, 
and from the fact that $i\beta^{-1}\alpha^{-1}$ maps $-x-y,y,x$ to $e_0,e_1,e_\infty$. 

(b) The commutativities of the left and right squares respectively follow from 
Def. \ref{def:ds} and Thm. \ref{thm:schneps:nu}. Let us show the commutativity of the top quadrilateral.  
The image of $\tilde f\in \mathbb L_{\mathrm{inert},>1}$ by the left composition is: 
\begin{equation}\label{pblv}
\tilde f\mapsto i(\tilde f)\mapsto (e_0\mapsto [i(\tilde f),e_0],e_\infty\mapsto [b_{i(\tilde f)},e_\infty],e_1\mapsto 0)
\in\mathfrak{sder}_{\{0,\infty\}},     
\end{equation}
while its image by the right composition is
\begin{align*}
    & \tilde f\stackrel{\nu}{\mapsto} D_{\beta\alpha(\tilde f),\beta\alpha(\tilde f)-\beta(v)}=(y\mapsto [y,\beta\alpha(\tilde f)],
x\mapsto [x,\beta\alpha(\tilde f)-\beta(v)],-x-y\mapsto 0)
    \\&
 \mapsto (e_1\mapsto[e_1,i(\tilde f)],e_\infty\mapsto[e_\infty,i(\tilde f)-i\alpha^{-1}(v)],e_0\mapsto 0)
 \in\mathfrak{sder}_{\{0,\infty\}}\hookrightarrow\underline{\mathfrak{ider}}_{\{0,1,\infty\}}, 
\end{align*}
where the computation of the first (resp. second) image follows from Lem. \ref{lem:sder}(b)
(resp. the fact that $i\beta^{-1}\alpha^{-1}$ maps $-x-y,y,x$ to $e_0,e_1,e_\infty$). 
The second image is equivalent modulo $\underline{\mathfrak{inn}}_{\{0,1,\infty\}}$ to 
the derivation 
$$
(e_1\mapsto 0,e_\infty\mapsto[e_\infty,-i\alpha^{-1}(v)],e_0\mapsto [e_0,-i(\tilde f)])
$$
(as both derivations differ by the inner derivation associated to $i(\tilde f)$), which coincides with the image of
\eqref{pblv}. 
\end{proof}

\begin{rem}
The Lie algebras  $\mathfrak{grt}_1$, $\mathfrak{dmr}_0$ and $\mathfrak{krv}$ may be viewed, via \eqref{diag:1er:mai}, as 
Lie subalgebras of $\underline{\mathfrak{ider}}_{\{0,1,\infty\}}/\underline{\mathfrak{inn}}_{\{0,1,\infty\}}$, 
satisfying the sequence of inclusions (the two first of which being equalities according to the conjecture of \cite{AT}) 
$$
\mathfrak{grt}_1\subset \mathfrak{dmr}_0\subset \mathfrak{krv}
\subset \underline{\mathfrak{ider}}_{\{0,1,\infty\}}/\underline{\mathfrak{inn}}_{\{0,1,\infty\}}. 
$$
The permutation group $\mathfrak S_3$ of the set $\{0,1,\infty\}$ acts on the Lie algebra
$\underline{\mathfrak{ider}}_{\{0,1,\infty\}}/\underline{\mathfrak{inn}}_{\{0,1,\infty\}}$, 
and this action preserves the Lie subalgebra $\mathfrak{krv}$ (personal communication of A. Alekseev). 
The involution $\tau$ of $\mathfrak{krv}$ from \cite{AT} a few lines before Prop. 8.7, which is induced by the exchange of 
$x$ and $y$, is  
induced by above inclusion and by the involution $s_{(1,\infty)}$ of $\underline{\mathfrak{ider}}_{\{0,1,\infty\}}/
\underline{\mathfrak{inn}}_{\{0,1,\infty\}}$. (The corresponding Lie subalgebra of fixed points, 
denoted $\mathfrak{krv}^{\mathrm{sym}}$, was recently studied in \cite{Ku}.) 
On the other hand,  by Cor. \ref{main:cor}(b), $\Theta$ from Lem. \ref{lem02:2506}(d) is 
an involution of $\mathfrak{dmr}_0$; this inclusion is induced by the above inclusion and by the involution 
$s_{(0,\infty)}$ of $\underline{\mathfrak{ider}}_{\{0,1,\infty\}}/\underline{\mathfrak{inn}}_{\{0,1,\infty\}}$; 
we do not know how to prove that $\mathfrak{dmr}_0$ is preserved by 
the action of $\mathfrak S_3$. 
Finally, $\mathfrak{grt}_1$ is contained in the Lie subalgebra of 
$\underline{\mathfrak{ider}}_{\{0,1,\infty\}}/\underline{\mathfrak{inn}}_{\{0,1,\infty\}}$ of invariants of the 
action of $\mathfrak S_3$, which is called in \cite{Dr} the Ihara Lie algebra. 
\end{rem}

\newpage
\appendix

\part*{Appendices}

The purpose of Appendices A and B is to give a proof of Prop. \ref{lemma:A:1403}, which is essentially equivalent to 
Prop. 2.2 in \cite{Sch1},  alternative to that of {\it loc. cit.} The objective of Appendix C is the proof of 
Prop. \ref{prop:Z}, which is a key ingredient in the proof of the inclusion
$\mathrm{Stab}_{\mathcal G}(\mathrm{GL}_3\hat V\bullet\rho_{\mathrm{DT}})\subset
\mathcal G_{\mathrm{inert}}$ (see Thm. \ref{thm:13:22:0205}).

\section{Poisson and Lie algebras}\label{app:poisson}

In this section, the base ring is the field $\mathbb Q$. 

\begin{defn}
(a) A Poisson algebra is the data of a commutative algebra $C$, together with a map 
$C^{\otimes 2}\to C$, $f\otimes g\mapsto \{f,g\}$, such that $(C,\{-,-\})$ is a Lie algebra, 
and that for any $g\in C$, the endomorphism $\{-,g\}$ of $C$ is a derivation of $C$. 

(b) A commutative algebra $C$ being given, a map $C^{\otimes 2}\to C$ satisfying these conditions
is called a Poisson bracket on $C$. 
\end{defn} 

The following result is well-known. 
\begin{lem}\label{lem:Poisson}
Let $A$ be an associative algebra, equipped with an increasing filtration $F_0A\subset F_1A\subset\cdots$, such that 
the associated graded algebra $\mathrm{gr}A$ is commutative. For $k,l\geq 0$ and $a,b$ in 
$\mathrm{gr}A$ of degrees $k,l$, one has $\tilde a\tilde b-\tilde b\tilde a\in F_{k+l-1}A$ for any 
lifts $\tilde a\in F_kA$, $\tilde b\in F_lB$, and the class of $\tilde a\tilde b-\tilde b\tilde a$
in $\mathrm{gr}_{k+l-1}A$ is independent of the choice of $\tilde a,\tilde b$, and denoted $\{a,b\}$.
The resulting map $F_kA\times F_lA\to F_{k+l-1}A$ is bilinear, and extends to a Poisson bracket on $\mathrm{gr}A$. 
\end{lem}

\begin{lem}\label{A3} (see \cite{Mei}, Thm 5.2 and Rem. 5.3)
Let $\mathfrak g$ be a $\mathbb Q$-Lie algebra. 

(a) The symmetric algebra $S(\mathfrak g)$ over $\mathfrak g$ is equipped with a Poisson bracket $\{-,-\}$, uniquely 
defined by the condition 
that the injection $\mathfrak g\hookrightarrow (S(\mathfrak g),\{-,-\})$ is a Lie algebra morphism. 

(b) $F_0U\mathfrak g=\mathbb Q$, $F_{k+1}U\mathfrak g=(\mathbb Q+\mathfrak g)\cdot F_kU\mathfrak g$ 
for $k\geq 0$ defines an increasing filtration of 
 $U(\mathfrak g)$, and there is an isomorphism $\mathrm{gr}U(\mathfrak g)\simeq S(\mathfrak g)$
 of graded algebras; in particular $\mathrm{gr}U(\mathfrak g)$ is commutative. 

 (c) The Poisson algebra structure on $S(\mathfrak g)$ associated by Lem. \ref{lem:Poisson} to the algebra filtration from (b) 
  is that of (a).    
\end{lem}

\section{Exactness of some complexes}

In this section, the base ring is the field $\mathbb Q$, and we will denote  $\mathfrak{lie}_{\{x,y\},\mathbb Q}$
by  $\mathfrak{lie}_{\{x,y\}}$. Set $L:=\mathbb Qy \oplus \mathfrak{lie}_{\{x,y\}}[\geq 2] \subset \mathfrak{lie}_{\{x,y\}}$. 
Then $S(L)\subset S(\mathfrak{lie}_{\{x,y\}})$ is a commutative subalgebra.

By Lem. \ref{A3}(a), the Lie bracket of $\mathfrak{lie}_{\{x,y\}}$ induces a Poisson structure on $S(\mathfrak{lie}_{\{x,y\}})$. 

\begin{lem}\label{lem:3:1603}
(a) $S(L)$ is stable under the endomorphism $\{x,-\}$ of 
$S(\mathfrak{lie}_{\{x,y\}})$. 

(b) The kernel of the endomorphism $f\mapsto \{x,f\}$ of $S(L)$ is $\mathbb Q$.
\end{lem}

\begin{proof}
(a) $\{x,-\}$ is the derivation of $S(\mathfrak{lie}_{\{x,y\}})$ (viewed as a commutative algebra) 
which extends the endomorphism $[x,-]$ of $\mathfrak{lie}_{\{x,y\}}$. 
The subspace $L\subset\mathfrak{lie}_{\{x,y\}}$ is stable under this endomorphism, 
therefore $\{x,-\}$ leaves stable the commutative subalgebra of $S(\mathfrak{lie}_{\{x,y\}})$
generated by $L$, which is $S(L)$. 

(b) $L$ is graded with finite dimensional components, and is stable under the endomorphism $[x,-]$ of 
$\mathfrak{lie}_{\{x,y\}}$, which is injective and of degree 1. Fix a graded complement $B$ of $[x,L]$ in $L$. 
Then the map 
\begin{equation}\label{iso:2115:2103}
\mathbb Q[t] \otimes B\to L
\end{equation}
induced by $t^d\otimes b\mapsto\mathrm{ad}_x^d(b)$ 
for $b\in B$, $d\geq 0$ is a morphism of graded $\mathbb Q$-vector spaces, which intertwines the endomorphisms $(t\cdot -)\otimes id_B$ 
and $[x,-]$. 

Let us prove that the map $\mathbb Q[t] \otimes B\to L$ is an isomorphism. Let $L=\oplus_{k\geq1}L_k$, $B=\oplus_{k\geq1}B_k$
be the homogeneous decompositions of $L,B$. Let us prove by induction on $k\geq1$ that $L_k=\sum_{0<l\leq k}
\mathrm{ad}_x^{k-l}(B_l)$ for any $k\geq 1$. One has obviously $L_1=B_1$, and the equality at step $k$ together with
$ L_{k+1}=B_{k+1}+[x,L_k]$ implies the equality at step $k+1$, which establishes the induction. This implies the 
surjectivity of $\mathbb Q[t] \otimes B\to L$. Let us prove its injectivity. An element of the kernel of this map 
has a unique expression as $\sum_{i\geq 0}t^i\otimes b_i$, where $b_i\in B$. Then $\sum_{i\geq 0}\mathrm{ad}_x^i(b_i)=0$
(equality in $L$). Let show inductively on $i_0\geq 0$ that $b_{i_0}=0$. One has $b_0=-\sum_{i>0}\mathrm{ad}_x^i(b_i)\in [x,L]$, 
hence $v_0\in B\cap [x,L]=0$ so $b_0=0$. Assume that $b_0=\ldots=b_{i_0-1}=0$. Then $\sum_{i\geq i_0}\mathrm{ad}_x^i(b_i)=0$. 
Since the restriction of $\mathrm{ad}_x=[x,-]$ to $L$ is injective, one obtains 
$\sum_{i\geq i_0}\mathrm{ad}_x^{i-i_0}(b_i)=0$,  which implies $b_{i_0}\in [x,L]$, hence 
$b_{i_0}\in [x,L]\cap B=0$, so $b_{i_0}=0$, establishing the induction, and therefore the injectivity of 
$\mathbb Q[t] \otimes B\to L$.

Let now $\mathcal B$ be a graded basis of $B$; it induces an isomorphism $\mathbb Q\mathcal B\to B$, which 
combined with the isomorphism \eqref{iso:2115:2103} gives rise to an isomorphism  
$\oplus_{b\in \mathcal B}\mathbb Q[t]\stackrel{\sim}{\to} L$, which induces an isomorphism 
\begin{equation}\label{iso:1:2103}
S(\oplus_{b\in \mathcal B}\mathbb Q[t])\stackrel{\sim}{\to} S(L), 
\end{equation}
which intertwines the derivation of $S(\oplus_{b\in \mathcal B}\mathbb Q[t])$ induced by the 
endomorphism of $\oplus_{b\in \mathcal B}\mathbb Q[t]$ given by multiplication of $t$, with $\{x,-\}$. 

The derivation $T$ of $S(\mathbb Q[t])$ induced by the endomorphism of $\mathbb Q[t]$ given by multiplication of $t$,
is compatible with the decomposition $S(\mathbb Q[t])=\oplus_{d\geq 0}S^d(\mathbb Q[t])$, denote by 
$T=\oplus_{d\geq 0}T^{(d)}$ the corresponding decomposition of $T$; one has $T^{(0)}=0$; moreover, 
one has an isomorphism $S^d(\mathbb Q[t])\simeq \mathbb Q[t_1,\dots,t_d]^{\mathfrak S_d}$, which intertwines 
$T^{(d)}$ with multiplication by $t_1+\cdots+t_d$. 

There is an isomorphism 
\begin{equation}\label{iso:2:2103}
\oplus_{\mathbf d\in \mathbb Z_{\geq 0}^{(\mathcal B)}}\otimes_{b\in \mathcal B}S^{\mathbf d(b)}(\mathbb Q[t])
\simeq S(\oplus_{b\in \mathcal B}\mathbb Q[t]),
\end{equation}
which interwines $\oplus_{\mathbf d\in \mathbb Z_{\geq 0}^{(\mathcal B)}}\sum_{b\in\mathcal B}T^{(\mathbf d(b))}$
with the derivation of $S(\oplus_{b\in \mathcal B}\mathbb Q[t])$ induced by the endomorphism of 
$\oplus_{b\in \mathcal B}\mathbb Q[t]$ given by multiplication of $t$. 

For each $\mathbf d\in \mathbb Z_{\geq 0}^{(\mathcal B)}$, the space $\otimes_{b\in\mathcal B}S^{\mathbf d(b)}(\mathbb Q[t])$
is isomorphic to the algebra of polynomials in the set of variables 
$t(b,i)$ indexed by the pairs $(b,i)$, where $b\in\mathcal B,i\in[\!\![1,\mathbf d(b)]\!\!]$, invariant under the action of 
$\prod_{b\in \mathcal B}\mathfrak S_{\mathbf d(v)}$, and the isomorphism intertwines $\otimes_{b\in\mathcal B}S^{\mathbf d(b)}(t\cdot -)$
with multiplication by $\sum_{b\in\mathcal B,i\in[\!\![1,\mathbf d(b)]\!\!]}t(b,i)$. Since polynomials form a domain, 
this multiplication is injective, which implies that $\oplus_{\mathbf d\in \mathbb Z_{\geq 0}^{(\mathcal B)}}
\otimes_{b\in \mathcal B}S^{\mathbf d(b)}(t\cdot -)$ is injective when $\mathbf d\neq 0$; on the other hand, this endomorphism 
is 0 if $\mathbf d=0$. As $\{x,-\}$ is the conjugation of the direct sum over $\mathbf d$ of the endomorphism indexed by 
$\mathbf d$ by the composition of \eqref{iso:1:2103} and \eqref{iso:2:2103}, the kernel of $\{x,-\}$ is 
to the image of $\otimes_{b\in \mathcal B}S^{\mathbf 0(b)}(\mathbb Q[t])$, which is $\mathbb Q$. 
\end{proof}

\begin{rem}
Let us sketch an alternative proof of Lem. \ref{lem:3:1603}(b), based on the following statements in the context of Lem. \ref{A3}: 
(a) the symmetrization map $S(\mathfrak g)\to U(\mathfrak g)$, defined to be the linear map such that 
$S^k(\mathfrak g)\ni x^k\mapsto x^k\in U(\mathfrak g)$ for any $k\geq 0$ and $x\in \mathfrak g$, is an linear 
isomorphism (PBW theorem); (b) a Lie algebra derivation $D$ of $\mathfrak g$ admits unique extensions both  
to a Poisson algebra derivation $D_{S(\mathfrak g)}$ of $S(\mathfrak g)$, and to an associative algebra 
derivation $D_{U(\mathfrak g)}$ of $U(\mathfrak g)$, which are moreover
intertwined by the symmetrization map. 

Recall that $[x,-]$ is a derivation of $L$. Then the endomorphism $\{x,-\}$ of $S(L)$ from the statement of 
Lem. \ref{lem:3:1603}(b) coincides with $[x,-]_{S(L)}$. By the above statements, the tensor product to the 
symmetrization maps of $\mathbb Qx$ and of $L$ induces an isomorphism 
$\mathbb Q[x] \otimes S(L)\to U(\mathbb Qx) \otimes U(L)$, which intertwines 
$id_{\mathbb Q[x]}\otimes [x,-]_{S(L)}$ and $id_{U(\mathbb Qx)}\otimes [x,-]_{U(L)}$. 
Since $\mathbb Qx\oplus L=\mathfrak{lie}_{\{x,y\}}$, the composition of product with the tensor 
product of natural inclusions induces an linear isomorphism $U(\mathbb Qx) \otimes U(L)\to 
U(\mathfrak{lie}_{\{x,y\}})$, which since $x$ commutes with $U(\mathbb Qx)$ intertwines $id_{U(\mathbb Qx)}\otimes [x,-]_{U(L)}$
with the inner derivation $[x,-]$ of $U(\mathfrak{lie}_{\{x,y\}})$. Finally,  
$U(\mathfrak{lie}_{\{x,y\}})$ is equal to the free algebra $\mathbb Q\langle x,y\rangle$ generated by $x,y$. 
Then $\mathrm{ker}([x,-])=\mathbb Q[x]$. It follows that the kernel of $id_{\mathbb Q[x]}\otimes \{x,-\}$, 
which is $\mathbb Q[x]\otimes\mathrm{ker}(\{x,-\})$
is the preimage of $\mathbb Q[x]$ in $\mathbb Q[x] \otimes S(L)$, which is $\mathbb Q[x]\otimes\mathbb Q$. 
This implies the statement of Lem. \ref{lem:3:1603}(b). 
\end{rem}

Denote by $S(\mathfrak{lie}_{\{x,y\}})_+$ the kernel of the projection map $S(\mathfrak{lie}_{\{x,y\}})
\to S^0(\mathfrak{lie}_{\{x,y\}})=\mathbb Q$. 

\begin{lem}\label{lem:poisson:1703}
(a) The sequence of maps $S(\mathfrak{lie}_{\{x,y\}})\to S(\mathfrak{lie}_{\{x,y\}})_+\oplus S(\mathfrak{lie}_{\{x,y\}})
\to S(\mathfrak{lie}_{\{x,y\}})$, where the first map is 
$c\mapsto (c\cdot x,\{c,x\})$ and the second map is $(a,b)\mapsto \{a,x\}-b\cdot x$, is an exact complex.  

(b) For each $k \geq 0$, the sequence of maps 
$S^k(\mathfrak{lie}_{\{x,y\}})\to S^{k+1}(\mathfrak{lie}_{\{x,y\}})\oplus S^k(\mathfrak{lie}_{\{x,y\}})
\to S(\mathfrak{lie}_{\{x,y\}})$ is an exact subcomplex of the complex from (a). 
\end{lem}

\begin{proof}
Since $L\oplus\mathbb Qx=\mathfrak{lie}_{\{x,y\}}$, the product induces an isomorphism 
$S(\mathbb Qx) \otimes S(L) \simeq S(\mathfrak{lie}_{\{x,y\}})$, which implies the direct sum 
decomposition $S(\mathfrak{lie}_{\{x,y\}}) \simeq S(L) \oplus x \cdot S(\mathfrak{lie}_{\{x,y\}})$, which in its turn
implies the direct sum decomposition
\begin{equation}\label{DIRECT:SUM:1603}
S(\mathfrak{lie}_{\{x,y\}})_+ \simeq S(L)_+ \oplus x \cdot S(\mathfrak{lie}_{\{x,y\}}),   
\end{equation}
$S(L)_+$ being the kernel of the projection $S(L)\to\mathbb Q$. 

(a) The said sequence of maps clearly forms a complex. Let  $a\in S(\mathfrak{lie}_{\{x,y\}})_+$, 
$b \in S(\mathfrak{lie}_{\{x,y\}})$ be such that $\{a,x\}=b\cdot x$. By \eqref{DIRECT:SUM:1603}, 
there exists $a_+ \in S(L)_+$ and $c\in S(\mathfrak{lie}_{\{x,y\}})$ such that $a=a_++x\cdot c$. Then 
one computes 
$$
\{a_+,x\}=\{a-x\cdot c,x\}
=\{a,x\}-x\cdot\{c,x\}-\{x,x\}\cdot c=b\cdot x-\{c,x\}\cdot x=(b-\{c,x\})\cdot x, 
$$
where the second equality follows from the Leibniz rule and the third equality follows from 
the assumption, the antisymmetry of the Poisson bracket and the commutativity of the symmetric algebra.  
This implies $\{a_+,x\} \in S(\mathfrak{lie}_{\{x,y\}})\cdot x$. 
On the other hand, since $L$ is stable under $[x,-]$, $S(L)_+$ is stable under $\{x,-\}$, therefore 
$\{a_+,x\} \in S(L)_+$. Then \eqref{DIRECT:SUM:1603} implies $\{a_+,x\}=0$, therefore by Lem. \ref{lem:3:1603}(b), 
$a_+ \in \mathbb Q$, which since $a_+\in S(L)_+$ implies $a_+=0$, therefore $a=c\cdot x$. This equality implies 
the second equality in $b\cdot x=\{a,x\}=\{c,x\}\cdot x$, which since $S(\mathfrak{lie}_{\{x,y\}})$ is a domain 
implies $b=\{c,x\}$. Therefore $(a,b)=(c\cdot x,\{c,x\})$. 

(b) follows from (a) and from the fact that the first map in the complex from (a) is the direct sum over $k\geq0$
of the maps from the statement of (b). 
\end{proof}

Denote by $\mathbb Q\langle x,y\rangle$ the free associative $\mathbb Q$-algebra with generators $x,y$,  
where $x,y$ have degree 1; it is equipped with a Hopf algebra structure for which $x,y$ are primitive. Then
$\mathcal P\mathbb Q\langle x,y\rangle=\mathfrak{lie}_{\{x,y\}}$. 
For $n\geq0$, denote by $\mathbb Q\langle x,y\rangle[n]$, $\mathfrak{lie}_{\{x,y\}}[n]$ the degree $n$ parts of 
$\mathbb Q\langle x,y\rangle$, $\mathfrak{lie}_{\{x,y\}}$. Denote by $\mathbb Q\langle x,y\rangle_+$ the kernel of the 
augmentation map $\mathbb Q\langle x,y\rangle\to\mathbb Q$. 

\begin{lem}\label{lemma:B:1403}
The sequence of maps 
$\mathbb Q\langle x,y\rangle\oplus \mathfrak{lie}_{\{x,y\}}\to
\mathbb Q\langle x,y\rangle_+^{\oplus2}\oplus\mathfrak{lie}_{\{x,y\}}\to \mathbb Q\langle x,y\rangle$, 
where the first map is $(c,u)\mapsto (xc+u,cx+u,[u,x])$ 
and the second map is $(a,b,z)\mapsto ax-xb-z$, is an exact complex of $\mathbb Q$-vector spaces.    
\end{lem}

\begin{proof} Let us denote by $A\to B\to \mathbb Q\langle x,y\rangle$ the said sequence of maps. It is clearly a complex of 
$\mathbb Q$-vector spaces, let us prove its exactness. 

Recall from Appendix \ref{app:poisson} that $\mathbb Q\langle x,y\rangle$ is equipped with an algebra filtration given by 
$F_0\mathbb Q\langle x,y\rangle=\mathbb Q$ and $F_k\mathbb Q\langle x,y\rangle=(\mathbb Q+\mathfrak{lie}_{\{x,y\}})
\cdot F_{k-1}\mathbb Q\langle x,y\rangle$ for $k>0$, inducing 
a graded algebra isomorphism $\mathrm{gr}\mathbb Q\langle x,y\rangle\simeq S(\mathfrak{lie}_{\{x,y\}})$, in particular
one has $\mathrm{gr}_k\mathbb Q\langle x,y\rangle\simeq S^k(\mathfrak{lie}_{\{x,y\}})$ for any $k\geq0$. 
The induced filtration on $\mathbb Q\langle x,y\rangle_+$ is defined by $F_k\mathbb Q\langle x,y\rangle_+:=
F_k\mathbb Q\langle x,y\rangle\cap\mathbb Q\langle x,y\rangle_+$. 
We will denote by $x\mapsto [x]_k$  the projection map $F_k\mathbb Q\langle x,y\rangle\to 
S^k(\mathfrak{lie}_{\{x,y\}})$. The algebra $S(\mathfrak{lie}_{\{x,y\}})$ is equipped with a Poisson 
algebra structure of degree $-1$ denoted $\{-,-\}$: if $a\in F_k\mathbb Q\langle x,y\rangle$ and 
$b\in F_l\mathbb Q\langle x,y\rangle$, then 
$a\cdot b\in F_{k+l}\mathbb Q\langle x,y\rangle$ and $[a\cdot b]_{k+l}=[a]_k\cdot[b]_l$, moreover
$a\cdot b-b\cdot a\in F_{k+l-1}\mathbb Q\langle x,y\rangle$ and $[a\cdot b-b\cdot a]_{k+l-1}
=\{[a]_k,[b]_l\}$. 

For $k\geq 0$, set $F_kA:=F_k\mathbb Q\langle x,y\rangle\oplus \mathfrak{lie}_{\{x,y\}}$ and 
$F_kB:=F_{k+1}\mathbb Q\langle x,y\rangle_+^{\oplus2}\oplus\mathfrak{lie}_{\{x,y\}}$. 
The complex $A\to B\to \mathbb Q\langle x,y\rangle$ induces a complex $F_kA\to F_kB\to \mathbb Q\langle x,y\rangle$. 

Let us prove inductively on $k\geq0$ the exactness of $F_kA\to F_kB\to \mathbb Q\langle x,y\rangle$. 

For $k=0$, this complex is 
$\mathbb Q\oplus \mathfrak{lie}_{\{x,y\}}\to 
\mathfrak{lie}_{\{x,y\}}^{\oplus 2}\oplus\mathfrak{lie}_{\{x,y\}}\to \mathbb Q\langle x,y\rangle$.  
Let $(a,b,z)\in\mathrm{ker}(\mathfrak{lie}_{\{x,y\}}^{\oplus 2}\oplus\mathfrak{lie}_{\{x,y\}}\to \mathbb Q\langle x,y\rangle)$. 
Then $ax-xb=z$, which implies that the image of the left-hand side under the projection map
$F_2\mathbb Q\langle x,y\rangle\to\mathrm{gr}_2\mathbb Q\langle x,y\rangle=S^2(\mathfrak{lie}_{\{x,y\}})$ is zero. 
Since this image is the product $(a-b)\cdot x$ and since $S(\mathfrak{lie}_{\{x,y\}})$ is a domain, this 
implies $a-b=0$ hence $a=b$. This implies that $(a,b,z)$ is the image of $(0,a)\in F_0A$ under $F_0A\to F_0B$. 

Assume $k\geq 0$ and that the complex $F_kA\to F_kB\to \mathbb Q\langle x,y\rangle$ is exact; let us prove the exactness of 
$F_{k+1}A\to F_{k+1}B\to \mathbb Q\langle x,y\rangle$.  
Let $(a,b,z)\in\mathrm{ker}(F_{k+1}B\to \mathbb Q\langle x,y\rangle)$. 
Then $(a,b,z)\in F_{k+2}\mathbb Q\langle x,y\rangle_+^{\oplus2}\oplus\mathfrak{lie}_{\{x,y\}}$ and $ax-xb=z$. 
This equality implies that the image of its left-hand side under the projection map
$F_{k+3}\mathbb Q\langle x,y\rangle\to\mathrm{gr}_{k+3}\mathbb Q\langle x,y\rangle=S^{k+3}(\mathfrak{lie}_{\{x,y\}})$ is zero. 
Since this image is the product $[a-b]_{k+2}\cdot x$ and since $S(\mathfrak{lie}_{\{x,y\}})$ is a domain, this 
implies $[a-b]_{k+2}=0$, therefore $e\in F_{k+1}\mathbb Q\langle x,y\rangle_+$, where $e:=b-a$.
Moreover, $ax-xb=z$ implies the equality $[a,x]-x\cdot e=z$. 
This equality implies that the image of its left-hand side under the projection map
$F_{k+2}\mathbb Q\langle x,y\rangle\to\mathrm{gr}_{k+2}\mathbb Q\langle x,y\rangle=S^{k+2}(\mathfrak{lie}_{\{x,y\}})$ is zero. 
This image is $\{[a]_{k+2},x\}-x\cdot [e]_{k+1}$, therefore $\{[a]_{k+2},x\}=x\cdot [e]_{k+1}$. Then Lem. \ref{lem:poisson:1703}(b)
implies the existence of $\gamma\in S^{k+1}(\mathfrak{lie}_{\{x,y\}})$, such that 
$[a]_{k+2}=\gamma\cdot x$ and $[e]_{k+1}=\{\gamma,x\}$. Let $c\in F_{k+1}\mathbb Q\langle x,y\rangle$ be such that 
$[c]_{k+1}=\gamma$, and define $(a',b'):=(a-xc,b-cx)$. Then 
$(a',b')\in F_{k+1}\mathbb Q\langle x,y\rangle_+^{\oplus2}$, which implies 
$(a',b',z)\in F_kB$. The equality $ax-xb=z$ implies $a'x-xb'=z$, therefore
$(a',b',z)\in \mathrm{ker}(F_kB\to C)$. By the induction assumption, there exists 
$(c',u)\in F_kA$, such that $(a',b',z)=\mathrm{im}((c',u)\in F_kA\to F_kB)$. 
Then $(c+c',u)\in F_{k+1}A$, and $(a,b,z)=\mathrm{im}((c+c',u)\in F_{k+1}A\to F_{k+1}B)$.
This implies the exactness of $F_{k+1}A\to F_{k+1}B\to \mathbb Q\langle x,y\rangle$.
\end{proof}

\begin{prop}\label{lemma:A:1403} (see \cite{Sch1}, prop. 2.2)
Let $n>0$, then the sequence of maps from Lem. \ref{lemma:B:1403} induces an exact complex
$\mathbb Q\langle x,y\rangle[n-1]\oplus \mathfrak{lie}_{\{x,y\}}[n]\to
\mathbb Q\langle x,y\rangle[n]^{\oplus2}\oplus\mathfrak{lie}_{\{x,y\}}[n+1]\to \mathbb Q\langle x,y\rangle$.  
\end{prop}

\begin{proof}
The complex from Lem. \ref{lemma:B:1403} is graded. Its exactness implies that all its graded components are exact, 
which implies the statement. 
\end{proof}

\section{Relations in free algebras}\label{appC}

In this section, we fix a commutative $\mathbb Q$-algebra $\mathbf k$. 
Denote by $\mathbf k\langle x,y\rangle$ the free associative $\mathbf k$-algebra with generators $x,y$,  
where $x,y$ have degree 1, and by $\mathbf k\langle\langle x,y\rangle\rangle$ its degree completion; 
these algebras are equipped with Hopf algebra structures (in the completed sense in the second case) 
for which $x,y$ are primitive. Then
$\mathcal P\mathbf k\langle x,y\rangle=\mathfrak{lie}_{\{x,y\}}$ and 
$\mathcal P\mathbf k\langle\langle x,y\rangle\rangle=\mathfrak{lie}_{\{x,y\}}^\wedge$, 
moreover $\mathcal G\mathbf k\langle\langle x,y\rangle\rangle=\mathrm{exp}(\mathfrak{lie}_{\{x,y\}}^\wedge)$. 
For $n\geq0$, denote by $\mathbf k\langle x,y\rangle[n]$, $\mathfrak{lie}_{\{x,y\}}[n]$ the degree $n$ parts of 
$\mathbf k\langle x,y\rangle$, $\mathfrak{lie}_{\{x,y\}}$ and by $\mathbf k\langle\langle x,y\rangle\rangle[\geq n]$, 
$\mathfrak{lie}_{\{x,y\}}^\wedge[\geq n]$ the degree $\geq n$ parts of 
$\mathbf k\langle\langle x,y\rangle\rangle$, $\mathfrak{lie}_{\{x,y\}}^\wedge$.  

\begin{defn}
(a) Let $X:=\{(a,b)\in (\mathbf k\langle\langle x,y\rangle\rangle^\times)^2|a\cdot x\cdot b^{-1}\in x+
\mathfrak{lie}_{\{x,y\}}^\wedge[\geq2]\}$. 

(b) Set $X^{(0)}:=X$ and for $n\geq 1$, define $X^{(n)}\subset X$ as the subset of pairs 
$(a,b)$ such that $a,b\in 1+\mathbf k\langle\langle x,y\rangle\rangle[\geq n]$.

(c) $\mathbb G$ is the group defined by 
$$
\mathbb G:=\mathcal G(\mathbf k\langle\langle x,y\rangle\rangle)\times 
(\mathbf k^\times\times\mathbf k\langle\langle x,y\rangle\rangle)^{\mathrm{op}}, 
$$
where $\mathbf k^\times\times\mathbf k\langle\langle x,y\rangle\rangle$ is equipped 
with the product $(\gamma,c)\cdot(\gamma',c'):=(\gamma\gamma',\gamma c'+c\gamma'+cxc')$. 

(d) Set $F^0\mathbb G:=\mathbb G$ and for $n\geq 1$, set $F^n\mathcal G(\mathbf k\langle\langle x,y\rangle\rangle)
:=\mathcal G(\mathbf k\langle\langle x,y\rangle\rangle)\cap (1+\mathbf k\langle\langle x,y\rangle\rangle[\geq n])$
and $F^n(\mathbf k^\times\times\mathbf k\langle\langle x,y\rangle\rangle):=\{1\}\times 
\mathbf k\langle\langle x,y\rangle\rangle[\geq n-1]$ and 
$$
F^n\mathbb G:=
F^n\mathcal G(\mathbf k\langle\langle x,y\rangle\rangle)\times F^n(\mathbf k^\times\times\mathbf k\langle\langle x,y\rangle\rangle). 
$$
\end{defn}

\begin{lem}\label{lem:ind:1803}
(a) The group $\mathbb G$ acts on the set $X$ by 
$$
(h,(\gamma,c))\cdot (a,b):=(h\cdot a\cdot (\gamma+xc),h\cdot b\cdot (\gamma+cx)). 
$$

(b) For any $n\geq 0$ and $x\in X^{(n)}$, there exists $g\in F^n\mathbb G$ such that $g\cdot x\in X^{(n+1)}$. 
\end{lem}

\begin{proof}
(a) is immediate. Let us prove (b). Let $x=(a,b)\in X$ and $g:=(1,(\epsilon(a)^{-1},0))\in\mathbb G$, where 
$\epsilon : \mathbf k\langle\langle x,y\rangle\rangle^\times\to\mathbf k^\times$ is the augmentation map. The relation 
$a\cdot x\cdot b^{-1}\in x+\mathfrak{lie}_{\{x,y\}}^\wedge[\geq2]$ implies $\epsilon(a)=\epsilon(b)$, which implies 
the second equality in $g\cdot x=(a/\epsilon(a),b/\epsilon(a))=(a/\epsilon(a),b/\epsilon(b))\in X^{(1)}$. This proves (b) when $n=0$. 
Let now $n\geq 1$ and $x=(a,b)\in X^{(n)}$. Let $a[n],b[n]$ be the components of $a,b$ in
$\mathbf k\langle x,y\rangle[n]$, then $a\equiv 1+a[n]$, $b\equiv 1+b[n]$ modulo  
$\mathbf k\langle\langle x,y\rangle\rangle[\geq n+1]$. Then $a\cdot x\cdot b^{-1}\equiv x+a[n]x-xb[n]$
modulo  $\mathbf k\langle\langle x,y\rangle\rangle[\geq n+2]$. It follows that $a[n]x-xb[n]\in 
\mathfrak{lie}_{\{x,y\}}[n+1]$. Then Prop. \ref{lemma:A:1403} implies that 
the sequence of maps from Lem. \ref{lemma:B:1403} induces an exact complex
$\mathbf k\langle x,y\rangle[n-1]\oplus \mathfrak{lie}_{\{x,y\}}[n]\to
\mathbf k\langle x,y\rangle[n]^{\oplus2}\oplus\mathfrak{lie}_{\{x,y\}}[n+1]\to \mathbf k\langle x,y\rangle$.
This implies 
the existence of a pair
$(c,u)$, where $c\in \mathbf k\langle x,y\rangle[n-1]$ and $u\in \mathfrak{lie}_{\{x,y\}}[n-1]$ such that 
\begin{equation}\label{rels:1803}
a[n]=xc+u,\quad  b[n]=cx+u. 
\end{equation}
Then $g:=(\mathrm{exp}(-u),(1,-c))\in F^n\mathbb G$, and the relations \eqref{rels:1803} imply $g\cdot x\in X^{(n+1)}$. 
This ends the proof of (b). 
\end{proof}

\begin{prop}\label{prop:Z}
Let $(a,b,z)$ be a triple where $a,b\in \mathbf k\langle\langle x,y\rangle\rangle^\times$, 
$z\in\mathfrak{lie}_{\{x,y\}}^\wedge[\geq2]$, and $a\cdot x=(x+z)\cdot b$. 

Then for some $(\gamma,c) \in \mathbf k^\times\times\mathbf k\langle\langle x,y\rangle\rangle$ and 
$h \in \mathcal G(\mathbf k\langle\langle x,y\rangle\rangle)$, 
one has $a=h\cdot (\gamma+xc)$, $b=h\cdot (\gamma+cx)$, $x+z=h\cdot x\cdot h^{-1}$. 
\end{prop}

\begin{proof}
The assumption on $(a,b,z)$ implies $(a,b)\in X$. Define inductively $x^{(n)}\in X^{(n)}$, 
$g^{(n)}\in F^n\mathbb G$ for $n\geq0$, by $x^{(0)}:=(a,b)$, and by the conditions that for any $n\geq 0$, $g^{(n)}\in F^n\mathbb G$
is such that $g^{(n)}\cdot x^{(n)}\in X^{(n+1)}$ (see Lem. \ref{lem:ind:1803}), 
and that $x^{(n+1)}:=g^{(n)}\cdot x^{(n)}$. The infinite product $\cdots \cdot g^{(1)}\cdot g^{(0)}$ defines an element
$g\in \mathbb G$. It follows from the sequence of inclusions $X=X^{(0)}\supset X^{(1)}\supset \cdots$ and from 
$\cap_{n\geq0}X^{(n)}=\{(1,1)\}$ that $g\cdot (a,b)=(1,1)$, therefore $(a,b)=g^{-1}\cdot (1,1)$. This implies the claimed conclusion, 
with $h\in\mathcal G(\mathbf k\langle\langle x,y\rangle\rangle)$
and $(\gamma,c)\in\mathbf k^\times\times\mathbf k\langle\langle x,y\rangle\rangle$ the elements such that 
$g^{-1}=(h,(\gamma,c))$. 
\end{proof}


\begin{thebibliography}{BGFr}

\bibitem[AT]{AT}
A. Alekseev, C. Torossian, {\it The Kashiwara-Vergne conjecture and Drinfeld’s associators}, 
Ann. of Math. (2) 175 (2012), no. 2, 415--463. 

\bibitem[AET]{AET}
A. Alekseev, B. Enriquez, C. Torossian,
{\it Drinfeld associators, Braid groups and explicit solutions of the Kashiwara–Vergne equations,}
Publications Mathématiques de l'IHES, Volume 112 (2010), pp. 143-189.


\bibitem[De]{De} P. Deligne, {\it Le groupe fondamental de la droite projective moins trois points,} Galois groups
over $\mathbb Q$ (Berkeley, CA, 1987), 79–297, Math. S. Res. Inst. Publ., 16, Springer, New York-Berlin,
1989.

\bibitem[DeG]{DeG} P. Deligne, and A. Goncharov, 
{\it Groupes fondamentaux motiviques de Tate mixte,} Ann. Sci. Ecole Norm. Sup. (4) 38 (2005), no. 1, 1–56.

\bibitem[DeT]{DT} P. Deligne, T. Terasoma, {\it Harmonic shuffle relation for associators,} preprint (2005).

\bibitem[Dr]{Dr} V. Drinfeld, {\it On quasitriangular quasi-Hopf algebras and on a group that is closely connected with 
$\mathrm{Gal}(\overline{\mathbb Q}/\mathbb Q)$}, Leningrad Math. J. 2 (1991), no. 4, 829--860.  

\bibitem[Ec1]{Ec1}
J. Ecalle, {\it A tale of three structures: the arithmetics of multizetas, the analysis of singularities,
the Lie algebra ARI,} in Differential Equations and the Stokes Phenomenon, Braaksma, Immink,
van der Put, Top eds., World Scientific 2002, 89–146.

\bibitem[Ec2]{Ec2} J. Ecalle, {\it The flexion structure and dimorphy: flexion units, singulators, generators, and
the enumeration of multizeta irreducibles,} With computational assistance from S. Carr. CRM
Series, 12, Asymptotics in dynamics, geometry and PDEs; generalized Borel summation. Vol. II,
27--211, Edizioni della Normale, Pisa, 2011.

\bibitem[EF0]{EF0} B. Enriquez, H. Furusho, {\it A stabilizer interpretation of double shuffle Lie algebras.} 
Int. Math. Res. Not. 2018, no. 22, 6870--6907. 

\bibitem[EF1]{EF1} B. Enriquez, H. Furusho, {\it The Betti side of the double shuffle theory. I. The harmonic coproduct.} 
Selecta Math. (N.S.) 27 (2021), no. 5, Paper No. 79, 106 pp.

\bibitem[EF2]{EF2} B. Enriquez, H. Furusho,
{\it The Betti side of the double shuffle theory. II. Double shuffle relations for associators.}
Selecta Math. (N.S.) 29 (2023), no. 1, Paper No. 3, 28 pp.


\bibitem[EF3]{EF3}  B. Enriquez, H. Furusho,
{\it The Betti side of the double shuffle theory. III. Bitorsor structures}, 
Selecta Math. (N.S.) 29 (2023), no. 2, Paper No. 27, 37 pp.

\bibitem[EF4]{EF4}  B. Enriquez, H. Furusho, 
{\it  The stabilizer bitorsors of the module and algebra harmonic coproducts are equal}, 
Adv. Math. 463 (2025), Paper No. 110128, 29 pp.

\bibitem[F1]{Fu:assandDS}
H. Furusho, {\it Double shuffle relation for associators,} Annals of Mathematics, Vol. 174 (2011), No. 1, 341-360.

\bibitem[F2]{Fu:confl}
H. Furusho, {\it The pentagon equation and the confluence relations,} Amer. J. Math. Vol 144, No 4, (2022) 873-894. 


\bibitem[FK]{FK} 
H.~Furusho and N.~Komiyama,
{\it Notes on Kashiwara-Vergne and double shuffle Lie algebras}, 
"Low Dimensional Topology and Number Theory", Springer Proc. Math. Stat., 456, 
Springer, Singapore, 2025, 63-80.

\bibitem[G]{Gi} J. Giraud, Cohomologie non abélienne. Die Grundlehren der mathematischen 
Wissenschaften, Band 179. Springer-Verlag, Berlin-New York, 1971.

\bibitem[HS]{HS}  M. Hirose, N. Sato, {\it Iterated integrals on $\mathbb P^1 \smallsetminus\{0,1,\infty,z\}$ 
and a class of relations among multiple zeta values,} Adv. Math. 348 (2019), 163–182.

\bibitem[IKZ]{IKZ}
K. Ihara, M. Kaneko, and D. Zagier, Derivation and double shuffle relations for multiple zeta
values, Compos. Math. 142 (2006), no. 2, 307–338.

\bibitem[Ka]{Kaw} H. Kawamura, {\it Ecalle's senary relation and dimorphic structures,}
preprint {\tt arXiv:2509.21252}. 

\bibitem[Ku]{Ku} Y. Kuno, {\it Emergent version of Drinfeld's associator equations,}
preprint {\tt arXiv:2504.02549}. 

\bibitem[LM]{LM} T.Q.T. Le and J. Murakami, {\it The universal Vassiliev-Kontsevich invariant for framed oriented
links,} Compositio Math. 102 (1996), no. 1, 41–64.

\bibitem[McL]{McL} S. Mac Lane, Categories for the Working Mathematician. Grad. Texts in Math., 5, 
Springer-Verlag, New York, 1998.

\bibitem[Me]{Mei} E. Meinrenken, {\it Clifford algebras and Lie theory,} 
Ergeb. Math. Grenzgeb. (3), 58, 
Springer, Heidelberg, 2013. 

\bibitem[R]{R} G. Racinet, {\it Doubles m\'elanges des polylogarithmes multiples aux racines de l'unit\'e.}
Publ. Math. Inst. Hautes \'Etudes Sci., No. 95 (2002), 185--231. 

\bibitem[S1]{Sch1} L. Schneps, {\it Double shuffle and Kashiwara-Vergne Lie algebras,} 
 J. Algebra  367 (2012), 54--74.

\bibitem[S2]{Sch2} L. Schneps, {\it The double shuffle Lie algebra injects into the Kashiwara-Vergne Lie algebra,} 
preprint {\tt arXiv:2504.14293}.

\bibitem[T]{T} T. Terasoma, {\it Geometry of multiple zeta values,}
International Congress of Mathematicians. Vol. II, 627–635.
European Mathematical Society (EMS), Z\"{u}rich, 2006. 

\end{thebibliography}
\end{document}